\documentclass[reqno,11pt]{amsart}
\usepackage{amsfonts}
\usepackage{graphicx}
\usepackage{amsfonts,amsmath, amssymb,bm}
\usepackage{bbm,mathrsfs,mathscinet}
\usepackage{cancel}
\usepackage{marginnote}
\usepackage{color}
\usepackage{textcomp}
\usepackage{tcolorbox}
\allowdisplaybreaks[4]

\numberwithin{equation}{section}

\newcommand{\R}{\mathbb{R}}

\newtheorem{theorem}{Theorem}[section]
\newtheorem{corollary}[theorem]{Corollary}
\newtheorem{lemma}[theorem]{Lemma}
\newtheorem{proposition}[theorem]{Proposition}
\newtheorem{remark}[theorem]{Remark}
\newtheorem{definition}[theorem]{Definition}

\def\v{\varepsilon}
\def\e{\epsilon}
\def\x{\xi}

\def\a{\alpha}
\def\b{\beta}
\def\g{\gamma}
\def\d{\delta}

\def\f{\frac}
\def\pa{\partial}
\def\dis{\displaystyle}
\def\la{\langle}
\def\ra{\rangle}
\def\mdiv{\mbox{\rm div}}

\begin{document}
	
\title[Prandtl Expansion for Steady  CNSF ]{Validity of Prandtl Expansion for Steady Compressible Navier-Stokes-Fourier Flows}

\author[Y. Guo]{Yan Guo}
\address[Y. Guo]{The Division of Applied Mathematics, Brown University, Providence, USA.}
\email{yan\_guo@brown.edu}

\author[Y. Wang]{Yong Wang}
\address[Y. Wang]{Academy of Mathematics and Systems Science, Chinese Academy of Sciences, Beijing 100190, China }
\email{yongwang@amss.ac.cn}

\begin{abstract}
	Assume no-slip boundary condition for the velocity field and either insulated or Dirichlet boundary conditions for the temperature field in a steady compressible fluid. In the inviscid limit $\varepsilon \rightarrow 0$, we develop a mathematical framework for the uniform-in-$\varepsilon$ remainder estimate for the linear steady compressible Navier-Stokes-Fourier equations around a Prandtl layer profile $(\rho_s, p_s, u_s, v_s, T_s)$ with both velocity and  thermal layers, which leads to the validity of the Prandtl layer expansion. We introduce the \textit{generalized stream function} $\phi$ for the remainder to solve the continuity equation, and solve the momentum equations $A_1 \equiv 0$, $A_2 \equiv 0$ by studying their div-curl system for $(A_1, A_2)$ with paired boundary conditions:
	\begin{align}\label{pair1}
		\begin{split}
			A_{1}|_{x=0} &  =\phi|_{x=0}=\phi_{xx}|_{x=0}=0,\\
			A_{2}|_{x=L} &  =\phi_{x}|_{x=L}=\phi_{xxx}|_{x=L}=0,  
		\end{split}
	\end{align}
	which exhibit an exact structural cancellation of the highest-order $\phi$ derivatives in $A_1$ and $A_2$. To solve the temperature $T$ equation, we introduce the paired \textit{pseudo entropy and its corresponding quotient} $q$:
	\begin{equation}\label{s}
		S:=T-\frac{1}{2\rho_{s}} p-\frac{1}{p_s}T_{sy}q\quad \text{and}\quad q:=\frac{T_{s}\phi}{u_{s}},
	\end{equation}
	which avoids an $\varepsilon$-loss through another exact cancellation in the bulk. To avoid degeneracy near $y = 0$ for the pseudo entropy $S$ under the Dirichlet boundary condition, we impose a boundary condition at $x = L$ to ensure structural positivity:
	\begin{equation}\label{0.3}
		\{2\rho_s u_s S_{x} - \kappa S_{yy}\}\big|_{x=L}=
		\f{\kappa }{p_s}T_{sy}q_{yy} + 2 \big(\f{\kappa}{p_s}T_{sy}\big)_{y}q_y + \big(\f{\kappa}{2\rho_s}\big)_{yy}p + \mathfrak{M}_0 (\mbox{nonlinear}).
	\end{equation}
	We construct $\phi$ from the curl of $(A_1, A_2)$ via an improved $H^4$ theory by Guo-Iyer \cite{Guo-Iyer-2024} and a new $H_{00}^{1/2}$ estimate, and solve for $p$ from the divergence of $(A_1, A_2)$ by imposing the viscous inflow boundary condition
	\begin{equation}\label{pressure}
		(\mu+\lambda) u_s \Delta_{\v}p\big|_{x=0}= -\sigma p_s \mathbf{d}_{11} B,
	\end{equation}
	with a {\it good unknown} $B$ to avoid the corner singularity and obtain structural positivity thanks to the exact cancellations from the paired boundary conditions.
\end{abstract}

\date{\today}
\maketitle

\tableofcontents
	
\thispagestyle{empty}

\section{Introduction}

\subsection{Background} Studying steady flows is crucial in fluid dynamics for  designing efficient systems (e.g., pipes, turbines, aircraft). Boundary layer study is fundamental in flying vessels 
surrounded by a  compressible fluid, for which the classical boundary layer theory for incompressible fluid becomes inadequate. 
In this paper, we prove the validity of boundary layer expansion for a steady compressible fluid. 
Despite its importance both from mathematical and physical standpoints, to
our knowledge, there are very limited mathematical results to solve the steady
compressible flows with both momentum and energy equations even for a finite Reynolds' number ($Re$)  in the presence of mixed non-slip and in-flow boundary condition.  
The main contribution of this paper is to provide a systematic mathematical framework to solve the  steady compressible Navier-Stokes-Fourier equations uniformly in $Re\gg1$. Both the velocity and thermal boundary layers are included. 

We consider the steady compressible Navier-Stokes-Fourier equations (CNSF) on the two dimensional domain $(x,Y)\in (0,L)\times \R_+=:\Omega$
{\small
	\begin{align}\label{1.1}
		\begin{cases}
			\dis \mbox{\bf div}(\rho^{NS} {U}^{NS})=0,\\[2mm]
			\dis  \rho^{NS} ( {U}^{NS}\cdot \mathbf{\nabla}) {U}^{NS} +   {\bf\nabla} P^{NS}   = \mu \, \v {\bf\Delta} {U}^{NS} + \lambda\, \v\, {\bf\nabla}\mbox{\bf div}{U}^{NS},\\
			\dis c_v \rho^{NS} (U^{NS}\cdot \nabla) T^{NS} + P^{NS} \mbox{\bf div}{U}^{NS} = \kappa \v {\bf\Delta} {T}^{NS} + 2\mu \v \big|\mathbb{S}(U^{NS})\big|^2 + (\lambda-\mu) \v  \big|\mbox{\bf div}{U}^{NS}\big|^2,
		\end{cases}
\end{align}}
where $\rho^{NS}$, ${U}^{NS}:=(u^{NS}, v^{NS})^t$, $T^{NS}$  are the density, velocity and temperature, respectively. Throughout this paper, we assume the ideal gas law, i.e., the pressure satisfies
\begin{align}\label{1.1-0}
	P^{NS}=\rho^{NS}T^{NS},
\end{align}
where the  specific gas constant has been taken as unity. Here $\v>0$ is the viscosity which is proportional to $1/{Re}$,  $\kappa>0$ is the \textit{normalized} heat conductive coefficient, and $\mu$ and $\lambda$ are the \textit{normalized} viscous coefficients. We assume that $\mu>0$ and $\lambda>0$ in this paper. Here $c_v>0$ is a positive constant, and we normalize it to be $1$ for simplicity.  We have denoted $\mbox{\bf div}=\pa_x + \pa_Y$, ${\bf\Delta}=\pa_{xx}+ \pa_{YY}$, and $\mathbb{S}(U)=\frac12 \big(\nabla U + \nabla U^T\big).$

The system \eqref{1.1} is supplemented by the no-slip boundary condition (BC) for velocity
\begin{align}\label{1.2}
	\big(u^{NS},\, v^{NS}\big)|_{Y=0}=(0,\,0),
\end{align}
and the following boundary condition for the temperature
\begin{align}\label{1.3-0}
	&\pa_YT^{NS}|_{Y=0}=0\,\, \big(\mbox{or}\,\,\, T^{NS}|_{Y=0}=T_b(x)\big).
\end{align}
This means we consider both Neumann (NT) and Dirichlet (DT) boundary conditions for the temperature.

\smallskip 

\noindent{\it The Euler background with mismatch:}  
Taking $\v=0$ in \eqref{1.1}, we obtain the compressible Euler equations
\begin{align}\label{Euler}
	\begin{cases}
		\dis \mbox{\bf div}(\rho_e U_e)=0,\\
		\dis \rho_e (U_e\cdot\nabla)U_e + \nabla P_e=0,\\
		\dis \rho_e (U_e\cdot \nabla) T_e + P_e \mbox{\bf div}{U}_e =0,\\
		\dis v_e|_{Y=0}=0.
	\end{cases}
\end{align}
We consider a far‑field Euler shear flow of the form
\begin{align}\label{1.3-1}
	(\rho_e^0,u_e^0,v_e^0,T_e^0)= (\rho_e^0(Y),u_e^0(Y),0,T_e^0(Y)),
\end{align} 
where the functions $\rho_e^0$, $u_e^0$, and $T_e^0$ are sufficiently smooth  satisfying
\begin{align}\label{1.4-0}
	p_e^0 = \rho_e^0(Y) T_e^0(Y) = 1.
\end{align}
Clearly, any such shear flow is a solution of the compressible Euler system \eqref{Euler}. For later use, we  assume 
\begin{align}\label{1.3}
	\begin{cases} 
		\dis 0<u_e^0(0)\leq u_e^0(Y) \leq a,   \\
		\dis  a:=\lim_{Y\rightarrow\infty} u_e^0(Y)>0, \\
		\dis |u_e^0(Y)-a|\leq a  e^{-Y},  \\
		\dis |\pa_{Y}^k u_e^0(Y)|\lesssim a e^{-Y},
	\end{cases}
	\& \qquad 
	\begin{cases} 
		\dis\f12\leq  T_e^0(Y) \leq 2, \\
		\dis  1:=\lim_{Y\to +\infty}T^0_e(Y)>0,\\
		\dis|T_e^0(Y)-1| \lesssim 
		a^2 e^{-Y},\\
		\dis |\pa_{Y}^k T_e^0(Y)|\lesssim a^2 e^{-Y}.
	\end{cases}
\end{align}
It is well known that, generically, such an Euler flow exhibits significant mismatches with the physical boundary conditions \eqref{1.2}-\eqref{1.3-0}:
\begin{align}\label{0.9}
	u_e^0(0)\neq 0\equiv u^{NS}(x,0)\quad\mbox{and}\quad  T_e^0(0)\neq T^{NS}(x,0).
\end{align}

\smallskip

\noindent{\it Incompressible Prandtl boundary layer:} To address the velocity mismatch, Prandtl introduced the boundary layer theory in 1904 for the simplified incompressible Navier-Stokes equations (formally, the Mach number $\mathscr{M}_a=0$ and $\rho^{NS}=1$ in $\eqref{1.1}_{1,2}$):
\begin{align}\label{0.9-0}
	\begin{cases}
		\mbox{\bf div} U^{NS}=0,\\
		({U}^{NS}\cdot \mathbf{\nabla}) {U}^{NS} +   {\bf\nabla} p^{NS}   = \mu \, \v {\bf\Delta} {U}^{NS}.
	\end{cases}
\end{align}
He derived the celebrated Prandtl equations:
\begin{align}\label{0.10}
	\begin{cases}
		\pa_x\bar{u}_p + \pa_y \bar{v}_p=0,\\
		\bar{u}_p \pa_x\bar{u}_p + \bar{v}_p \pa_y \bar{u}_p  + \pa_xp_{e}(x,0) =\mu \pa_{yy} \bar{u}_p,
	\end{cases}
\end{align}
where the stretched variable $y:=\f{Y}{\sqrt{\v}}$. Such a Prandtl boundary layer characterizes a sharp transition of the velocity field within a thin layer of $\sqrt{\v}$ near the boundary $\{Y=0\}$. The Prandtl theory states that
the behavior near the boundary is described by the Prandtl boundary layer equations, while the region away from the boundary is approximated by the inviscid Euler equations, i.e.,
\begin{align}\label{0.11}
	\mbox{Incompressible NS}=\mbox{Euler} + \mbox{Prandtl layer} + O(\sqrt{\v}).
\end{align}

\noindent{\it Thermal boundary layer:}  
Clearly, the classical Prandtl layer fails to address the mismatch of the temperature in \eqref{0.9}, i.e., $T_e^0(0)\neq T^{NS}(x,0)$. Consequently, the classical Prandtl boundary layer theory for incompressible flows (formally, the Mach number $\mathscr{M}_a=0$ and $\rho^{NS}=1$ in $\eqref{1.1}_{1,2}$) with large Reynolds number
becomes insufficient for modeling compressible fluids with $\mathscr{M}_a\neq0$ that exhibit sharp temperature variations within a thin layer adjacent to a solid surface (the thermal boundary layer) \cite{Schlichting}, at which the temperature
satisfies naturally either Neumann (isothermal body) or Dirichlet boundary
conditions. It is well-known that such a thermal layer interacts strongly and nonlinearly with the classical velocity Prandtl layer through CNSF,  playing a crucial role in many important applications in Aeronautics  and Astronautics.


\smallskip

Due to the appearance of Prandtl boundary layer, we define 
\begin{align}
	\begin{split}
		&y:=\f{Y}{\sqrt{\v}},\quad \nabla_{\v}=(\sqrt{\v}\pa_x, \pa_y),\quad   \Delta_{\v}=\v \pa_{xx} + \pa_{yy} = \mbox{div}_{\v} \nabla_{\v},\\
		&u^\v(x,y):=u^{NS}(x,Y),\quad v^{\v}(x,y):= \f{1}{\sqrt{\v}}v^{NS}(x,Y),\\
		&\rho^{\v}(x,y):=\rho^{NS}(x,Y), \quad T^{\v}(x,y)=T^{NS}(x,Y),\quad p^{\v}(x,y)=P^{\v}(x,Y).
	\end{split}
\end{align} 
Then the steady CNSF \eqref{1.1} is rewritten as
\begin{align}\label{1.4}
	\begin{cases}
		\dis \pa_x(\rho^{\v}u^\v) + \pa_y(\rho^{\v} v^{\v})=0,\\[2mm]
		\dis \rho^{\v} [u^{\v}\pa_x u^{\v} + v^{\v} \pa_y u^{\v}] + \pa_x p^\v
		=\mu \Delta_{\v} u^{\v} + \lambda\,\v \pa_x\mdiv{U}^\v,\\[2mm]
		\dis \v \rho^{\v} [u^{\v}\pa_x v^{\v} + v^{\v} \pa_y v^{\v}] + \pa_y p^\v
		=\mu \v \Delta_{\v} v^{\v} + \lambda\,\v \pa_y\mdiv{U}^\v,\\[2mm]
		\dis 2\rho^{\v} [u^{\v}\pa_x T^{\v} + v^{\v} \pa_y T^{\v}] -  ( U^{\v}\cdot\nabla) p^{\v} = \kappa \Delta_{\v}T^{\v} + 2\mu   |\mathbb{S}^{\v}(U^{\v})|^2 + (\lambda-\mu) \v |\mbox{div} U^{\v}|^2,
	\end{cases}
\end{align}
where we have denoted ${U}^\v= (u^\v,v^\v)^{t}$, $\mdiv{U}^\v=\pa_x u^\v+\pa_y v^\v$, and 
\begin{align*}
	\mathbb{S}^{\v}(U^{\v})=  \left( 
	\begin{array}{cc}
		\sqrt{\v} \pa_x u^{\v}, & \f12 (\pa_{y}u^{\v} + \sqrt{\v} \sqrt{\v}\pa_x v^{\v})\\
		\f12 (\pa_{y}u^{\v} + \sqrt{\v} \sqrt{\v}\pa_x v^{\v}), & \sqrt{\v}\pa_y v^{\v}
	\end{array}
	\right),\quad \mbox{\bf div}U^{NS}\equiv \mbox{div} U^{\v}.
\end{align*}
The boundary conditions of velocity and temperature at $\{y=0\}$ become:
\begin{align}\label{1.4-1}
	\begin{cases}
		[u^{\v},\, v^{\v}]|_{y=0}=[0,\,0],\\[1.5mm]
		\pa_yT^{\v}|_{y=0}=0\,\,\,\,\, \mbox{for NT},\\
		\mbox{or}\,\, T^{\v}|_{y=0}=T_b 
		\,\,\,\,\, \mbox{for DT}.
	\end{cases}
\end{align}

We aim to construct the solutions of steady CNSF \eqref{1.4}-\eqref{1.4-1} in the following form ($n, N_0>1$ are suitably large number)
\begin{align}\label{1.5}
	\begin{split}
		\rho^\v&=\rho_e^0 + \rho^0_p + \sum_{i=1}^n \sqrt{\v}^i \big[\rho_e^i(x,Y) + \rho_p^i(x,y)\big]  + \v^{N_0}\rho^{(\v)}=:\rho_s + \v^{N_0}\rho^{(\v)},\\
		u^\v&= u_e^0 + u_p^0+ \sum_{i=1}^n \sqrt{\v}^i \big[u_e^i(x,Y) + u_p^i(x,y)\big] + \v^{N_0} u^{(\v)}=:u_s + \v^{N_0} u^{(\v)},\\
		v^\v&=v_e^1+ v_p^0+ \sum_{i=1}^{n-1} \sqrt{\v}^i \big[v_e^{i+1}(x,Y) + v_i(x,y)\big]  + \sqrt{\v}^n v_p^n + \v^{N_0} v^{(\v)}
		=:v_s + \v^{N_0} v^{(\v)},\\
		T^\v&=T_e^{0}+ T^{0}_p+ \sum_{i=1}^n \sqrt{\v}^i \big[T_e^i(x,Y) + T_p^i(x,y)\big]  + \v^{N_0}T^{(\v)}=:T_s + \v^{N_0}T^{(\v)}.
	\end{split}
\end{align}
Here $(\rho^i_e, u^i_e, v^i_e, T^i_e), i=1,\cdots, n$ are the high-order corrections to the Euler equations; $(\rho^0_p, u^0_p, v^0_p, T^0_p)$ is the leading-order Prandtl boundary layer, while $(\rho^i_p, u^i_p, v^i_p, T^i_p), i=1,\cdots, n$ are its high-order Prandtl corrections. The remainder $(\rho^{(\v)},u^{(\v)},v^{(\v)}, T^{(\v)})$ now satisfy a nonlinear CNSF system.

For the leading order Prandtl boundary layer, we denote 
\begin{align}\label{3.5-0}
	\begin{cases}
		\bar{\rho}^0_p:= \rho_e^0(0) + \rho^0_p(x,y),\\
		\bar{u}^0_p:= u_e^0(0) + u^0_p(x,y),\\
		\bar{v}^0_p:= v_e^1(x,0) + v^0_p(x,y),\\
		\bar{T}^0_p:= T_e^0(0) + T^0_p(x,y),
	\end{cases}
\end{align}
which satisfy the following steady compressible Prandtl equations
\begin{align}\label{3.6-0}
	\begin{cases}
		\pa_y\big(\bar{\rho}^0_p \bar{T}^0_p\big)=0,\\
		(\bar{\rho}^0_p \bar{u}^0_p)_x + (\bar{\rho}^0_p \bar{v}^0_p)_y =0,\\
		\bar{\rho}^0_p[\bar{u}^0_p \, \pa_x + \bar{v}^0_p \, \pa_y] \bar{u}^0_p + (\bar{\rho}^0_p \bar{T}^0_p)_x = \mu \pa_{yy} \bar{u}^0_p, \\
		\bar{\rho}^0_p[\bar{u}^0_p \, \pa_x + \bar{v}^0_p \, \pa_y] \bar{T}^0_p + \bar{\rho}^0_p \bar{T}^0_p [\pa_x \bar{u}^0_p + \pa_y \bar{v}^0_p] = \kappa \pa_{yy}\bar{T}^0_{p} + \mu |\pa_y \bar{u}^0_p|^2.
	\end{cases}
\end{align}
Guo-Wang  \cite{Guo-Wang} have already constructed the local-in-$x$ solutions for the Prandtl boundary layer $(\rho^0_p, u^0_p, v^0_p, T^0_p)$ and higher order corrections $(\rho^i_e, u^i_e, v^i_e, T^i_e)$ \& $(\rho^i_p, u^i_p, v^i_p, T^i_p), i=1,\cdots, n$, and proved that  $(\rho_s, u_s, v_s,T_s)$ satisfy the following approximate steady CNSF
\begin{align}\label{1.7}
	\begin{cases}
		\dis \pa_x(\rho_s u_s) + \pa_y(\rho_s v_s)=\pa_xA_{s1}+\pa_y A_{s2},\\
		\dis \rho_s [u_s\pa_x u_s + v_s \pa_y u_s] + \pa_x p_s
		=\mu  \Delta_{\v} u_s + \lambda \v \pa_x\mdiv{U_s}  + B_{s1},\\
		\dis \v \rho_s [u_s\pa_xv_s+ v_s \pa_y v_s] + \pa_y p_s
		=\mu \v \Delta_{\v} v_s + \lambda \v \pa_y\mdiv{U_s}  + \v B_{s2},\\
		\dis 2\rho_s [u_s\pa_x T_s + v_s \pa_y T_s] -(U_s\cdot \nabla)p_s
		=\kappa  \Delta_{\v} T_s + 2\mu |\mathbb{S}^{\v}(U_s)|^2 + (\lambda-\mu) \v  |\mbox{div}U_s|^2 + B_{s3},
	\end{cases}
\end{align}
with $p_s= \rho_s T_s$ and $U_s=(u_s, v_s)^{t}$. Also it holds that $(u_s, v_s)|_{y=0}=(0,0)$, and $T_{sy}|_{y=0}=0$ for NT (or $T_s|_{y=0}=T_b(x)$ for DT). Here $A_{s1},A_{s2}$ and $B_{s1},B_{s2}, B_{s3}$ are some given smooth functions satisfying
\begin{align}\label{1.7-1}
	\begin{split}
		&|\pa_x^i\pa_y^j(A_{s1},A_{s2},B_{s1},B_{s2},B_{s3})|\lesssim O(1) \v^{n}\big\{\la y\ra^{-\mathfrak{l}_0} + \la \sqrt{\v}y\ra^{-\mathfrak{l}_0}\big\},\quad \mbox{for}\,\,\, i+j\geq 0,\,\,  y\geq 0,\\
		& A_{s1}=O(\v^n) u_s, \quad A_{s2}=O(\v^n) u_s^2,\quad \mbox{for}\,\, y\in [0,1],
	\end{split}
\end{align}
where $\mathfrak{l}_0>1$ is a suitably large constant, see \cite{Guo-Wang} for more details. 
Due to $\eqref{1.7}_1$ and \eqref{1.7-1}, we also have 
\begin{align}\label{1.7-2}
	(u_{sx}+v_{sy})|_{y=0}=0 \,\,\, \Longrightarrow \,\,\,   v_{sy}|_{y=0}=0 \quad \Longrightarrow|v_s|\lesssim ay^2\,\,\, \mbox{for}\,\, y\in[0,1].
\end{align}
Also it follows from Remark 1.5 in \cite{Guo-Wang} that 
\begin{align}\label{1.7-3}
	\begin{split}
		&u_s\cong ay \,\,\mbox{for}\,\, y\in[0,1]\quad \& \quad 
		u_s\cong a>0\,\, \, \mbox{for}\,\, y\geq 1,\\
		& T_s \cong 1 \,\,\mbox{for}\,\, y\geq 0 (\mbox{for DT and NT}) \quad \& \quad  T_{sy}\cong |u_{sy}|^2\cong a^2\,\,\, \mbox{for}\,\, y\in[0,1] \,(\mbox{only for DT}),\\
		&|(\pa^l_{x}\pa^{j}_{y} \rho^{i}_{p}, \pa^l_{x}\pa^{j}_{y} T^{i}_{p},  \pa^l_{x}\pa^{j}_{y} u^{i}_{p}, \pa^{l}_x\pa^{j+1}_{y} v^{i}_{p})|\lesssim a (1+y)^{-\mathfrak{l}_0} \quad \mbox{for}\quad l+j\geq 0,\,\, i=0,1,\cdots, n,\\
		&|(\pa^l_{x}\pa^{j}_{Y} \rho^{i}_{e}, \pa^l_{x}\pa^{j}_{Y} T^{i}_{e}, \pa^l_{x}\pa^{j}_{Y} u^{i}_{e},\pa^{l}_x\pa^{j}_{Y} v^{i}_{e})|\lesssim a (1+\sqrt{\v}y)^{-\mathfrak{l}_0} \quad \mbox{for}\quad l+j\geq 0,\,\, i=1,\cdots, n.
	\end{split}
\end{align}
For later use, we define the background Mach number
\begin{align}
	&\dis \mathscr{M}_a:=\sup_{(x,y)\in\Omega} \Big\{\f{|U_s|}{\sqrt{2T_s}}\Big\} \cong a.
\end{align}

By \eqref{1.4-0}, we note  $p_{sy}=O(\sqrt{\v})$, so that 
we denote 
\begin{align}\label{1.19}
	p^{\v}:=p_s + \v^{N_0} p^{(\v)}.
\end{align}
Then, by \eqref{1.1-0}, it is clear to know that 
\begin{align}\label{1.19-1}
	\begin{split}
		\rho^{(\v)} &
		= \f{1}{T^{\v}} \big[p^{(\v)} - \rho_s T^{(\v)}\big] 
		=\underbrace{\f{1}{T_s} p^{(\v)} - \f{\rho_s}{T_s} T^{(\v)}}_{:=\bar{\rho}^{(\v)}} + \underbrace{(\f{1}{T^{\v}} - \f{1}{T_s}) \big[p^{(\v)} - \rho_s T^{(\v)}\big]}_{:=N_{\rho}=O(1)\v^{N_0}|(p^{(\v)}, T^{(\v)})|^2}.
	\end{split}
\end{align}
With the help of \eqref{1.5} and \eqref{1.7}, we can obtain the remainder equations of $\big(\rho^{(\v)},u^{(\v)},v^{(\v)}, T^{(\v)}\big)$:
\begin{align}\label{1.8-1}
	\begin{cases}
		\dis \mdiv\big[\rho_s U^{(\v)} +   U_s \rho^{(\v)} +  \v^{N_0} \rho^{(\v)} U^{(\v)} + \v^{-N_0} A_s\big]=0,\\[2mm]
		\dis \rho^{\v} (U_s\cdot \nabla) u^{(\v)} + \rho^{\v} (U^{(\v)}\cdot\nabla)u_s  +  (U_s\cdot\nabla)u_s\, \rho^{(\v)} +  p_x^{(\v)}\\
		\qquad\qquad\qquad\qquad\qquad\qquad\qquad\qquad = \mu \Delta_{\v} u^{(\v)}  + \lambda \v \pa_x\mdiv{U}^{(\v)} + R_{u}^{(\v)}  - \v^{-N_0}B_{s1},\\
		\dis \v \rho^{\v} (U_s\cdot\nabla)v^{(\v)} + \v \rho^{\v} (U^{(\v)}\cdot\nabla)v_s + \v\, \rho^{(\v)} (U_s\cdot\nabla)v_s +  p^{(\v)}_y \\
		\qquad\qquad\qquad\qquad\qquad\qquad\qquad\qquad = \mu \v \Delta_{\v} v^{(\v)} + \lambda \v \pa_y\mdiv{U}^{(\v)}  + \v R_{v}^{(\v)}  - \v^{1-N_0}B_{s2},\\
		\dis 2 \rho_s (U_s\cdot\nabla)T^{(\v)} + 2\rho_s (U^{(\v)}\cdot\nabla)T_s
		+ 2(U_s\cdot\nabla)T_s \,  \rho^{(\v)}  - (U_s\cdot\nabla)p^{(\v)}  - (U^{(\v)}\cdot\nabla)p_s\\
		\dis \qquad   =\kappa \Delta_{\v}T^{(\v)} + 4\mu \mathbb{S}^{\v}(U_s):\mathbb{S}^{\v}(U^{(\v)}) + 2 (\lambda-\mu) \v \mbox{div}U_s\cdot \mbox{div}U^{(\v)} + R^{(\v)}_{T} - \v^{-N_0}B_{s3},
	\end{cases}
\end{align}
where $U^{(\v)} := (u^{(\v)},v^{(\v)})^{t}$ and
\begin{align}\label{1.9}
	\begin{split}
		R_u^{(\v)}&:=-\v^{N_0} \rho_s(U^{(\v)}\cdot \nabla)u^{(\v)}   -\v^{2N_0} \rho^{(\v)}(U^{(\v)}\cdot\nabla)u^{(\v)}, \\
		R_v^{(\v)}&:=-\v^{N_0} \rho_s(U^{(\v)}\cdot \nabla)v^{(\v)}   -\v^{2N_0} \rho^{(\v)}(U^{(\v)}\cdot\nabla)v^{(\v)} ,\\
		R^{(\v)}_{T}&:=- 2\v^{N_0}\big[\rho_s  U^{(\v)} + \rho^{(\v)} U_s + \v^{N_0}\rho^{(\v)} U^{(\v)}\big]\cdot \nabla T^{(\v)} - 2\rho^{(\v)} (U^{(\v)}\cdot\nabla)T_s \\
		&\quad  + \v^{N_0} ( U^{(\v)}\cdot\nabla) p^{(\v)}  + 2\mu \v^{N_0}  |\mathbb{S}^{\v}(U^{(\v)})|^2 + (\lambda-\mu) \v^{1+N_0} |\mbox{div} U^{(\v)}|^2.
	\end{split}
\end{align}
The boundary conditions of $(u^{(\v)}, v^{(\v)}, T^{(\v)})$ are
\begin{align}\label{1.8-2}
	\begin{cases}
		\dis (u^{(\v)}, v^{(\v)})|_{y=0}=(0,0),\\
		\dis T^{(\v)}_y|_{y=0}=0 \,\,\,\,\, \mbox{for NT},\\ 
		\dis \mbox{or}\,\, T^{(\v)}|_{y=0}=0 \,\,\,\,\, \mbox{for DT}.
	\end{cases}
\end{align}

We note that the Prandtl layer expansion $(\rho_s, u_s, v_s, T_s)$ is not periodic-in-$x$, and 
there is no natural physical BC {\it for the remainder equations \eqref{1.8-1}} at
$x=0,\, L$ in a channel. A rigorous justification of the Prandtl expansion \eqref{1.5} can
be formulated as the following PDE problem: assume the remainder $(u^{(\v)
},v^{(\v)}, p^{(\v)},T^{(\v)})$ is bounded (in a suitable
sense) at $x=0,\, L$, then prove $(u^{(\v)
},v^{(\v)}, p^{(\v)},T^{(\v)})$ is bounded in the bulk.

It is well-known that the basic energy estimate fails due to the boundary layer profile from the convection term 
\begin{align}\label{1.28-v}
	\la u_{sy}u^{(\v)}, \, v^{(\v)}\ra,
\end{align}
for which $v^{(\v)}$ is $\v^{-1/2}$ out of control from the natural dissipation $\v^{1/2}\|\nabla_\v v^{(\v)}\|$.   Our main goal is to develop a PDE framework to obtain $\varepsilon$-uniform
solvability and bounds for the remainder $(u^{(\v)
},v^{(\v)}, p^{(\v)},T^{(\v)})$ subject to the natural physical BC \eqref{1.8-2}
at $y=0$, with \textit{suitably imposed well-posed} mathematical BC
at $x=0,\, L$. These lateral BCs (and their non-homogeneous versions)
characterize mathematically the boundedness of the remainder $(u^{(\v)
},v^{(\v)}, p^{(\v)},T^{(\v)})$ at $x=0,\, L$.

In comparison to the incompressible case with only the velocity field
$(u^{(\v)}, v^{(\v)})$, we adopt the quotient estimates in Guo-Iyer \cite{Guo-Iyer-2024} 
to overcome the loss of $\v^{1/2}$ in \eqref{1.28-v}. The new mathematical challenge is to \textit{design compatible} BCs for $(p^{(\v)},T^{(\v)})$ at $x=0,\, L$ to
ensure their   solvability and uniform-in-$\v$ bounds, particularly along the boundaries $x=0,L$ and $y=0$.

\subsection{Highlights of framework and main theorem}\label{Sec1.2}
In the following, we explain the key ideas and highlights of our framework. From now on, we drop the superscript of $(\rho^{(\v)}, p^{(\v)}, T^{(\v)}, u^{(\v)}, v^{(\v)})$ when no confusion arises.

\smallskip

\noindent $\bullet $ {\it Generalized stream function  $\phi $.}

We introduce the {\it stream function }
$\phi $ such that 
\begin{align}\label{7.1-1}
\rho^{\v} U  = \nabla^{\perp} \phi - U_s \rho  - \v^{-N_0} A_{s} \quad\mbox{with}\,\,\,\nabla^{\perp} :=(\pa_y, -\pa_x)^t,
\end{align}
which solves the mass continuity equation $\eqref{1.8-1}_1$. Then, by noting $U_s|_{y=0}=A_s|_{y=0}=0$,  we can reformulate the
non-slip BC $\eqref{1.8-2}_1$ equivalently as:
\begin{align}\label{2.17-0}
\phi|_{y=0}&=\phi_{y}|_{y=0}=0.
\end{align}
Using \eqref{7.1-1} and  \eqref{1.19-1}, after tedious calculations, one can rewrite the momentum equations $\eqref{1.8-1}_{2,3}$ in the equivalent form \eqref{7.36}, see more details in Appendix A.

\medskip

\noindent $\bullet $ \textit{Div-curl reformulation and paired boundary conditions \eqref{pair1}. }

We reformulate the problem in $[\phi,p,T]$.
We solve equivalently the momentum equations by taking the divergence and curl of momentum equations \eqref{7.36} (i.e.,  $A_1=0,\,\, A_2=0$ in the bulk), roughly as
\begin{align}\label{1.33}
\big(u_s^2 q_{xy}\big)_y  + \v \big(u_s^2 q_{xx}\big)_x   - \mu \Delta_{\v}^2\phi
&= \big(\f{\rho_s}{T_s} u_s^2 T_x\big)_y +  \cdots,\\
\Delta_{\v}p +  (\mu+\lambda)\v \frac{1}{p^{\v}} \big(U^{\v}\cdot \nabla\big) \Delta_{\v}p
&=O(\v)\rho_{sy}\big\{\Delta_{\v}\phi_y +  \Delta_{\v}\phi_x\big\} +\cdots.\label{1.34}
\end{align}

One of the main foci of the paper is the solvability of \eqref{1.33}-\eqref{1.34}.
In order to recover original momentum equations $A_{1}= A_{2}=0,$ we
must impose vanishing condition for $A_{1}$ or $A_{2}$ at the boundary by the
{\it div-curl} inversion Lemma \ref{lemR.8}.  The paired BC \eqref{pair1}  precisely
eliminate the top 3rd order derivatives $\nabla_{\varepsilon}^{3}\phi$ at $A_{1}$ and
$A_{2}$ at $x=0,\, L$  connecting naturally to Guo-Iyer \cite{Guo-Iyer-2024} $H^{4}$
theory. Together with $A_1|_{y=0}=0$, we obtain  the following desired estimates for closure:
\begin{align}
\|\phi\|_{H^{4}} &  \lesssim \|\nabla p\|_{H^{1}}+ \|u_{s}p_{y}\|_{H^{2}}+\cdots,\\
\|\nabla p\|_{H^{1}} &  \lesssim \|\phi_{yyy}\|_{H_{y=0}^{1/2}}+\cdots,\\
\|u_{s}p_{y}\|_{H^{2}} &  \lesssim \|u_{s}\nabla_{\varepsilon}^{3}\phi\|_{H_{x=0}^{1/2}}+|u_{s}\nabla_{\varepsilon}^{3}\phi|_{H_{x=L}^{1/2}}+\cdots.
\end{align}
It is important to note that for other types of BCs for $\phi$ (e.g.,
the classical no-slip BC $\phi=\phi_{x}=0$ at $x=0,\, L$)$,$ 3rd
order derivatives $\nabla_{\varepsilon}^{3}\phi$ can not be eliminated, which
leads to severe loss of derivative\textbf{ }for the pressure estimate.

The design of such a paired BC \eqref{pair1} is the first technical highlight with structural
cancellation and the starting point of our analysis. 


 \medskip
\noindent $\bullet $ \textit{Paired pseudo entropy $S$ and quotient $q$ \eqref{s}, modified div-curl system} 

We note from the temperature equation $\eqref{1.8-1}_4$: 
\begin{align}
2\rho_{s}u_{s}T_{x}-\kappa\Delta_{\v}T=[2T_{sy}- \rho_s^{-1} p_{sy} 
]\phi_{x}+u_{s}p_{x}+\cdots.
\end{align}
A rough $H^{2}$ estimate via $T_{xx}$ multiplier yields natural bound
\begin{align}
\|T_{xy}\| \backsim \v^{-\f12} \left\| [2T_{sy}-\rho_s^{-1} p_{sy} ]\phi
_{x}+u_{s}p_{x}\right\|,
\end{align}
which leads to a $\varepsilon^{-1/2}$ loss in $ \big(\f{\rho_s}{T_s} u_s^2 T_x\big)_y$ for
the closure for $\phi$ from \eqref{1.33}!

We introduce \textit{new good unknowns} to avoid such a loss of $\varepsilon
^{1/2}.$ To make use of the exact structure of the problematic terms,  we define the {\it pseudo entropy $S$ and its corresponding quotient $q$} in (\ref{s}).  We note that 
the  pseudo entropy $S$ now satisfies  (see Appendix \ref{subsec10.3} for  details of calculations and notations)
\begin{align}\label{T.1}
	2 \rho_s u_s S_x + 2\rho_s v_s S_y - \kappa \Delta_{\v}S
	&= \kappa \Delta_{\v}\big(\frac{1}{2\rho_s}p\big) + \kappa \Delta_{\v}\big(\f{1}{p_s}T_{sy}q\big) + \mathfrak{J}(\phi,p,T).
\end{align}
For the multiplier $S_{xx}$ of \eqref{T.1}, with possible integration by parts in $x$ or $y$ in $\big\la T_{sy}q_{yy},\, S_{xx}\big\ra,$
one  gets
\[
  \|S_{xy}\| \backsim \v^{-\f12}\|\Delta_{\v}p\| + \|q_{xy}\|
\]
which could be exactly controlled in the expected $\phi$ estimate and $p$ estimate.

On the other hand, our precise design of quotient $q$ (distinct from Guo-Iyer \cite{Guo-Iyer-2024}) allows an exact cancellation of the bad term $u_s^2 T_s^{-2}T_{sy} q_x$ created by the new unknown $S$ in the 1st momentum equation $\eqref{7.36}_1$, i.e., 
\begin{align}\label{3.4-0}
\big[u_s \phi_{xy} - u_{sy} \phi_x\big] + u_s^2T_s^{-2}T_{sy} q_x
&= m u_s^2 q_{xy} + u_s\bar{u}_{sx} q_y + [u_s \bar{u}_{sxy} - u_{sy} \bar{u}_{sx}]q\nonumber\\
&\quad  +  \cancel{u_s^2 (T_s^{-1})_y q_x}  + \cancel{u_s^2 T_s^{-2}T_{sy} q_x},
\end{align}
where we have used the following notations
\begin{align}
m:=T_s^{-1} \quad \mbox{and}\quad \bar{u}_s:=m u_s \equiv T_s^{-1} u_s.
\end{align}
Then, based on \eqref{3.4-0}, after tedious calculations, one can further rewrite the momentum equations  $\eqref{1.8-1}_{2,3}$ as
\begin{align}\label{2.92}
	\begin{cases}
		\dis A_1:=m u_s^2 q_{xy}  + \mathbf{d}_{11}\, p_x+ \hat{d}_{11} S_x  + g_1(p,\phi, T(S,p,q))  + \mathfrak{g}_1(q)  -\frac{\mu}{\rho^{\v}} \Delta_{\v}\phi_y \\
		\dis \qquad\,\,\,  + \frac{\mu+\lambda}{p^{\v}} u^{\v} \Delta_{\v}p - \frac{\lambda}{p^{\v}} \big[1- \f{\rho^{\v}}{2\rho_s}\big] u^{\v} p_{yy} + \frac{\lambda\v}{p^{\v}} v^{\v} p_{xy}  + \frac{\lambda}{T^{\v}} u^{\v} S_{yy} - G_1 - \mathcal{N}_1=0,\\[3mm]
		\dis A_2:= -\v m u_s^2 q_{xx} +   \tilde{d}_{22}\, p_y + g_2(p,\phi, T(S,p,q)) + \v \mathfrak{g}_2(q) +\frac{\mu\v }{\rho^{\v}}   \Delta_{\v}\phi_x + \frac{\mu\v}{p^{\v}} v^{\v} \Delta_{\v}p \\
		\dis \qquad\,\,\,  + \frac{\lambda\v}{p^{\v}} \big[1-\f{\rho^{\v}}{2\rho_s}\big] u^{\v} p_{xy} + \frac{\lambda\v}{p^{\v}} v^{\v} p_{yy} - \frac{\lambda\v}{T^{\v}} u^{\v} S_{xy}  - \v G_2 - \v \mathcal{N}_2=0.
	\end{cases}
\end{align}
We refer the readers to Appendix \ref{subsec10.1} for further details of calculations and notations of \eqref{2.92}.

We note that   \eqref{T.1} and \eqref{2.92} are now
compatible with desired estimates for $[\phi ,p].$ Moreover, the
introduction of pseudo entropy $S\,$\ entails a natural structure of the
generalized Rayleigh operator for $\phi $ equation and the correct subsonic
(ellipticity) condition for solving $p$%
\begin{equation}\label{1.28}
\mathbf{d}_{11}= 1-\f1{2T_s} u_{s}^{2} + O(1)\v u_s >0.
\end{equation}

By considering the curl of momentum equations \eqref{2.92}, i.e., $	\pa_y\big(b_1A_1\big) - \pa_{x} \big(b_2 A_2\big) =0$, to obtain
\begin{align}\label{1.33-0} 
	\big(b_1m u_s^2 q_{xy}\big)_y  + \v \big(b_2mu_s^2 q_{xx}\big)_x   - \mu\Big(\frac{b_1}{\rho^{\v}} \Delta_{\v}\phi_y\Big)_y - \mu\v \Big(\frac{b_2}{\rho^{\v}}   \Delta_{\v}\phi_x\Big)_x 
	=\big(\f{\rho_s}{T_s}b_1 u_s^2 S_x\big)_y   + \cdots,
\end{align}
where we have denoted 
\begin{align}\label{84}
	b_1:=\mathbf{d}_{11}^{-1},\quad b_2=\tilde{d}_{22}^{-1}.
\end{align}

The introduction of paired  pseudo entropy $S$ and the corresponding quotient $q$ marks the second technical highlight with structural cancellation. 

\smallskip

\noindent $\bullet $ $H^{4}$ \textit{theory for the stream function  $\phi$ with $H^{1/2}_{00}$}.

We follow Guo-Iyer \cite{Guo-Iyer-2024} quotient estimate closely. For $q_{xx}$ multiplier, even  with better $S_{xy}$ estimate, the most
delicate part is to control the interaction with the pseudo entropy, i.e., $\big\la ( \f{b_1}{T_s} \rho_s u_s^2 S_x )_y, \,  q_{xx} \big\ra$. In fact, by noting the spacial weight $u_s^2$, we integrate by parts in $y$ and use the pseudo equation \eqref{T.1} to obtain
\begin{align}
	\big\la ( \f{b_1}{T_s} \rho_s u_s^2 S_x )_y, \,  q_{xx} \big\ra &= -\big\la \f{b_1}{T_s} \rho_s u_s^2 S_x , \,  q_{xxy} \big\ra \nonumber\\
	&\backsim -  \big\la  u_s\Delta_{\v}S ,\, q_{xxy} \big\ra  - \big\la u_s \Delta_{\v} p,\, q_{xxy} \big\ra - \big\la u_s T_{sy} \Delta_{\v} q,\, q_{xxy} \big\ra + \cdots.
\end{align}
Then integrating by parts in $x,y$, we can control terms $\big\la  u_s S_{yy} ,\, q_{xxy} \big\ra$ and $\big\la u_s T_{sy} q_{yy},\, q_{xxy} \big\ra$, see step 4.1 in the proof of Lemma \ref{lem3.1} for more details. This also emphasizes the strong coupling of the velocity and thermal layers.

One of the fundamental {\it new} difficulties is the estimate for $\|\nabla^{2}p\|$. To this end, we must control the trace
$\|\phi_{yyy}\|_{H_{y=0}^{1/2}}$ from $A_1=0$ at $y=0$. Consequently, we must control $\|\phi_{xyyy}\|$  without the
weight $u_{s}$, unlike the weighted norm  $\|\sqrt{u_s}\phi_{xyyy}\|$ estimated in Guo-Iyer \cite{Guo-Iyer-2024}.  
For this purpose,  we must control (see \eqref{3.Tr152} and \eqref{3.Tr154-1})
\begin{equation*}
\v	\int_{0}^{L}\phi _{xyyy}(x,0)\phi _{xyy}(x,0)dx.
\end{equation*}
This is a typical difficulty in high energy estimate for which $\phi_{xyyy}|_{y=0}$ itself is beyond the desired $H^{4}$ norm of $\phi$. However, this can be overcome if one can make a formal integration by part of $1/2$ derivative $\partial _{x}$ for the boundary contribution, which is valid for a periodic function but invalid in general interval $[0,L]$.

Thanks to the boundary conditions $\phi_{xyy}|_{x=L}=\phi_{yyy}|_{x=0}=0$ from \eqref{pair1}, we are able to employ a subtle $H_{00}^{1/2}$ theory (see Lemmas \ref{lemA.9} \& \ref{lemAH.1}) for an integration by part of $1/2$ order of $\partial _{x}$  (see \eqref{3.Tr154-1}--\eqref{3.Tr27})
\begin{align}
\v \int_0^L\phi_{xyyy}  \phi_{xyy} dx \Big|_{y=0} 
&= \v \int_0^L \phi_{xyyy}  \phi_{xyy} \chi(\frac{4x}{L}) dx \Big|_{y=0} +  \v \int_0^L \phi_{xyyy}  \phi_{xyy} \bar{\chi}(\frac{4x}{L}) dx \Big|_{y=0} \nonumber\\
&\lesssim \f1N \|(\v \phi_{xxyy}, \sqrt{\v}\phi_{xyyy})\|^2 + N \|\phi_{yyyy}\|^2 + L^{\f14} \sqrt{\v}[\phi]_{3,1},
\end{align}
where $\chi(\cdot)$ is a smooth cut-off function defined in \eqref{6.0.1}.

Multiple applications of subtle $H_{00}^{1/2}$ theory mark the third highlight of our framework, which enable us to make crucial integration by	parts of $1/2$ derivative at the boundary. We remark that the need of moving	$1/2$ derivative to match the trace theorem at the boundary is a generic	difficulty in energy method in many PDE problems.

\medskip

\noindent $\bullet $ \textit{Pressure estimate with boundary condition \eqref{pressure}.}

One of the main challenges is to solve the pressure $p$ of \eqref{1.34} uniformly in
$\varepsilon$ with
\begin{align}\label{1.50}
A_{1}|_{y=0}=A_{1}|_{x=0}=A_{2}|_{x=L}=0.
\end{align}
 This is creafully chosen so that with $\Delta_{\varepsilon}%
\phi|_{x=0}=\Delta_{\varepsilon}\phi_{x}|_{x=L}=0,$ we may take $\partial_{y}$
derivative on $A_{1}$ and $A_{2}.$ We may  view \eqref{1.34} as a
{\it transport} equation for $\Delta_{\varepsilon}p$,  while  \eqref{1.50} provides
non-classical BC to further solve $p.$ We therefore can impose an inflow BC
for $\Delta_{\varepsilon}p|_{x=0}$ lead to positive contribution at $x=0, L$.


We explain the key idea on the  positive contribution at $x=0, L$. Multiplying \eqref{1.34} by $\Delta_{\v}p$, we have
\begin{align}\label{1.42}
&\|(\v p_{xx}, p_{yy})\|^2 + 2\|\sqrt{\v}p_{xy}\|^2 +  (\mu+\lambda) \v \int_0^\infty \frac{u_s}{2p_s} |\Delta_{\v}p|^2 dy \Big|_{x=L} + 2\v \int_{0}^\infty p_{x}p_{yy} dy\Big|_{x=L} \nonumber\\
&\leq  (\mu+\lambda) \v \int_0^\infty \frac{u_s}{2p_s} |\Delta_{\v}p|^2 dy\Big|_{x=0} +  2\v \int_{0}^\infty p_{x}p_{yy} dy\Big|_{x=0} + (\cdots),
\end{align}
where we only keep the linear leading terms, and $(\cdots)$ contains the corner terms and nonlinear terms. Combining  both boundary conditions $A_1\big|_{x=0}=0$ and \eqref{pressure},
we have 
\begin{align*}
 \frac{1}{2} \lambda  u_s p_{yy} \big|_{x=0}  &= p_s\mathbf{d}_{11}\, p_x + (\mu+\lambda) u_s \Delta_{\v}p + \cdots = - (\sigma-1)  p_s\mathbf{d}_{11}\, p_x + \cdots,
\end{align*}
which yields a positive contribution
\begin{align}
\mbox{Terms at $x=0$ of } \eqref{1.42} 
&=\v \int_0^\infty\Big\{\f{\sigma^2\mathbf{d}_{11}}{2(\mu+\lambda)} - \f{4}{\lambda}(\sigma-1)\Big\} \f{p_s}{u_s}\mathbf{d}_{11} p_x^2 dy \Big|_{x=0}\nonumber\\
&\lesssim -c_0\int_0^\infty \f{p_x^2}{u_s}  dy \Big|_{x=0},
\end{align}
where we have optimized $\sigma>1$ from 
\begin{align}\label{1.44}
\f{\sigma^2\mathbf{d}_{11}}{2(\mu+\lambda)} - \f{4}{\lambda}(\sigma-1)
&=\f{1}{2(\mu+\lambda)} \Big\{\Big(\sigma\sqrt{\mathbf{d}_{11}} - \f{4(1+\f{\mu}{\lambda})}{\sqrt{\mathbf{d}_{11}}}\Big)^2 -\f{8}{\mathbf{d}_{11}}\Big(1+\f{\mu}{\lambda}\Big) \Big(2(1+\f{\mu}{\lambda})-\mathbf{d}_{11}\Big) \Big\}\nonumber\\
&\leq -c_0<0.
\end{align}
 In \eqref{1.44}, the smallness of Mach number is not needed, only  the subsonic condition is enough.

On the other hand, we have from $A_2|_{x=L}=0$ that 
\begin{align*}
p_y\big|_{x=L}=-\lambda \v u_s p_{xy} + \cdots
\end{align*}
which also yields a positive contribution at $x=L$ that 
\begin{align}
2\v \int_{0}^\infty p_{x}p_{yy} dy\Big|_{x=L} = - 2\v \int_{0}^\infty p_{xy}p_{y} dy\Big|_{x=L} + \cdots= \f{2}{\lambda} \int_0^\infty u_s |p_{xy}|^22 dy \Big|_{x=L} + \cdots.
\end{align}

The optimal choice of $\sigma$ is more or less unique for arbitrary $\lambda$ and $\mu$ in \eqref{1.44}. The discovery of BC \eqref{pressure} marks the fourth highlight of our framework, which leads to intricate positivity at the boundaries
$x=0, L$ and solvability of the pressure in $H^{2}.$

\medskip

\noindent $\bullet $ \textit{Estimate of the pseudo entropy  $S$ with \eqref{0.3}}.

To solve the pseudo entropy equation \eqref{T.1}, we impose the following boundary conditions at $x=0,L$
\begin{align}
	S \big|_{x=0}= -  \frac{\chi}{2\rho_s}p \quad \mbox{and}\quad  S_{x}\big|_{x=L}=0 
	\quad \mbox{for NT},\label{T.2}\\
 S \big|_{x=0}=0 \quad \mbox{and}\,\,\, \eqref{0.3}\,\,
		\,\, \mbox{for DT},\label{T.2-D}
\end{align}
Noting \eqref{2.17-0}, $T_{sy}|_{y=0}=0$ for NT and $\eqref{1.8-2}_{2,3}$, one derives
\begin{align}\label{T.3}
	\begin{split}
		S_y \big|_{y=0} = - \big(\frac{1}{2\rho_s}p\big)_y \,\,\,\,\, \mbox{for NT}\, \Big(\mbox{or}\,\, S \big|_{y=0}= -\f{1}{2\rho_s}p \,\,\,\, \mbox{for DT}\Big).
	\end{split}
\end{align}
We also impose
\begin{align}\label{3.61-11}
	S(0,0)=0,
\end{align}
which yields  $S(x,+\infty)\neq 0$  in general. 

Our uniform-in-$\v$ estimate  for the pseudo entropy $S$ depends on
the multiplier $S_{xx}$: 
\begin{align}
\big\la S_{yy},\, S_{xx} \big\ra 
&=  \|S_{xy}\|^2 - \la S_{y},\, S_{xy} \ra_{x=L} +  \la S_{y},\, S_{xy} \ra_{x=0} - \la S_{y},\, S_{xx}  \ra_{y=0}.
\end{align}
For the case of NT, thanks to the insulated BC $T_{y} |_{y=0} =0$ and \eqref{T.2}, we express
the key mixed boundary term at $y=0$ as a solid integral in the bulk, then a
side integral at $x=0$  
\begin{align}
\big\la S_{y},\, S_{xx} \big\ra_{y=0}&=-\Big\la \big(\frac{\chi }{2\rho _{s}}p\big)_{y},S_{xx}\Big\ra_{y=0} 
=\Big\la \big(\frac{\chi }{2\rho _{s}}p\big)_{yy},\, S_{xx}\Big\ra + \Big\la \big(\frac{\chi }{2\rho _{s}} 
p\big)_{y},\, S_{xxy}\Big\ra \nonumber\\
&=-\Big\la \big(\frac{\chi }{2\rho _{s}}p\big)_{y},\, S_{xy}\Big\ra_{x=0}+\cdots.
\end{align}
There is an {\it exact cancellation} of $\dis -\Big\la \big(\frac{\chi }{2\rho _{s}}p\big)_{y},\, S_{xy}\Big\ra_{x=0}$ and $\la S_{y},\, S_{xy} \ra_{x=0}$ in the process leading to desired uniform estimate, see the proof of Lemma \ref{lemTH6.3} for details. For the uniform estimate of third order derivative $[[[S]]]_{3,w}$, the $H^{1/2}_{00}$ argument also plays a key role, see \eqref{1TC.76}. 

For the case of DT with $S_{xx}$ multiplier, we note  
\begin{align}\label{A1.56}
A_1|_{y=0}=0 \quad  \Longleftrightarrow \quad p_x|_{y=0} \backsim \phi_{yyy}(x,0) + \cdots,
\end{align}
the term $\big\la S_{y},\, S_{xx} \big\ra_{y=0}$ is more singular and out of control.  
Instead, by imposing the new BC \eqref{0.3}, after delicate integration by parts, we focus on $x=L$ and $y=0$ terms
\begin{align}
&\Big\la S_{x}, \, \kappa S_{yy}- \rho_s u_s S_x +\f{\kappa }{p_s}T_{sy}q_{yy} + 2 \big(\f{\kappa}{p_s}T_{sy}\big)_{y}q_y + \big(\f{\kappa}{2\rho_s}\big)_{yy}p + \mathfrak{M}_0\Big\ra_{x=L} + \big\la S_x, S_{xy}\big\ra_{y=0}  \nonumber\\
&\backsim \|\sqrt{u_s}S_x\|^2_{x=L} - \big\la p_x, S_{xy}\big\ra_{y=0} + \cdots \nonumber\\
&\backsim  \|\sqrt{u_s}S_x\|^2_{x=L} - \big\la \phi_{yyy}, S_{xy}\big\ra_{y=0} + \cdots\nonumber\\
&\backsim \|\sqrt{u_s}S_x\|^2_{x=L} + \big\la \phi_{yyyy}, S_{xy}\big\ra + \big\la \phi_{yyy}, S_{yy}\big\ra_{x=L} + \cdots.
\end{align}
Then, again with the help of \eqref{0.3},  we can obtain 
\begin{align}
\|S_{yy}\|^2_{x=L}\lesssim a \|\sqrt{u_s}S_x\|^2_{x=L} + \cdots.
\end{align}
 Moreover, a subtle improvement to estimate  $\|u_s\Delta_{\varepsilon}p_{y}\|$ is needed to obtain a small constant for closure, see Lemma \ref{CBlem7.6} for details.

 We remark that \eqref{0.3} is almost unique with respect to $p,\phi$ (out of control otherwise).
The discovery of a BC \eqref{0.3} at $x=L$ marks the fifth highlight of our framework with 
\textit{structural cancellation}, demonstrating the intricate interaction
between the temperature $T$ with the velocity and pressure.

\smallskip

\noindent $\bullet $ {\it PDE construction with non-classical BC.}

Even though these new imposed BCs lead to uniform {\it a priori} estimates,
construction of such solutions requires delicate analysis due to their
non-standard nature. A major effort is devoted to the careful study of such
well-posedness question througout the paper. We illustrate the precise construction of pressure as an example.

We remark that it is not standard to solve the pressure equation in the
present of the new viscous-inflow BC (\ref{pressure}) with 
\begin{align}\label{3.61-1}
	p(0,0) =0,
\end{align}
which yields that $p(x,+\infty)\neq 0$. Furthermore, due to
subtle difficulty at the corners $(0,0)$ and $(0,L),$ we must avoid corner
values  $\nabla _{\varepsilon }p(0,0),$ $\nabla _{\varepsilon }p(0,L)$,  which are exactly beyond of our working functional space$.$ To
overcome these difficulties, it is important for us to reformulate the
pressure equation as a \textit{coupled} system in terms of $\mathbf{w}%
=\Delta _{\varepsilon }p$ and a new principle {\it good unknown}: 
\begin{equation*}
	B:=  p_x + \f{1}{\mathbf{d}_{11}} \big\{g_1(p,\phi, T(S,p,q))  - \mathfrak{h}_1 - \mathfrak{N}_{11}\big\},
\end{equation*}%
where $\mathfrak{N}_{11}, \mathfrak{h}_1$ are nonlinear terms defined in \eqref{17.14-1} and \eqref{7.33-11} respectively.
The equations  of $\mathbf{w}$ and $B$ are 
\begin{align} \label{w1}
	\begin{cases}
		\dis \mathbf{w}+(\mu +\lambda )\varepsilon \frac{1}{p^{\v}}(U^{\v}\cdot \nabla ) 
		\mathbf{w}  = \cdots,  \\
		\dis (\mu +\lambda )u_{s}\mathbf{w} + 2(\frac{\mu }{\lambda }+1)p_{s} 
		\mathbf{d}_{11}B \big|_{x=0}  =0,
	\end{cases}
\end{align}
and
\begin{align}
	\begin{cases}
		\dis \Delta _{\varepsilon } B =\mathbf{w}_{x} + \Delta_{\v}\Big( \f{1}{\mathbf{d}_{11}} \big\{g_1(p,\phi, T(S,p,q))  - \mathfrak{h}_1 - \mathfrak{N}_{11}\big\}\Big),  \label{B1} \\[3mm]
		\dis (2\frac{\mu }{\lambda }+1)p_{s}\mathbf{d}_{11}B -\varepsilon \frac{1}{2} 
		\lambda (1+\chi )u_{s}B_{x}\big|_{x=0} =-\frac{1}{2}\lambda (1+\chi )u_{s}\mathbf{w} + \cdots,\\[3mm]
		\dis \v B_{x} - \lambda\v\frac{b_2}{2p_s} u_s  B_{yy} \big|_{x=L}  = \mathbf{w} + \cdots,\\[3mm]
		\dis 	B\big|_{y=0}=\frac{\mu }{\rho^{\v}}\phi _{yyy}.
	\end{cases}
\end{align}
The boundary condition in \eqref{B1} is in fact the crucial necessary boundary condition \eqref{1.50}.
An auxiliary boundary correction $\theta $ is constructed in Lemma \ref{lemAF} to
remove the boundary contribution 
\begin{equation*}
	B\big|_{y=0}=\frac{\mu }{\rho^{\v}}\phi _{yyy},
\end{equation*}%
and remove the dependence of corner value of $B$ at $(0,0)$ and $(0,L)$.
Furthermore, further subtle $\delta -$approximation of \eqref{w1}-\eqref{B1}
are also needed to carry out the construction, see section \ref{sec7} for more details.

\smallskip

Our framework leads to the following main theorem:
\begin{theorem}[Main Theorem]\label{thm1.1}
	Assume the background Mach number $\mathscr{M}_{a}$ is suitably small and $0<\v\ll L\ll\mathscr{M}_a$. Under the boundary conditions \eqref{pair1}, \eqref{0.3}, \eqref{pressure}, \eqref{2.17-0}, \eqref{T.2}-\eqref{3.61-11} and \eqref{3.61-1}, there exist  strong solutions $(\phi, p, S)$ to \eqref{T.1} and \eqref{2.92} satisfying
	\begin{align}\label{1.51}
		\|(\phi,p,S) \hat{w}^{\v}_0\|^2_{\mathbf{X}^{\v}}\lesssim \v^{N_0}.
	\end{align}
	Consequently,  we  construct strong solutions $(\rho^{\v}, u^{\v}, v^{\v}, T^{\v})$ of the steady  compressible Navier-Stokes-Fourier equations \eqref{1.4}-\eqref{1.4-1} in the form of \eqref{1.5} such that
	\begin{align}\label{1.52}
		\|(\rho^\v- \rho^0_e-\rho^0_p, u^{\v}-u^0_e-u^0_p, v^{\v}-v_e^1-v^0_p, T^{\v}-T_e^0-T^0_p)\|_{L^\infty}\lesssim \sqrt{\v}.
	\end{align}
	Therefore we justify the validity of Prandtl layer expansion of steady compressible Navier-Stokes-Fourier in the inviscid limit.
\end{theorem}

\begin{remark}
	All the boundary conditions can be inhomogeneous of high order $\v$ which ensure the boundedness for the remainder at the boundaries.
\end{remark}

\begin{remark}
	We note that $L\ll 1$ is required for the solvability of Prandtl
	equations \eqref{3.6-0} in \cite{Guo-Wang} for general data, and it is also required for the validity of
	Prandtl layer expansion for the incompressible case \cite{Guo-Iyer-2024}. As in \cite{GZ,Iyer-M} , we
	expect this can be improved for more special Prandtl profiles, which is left
	for a future study.
\end{remark}

\begin{remark}
	If singularities appear in the Euler solution $\rm($see, e.g.,  Chen et al. \cite{BCF2009,CHWX2019,ChenF2003} and   Xin et al. \cite{LXY2009}$\rm)$, it is an interesting problem to study  the interaction between Prandtl layer and shock wave $\rm($or contact discontinuity$\rm)$. Such problems will be addressed in future work.
\end{remark}

\subsection{Literature review} Due to the importance of inviscid limit in fluids, there have been lots of results related to our current study.

\noindent{$\bullet$ \it Steady compressible Navier-Stokes-Fourier flows with $\v\gtrsim1$}: 
There are numerous results on the well-posedness of the steady compressible Navier-Stokes equations ($\v\backsim 1$). Here we review only those most closely related to our work;  readers are referred to the cited literature for a more comprehensive bibliography. For general existence theory of steady isentropic compressible Navier-Stokes equations, Lions \cite{Lions-1998} first established the existence of weak solutions  under certain restrictions on the adiabatic exponent, later sharpened by Novotny-Straskraba \cite{Novotny2004}, we also refer the reader to Jiang-Zhou \cite{JZ2018}, Kreml-Mucha-Pokorny \cite{KMP2018}, as well as the references therein.
For the general existence theory of steady CNSF, Mucha-Pokorny-Zatorska \cite{Mucha2016} proved the existence of weak solutions and variational weak solutions in smooth bounded domains; we note that the notion of variational weak solutions for CNSF was originally introduced by Feireisl \cite{Feireisl-2004}. As for the case of strong solutions, Kweon-Kellogg \cite{KK2004} established the existence and regularity of solutions to CNSF with nonzero boundary conditions in a two-dimensional polygon under sufficiently large viscosity and heat conductivity. Plotnikov-Ruban-Sokolowski \cite{Plotnikov} proved the existence of strong solutions to steady  CNSF in a smooth bounded three-dimensional domain with inhomogeneous boundary data for $\mathscr{M}_a\ll1$ and $Re\ll1$.  Dou-Jiang-Jiang-Yang \cite{Dou2015} showed the existence of strong solutions to the Dirichlet boundary value problem of steady CNSF with large external forces in a smooth bounded domain, provided $\mathscr{M}_a\ll1$; see also \cite{HWW2024} for the corresponding analysis of thermally driven flows. We also refer to \cite{GJZ}, which established the existence of strong solutions for steady isentropic compressible Navier-Stokes equations with inflow condition in a two-dimensional channel without smallness on Reynolds number or Mach number.

\noindent{$\bullet$ \it Steady incompressible Navier-Stokes flows with $\v\to0+$}: A fundamental question is to describe the asymptotic behavior of solutions of Navier-Stokes equations as  $\v\to0+$. Generically, there is a mismatch for the boundary conditions at the boundary $\{Y=0\}$ between the viscous flows and inviscid flows. In 1904, L. Prandtl \cite{Prandtl} proposed the famous boundary layer theory, which are already applied successfully in the engineering, to explain such a mismatch. We remark that the original Prandtl's theory was made for two-dimensional steady problem.

For the steady incompressible Prandtl equations \eqref{0.10}, the foundational result on the existence and uniqueness of strong solutions via the von Mises transformation and maximum principle is due to Oleinik  \cite{Oleinik1}(see also the monograph \cite{Oleinik-1999}), who also showed that favorable pressure gradients yield solutions global-in-$x$. Guo-Iyer \cite{Guo-Iyer-CMP} established the higher regularity of Prandtl solutions via quotient estimates and energy methods, leading to a Prandtl layer expansion up to any order, while Iyer-Masmoudi \cite{Iyer-M0} established global-in-$x$ regularity of the Prandtl system without invoking the classical von Mises transformation.  Wang-Zhang \cite{WZ} obtained the global $C^\infty$-regularity under favorable pressure gradients by using the maximum principle method. 
Serrin \cite{Serrin1967}  proved that the Blasius profile is indeed a self-similar attractor of Prandtl  equations  by using maximum principle techniques,  recently Iyer \cite{Iyer4} proved the asymptotic stability with explicit decay rate near the self-similar Blasius profile by using a new energy method, later Wang-Zhang \cite{WZ1} extended this result to  general initial data and Jia-Lei-Yuan \cite{JLY2025} established the optimal decay rate. 
Dalibard-Masmoudi \cite{DM} showed that separation occurs for the steady Prandtl equations  under an adverse pressure gradient $p'(x)=1$. Later, Shen-Wang-Zhang \cite{Wang-Zhang2021} extended this result to more general adverse pressure gradients satisfying $p'(x)\geq c>0$. Most recently, Iyer-Masmoudi \cite{Iyer-M2022} described the solution of steady Prandtl equations after separation for the first time.

It is a challenging problem to rigorously justify Prandtl's boundary layer expansion since it is indeed a multi-scaled problem. Recently, significant progress has been made for the steady incompressible Navier-Stokes equations. Guo-Nguyen \cite{Guo-N} justified the Prandtl expansion for  incompressible Navier-Stokes flows with moving boundary, and Iyer \cite{Iyer2} established the global-in-$x$ validity  subsequently, see also \cite{Iyer} for a nonshear Euler flows. 

For a motionless boundary, Gerard-Varet-Maekawa \cite{GM2} provided a  sobolev stability of Prandtl expansion for the case of external force and periodic-in-$x$. Guo-Iyer \cite{Guo-Iyer-2024}  justified the validity of the steady Prandtl layer expansion, including the celebrated Blasius boundary layer. Subsequently, Gao-Zhang \cite{GZ} established the expansion over a large interval in $x$ for Prandtl layer with a favorable sign. Iyer-Masmoudi \cite{Iyer-M} proved the global-in-$x$ stability of steady Prandtl expansions for two-dimensional incompressible Navier-Stokes equations. 

Recently, Iyer \cite{Iyer1} established the validity of the Prandtl expansion for incompressible Navier-Stokes flows over a rotating disk. Fei-Gao-Lin-Tao \cite{FGLT2} proved the existence of Prandtl-Batchelor flows on a disk with a wall velocity slightly different from rigid rotation by using the Prandtl expansion; see also \cite{FGLT1} for flows on an annulus and \cite{CFLZ} for Prandtl-Batchelor flows with a point vortex on a disk. Fei-Pan-Zhao \cite{FPZ} established the vanishing viscosity limit in a horizontally periodic strip.

\noindent{$\bullet$ \it Steady compressible Navier-Stokes-Fourier flows with $\v\to0+$}: 
Recently, Guo-Wang \cite{Guo-Wang} proved local-in-$x$ well-posedness of steady Prandtl layer and constructed Prandtl expansion $(\rho_s, u_s, v_s, T_s)$ of CNSF in strong Sobolev norms, a key foundation for the validity of Prandtl expansion in the present work.

For the isentropic compressible Navier-Stokes equations, Li-Yang-Zhang \cite{LYZ} established a Prandtl expansion in an $x$-periodic setting with external forcing, Jiang-Zhou \cite{JZ} proved the inviscid limit of steady supersonic shear flow with inflow boundary conditions.

\noindent{$\bullet$ \it Unsteady flows with $\v\to0+$}: There have been extensive results on this topic, we could only select those closely related to our current study. We recall the celebrated  Kato's criterion \cite{Kato,E} in $L^2$ convergence for the incompressible Navier-Stokes equations, we refer subsequent results in this direction \cite{Kelliher1,Kelliher,Temam,WangT1},  among others.
Mazzucato-Taylor \cite{MT} studied the vanishing viscosity for plane parallel channel flow. Constantin-Kukavica-Vicol \cite{Constantin} proved the inviscid limit under the conditions of the Oleinik  no back-flow and a vorticity lower bound.

For the unsteady incompressible Navier-Stokes equations with no-slip boundary conditions, the validity of the Prandtl expansion has been established for analytic data by Sammartino-Caflisch \cite{SC1,SC2}; further extensions were provided by Wang-Wang-Zhang \cite{WWZ} and Nguyen-Nguyen \cite{NN}. Maekawa \cite{Maekawa} justified the inviscid limit for initial vorticity that vanishes near the boundary in $\R_+^2$, and Fei-Tao-Zhang \cite{FTZ} later generalized Maekawa’s result to $\R_+^3$ using a direct energy method. Gerard-Varet-Maekawa-Masmoudi \cite{GMM} treated data in the Gevrey class. We also refer the reader to \cite{MM} for a comprehensive review on the inviscid limit and boundary layers.

In a stark contrast to the steady flow, based on \cite{Grenier2000},  Grenier-Nguyen \cite{GreN2} proved that the Prandtl limit is invalid in Sobolev space in general. Then, based on \cite{GGN1}, Grenier-Nguyen \cite{GN2024,GreN1} further proved the  Prandtl expansion is invalid. On the other hand, in the present of magnetic effect, Liu-Yang-Xie \cite{LXY} proved the validity of Prandtl expansion in Sobolev space.


Regarding CNSF, Liu-Wang-Yang \cite{LWY} established the local-in-time existence of  unsteady Prandtl layer. Recently, Yang-Zhang  \cite{YZ} proved the spectral instability
of subsonic boundary layers for isentropic Navier-Stokes equations with Mach number less than $1/\sqrt{3}$, then  Masmoudi-Wang-Wu-Zhang \cite{MWWZ2024} removed the restriction.

\noindent{\it Organization of the paper:} In Section \ref{sec3}, we introduce the $H^{1/2}_{00}$-estimate, div-curl inversion lemma, function spaces, norms, and some useful estimates that will be used frequently. Section \ref{sec2} presents the approximation systems and schemes. The existence and uniform estimates of the stream function are established in Section \ref{sec4-0}. Section \ref{sec7} proves the existence and uniform estimates of the pressure. In Section \ref{sec6}, we establish the existence and uniform estimates of the pseudo entropy for both Neumann and Dirichlet boundary value problems (BVPs). Finally, in Section \ref{Sec7}, we construct solutions to systems \eqref{T.1} and \eqref{2.92} and complete the proof of Theorem \ref{thm1.1}. Appendices \ref{SecB} and \ref{secC} establish the well-posedness of the temperature equation with Neumann boundary condition and the pseudo entropy equation with Dirichlet boundary condition, respectively, which are required for the existence of  pseudo entropy.

{\it Due to complexity of the paper, we provide details of reformulation for the steady CNSF, and include all its definitions and notations in  Appendix \ref{sec:appendixA} for the convenience of the readers.}

\section{Preliminary: Function Spaces, Norms and Useful Estimates}\label{sec3}

\subsection{$H^{1/2}_{00}$-space and useful inequalities}
In this part, we introduce the  $H^{1/2}_{00}$-space and some useful inequalities which are crucial in $H^4$-estimate of $\phi$.

\begin{definition}[\cite{Lions}]
We denote
\begin{align}\label{D.7}
	t^{\ast}:=
	\begin{cases}
		t,\quad t\in [0,1],\\
		1,\quad t\in (1,\infty).
	\end{cases}
\end{align}
Define 
\begin{align*}
	H^{1/2}_{00}(\R_+):=\Big\{ f ~\big|~ f\in H^{1/2}(\R_+) ~~ \& ~~ (t^{\ast})^{-\f12} f \in L^2(\R_+) \Big\} \subsetneq H^{1/2}(\R_+),
\end{align*}
with 
\begin{align}\label{D.6}
\|f\|^2_{H^{1/2}_{00}}:=\|f\|^2_{H^{1/2}} + \int_0^\infty (t^{\ast})^{-1} |f(t)|^2 dt.
\end{align}

\end{definition}

\begin{lemma}[$H^{1/2}_{00}$-estimate]\label{lemA.9}
	Assume $f(t)\in H^{1/2}_{00}(\R_+)$ and  $g(t)\in H^{1/2}(\R_+)$ , then it holds that 
	\begin{align}\label{D.8}
		\Big|\int_0^{\infty} \pa_t f\cdot g\,dt\Big|=\Big|\int_0^{\infty} f\cdot \pa_tg\,dt\Big|\lesssim \|f\|_{H^{1/2}_{00}}\cdot \|g\|_{H^{1/2}}.
	\end{align}
\end{lemma}

\noindent{\bf Proof.} We define the bilinear map
\begin{align*}
\mathscr{B}(f,g) = \int_0^\infty f' g dy = -\int_0^\infty f g' dy\quad \forall \, f,\, g\in C^{\infty}_c(\R_+).
\end{align*}
Applying the interpolation result in Sections 11 $\&$ 12 of \cite{Lions}, then the bilinear map
\begin{align*}
\mathscr{B}(f,g) : C^\infty_c((0,\infty))  \times C^\infty_c((0,\infty)) \to \R,
\end{align*}
can be uniquely extended to a bounded bilinear map
\begin{align*}
\mathscr{B}(\cdot,\cdot) : H^{1/2}_{00}(\R_+) \times H^{1/2}(\R_+) \to \R.
\end{align*}
Therefore the proof of Lemma \ref{lemA.9} is completed. $\hfill\Box$

\begin{lemma}\label{lemA.01}
	Let $f\in H^1(\R_+\times [0,b])$, it  hods that
	\begin{align}\label{Z.11-2}
		\|f(\cdot,0)\|^2_{H^{1/2}} 
		&\lesssim \|(f_x, f)\| \big\{ \|f_y\| + b^{-1}\|f\|\big\}.
	\end{align}
\end{lemma}

\noindent{\bf Proof.} 
Firstly we consider the case $h(x,y)\in H^1(\mathbb{R}\times[0,b])$ We note that
\begin{align}\label{Z.10}
	\int_{\mathbb{R}}  (1 + |\xi|) |\hat{h}(\xi,0)|^2 d\xi = \int_{\mathbb{R}}(1 + |\xi|) |\hat{h}(\xi,y)|^2 d\xi - \int_0^y \int_{\mathbb{R}} 2 (1 + |\xi|) \hat{h}(\xi,z) \hat{h}_y(\xi,z)  d\xi dz.
\end{align}
Integrating \eqref{Z.10} over $y\in [0,b]$, one gets
\begin{align}\label{Z.11-1}
	\|h(\cdot,0)\|^2_{H^{1/2}(\mathbb{R})} 
	&\lesssim \|(h_x, h)\| \big\{ \|h_y\| + \f{1}{b}\|h\|\big\}.
\end{align}

Let $\tilde{f}$ be the even extension of $f$ in $x$, then we know that $\tilde{f}\in H^1(\R\times [0,b])$. Applying \eqref{Z.11-1} to $\tilde{f}$, one concludes \eqref{Z.11-2}.
Therefore the proof of Lemma \ref{lemA.01} is completed. $\hfill\Box$


\begin{lemma}\label{lemAH.1}
Let $f(x,y) \in H^1(\R_+\times [0,b])$ and $f|_{x=0}=0$, then it holds that $f(\cdot, 0) \in H^{1/2}_{00}(\R_+)$ with
\begin{align}\label{AH.0}
\|f(\cdot, 0)\|^2_{H^{1/2}_{00}} \lesssim \|(f_x, f)\| \big\{b^{-1}\|f\| + \|f_y\|\big\}.
\end{align}
\end{lemma}

\noindent{\bf Proof.} Firstly it follows from Lemma \ref{lemA.01} that 
\begin{align}\label{AH.2}
	\|f(\cdot,0)\|^2_{H^{1/2}}\lesssim \|(f_x, f)\| \big\{ b^{-1} \|f\| + \|f_y\|\big\}.
\end{align}

On the other hand, noting  $f|_{x=0}=0$, we apply the Hardy inequality to get
\begin{align}\label{AH.1}
\|x^{-1}f\| \lesssim \|f_x\|.
\end{align}
Then it is clear to know that  
\begin{align}\label{AH.3}
\int_0^\infty \f{1}{x^{\ast}} |f(x,0)|^2
&\leq b^{-1} \int_0^b\int_0^\infty (x^{\ast})^{-1} |f(x,y)|^2dxdy + \int_0^b\int_0^\infty 2(x^{\ast})^{-1} |f(x,y) f_y(x,y)| dxdy \nonumber\\
&\lesssim \|(x^{\ast})^{-1}f\| \big\{ b^{-1} \|f\| +  \cdot \|f_y\|\big\} \lesssim \big\{\|f\| + \|x^{-1}f\|\big\} \big\{b^{-1} \|f\| + \|f_y\|\big\} \nonumber\\
&\lesssim \|(f_x, f)\| \big\{b^{-1}\|f\| + \|f_y\|\big\},
\end{align}
where we have used \eqref{AH.1}.

Finally, by noting the definition of $H^{1/2}_{00}(\R_+)$ in \eqref{D.6}, we conclude  \eqref{AH.0} from \eqref{AH.2} and \eqref{AH.3}. Therefore the proof of Lemma \ref{lemAH.1} is completed. $\hfill\Box$


\begin{lemma}\label{lemB.5}
	It holds that 
	\begin{align}\label{B.14}
		\|f\|  
		&\lesssim 
		a^{-1} \|u_sf\| + \|f_{y}\|.
	\end{align}
\end{lemma}

\noindent{\bf Proof.} Let $y_1\in [1,2]$, $y\in [0,1]$, it is clear that 
\begin{align*}
	f(x,y) = f(x,y_1) -  \int_y^{y_1} f_{y}(x,s) ds,
\end{align*}
which implies that 
\begin{align*}
	\int_0^L \int_0^1 |f(x,y)|^2 dy \lesssim \int_0^L |f(x,y_1)|^2 dx + \int_0^L \int_0^2 |f_{y}(x,s)|^2 ds dx.
\end{align*}
Integrating above inequality over $y_1\in [1,2]$, one gets that 
\begin{align*}
	\|f \mathbf{1}_{\{y\in[0,1]\}}\| \lesssim \|f\mathbf{1}_{\{y\in[1,2]\}}\| + \|f_{y} \mathbf{1}_{\{y\in[0,2]\}}\|,
\end{align*}
which yields that 
\begin{align*} 
	\|f\|  &\leq  \|f\mathbf{1}_{\{y\geq1\}}\| + \|f\mathbf{1}_{\{y\in[0,1]\}}\|  \lesssim  \|f\mathbf{1}_{\{y\geq1\}}\| + \|f_{y} \mathbf{1}_{\{y\in[0,2]\}}\| \lesssim  a^{-1} \|u_sf\| + \|f_{y}\|.
\end{align*}
Therefore the proof of Lemma \ref{lemB.5} is completed. $\hfill\Box$

\begin{lemma}\label{lemA.8}
	Let $f$ be a function with $\dis \int_0^\infty y|f(y)|^2dy <\infty$, and denote $f^{\d}(y):=f(y+\d)$. Then it holds that $f^{\d} \in L^2(\R_+)$ with
	\begin{align*}
		\int_0^{\infty} (y+\d) |f^{\d}(y)|^2<\int_0^\infty y|f(y)|^2 dy <\infty, \quad \forall \,\, \d\in [0,1],\\
		\|\sqrt{y}(f^{\d}-f)\|_{L^2}\to 0,\quad \mbox{as}\,\,\, \delta\to0.
	\end{align*}
\end{lemma}

\noindent{\bf Proof.} The proof is direct, we omit the details for brevity. $\hfill\Box$

\subsection{Div-curl inversion lemma} 
To recover the momentum equations \eqref{2.109} from the div-curl equations \eqref{1.33}-\eqref{1.34},
we need the following generalized div-curl inversion lemma.
\begin{lemma}\label{lemR.8}
	Let $\mathscr{A}_1,\mathscr{A}_2 \in H^1$  satisfy the following generalized div-curl system
	\begin{align}\label{B.28-1}
		\begin{cases}
			\dis	\pa_x (\v \mathscr{A}_1) + \pa_y \mathscr{A}_2 =0,\\
			\dis	\pa_y (b_1 \mathscr{A}_1) - \pa_x (b_2 \mathscr{A}_2)=0,\\
			\dis	\mathscr{A}_1|_{x=0}=\mathscr{A}_2|_{x=L}=\mathscr{A}_1|_{y=0}=0.
		\end{cases}
	\end{align}
	Then it holds  $(\mathscr{A}_1,\mathscr{A}_2)=(0,0)$.
\end{lemma}
\noindent{\bf Proof.}  Noting $\eqref{B.28-1} $, we introduce the stream function $\psi$ such that
\begin{align}\label{B.6}
	\v \mathscr{A}_1= \psi_y,\qquad \mathscr{A}_2=-\psi_x,
\end{align}
with 
\begin{align}\label{B.28-2}
	\psi_y \big|_{x=0}=0,\quad \psi_y\big|_{y=0}  =0,\quad \psi_x\big|_{x=L}=0. 
\end{align}
Due to  $\psi_y |_{x=0}=0$,  we take 
\begin{align}\label{B.28-3}
	\psi(0,y)=0.
\end{align}
Then it holds that 
\begin{align*}
	\psi(x,y)=\int_0^x \psi_x(s,y)ds\quad \Longrightarrow\quad \psi \in L^2([0,L]\times \mathbb{R}_+).
\end{align*}
Thus we have $\psi\in H^2([0,L]\times \mathbb{R}_+)$.

We  rewrite $\eqref{B.28-1}_2$ as
\begin{align}\label{R.33-0}
	(b_1 \psi_y)_y + \v (b_2 \psi_x)_x=0.
\end{align}
Multiplying \eqref{R.33-0}  by  $\psi$, one obtains 
\begin{align}\label{B.33-1}
	\int_0^L\int_{0}^\infty b_1 |\psi_y|^2 + \v b_2 |\psi_x|^2 dydx =\int_0^L b_1 \psi_y \psi dx \Big|_{y=0} + \v \int_0^L b_2 \psi_x \psi d y \Big|_{x=0}^{x=L} =0,
\end{align}
where we have used the boundary conditions \eqref{B.28-2}-\eqref{B.28-3}. Hence it follows from \eqref{B.33-1} and \eqref{B.6} that 
\begin{align*}
	\psi\equiv 0\quad \Longrightarrow\quad   (\mathscr{A}_1,\mathscr{A}_2)=(0,0).
\end{align*}
Therefore the proof of Lemma \ref{lemR.8} is completed. $\hfill\Box$

\subsection{Function spaces, norms and notations}\label{sec3.1}
For later use,  we define the weight functions
\begin{align}\label{6.0-0}
	w=\mbox{either}\, 1 \, \mbox{or}\,\, w_0 \quad \mbox{and} \quad  \hat{w}:=1+w,\quad \hat{w}^{\v}:=1+\sqrt{\v}w,
\end{align}
with $w_0=y(1+Y)^{l}=y (1 + \sqrt{\v}y)^{l}$ for $\mathfrak{l}_0\gg l\gg1$. It is clear that 
\begin{align}\label{6.0}
	|\pa_yw_0|\lesssim \sqrt{\v}w_0 + 1\equiv \hat{w}^{\v}_0.
\end{align}
Throughout this paper, unless stated otherwise, we always take the weight function as in \eqref{6.0-0}.

We fix a cut-off function
\begin{align}\label{6.0.1}
	\chi(y):=
	\begin{cases}
		1,\quad |y|\in [0,1),\\
		0,\quad |y|\in (2,\infty),
	\end{cases}
	\quad \mbox{and}\quad \bar{\chi}(y):=1-\chi(y),
\end{align}
with $\chi'(y)\leq 0$ for $y\geq0$. We also denote 
\begin{align}
\mathbf{1}_{\bf N}:=
\begin{cases}
	1,\quad \mbox{for NT},\\
	0,\quad \mbox{for DT},
\end{cases}
\quad \mbox{and}\quad 
\mathbf{1}_{\bf D}:=
\begin{cases}
	1,\quad \mbox{for DT},\\
	0,\quad \mbox{for NT}.
\end{cases}
\end{align}

\medskip

In this paper, we always assume $1\gg L \gg \v >0$. From now on, we use the following  notations 
\begin{align}\label{A.1}
	\begin{split}
		\||q|\|_{w}&:= \|\sqrt{u_s}(q_{xyy},\sqrt{\v}q_{xxy}, \v q_{xxx})w\| + L^{-\f18}\|\sqrt{u_s}q_{yyy} w\|,\\
		\|||\phi||\|_{w}&:= \|(\sqrt{\v}\phi_{xyyy},\v \phi_{xxyy}, \v^{\f32}\phi_{xxxy}, \v^2\phi_{xxxx}, \phi_{yyyy})w\|.
	\end{split}
\end{align}
We define
\begin{align}\label{6.1}
	\begin{split}
		[\phi]_{3,w}&:=\||q|\|_{w}^2 + L^{-\f14}\|u_s\nabla_{\v}q_xw\|^2 + |q|_{\pa,2,w}^2,\\ 
		[\phi]_{4,w}&:=\|||\phi||\|_{w}^2 + \|\sqrt{u_s}\big(\v^{\f32}\phi_{xxx},\sqrt{\v}\phi_{xyy}\big)w\|_{x=0}^2 + \|\sqrt{u_s} \v \phi_{xxy}w\|_{x=L}^2,
	\end{split}
\end{align}
where $\dis q= \phi/\bar{u}_s$ and 
\begin{align}\label{6.2}
	\begin{split}
		|q|_{\pa,2,w}&:=\|u_s q_{xy} w\|_{x=0} + \|\sqrt{\v}u_s q_{xx}w\|_{x=L}  + \|\sqrt{u_{sy}}q_{xy}\|_{y=0}  \\
		&\quad  + L^{-\f18}\|\sqrt{u_{sy}}q_{yy}\|_{y=0}  + L^{-\f18}\|u_sq_{yy}w\|_{x=L}.
	\end{split}
\end{align}

For the pressure $p$, we define
\begin{align} \label{6.3-1}
\begin{split}
[[p]]_{2,w}
:&=\|(p_{y},\sqrt{\v}p_{x})w\|^2  + \v \|\frac{B}{\sqrt{u_s}}w\|^2_{x=0} + \mathfrak{K}_{2,w}(p),\\
[[p]]_{3,w}
:&= \|u_s \nabla^2_{\v}p_y\, w\|^2 + \v \|\sqrt{u_s^3}\Delta_{\v}p_yw\|_{x=L}^2 +\|\sqrt{u_s}p_{yy}w\|_{x=L}^2 + \|\sqrt{u_s}\sqrt{\v}p_{xy}w\|_{x=0}^2, 
\end{split}
\end{align}
where 
\begin{align}\label{6.3-2}
\mathfrak{K}_{2,w}(p):&=\|\nabla^2_{\v}p\,w\|^2 +  \v \|\sqrt{u_s}\Delta_{\v}pw\|_{x=L}^2  + \v \|\sqrt{u_s}\sqrt{\v}p_{xy}w\|^2_{x=L} + \v\|\sqrt{u_s}p_{yy}w\|^2_{x=0}.
\end{align}

For the pseudo entropy, we define
\begin{align}\label{3.0-0S}
\begin{split}
[[[S]]]_{2,w}
:&=\|(\sqrt{\v}S_{xx}, S_{xy})w\|^2 + L^{-\f12}\|S_{yy}w\|^2 + L^{-\f12}\|u_s S_xw\|^2 +   L^{-2}\|S_yw\|^2 \\
&\quad  + \|\sqrt{u_s}S_xw\|^2_{x=0, L} 
+  L^{-1}\|S_yw\|_{x=0,L}^2,\\
[[[S]]]_{3,w}
:&=\|u_s(\v S_{xxy}, \sqrt{\v} S_{xyy},   S_{yyy})w\|^2,\\
[[[S]]]_{b,w}
:&= \|\sqrt{u_s}\sqrt{\v} S_{xy}w\|^2_{x=0} + \mathbf{1_{N}}\|\sqrt{u_s} S_{yy}w\|^2_{x=L} + \mathbf{1_{D}}\|S_{yy}w\|^2_{x=L}.
\end{split}
\end{align}


For later use, we also denote the combination norm
\begin{align}\label{3.0-0}
	\|(\phi, p, S)w\|^2_{\mathbf{X}^{\v}}
	:&=[\phi]_{3,w} + [\phi]_{4,w} + \f{1}{\v} [[p]]_{2,w} + [[p]]_{3,w} + \|u_s \sqrt{\v}\Delta_{\v}p_x w\|^2   \nonumber\\
	&\quad + [[[S]]]_{2,w} + [[[S]]]_{b,w} + [[[S]]]_{3,w}.
\end{align}

As pointed out previous, it is hard to control $\|\nabla_{\v}p w\|$ for the approximate solutions, so we introduce the following notation (see \eqref{3.63-1} for the definition of $\mathbf{d}_2$)
\begin{align}\label{3.0-1}
	\begin{split}
		\mathscr{K}(p):&=\|(\mathbf{d}_2 p_x, p_x \mathbf{1}_{\{y\in[0,2]\}})\|^2 + \v^{-1+\beta}\|(\mathbf{d}_2p_y, p_y \mathbf{1}_{\{y\in[0,2]\}})\|^2,\quad \mbox{for}\,\beta \in (\f12,1),
	\end{split}
\end{align} 
which will help us to deal with the first order derivative of pressure.
For the higher order derivative of pressure in the procedure of $\delta$-approximation, we also need the following $\delta$-dependent  notations 
\begin{align}\label{3.0}
	\begin{split}
		\mathcal{K}^{\delta}_{2,w}(p):&= \| \nabla_{\v}^2 p w\|^2   + \v  \|\sqrt{u_s^{\delta}}\Delta_{\v}p w\|^2_{x=L} + \v \|\sqrt{u_s^{\delta}}\sqrt{\v}p_{xy} w\|^2_{x=L},\\
		\mathcal{K}^{\delta}_{3,w}(p):&=\|u_s^{\delta}\nabla^{2}_{\v}p_{y}w\|^2 + \|\sqrt{u_s^{\delta}}p_{yy}w\|^2_{x=L}  + \v \|\sqrt{u_s^{\delta}}p_{xy}w\|^2_{x=0} + \v \|u^{\delta}_s\sqrt{u_s^{\delta}} \Delta_{\v}p_y w\|^2_{x=L},\\
	 \bar{\mathcal{K}}^{\delta}_{3,w}(p) :&= \|u^{\delta}_s \sqrt{\v}\Delta_{\v}p_x w\|^2,
	\end{split}
\end{align}
where $u_s^{\delta}:=\delta + u_s$.
We also define the $\delta$-dependent combination norm
\begin{align}\label{3.0-2}
	\|(\phi, p, S)w\|^2_{\mathbf{Y}^{\v}_{\delta}}
	:&=[\phi]_{3,w} + [\phi]_{4,w} + \mathscr{K}(p)  + \f{1}{\v} \mathcal{K}^{\delta}_{2,w}(p) + \mathcal{K}^{\delta}_{3,w}(p) +  \bar{\mathcal{K}}^{\delta}_{3,w}(p)  \nonumber\\
	&\quad + [[[S]]]_{2,w} + [[[S]]]_{b,w} + [[[S]]]_{3,w}.
\end{align}
If $\delta=0$ in \eqref{3.0}-\eqref{3.0-2}, we shall simply denote them as $\mathcal{K}_{2,w}$, $\mathcal{K}_{3,w}$, $\bar{\mathcal{K}}_{3,w}$ and $\|(\cdots)w\|^2_{\mathbf{Y}^{\v}}$, respectively.

\subsection{Some estimates on $\phi$} 
For later use, we first recall a useful Hardy-type inequality.
\begin{lemma}[Guo-Iyer\cite{Guo-Iyer-2024}]\label{lemA.0}
	Let $w=1$ or $w_0$. We assume that $f$ satisfy $f|_{y=0}=0$ and $f\rightarrow 0$ as $y\rightarrow \infty$, then it holds that 
	\begin{align}\label{A.5}
		\begin{split}
			\|\frac{f}{y}w\|&\lesssim \|f_yw\|+ \|\sqrt{\v}fw\|,\\
			\|\frac{f}{y}w\|&\lesssim \|f_yw\|+ \|fw_y\|.
		\end{split}
	\end{align}
\end{lemma}

Using \eqref{pair1}, one has 
\begin{align}\label{A.6}
	\pa_x q:=\frac{\bar{u}_s \phi_x- \bar{u}_{sx}\phi}{\bar{u}_s^2} \quad \Longrightarrow	\quad \pa_x q|_{x=L}=-\frac{\bar{u}_{sx}}{\bar{u}_s} q |_{x=L}.
\end{align}

\begin{lemma}[Guo-Iyer \cite{Guo-Iyer-2024}]\label{lemA.2}
	Let $L\ll 1$, it holds that
	\begin{align}
		&\|\pa_y^{j}\pa_x^k \phi\, w\|\lesssim L \|\pa_y^{j}\pa_x^{k+1}\phi\, w\|,\quad \mbox{for}\,\, k=0,1,2,3,\label{A.7}\\
		&\|u_s^k \pa_y^j q\, w\|\lesssim L \|u_s^k \pa_y^j \pa_x q w\|,\quad \mbox{for}\,\, k\geq0,\label{A.8}\\
		&\|q_x w\|\lesssim L \|q_{xx}w\|.\label{A.9}
	\end{align}
\end{lemma}

\begin{lemma}\label{lemA.3}
	Let $L\ll a<1$ and $\alpha\in (0,\f1{12}]$, then it holds that 
	\begin{align}
		&\|q_{yy}w\| 
		\lesssim L^{\f12-} [\phi]_{3,\hat{w}}^{\f12},\label{A.10}\\
		&\|\nabla_{\v}q_xw\| + \|\frac{q_x}{y}w\| 
		\lesssim  L^{\f{1}{2}\alpha-}[\phi]_{3,\hat{w}}^{\f12},\label{A.11}\\
		&\|\nabla_{\v}q w\| 
		\lesssim L[\phi]_{3,\hat{w}}^{\f12}\label{A.12} 
	\end{align}
	and 
	\begin{align}
		&\|\phi_{yyy}w\| + \|\sqrt{\v}\phi_{xx}w\|
		\lesssim L^{\f18} [\phi]_{3,\hat{w}}^{\f12},\label{A.12-2}\\
		&\|(\phi_{xyy},\sqrt{\v}\phi_{xxy}, \v\phi_{xxx})w\|
		\lesssim \sqrt{a} [\phi]_{3,\hat{w}}^{\f12},\label{A.12-4}\\
		&\|\phi_{yy}w\| + \|\nabla_{\v}\phi w\| \lesssim   L[\phi]_{3,\hat{w}}^{\f12},\label{A.12-5}\\
		&\|\phi_{xy}w\| \lesssim L^{\f12 \alpha-} [\phi]_{3,\hat{w}}^{\f12}.\label{A.12-6}
	\end{align}
	We remark that the above constants are independent of $\alpha\in(0,\f1{12}]$.
\end{lemma}

\noindent{\bf Proof.} The idea of proof is similar as \cite{Guo-Iyer-2024}  but with more delicate  coefficients. We divide it into two steps  since it is very long.

\noindent{\it Step 1.} It is clear that 
\begin{align}\label{A.13}
\|q_{yy}w\| \leq \|q_{yy}w \chi(y)\| + \|q_{yy}w \bar{\chi}(y)\|.
\end{align}
Integrating by parts, we have 
\begin{align}\label{A.14}
	&\|q_{yy}w \chi \|^2
	=\int_0^L\int_0^\infty \pa_yy\, q_{yy}^2 w^2 \,\chi^2(y) \, dydx\nonumber\\
	&=-2\int_0^L\int_0^\infty  y \Big[q_{yy} q_{yyy} w^2 \,\chi^2  + q_{yy}^2 w w_y \,\chi^2  + q_{yy}^2 w^2 \chi' \chi \Big]  \, dydx\nonumber\\
	&\lesssim \|\sqrt{y}q_{yy}w\chi\|\cdot \|\sqrt{y}q_{yyy}w\chi\| 
	+  \|\sqrt{y}q_{yy}w\, \mathbf{1}_{\{1\leq y\leq 2\}}\|^2   +  \|\sqrt{y}q_{yy}w_y\chi\|^2\nonumber\\
	&\lesssim \frac{L}{a}  \|\sqrt{u_s}q_{xyy}w\chi\|\cdot \|\sqrt{u_s}q_{yyy}w\chi\| + L^2 \frac{1}{a}  \|\sqrt{u_s}q_{xyy}\hat{w}\chi\|^2  
	 \lesssim  \frac{L}{a} \||q|\|_{\hat{w}}^2.
\end{align}
A direct calculation shows that 
\begin{align*}
	\|q_{yy}w \bar{\chi}(y)\|
	&\lesssim L \|q_{xyy}w \bar{\chi}(y)\|
	\lesssim \frac{L}{\sqrt{a}}\||q|\|_{w}, 
\end{align*}
which, together with \eqref{A.13}-\eqref{A.14}, yields that 
\begin{align}\label{A.15}
	\|q_{yy}w\|  \lesssim  \frac{\sqrt{L}}{\sqrt{a}} \||q|\|_{\hat{w}} \lesssim \frac{\sqrt{L}}{\sqrt{a}} [\phi]_{3,\hat{w}}^{\f12}.
\end{align}

\smallskip

For later use, we consider a general inequality. Let $g$ be a generic function now. It is clear that 
\begin{align}\label{A.19}
	\|gw\bar{\chi}(\frac{ay}{\xi})\| 
	\lesssim \frac{1}{\xi} \|u_sgw \, \mathbf{1}_{\{ay\geq \xi\}}\|.
\end{align}
where $\xi\in (0,1]$ is a small parameter determined later. 
Integrating by parts, one has
\begin{align*}
	&\|gw\chi(\frac{ay}{\xi})\|^2=  \int_0^L\int_0^\infty \pa_yy\, g^2 w^2 \chi^2(\frac{ay}{\xi}) \, dydx\nonumber\\
	&=-2\int_0^L\int_0^\infty  y \Big[g g_y w^2 \,\chi^2(\frac{ay}{\xi}) + g^2 w w_y \,\chi^2(\frac{ay}{\xi}) + g^2 w^2 \,\frac{a}{\xi}\chi'(\frac{ay}{\xi})\chi(\frac{ay}{\xi})\Big]  \, dydx\nonumber\\
	&\lesssim \|gw\chi\| \Big(\|yg_y w\chi\| + \|ygw_y\chi\| 
	+ \frac{1}{\xi}\|aygw\, \mathbf{1}_{\{\xi\leq ay \leq 2\xi\}}\|\Big),
\end{align*}
which yields that 
\begin{align}\label{A.20}
	\|gw\chi(\frac{ay}{\xi})\| 
	&\lesssim \frac{\sqrt{\xi}}{a} \|\sqrt{u_s} g_yw\chi\| + \big(\frac{1}{a}+\f{1}{\xi}\big) \|u_s g \hat{w}\chi\|,
\end{align}
Combining \eqref{A.19} and \eqref{A.20}, one has that 
\begin{align}\label{A.21}
	\|gw\|&\lesssim  \frac{\sqrt{\xi}}{a}  \|\sqrt{u_s} g_yw\chi\| +  \big(\frac{1}{a}+\f{1}{\xi}\big)   \|u_s g \hat{w}\|.
\end{align}

Taking $g=q_{xy},\sqrt{\v}q_{xx}$ in \eqref{A.21}, one can obtain
\begin{align}
	\|q_{xy}w\|
	&\lesssim  \frac{\sqrt{\xi}}{a}  \||q|\|_{w} + \big(\frac{1}{a} + \f{1}{\xi}\big)   \|u_s q_{xy} \hat{w}\| \lesssim \frac{L^{\f{\alpha}{2}}}{a}[\phi]_{3,\hat{w}}^{\f12},\label{A.23} \\
	\|\sqrt{\v}q_{xx}w\|
	&\lesssim \frac{\sqrt{\xi}}{a}  \||q|\|_{w} +  \big(\frac{1}{a}+\f{1}{\xi}\big)  \|u_s \sqrt{\v}q_{xx} \hat{w}\|\lesssim \frac{L^{\f{\alpha}{2}}}{a}[\phi]_{3,\hat{w}}^{\f12}.\label{A.24}
\end{align}
where we have taken $\xi=L^{\alpha}$ in above two estimates. 

It follows from \eqref{A.8}-\eqref{A.9} and \eqref{A.23}-\eqref{A.24}  that 
\begin{align}\label{A.24-1}
	\|(q_y,\sqrt{\v}q_x)w\|\lesssim L \|(q_{xy}, \sqrt{\v}q_{xx})w\|
	\lesssim L [\phi]_{3,\hat{w}}^{\f12}.
\end{align}
Also we have from $\eqref{A.5}_1$ and \eqref{A.23} that
\begin{align}\label{A.27}
	\|\frac{q_x}{y}w\|&\lesssim \|q_{xy}w\|+ \|\sqrt{\v}q_x w\| 
	\lesssim \frac{L^{\f{\alpha}{2}}}{a}[\phi]_{3,\hat{w}}^{\f12}.
\end{align}

\smallskip

\noindent{\it Step 2.} A direct calculation shows that 
\begin{align}\label{A.34}
	\begin{split}
		&\phi_x=\bar{u}_s q_x + \bar{u}_{sx}q,\qquad \phi_y = \bar{u}_sq_y + \bar{u}_{sy} q,\\
		&\phi_{yy}=\bar{u}_{s} q_{yy} +2  \bar{u}_{sy}q_y+  \bar{u}_{syy} q, \\
		&\phi_{xx}=\bar{u}_{s} q_{xx} + 2\bar{u}_{sx}q_x + \bar{u}_{sxx} q,\\
		&\phi_{xy}=\bar{u}_{s} q_{xy} + \bar{u}_{sx}q_y+ \bar{u}_{sy} q_x + \bar{u}_{sxy} q,\\
		&\phi_{yyy}=\bar{u}_{s} q_{yyy} + 3\bar{u}_{sy}q_{yy} + 3 \bar{u}_{syy} q_y + \bar{u}_{syyy} q,\\
		&\phi_{xxx}=\bar{u}_{s} q_{xxx} + 3\bar{u}_{sx}q_{xx} + 3 \bar{u}_{sxx} q_x + \bar{u}_{sxxx} q,\\
		&\phi_{xyy}=\bar{u}_{s} q_{xyy} + \bar{u}_{sxyy} q + \bar{u}_{syy}q_{x} + \bar{u}_{sx} q_{yy} +2 \bar{u}_{sxy} q_y + 2\bar{u}_{sy} q_{xy} ,\\
		&\phi_{xxy}=\bar{u}_{s} q_{xxy} + \bar{u}_{sxxy} q + \bar{u}_{sxx}q_{y} + \bar{u}_{sy} q_{xx} +2 \bar{u}_{sxy} q_x + 2\bar{u}_{sx} q_{xy}.
	\end{split}
\end{align}
For $\phi_{yyy}$, using \eqref{A.15} and \eqref{A.24-1}, we have
\begin{align}\label{A.38}
	\begin{split}
		&\|\bar{u}q_{yyy}w\|\lesssim \sqrt{a}  \|\sqrt{u_s}q_{yyy}w\| \lesssim L^{\f18}\sqrt{a}\||q|\|_{w},\\
		&\|\bar{u}_{sy}q_{yy}w\|\lesssim a\|q_{yy}w\|\lesssim \sqrt{aL}\||q|\|_{\hat{w}},\\
		&\|\bar{u}_{syy}q_{y}w\|\lesssim a\|q_{y}w\| \lesssim  L [\phi]_{3,\hat{w}}^{\f12}.
	\end{split}
\end{align}
For term involving $q$, it follows from \eqref{A.24-1} and \eqref{A.27} that 
\begin{align}\label{A.39}
	\|\bar{u}_{syyy}qw\|
	&\lesssim L \| \frac{q_x}{y} w\| + L \|\sqrt{\v}q_xw\| \lesssim L [\phi]_{3,\hat{w}}^{\f12},
\end{align}
which, together with \eqref{A.39}, yields that 
\begin{align}\label{A.40}
	\|\phi_{yyy}w\|
	\lesssim L^{\f18}  [\phi]_{3,\hat{w}}^{\f12}.
\end{align}

For $\phi_{xyy}$, using \eqref{A.15}, \eqref{A.23} and \eqref{A.24-1}, one can obtain 
\begin{align}\label{A.41}
	\begin{split}
		&\|\bar{u}_s q_{xyy}w\| \lesssim \sqrt{a} \||q|\|_{w},\\
		&\|\bar{u}_{sx}q_{yy}w\|\lesssim a \|q_{yy}w\|\lesssim \sqrt{aL} [\phi]_{3,\hat{w}}^{\f12},\\
		&\|\bar{u}_{sy}q_{xy}w\|\lesssim a \|q_{xy}w\|
		\lesssim L^{\f12\alpha}  [\phi]_{3,\hat{w}}^{\f12},\\
		&\|(\bar{u}_{syy}q_x, \bar{u}_{sxy}q_y)w\|  \lesssim a\|\frac{q_{x}}{y}w\| + \|(q_y, \sqrt{\v}q_x) w\| \lesssim L^{\f12\alpha}  [\phi]_{3,\hat{w}}^{\f12},\\
		&\|\bar{u}_{sxyy}qw\|\lesssim L\|\frac{q_x}{y}w\| + L \|\sqrt{\v}q_xw\| \lesssim L [\phi]_{3,\hat{w}}^{\f12},
	\end{split}
\end{align}
which  immediately  yields that
\begin{align}\label{A.44}
	\|\phi_{xyy}w\|\lesssim \sqrt{a} [\phi]_{3,\hat{w}}^{\f12}.
\end{align}

For $(\sqrt{\v}\phi_{xxy}, \v \phi_{xxx})$, using \eqref{A.23}-\eqref{A.24} and \eqref{A.24-1}, one can similarly get
\begin{align}\label{A.45}
	\|(\sqrt{\v}\phi_{xxy}, \v \phi_{xxx})w\|\lesssim  \sqrt{a} [\phi]_{3,\hat{w}}^{\f12}.
\end{align}


For $\phi_{yy}$, it follows from \eqref{A.44} that 
\begin{align}\label{A.47}
	\|\phi_{yy}w\|\lesssim L \|\phi_{xyy}w\|\lesssim L[\phi]_{3,\hat{w}}^{\f12}.
\end{align}
For $\phi_{xy},\phi_{xx}$, we have from \eqref{A.24-1} and \eqref{A.27} that
\begin{align}\label{A.48}
\begin{split}
\|\phi_{xy}w\|
&\lesssim  \|u_{s}q_{xy}w\| +  \|\nabla_{\v}qw\| + a \|\frac{q_x}{y}w\|\lesssim L^{\f12 \alpha} [\phi]_{3,\hat{w}}^{\f12},\\
\|\sqrt{\v}\phi_{xx}w\|&\lesssim \|\sqrt{\v}u_s q_{xx}w\| +  \|\sqrt{\v}q_xw\| \lesssim L^{\f18} [\phi]_{3,\hat{w}}^{\f12}. 
\end{split}
\end{align}
For $\phi_y, \phi_x$, it follows from \eqref{A.48}  that 
\begin{align*}
	\|(\phi_y, \sqrt{\v}\phi_x)w\|&\lesssim L\|(\phi_{xy}, \sqrt{\v}\phi_{xx})w\| \lesssim L  [\phi]_{3,\hat{w}}^{\f12}.
\end{align*}
Therefore the proof of Lemma \ref{lemA.3} is completed. $\hfill\Box$

\begin{lemma}[Trace estimate]\label{lemA.4}
	It holds that 
	\begin{align}
		\|(q_x,q_y) w\|_{x=L}+ \|\phi_{yy}w\|_{x=L}+ \|\v\sqrt{u_s}q_{xx}w\|_{x=L} 
		+ \|q_{xy}w\|_{x=L}
		&\lesssim L^{\f12-}[\phi]_{3,\hat{w}}^{\f12},\label{A.50}\\
		\|\sqrt{\v}q_x w\|_{x=0}+ \|\sqrt{\v}q_{xx} w\|_{x=0}
		&\lesssim L^{\f12-}[\phi]_{3,\hat{w}}^{\f12},\label{A.50-1}\\
		\|\v^{\f14}\frac{q_x}{\la y\ra} w\|_{x=0}+ \|\v^{\f14}q_{xx} w\|_{x=0}
		&\lesssim L^{\f12\alpha-}[\phi]_{3,\hat{w}}^{\f12},\label{A.50-2}\\
		\|\sqrt{u_s}q_{xyy}w\|_{x=L}
		&\lesssim L^{\f12-}[\phi]_{3,\hat{w}}^{\f12}.\label{A.50-3}
	\end{align}
\end{lemma}
\noindent{\bf Proof.} The lemma can be proved by similar arguments as in \cite{Guo-Iyer-2024} and Lemma \ref{lemA.3}, we omit the details of proof for simplicity of presentation. $\hfill\Box$

\medskip

Next we present a new improved trace estimate which is important for the uniform estimate of pseudo entropy of DT case  in the Section \ref{SecD7}.
\begin{lemma}\label{lemAD.1}
Let $\alpha\in (0,\f1{12}]$. It holds that 
	\begin{align}\label{AD.1}
		\|q_{yy}w\|^2_{x=L} + \|\sqrt{\v}q_{xy}w\|^2_{x=0}&\lesssim L^{\f14 \alpha-} [\phi]_{3,\hat{w}} + L^{\f14 \alpha-} [\phi]_{4,1}.
	\end{align}
\end{lemma}

\noindent{\bf Proof.} Let $\g \in (0,\f12)$ be a suitably small positive constant determined later. We take $\xi=L^{\alpha}$.

 It is noted that 
\begin{align}\label{D3.27}
	&	\|q_{yy}w\|^2_{x=L}
=2\int_0^L\int_0^\infty q_{yy} q_{xyy} \big(w\chi(\f{y}{\xi})\big)^2  dydx  + O(1)L^{\f12} [\phi]_{3,w} \nonumber\\
	&= 2\int_0^L\int_0^\infty y^{-1+\g}\big\{q_{yy}(x,y)-q_{yy}(x,0)\big\}\cdot  y^{1-\g} q_{xyy}  \big(w\chi(\f{y}{\xi})\big)^2 dydx \nonumber\\
	&\quad + 2\int_0^L\int_0^\infty q_{yy}(x,0)  q_{xyy}  \big(w\chi(\f{y}{\xi})\big)^2 dydx + O(1)L^{\f12} [\phi]_{3,w}\nonumber\\
	&\lesssim  L^{\f12} [\phi]_{3,w} + C_{\g} \|y^{1-\g}q_{xyy} \chi(\f{y}{\xi})\| \big\{ \|y^{\g} q_{yyy}\chi(\f{y}{\xi})\| + L^{-\a}\|q_{yy}\| + L^{-\a}\|q_{yy}\|_{y=0}\big\}  \nonumber\\
	&\quad  +  \Big|\int_0^L \big(q_{yy} q_{xy}w\big)(x,0) dx\Big| + \Big|\int_0^L\int_0^\infty \big(\chi^2(\f{y}{\xi})w^2\big)_y  q_{yy}(x,0)  q_{xy}(x,y) dydx\Big| \nonumber\\
	&\lesssim  L^{\f1{4}\alpha} [\phi]_{3,\hat{w}} + L^{\f14\alpha} \|\sqrt{y}q_{xyy} \chi(\f{y}{\xi})\|\cdot \|y^{\gamma} q_{yyy}\chi(\f{y}{\xi})\|.
\end{align}
For $\|y^{\g} q_{yyy}\chi(\f{y}{\xi})\|$ with $\gamma>0$ suitably small, we have 
\begin{align}\label{D3.28}
	\|y^{\g}\phi_{yyyy}\chi(\f{y}{\xi})\|^2
	&\geq \f34 \|y^{\g} \chi(\f{y}{\xi}) \big[u_s q_{yyyy} + 4 u_{sy}q_{yyy}\big]\|^2 - L^{\f12}[\phi]_{3,1}\nonumber\\
	&\geq \f34 \big\{\|y^{\g} u_s q_{yyyy} \chi(\f{y}{\xi})\|^2 + 10\|y^{\g} u_{sy}q_{yyy}\chi(\f{y}{\xi})\|^2\big\} - L^{\f12}[\phi]_{3,1},
\end{align}
where we have used the fact
\begin{align*}
 8\big\la y^{2\g} u_su_{sy} q_{yyyy},\, q_{yyy} \chi^2(\f{y}{\xi})\big\ra 
	&=-4\big\la \big(y^{2\g}u_s u_{sy} \chi^2(\f{y}{\xi})\big)_y,\, |q_{yyy}|^2 \big\ra - \cancel{4\big\la y^{2\g} u_s u_{sy},\, |q_{yyy}|^2 \big\ra_{y=0}}\nonumber\\
	&= -4 [1\pm O(1)\g]\|y^{\g}u_{sy}q_{yyy}\chi(\f{y}{\xi})\|^2  + L^{\f1{12}-} [\phi]_{3,1}.
\end{align*}
Thus we have from \eqref{D3.27}-\eqref{D3.28} that 
\begin{align}\label{AD.4}
	\|q_{yy}w\|^2_{x=L}
	&\lesssim L^{\f14 \a-} [\phi]_{3,\hat{w}} +  L^{\f14 \a-} [\phi]_{4,1}.
\end{align}

\medskip

Similar calculations as in \eqref{D3.27}-\eqref{AD.4}, we can also derive  
\begin{align}\label{AD.7}
	\|\sqrt{\v}q_{xy}\|^2_{x=0}&\lesssim L^{\f14 \alpha-} [\phi]_{3,\hat{w}} + L^{\f14 \alpha-} [\phi]_{4,1}.
\end{align}
The details are omitted for simplicity of presentation.

Hence we conclude \eqref{AD.1} from \eqref{AD.4} and \eqref{AD.7}.  Therefore the proof of Lemma \ref{lemAD.1} is completed. $\hfill\Box$

\subsection{Some estimates on pressure and pseudo entropy} 
For later use, we denote
\begin{align}\label{3.63-1}
	\begin{split}
		\mathbf{d}_i:&=a \big\{ \la y\ra^{-\f12\mathfrak{l}_0} + \sqrt{\v}^{i-1} \la \sqrt{\v}y\ra^{-\f12\mathfrak{l}_0}\big\},\quad i=1,2,3,
	\end{split}
\end{align}
Then, from \eqref{1.7-3}, we have 
\begin{align}\label{3.64-1}
\begin{split}
|(\rho_{sx}, \rho_{sy},  \rho_{sxx},T_{sx}, T_{sy}, T_{sxx})|
&\cong \mathbf{d}_2, \\
|(\rho_{syy}, \rho_{sxy}, T_{syy},T_{sxy})|
&\cong \mathbf{d}_3,
\end{split}
\end{align}
which, together with \eqref{A.4}-\eqref{7.33-1}, yields that 
\begin{align}\label{3.65-1}
\begin{split}
&|(\tilde{d}_{13}, \hat{d}_{13})| + |\pa_x(\tilde{d}_{13}, \hat{d}_{13})|\lesssim  \mathbf{d}_2  \quad\&\quad  |\pa_y(\tilde{d}_{13}, \hat{d}_{13})|\lesssim \mathbf{d}_3,\\
&|(\tilde{d}_{23}, \hat{d}_{23})| + |\pa_x(\tilde{d}_{23}, \hat{d}_{23})|\lesssim  \mathbf{d}_1  \quad\&\quad  |\pa_y(\tilde{d}_{23}, \hat{d}_{23})|\lesssim \mathbf{d}_2.
\end{split}
\end{align}

\begin{lemma}\label{lemA.5-1}
	Assume \eqref{3.61-1} and \eqref{3.61-11} hold. For $w=1$ or $w_0$, we have 
	\begin{align}\label{3.57-1}
		\begin{split}
			\|(\sqrt{\v}p_x,p_y)w\|^2_{x=0,L}&  \lesssim \f{1}{\sqrt{\v}} [[p]]_{2,w},
		\end{split}
	\end{align}
	\begin{align}\label{3.58-1}
		\begin{split}
			\|\mathbf{d}_2p\|_{x=0}^2 + \v \|\mathbf{d}_2pw\|^2_{x=0}   &\lesssim  \f{a^2}{\sqrt{\v}}[[p]]_{2,1},\\
			\|\mathbf{d}_2 p\|^2_{x=L} + \v \|\mathbf{d}_2 pw\|^2_{x=L}  & \lesssim  a^2 L \frac{1}{\v} [[p]]_{2,1},\\
			\|\mathbf{d}_2 S\|_{x=0,L}^2  + \v \|\mathbf{d}_2 S w\|_{x=0,L}^2  &\lesssim  a^2 L  [[[S]]]_{2,1},
		\end{split}
	\end{align}
	and 
	\begin{align}\label{6.174}
		\begin{split}
			\|\mathbf{d}_2p\|^2 + \v\|\mathbf{d}_2pw\|^2 &\lesssim  a^2 L^2\frac{1}{\v}[[p]]_{2,1},\\ 
			\|\mathbf{d}_2S\|^2 + \v\|\mathbf{d}_2S w\|^2 &\lesssim a^2 L^2 [[[S]]]_{2,1},\\
			\|p \mathbf{1}_{\{y\in[0,2]\}}\|^2&\lesssim L^2 \mathscr{K}(p) + L^2 \f{1}{\v}\mathfrak{K}_{2,1}(p),\\
			\|S \mathbf{1}_{\{y\in[0,2]\}}\|^2&\lesssim L^2 [[[S]]]_{2,1}.
		\end{split}
	\end{align}
\end{lemma}

\noindent{\bf Proof.} From Sobolev inequality, it is clear
\begin{align}\label{5.40-2}
	\begin{split}
		\|p_x\|_{y=0}^2&\lesssim \|p_{x}\|^2 + \|p_{xy}\|^2 \lesssim \frac{1}{\v} [[p]]_{2,1},\\
		\|p_y\|_{y=0}^2&\lesssim \|p_{y}\|^2 + \|p_{yy}\|^2 \lesssim  [[p]]_{2,1},\\
		\|S_x\|_{y=0}^2&\lesssim \|S_{x}\|^2 + \|S_{xy}\|^2 \lesssim [[[S]]]_{2,1},
	\end{split}
\end{align}
where we have used \eqref{B.14} to derive the following estimate
\begin{align}\label{5.40-3}
	\|S_x\|\lesssim a^{-1}\|u_s S_x\| + \|S_{xy}\|\lesssim [[[S]]]_{2,1}.
\end{align}
It is direct to get
\begin{align}\label{6.92}
	\|(\sqrt{\v}p_x,p_y)w\|^2_{x=0,L} \lesssim \f1L \|(\sqrt{\v}p_x,p_y)w\|^2 + \|(\sqrt{\v}p_{xx},p_{xy})w\|\cdot\|(\sqrt{\v}p_{x},p_{y})w\|  \lesssim \f1{\sqrt{\v}} [[p]]_{2,w}.
\end{align}

Noting \eqref{3.61-1} and \eqref{3.61-11}, we have that 
\begin{align}\label{6.59}
	\|\mathbf{d}_2 p\|_{x=0}^2 + \v \|\mathbf{d}_2pw\|^2_{x=0} &\lesssim  \|\mathbf{d}_2p\|^2_{x=0}   
	\lesssim  a^2 \|p_y\|_{x=0}^2\int_0^\infty y\big\{\la y\ra^{-\f12\mathfrak{l}_0} + \v \la \sqrt{\v}y\ra^{-\f12\mathfrak{l}_0}\big\} dy \nonumber\\ 
	&\lesssim a^2 \|p_y\|^2_{x=0} \lesssim  \f{a^2}{\sqrt{\v}}[[p]]_{2,1},\\
	\|\mathbf{d}_2 S\|_{x=0}^2 +  \v \|\mathbf{d}_2 S w\|_{x=0}^2 &\lesssim \|\mathbf{d}_2 S\|_{x=0}^2 \lesssim  a^2 \|S_y\|^2_{x=0} \lesssim a^2 L [[[S]]]_{2,1},
\end{align}
where we have used  \eqref{6.92} and \eqref{3.0-0S}.

Also it clear to have that 
\begin{align}\label{6.93}
	&\|\mathbf{d}_2 p\|^2_{x=L} + \v \|\mathbf{d}_2 pw\|^2_{x=L} \lesssim 	\|\mathbf{d}_2 p\|^2_{x=L} \nonumber\\
	&\lesssim a^2 \int_0^\infty \big\{\la y\ra^{-\f12\mathfrak{l}_0} + \v \la \sqrt{\v}y\ra^{-\f12\mathfrak{l}_0}\big\} \Big(|p(L,0)|^2 + |p(L,y)-p(L,0)|^2\Big) dy\nonumber\\
	&\lesssim La^2 \|p_{x}\|^2_{y=0} + a^2 \|p_y\|^2_{x=L} \lesssim a^2 L \frac{1}{\v} [[p]]_{2,1},
\end{align}
and
\begin{align}\label{6.101}
	\begin{split}
		&\|\mathbf{d}_2 S\|_{x=L}^2  + \v \|\mathbf{d}_2 S w\|_{x=L}^2  \lesssim  \|\mathbf{d}_2 S\|_{x=L}^2  \\
		&\lesssim a^2 \int_0^\infty \big\{\la y\ra^{-\f12\mathfrak{l}_0} + \v \la \sqrt{\v}y\ra^{-\f12\mathfrak{l}_0}\big\} \Big(|S(L,y)-S(L,0)|^2 + |S(L,0)|^2\Big)  dy \\
		&\lesssim a^2 \|S_y\|^2_{x=L} + La^2 \|S_x\|^2_{y=0}  \lesssim a^2L [[[S]]]_{2,1},
	\end{split}
\end{align}
where we have used \eqref{5.40-2} and \eqref{6.92}.

It follows from \eqref{5.40-2} that 
\begin{align}\label{3.67-1}
	& \|\mathbf{d}_2p\|^2 + \v\|\mathbf{d}_2pw\|^2\lesssim \|\mathbf{d}_2p\|^2 \nonumber\\
	&\lesssim a^2  \int_0^L\int_0^\infty \big\{\la y\ra^{-\f12\mathfrak{l}_0} + \v \la \sqrt{\v}y\ra^{-\f12\mathfrak{l}_0}\big\} \big([p(x,y)-p(x,0)]^2 + |p(x,0)|^2\big)  dydx  \nonumber\\
	&\lesssim a^2 \|p_y\|^2 + a^2 L^{2} \|p_x\|^2_{y=0}
	\lesssim  L^2 \f{1}{\v}[[p]]_{2,1},
\end{align}
and
\begin{align}\label{3.119-03}
	&\|\mathbf{d}_2 S\|^2 +  \v \|\mathbf{d}_2 S w\|^2 \lesssim \|\mathbf{d}_2 S\|^2  \nonumber\\
	&\lesssim a^2 \int_0^L\int_0^\infty \big\{\la y\ra^{-\f12\mathfrak{l}_0} + \v \la \sqrt{\v}y\ra^{-\f12\mathfrak{l}_0}\big\} \big([S(x,y)-S(x,0)]^2 + |S(x,0)|^2\big)  dydx  \nonumber\\
	&
	\lesssim a^2 \|S_y\|^2 + a^2 L^2 \|S_x\|^2_{y=0} \lesssim a^2 L^2   [[[S]]]_{2,1}.
\end{align}

Finally, we have 
\begin{align}
	&\|p \mathbf{1}_{\{y\in[0,2]\}}\|^2\lesssim \int_0^L\int_0^2  \big([p(x,y)-p(x,0)]^2 + |p(x,0)|^2\big)  dydx  \nonumber\\
	&\lesssim \|p_y\mathbf{1}_{\{y\in[0,2]\}}\|^2 + L^2 \|p_x\|^2_{y=0}  \lesssim \|p_y\mathbf{1}_{\{y\in[0,2]\}}\|^2 + L^2 \|(p_x,p_{xy})\mathbf{1}_{\{y\in[0,2]\}}\|^2\nonumber\\
	&\lesssim L^2 \mathscr{K}(p) + L^2 \f{1}{\v}\mathfrak{K}_{2,1}(p),
\end{align}
and
\begin{align}
	\|S \mathbf{1}_{\{y\in[0,2]\}}\|^2 &\lesssim \int_0^L\int_0^2  \big([S(x,y)-S(x,0)]^2 + |S(x,0)|^2\big)  dydx  \nonumber\\
	&\lesssim \|S_y\mathbf{1}_{\{y\in[0,2]\}}\|^2 + L^2 \|S_x\|^2_{y=0}  
	\lesssim L^2[[[S]]]_{2,1}.
\end{align}
Therefore the proof of Lemma \ref{lemA.5-1} is completed. $\hfill\Box$

\begin{lemma}\label{lemA.6-0}
	It holds that 
	\begin{align}\label{2TC.3}
		\begin{split}
			\|S_x w\|^2&\lesssim [[[S]]]_{2,\hat{w}},\\
			\|\sqrt{\v}S_xw\|^2_{x=0,L}&\lesssim \sqrt{\v} [[[S]]]_{2,\hat{w}}.
		\end{split}
	\end{align}
	and
	\begin{align}\label{D3.30-0}
		\v \|\sqrt{u_s} S_{xy}w\|^2_{x=0, L}
		&\lesssim \v^{1-}\|S_{xy}w\|^2 + \|\v u_s   S_{xxy}w\|\cdot \|S_{xy}w\|.
	\end{align}
\end{lemma}

\noindent{\bf Proof.} Noting \eqref{B.14}, we have 
\begin{align}\label{2TC.3-1}
	\|S_x w\|^2\lesssim  a^{-2}\|u_s S_x w\|^2 + \|S_x\|^2 \lesssim  a^{-2}\|u_s S_x \hat{w}\|^2 + \|S_{xy}\|^2,
\end{align}
which yields immediately $\eqref{2TC.3}_1$.

It follows from the Sobolev inequality and \eqref{2TC.3-1} that
\begin{align*}
	\|\sqrt{\v}S_xw\|^2_{x=0,L}&\lesssim \|\sqrt{\v}S_{xx}w\|\cdot \|\sqrt{\v}S_xw\| + L^{-1}\|\sqrt{\v}S_xw\|^2  \lesssim  \sqrt{\v} [[[S]]]_{2,\hat{w}}.
\end{align*}
It is clear to have 
\begin{align*}
	\v \|\sqrt{u_s} S_{xy}w\|^2_{x=0, L}&\lesssim \v^{1-}\|S_{xy}w\|^2 + \|\v u_s S_{xxy}w\|\cdot \|S_{xy}w\|.
\end{align*}
Therefore the proof of Lemma \ref{lemA.6-0} is completed. $\hfill\Box$


\begin{lemma}\label{lemA.6}
	Recall $T= S + \frac{1}{2\rho_s}p + \f{1}{p_s}T_{sy}q$. It holds that 
	\begin{align}\label{6.175}
		\begin{split}
			\|T_x w\|^2 &\lesssim L^{\f12\alpha-} [\phi]_{3,\hat{w}} + [[[S]]]_{2,w} + \f{1}{\v} [[p]]_{2,w}+ L\f{1}{\v^2} [[p]]_{2,1}  \mathbf{1}_{\{w=w_0\}},\\
			\|T_y w\|^2 &\lesssim L  [\phi]_{3,\hat{w}} + L[[[S]]]_{2,w} + L\f{1}{\v} [[p]]_{2,w}+ L\f{1}{\v^2} [[p]]_{2,1}  \mathbf{1}_{\{w=w_0\}},
		\end{split}
	\end{align}
	and
	\begin{align}\label{6.176}
		\begin{split}
			\|(\sqrt{\v}T_{xx},T_{xy})w\|^2&  \lesssim   [[[S]]]_{2,w} + \frac{1}{\v} [[p]]_{2,\hat{w}}   +  L^{\f12\alpha-} [\phi]_{3,1},\\
			\|(T_{yy},\Delta_{\v}T)w\|^2 &  \lesssim L^{\f12}[[[S]]]_{2,w} +   L\f{1}{\v}[[p]]_{2,\hat{w}}  + L^{1-} [\phi]_{3,1}.
		\end{split}
	\end{align}
\end{lemma}

\noindent{\bf Proof.}  
A direct calculation shows that 
\begin{align}\label{5.6-1}
	\begin{split}
		T_{x}&= S_{x} + \frac{1}{2\rho_s} p_x + (\frac{1}{2\rho_s})_x p + \big(\f{1}{p_s}T_{sy}q\big)_x,\\
		T_{y}&=S_{y} + \frac{1}{2\rho_s} p_y + (\frac{1}{2\rho_s})_y p + \big(\f{1}{p_s}T_{sy}q\big)_y,
	\end{split}
\end{align}
and
\begin{align}\label{5.6-2}
	\begin{split}
		T_{xx} &= S_{xx} + \frac{1}{2\rho_s}p_{xx} + (\frac{1}{\rho_s})_x p_x +(\frac{1}{2\rho_s})_{xx}p + \big(\f{1}{p_s}T_{sy}q\big)_{xx} ,\\
		T_{yy} &= S_{yy} + \frac{1}{2\rho_s}p_{yy} + (\frac{1}{\rho_s})_y p_y +(\frac{1}{2\rho_s})_{yy}p + \big(\f{1}{p_s}T_{sy}q\big)_{yy},\\
		T_{xy} &= S_{xy} + \frac{1}{2\rho_s}p_{xy} + (\frac{1}{2\rho_s})_x p_y + (\frac{1}{2\rho_s})_y p_x +(\frac{1}{2\rho_s})_{xy}p + \big(\f{1}{p_s}T_{sy}q\big)_{xy},\\
		\Delta_{\v}T&= \Delta_{\v}S + \frac{1}{2\rho_s}\Delta_{\v}p  + [\Delta_{\v},\frac{1}{2\rho_s}]p + \Delta_{\v}\big(\f{1}{p_s}T_{sy}q\big).
	\end{split}
\end{align}

Noting \eqref{5.6-1}, then applying $\eqref{2TC.3}_1$, \eqref{6.174} and \eqref{A.11}-\eqref{A.12},  we  conclude \eqref{6.175} directly. 

Noting \eqref{5.6-2}, then applying  \eqref{6.174} and   \eqref{A.10}-\eqref{A.12}, one shows that 
\begin{align}
	\|\sqrt{\v}T_{xx}w\|^2&\lesssim \|\sqrt{\v} S_{xx}w\|^2 + \|\sqrt{\v} p_{xx}w\|^2 + \|\sqrt{\v} \rho_{sx} p_xw \|^2 + \|\sqrt{\v} \rho_{sxx} p w\|^2 + \|\sqrt{\v}\big(\f{T_{sy}}{p_s}q\big)_{xx}w\|^2 \nonumber\\
	&\lesssim [[[S]]]_{2,w} +  \frac{1}{\v} [[p]]_{2,w}+ L\frac{1}{\v} [[p]]_{2,1} +  L^{\f12\alpha} [\phi]_{3,1}, \label{6.177}\\
	\|T_{xy}w\|^2 &\lesssim \|S_{xy}w\|^2 + \|p_{xy}w\|^2 + \|\rho_{sx} p_y w\|^2 + \|\rho_{sy} p_x w\|^2 + \|\mathbf{d}_3 p w\|^2 +  \|\big(\f{T_{sy}}{p_s}q\big)_{xy}w\|^2\nonumber\\
	&\lesssim [[[S]]]_{2,w} + \frac{1}{\v} [[p]]_{2,w}+ L\frac{1}{\v} [[p]]_{2,1}  +  L^{\f12\alpha} [\phi]_{3,1},	\label{6.178}\\
	\|T_{yy}w\|^2  &\lesssim \|S_{yy}w\|^2 + \|p_{yy}w\|^2 + \|\rho_{sy} p_yw\|^2 + \|\rho_{syy} pw\|^2 + \|\big(\f{T_{sy}}{p_s}q\big)_{yy}w\|^2 \nonumber\\
	&\lesssim L^{\f12}[[[S]]]_{2,w} + [[p]]_{2,w}  + L\f{1}{\v}[[p]]_{2,1}  + L^{1-} [\phi]_{3,1},\label{6.179}\\
	\|\Delta_{\v}Tw\|^2 &\lesssim L^{\f12}[[[S]]]_{2,w} + \f{1}{\sqrt{\v}}[[p]]_{2,w}  + L\f{1}{\v}[[p]]_{2,1}  + L^{1-} [\phi]_{3,1},\label{6.179-1}
\end{align}
which conclude \eqref{6.176}. Therefore the proof of Lemma \ref{lemA.6} is completed. $\hfill\Box$

\subsection{Some useful estimates}\label{Sec2.4}
Under the assumption $\|(\phi, p, S)\hat{w}^{\v}_0\|^2_{\mathbf{X}^{\v}}\lesssim 1$ and $p(0,0)=S(0,0)=0$, we shall derive some $L^\infty$-estimates for $(p, T, \rho)$ and $(u, v)$ for later use.

1. We have  
\begin{align}
	|p(x,y)|^2&\lesssim |p(x,y) - p(0,y)|^2 + |p(0,y)|^2 
	\lesssim L \|p_x(\cdot,y)\|^2_{L^2_x} + \|\la y\ra p_y\|^2_{x=0} \nonumber\\
	&\lesssim L \f{1}{\v} [[p]]_{2,1} + \f{1}{\sqrt{\v}} [[p]]_{2,w_0} \lesssim L \f{1}{\v} [[p]]_{2,\hat{w}_0},\label{3.87}\\
	|S(x,y)|^2
	&\lesssim |S(x,y) - S(0,y)|^2 + |S(0,y)|^2 
	\lesssim L \|S_x(\cdot,y)\|^2_{L^2_x} + \|\la y\ra S_y\|^2_{x=0} \nonumber\\
	&\lesssim L [[[S]]]_{2,1} +   L [[[S]]]_{2,w_0}\lesssim L [[[S]]]_{2,\hat{w}_0}.\label{3.88}
\end{align}
Noting  \eqref{s}, using \eqref{3.87}-\eqref{3.88} and  , one obtains that
\begin{align}\label{3.90}
	|T(x,y)|^2 
	&\lesssim L \f{1}{\v} [[p]]_{2,\hat{w}_0} + L [[[S]]]_{2,\hat{w}_0} + L \f{1}{\sqrt{\v}} [\phi]_{3,1}.
\end{align}
Recall $\bar{\rho}, N_{\rho}$ in \eqref{1.19-1}, we have from \eqref{3.87} and \eqref{3.90} that 
\begin{align}\label{3.91}
	\begin{split}
		&\v^{N_0} |\bar{\rho}| \lesssim \v^{N_0} \Big(L \f{1}{\v} [[p]]_{2,\hat{w}_0} + L [[[S]]]_{2,\hat{w}_0} + L \f{1}{\sqrt{\v}} [\phi]_{3,1}\Big)^{\f12} \lesssim \v^{N_0-\f14}\ll 1,\\
		& \v^{N_0}N_{\rho} 
		\lesssim \v^{2N_0} \Big(L \f{1}{\v} [[p]]_{2,\hat{w}_0} + L [[[S]]]_{2,\hat{w}_0} + L \f{1}{\sqrt{\v}} [\phi]_{3,1}\Big)  \lesssim \v^{2N_0-1}\ll 1,
	\end{split}
\end{align}
Then we have from \eqref{3.91} that 
\begin{align}\label{3.91-0}
	|\v^{N_0}\rho|=|\v^{N_0}\bar{\rho} + \v^{N_0} N_{\rho}|\lesssim \v^{N_0-\f14}.
\end{align}
Hence it follows from \eqref{3.90}-\eqref{3.91-0} that
\begin{align}\label{3.92}
	\begin{split}
		\rho^{\v} = \rho_s + O(1) \v^{N_0-\f14} \cong \rho_s \quad \mbox{and}\quad 
		T^{\v} = T_s + O(1) \v^{N_0-\f14}\cong T_s.
	\end{split}
\end{align}

\smallskip

2. For $u^{\v}, v^{\v}$, it follows from \eqref{7.1-1}, \eqref{3.91-0}-\eqref{3.92} and \eqref{2.122-0} that 
\begin{align}\label{3.93}
	\begin{split}
		u^{\v}:&=u_s + O(1)[\v^{N_0}|\phi_y|+u_s |\v^{N_0}\rho| + A_{s1}] = u_s + O(1)\v^{N_0-1} u_s,\\
		v^{\v}:&=v_s + O(1)[\v^{N_0}|\phi_x|+v_s |\v^{N_0}\rho| + A_{s2}]= v_s + O(1) \v^{N_0-\f32} u_s^{2-\f1k},
	\end{split}
\end{align}
for any $k\geq 2$, where we have used 
\begin{align*}
	|\phi_x(x,y)|&=\Big|\int_0^y dz \int_0^z \phi_{xyy}(x,s) ds \Big|\leq y^{2-\f1k} \|\phi_{xyy}(x,\cdot)\|_{L^k_y} \nonumber\\
	&\lesssim C(\f1L,k) y^{2-\f1k} \|\phi_{xyy}\|_{H^1} \lesssim \v^{-\f54} y^{2-\f1k}.
\end{align*}
Similarly, for any $k \geq 2$, one has
\begin{align}\label{3.94}
	\begin{split}
		\v^{N_0}|u_s(\nabla\rho, \nabla T)| &\lesssim \v^{N_0-1} u_s  + \v^{N_0} |u_s (\nabla p, \nabla T)|\lesssim \v^{N_0-2}u_s^{1-\f1k},\\
		\v^{N_0}|(u_x, \nabla v)|&\lesssim \v^{N_0-1}u_s  + \v^{N_0} |(u_s \nabla p, u_s \nabla T, \nabla\phi_x)|\lesssim \v^{N_0-2} u_s^{1-\f1k},\\
		\v^{N_0}|u_y|&\lesssim \v^{N_0-1}  + \v^{N_0} |(u_s \nabla p, u_s \nabla T,  \phi_{yy})|\lesssim \v^{N_0-2}.
	\end{split}
\end{align}

\medskip

\section{Approximation Systems and Iteration Schemes}\label{sec2}
Given $(\phi, p, S)$,  then we  regard
\begin{align}\label{2.122-1}
	g_i=g_i(\phi,p,T(S,p,q)), \quad \mathfrak{g}_i=\mathfrak{g}_i(q),\quad \mathcal{N}_i= \mathcal{N}_{i}(\phi,p,T(S,p,q)), \quad i=1,2,
\end{align}
as known functions, where  $T(S,p,q)=S + \f{1}{2\rho_s}p + \f{1}{p_s}T_{sy} q$. We  also regard 
\begin{align}
	(\rho^{\v}, T^{\v}, u^{\v}, v^{\v})=(\rho_s+\v^{N}\rho, \,  T_s + \v^{N_0}T,\,  u_s + \v^{N_0} u,\,  v_s + \v^{N_0} v),
\end{align}
as given functions with 
\begin{align}
	\begin{split}
		u:=\f{1}{\rho^{\v}} \big[\phi_y - u_s\, \rho - \v^{-N_0}A_{s1} \big]\,\,\,\mbox{and}\,\,\,
		v:=\f{1}{\rho^{\v}} \big[-\phi_x - v_s\, \rho - \v^{-N_0}A_{s2} \big].
	\end{split}
\end{align}

\smallskip

For above given $(\phi,p, S)$, we assume 
\begin{align}\label{2.122-0}
	\begin{split}
		\|(\phi, p, S)\hat{w}^{\v}_0\|^2_{\mathbf{X}^{\v}}\leq \v,
	\end{split}
\end{align}
and the compatibility conditions
\begin{align}\label{2.123}
\begin{cases}
\dis 	\phi|_{y=0}=\phi_y|_{y=0}=0,\quad \phi|_{x=0}=\phi_{xx} |_{x=0}=\phi_{x}|_{x=L}=\phi_{xxx} |_{x=L}=0,\\[1.5mm]
\dis 	S_y \big|_{y=0} 
= - \big(\frac{1}{2\rho_s}p\big)_y ,\quad  S \big|_{x=0}= -  \frac{\chi}{2\rho_s}p \quad \mbox{and}\quad  S_{x}\big|_{x=L}=0 \quad \mbox{for NT}, \\[1.5mm]
 S \big|_{y=0}= -\f{1}{2\rho_s}p \quad \mbox{and}\quad  S \big|_{x=0}=0\,\,\, 
\,\, \mbox{for DT}. 
\end{cases}
\end{align}

It is  hard to construct solution for system  \eqref{T.1} and \eqref{2.92} with \eqref{pair1},  \eqref{pressure}, \eqref{2.17-0}, \eqref{T.2}-\eqref{3.61-11},  \eqref{A1.56} and \eqref{3.61-1}. To archive our target, we need several approximations. In fact, for the unknown functions $(\Phi, P, \mathbf{S})$,  we consider the {\it modified} momentum system 
\begin{align}\label{2.109}
	\begin{cases}
		\dis \mathbf{A}_1:= m u_s^2 \mathbf{q}_{xy}  + \mathbf{d}_{11}\, P_x -\frac{\mu}{\rho^{\v}} \Delta_{\v}\Phi_y + \frac{\mu+\lambda}{p^{\v}} u^{\v} \Delta_{\v}P  - \frac{\lambda}{p^{\v}}\big[1- \f{\rho^{\v}}{2\rho_s}\big] u^{\v} P_{yy}\\
		\dis \qquad\,\,  + \frac{\lambda\v}{p^{\v}} v^{\v} P_{xy}  + \frac{\lambda}{T^{\v}} u^{\v} \mathbf{S}_{yy} + \hat{d}_{11} \zeta_x  + g_1(p,\phi, T(S,p,q))  + \mathfrak{g}_1(q)  \\
		\dis \qquad\,\, - G_{11}(\tilde{T}(\zeta,p,q)) - G_{12}(p,q) - \mathcal{N}_1=0,\\
		\dis\mathbf{A}_2:=  -\v m u_s^2 \mathbf{q}_{xx} +   \tilde{d}_{22}\, P_y +\frac{\mu\v }{\rho^{\v}}   \Delta_{\v}\Phi_x + \frac{\mu\v}{p^{\v}} v^{\v} \Delta_{\v}P  + \frac{\lambda\v}{p^{\v}} \big[1- \f{\rho^{\v}}{2\rho_s}\big] u^{\v} P_{xy}\\
		\dis\qquad\,\,\, + \frac{\lambda\v}{p^{\v}} v^{\v} P_{yy} - \frac{\lambda\v}{T^{\v}} u^{\v} \mathbf{S}_{xy}  +  g_2(p,\phi, T(S,p,q)) + \v \mathfrak{g}_2(q) \\
		\dis\qquad\,\,\,  -  \v G_{21}(\tilde{T}(\zeta,p,q))  - \v G_{22}(p,q) - \v \mathcal{N}_2=0,
	\end{cases}
\end{align}
where  $\mathbf{q}:= \f{\Phi}{\bar{u}_s}$ and  $\tilde{T}(\zeta, p,q)=\zeta + \f{1}{2\rho_s}p + \f{1}{p_s}T_{sy} q$. We remark that all the notations in \eqref{2.109} can be found in Appendix \ref{subsec10.1}. For  $\Phi$, we impose the following boundary conditions 
\begin{align}\label{2.113} 
\Phi|_{y=0}&=\Phi_y|_{y=0}=0,\quad  \Phi|_{x=0}=\Phi_{xx}|_{x=0}=\Phi_{x}|_{x=L}=\Phi_{xxx}|_{x=L}=0.
\end{align}
For the pseudo entropy  $\mathbf{S}$ in the case of NT, we consider the following {\it modified} boundary value problem (BVP)
\begin{align} \label{2.111}
\begin{cases}
\dis 2 \rho_s u_s \mathbf{S}_x + 2\rho_s v_s \mathbf{S}_y - \kappa \Delta_{\v}\mathbf{S} \\
\dis \qquad \qquad\qquad  = \kappa \Delta_{\v}\Big(\frac{1}{2\rho_s}[\chi P + \bar{\chi}p]\Big) + \kappa \Delta_{\v}\Big(\f{1}{p_s}T_{sy}q\Big)  + \mathfrak{J}\big(\phi,p,T(S,p,q)\big),\\
\dis \mathbf{S}_y \big|_{y=0} = - \big(\frac{1}{2\rho_s}P\big)_y, \quad \mathbf{S} \big|_{x=0}= -  \frac{\chi}{2\rho_s}P, \,\, \, \mbox{and} \,\, \, \mathbf{S}_{x}\big|_{x=L}=0.
\end{cases} 
\end{align}
While, for the case of DT, the {\it modified} BVP of pseudo entropy $\mathbf{S}$ is
\begin{align}\label{2.111-6D}
\begin{cases}
\dis	2 \rho_s u_s \mathbf{S}_x + 2\rho_s v_s \mathbf{S}_y - \kappa \Delta_{\v}\mathbf{S}  = \kappa \Delta_{\v}\big(\frac{1}{2\rho_s} p \big) + \kappa \Delta_{\v}\big(\f{1}{p_s}T_{sy}q\big)  + \mathfrak{J}(\phi,p,T(p,S,q)),\\
\dis \mathbf{S}|_{y=0}=-\f1{2\rho_s}\mathfrak{P}(x),\quad \mathbf{S}|_{x=0}=0 \,\,\,\, \mbox{and}\,\,\, \big\{2\rho_s u_s \mathbf{S}_x-\kappa \mathbf{S}_{yy}\big\}|_{x=L}= \mathfrak{G}(L,y),
\end{cases}
\end{align}
where the functions $\mathfrak{P}(x)$ and $\mathfrak{G}(x,y)$ are given by
\begin{align}\label{D3.1-1}
	\begin{split}
	\dis \mathfrak{P}_x:&=\big(\mu\rho_s^{-1} \phi_{yyy} - g_1 + \mathfrak{h}_1 + \mathfrak{N}_{11} - \mu \v^{N_0} \rho (\rho^{\v}\rho_s)^{-1}  \phi_{yyy}\big)(x,0)\,\,\,\mbox{with}\,\,\,
	\mathfrak{P}(0)=0,\\
	\dis \mathfrak{G}:&= \kappa p_s^{-1}T_{sy}q_{yy} + 2\kappa \big(p_s^{-1}T_{sy}\big)_{y}q_y + \kappa\big((2\rho_s)^{-1}\big)_{yy}p + \mathfrak{M}_0.
	\end{split}
\end{align}
All the notations in \eqref{2.111}-\eqref{D3.1-1} can be found in Appendix \ref{subsec10.1}.

We point out that  $\eqref{D3.1-1}_1$ arises from the necessary boundary condition $A_1|_{y=0}=0$. It is crucial to use $\mathfrak{P}$ rather than $p(x,0)$, since better trace estimate can be derived by using $\mathfrak{P}$.
Furthermore, we can prove that $\mathfrak{P}(x)=p(x,0)$ when a fixed point is established for the iteration $(\phi, p, S) \to (\Phi, P, \mathbf{S})$.
We also remark that the equations of pseudo entropy $\mathbf{S}$ differ slightly between the Neumann and Dirichlet boundary value problems. The main reason is that the boundary condition $\mathbf{S} \big|_{x=0}= -  \frac{\chi}{2\rho_s}P$ in   \eqref{2.111}  plays an important role when we derive boundary condition of pressure $P$ from the first momentum equation $\eqref{2.109}_1$, see Lemma \ref{lemAF-1} for details.

For $\zeta$ in \eqref{2.109}, it satisfies the following BVP:
\begin{align} \label{2.111-1}
	\begin{cases}
		\dis 2 \rho_s u_s \zeta_x + 2\rho_s v_s \zeta_y - \kappa \Delta_{\v}\zeta  
		= \kappa \Delta_{\v}\Big(\frac{1}{2\rho_s}p\Big) + \kappa \Delta_{\v}\Big(\f{1}{p_s}T_{sy}q\Big) + \mathfrak{J}\big(\phi,p,T(S,p,q)\big),\\[3mm]
		\dis \zeta_y \big|_{y=0} = - \big(\frac{1}{2\rho_s}p\big)_y , \quad \zeta \big|_{x=0}= -  \frac{\chi}{2\rho_s}p,  \,\, \,\mbox{and}\,\,\, \zeta_{x}\big|_{x=L}=0,
		\quad \mbox{for NT},\\[3mm]
		\dis \zeta|_{y=0}=-\f1{2\rho_s}\mathfrak{P}(x),\quad \zeta|_{x=0}=0 \,\,\,\, \mbox{and}\,\,\, \big\{2\rho_s u_s \zeta_x-\kappa \zeta_{yy}\big\}|_{x=L}= \mathfrak{G}  \quad \mbox{for DT}.
	\end{cases}
\end{align}
The introduction of auxiliary function $\zeta$ in \eqref{2.109}  is necessary exclusively for the pseudo entropy of NT case, it  plays an important role in the multiplier $\mathbf{q}_{xx}$, see \eqref{3.50-1}-\eqref{3.64-2} for details.
For the pseudo entropy of DT case, we point out that \eqref{2.111-6D} is structurally identical to \eqref{2.111-1}, that means we always have $\mathbf{S}\equiv \zeta$. 
To maintain notational consistency in constructing solutions for both NT and DT, we retain the notation $\zeta$ throughout our treatment of DT case. 
Also, since \eqref{2.111-6D} and \eqref{2.111-1} are independent of the variables $\Phi$ and $P$ defined in \eqref{2.109},  then we  can solve   \eqref{2.111-6D} and \eqref{2.111-1} independently as the first step.


\subsection{Deriving the curl and divergence equations from \eqref{2.109}}
We consider the curl of momentum equations \eqref{2.109}, i.e.,
\begin{align*}
	\pa_y\big(b_1\mathbf{A}_1\big) - \pa_{x} \big(b_2 \mathbf{A}_2\big) =0,
\end{align*}
which can be represented as  
\begin{align}\label{2.112}
	&\big(b_1m u_s^2 \mathbf{q}_{xy}\big)_y  + \v \big(b_2mu_s^2 \mathbf{q}_{xx}\big)_x   - \mu\Big(\frac{b_1}{\rho^{\v}} \Delta_{\v}\Phi_y\Big)_y - \mu\v \Big(\frac{b_2}{\rho^{\v}}   \Delta_{\v}\Phi_x\Big)_x\nonumber\\
	&= \mathcal{F}_{R}(P,\mathbf{S}, \zeta, \phi, p, S, ) =: \sum_{i=1}^8 \mathcal{F}_{Ri},
\end{align}
where the terms $\mathcal{F}_{Ri}, i=1,\cdots, 8$ are defined in Appendix \ref{appA.3}.  
Here the term  $\big(b_1m u_s^2 {\bf q}_{xy}\big)_y  + \v \big(b_2mu_s^2 {\bf q}_{xx}\big)_x$ is regarded as a generalized Rayleigh operator.

\smallskip

We consider the divergence of momentum equations \eqref{2.109}, i.e.,
\begin{align*}
	\v \pa_x \mathbf{A}_1 + \pa_y \mathbf{A}_2 =0,
\end{align*}
which is rewritten as
\begin{align}\label{2.115}
	\Delta_{\v}P + \v (\mu+\lambda) \frac{1}{p^{\v}} \big(U^{\v}\cdot \nabla\big) \Delta_{\v}P 
	&=\mathcal{G}(P,\mathbf{S}, \Phi, p, S, \zeta, \phi)=:\sum_{i=1}^8 \mathcal{G}_i,
\end{align}
where  the terms $\mathcal{G}_{i}, i=1,\cdots, 8$ are also given in Appendix \ref{appA.3}.

To solve \eqref{2.115}, we supplement it with the boundary condition
\begin{align}\label{17.15}
	(\mu+\lambda) u_s \Delta_{\v}P\big|_{x=0}= -\sigma p_s \mathbf{d}_{11} B,
\end{align}
with $B$ defined as
\begin{align}\label{17.14}
	B:&= P_x + \f{1}{\mathbf{d}_{11}}\big\{g_1(p,\phi, T(S,p,q))  - \mathfrak{h}_1 - \mathfrak{N}_{11}\big\} 
	\equiv P_x + \bar{g}_1(p,\phi, T(S,p,q)),
\end{align}
where we have used the notation
\begin{align}\label{17.14-1}
	\mathfrak{N}_{11}:= - \v^{N_0} \f{2\mu}{\rho^{\v}}\rho_y\cdot u_y ,
\end{align}
which is indeed a nonlinear term in \eqref{9.26}, see also  $\mathfrak{h}_1$ in \eqref{7.33-11}.

To recover the momentum equations \eqref{2.109} from \eqref{2.112} and \eqref{2.115}, we also need  
\begin{align}\label{2.104}
	\mathbf{A}_1\big|_{y=0}=\mathbf{A}_1\big|_{x=0}=\mathbf{A}_2\big|_{x=L}=0,
\end{align}
when we apply the div-curl arguments, see Lemma \ref{lem9.1}.

\subsection{Schemes} 
Even for the div-curl equations \eqref{2.112} and \eqref{2.115}, we are still not clear how to solve them, so we need  further approximations. We explain our iterative schemes.

1) For given $(\phi, p, S)$ with  \eqref{2.122-0}-\eqref{2.123}, we first solve  $\zeta$ of \eqref{2.111-1}, and $\mathbf{S}$ of \eqref{2.111-6D}. The details are presented in Sections \ref{secT.2} \& \ref{SecD7}, respectively.

2)   Given $(\hat{P},\hat{\mathbf{S}})$. For the case of NT, we need $(\hat{P},\hat{\mathbf{S}})$ to satisfy the following compatibility conditions
\begin{align}\label{2.125}
	\hat{\mathbf{S}} \big|_{x=0}= -  \frac{\chi}{2\rho_s}\hat{P},\,\,\, \hat{\mathbf{S}}_{x}\big|_{x=L}=0 \,\, \,\mbox{and}\,\,\,
	\hat{\mathbf{S}}_y \big|_{y=0} = - \big(\frac{1}{2\rho_s}\hat{P}\big)_y.
\end{align}
For the case of DT, we take $\hat{\mathbf{S}}=\mathbf{S}$ since  $\mathbf{S}$ has already been solved in 1). 

With above preparations,  we  consider the approximate problem for the curl equation \eqref{2.112} 
\begin{align}\label{2.121}
	\begin{cases}
		\dis \big(b_1 m u_s^2 \mathbf{q}_{xy}\big)_y  + \v \big(b_2 mu_s^2 \mathbf{q}_{xx}\big)_x   - \mu\Big(\frac{b_1}{\rho^{\v}} \Delta_{\v}\Phi_y\Big)_y-\mu\v \Big(\frac{b_2}{\rho^{\v}}   \Delta_{\v}\Phi_x\Big)_x  =\mathcal{F}_{R}(\hat{P}, \hat{\mathbf{S}},  p, S,  \zeta, \phi),\\
		\dis \Phi|_{x=0}=\Phi_{xx} |_{x=0}=\Phi_{x}|_{x=L}=\Phi_{xxx}|_{x=L} = \Phi|_{y=0}=\Phi_y|_{y=0}=0,
	\end{cases}
\end{align}
see Section \ref{sec4-0}  for details.

\smallskip

3) For the pressure $P$,  it is hard to solve the original BVP problem \eqref{2.115}-\eqref{17.15} and \eqref{2.104} with enough regularity in  $(0,L)\times \mathbb{R}_+$. Instead, we consider the following  approximation problem 
\begin{align}\label{2.122}
	\begin{cases}
		\dis \Delta_{\v}P +  \v (\mu+\lambda) \f{1}{p^{\v}}(U^{\v}\cdot\nabla)\Delta_{\v}P =G(\hat{P}, \hat{\mathbf{S}}, \Phi, p, S, \zeta, \phi) ,\,\, (x,y)\in (0,L)\times \mathbb{R}_+,\\
		\dis (\mu+\lambda) \f{1}{p_s} u_s \mathbf{w} \big|_{x=0}=-\sigma \mathbf{d}_{11} B ,  
	\end{cases}
\end{align}
where $B$ is the one defined in \eqref{17.14}. A modified version of \eqref{2.104} is also required; see \eqref{18.39-2} for details. As mentioned in Section \ref{Sec1.2}, to solve \eqref{2.122}, we shall first establish the existence and uniform estimate of $\mathbf{w}$ for a slightly modified equation of \eqref{w1}, then recover $P$ by solving $B$ for a modified equation of \eqref{B1}. It is not a good idea to solve $P$ directly from $\mathbf{w}$ due to the bad influence of corner points  $(0,0)$ \& $(L,0)$. In fact, by using $B$, we can extend \eqref{B1} in the corner domain $(0,L)\times \mathbb{R}_+$ into $(0,L)\times \mathbb{R}$, then it will help us to derive the uniform  estimate  of higher order derivatives.  We remark that, for strictness,   a further viscous $\delta$-approximation argument for \eqref{w1} is needed, see Section \ref{sec7} for details.

\smallskip

4) Let  $P$ be the one established in 3). For the pseudo entropy BVP \eqref{2.111} of NT, we study
\begin{align} \label{2.111-0}
	\begin{cases}
		\dis 2 \rho_s u_s \mathbf{S}_x + 2\rho_s v_s \mathbf{S}_y - \kappa \Delta_{\v}\mathbf{S} 
		= \kappa \Delta_{\v}\big(\frac{1}{2\rho_s}[\chi P + \bar{\chi} p]\big) + \kappa \Delta_{\v}\big(\f{1}{p_s}T_{sy}q\big)  + \mathfrak{J},\\[2mm]
		\dis \mathbf{S}_y \big|_{y=0} = - \big(\frac{1}{2\rho_s}P\big)_y,\quad \mathbf{S} \big|_{x=0}= -  \frac{\chi}{2\rho_s}P,  \,\, \,\mbox{and}\,\,\, \mathbf{S}_{x}\big|_{x=L}=0,
	\end{cases}
\end{align} 
see Section \ref{secT.1} for details.

\medskip

5) For the mapping $(\hat{P},\hat{\mathbf{S}}) \to (\Phi, P, \mathbf{S})$, we can prove that it is indeed a contraction.  We emphasize that it is hard  to get uniform-in-$\v$ estimate for $\|(P_x, \f{1}{\sqrt{\v}}P_y)\|$. Fortunately, we have some slightly good estimate for $\|\mathbf{d}_2(P_x, P_y)\|$ and $\|(P_x, P_y) \mathbf{1}_{\{y\in[0,2]\}}\|$ which are enough at least in this step.  However, when we consider the mapping 
$$(\phi,p, S) \to (\Phi, P, \mathbf{S}),$$
better uniform-in-$\v$ estimate of $\|(P_x, \f{1}{\sqrt{\v}}P_y)\|$ is needed which can be established by {\it div-curl} argument. Then we can prove that there exists a fixed point in above mapping, i.e., $(\phi,p, S)=(\Phi, P, \mathbf{S})$, and hence $\zeta=\mathbf{S}$. Finally, we point out that it is impossible to apply the {\it div-curl} argument for solutions established in 2) \& 3) because the {\it div-curl} structure is destroyed.

\smallskip

\section{Existence and Uniform Estimate for Stream Function}\label{sec4-0}
In this section, we aim to establish the existence and uniform estimates for the solution $\Phi$ to \eqref{2.121}.
\begin{theorem}\label{thm4.1}
	There exists a unique solution $\Phi$ to \eqref{2.121} satisfying
	\begin{align}
		[\Phi]_{3,\hat{w}_0^{\v}}
		&\lesssim L^{\f14-} [\Phi]_{4,\hat{w}_0^{\v}} + (\f{1}{N} + N a ) \|(\phi, p, S(\zeta))\hat{w}_0^{\v}\|^2_{\mathbf{X}^{\v}}  + Na  \|(0, \hat{P}, \hat{\mathbf{S}})\hat{w}_0^{\v}\|^2_{\mathbf{Y}^{\v}} \nonumber\\
		&\quad +  L^{\f12\alpha-} \|u_s\sqrt{\v}\Delta_{\v}\tilde{T}_{x}\hat{w}_0^{\v}\|^2  + \f{1}{\v} \|(\mathfrak{N}_3, \mathcal{F}_{R8}) \hat{w}_0^{\v}\|^2  +  \v^{N_0 },\label{18.3}\\
		[\Phi]_{4,\hat{w}_0^{\v}}
		&\lesssim a^3 [\Phi]_{3,\hat{w}_0^{\v}} +  [[p]]_{3,\hat{w}_0^{\v}} + [[[\zeta]]] _{3,\hat{w}_0^{\v}}  + \mathcal{K}_{3,\hat{w}_0^{\v}}(\hat{P})  + [[[\hat{\mathbf{S}}]]]_{3,\hat{w}_0^{\v}}  +  a^2  \|(\phi, p, S(\zeta)) \hat{w}_0^{\v}\|^2_{\mathbf{X}^{\v}}  \nonumber\\
		&\quad +  L^{\f12}  \|(0, \hat{P}, \hat{\mathbf{S}}) \hat{w}_0^{\v}\|^2_{\mathbf{Y}^{\v}} + \v^{1-} \|u_s\sqrt{\v} \Delta_{\v}\tilde{T}_x \hat{w}_0^{\v}\|^2 + \|\mathcal{F}_{R8}\hat{w}_0^{\v}\|^2.\label{18.4}
	\end{align}
\end{theorem}

Recall the given functions $(\phi, p, S)$ and $(\hat{P}, \hat{\mathbf{S}})$ in Section \ref{sec2}. Noting the smallness of $u_s$, the existence of a unique solution $\Phi$ to \eqref{2.121} can be established by arguments similar  to those in \cite{Guo-Iyer-2024}; we omit the details for simplicity. Thus, it remains to establish the uniform estimates \eqref{18.3}-\eqref{18.4}.

\subsection{Quotient estimate}\label{sec4}
In this subsection, we concerns on the quotient estimate $\||\mathbf{q}|\|_w$. The coupling with the pseudo entropy makes this problem considerably more challenging than the case for the incompressible Navier-Stokes equations \cite{Guo-Iyer-2024}.
\begin{lemma}\label{lem3.1}
It holds that 
\begin{align}\label{3.1}
&\|\sqrt{u_s} \big(\mathbf{q}_{xyy},\, \sqrt{\v} \mathbf{q}_{xxy},\, \v \mathbf{q}_{xxx}\big)w\|^2 
+ \|u_s \mathbf{q}_{xy} w\|_{x=0}^2 + \|\sqrt{\v}u_s \mathbf{q}_{xx}w\|^2_{x=L}  + \|\sqrt{u_{sy}}\mathbf{q}_{xy}\|_{y=0}^2 \nonumber\\
&\lesssim   L^{\f18} \|\mathbf{q}_{xx}w_y\|^2 + L^{\f12\alpha-} [\Phi]_{3,\hat{w}} + L^{\f12-}[\Phi]_{4,\hat{w}} + \big|\big\la \mathcal{F}_{R}(\hat{P}, \hat{\mathbf{S}}, \zeta, p, \phi, S),\, \mathbf{q}_{xx}w^2\big\ra\big|,
\end{align}
with 
\begin{align}\label{3.1-0}
&\big|\big\la \mathcal{F}_{R}(\hat{P}, \hat{\mathbf{S}}, \zeta, p, \phi, S),\, \mathbf{q}_{xx}w^2\big\ra\big|\nonumber\\
&\lesssim \f1N \big\{\|\mathbf{q}_{xx}w_y\|^2  + [\Phi]_{3,\hat{w}} \big\} + \big(\f{1}{N} + N a\big) \|(\phi, p, S(\zeta))\hat{w}\|^2_{\mathbf{X}^{\v}}   + Na  \|(0, \hat{P}, \hat{\mathbf{S}})\hat{w}\|^2_{\mathbf{Y}^{\v}} \nonumber\\
&\,\,\, + L^{\f12} \f{1}{\v}\|(0, p, S)\|^2_{\mathbf{X}^{\v}} \mathbf{1}_{\{w=w_0\}}+  L^{\f12\alpha-} \|u_s\sqrt{\v}\Delta_{\v}\tilde{T}_{x}w\|^2   + \f{1}{\v} \|(\sqrt{u_s}\mathfrak{N}_3, \mathcal{F}_{R8}) w\|^2 
+  \v^{N_0 }. 
\end{align}
\end{lemma}

\noindent{\bf Proof.} Multiplying \eqref{2.121} by $\mathbf{q}_{xx}w^2$, one has 
\begin{align}\label{3.2}
&\big\la\big(b_1m u_s^2 \mathbf{q}_{xy}\big)_y,\, \mathbf{q}_{xx}w^2\big\ra + \big\la \v \big(b_2mu_s^2 \mathbf{q}_{xx}\big)_x ,\, \mathbf{q}_{xx}w^2\big\ra  \nonumber\\
&-\Big\la\mu\,\Big(\frac{b_1}{\rho^{\v}} \Delta_{\v}\Phi_y\Big)_y,\, \mathbf{q}_{xx}w^2\Big\ra
-\Big\la\mu\,\v \Big(\frac{b_2}{\rho^{\v}}   \Delta_{\v}\Phi_x\Big)_x,\, \mathbf{q}_{xx}w^2\Big\ra \nonumber\\
&= \big\la \mathcal{F}_{R}(\hat{P}, \hat{\mathbf{S}}, \zeta, p, \phi, S),\, \mathbf{q}_{xx}w^2\big\ra.
\end{align}

\noindent{\it Step 1. Estimate on Rayleigh term.} 
 Integrating by parts with respect to $y$, one can obtain 
\begin{align}\label{3.12}
\big\la\big(b_1m u_s^2 \mathbf{q}_{xy}\big)_y,\, \mathbf{q}_{xx}w^2\big\ra 
&=-\Big\la b_1m u_s^2 \mathbf{q}_{xy},\, \mathbf{q}_{xxy}w^2\Big\ra - \big\la b_1m u_s^2 \mathbf{q}_{xy},\, 2\mathbf{q}_{xx}w w_y\big\ra\nonumber\\
&\geq \f12 \|\sqrt{mb_1}u_s\mathbf{q}_{xy}w\|_{x=0}^2 - C L^{\f18} [\Phi]_{3,\hat{w}} - CL^{\f18} \|\mathbf{q}_{xx}w_y\|^2,
\end{align}
where we have used the following two estimates
\begin{align}\label{3.5}
&-\Big\la b_1m u_s^2 \mathbf{q}_{xy},\, \mathbf{q}_{xxy}w^2\Big\ra \nonumber\\
&=\f12 \|\sqrt{mb_1}u_s\mathbf{q}_{xy}w\|_{x=0}^2 -  \f12 \|\sqrt{mb_1}u_s\mathbf{q}_{xy}w\|_{x=L}^2 
+ \f12 \Big\la (m b_1u_s^2)_x, \, |\mathbf{q}_{xy}|^2w^2\Big\ra \nonumber\\
&=\f12 \|\sqrt{mb_1}u_s\mathbf{q}_{xy}w\|_{x=0}^2 + O(1) L^{\f14-} [\Phi]_{3,\hat{w}},
\end{align}
and
\begin{align}\label{3.11}
\big|\big\la b_1m u_s^2 \mathbf{q}_{xy},\, 2\mathbf{q}_{xx}w w_y\big\ra\big| 
&\lesssim a \|\mathbf{q}_{xx}w_y\|\cdot  \|u_s\mathbf{q}_{xy}w\|   \lesssim   L^{\f18} \|\mathbf{q}_{xx}w_y\|^2 + L^{\f18} [\Phi]_{3,w}.
\end{align}

A direct calculation shows that 
\begin{align}\label{3.14}
\big\la \v \big(b_2\,mu_s^2 \mathbf{q}_{xx}\big)_x ,\, \mathbf{q}_{xx}w^2\big\ra 
&=\Big\la \v m b_2 u_s^2 \mathbf{q}_{xxx},\, \mathbf{q}_{xx}w^2\Big\ra + \big\la \v \big(b_2\,mu_s^2\big)_x \mathbf{q}_{xx} ,\, \mathbf{q}_{xx}w^2\big\ra \nonumber\\
&=\f12 \|\sqrt{mb_2}\sqrt{\v}u_s\mathbf{q}_{xx}w\|_{x=L}^2 - \f12 \|\sqrt{mb_2}\sqrt{\v}u_s\mathbf{q}_{xx}w\|_{x=0}^2\nonumber\\
&\quad  + \f12 \big\la \v (m b_2 u_s^2)_x \mathbf{q}_{xx}^2,\, w^2\big\ra 
\nonumber\\
&\geq  \f12 \|\sqrt{b_2}\sqrt{\v}u_s\mathbf{q}_{xx}w\|_{x=L}^2 - C L^{\f14} [\Phi]_{3,\hat{w}}.
\end{align}

\noindent{\it Step 2. Estimate on viscous term.}  Our aim is to establish the following estimate
\begin{align}\label{3.18}
&\Big\la\mu\Big(\frac{b_1}{\rho^{\v}} \Delta_{\v}\Phi_y\Big)_y,\, \mathbf{q}_{xx}w^2\Big\ra + \Big\la\mu\v \Big(\frac{b_2}{\rho^{\v}}   \Delta_{\v}\Phi_x\Big)_x,\, \mathbf{q}_{xx}w^2\Big\ra\nonumber\\
&\lesssim  - \|\sqrt{u_s}\big(\mathbf{q}_{xyy}, \sqrt{\v}\mathbf{q}_{xxy}, \v \mathbf{q}_{xxx}\big)w\|^2 - \|\sqrt{u_{sy}}\mathbf{q}_{xy}w\|^2_{y=0} + L \, \|\mathbf{q}_{xx}w_y\|^2\nonumber\\
&\quad  + L^{\f12 \alpha} [\Phi]_{3,\hat{w}} + L^{\f14}[\Phi]_{4,\hat{w}} .
\end{align}

It is clear to note that
\begin{align}\label{3.18-1}
&\Big\la\mu\Big(\frac{b_1}{\rho^{\v}} \Delta_{\v}\Phi_y\Big)_y,\, \mathbf{q}_{xx}w^2\Big\ra + \Big\la\mu\v \Big(\frac{b_2}{\rho^{\v}}   \Delta_{\v}\Phi_x\Big)_x,\, \mathbf{q}_{xx}w^2\Big\ra\nonumber\\
&=\Big\la\mu\Big(\frac{b_1}{\rho_s} \Delta_{\v}\Phi_y\Big)_y,\, \mathbf{q}_{xx}w^2\Big\ra + \Big\la\mu\v \Big(\frac{b_2}{\rho_s}   \Delta_{\v}\Phi_x\Big)_x,\, \mathbf{q}_{xx}w^2\Big\ra \nonumber\\
&\quad + \v^{\f12 N_0} \big\{[\Phi]_{3,\hat{w}} + [\Phi]_{4,\hat{w}}\big\},
\end{align}
where we have used \eqref{2.122-0} and the following fact
\begin{align}\label{3.18-2}
\v	\|\nabla (p, \rho, S) w\|_{L^4} \lesssim C(\f{1}{L}) \v \|\nabla (p, \rho, S) w\|_{H^1}
\lesssim C(\f1L) \|(\phi, p, S) \hat{w}^{\v}_0\| \lesssim C(\f{1}{L}).
\end{align}
Then we need only to control 1th and 2th terms on RHS of \eqref{3.18-1}. Since the proof is very  delicate, we divide it into several substeps.

\noindent{\it Step 2.1. Estimate on $\Big\la\Big(\frac{b_1}{\rho_s} \Phi_{yyy}\Big)_y,\, \mathbf{q}_{xx}w^2\Big\ra$.} Integrating by parts in $y$ and then in $x$, one has 
\begin{align}\label{3.19}
\Big\la\Big(\frac{b_1}{\rho_s} \Phi_{yyy}\Big)_y,\, \mathbf{q}_{xx}w^2\Big\ra 
&=\Big\la \frac{b_1}{\rho_s} \Phi_{xyyy} ,\, \mathbf{q}_{xy}w^2\Big\ra + \Big\la \big(\frac{b_1}{\rho_s}\big)_x \Phi_{yyy} ,\, \mathbf{q}_{xy}w^2\Big\ra \nonumber\\
&\quad - \Big\la \frac{b_1}{\rho_s} \Phi_{yyy} ,\, \mathbf{q}_{xy}w^2\Big\ra_{x=L} + \Big\la  \frac{b_1}{\rho_s} \Phi_{yy} ,\, 2 \mathbf{q}_{xxy}w w_y\Big\ra\nonumber\\
&\quad  + \Big\la  \frac{b_1}{\rho_s} \Phi_{yy} ,\, 2 \mathbf{q}_{xx}(w w_y)_y\Big\ra + \Big\la  \big(\frac{b_1}{\rho_s}\big)_y \Phi_{yy} ,\, 2 \mathbf{q}_{xx}w w_y\Big\ra
\end{align}
Integrating by parts in $y$, we get
\begin{align}\label{3.20}
	\eqref{3.19}_1
	&=-\Big\la \frac{b_1}{\rho_s} \Phi_{xyy} ,\, \mathbf{q}_{xyy}w^2\Big\ra - \Big\la \frac{b_1}{\rho_s} \Phi_{xyy} ,\, 2 \mathbf{q}_{xy}w w_y\Big\ra - \Big\la \big(\frac{b_1}{\rho_s}\big)_y \Phi_{xyy} ,\, \mathbf{q}_{xy}w^2\Big\ra\nonumber\\
	&\quad - \Big\la \frac{b_1}{\rho_s} \Phi_{xyy} ,\, \mathbf{q}_{xy}w^2\Big\ra_{y=0}.
\end{align}
Noting $u_s|_{y=0}=u_{sx}|_{y=0}=0$ and $\mathbf{q}|_{y=0}=\mathbf{q}_{x}|_{y=0}=0$, we have from \eqref{A.34} that 
\begin{align}\label{3.13}
\Phi_{xyy}\big|_{y=0}
&=\big( 2[m_x u_{sy} + m u_{sxy}] \mathbf{q}_y  + 2 m u_{sy} \mathbf{q}_{xy} \big)(x,0).
\end{align}
It holds that 
\begin{align}\label{3.23}
	u_{sy}&=u_{py}^0 + u_{ey}^0 + \sqrt{\v} (\cdots)_y 
	=u^0_{py} + \sqrt{\v} u_{eY}^0 + \sqrt{\v} (\cdots)_y,
\end{align}
which, with the smallness of $\v>0$, yields that
\begin{align}\label{3.24}
	u_{sy} \big|_{y=0}= \bar{u}_{py}\big|_{y=0}  + O(1)\sqrt{\v}
	=a+ O(1)\sqrt{\v}>0,
\end{align}
Then it follows from \eqref{3.13} and \eqref{3.24} that
\begin{align}\label{3.22}
\eqref{3.20}_4
&= - (2-O(1)L) \|\sqrt{m  b_1 (\rho_s)^{-1} u_{sy}}\mathbf{q}_{xy}w\|_{y=0}^2.
\end{align}

We have from  \eqref{A.34} that 
\begin{align}\label{3.25}
\eqref{3.20}_1
&=   -\Big\la \frac{b_1}{\rho_s} [\bar{u}_{sxyy}\mathbf{q}+ \bar{u}_{sx}\mathbf{q}_{yy}+ \bar{u}_{syy}\mathbf{q}_x+ 2 \bar{u}_{sxy}\mathbf{q}_y + 2\bar{u}_{sy}\mathbf{q}_{xy}] ,\, \mathbf{q}_{xyy}w^2\Big\ra\nonumber\\
&\quad -\|\sqrt{ b_1 (\rho_s)^{-1} \bar{u}_s} \mathbf{q}_{xyy}w\|^2.
\end{align}
Similar as in \eqref{3.23}, one can obtain
\begin{align*}
	u_{sxyy}=u_{pxyy}^0 + \sqrt{\v}^{2} u^0_{exYY} + \sqrt{\v}  (\cdots)_y,
\end{align*}
which, together with Lemma \ref{lemA.3},  yields that 
\begin{align}
	-\Big\la \frac{b_1}{\rho_s} \bar{u}_{sxyy}\mathbf{q} ,\, \mathbf{q}_{xyy}w^2\Big\ra
	&\lesssim L \|\sqrt{u_s}\mathbf{q}_{xyy}w\|\cdot \|\frac{\mathbf{q}_x}{y}w\|
	 \lesssim  L^{-} [\Phi]_{3,\hat{w}},  \label{3.26}\\
-\Big\la \frac{b_1}{\rho_s} \bar{u}_{syy}\mathbf{q}_x ,\, \mathbf{q}_{xyy}w^2\Big\ra 
&\lesssim a \|\sqrt{u_s}\mathbf{q}_{xyy}w\|\cdot \|\frac{\mathbf{q}_x}{y}w\|
\lesssim  L^{\f12 \alpha} [\Phi]_{3,\hat{w}},\label{3.28}\\
-\Big\la \frac{b_1}{\rho_s} \bar{u}_{sx}\mathbf{q}_{yy} ,\, \mathbf{q}_{xyy}w^2\Big\ra 
&\lesssim \|\sqrt{u_s}\mathbf{q}_{yy}w\|\cdot \|\sqrt{u_s}\mathbf{q}_{xyy}w\|  \lesssim L [\Phi]_{3,w}.\label{3.29}
\end{align}
 For the remaining two terms on RHS of \eqref{3.25}, integrating by parts in $y$ and using Lemma \ref{lemA.3}, one has
\begin{align}\label{3.30}
	&-\Big\la \frac{b_1}{\rho_s} 2 \bar{u}_{sxy}\mathbf{q}_y,\, \mathbf{q}_{xyy}w^2\Big\ra\nonumber\\
	&=  \Big\la \frac{b_1}{\rho_s} 2 \bar{u}_{sxy}\mathbf{q}_y,\, \mathbf{q}_{xy}w^2\Big\ra_{y=0} +  \Big\la \frac{b_1}{\rho_s} 2 \bar{u}_{sxy}\mathbf{q}_{yy},\, \mathbf{q}_{xy}w^2\Big\ra + \Big\la \big(\frac{b_1}{\rho_s} 2 \bar{u}_{sxy} w^2\big)_y\mathbf{q}_y,\, \mathbf{q}_{xy}w^2\Big\ra\nonumber\\
	&\lesssim L^-\|\sqrt{u_{sy}}\mathbf{q}_{xy}\|^2_{y=0} + L^{\f12-}\||\mathbf{q}|\|_{\hat{w}} \|\mathbf{q}_{xy}w\| + L\|\mathbf{q}_{xy}w\|^2 
	\lesssim L^{\f12-} [\Phi]_{3,\hat{w}},
\end{align}
and 
\begin{align}\label{3.31}
	&-\Big\la \frac{b_1}{\rho_s} 2 \bar{u}_{sy}\mathbf{q}_{xy},\, \mathbf{q}_{xyy}w^2\Big\ra\nonumber\\
	&= \Big\la \frac{b_1}{\rho_s} \bar{u}_{sy}\mathbf{q}_{xy},\, \mathbf{q}_{xy}w^2\Big\ra_{y=0} + \Big\la (\frac{b_1}{\rho_s} \bar{u}_{sy})_y \mathbf{q}_{xy},\, \mathbf{q}_{xy}w^2\Big\ra + \Big\la \frac{b_1}{\rho_s} \bar{u}_{sy}\mathbf{q}_{xy},\, 2\mathbf{q}_{xy}w w_y\Big\ra\nonumber\\
	&= \|\sqrt{ b_1 m (\rho_s)^{-1}  u_{sy}} \mathbf{q}_{xy}w\|^2_{y=0} + O(1) L^{\alpha-} [\Phi]_{3,\hat{w}}.
\end{align}
Combining \eqref{3.25}-\eqref{3.31}, one gets
\begin{align}\label{3.32}
	\eqref{3.20}_1
	&\leq -\|\sqrt{(\rho_s)^{-1}  b_1 m u_s} \mathbf{q}_{xyy}w\|^2 + [1+ L^-] \|\sqrt{(\rho_s)^{-1}  b_1 mu_{sy}} \mathbf{q}_{xy}\|^2_{y=0} +  L^{\f12 \alpha-} [\Phi]_{3,\hat{w}}.
\end{align}
 
For the 2th \& 3th terms on RHS of \eqref{3.20}, using Lemma \ref{lemA.3}, we have 
\begin{align}\label{3.33}
\eqref{3.20}_{2,3}
&\lesssim \|\Phi_{xyy}w\|\cdot \|\mathbf{q}_{xy}\hat{w}\| 
 \lesssim L^{\f12 \alpha-} [\Phi]_{3,\hat{w}} 
\end{align}
which, together with \eqref{3.20}, \eqref{3.22} and \eqref{3.32}, yields that 
\begin{align}\label{3.34}
\eqref{3.19}_1
&\leq -\|\sqrt{(\rho_s)^{-1}  b_1 m u_s} \mathbf{q}_{xyy}w\|^2  -(1-L^{-})\|\sqrt{(\rho_s)^{-1}  b_1 m u_{sy}}\mathbf{q}_{xy}w\|_{y=0}^2   + L^{\f12 \alpha-} [\Phi]_{3,\hat{w}}.
\end{align}

For the 2th, 5th and 6th terms on RHS of \eqref{3.19}, similarly as \eqref{3.33},  we obtain that 
\begin{align}\label{3.35}
\eqref{3.19}_{2,5,6}	&\lesssim \|\Phi_{yyy}w\|\cdot \|\mathbf{q}_{xy}w\| + \|\Phi_{yy}\hat{w}\|\cdot \|\mathbf{q}_{xx}w_y\|  \lesssim  L^{\f18} [\Phi]_{3,w}  +  L \|\mathbf{q}_{xx}w_y\|^2 .
\end{align}

For the 3th term on RHS of \eqref{3.19}, we note from \eqref{A.6} that 
\begin{align}\label{3.36}
	&- \Big\la \frac{b_1}{\rho_s} \Phi_{yyy} ,\, \mathbf{q}_{xy}w^2\Big\ra_{x=L}
	=\Big\la \frac{b_1}{\rho_s} \Phi_{yyy} ,\, \Big(\frac{\bar{u}_{sx}}{\bar{u}_s}\mathbf{q}\Big)_y w^2\Big\ra_{x=L}\nonumber\\
	&=\Big\la \frac{b_1}{\rho_s} \Phi_{xyyy} ,\, \Big(\frac{\bar{u}_{sx}}{\bar{u}_s}\mathbf{q}\Big)_y w^2\Big\ra + \Big\la (\frac{b_1}{\rho_s})_x \Phi_{yyy} ,\, \Big(\frac{\bar{u}_{sx}}{\bar{u}_s}\mathbf{q}\Big)_y w^2\Big\ra + \Big\la \frac{b_1}{\rho_s} \Phi_{yyy} ,\, \Big(\frac{\bar{u}_{sx}}{\bar{u}_s}\mathbf{q}\Big)_{xy} w^2\Big\ra\nonumber\\
	&=-\Big\la \frac{b_1}{\rho_s} \Phi_{xyy} ,\, \Big(\frac{\bar{u}_{sx}}{\bar{u}_s}\mathbf{q}\Big)_{yy} w^2\Big\ra - \Big\la \big(\frac{b_1}{\rho_s}\big)_y \Phi_{xyy} ,\,  \Big(\frac{\bar{u}_{sx}}{\bar{u}_s}\mathbf{q}\Big)_y w^2\Big\ra - \Big\la \frac{b_1}{\rho_s} \Phi_{xyy} ,\,  \Big(\frac{\bar{u}_{sx}}{\bar{u}_s} \mathbf{q}\Big)_y w w_y\Big\ra \nonumber\\
	&\quad  + \Big\la \big(\frac{b_1}{\rho_s}\big)_x \Phi_{yyy} ,\,  \Big(\frac{\bar{u}_{sx}}{\bar{u}_s} \mathbf{q}\Big)_y w^2\Big\ra + \Big\la \frac{b_1}{\rho_s} \Phi_{yyy} ,\,  \Big(\frac{\bar{u}_{sx}}{\bar{u}_s} \mathbf{q}\Big)_{xy}w^2\Big\ra - \Big\la \frac{b_1}{\rho_s} \Phi_{xyy} ,\,  \Big(\frac{\bar{u}_{sx}}{\bar{u}_s} \mathbf{q}\Big)_y w^2\Big\ra_{y=0}.
\end{align}
Noting, for $y\in [0,1]$
\begin{align*}
	\big(\frac{\bar{u}_{sx}}{\bar{u}_s}\big)(x,y)
	&=\frac{\bar{u}_{sxy}(x,0)+ \cdots + \frac{1}{k!}\pa_{y}^k\bar{u}_{sx}(x,0)\,y^{k-1} + \frac{y^{k}}{k!} \int_0^1 \pa_{y}^{k+1} \bar{u}_{sx}(x,\theta y)(1-\theta)^{k-1}d\theta}{\bar{u}_{sy}(x,0)+ \cdots + \frac{1}{k!}\pa_{y}^k \bar{u}_{s}(x,0)\,y^{k-1} + \frac{y^{k}}{k!} \int_0^1 \pa_{y}^{k+1} \bar{u}_{s}(x,\theta y)(1-\theta)^{k-1}d\theta},
\end{align*}
which yields that $\dis \pa_{y}^l\big(\frac{\bar{u}_{sx}}{\bar{u}_s}\big)$ does not create singularity near $y=0$. Hence we have 
\begin{align}\label{3.38}
\begin{split}
\eqref{3.36}_{1,2,3}	&\lesssim  \|\Phi_{xyy}w\|\Big\{\|\mathbf{q}_{yy}w\|+\|\mathbf{q}_y\hat{w}\| + L \|\frac{\mathbf{q}_x}{y}\hat{w}\| + L \|\sqrt{\v}\mathbf{q}_x\hat{w}\|\Big\} 
\lesssim L^{\f12-} [\Phi]_{3,\hat{w}},\\
\eqref{3.36}_{4,5}	&\lesssim  \|\Phi_{yyy}w\| \Big\{\|\mathbf{q}_{xy}w\|+\|\mathbf{q}_y\hat{w}\|+  \|\frac{\mathbf{q}_x}{y}\hat{w}\| + \|\sqrt{\v}\mathbf{q}_x\hat{w}\| \Big\} \lesssim  L^{\f18} [\Phi]_{3,\hat{w}}, 
\end{split}
\end{align}
and
\begin{align}\label{3.40}
\eqref{3.36}_6 
	&=-\Big\la \frac{b_1}{\rho_s}  [2\bar{u}_{sxy}\mathbf{q}_y + 2\bar{u}_{sy}\mathbf{q}_{xy}] ,\, \frac{\bar{u}_{sx}}{\bar{u}_s}\, \mathbf{q}_y\, w^2\Big\ra_{y=0}
	 \lesssim L^{-} \|\sqrt{u_{sy}}\mathbf{q}_{xy}w\|_{y=0}^2 \lesssim L [\Phi]_{3,w}.
\end{align}
Thus it follows from  \eqref{3.36}-\eqref{3.40} that
\begin{align}\label{3.41}
\eqref{3.19}_3
&\lesssim  L^{\f18} [\Phi]_{3,\hat{w}}.
\end{align}

For the 4th term on RHS of \eqref{3.19}, integrating by parts in $x$, one gets
\begin{align}\label{3.42}
\eqref{3.19}_4
&= - \Big\la  \frac{b_1}{\rho_s} \Phi_{xyy} ,\, 2 \mathbf{q}_{xy}w w_y\Big\ra - \Big\la  (\frac{b_1}{\rho_s})_x \Phi_{yy} ,\, 2 \mathbf{q}_{xy}w w_y\Big\ra  + \Big\la  \frac{b_1}{\rho_s} \Phi_{yy} ,\, 2 \mathbf{q}_{xy}w w_y\Big\ra_{x=L}\nonumber\\
&\lesssim \|\Phi_{xyy}w\|\cdot \|\mathbf{q}_{xy}\hat{w}\| + \|\Phi_{yy}w\|_{x=L}\cdot\|\mathbf{q}_{xy} w_y\|_{x=L} 
\lesssim L^{\f12 \alpha} [\Phi]_{3,\hat{w}}.
\end{align}

Substituting \eqref{3.34}-\eqref{3.35} and  \eqref{3.41}-\eqref{3.42} into \eqref{3.19}, we can obtain that 
\begin{align}\label{3.43}
\Big\la\Big(\frac{b_1}{\rho_s} \Phi_{yyy}\Big)_y,\, \mathbf{q}_{xx}w^2\Big\ra
&\leq -\|\sqrt{(\rho_s)^{-1}  b_1 m u_s} \mathbf{q}_{xyy}w\|^2  -(1- L^{-})\|\sqrt{ (\rho_s)^{-1}  b_1 m u_{sy}}\mathbf{q}_{xy}w\|_{y=0}^2  \nonumber\\
&\quad  +  L^{-} \|\mathbf{q}_{xx}w_y\|^2 + L^{\f12 \alpha-} [\Phi]_{3,\hat{w}} . 
\end{align}

\noindent{\it Step 2.2. Estimate on $\Big\la \v \Big(\frac{b_1}{\rho_s} \Phi_{xxy}\Big)_y,\, \mathbf{q}_{xx}w^2\Big\ra$ and $\dis \Big\la \v \Big(\frac{b_2}{\rho_s} \Phi_{xyy}\Big)_x,\, \mathbf{q}_{xx}w^2\Big\ra$.} Integrating by parts in $y$, one has 
\begin{align}\label{3.44}
&\Big\la \v \Big(\frac{b_1}{\rho_s} \Phi_{xxy}\Big)_y,\, \mathbf{q}_{xx}w^2\Big\ra =-\Big\la \v \frac{b_1}{\rho_s} \Phi_{xxy},\, \mathbf{q}_{xxy}w^2\Big\ra - \Big\la \v \frac{b_1}{\rho_s} \Phi_{xxy},\, \mathbf{q}_{xx}2w w_y \Big\ra\nonumber\\
&=- \|\sqrt{ (\rho_s)^{-1} b_1m u_s} \sqrt{\v}\mathbf{q}_{xxy}w\|^2 - \Big\la \v \frac{b_1}{\rho_s} \Phi_{xxy},\, \mathbf{q}_{xx}2w w_y \Big\ra \nonumber\\
&\quad  -\Big\la \v \frac{b_1}{\rho_s} \big[ \bar{u}_{sxxy} \mathbf{q} + \bar{u}_{sxx}\mathbf{q}_{y} + \bar{u}_{sy} \mathbf{q}_{xx} +2 \bar{u}_{sxy} \mathbf{q}_x + 2\bar{u}_{sx} \mathbf{q}_{xy}\big],\, \mathbf{q}_{xxy}w^2\Big\ra. 
\end{align}
A direct calculation shows that 
\begin{align}\label{3.48}
	\begin{split}
		- \Big\la \v \frac{b_1}{\rho_s} \Phi_{xxy},\, \mathbf{q}_{xx}2w w_y \Big\ra 
		&\lesssim \sqrt{\v} \|\mathbf{q}_{xx}w_y\|^2 + \sqrt{\v} [\Phi]_{3,w},\\
		-\Big\la \v \frac{b_1}{\rho_s} 2 \bar{u}_{sxy} \mathbf{q}_x ,\, \mathbf{q}_{xxy}w^2\Big\ra
		&\lesssim \sqrt{\v}\|\sqrt{\v u_s}\mathbf{q}_{xxy}w\| \cdot \|\frac{\mathbf{q}_x}{y}w\| \lesssim \sqrt{\v} [\Phi]_{3,w},\\
	\Big\la \v \frac{b_1}{\rho_s}\bar{u}_{sxxy} \mathbf{q} ,\, \mathbf{q}_{xxy}w^2\Big\ra
& \lesssim  \sqrt{\v}\|\sqrt{\v u_s}\mathbf{q}_{xxy}w\|\cdot   L\|\frac{\mathbf{q}_x}{y}w\|
\lesssim \sqrt{L\v} [\Phi]_{3,w},\\
-\Big\la \v \frac{b_1}{\rho_s} \big[\bar{u}_{sxx}\mathbf{q}_{y} + 2\bar{u}_{sx} \mathbf{q}_{xy}\big],\, \mathbf{q}_{xxy}w^2\Big\ra & \lesssim \sqrt{\v}\|\sqrt{\v u_s}\mathbf{q}_{xxy}w\| \cdot \|\mathbf{q}_{xy}w\| \lesssim \sqrt{\v} [\Phi]_{3,w}.
\end{split}
\end{align}
Integrating by parts in $y$, we have 
\begin{align}\label{3.50}
	-\Big\la \v \frac{b_1}{\rho_s}  \bar{u}_{sy} \mathbf{q}_{xx},\, \mathbf{q}_{xxy}w^2\Big\ra
	&= \f12 \Big\la \v \big( \frac{b_1}{\rho_s} \bar{u}_{sy} w^2\big)_y \mathbf{q}_{xx}^2,\, w^2\Big\ra 
	\lesssim  L^{\alpha-}  [\Phi]_{3,\hat{w}}.
\end{align}

Combining \eqref{3.44}-\eqref{3.50}, we obtain
\begin{align}\label{3.51}
\Big\la \v \Big(\frac{b_1}{\rho_s} \Phi_{xxy}\Big)_y,\, \mathbf{q}_{xx}w^2\Big\ra
&\leq - \|\sqrt{(\rho_s)^{-1} b_1m u_s} \sqrt{\v}\mathbf{q}_{xxy}w\|^2 
 + \sqrt{\v} \|\mathbf{q}_{xx}w_y\|^2 +  L^{\alpha-} [\Phi]_{3,\hat{w}}.
\end{align}
 
Similar as \eqref{3.51}, one also has
\begin{align}\label{3.52}
 \Big\la \v \Big(\frac{b_2}{\rho_s} \Phi_{xyy}\Big)_x,\, \mathbf{q}_{xx}w^2\Big\ra
&\leq - \|\sqrt{(\rho_s)^{-1} b_2m u_s}\sqrt{\v} \mathbf{q}_{xxy}w\|^2 + \sqrt{\v} \|\mathbf{q}_{xx}w_y\|^2 +  L^{\alpha-} [\Phi]_{3,\hat{w}}.
\end{align}

\smallskip 

\noindent{\it Step 2.3. Estimate on $\Big\la \v^2 \Big(\frac{b_2}{\rho_s} \Phi_{xxx}\Big)_x,\, \mathbf{q}_{xx}w^2\Big\ra$.} Integrating by parts in $x$, one has 
\begin{align}\label{3.53}
	&\Big\la \v^2 \Big(\frac{b_2}{\rho_s} \Phi_{xxx}\Big)_x,\, \mathbf{q}_{xx}w^2\Big\ra
	=-\Big\la \v^2  \frac{b_2}{\rho_s} \Phi_{xxx},\, \mathbf{q}_{xxx}w^2\Big\ra
	- \Big\la \v^2  \frac{b_2}{\rho_s} \Phi_{xxx},\, \mathbf{q}_{xx}w^2\Big\ra_{x=0}.
\end{align}
It is clear that 
\begin{align}\label{3.54}
\eqref{3.53}_1	&=-\Big\la \v^2  \frac{b_2}{\rho_s} \big[\bar{u}_s \mathbf{q}_{xxx} + 3 \bar{u}_{sx}\mathbf{q}_{xx} + 3\bar{u}_{sxx}\mathbf{q}_x + \bar{u}_{sxxx}\mathbf{q}\big],\, \mathbf{q}_{xxx}w^2\Big\ra\nonumber\\
	&=  - \|\sqrt{(\rho_s)^{-1} b_2m u_s} \v \mathbf{q}_{xxx}w\|^2 + O(1)\v^{\f12-} [\Phi]_{3,\hat{w}}.
\end{align}
For  2th term on RHS of \eqref{3.53}, we note from \eqref{2.121} that 
\begin{align*}
	\mathbf{q}_{xx}\big|_{x=0} = -2 \frac{\bar{u}_{sx}}{\bar{u}_s} \mathbf{q}_x\big|_{x=0},  
\end{align*}
which, together with Lemma \ref{lemA.4}, yields that 
\begin{align}\label{3.56} 
\eqref{3.53}_2	&\lesssim \|\v^{\f32}\Phi_{xxx}w\|_{x=0}\cdot \|\sqrt{\v}\mathbf{q}_xw\|_{x=0}
	\lesssim L^{\f12-}  [\Phi]_{3,w} +  L^{\f12-}   \|\v^2 \Phi_{xxxx}w\|.
\end{align}
Hence, substituting \eqref{3.54} and \eqref{3.56} into \eqref{3.53}, we obtain
\begin{align}\label{3.57}
\Big\la \v^2 \Big(\frac{b_2}{\rho_s} \Phi_{xxx}\Big)_x,\, \mathbf{q}_{xx}w^2\Big\ra
&\leq - \|\sqrt{(\rho_s)^{-1} b_2m u_s} \, \v \mathbf{q}_{xxx}w\|^2  + L^{\f12-} [\Phi]_{3,\hat{w}}    + 
L^{\f12-} \|\v^2\Phi_{xxxx}w\|^2.
\end{align}

\medskip

\noindent{\it Step 2.4.} Combining \eqref{3.57}, \eqref{3.51}-\eqref{3.52} and \eqref{3.43}, we can conclude \eqref{3.18}. 

\medskip

\noindent{\it Step 3.} Substituting \eqref{3.12}, \eqref{3.14} and \eqref{3.18} into \eqref{3.2},  we conclude \eqref{3.1}. 

\medskip

\noindent{\it Step 4. Estimate on the $\big\la \mathcal{F}_{R}(\hat{P}, \hat{\mathbf{S}}, \zeta, p, \phi, S),\, \mathbf{q}_{xx}w^2\big\ra$.}  We still need to control the remainder term on RHS of \eqref{3.1}. Since the proof is very complicate, we divide it into several substeps.

\smallskip

\noindent{\it Step 4.1. Estimate on $\big\la \mathcal{F}_{R1},\, \mathbf{q}_{xx}w^2  \big\ra$.} Recall $\mathcal{F}_{R1}$ in $\eqref{2.98}$. Integrating by parts in $y$, one gets that 
\begin{align}\label{3.50-1}
\big\la \mathcal{F}_{R1},\, \mathbf{q}_{xx}w^2  \big\ra 
&=\big\la  \f{\rho_s}{T_s}b_1  [u_s^2 + O(1) u_s \v ]\zeta_x ,\, \mathbf{q}_{xxy}w^2  \big\ra + \big\la  b_1 \f{\rho_s}{T_s} [u_s^2 + O(1) u_s \v ]\zeta_x ,\, \mathbf{q}_{xx}2ww_y  \big\ra\nonumber\\
&=\big\la \f{b_1}{T_s} \rho_s u_s^2 \zeta_x ,\, \mathbf{q}_{xxy}w^2  \big\ra   + L^{\f14} \|\mathbf{q}_{xx}w_y\|^2  +  L^{\f15} \big\{[[[\zeta]]]_{2,\hat{w}} +  [\Phi]_{3, \hat{w}}\big\}.
\end{align} 

For 1th term on RHS of \eqref{3.50-1}, we use  \eqref{2.111-1} to have
\begin{align}\label{3.50-2}
\Big\la \f{b_1}{T_s}  \rho_s u_s^2 \zeta_x ,\, \mathbf{q}_{xxy}w^2  \Big\ra 
&= \f{\kappa}2 \Big\la \v  \f{b_1}{T_s} u_s \zeta_{xx} ,\, \mathbf{q}_{xxy}w^2  \Big\ra  + \f{\kappa}2 \Big\la \f{b_1}{T_s} u_s \zeta_{yy} ,\, \mathbf{q}_{xxy}w^2  \Big\ra  \nonumber\\
&\quad +  \f{\kappa}2 \Big\la \f{b_1 u_s}{T_s} \Delta_{\v}\big(\f{1}{4\rho_s}p\big),\, \mathbf{q}_{xxy}w^2  \Big\ra + \f{\kappa}2 \Big\la \f{b_1u_s}{T_s}  \Delta_{\v} \big(\f{1}{p_s}T_{sy}q\big),\, \mathbf{q}_{xxy}w^2  \Big\ra  \nonumber\\
&\quad  - \Big\la \f{b_1}{T_s}  \rho_s u_s v_s \zeta_y ,\, \mathbf{q}_{xxy}w^2  \Big\ra + \f12 \Big\la \f{b_1}{T_s}  u_s \mathfrak{J} ,\, \mathbf{q}_{xxy}w^2  \Big\ra .
\end{align}
It is clear that
\begin{align}\label{3.52-1}
\eqref{3.50-2}_1	&\lesssim \sqrt{a}\|\sqrt{\v}\zeta_{xx}w\|\cdot \||\mathbf{q}|\|_w \lesssim \f{1}{N} [\Phi]_{3,w} + N a [[[\zeta]]]_{2,w}.
\end{align} 
By noting $\dis\zeta|_{x=0}=-\f{1}{2\rho_s}\chi p$ for NT and $\zeta|_{x=0}=0$ for DT, we integrate by parts to obtain 
\begin{align}\label{3.53-1}
\eqref{3.50-2}_2&=-\Big\la \f{\kappa b_1}{2T_s} u_s \zeta_{xyy} ,\, \mathbf{q}_{xy}w^2 \Big\ra - \Big\la \big(\f{\kappa b_1}{2T_s} u_s\big)_x \zeta_{yy} ,\, \mathbf{q}_{xy}w^2 \Big\ra + \Big\la \f{\kappa b_1}{2T_s} u_s \zeta_{yy} ,\, \mathbf{q}_{xy}w^2 \Big\ra_{x=L} \nonumber\\
&\quad - \Big\la \f{\kappa b_1}{2T_s} u_s \zeta_{yy} ,\, \mathbf{q}_{xy}w^2 \Big\ra_{x=0}\nonumber\\
&=\Big\la \f{\kappa b_1}{2T_s} u_s \zeta_{xy} ,\, \mathbf{q}_{xyy}w^2 \Big\ra + \Big\la \big(\f{\kappa b_1}{2T_s} u_s w^2\big)_y \zeta_{xy} ,\, \mathbf{q}_{xy}w^2 \Big\ra  - \Big\la \big(\f{\kappa b_1}{2T_s}u_s\big)_x \zeta_{yy} ,\, \mathbf{q}_{xy}w^2 \Big\ra  \nonumber\\
&\quad + \Big\la \f{\kappa b_1}{2T_s} u_s \zeta_{yy} ,\, \mathbf{q}_{xy}w^2\Big\ra_{x=L} - \Big\la \f{\kappa b_1}{2T_s} u_s \zeta_{yy} ,\, \mathbf{q}_{xy}w^2 \Big\ra_{x=0}\nonumber\\
&\lesssim \f{1}{N} [\Phi]_{3,\hat{w}} + Na \f{1}{\v} [[p]]_{2,1} + L^{\f12-} [[p]]_{3,1} + N a [[[\zeta]]]_{2,w} +  L^{\f14-} [[[\zeta]]]_{b,w},
\end{align}
where we have used the following estimate for the case of NT
\begin{align}\label{3.53-2}
&-\mathbf{1_N}\, \Big\la \f{b_1}{T_s} u_s \zeta_{yy} ,\, \mathbf{q}_{xy}w^2  \Big\ra_{x=0}
=\Big\la \f{b_1}{T_s} u_s \big(\f{\chi}{2\rho_s}p\big)_{yy} ,\, \mathbf{q}_{xy}w^2  \Big\ra_{x=0} \nonumber\\
&=\Big\la \f{b_1}{2p_s} u_sp_{yy} \chi,\, \mathbf{q}_{xy}w^2  \Big\ra_{x=0} + \Big\la \f{b_1}{T_s}\big(\f{\chi}{2\rho_s}\big)_y u_s  p_{y} ,\, \mathbf{q}_{xy}w^2  \Big\ra_{x=0} + \Big\la \f{b_1}{T_s}\big(\f{\chi}{2\rho_s}\big)_{yy}  u_s  p ,\, \mathbf{q}_{xy}w^2  \Big\ra_{x=0}\nonumber\\
&=-\Big\la \f{b_1}{2p_s} u_sp_{xyy} \chi,\, \mathbf{q}_{xy}w^2  \Big\ra - \Big\la \f{b_1}{2p_s} u_sp_{yy} \chi,\, \mathbf{q}_{xxy}w^2  \Big\ra - \Big\la \big(\f{b_1}{2p_s} u_s\big)_x p_{yy} \chi,\, \mathbf{q}_{xy}w^2  \Big\ra \nonumber\\
&\quad + \Big\la \f{b_1}{2p_s} u_sp_{yy} \chi,\, \mathbf{q}_{xy}w^2  \Big\ra_{x=L} + O(1) \|u_s\mathbf{q}_{xy}w\|_{x=0}\cdot \|\chi p_y\|_{x=0}\nonumber\\
&=\Big\la \f{b_1}{2p_s} u_sp_{xy} \chi,\, \mathbf{q}_{xyy}w^2  \Big\ra + \Big\la \big(\f{b_1}{2p_s} u_s \chi\big)_y p_{xy},\, \mathbf{q}_{xy}w^2  \Big\ra + \Big\la \f{b_1}{2p_s} u_sp_{xy} \chi,\, 2\mathbf{q}_{xy}ww_y  \Big\ra\nonumber\\
&\quad + O(1)\sqrt{a}\big(\|\sqrt{u_s}\sqrt{\v}\mathbf{q}_{xxy}\| + \sqrt{\v}\|\mathbf{q}_{xy}\|\big) \cdot \|\f{1}{\sqrt{\v}}p_{yy}w\| +  O(1) \|u_s\mathbf{q}_{xy}w\|_{x=0}\cdot \|\chi p_y\|_{x=0} \nonumber\\
&\quad + O(1) \sqrt{a}\|\sqrt{u_s}p_{yy}w\|_{x=L}\cdot \|\mathbf{q}_{xy}w\|_{x=L}\nonumber\\
&\lesssim \f{1}{N} [\Phi]_{3,\hat{w}}  + Na \f{1}{\v} [[p]]_{2,1} + L^{\f12-} [[p]]_{3,1}.
\end{align}

A direct calculation shows that 
\begin{align}\label{3.53-3} 
\eqref{3.50-2}_3
&=\Big\la \f{\kappa b_1}{8p_s} u_s  \Delta_{\v}p,\, \mathbf{q}_{xxy}w^2  \Big\ra +  \Big\la \f{\kappa b_1}{2T_s} u_s \nabla_{\v}(\f{1}{2\rho_s})\cdot \nabla_{\v}p,\, \mathbf{q}_{xxy}w^2  \Big\ra \nonumber\\
&\quad  +  \Big\la \f{\kappa b_1}{2T_s}\Delta_{\v}(\f{1}{4\rho_s})\, u_s p,\, \mathbf{q}_{xxy}w^2  \Big\ra \nonumber\\
&\lesssim  \f1N [\Phi]_{3,w} + Na\f{1}{\v} [[p]]_{2,\hat{w}},
\end{align}
where we have used Lemmas \ref{lemA.4}-\ref{lemA.5-1} to derive the following fact
\begin{align}\label{3.54-2}
\Big\la \f{b_1}{T_s}\Delta_{\v}(\f{1}{4\rho_s})\, u_s p,\, \mathbf{q}_{xxy}w^2  \Big\ra 
&=-\Big\la \f{b_1}{T_s}\Delta_{\v}(\f{1}{4\rho_s})\, u_s p_x,\, \mathbf{q}_{xy}w^2  \Big\ra -  \Big\la \Big(\f{b_1}{T_s}\Delta_{\v}(\f{1}{4\rho_s})\, u_s\Big)_x p,\, \mathbf{q}_{xy}w^2  \Big\ra \nonumber\\
&\quad +  \Big\la \f{b_1}{T_s}\Delta_{\v}(\f{1}{4\rho_s})\, u_s p,\, \mathbf{q}_{xy}w^2  \Big\ra_{x=L} -  \Big\la \f{b_1}{T_s}\Delta_{\v}(\f{1}{4\rho_s})\, u_s p,\, \mathbf{q}_{xy}w^2  \Big\ra_{x=0}\nonumber\\
&\lesssim  \|u_s\mathbf{q}_{xy}w\| \big\{\|\mathbf{d}_3 p_xw\| + \|\mathbf{d}_3 p w\|\big\} +  \|\mathbf{q}_{xy}w\|_{x=L} \|\mathbf{d}_3 p w\|_{x=L} \nonumber\\
&\quad + \|u_s\mathbf{q}_{xy}w\|_{x=0} \|\mathbf{d}_3 p w\|_{x=0} \nonumber\\
&\lesssim L^{\f18} [\Phi]_{3,\hat{w}} + L^{\f18} \f{1}{\v} [[p]]_{2,1}. 
\end{align}
We point out that the following estimate has been used in \eqref{3.54-2}
\begin{align}\label{3.72-2}
\|\mathbf{d}_3\, p(1+w_0)\|^2 \lesssim \|\mathbf{d}_2\, p\|^2 \lesssim L^2 \f{1}{\v} [[p]]_{2,1}.
\end{align}

It holds that 
\begin{align}\label{3.72-3}
\eqref{3.50-2}_4&= \Big\la \f{\kappa b_1}{2T_s} u_s \big(\f{1}{p_s}T_{sy}q\big)_{yy},\, \mathbf{q}_{xxy}w^2  \Big\ra + L^{\f12 \alpha-} [\Phi]_{3,w} + L^{\f12 \alpha-} [\phi]_{3,w} \nonumber\\
&=-\Big\la \f{\kappa b_1}{2T_s} u_s \big(\f{1}{p_s}T_{sy}q\big)_{xyy},\, \mathbf{q}_{xy}w^2  \Big\ra - \Big\la \big(\f{\kappa b_1}{2T_s} u_s\big)_x \big(\f{1}{p_s}T_{sy}q\big)_{yy},\, \mathbf{q}_{xy}w^2  \Big\ra \nonumber\\
&\quad + \Big\la \f{\kappa b_1}{2T_s} u_s \big(\f{1}{p_s}T_{sy}q\big)_{yy},\, \mathbf{q}_{xy}w^2  \Big\ra_{x=L} + L^{\f12 \alpha-} [\Phi]_{3,\hat{w}} + L^{\f12 \alpha-} [\phi]_{3, \hat{w}} \nonumber\\
&\lesssim L^{\f12 \alpha-} [\Phi]_{3,\hat{w}} + L^{\f12 \alpha-} [\phi]_{3, \hat{w}},
\end{align}
and
\begin{align}\label{3.54-1}
\eqref{3.50-2}_5&=  \Big\la \f{\kappa b_1}{2T_s}  \rho_s u_s v_s \zeta_{xy} ,\, \mathbf{q}_{xy}w^2  \big\ra + \Big\la \big(\f{\kappa b_1}{2T_s}  \rho_s u_s v_s \big)_x \zeta_y ,\, \mathbf{q}_{xy}w^2  \Big\ra \nonumber\\
&\quad - \Big\la \f{\kappa b_1}{2T_s}  \rho_s u_s v_s \zeta_y ,\, \mathbf{q}_{xy}w^2  \Big\ra_{x=L} + \Big\la \f{b_1}{T_s}  \rho_s u_s v_s \zeta_y ,\, \mathbf{q}_{xy}w^2  \Big\ra_{x=0}\nonumber\\
&\lesssim  L^{\f18} [\Phi]_{3,\hat{w}} + L^{\f18}[[[\zeta]]]_{2,\hat{w}}.
\end{align}

For term involving $\mathfrak{J}$, we note that 
\begin{align}\label{3.56-1}
&-\Big\la \f{b_1}{T_s}  u_s \Big\{2\rho_sv_s \big(\f{1}{p_s}T_{sy}q\big)_y +  2\big[\rho_su_s \big(\f{1}{p_s}T_{sy}\big)_x-  T_{sy} \bar{u}_{sx}\big]q\Big\} ,\, \mathbf{q}_{xxy}w^2  \Big\ra\nonumber\\
&=\Big\la \pa_x\Big\{\f{b_1}{T_s}  u_s \Big(2\rho_sv_s \big(\f{1}{p_s}T_{sy}q\big)_y + 2\big[\rho_su_s \big(\f{1}{p_s}T_{sy}\big)_x-  T_{sy} \bar{u}_{sx}\big]q\Big)\Big\} ,\, \mathbf{q}_{xy}w^2  \Big\ra \nonumber\\
&\quad - \Big\la \f{b_1}{T_s}  u_s \Big(2\rho_sv_s \big(\f{1}{p_s}T_{sy}q\big)_y +  2\big[\rho_su_s \big(\f{1}{p_s}T_{sy}\big)_x-  T_{sy} \bar{u}_{sx}\big]q\Big) ,\, \mathbf{q}_{xy}w^2  \Big\ra_{x=L}\nonumber\\
&\lesssim L^{\alpha} [\Phi]_{3,\hat{w}} + L^{\alpha} [\phi]_{3,\hat{w}},
\end{align}
and
\begin{align}\label{3.57-2}
\Big\la \f{b_1}{p_s}  u_s p_{sy}  \phi_x ,\, \mathbf{q}_{xxy}w^2  \Big\ra 
&=-\Big\la \f{b_1}{p_s}  u_s p_{sy}  \phi_{xx} ,\, \mathbf{q}_{xy}w^2  \Big\ra - \Big\la (\f{b_1}{p_s}  u_s p_{sy})_x  \phi_x ,\, \mathbf{q}_{xy}w^2  \Big\ra\nonumber\\
&\quad  - \Big\la \f{b_1}{p_s}  u_s p_{sy}  \phi_x ,\, \mathbf{q}_{xy}w^2  \Big\ra_{x=0} \nonumber\\
&\lesssim \sqrt{L} [\Phi]_{3,\hat{w}} + \sqrt{L} [\phi]_{3,\hat{w}},
\end{align}
where we have used the fact $p_{sy}\cong \sqrt{\v}$.

Integrating by parts in $x$, one also has
\begin{align}\label{3.58-2}
&\Big\la  \f{b_1}{T_s}\rho_s u_s (U_s\cdot \nabla) (\frac{1}{\rho_s})\, p ,\, \mathbf{q}_{xxy}w^2  \Big\ra \nonumber\\
&=-\Big\la  \f{b_1}{T_s}  u_s \rho_s (U_s\cdot \nabla) (\frac{1}{\rho_s})\, p_x ,\, \mathbf{q}_{xy}w^2  \Big\ra - \Big\la  \big\{\f{b_1}{T_s}  \rho_s u_s (U_s\cdot \nabla) (\frac{1}{\rho_s})\big\}_x\, p ,\, \mathbf{q}_{xy}w^2  \Big\ra\nonumber\\
&\quad + \Big\la  \f{b_1}{T_s} \rho_s  u_s (U_s\cdot \nabla) (\frac{1}{\rho_s})\, p ,\, \mathbf{q}_{xy}w^2  \Big\ra_{x=L} - \Big\la  \f{b_1}{T_s}  \rho_s  u_s(U_s\cdot \nabla) (\frac{1}{\rho_s})\, p ,\, \mathbf{q}_{xy}w^2  \Big\ra_{x=0}\nonumber\\
&\lesssim \|u_s\mathbf{q}_{xy}w\| \big(\|p_x\| + \|\mathbf{d}_2pw\|\big) + \|\mathbf{q}_{xy}w\|_{x=L} \|\mathbf{d}_2 pw\|_{x=L} + \|u_s \mathbf{q}_{xy}w\|_{x=0} \|\mathbf{d}_2 pw\|_{x=0} \nonumber\\
&\lesssim L^{\f18} [\Phi]_{3,\hat{w}} + L^{\f18} \f{1}{\v} [[p]]_{2,\hat{w}}+   L^{2} \f{1}{\v^2} [[p]]_{2,1}\cdot \mathbf{1}_{\{w=w_0\}}.
\end{align}

It follows from  \eqref{7.47} that (for both NT \& DT)
\begin{align}\label{1C.43-1}
	\begin{split}
		|\mathfrak{M}_0|
		&\lesssim \mathbf{d}_2 u_s |(p_y,S_y)| + a \mathbf{d}_2 |(p,S)| + \mathbf{d}_2 |(\phi_{yy}, \phi_y)| + a \mathbf{d}^2_2  \big|(q_y,q)\big|,\\
		|\pa_x\mathfrak{M}_0|
		&\lesssim \mathbf{d}_2u_s |(p_{xy}, S_{xy})| + a \mathbf{d}_2 |(p_x, p_y, S_x,S_y)| + a\mathbf{d}_2|(p,S)| \\
		&\quad   + \mathbf{d}_2 |(\phi_{xyy}, \phi_{xy}, \phi_{yy}, \phi_y)| + u_s \mathbf{d}_2^2 |q_{xy}|+ a \mathbf{d}_2^2 |(q_x,q_y,q)|,
	\end{split}
\end{align}
which, together with \eqref{2.123}, yields that(for both NT \& DT)
\begin{align}\label{3.59-1}
\begin{split}
\|\mathfrak{M}_0w\|^2_{x=0} 
&\lesssim \|p_y \|^2_{x=0} +  \|\mathbf{d}_2 pw\|^2_{x=0}  \lesssim \f{1}{\sqrt{\v}} [[p]]_{2,1} +  \f{1}{\v^{\f32}} [[p]]_{2,1} \cdot \mathbf{1}_{\{w=w_0\}},\\
\|\mathfrak{M}_0w\|^2_{x=L}&\lesssim \|\mathfrak{M}_0w\|^2_{x=0}  + \|(\pa_x\mathfrak{M}_0, \mathfrak{M}_0)w\|^2,
\end{split}
\end{align}
and
\begin{align}\label{1TC.44} 
	\begin{split}
		\|\mathfrak{M}_0w\|^2&\lesssim L [\phi]_{3,1} +  L^2 \f{1}{\v} [[p]]_{2,1} + L [[[S]]]_{2,1} +  L\Big\{\f{1}{\v^2} [[p]]_{2,1} + \f{1}{\v} [[[S]]]_{2,1}\Big\} \mathbf{1}_{\{w=w_0\}}, \\
		\|\pa_x \mathfrak{M}_0w\|^2& \lesssim a^2 \f{1}{\v} [[p]]_{2,1} + a^2 [[[S]]]_{2,1}  
		+ a^2 [\phi]_{3,\hat{w}} +  L \Big\{\f{1}{\v^2} [[p]]_{2,1} + \f{1}{\v} [[[S]]]_{2,1}\Big\} \mathbf{1}_{\{w=w_0\}}.
	\end{split}
\end{align}

We note from \eqref{7.35-2} and \eqref{7.36-2} that (for both NT \& DT)
\begin{align}
	\begin{split}\label{1TC.47}
		|\mathfrak{M}_1|&\lesssim u_s \mathbf{d}_1|(p_y,S_y)| + u_s\mathbf{d}_2 |(p_x,S_x)| + u_s\mathbf{d}_2|(p,S)| + \mathbf{d}_1 |(\nabla^2\phi , \nabla\phi)|  + au_s \mathbf{d}_2 |(q_y,q_x,q)|,\\
		|\mathfrak{M}_2|&\lesssim u_s \mathbf{d}_2|(p_x,p_y,S_x,S_y)| + u_s\mathbf{d}_2^2|(p,S)| + \mathbf{d}_2 |\nabla\phi| + a u_s \mathbf{d}_2 |(q_x,q_y,q)|,
	\end{split}
\end{align}
which implies that (for both NT \& DT)
\begin{align}\label{1TC.48}
	\begin{split}
		\|\sqrt{\v}(\mathfrak{M}_1,\mathfrak{M}_2)w\|^2&\lesssim   L^2 \f{1}{\v} [[p]]_{2,1} + L [[[S]]]_{2,1} + L[\phi]_{3,\hat{w}}.
	\end{split}
\end{align} 

Using \eqref{3.59-1}-\eqref{1TC.44} and \eqref{1TC.48}, we obtains
\begin{align}\label{3.60}
&\Big\la \f{b_1}{T_s}  u_s\mathfrak{M} ,\, \mathbf{q}_{xxy}w^2  \Big\ra 
 \nonumber\\
&=- \Big\la \f{b_1}{T_s}  u_s \pa_x \mathfrak{M}_0 ,\, \mathbf{q}_{xy}w^2  \Big\ra - \Big\la (\f{b_1}{T_s}  u_s)_x \mathfrak{M}_0 ,\, \mathbf{q}_{xy}w^2  \Big\ra  +  \Big\la \f{b_1}{T_s}  u_s \mathfrak{M}_0 ,\, \mathbf{q}_{xy}w^2  \Big\ra_{x=L} \nonumber\\
&\quad  -  \Big\la \mathfrak{M}_0 ,\, \f{b_1}{T_s}  u_s \mathbf{q}_{xy}w^2  \Big\ra_{x=0} + \v \Big\la \mathfrak{M}_1+ \mathfrak{M}_2 ,\, \f{b_1}{T_s}  u_s \mathbf{q}_{xxy}w^2  \Big\ra \nonumber\\
&\lesssim\|(\pa_x \mathfrak{M}_0, \mathfrak{M}_0) \, w\|\cdot \|u_s\mathbf{q}_{xy}w\| +  L^{\f12-} [\Phi]_{3,w}^{\f12} \|\mathfrak{M}_0w\|_{x=L}  + \sqrt{\v} [\Phi]_{3,w}^{\f12}\|(\mathfrak{M}_{1}, \mathfrak{M}_2)w\|  \nonumber\\
&\quad + \|u_s\mathbf{q}_{xy}w\|_{x=0} \cdot \|\mathfrak{M}_0w\|_{x=0} \nonumber\\
&\lesssim L^{\f18} [\Phi]_{3,\hat{w}}  +  L^{\f18} \f{1}{\v} [[p]]_{2,\hat{w}} + L^{\f18} [[[S]]]_{2,\hat{w}} + \f{1}{\v^{\f32}} [[p]]_{2,1} \cdot \mathbf{1}_{\{w=w_0\}} + L^{\f18} [\phi]_{3,\hat{w}}.
\end{align}
Also it is clear that 
\begin{align}\label{3.61-2}
\Big\la \mathfrak{N}_3,\, \f{b_1}{T_s}  u_s \mathbf{q}_{xxy}w^2  \Big\ra  
\lesssim \f{1}{N} [\Phi]_{3,w} + N \f{1}{\v} \|\sqrt{u_s}\mathfrak{N}_3 w\|^2. 
\end{align}
Hence, noting the definition of $\mathfrak{J}$ in \eqref{2.84},  we have from \eqref{3.56-1}-\eqref{3.61-2} that 
\begin{align}\label{3.62-1}
\eqref{3.50-2}_6
& \lesssim \f{1}{N} [\Phi]_{3,\hat{w}}  + L^{\alpha} [\phi]_{3,\hat{w}} + L^{\f18} \f{1}{\v} [[p]]_{2,\hat{w}}+ L^{\f18} [[[S]]]_{2,\hat{w}} + L \f{1}{\v^2} [[p]]_{2,1}\cdot  \mathbf{1}_{\{w=w_0\}} \nonumber\\
&\quad  + N \f{1}{\v} \|\sqrt{u_s}\mathfrak{N}_3 w\|^2 . 
\end{align}

Substituting \eqref{3.52-1}-\eqref{3.54-1} and \eqref{3.62-1} into \eqref{3.50-2}, one obtains that 
\begin{align}\label{3.63-2}
\big\la \f{b_1}{T_s}  \rho_s u_s^2 S_x ,\, \mathbf{q}_{xxy}w^2  \big\ra 
&\lesssim \f{1}{N} [\Phi]_{3,\hat{w}}  + Na [[[\zeta]]]_{2,\hat{w}}+ L^{\f14-} [[[\zeta]]]_{b,w} + Na \f{1}{\v} [[p]]_{2,\hat{w}} + L^{\f12-} [[p]]_{3,1} \nonumber\\
&\,\,\,+  L\f{1}{\v^2} [[p]]_{2,1} \mathbf{1}_{\{w=w_0\}}  + L^{\f18}  [[[S]]]_{2,\hat{w}} + L^{\f12\alpha-}[\phi]_{3,\hat{w}}  +  \f{1}{\v} \|\sqrt{u_s}\mathfrak{N}_3 w\|^2,
\end{align}
which, together with \eqref{3.50-1}, yields that 
\begin{align}\label{3.64-2}
\big\la \mathcal{F}_{R1},\, \mathbf{q}_{xx}w^2  \big\ra 
&\lesssim   L^{\f14} \|\mathbf{q}_{xx}w_y\|^2 + \f{1}{N} [\Phi]_{3,\hat{w}} + Na \|(\phi, p, S\,(\zeta)) \hat{w}\|^2_{\mathbf{X}^{\v}} +  L\f{1}{\v^2} [[p]]_{2,1}\, \mathbf{1}_{\{w=w_0\}} \nonumber\\
&\quad  +  \f{1}{\v} \|\sqrt{u_s}\mathfrak{N}_3 w\|^2.
\end{align}

\medskip

\noindent{\it Step 4.2. Estimate on $\big\la \mathcal{F}_{R2}, \, \mathbf{q}_{xx}w^2 \big\ra $.} Recall $\mathcal{F}_{R2}$ \eqref{2.98} and $\mathfrak{g}_1, \mathfrak{g}_2$ in \eqref{2.93}. It is direct to have that 
\begin{align}\label{3.65-2}
\big\la \v (b_2 \mathfrak{g}_2)_x, \, \mathbf{q}_{xx}w^2 \big\ra \lesssim  L^{\f14}  [\Phi]_{3,\hat{w}} + L^{\f14}  [\phi]_{3,\hat{w}}.
\end{align}
Integrating by parts in $y,x$, one has
\begin{align}\label{3.66-2}
\big\la  (b_1 \mathfrak{g}_1)_y , \, \mathbf{q}_{xx}w^2 \big\ra 
&= \big\la (b_1 \mathfrak{g}_1)_x  , \, \mathbf{q}_{xy}w^2 \big\ra - \big\la b_1 \mathfrak{g}_1  , \, \mathbf{q}_{xx}2w w_y \big\ra  - \big\la b_1 \mathfrak{g}_1  , \, \mathbf{q}_{xy}w^2 \big\ra_{x=L} \nonumber\\
&\quad + \big\la b_1 \mathfrak{g}_1  , \, \mathbf{q}_{xy}w^2 \big\ra_{x=0}\nonumber\\
&\lesssim  L^{\f18} [\Phi]_{3,\hat{w}} + L^{\f12}\|\mathbf{q}_{xx}w_y\|^2 +  L^{\f18} [\phi]_{3,\hat{w}}.
\end{align}
Thus it follows from \eqref{3.65-2}-\eqref{3.66-2} that 
\begin{align}\label{3.67-2}
	\big\la \mathcal{F}_{R2}, \, q_{xx}w^2 \big\ra \lesssim  L^{\f18} [\Phi]_{3,\hat{w}} + L^{\f12}\|\mathbf{q}_{xx}w_y\|^2 +  L^{\f18} [\phi]_{3,\hat{w}}.
\end{align}

\medskip

\noindent{\it Step 4.3. Estimate on $\big\la \mathcal{F}_{R3}, \, \mathbf{q}_{xx}w^2 \big\ra $.} Recall $\mathcal{F}_{R3}$ in \eqref{2.82} and $\dis T=S + \f{1}{2\rho_s}p + \f{1}{p_s}T_{sy} q$, we can rewrite $g_1, g_2$ as
\begin{align}\label{3.66-0}
g_1(p,T,\phi)
&=\Big[\tilde{d}_{12} + \f{1}{2\rho_s}\hat{d}_{12}\Big] p_y + \Big[\tilde{d}_{13} + \f{1}{2\rho_s}\hat{d}_{13} + \big(\f{1}{2\rho_s}\hat{d}_{11}\big)_x + \hat{d}_{12}\,(\f{1}{2\rho_s})_y \Big]\, p\nonumber\\
&\quad + \hat{d}_{12}\, S_y  + \hat{d}_{13}\, S + \hat{d}_{12}\big(\f{1}{p_s}T_{sy} q\big)_y + \hat{d}_{13}\f{1}{p_s}T_{sy} q  + I_u(\phi)   - \mathfrak{R}_{s1},
\end{align}
and
\begin{align}\label{3.66-1}
g_2(p,T,\phi)
&=\v \Big[\tilde{d}_{21} + \f{1}{2\rho_s}\hat{d}_{21} \Big] p_x + \v \hat{d}_{22} p_y + \v \Big\{\tilde{d}_{23} + \f{1}{2\rho_s}\hat{d}_{23} + (\f{1}{2\rho_s})_x \hat{d}_{21} + (\f{1}{2\rho_s})_y \hat{d}_{22}\Big\} p \nonumber\\
&\quad + \v \hat{d}_{21} S_x + \v \hat{d}_{22}\, S_y + \v \hat{d}_{23} S 
 +  \v \hat{d}_{21}\, \big(\f{1}{p_s}T_{sy} q\big)_x  +  \v \hat{d}_{22}\, \big(\f{1}{p_s}T_{sy} q\big)_y \nonumber\\
&\quad +  \v \hat{d}_{23}\, \f{1}{p_s}T_{sy} q  + \v I_{v}(\phi) - \v \mathfrak{R}_{s2}.
\end{align}

Noting \eqref{7.12-1}-\eqref{7.13-1}, we rewrite $I_u(\phi)$ and $I_v(\phi)$ as
\begin{align}\label{3.67-0}
\begin{split}
I_u(\phi)&= \big[v_s + \frac{2\mu}{\rho_s}\rho_{sy}\big]\phi_{yy}  +   [\cdots] \phi_{y}
+ \v  [\cdots] \phi_{xy}  + \v  [\cdots] |(\phi_{xx},\phi_x)|, \\
I_{v}(\phi)&=  \big[\frac{2\mu+\lambda }{\rho_s^2}\rho_{sy} -v_s\big] \phi_{xy} + [\cdots] \v \phi_{xx} + [\cdots] \phi_x + [\cdots] \phi_y.
\end{split}
\end{align}
Integrating by parts in $y$ and then in $x$, one has
\begin{align}\label{3.58}
\big\la(b_1[\cdots]\phi_{yy})_y, \, \mathbf{q}_{xx}w^2 \big\ra
&=-\big\la b_1[\cdots] \phi_{yy}, \, \mathbf{q}_{xy}w^2 \big\ra_{x=L} + \big\la b_1[\cdots] \phi_{xyy}, \, \mathbf{q}_{xy}w^2 \big\ra \nonumber\\
&\quad  + \big\la (b_1[\cdots])_x \phi_{yy}, \, \mathbf{q}_{xy}w^2 \big\ra
- \big\la b_1[\cdots] \phi_{yy}, \, \mathbf{q}_{xx} 2w w_y \big\ra\nonumber\\
&\lesssim \|\phi_{yy}w\|_{x=L}\cdot \|\mathbf{q}_{xy}w\|_{x=L} + \|\phi_{xyy}w\| \big\{\|\mathbf{q}_{xy}w\| + L\|\mathbf{q}_{xx}w_y\|\big\}\nonumber\\
&\lesssim L^{\f12\alpha-} [\Phi]_{3,\hat{w}} + L \|\mathbf{q}_{xx}w_y\|^2 + L^{\f12\alpha-} [\phi]_{3,\hat{w}},
\end{align}
and
\begin{align}\label{3.59-2}
  \big\la(b_1 [\cdots] \phi_{y})_y, \, \mathbf{q}_{xx}w^2 \big\ra 
&= - \big\la b_1 [\cdots] \phi_{y} , \, \mathbf{q}_{xx} 2ww_y \big\ra + \big\la b_1 [\cdots] \phi_{xy} , \, \mathbf{q}_{xy}w^2 \big\ra  \nonumber\\
&\quad + \big\la (b_1 [\cdots])_x \phi_{y} , \, \mathbf{q}_{xxy}w^2 \big\ra - \big\la b_1 [\cdots] \phi_{y} , \, \mathbf{q}_{xy}w^2 \big\ra_{x=L}\nonumber\\
&\lesssim L^{\f12\alpha-} [\Phi]_{3,\hat{w}} + L \|\mathbf{q}_{xx}w_y\|^2 + L^{\f12\alpha-} [\phi]_{3,\hat{w}}.
\end{align}

Also 
\begin{align*}
	&\v\, \Big\la \big(b_1[\cdots] \phi_{xy}\big)_y + \big(b_1[\cdots]\phi_{xx}\big)_y + \big(b_1[\cdots] \phi_x\big)_y, \,\, \mathbf{q}_{xx}w^2 \Big\ra\nonumber\\
	&\lesssim \|\sqrt{\v}\mathbf{q}_{xx}w\| \Big\{\sqrt{\v}\|\phi_{xxy}w\| + \sqrt{\v}\|\phi_{xyy}w\| + \sqrt{\v}\|\phi_{xy}\| + \|(\sqrt{\v}\phi_{xx},\sqrt{\v}\phi_x)w\|\Big\} \nonumber\\
	&\lesssim   L^{\f12\alpha-} [\Phi]_{3,\hat{w}}  + L^{\f12\alpha-} [\phi]_{3,\hat{w}},
\end{align*}
which, together with \eqref{7.12-1},  \eqref{3.58}-\eqref{3.59-2}, yields that 
\begin{align}\label{3.61}
\big\la(b_1 I_{u}(\phi))_y, \, \mathbf{q}_{xx}w^2 \big\ra
&\lesssim L^{\f12\alpha-} [\Phi]_{3,\hat{w}} + L \|\mathbf{q}_{xx}w_y\|^2 + L^{\f12\alpha-} [\phi]_{3,\hat{w}}.
\end{align}

Noting $\eqref{3.67-0}_2$, we have
\begin{align}\label{3.62}
\big\la \v \big( b_2 I_{v}(\phi)\big)_x,\, \mathbf{q}_{xx}w^2\big\ra 
&\lesssim \|\sqrt{\v}\mathbf{q}_{xx}w\|\Big\{ \sqrt{\v} \|(\phi_y,\phi_x)w\| + \sqrt{\v}\|(\phi_{xx},\phi_{xy},\phi_{yy})w\| \nonumber\\
&\qquad\qquad\qquad  + \sqrt{\v}\|(\v \phi_{xxx}, \phi_{xxy}, \phi_{xyy})w\| \Big\} \nonumber\\
&\lesssim L^{\f12\alpha-} [\Phi]_{3,\hat{w}} + L^{\f12\alpha-} [\phi]_{3,\hat{w}}.
\end{align}
 
\smallskip

Noting \eqref{3.64-1}, then integrating by parts in $x$, one gets that 
\begin{align}\label{3.62-2}
	&\Big\la \Big(b_1\Big[\hat{d}_{12}\big(\f{T_{sy}}{p_s} q\big)_y + \hat{d}_{13}\f{T_{sy}}{p_s} q\Big]\Big)_y ,\,  \mathbf{q}_{xx}w^2\Big\ra\nonumber\\
	&=\Big\la \Big(b_1\Big[\hat{d}_{12}\big(\f{T_{sy}}{p_s} q\big)_y + \hat{d}_{13}\f{T_{sy}}{p_s} q\Big]\Big)_x,\,  \mathbf{q}_{xy}w^2\Big\ra - \Big\la b_1\Big[\hat{d}_{12}\big(\f{T_{sy}}{p_s} q\big)_y + \hat{d}_{13}\f{T_{sy}}{p_s} q\Big],\,  \mathbf{q}_{xy}w^2\Big\ra_{x=L} \nonumber\\
	&\quad -\Big\la b_1\Big[\hat{d}_{12}\big(\f{T_{sy}}{p_s} q\big)_y + \hat{d}_{13}\f{T_{sy}}{p_s} q\Big]  ,\,  \mathbf{q}_{xx}2ww_y\Big\ra \nonumber\\
	& \lesssim L^{\f12\alpha-} [\Phi]_{3,\hat{w}} + L \|\mathbf{q}_{xx}w_y\|^2 + L^{\f12\alpha-} [\phi]_{3,\hat{w}}.
\end{align}
It is direct to have 
\begin{align}\label{3.62-3}
	&\Big\la  \Big(b_2\Big[\v \hat{d}_{21}\, \big(\f{T_{sy}}{p_s}q\big)_x  +  \v \hat{d}_{22}\, (\f{T_{sy}}{p_s} q)_y +  \v \hat{d}_{23}\, \f{T_{sy}}{p_s} q\Big]\Big)_x ,\,  \mathbf{q}_{xx}w^2\Big\ra  
 \lesssim L^{\f12\alpha-} [\Phi]_{3,\hat{w}} + L^{\f12\alpha-} [\phi]_{3,\hat{w}}.
\end{align}

\smallskip

Using   Lemma \ref{lemA.3},  we have 
\begin{align}
\big\la (b_1[\tilde{d}_{12}+\cdots] p_y)_y,\, \mathbf{q}_{xx}w^2 \big\ra
&=\big\la b_1 [\tilde{d}_{12}+\cdots]\,p_{yy},\, \mathbf{q}_{xx}w^2 \big\ra + \big\la  (b_1 [\tilde{d}_{12}+\cdots])_y\,p_y,\, \mathbf{q}_{xx}w^2 \big\ra\nonumber\\
& \lesssim L^{\f12\alpha-} [\Phi]_{3,\hat{w}} + L^{\f12\alpha-} \frac{1}{\v} [[p]]_{2,w},\label{3.68}\\
\big\la \v (b_2 [\tilde{d}_{21}+\cdots]\, p_x)_x,\, \mathbf{q}_{xx}w^2 \big\ra & =\big\la \v  b_2[\tilde{d}_{21}+\cdots]\, p_{xx},\, \mathbf{q}_{xx}w^2 \big\ra + \big\la \v (b_2[\tilde{d}_{21}+\cdots])_x \, p_x,\, \mathbf{q}_{xx}w^2 \big\ra \nonumber\\
&\lesssim L^{\f12\alpha-} [\Phi]_{3,\hat{w}} + L^{\f12\alpha-} \frac{1}{\v} [[p]]_{2,w},\label{3.63},\\
\big\la \v (b_2 \hat{d}_{22}\, p_y)_x,\, \mathbf{q}_{xx}w^2 \big\ra&\lesssim  \sqrt{\v} [\Phi]_{3,\hat{w}} + \frac{1}{\sqrt{\v}} [[p]]_{2,w},\label{3.63-3}
\end{align}
and
\begin{align}\label{3.64}
	\big\la \v (b_2 [\tilde{d}_{23}+\cdots] p)_x,\, \mathbf{q}_{xx}w^2 \big\ra
	&\lesssim  L^{\f12\alpha-} [\Phi]_{3,w}^{\f12} \Big(\sqrt{\v}\|p_xw\| + \|\mathbf{d}_2  pw\|\Big) \nonumber\\
	&\lesssim L^{\f14\alpha} [\Phi]_{3,\hat{w}} + L^2 \f{1}{\v}[[p]]_{2,w} + L \f{1}{\v^2} [[p]]_{2,1}\cdot \mathbf{1}_{\{w=w_0\}} .
\end{align}

Noting \eqref{3.65-1}, then using \eqref{3.72-2} and Lemma \ref{lemA.3}, we obtain
\begin{align}\label{3.72-1}
\big\la (b_1 [\tilde{d}_{13} + \cdots]\, p)_y,\, \mathbf{q}_{xx}w_0^2 \big\ra  
&\lesssim \|\sqrt{\v}\mathbf{q}_{xx}w_0\|\bigg\{\|\frac{1}{\sqrt{\v}} p_yw_0\| + \frac{1}{\sqrt{\v}} \|\mathbf{d}_3pw_0\| \bigg\}\nonumber\\
&\lesssim L^{\f12\alpha-} [\Phi]_{3,\hat{w}_0} +  L^{\f14\alpha-} \f1{\v}[[p]]_{2,w_0}  + L \f1{\v^2}[[p]]_{2,1} \cdot \mathbf{1}_{\{w=w_0\}},
\end{align}
where we have used \eqref{3.72-2}.

The method used in \eqref{3.72-1} does not apply to the case 
 $w=1$, as it requires an additional factor of $\sqrt{\v}$.  A direct calculation shows that 
\begin{align}\label{3.69}
	\big\la (b_1  [\tilde{d}_{13} + \cdots] p)_y,\, \mathbf{q}_{xx} \big\ra 
	&= \big\la (b_1 [\tilde{d}_{13} + \cdots])_y\, p,\, \mathbf{q}_{xx}  \big\ra +  L^{\f12\alpha-} [\Phi]_{3,1} + L^{\f12\alpha-} \frac{1}{\v}[[p]]_{2,1}.
\end{align}
The first term on the right-hand side of \eqref{3.69} is difficult to control directly and requires careful treatment. Applying integration by parts in $x$, we obtain
\begin{align}\label{3.69-1}
	\big\la (b_1  [\tilde{d}_{13} + \cdots])_y\,p,\, \mathbf{q}_{xx}\big\ra
	&= \big\la (b_1  [\tilde{d}_{13} + \cdots])_y\,p,\, \mathbf{q}_{x} \big\ra \Big|_{x=0}^{x=L}  - \big\la (b_1  [\tilde{d}_{13} + \cdots])_y\,p_x,\, \mathbf{q}_{x} \big\ra \nonumber\\
	&\quad  - \big\la (b_1 [\tilde{d}_{13} + \cdots])_{xy}\, p,\, \mathbf{q}_{x} \big\ra.
\end{align}
For the last term on RHS of \eqref{3.69-1},   one has from Lemmas \ref{lemA.3} $\&$ \ref{lemA.5-1} that  
\begin{align}\label{3.69-3}
	&|\big\la (b_1  [\tilde{d}_{13} + \cdots])_{xy}\, p,\, \mathbf{q}_{x}\big\ra|\nonumber\\
	& \leq 
	\big|\big\la (b_1 [\tilde{d}_{13} + \cdots])_{xy}\, [p(x,y)-p(x,0)],\, \mathbf{q}_{x}\big\ra\big| + \big|\big\la (b_1  [\tilde{d}_{13} + \cdots])_{xy}\, p(x,0),\, \mathbf{q}_{x}\big\ra\big|\nonumber\\
	&\lesssim a^{2}\iint \big\{\la y\ra^{-\f12\mathfrak{l}_0} + \v \la \sqrt{\v}y\ra^{-\f12\mathfrak{l}_0}\big\} \big\{\sqrt{y} |\mathbf{q}_x|  \cdot \|p_y(x,\cdot)\|_{L^2_y} + |\mathbf{q}_x|\cdot |p(x,0)| \big\}dydx 
	\nonumber\\
	&\lesssim \|\frac{\mathbf{q}_x}{y}\|\cdot \|p_y\|   +   \|\frac{\mathbf{q}_x}{y}\|\cdot  \|p\|_{y=0}  
	\lesssim L [\Phi]_{3,1} + L \f1{\v}[[p]]_{2,1}.
\end{align}
For the second term on RHS of \eqref{3.69-1}, it follows from Lemma \ref{lemA.3}  that 
\begin{align}\label{3.69-5}
|\big\la (b_1 [\tilde{d}_{13} + \cdots])_y\, p_x,\, \mathbf{q}_{x} \big\ra| 
&\lesssim a^{2} \iint \big\{\la y\ra^{-\f12\mathfrak{l}_0} + \v \la \sqrt{\v}y\ra^{-\f12\mathfrak{l}_0}\big\} |p_x \mathbf{q}_x| dydx \lesssim \|\frac{\mathbf{q}_x}{y}\| \|p_x\| \nonumber\\
& \lesssim L^{\f12\alpha-} [\Phi]_{3,1} + L^{\f14\alpha-}\frac{1}{\v} [[p]]_{2,1}.
\end{align}
For the boundary term of \eqref{3.69-1}  at $x=L$, it follows from \eqref{3.65-1} and  Lemmas \ref{lemA.4} \& \ref{lemA.5-1}  that
\begin{align}\label{3.69-6}
	\big\la (b_1 [\tilde{d}_{13} + \cdots])_y\, p,\, \mathbf{q}_{x} \big\ra_{x=L} &\lesssim \|\mathbf{q}_x\|_{x=L}\cdot \|\mathbf{d}_3 p\|_{x=L}  \lesssim L^{1-} [\Phi]_{3,1} + L\f1{\v} [[p]]_{2,1},
\end{align}
and
\begin{align}\label{3.69-7}
	\big\la (b_1 [\tilde{d}_{13} + \cdots])_y\, p,\, \mathbf{q}_{x} \big\ra_{x=0} 
	&\lesssim a^{2} \int_0^\infty \big\{\la y\ra^{-\f12\mathfrak{l}_0} + \v \la \sqrt{\v}y\ra^{-\f12\mathfrak{l}_0}\big\} |\mathbf{q}_x p| dy \Big|_{x=0}   \nonumber\\
	&\lesssim a^{2} \|p_y\|_{x=0}  \int_0^\infty \sqrt{y} \big\{\la y\ra^{-\f12\mathfrak{l}_0} + \v \la \sqrt{\v}y\ra^{-\f12\mathfrak{l}_0}\big\} |\mathbf{q}_x(0,y)| dy \nonumber\\
	&\lesssim \|p_y\|_{x=0} \cdot \|\frac{\mathbf{q}_x}{\la y\ra}\|_{x=0} + \|p_y\|_{x=0} \cdot \|\sqrt{\v}\mathbf{q}_x\|_{x=0} \nonumber\\
	&\lesssim    L^{\f12\alpha-} [\Phi]_{3,1} +  L^{\f12\alpha-} \f1{\v} [[p]]_{2,1}.
\end{align}
Substituting \eqref{3.69-3}-\eqref{3.69-7} into \eqref{3.69-1}, one has that 
\begin{align}\label{3.69-9}
	|\big\la (b_1 [\tilde{d}_{13} + \cdots])_y\, p,\, \mathbf{q}_{xx}\big\ra| 
	&\lesssim L^{\f12\alpha-} [\Phi]_{3,1} +  L^{\f12\alpha-} \f1{\v} [[p]]_{2,1}.
\end{align}
Hence, combining \eqref{3.69-9} and \eqref{3.69}, one has
\begin{align}\label{3.69-10}
	\big\la (b_1 [\tilde{d}_{13} + \cdots]\, p)_y,\, \mathbf{q}_{xx} \big\ra   
	&\lesssim  L^{\f12\alpha-} [\Phi]_{3,1} +  L^{\f12\alpha-} \f1{\v} [[p]]_{2,1}.
\end{align}

\smallskip

Using \eqref{6.174}, we get 
\begin{align}\label{3.104-1}
	&\big\la \v (b_2\hat{d}_{21} S_x)_x +\v (b_2\hat{d}_{22} S_y)_x + \v (b_2\hat{d}_{23} S)_x,\, \mathbf{q}_{xx}w^2 \big\ra \nonumber\\
	&\lesssim L^{\f12\alpha-} [\Phi]_{3,\hat{w}} + L^{\f12\alpha-}  [[[S]]]_{2,\hat{w}} + L\f{1}{\v} [[[S]]]_{2,1}  \mathbf{1}_{\{w=w_0\}}.
\end{align}
 Integrating by parts, one can obtain
\begin{align}\label{3.112R.5}
&-\big\la (b_1\hat{d}_{12} S_{y})_y ,\, \mathbf{q}_{xx}w^2 \big\ra \nonumber\\
&= - \big\la b_1\hat{d}_{12} S_{xy}  ,\, \mathbf{q}_{xy}w^2 \big\ra  - \big\la (b_1\hat{d}_{12})_x S_{y}  ,\, \mathbf{q}_{xy}w^2 \big\ra   + \big\la b_1\hat{d}_{12} S_{y}  ,\, \mathbf{q}_{xx}2ww_y \big\ra\nonumber\\
&\quad  +   \big\la b_1\hat{d}_{12} S_{y}  ,\, \mathbf{q}_{xy}w^2 \big\ra_{x=L} -  \big\la b_1\hat{d}_{12} S_{y}  ,\, \mathbf{q}_{xy}w^2 \big\ra_{x=0}  \nonumber\\
&\lesssim L  \|\mathbf{q}_{xx}w_y\|^2 + L^{\f12\alpha-} [\Phi]_{3,\hat{w}}  + L^{\f12\alpha-} [[[S]]]_{2,\hat{w}} + L\f{1}{\v} [[p]]_{2,1},
\end{align}
where we have used $ S|_{y=0} =0$ for case of DT, and the following estimation for the case of NT
\begin{align}\label{3.112R.5-1}
&\big\la b_1\hat{d}_{12} S_{y}  ,\, \mathbf{q}_{xy}w^2 \big\ra_{x=0} \nonumber\\
&=\Big\la b_1\hat{d}_{12} \big(\f{\chi}{2\rho_s}p\big)_{y}  ,\, \mathbf{q}_{xxy}w^2 \Big\ra + \Big\la \big(b_1\hat{d}_{12} \big(\f{\chi}{2\rho_s}p\big)_{y}\big)_x,\, \mathbf{q}_{xy}w^2 \Big\ra   - \Big\la b_1\hat{d}_{12} \big(\f{\chi}{2\rho_s}p\big)_{y},\, \mathbf{q}_{xy}w^2 \Big\ra_{x=L} \nonumber\\
&=-\Big\la \big(b_1\hat{d}_{12} \big(\f{\chi}{2\rho_s}p\big)_{y}\big)_y  ,\, \mathbf{q}_{xx}w^2 \Big\ra - \Big\la b_1\hat{d}_{12} \big(\f{\chi}{2\rho_s}p\big)_{y}  ,\, \mathbf{q}_{xx}2ww_y \Big\ra  + L^{\f12\alpha-} [\Phi]_{3,1} + L^{\f12\alpha-}\f{1}{\v} [[p]]_{2,1} \nonumber\\
&=-\Big\la \big(b_1\hat{d}_{12} \big(\f{\chi}{2\rho_s}\big)_{y}\big)_y p  ,\, \mathbf{q}_{xx}w^2 \Big\ra - \Big\la b_1\hat{d}_{12} \big(\f{\chi}{2\rho_s}\big)_{y}  p ,\, \mathbf{q}_{xx}2ww_y \Big\ra  + L^{\f12\alpha-} \f{1}{\v} [[p]]_{2,1} + L^{\f12\alpha-} [\Phi]_{3,1} \nonumber\\
&\lesssim L^{\f12\alpha-} \f{1}{\v} [[p]]_{2,1} + L^{\f12\alpha-} [\Phi]_{3,1}.
\end{align}
Here we remark that similar estimate as in \eqref{3.69-9} are applied in the last inequality of \eqref{3.112R.5-1}.

\smallskip

It is direct to have 
\begin{align}\label{3.106-1} 
&-\big\la (b_1\hat{d}_{13} S)_y,\, \mathbf{q}_{xx}w^2 \big\ra\nonumber\\
&= \big\la b_1\hat{d}_{13} S,\, \mathbf{q}_{xx} 2ww_y \big\ra - \big\la b_1\hat{d}_{13} S_x,\, \mathbf{q}_{xy} w^2 \big\ra - \big\la (b_1\hat{d}_{13})_x S,\, \mathbf{q}_{xy} w^2 \big\ra \nonumber\\
&\quad + \big\la b_1\hat{d}_{13} S,\, \mathbf{q}_{xy}w^2 \big\ra_{x=L} + \big\la b_1\hat{d}_{13}S,\, \mathbf{q}_{xy}w^2 \big\ra_{x=0}\nonumber\\
&\lesssim \sqrt{L}\|\mathbf{q}_{xx}w_y\|^2 + L^{\f12\alpha-} [\Phi]_{3,\hat{w}} + L^{\f12\alpha-}  [[[S]]]_{2,\hat{w}} + \sqrt{L}\f{1}{\v} [[[S]]]_{2,1} \mathbf{1}_{\{w=w_0\}},
\end{align}
where we have used $S|_{x=0}=0$ for the case of DT, and the following estimation for the case of NT
\begin{align}\label{3.103-10}
&\big\la b_1\hat{d}_{13} S,\, \mathbf{q}_{xy}w^2 \big\ra_{x=0} = \Big\la \hat{d}_{13}\f{b_1\chi}{2\rho_s} p,\, \mathbf{q}_{xy}w^2 \Big\ra_{x=0} \nonumber\\
&=\Big\la \hat{d}_{13}\f{b_1\chi}{2\rho_s} p,\, \mathbf{q}_{xxy}w^2 \Big\ra   + \Big\la \big(\hat{d}_{13}\f{b_1\chi}{2\rho_s} p\big)_x ,\, \mathbf{q}_{xy}w^2 \Big\ra - \Big\la \hat{d}_{13}\f{b_1\chi}{2\rho_s} p,\, \mathbf{q}_{xy}w^2 \Big\ra_{x=L} \nonumber\\
&=-\Big\la \hat{d}_{13}\f{b_1\chi}{2\rho_s} p_y,\, \mathbf{q}_{xx}w^2 \Big\ra - \Big\la \big(\hat{d}_{13}\f{b_1\chi}{2\rho_s}w^2\big)_y p,\, \mathbf{q}_{xx} \Big\ra  +  L^{\f12\alpha-} [\Phi]_{3,1} + L^{\f12\alpha-} \f{1}{\v} [[p]]_{2,1}\nonumber\\
&\lesssim  L^{\f12\alpha-} [\Phi]_{3,1} + L^{\f12\alpha-} \f{1}{\v} [[p]]_{2,1}.
\end{align}
Also we have used similar estimate as in \eqref{3.69-9}  to control $\dis \Big\la \big(\hat{d}_{13}\f{b_1\chi}{2\rho_s}w^2\big)_y p,\, \mathbf{q}_{xx} \Big\ra$ in \eqref{3.103-10}.

Noting $\eqref{2.98}_3$ and \eqref{3.66-0}-\eqref{3.66-1}, one has from \eqref{3.61}-\eqref{3.62-3}, \eqref{3.68}-\eqref{3.72-1} and \eqref{3.69-10}-\eqref{3.106-1} that
\begin{align}\label{3.91-1}
\big\la \mathcal{F}_{R3}, \, \mathbf{q}_{xx}w^2 \big\ra  
&\lesssim L^{\f12} \|\mathbf{q}_{xx}w_y\|^2 + L^{\f12\alpha} [\Phi]_{3,\hat{w}}  + L^{\f12\alpha-} \f{1}{\v} [[p]]_{2,\hat{w}} + L^{\f12\alpha-}  [[[S]]]_{2,\hat{w}}\nonumber\\
&\quad  + L^{\f12} \Big\{\f{1}{\v} [[[S]]]_{2,1} + \f1{\v^2}[[p]]_{2,1} \Big\} \mathbf{1}_{\{w=w_0\}}  +  L^{\f12\alpha-} [\phi]_{3,\hat{w}} + \v^{2N_0-1}.
\end{align}

\medskip

{\it Step 4.4. Estimate on $\big\la \mathcal{F}_{R4},\, \mathbf{q}_{xx}w^2\big\ra$.} Recall $\eqref{2.98}_4$, \eqref{2.83} and $\tilde{T}=\zeta + \f{1}{2\rho_s} p + \f{1}{p_s}T_{sy} q$, we have 
\begin{align}\label{3.92-1}
\mathcal{F}_{R4}
&=(\mu+\lambda)\Big(\frac{b_1 }{T^{\v}} u^{\v} \Delta_{\v}\tilde{T}\Big)_y  - \mu \v \Big(\frac{b_2}{T^{\v}} v^{\v} \Delta_{\v}\tilde{T}\Big)_x  + \lambda\v\frac{b_1-b_2}{T^{\v}} v^{\v} \tilde{T}_{xyy} \nonumber\\
&\quad  + \lambda\v \big(\frac{b_1}{T^{\v}} v^{\v}\big)_y \tilde{T}_{xy} - \lambda\v\big(\frac{b_2}{T^{\v}}v^{\v}\big)_x \tilde{T}_{yy}.
\end{align}
  From \eqref{84}, we note  $b_1-b_2=u_s^2 + O(1)\v u_s$.

Integrating by parts in $y$, then using \eqref{2.122-0}, \eqref{3.92}-\eqref{3.94}, \eqref{6.176}  one has that  
\begin{align}
&\v \Big\la  \frac{b_1-b_2}{T^{\v}}  v^{\v} \tilde{T}_{xyy} ,\, \mathbf{q}_{xx}w^2\Big\ra\nonumber\\
&=- \v \Big\la  \frac{b_1-b_2}{T^{\v}}  v^{\v} \tilde{T}_{xy} ,\, \mathbf{q}_{xxy}w^2\Big\ra - \v \Big\la  \frac{b_1-b_2}{T^{\v}} v^{\v} \tilde{T}_{xy} ,\, \mathbf{q}_{xx} 2ww_y\Big\ra \nonumber\\
&\quad - \v \Big\la  \big(\frac{b_1-b_2}{T^{\v}}\big)_y  v^{\v} \tilde{T}_{xy} ,\, \mathbf{q}_{xx}w^2\Big\ra  - \v \Big\la  \frac{b_1-b_2}{T^{\v}} v^{\v}_{y} \tilde{T}_{xy} ,\, \mathbf{q}_{xx}w^2\Big\ra\nonumber\\
&\lesssim \f1N [\Phi]_{3,\hat{w}} + Na^4 \|(\sqrt{\v}\tilde{T}_{xx}, \sqrt{\v}\tilde{T}_{xy})w\|^2 + \v^{N_0} \|\sqrt{\v}\tilde{T}_{xy}w\|^2\cdot \|(u_s T_y, v_y)\|_{L^\infty} \nonumber\\
&\lesssim \f1N [\Phi]_{3,\hat{w}}  + Na^4 [[[\zeta]]]_{2,w} + Na^4\frac{1}{\v} [[p]]_{2,\hat{w}}  +  L^{\f12\alpha-} [\phi]_{3,1},
\end{align}
and
\begin{align}
&\Big\la \Big(\frac{b_1 }{T^{\v}} u^{\v} \Delta_{\v}\tilde{T}\Big)_y,\, \mathbf{q}_{xx}w^2\Big\ra\nonumber\\
&= - \Big\la  \frac{b_1u_s }{T_s}  \Delta_{\v}\tilde{T} ,\, \mathbf{q}_{xxy}w^2\Big\ra - \Big\la \frac{b_1 u_s}{T_s}  \Delta_{\v}\tilde{T},\, \mathbf{q}_{xx}2ww_y\Big\ra   + C\v [\Phi]_{3,\hat{w}} + C\v \|\Delta_{\v}\tilde{T} w\|^2 \nonumber\\
&= - \Big\la  \frac{b_1u_s }{T_s} \Delta_{\v}\zeta,\, \mathbf{q}_{xxy}w^2\Big\ra - \Big\la  \frac{b_1u_s }{T_s} \Delta_{\v}\big(\f{p}{2\rho_s}\big) ,\, \mathbf{q}_{xxy}w^2\Big\ra - \Big\la \frac{b_1u_s }{T_s} \Delta_{\v} \big(\f{T_{sy}}{p_s}q\big) ,\, \mathbf{q}_{xxy}w^2\Big\ra  \nonumber\\
&\quad  + Ca\|\mathbf{q}_{xx}w_y\|\cdot \|\Delta_{\v}\tilde{T}w\|  + C\v [\Phi]_{3,\hat{w}} + C \v \|\Delta_{\v}\tilde{T} w\|^2 \nonumber\\
&\lesssim \f{1}{N} [\Phi]_{3,\hat{w}} + \f{1}{N} \|\mathbf{q}_{xx}w_y\|^2 + Na [[[\zeta]]]_{2,w}  + L^{\f14-} [[[\zeta]]]_{b,w} + Na \frac{1}{\v} [[p]]_{2,\hat{w}} \nonumber\\
&\quad   + L^{\f12-} [[p]]_{3,1}  +  L^{\f12\alpha-} [\phi]_{3,\hat{w}},
\end{align}
where we have used similar arguments as in \eqref{3.52-1}-\eqref{3.72-3}. 

Similarly one can obtain
\begin{align}
&\Big\la \v \Big(\frac{b_2}{T^{\v}} v^{\v} \Delta_{\v}\tilde{T}\Big)_x,\, \mathbf{q}_{xx}w^2\Big\ra
\lesssim  L^{\f12\alpha-} [\Phi]^{\f12}_{3,\hat{w}} \Big\{\|u_s\sqrt{\v}\Delta_{\v}\tilde{T}_{x}w\| + \|u_s \Delta_{\v}\tilde{T}w\|\Big\} \nonumber\\
&\lesssim L^{\f12\alpha-} [\Phi]_{3,\hat{w}}  +  L^{\f12}[[[\zeta]]]_{2,w}  + L\f{1}{\v}[[p]]_{2,\hat{w}}  + L^{1-} [\phi]_{3,1}   +  L^{\f12\alpha-} \|u_s\sqrt{\v}\Delta_{\v}\tilde{T}_{x}w\|^2,
\end{align}
and
\begin{align}\label{3.92-10}
&\Big\la  \lambda\v \big(\frac{b_1}{T^{\v}} v^{\v}\big)_y \tilde{T}_{xy} - \lambda\v\big(\frac{b_2}{T^{\v}}v^{\v}\big)_x \tilde{T}_{yy},\, \mathbf{q}_{xx}w^2\Big\ra  \lesssim \sqrt{\v} [\Phi]_{3,\hat{w}} + \sqrt{\v}\|(\tilde{T}_{xy}, \tilde{T}_{yy})w\|^2\nonumber\\
&\lesssim  \sqrt{\v} [\Phi]_{3,\hat{w}}  + \sqrt{\v}[[[\zeta]]]_{2,\hat{w}}  +  \f{1}{\sqrt{\v}}[[p]]_{2,1}  + \sqrt{\v} [\phi]_{3,1}.
\end{align}
Hence, combining \eqref{3.92-1}-\eqref{3.92-10}, one obtains 
\begin{align}\label{3.110-1}
	\big\la \mathcal{F}_{R4},\, q_{xx}w^2\big\ra
	&\lesssim \f{1}{N} [\Phi]_{3,\hat{w}} + \f{1}{N} \|\mathbf{q}_{xx}w_y\|^2 + Na [[[\zeta]]]_{2,w}  +L^{\f14-} [[[\zeta]]]_{b,w} + Na \frac{1}{\v} [[p]]_{2,\hat{w}} \nonumber\\
	&\quad   + L^{\f12-} [[p]]_{3,1}  +  L^{\f12\alpha-} [\phi]_{3,\hat{w}}  +  L^{\f12\alpha-} \|u_s\sqrt{\v}\Delta_{\v}\tilde{T}_{x}w\|^2.
\end{align}

\medskip

{\it Step 4.5. Estimate on $\big\la \mathcal{F}_{R5},\, \mathbf{q}_{xx}w^2\big\ra$.}  Noting the estimates in Section \ref{Sec2.4}, by similar arguments as in \eqref{3.72-1} and \eqref{3.69-10}, one has 
\begin{align}\label{3.111-1}
\Big\la \Big(\f{b_1}{T^{\v}}u^{\v} (\f{1}{2\rho_s})_{yy}  p\Big)_y,\, \mathbf{q}_{xx}w^2 \Big\ra
&\lesssim  L^{\f12\alpha-} [\Phi]_{3,\hat{w}} +  L^{\f12\alpha-} \f1{\v}[[p]]_{2,\hat{w}}  + L \f1{\v^2}[[p]]_{2,1}  \mathbf{1}_{\{w=w_0\}}.
\end{align}
Also it holds that 
\begin{align}
\Big\la \Big(\f{b_1}{T^{\v}} u^{\v} (\f{1}{\rho_s})_{y}p_y\Big)_y,\, \mathbf{q}_{xx}w^2 \Big\ra \lesssim  L^{\f12\alpha-} [\Phi]_{3,\hat{w}} +  L^{\f12\alpha-} \f1{\v}[[p]]_{2,\hat{w}},
\end{align}
and
\begin{align}
&\Big\la \Big(\f{b_1}{T^{\v}}u^{\v} (\f{T_{sy}}{p_s}q)_{yy}\Big)_y,\, \mathbf{q}_{xx}w^2 \Big\ra\nonumber\\
&=-\Big\la \f{b_1u_s}{T_s} (\f{T_{sy}}{p_s}q)_{xyy} ,\, \mathbf{q}_{xy}w^2 \Big\ra - \Big\la (\f{b_1u_s}{T_s})_x (\f{T_{sy}}{p_s}q)_{yy} ,\, \mathbf{q}_{xy}w^2 \Big\ra + \v[\Phi]_{3,\hat{w}} +  \v [\phi]_{3,\hat{w}}\nonumber\\
&\quad + \Big\la \f{b_1u_s}{T_s} (\f{T_{sy}}{p_s}q)_{yy} ,\, \mathbf{q}_{xx} 2ww_y \Big\ra + \Big\la \f{b_1u_s}{T_s} (\f{T_{sy}}{p_s}q)_{yy} ,\, \mathbf{q}_{xy}w^2 \Big\ra_{x=L} \nonumber\\
&\lesssim L^{\f12\alpha-} [\Phi]_{3,\hat{w}} + L^{\f12\alpha-} [\phi]_{3,\hat{w}} + L^{\f18} \|\mathbf{q}_{xx}w_y\|^2 
.
\end{align}

A direct calculation shows that 
\begin{align}\label{3.114}
\big\la \v \big(b_2G_{22}\big)_x, \mathbf{q}_{xx}w^2\big\ra
&\lesssim L^{\f12\alpha} [\Phi]_{3,\hat{w}} + L^{\f12\alpha-} [\phi]_{3,\hat{w}} + [[p]]_{2,\hat{w}}.
\end{align}

Recall $\eqref{2.98}_5$ and \eqref{2.96-1}, then we have from  \eqref{3.111-1}-\eqref{3.114}, one obtains 
\begin{align}\label{3.115}
\big\la \mathcal{F}_{R5},\, \mathbf{q}_{xx}w^2\big\ra
&\lesssim  L^{\f18} \|\mathbf{q}_{xx}w_y\|^2  + L^{\f12\alpha-} [\Phi]_{3,\hat{w}}  + L^{\f12\alpha-} [\phi]_{3,\hat{w}} \nonumber\\
&\quad  +  L^{\f12\alpha-} \f1{\v}[[p]]_{2,\hat{w}}  + L \f1{\v^2}[[p]]_{2,1}  \mathbf{1}_{\{w=w_0\}}.
\end{align}

\smallskip

{\it Step 4.6. Estimate on $\big\la \mathcal{F}_{R6}(\hat{\mathbf{S}}),\, \mathbf{q}_{xx}w^2\big\ra$.} Integrating by parts, one has
\begin{align}\label{3.116}
&\lambda\Big\la \big(\f{b_1u^{\v}}{T^{\v}} \hat{\mathbf{S}}_{yy}\big)_y ,\, \mathbf{q}_{xx}w^2\Big\ra 
\nonumber\\
&=-\lambda\Big\la \f{b_1u_s}{T_s} \hat{\mathbf{S}}_{xy} ,\, \mathbf{q}_{xyy}w^2\Big\ra - \lambda\Big\la \big(\f{b_1u_s}{T_s}w^2\big)_y \hat{\mathbf{S}}_{xy} ,\, \mathbf{q}_{xy}w^2\Big\ra  +\lambda \Big\la \big(\f{b_1u_s}{T_s}\big)_x \hat{\mathbf{S}}_{yy} ,\, \mathbf{q}_{xy}w^2\Big\ra  \nonumber\\
&\quad- \lambda\Big\la \f{b_1u_s }{T_s}\hat{\mathbf{S}}_{yy} ,\, \mathbf{q}_{xx} 2ww_y\Big\ra - \lambda\Big\la \f{b_1u_s}{T_s} \hat{\mathbf{S}}_{yy} ,\, \mathbf{q}_{xy}w^2\Big\ra_{x=L}  + \lambda\Big\la \f{b_1}{T_s}u_s \hat{\mathbf{S}}_{yy} ,\, \mathbf{q}_{xy}w^2\Big\ra_{x=0} \nonumber\\ 
&\quad + C\v  [\Phi]_{3,\hat{w}}+ C\v [[[\hat{\mathbf{S}}]]]_{2,w} \nonumber\\
&\lesssim  \f{1}{N} [\Phi]_{3,\hat{w}} + L^{\f14} \|\mathbf{q}_{xx}w_y\|^2  + Na \f{1}{\v} \mathcal{K}_{2,1}(\hat{P}) + L^{\f12-}  \mathcal{K}_{3,1}(\hat{P}) + \v^{\f{1-\beta}{2}} \mathscr{K}(\hat{P}) + N a [[[\hat{\mathbf{S}}]]]_{2,w} \nonumber\\
&\quad +L^{\f12-} \|\sqrt{u_s}\hat{\mathbf{S}}_{yy}w\|^2_{x=L},
\end{align}
where we have used  $\hat{\mathbf{S}}|_{x=0}=0$ for the case of DT, and used \eqref{2.125} to derive the following estimation for the case of NT
\begin{align*} 
	&-\Big\la \f{b_1u_s}{T_s}  \hat{\mathbf{S}}_{yy} ,\, \mathbf{q}_{xy}w^2  \Big\ra_{x=0}
	=\Big\la \f{b_1u_s}{T_s}  \big(\f{\chi}{2\rho_s}\hat{P}\big)_{yy} ,\, \mathbf{q}_{xy}w^2  \Big\ra_{x=0} \nonumber\\
	&=\Big\la \f{b_1u_s}{2p_s} \hat{P}_{yy} \chi,\, \mathbf{q}_{xy}w^2  \Big\ra_{x=0} + \Big\la \f{b_1u_s}{T_s}\big(\f{\chi}{\rho_s}\big)_y \hat{P}_{y} ,\, \mathbf{q}_{xy}w^2  \Big\ra_{x=0} + \Big\la \f{b_1u_s}{T_s}\big(\f{\chi}{2\rho_s}\big)_{yy} \hat{P} ,\, \mathbf{q}_{xy}w^2  \Big\ra_{x=0}\nonumber\\
	&=-\Big\la \f{b_1u_s}{2p_s}\hat{P}_{xyy} \chi,\, \mathbf{q}_{xy}w^2  \Big\ra - \Big\la \f{b_1u_s}{2p_s} \hat{P}_{yy} \chi,\, \mathbf{q}_{xxy}w^2  \Big\ra - \Big\la \big(\f{b_1u_s}{2p_s}\big)_x \hat{P}_{yy} \chi,\, \mathbf{q}_{xy}w^2  \Big\ra \nonumber\\
	&\quad + \Big\la \f{b_1u_s}{2p_s}\hat{P}_{yy} \chi,\, \mathbf{q}_{xy}w^2  \Big\ra_{x=L} + C \|u_s\mathbf{q}_{xy}w\chi\|_{x=0}\cdot \|\chi \hat{P}_y\|_{x=0}\nonumber\\
	&=\Big\la \f{b_1u_s}{2p_s}  \hat{P}_{xy} \chi,\, \mathbf{q}_{xyy}w^2  \Big\ra + \Big\la \big(\f{b_1u_s}{2p_s}  \chi w^2\big)_y \hat{P}_{xy},\, \mathbf{q}_{xy}w^2  \Big\ra  + C\sqrt{a} [\Phi]_{3,w}^{\f12} \cdot \|\f{1}{\sqrt{\v}}\hat{P}_{yy}w\| \nonumber\\
	&\quad  + C\sqrt{a}\|\sqrt{u_s}\hat{P}_{yy}w\chi\|_{x=L}\cdot \|\mathbf{q}_{xy}w\|_{x=L} +  C\|u_s\mathbf{q}_{xy}w\|_{x=0}\cdot \|\chi \hat{P}_y\|_{x=0} \nonumber\\
	&\lesssim \f{1}{N} [\Phi]_{3,1}  + Na \f{1}{\v} \mathcal{K}_{2,1}(\hat{P}) + L^{\f12-} \mathcal{K}_{3,1}(\hat{P}) + \v^{\f{1-\beta}{2}} \mathscr{K}(\hat{P}).
\end{align*}
Here we remark that we can only obtain a uniform bound for $\dis \v^{-1+\beta}\|\hat{P}_y\chi\|^2$, which is weaker than the corresponding estimate for $p$. Fortunately, this still suffices to close the estimate in our argument.

Integrating by parts, then using \eqref{3.94},  one has 
\begin{align*}
\lambda\v \Big\la \big(\f{b_1}{T^{\v}}u^{\v} \hat{\mathbf{S}}_{xy}\big)_x ,\, \mathbf{q}_{xx}w^2\Big\ra 
&=-\lambda\v \Big\la \f{b_1}{T^{\v}}u^{\v} \hat{\mathbf{S}}_{xx}  ,\, \mathbf{q}_{xxy}w^2\Big\ra - \lambda\v \Big\la \f{b_1}{T^{\v}}u^{\v} \hat{\mathbf{S}}_{xx}  ,\, \mathbf{q}_{xx} 2ww_y\Big\ra \nonumber\\
&\quad  + \lambda\v \Big\la \big(\f{b_1}{T^{\v}}u^{\v}\big)_x  \hat{\mathbf{S}}_{xy},\, \mathbf{q}_{xx}w^2\Big\ra - \lambda\v \Big\la (\f{b_1}{T^{\v}}u^{\v})_y \hat{\mathbf{S}}_{xx}  ,\, \mathbf{q}_{xx}w^2\Big\ra\nonumber\\
&\lesssim \f{1}{N} [\Phi]_{3,\hat{w}} + Na [[[\hat{\mathbf{S}}]]]_{2,w},
\end{align*}
which, together with $\eqref{2.98-1}_1$ and \eqref{3.116}, yields that 
\begin{align}\label{3.119}
\big\la \mathcal{F}_{R6}(\hat{\mathbf{S}}),\, \mathbf{q}_{xx}w^2\big\ra 
 &\lesssim \f{1}{N} [\Phi]_{3,\hat{w}} + L^{\f14} \|\mathbf{q}_{xx}w_y\|^2  + Na \f{1}{\v} \mathcal{K}_{2,1}(\hat{P}) + L^{\f12-}  \mathcal{K}_{3,1}(\hat{P}) + \v^{\f{1-\beta}{2}} \mathscr{K}(\hat{P}) \nonumber\\
 &\quad  + N a [[[\hat{\mathbf{S}}]]]_{2,w} + L^{\f12-} \|\sqrt{u_s}\hat{\mathbf{S}}_{yy}w\|^2_{x=L}.
\end{align}

\medskip

{\it Step 4.7. Estimate on  $\big\la \mathcal{F}_{R7}(\hat{P}),\, \mathbf{q}_{xx}w^2\big\ra$.} 
Recall $\eqref{2.98-1}_2$, we rewrite $\mathcal{F}_{R7}(\hat{P})$ as
\begin{align}\label{3.98-1}
\mathcal{F}_{R7}(\hat{P})
&=\mu\v \Big(\frac{b_2}{p^{\v}} v^{\v} \Delta_{\v}\hat{P}\Big)_x - \mu \Big(\frac{b_1}{p^{\v}} u^{\v} \Delta_{\v}\hat{P}\Big)_y - \lambda \Big(\frac{b_1\rho^{\v}}{2\rho_sp^{\v}} u^{\v} \hat{P}_{yy}\Big)_y + \lambda\v\Big(\frac{b_2-b_1}{p^{\v}} (U^{\v}\cdot\nabla)\hat{P}_x \Big)_y \nonumber\\
&\quad + \lambda \v  \Big(\frac{b_2\rho^{\v}}{2\rho_sp^{\v}} u^{\v} \hat{P}_{xx}\Big)_y  + \lambda \v \big(\frac{b_2}{p^{\v}} U^{\v} \big)_x  \cdot\nabla\hat{P}_y - \lambda \v \big(\frac{b_2}{p^{\v}} U^{\v}\big)_y \cdot\nabla\hat{P}_x +  \lambda \v \big(\frac{b_2\rho^{\v}}{2\rho_s p^{\v}} u^{\v}\big)_x \hat{P}_{xy}\nonumber\\
&\quad - \lambda \v  \Big(\frac{b_2\rho^{\v}}{2\rho_sp^{\v}} u^{\v}\Big)_y \hat{P}_{xx}. 
\end{align}
Integrating by parts and using \eqref{3.87}-\eqref{3.94}, one has 
\begin{align}
&\Big\la \v \Big(\frac{b_2}{p^{\v}} v^{\v} \Delta_{\v}\hat{P}\Big)_x,\, \mathbf{q}_{xx}w^2\Big\ra \nonumber\\
&=-\Big\la \v \frac{b_2}{p^{\v}} v^{\v} \Delta_{\v}\hat{P} ,\, \mathbf{q}_{xxx}w^2\Big\ra + \Big\la \v \frac{b_2}{p^{\v}} v^{\v} \Delta_{\v}\hat{P} ,\, \mathbf{q}_{xx}w^2\Big\ra_{x=L} - \Big\la \v \frac{b_2}{p^{\v}} v^{\v} \Delta_{\v}\hat{P} ,\, \mathbf{q}_{xx}w^2\Big\ra_{x=0} \nonumber\\
&\lesssim   \|\sqrt{u_s}\v \mathbf{q}_{xxx}w\|\cdot \|\Delta_{\v}\hat{P}w\| +  \v^{\f12-} \|u_s\,\sqrt{\v} \mathbf{q}_{xx}w\|_{x=L}\cdot \|\sqrt{u_s}\Delta_{\v}\hat{P} w\|_{x=L} \nonumber\\
&\quad  + \|\sqrt{\v}\mathbf{q}_{xx}w\|_{x=0}\cdot \sqrt{\v}\|v_s \Delta_{\v}\hat{P} w\|_{x=0} \nonumber\\
&\lesssim \v^{\f12-} [\Phi]_{3,\hat{w}} + \v^{\f12-}\f{1}{\v} \mathcal{K}_{2,w}(\hat{P})  + \sqrt{\v} \bar{\mathcal{K}}_{3,w}(\hat{P}),
\end{align}
and
\begin{align}
&\Big\la \Big\{- \mu \frac{b_1}{p^{\v}} u^{\v} \Delta_{\v}\hat{P}  - \lambda \frac{b_1\rho^{\v}}{2\rho_sp^{\v}} u^{\v} \hat{P}_{yy}  + \lambda\v \frac{b_2-b_1}{p^{\v}} (U^{\v}\cdot\nabla)\hat{P}_x -  \lambda \v \frac{b_2\rho^{\v}}{2\rho_sp^{\v}} u^{\v} \hat{P}_{xx}\Big\}_y ,\, \mathbf{q}_{xx}w^2\Big\ra \nonumber\\
&= \Big\la \mu \frac{b_1}{p^{\v}} u^{\v} \Delta_{\v}\hat{P} + \lambda \frac{b_1\rho^{\v}}{2\rho_sp^{\v}} u^{\v} \hat{P}_{yy} - \lambda\v \frac{b_2-b_1}{p^{\v}} (U^{\v}\cdot\nabla)\hat{P}_x + \lambda \v \frac{b_2\rho^{\v}}{2\rho_sp^{\v}} u^{\v} \hat{P}_{xx} ,\, \mathbf{q}_{xxy}w^2\Big\ra  \nonumber\\
&\,\,\, +\Big\la \mu \frac{b_1}{p^{\v}} u^{\v} \Delta_{\v}\hat{P} + \lambda \frac{b_1\rho^{\v}}{2\rho_sp^{\v}} u^{\v} \hat{P}_{yy} - \lambda\v \frac{b_2-b_1}{p^{\v}} (U^{\v}\cdot\nabla)\hat{P}_x + \lambda \v \frac{b_2\rho^{\v}}{2\rho_sp^{\v}} u^{\v} \hat{P}_{xx} ,\, \mathbf{q}_{xx} 2ww_y\Big\ra \nonumber\\
&\lesssim \f{1}{N} [\Phi]_{3,\hat{w}} + \sqrt{\v}\|\mathbf{q}_{xx}w_y\|^2  + Na \f{1}{\v} \mathcal{K}_{2,w}(\hat{P}),
\end{align}
where we have used 
\begin{align}\label{3.122}
\sqrt{\v} \|\sqrt{u_s} \Delta_{\v}\hat{P} w\|^2_{x=0} 
&\lesssim  \sqrt{\v} \|\sqrt{u_s}\Delta_{\v}\hat{P}w\|^2_{x=L} + \|u_s\sqrt{\v}\Delta_{\v}\hat{P}_x w\|\cdot \|\Delta_{\v}\hat{P} w\| \nonumber\\
&\lesssim \f{1}{\sqrt{\v}} \mathcal{K}_{2,w}(\hat{P}) + \sqrt{\v} \bar{\mathcal{K}}_{3,w}(\hat{P}).
\end{align}

Using \eqref{3.87}-\eqref{3.94}, one also gets that
\begin{align}\label{3.123}
&\v \Big\la \big(\frac{b_2}{p^{\v}} U^{\v} \big)_x  \cdot\nabla\hat{P}_y - \big(\frac{b_2}{p^{\v}} U^{\v}\big)_y \cdot\nabla\hat{P}_x +  \big(\frac{b_2\rho^{\v}u^{\v}}{2\rho_s p^{\v}} \big)_x \hat{P}_{xy}  -  \Big(\frac{b_2\rho^{\v}u^{\v}}{2\rho_sp^{\v}} \Big)_y \hat{P}_{xx},\, \mathbf{q}_{xx}w^2\Big\ra\nonumber\\
& \lesssim  L^{\f12\alpha-} [\Phi]_{3,w} + L^{\f12\alpha-}\f{1}{\v} \mathcal{K}_{2,w}(\hat{P}).
\end{align}

Combining  \eqref{3.98-1}-\eqref{3.123}, one obtains
\begin{align}
\big\la \mathcal{F}_{R7}(\hat{P}),\, \mathbf{q}_{xx}w^2\big\ra
&\lesssim  \f{1}{N} [\Phi]_{3,\hat{w}} + \sqrt{\v}\|\mathbf{q}_{xx}w_y\|^2  + Na \f{1}{\v} \mathcal{K}_{2,w}(\hat{P}) +   \sqrt{\v} \bar{\mathcal{K}}_{3,w}(\hat{P}).
\end{align}

\smallskip

{\it Step 4.8.} Combining above estimates in steps 4.1-4.7, 
one concludes \eqref{3.1-0}. Therefore the proof of Lemma \ref{lem3.1} is completed. $\hfill\Box$

\medskip

\begin{lemma}\label{lem3.2}
It holds that 
\begin{align}\label{3.103}
&L^{-\f{1}{4}}\Big\{\|u_s \mathbf{q}_{yy}w\|_{x=L}^2 + \|\sqrt{u_{sy}}\mathbf{q}_{yy}\|_{y=0}^2 + \|\sqrt{u_s}\big(\mathbf{q}_{yyy},\sqrt{\v}\mathbf{q}_{xyy}, \v \mathbf{q}_{xxy}\big)w\|^2\Big\}  \nonumber\\
&\lesssim   L^{\f14-}  [\Phi]_{3,\hat{w}} + \v^{1-} [\Phi]_{4,w}
  + L^{\f14-}\|\mathcal{F}_{R}w\|^2,
\end{align}
with
\begin{align}\label{3.103-0}
\|\mathcal{F}_{R}w\|^2&\lesssim [[p]]_{3,w} + [[[\zeta]]] _{3,w}  + \mathcal{K}_{3,w}(\hat{P})  + [[[\hat{\mathbf{S}}]]]_{3,w}  +  a^2  \|(\phi, p, S(\zeta)) \hat{w}\|^2_{\mathbf{X}^{\v}}  \nonumber\\
&\quad +  L^{\f12}  \|(0, \hat{P}, \hat{\mathbf{S}})\hat{w}\|^2_{\mathbf{Y}^{\v}} + \v^{1-} \|u_s\sqrt{\v} \Delta_{\v}\tilde{T}_x w\|^2 + \|\mathcal{F}_{R8}w\|^2.
\end{align}
\end{lemma}

\noindent{\bf Proof.} Multiplying \eqref{2.121} by $\mathbf{q}_{yy}w^2$, one has 
\begin{align}\label{3.104}
	&\big\la\big(b_1 m u_s^2 \mathbf{q}_{xy}\big)_y,\, \mathbf{q}_{yy}w^2\big\ra + \big\la \v \big(b_2 mu_s^2 \mathbf{q}_{xx}\big)_x ,\, \mathbf{q}_{yy}w^2\big\ra  \nonumber\\
	&-\mu\Big\la \Big(\frac{b_1}{\rho_s} \Delta_{\v}\Phi_y\Big)_y,\, \mathbf{q}_{yy}w^2\Big\ra
	- \mu \Big\la \v \Big(\frac{b_2}{\rho_s}   \Delta_{\v}\Phi_x\Big)_x,\, \mathbf{q}_{yy}w^2\Big\ra\nonumber\\
	& = \big\la \mathcal{F}_{R}(\hat{P}, \hat{\mathbf{S}}, \hat{\Phi}, \zeta, p, \phi, S),\, \mathbf{q}_{yy}w^2\big\ra + O(1)   \v^{\f12 N_0} \big\{[\Phi]_{3,\hat{w}} + [\Phi]_{4,\hat{w}}\big\},
\end{align}
where we have used similar arguments as in \eqref{3.18-1}.
We divide the rest estimate of \eqref{3.104}  into several steps.

\smallskip

\noindent{\it Step 1. Estimate on the Rayleigh term.} Integrating by parts and using Lemma \ref{lemA.3}, one obtains that 
\begin{align}\label{3.85}
&\Big\la (b_1 m u_s^2 \mathbf{q}_{xy})_y ,\, \mathbf{q}_{yy}w^2\Big\ra\nonumber\\
&=\Big\la  b_1m u_s^2 \mathbf{q}_{xyy} ,\, \mathbf{q}_{yy}w^2\Big\ra + \Big\la (b_1m)_y u_s^2 \mathbf{q}_{xy} ,\, \mathbf{q}_{yy}w^2\Big\ra + \Big\la 2b_1 m u_su_{sy} \mathbf{q}_{xy} ,\, \mathbf{q}_{yy}w^2\Big\ra\nonumber\\
&=\f12 \Big\la  b_1m u_s^2 \mathbf{q}_{yy} ,\, \mathbf{q}_{yy}w^2\Big\ra_{x=L} - \f12 \Big\la  (b_1m u_s^2)_x \mathbf{q}_{yy} ,\, \mathbf{q}_{yy}w^2\Big\ra  + \Big\la (b_1m)_y u_s^2 \mathbf{q}_{xy} ,\, \mathbf{q}_{yy}w^2\Big\ra \nonumber\\
&\quad + \Big\la 2b_1 m u_su_{sy} \mathbf{q}_{xy} ,\, \mathbf{q}_{yy}w^2\Big\ra\nonumber\\
&= \f12 \|\sqrt{b_1m}u_s \mathbf{q}_{yy}w\|_{x=L}^2  +  O(1) \sqrt{L} [\Phi]_{3,\hat{w}}.
\end{align}
It is clear that 
\begin{align}\label{3.85-1}
\big\la \v \big(b_2 mu_s^2 \mathbf{q}_{xx}\big)_x ,\, \mathbf{q}_{yy}w^2\big\ra 
&=\big\la \v b_2 mu_s^2 \mathbf{q}_{xxx} + \v \big(b_2 mu_s^2\big)_x \mathbf{q}_{xx} ,\, \mathbf{q}_{yy}w^2\big\ra  \lesssim \sqrt{L} [\Phi]_{3,\hat{w}}.
\end{align}

\medskip

\noindent{\it Step 2. Estimate on the viscous term.} We aim to prove the following facts
\begin{align}\label{3.108-0}
	&-\Big\la \Big(\frac{b_1}{\rho_s} \Delta_{\v}\Phi_y\Big)_y,\, \mathbf{q}_{yy}w^2\Big\ra
	-\Big\la \v \Big(\frac{b_2}{\rho_s}   \Delta_{\v}\Phi_x\Big)_x,\, \mathbf{q}_{yy}w^2\Big\ra\nonumber\\
	&\lesssim  - \|\sqrt{u_s}\big(\mathbf{q}_{yyy},\, \sqrt{\v}\mathbf{q}_{xyy}, \v \mathbf{q}_{xxy}\big)w\|^2 - \|\sqrt{u_{sy}}\mathbf{q}_{yy}\|_{y=0}^2 + L^{1-}[\Phi]_{3,\hat{w}} + \v [\Phi]_{4,w}.
\end{align}

\noindent{\it Step 2.1. Estimate on $\dis \Big\la \Big(\frac{b_1}{\rho_s} \Delta_{\v}\Phi_y\Big)_y,\, \mathbf{q}_{yy}w^2\Big\ra$.} It is clear that 
\begin{align}\label{3.108}
	\Big\la \Big(\frac{b_1}{\rho_s} \Delta_{\v}\Phi_y\Big)_y,\, \mathbf{q}_{yy}w^2\Big\ra
	&=\Big\la \Big(\frac{b_1}{\rho_s} \Phi_{yyy}\Big)_y,\, \mathbf{q}_{yy}w^2\Big\ra + \Big\la  \frac{b_1}{\rho_s} \v \Phi_{xxyy} ,\, \mathbf{q}_{yy}w^2\Big\ra 
	+ O(1)\sqrt{\v}  [\Phi]_{3,\hat{w}}.
\end{align}
For 1th term on RHS of \eqref{3.108}, integrating by parts in $y$, one gets that 
\begin{align}\label{3.109}
\eqref{3.108}_1&=-\Big\la \frac{b_1}{\rho_s}  \Phi_{yyy} ,\, \mathbf{q}_{yyy}w^2\Big\ra  
	-\Big\la  \frac{b_1}{\rho_s}  \Phi_{yyy} ,\, 2\mathbf{q}_{yy}ww_y\Big\ra  -\Big\la \frac{b_1}{\rho_s}  \Phi_{yyy} ,\, \mathbf{q}_{yy}w^2\Big\ra_{y=0}.
\end{align}
Using \eqref{A.34}, we have 
\begin{align}\label{3.110}
\eqref{3.109}_1
	&= - \|\sqrt{b_1\rho_s^{-1}m u_s}\mathbf{q}_{yyy}w\|^2 - \Big\la 3  \frac{b_1}{\rho_s}  \bar{u}_{sy}\mathbf{q}_{yy},\, \mathbf{q}_{yyy}w^2\Big\ra - \Big\la 3 \frac{b_1}{\rho_s}  \bar{u}_{syy}\mathbf{q}_{y},\, \mathbf{q}_{yyy}w^2\Big\ra\nonumber\\
	&\quad - \Big\la  \frac{b_1}{\rho_s} \bar{u}_{syyy}\mathbf{q},\, \mathbf{q}_{yyy}w^2\Big\ra.
\end{align}
Integrating in $y$, one obtains
\begin{align}\label{3.111}
\eqref{3.110}_2
&= \f32  \Big\la \frac{b_1}{\rho_s} \bar{u}_{sy}\mathbf{q}_{yy},\, \mathbf{q}_{yy}w^2\Big\ra_{y=0} + \f32 \Big\la \big(\frac{b_1}{\rho_s} \bar{u}_{sy} w^2\big)_y \mathbf{q}_{yy},\, \mathbf{q}_{yy}w^2\Big\ra\nonumber\\
&= \f32 \|\sqrt{b_1\rho_s^{-1}m u_{sy}}\mathbf{q}_{yy}\|^2_{y=0} + O(1) L^{1-} [\Phi]_{3,\hat{w}},
\end{align}
\begin{align}\label{3.112}
\eqref{3.110}_3	&=\Big\la 3\frac{b_1}{\rho_s}  \bar{u}_{syy}\mathbf{q}_{y},\, \mathbf{q}_{yy}w^2\Big\ra_{y=0}
	 +  \Big\la 3\big(\frac{b_1}{\rho_s}  \bar{u}_{syy} w^2\big)_y \mathbf{q}_{y},\, \mathbf{q}_{yy}w^2\Big\ra \nonumber\\
	&\quad  +  \Big\la 3 \frac{b_1}{\rho_s}  \bar{u}_{syy} ,\, |\mathbf{q}_{yy}|^2w^2\Big\ra \lesssim  L [\Phi]_{3,\hat{w}},  
\end{align}
and 
\begin{align}\label{3.113}
\eqref{3.110}_4	&=\Big\la \frac{b_1}{\rho_s} \bar{u}_{syyy}\mathbf{q}_y,\, \mathbf{q}_{yy}w^2\Big\ra + \Big\la \big(\frac{b_1}{\rho_s} \bar{u}_{syyy}w^2\big)_y \mathbf{q},\, \mathbf{q}_{yy}\Big\ra 
	\lesssim L [\Phi]_{3,\hat{w}}.
\end{align}
Then we have from \eqref{3.110}-\eqref{3.113} that 
\begin{align}\label{3.113-1}
\eqref{3.109}_1 &\leq  - \|\sqrt{b_1\rho_s^{-1}m u_s}\mathbf{q}_{yyy}w\|^2 +\f32 \|\sqrt{b_1\rho_s^{-1}m u_{sy}}\mathbf{q}_{yy}\|^2_{y=0} + L^{1-} [\Phi]_{3,\hat{w}}.
\end{align}

\smallskip

It is clear that 
\begin{align}\label{3.98}
\eqref{3.109}_2 \lesssim \|\Phi_{yyy}w\| \cdot \|\mathbf{q}_{yy}w_y\| \lesssim L^{\f12-} [\Phi]_{3,\hat{w}},
\end{align}
and
\begin{align}\label{3.99-1}
\eqref{3.109}_3 
	&= -\Big\la \frac{b_1}{\rho_s} [3\bar{u}_{sy}\mathbf{q}_{yy} + 3 \bar{u}_{syy} \mathbf{q}_y] ,\, \mathbf{q}_{yy}w^2\Big\ra_{y=0}  
	 =-3\|\sqrt{b_1\rho_s^{-1}m u_{sy}}\mathbf{q}_{yy}\|^2_{y=0} +  L  [\Phi]_{3,1}. 
\end{align}
Hence it follows from \eqref{3.109} and \eqref{3.113-1}-\eqref{3.99-1} that 
\begin{align}\label{3.117}
\eqref{3.108}_1 &\leq   - \|\sqrt{b_1\rho_s^{-1} m u_s}\mathbf{q}_{yyy}w\|^2 - \f32\|\sqrt{mb_1\rho_s^{-1}} \sqrt{u_{sy}}\mathbf{q}_{yy}w\|^2_{y=0}  +  CL^{1-} [\Phi]_{3,\hat{w}}.
\end{align}

\medskip

For 2th term on RHS \eqref{3.108}, integrating by parts in $x$, we have 
\begin{align}\label{3.101}
\eqref{3.108}_2
	&= - \Big\la \v  \frac{b_1}{\rho_s}\big[\bar{u}_{sxyy} \mathbf{q} + \bar{u}_{syy}\mathbf{q}_{x} 
	+2\bar{u}_{sxy} \mathbf{q}_y + 2\bar{u}_{sy} \mathbf{q}_{xy}\big] ,\, \mathbf{q}_{xyy}w^2\Big\ra   + C \v  [\Phi]_{3,\hat{w}}\nonumber\\
	&\quad 
	 - \|\sqrt{m b_1\rho_s^{-1}u_s}\sqrt{\v}\mathbf{q}_{xyy}w\|^2.
\end{align}
Integrating by parts in $y$, one has 
\begin{align}\label{3.103-1}
	& - \Big\la \frac{b_1}{\rho_s} \v \big[\bar{u}_{sxyy} \mathbf{q} + \bar{u}_{syy}\mathbf{q}_{x} +2 \bar{u}_{sxy} \mathbf{q}_y\big] ,\, \mathbf{q}_{xyy}w^2\Big\ra \nonumber\\
	&=2 \v \Big\la  \frac{b_1}{\rho_s} \bar{u}_{sxy} \mathbf{q}_y ,\, \mathbf{q}_{xy}w^2\Big\ra_{y=0}
	+\Big\la \frac{b_1}{\rho_s} \v \big[\bar{u}_{sxyy} \mathbf{q}_y + \bar{u}_{syy}\mathbf{q}_{xy} +2\bar{u}_{sxy} \mathbf{q}_{yy}\big] ,\, \mathbf{q}_{xy}w^2\Big\ra \nonumber\\
	&\quad + \v  \Big\la \big[(\frac{b_1}{\rho_s}\bar{u}_{sxyy}w^2)_y \mathbf{q} + (\frac{b_1}{\rho_s}\bar{u}_{syy}w^2)_y\mathbf{q}_{x} +2 (\frac{b_1}{\rho_s}\bar{u}_{sxy}w^2)_y \mathbf{q}_y\big] ,\, \mathbf{q}_{xy}\Big\ra \lesssim \v [\Phi]_{3,w},
\end{align}
and
\begin{align}\label{3.121}
	- \Big\la 2\v\frac{b_1}{\rho_s} \bar{u}_{sy} \mathbf{q}_{xy}  ,\, \mathbf{q}_{xyy}w^2\Big\ra 
	&\lesssim  \v \|\sqrt{u_{sy}}\mathbf{q}_{xy}\|_{y=0}^2
	+  \v  \Big\la \big(\frac{b_1}{\rho_s}\bar{u}_{sy} w^2\big)_y \mathbf{q}_{xy}  ,\, \mathbf{q}_{xy}w^2\Big\ra \lesssim \v  [\Phi]_{3,\hat{w}}.
\end{align}
Substituting \eqref{3.103-1}-\eqref{3.121} into \eqref{3.101}, one obtains 
\begin{align}\label{3.105}
\eqref{3.108}_2 \leq -\mu \|\sqrt{b_1\rho_s^{-1}mu_s}\sqrt{\v}\mathbf{q}_{xyy}w\|^2  +C \v  [\Phi]_{3,\hat{w}}.
\end{align}


Combining \eqref{3.108}, \eqref{3.117} and \eqref{3.105}, one has
\begin{align}\label{3.107}
	\Big\la\Big(\frac{b_1}{\rho_s} \Delta_{\v}\Phi_y\Big)_y,\, \mathbf{q}_{yy}w^2\Big\ra
	&\leq  -\|\sqrt{b_1\rho_s^{-1}m u_s}(\mathbf{q}_{yyy},\sqrt{\v}\mathbf{q}_{xyy})w\|^2 - \|\sqrt{m b_1\rho_s^{-1}} \sqrt{u_{sy}}\mathbf{q}_{yy}w\|^2_{y=0}  \nonumber\\
	&\quad  + L^{1-}[\Phi]_{3,\hat{w}} 
\end{align}

\noindent{\it Step 2.2. Estimate on $\dis \Big\la \v \Big(\frac{b_2}{\rho_s}   \Delta_{\v}\Phi_x\Big)_x,\, \mathbf{q}_{yy}w^2\Big\ra$.} It is clear that 
\begin{align}\label{3.125}
	\Big\la \v \Big(\frac{b_2}{\rho_s}   \Delta_{\v}\Phi_x\Big)_x,\, \mathbf{q}_{yy}w^2\Big\ra
	&=\Big\la  \v  \frac{b_2}{\rho_s}\Phi_{xxyy} ,\, \mathbf{q}_{yy}w^2\Big\ra + \Big\la \v^2 \Big(\frac{b_2}{\rho_s} \Phi_{xxx}\Big)_x,\, \mathbf{q}_{yy}w^2\Big\ra +  C \v  [\Phi]_{3,\hat{w}} 
\end{align}
Similar as \eqref{3.101}-\eqref{3.105}, it holds that 
\begin{align}\label{3.126}
\eqref{3.125}_1	&\leq - \|\sqrt{b_2\rho_s^{-1}m u_s}\sqrt{\v}\mathbf{q}_{xyy}w\|^2  + C \v  [\Phi]_{3,\hat{w}} .
\end{align}

For 2th term on RHS of \eqref{3.125}, integrating by parts in $x$ and $y$, one obtains that 
\begin{align}\label{3.127}
\eqref{3.125}_2
&=\v^2 \Big\la\frac{b_2}{\rho_s} \Phi_{xxy},\, \mathbf{q}_{xy}w^2\Big\ra_{x=L}  -\v^2 \Big\la\frac{b_2}{\rho_s} \Phi_{xxy},\, \mathbf{q}_{xxy}w^2\Big\ra  + C \v  [\Phi]_{3,\hat{w}} \nonumber\\
&=- \|\sqrt{b_2\rho_s^{-1}m u_s}\v \mathbf{q}_{xxy}w\|^2 + C\v  [\Phi]_{4,w} + C\v [\Phi]_{3,\hat{w}}
\end{align}
where we have used the following two estimates
\begin{align*}
\v^2 \Big\la\frac{b_2}{\rho_s} \Phi_{xxy},\, \mathbf{q}_{xy}w^2\Big\ra_{x=L} 
&\lesssim \v \|\v^{\f32}\Phi_{xxxy}w\|^{\f12}\cdot \|\sqrt{\v}\Phi_{xxy}w\|^{\f12} \cdot \|\mathbf{q}_{xy}w\|_{x=L} \nonumber\\
&\lesssim \v  [\Phi]_{4,w} + \v [\Phi]_{3,\hat{w}},
\end{align*}
and
\begin{align*}
- \v^2 \Big\la\frac{b_2}{\rho_s} \Phi_{xxy},\, \mathbf{q}_{xxy}w^2\Big\ra
	&=- \|\sqrt{b_2\rho_s^{-1}m u_s}\v \mathbf{q}_{xxy}w\|^2 + \v  [\Phi]_{3,w}.
\end{align*}

Hence it follows from \eqref{3.125}-\eqref{3.127} that 
\begin{equation}\label{3.130}
\Big\la \v \Big(\frac{b_2}{\rho_s}   \Delta_{\v}\Phi_x\Big)_x,\, \mathbf{q}_{yy}w^2\Big\ra 
\leq - \|\sqrt{b_2\rho_s^{-1}m u_s}\big(\sqrt{\v}\mathbf{q}_{xyy},\v \mathbf{q}_{xxy}\big)w\|^2 + C\v \big\{  [\Phi]_{3,\hat{w}} +   [\Phi]_{4,w} \big\}.
\end{equation}
 Combining \eqref{3.107} and \eqref{3.130}, one obtains \eqref{3.108-0}.

\smallskip

\noindent{\it Step 3.} Substituting \eqref{3.85}-\eqref{3.85-1} and \eqref{3.108-0} into \eqref{3.2},  we have
\begin{align*}
	&\|\sqrt{u_s}\big(\mathbf{q}_{yyy},\sqrt{\v}\mathbf{q}_{xyy}, \v \mathbf{q}_{xxy}\big)w\|^2 + \|\sqrt{u_{sy}}\mathbf{q}_{yy}\|_{y=0}^2 + \|u_s \mathbf{q}_{yy}w\|_{x=L}^2 \nonumber\\
	&\lesssim L^{\f12-}[\Phi]_{3,\hat{w}} + \v [\Phi]_{4,w} 
	+  L^{\f12-}\|\mathcal{F}_{R}w\|^2,
\end{align*}
which  concludes \eqref{3.103}.  

\smallskip

\noindent{\it Step 4. Estimate on the remainder term $\|\mathcal{F}_{R}w\|^2$.} It is clear that 
\begin{align}\label{3.133}
\begin{split}
\|\mathcal{F}_{R1}w\|^2 &\lesssim \|\hat{d}_{11}\zeta_{xy}w\|^2 + \|\big(b_1\hat{d}_{11}\big)_y \zeta_x w\|^2 \lesssim a^4[[[\zeta]]]_{2,\hat{w}},\\ 
\|\mathcal{F}_{R2}w\|^2 &\lesssim \|\v \big(b_2\, \mathfrak{g}_2\big)_xw\|^2 +  \|\big(b_1\mathfrak{g}_1\big)_yw\|^2 \lesssim L^{\f12\alpha-} [\phi]_{3,\hat{w}}.
\end{split}
\end{align}

\smallskip

For $\mathcal{F}_{R3}$, we note that 
\begin{align}\label{3.134}
\begin{split}
&\|(b_1 [\tilde{d}_{12}+\cdots]p_y)_y w\|^2 + \|(\v b_2\hat{d}_{22} p_y)_xw\|^2 + \|(\v b_2 [\tilde{d}_{21}+\cdots]p_x)_x w\|^2\lesssim a^2 [[p]]_{2,w},\\
&\|(b_1 \hat{d}_{12}S_y)_y w\|^2\lesssim a^2 \|(S_{yy}, S_y)w\|^2 \lesssim L^{\f12}[[[S]]]_{2,w},\\
&\|(\v b_2 \hat{d}_{21}S_x)_x w\|^2 \lesssim \v \|(\sqrt{\v}S_{xx}, \sqrt{\v}S_x)w\|^2 \lesssim \v [[[S]]]_{2,\hat{w}},\\ 
&\|(\v b_2 \hat{d}_{22}S_y)_x w\|^2 \lesssim \v^2 \|(S_{xy}, S_y)w\|^2\lesssim \v^2  [[[S]]]_{2,w}.
\end{split}
\end{align}

Using Lemma \ref{lemA.5-1}, one has
\begin{align} 
\begin{split}
\|\big(b_1[\tilde{d}_{13}+\cdots] p\big)_yw\|^2 +\|\v \big(b_2[\tilde{d}_{23}+\cdots] p\big)_xw\|^2 
&\lesssim L^2 \f{1}{\v} [[p]]_{2,\hat{w}},\\
\|\big(b_1\hat{d}_{13} S\big)_y w\|^2 + \|\v \big(b_2\hat{d}_{23} S\big)_x w\|^2 
&\lesssim  L [[[S]]]_{2,\hat{w}}. 
\end{split}
\end{align}
It follows from \eqref{3.67-0} that 
\begin{align}\label{3.135-2}
&\|(b_1 I_u(\phi))_yw\|^2 + \|(\v b_2 I_v(\phi))_x w\|^2 + 
\|\big\{b_1\hat{d}_{12}\big(\f{T_{sy}}{p_s} q\big)_y\big\}_y w\|^2 + \|\big\{b_1\hat{d}_{13}\f{T_{sy}}{p_s} q\big\}_y w\|^2  \nonumber\\
&+ \|\v \big\{b_1\hat{d}_{21}\big(\f{T_{sy}}{p_s} q\big)_x\big\}_x w\|^2 + \|\v\big\{b_2\hat{d}_{22}\big(\f{T_{sy}}{p_s} q\big)_y\big\}_x w\|^2  + \|\v\big\{b_2\hat{d}_{23}\f{T_{sy}}{p_s} q\big\}_x w\|^2  
\lesssim a^2 [\phi]_{3,\hat{w}}.
\end{align}

Hence, noting $\eqref{2.98}_3$ and \eqref{3.66-0}-\eqref{3.66-1}, then we have from \eqref{3.134}-\eqref{3.135-2} that 
\begin{align}\label{3.136} 
\|\mathcal{F}_{R3}w\|^2 
&\lesssim  a^2 [\phi]_{3,\hat{w}} + L^2 \f{1}{\v} [[p]]_{2,\hat{w}} + L^{\f12} [[[S]]]_{2,\hat{w}}.
\end{align}

\medskip

For $\mathcal{F}_{R4}$, noting \eqref{3.94}, we have from \eqref{3.92-1} that 
\begin{align}\label{3.137}
\|\mathcal{F}_{R4}  w\|^2 
&\lesssim 
\v^{1-} \|u_s \sqrt{\v} \Delta_{\v}\tilde{T}_x w\|^2 + \|u_s \nabla^2_{\v}\tilde{T}_y w\|^2 +a^2\|(\Delta_{\v}\tilde{T}, \v\tilde{T}_{xy}, \v\tilde{T}_{yy})w\|^2 \nonumber\\
&\lesssim \v^{1-} \|u_s\sqrt{\v} \Delta_{\v}\tilde{T}_x w\|^2  + [[[\zeta]]] _{3,w} + [[p]]_{3,w}  + L\f{1}{\v}[[p]]_{2,\hat{w}} \nonumber\\
&\quad + a^2 [\phi]_{3,\hat{w}} + L^{\f12}[[[S]]]_{2,w},
\end{align}
where we have used Lemma \ref{lemA.6} and the following estimate
\begin{align}\label{3.137-1}
\|u_s \nabla^2_{\v}\tilde{T}_y w\|^2 \lesssim [[[\zeta]]] _{3,w}  + [[p]]_{3,w} + L\f{1}{\v}[[p]]_{2,\hat{w}} + a^2 [\phi]_{3,\hat{w}}.
\end{align}

Using \eqref{3.94}, one has from $\eqref{2.98}_5$ and \eqref{2.96-1} that 
\begin{align}
\|\mathcal{F}_{R5} w\|^2 
&\lesssim \|\big(b_1 G_{12}\big)_y w\|^2 + \v^2 \|\big(b_2 G_{22}\big)_x w\|^2  \lesssim L^2 \f{1}{\v} [[p]]_{2,\hat{w}} + a^2 [\phi]_{3,\hat{w}}.
\end{align}
It follows from $\eqref{2.98-1}_1$, \eqref{3.98-1} and \eqref{3.94} that 
and 
\begin{align}
\|\mathcal{F}_{R6}w\|^2&\lesssim \|u_s \hat{\mathbf{S}}_{yyy}w\|^2 + \|\v u_s \hat{\mathbf{S}}_{xxy}w\|^2 + a^2 \|(\hat{\mathbf{S}}_{yy}, \v \hat{\mathbf{S}}_{xy})w\|^2 \nonumber\\
&\lesssim [[[\hat{\mathbf{S}}]]]_{3,w} + L^{\f12} [[[\hat{\mathbf{S}}]]]_{2,w},\\
\|\mathcal{F}_{R7}(\hat{P}) w\|^2 &\lesssim \v \|u_s \sqrt{\v} \Delta_{\v}\hat{P}_x w\|^2 + \|u_s \nabla^2_{\v}\hat{P}_y w\|^2 +  a^2 \|\nabla^2_{\v}\hat{P} w\|^2\nonumber\\
&\lesssim  \mathcal{K}_{3,w}(\hat{P})  + \v^{1-} \bar{\mathcal{K}}_{3,w}(\hat{P})  +  a^2 \mathcal{K}_{2,w}(\hat{P}).\label{3.138}
\end{align}

Finally, combing  \eqref{3.133}-\eqref{3.138}, we conclude \eqref{3.103-0}. 
Therefore the proof of Lemma \ref{lem3.2} is completed. $\hfill\Box$

\medskip

\begin{lemma}\label{lem3.3}
It holds that 
\begin{align}\label{3.138-0}
\|u_s\nabla_{\v}\mathbf{q}_xw\|^2  
&\lesssim  L^{\f12-} \Big\{\|\mathbf{q}_{xx}w_y\|^2 +  [\Phi]_{3,w} + [\Phi]_{4,w} 
 + \|(\phi, p, S \, (\zeta)) \hat{w}\|^2_{\mathbf{X}^{\v}}  + \|(0, \hat{P}, \hat{\mathbf{S}})w\|^2_{\mathbf{Y}^{\v}} \nonumber\\
&\quad  + \f{1}{\v} \|(0, p, S )\|^2_{\mathbf{X}^{\v}}\mathbf{1}_{\{w=w_0\}}  +  \|u_s \sqrt{\v} \Delta_{\v}\tilde{T}_xw\|^2  + \f{1}{\v}\|\mathcal{F}_{R8}w\|^2  + \v^{N_0}\Big\} .
\end{align}
\end{lemma}

\noindent{\bf Proof.} Multiplying \eqref{2.121} by $\mathbf{q}_{x}w^2$, one has 
\begin{align}\label{3.139}
&\big\la\big(b_1 m u_s^2 \mathbf{q}_{xy}\big)_y,\, \mathbf{q}_{x}w^2\big\ra + \big\la \v \big(b_2 mu_s^2 \mathbf{q}_{xx}\big)_x ,\, \mathbf{q}_{x} w^2\big\ra  \nonumber\\
&-\Big\la\mu\,\Big(\frac{b_1}{\rho_s} \Delta_{\v}\Phi_y\Big)_y,\, \mathbf{q}_{x}w^2\Big\ra-\Big\la\mu\,\v \Big(\frac{b_2}{\rho_s}   \Delta_{\v}\Phi_x\Big)_x,\, \mathbf{q}_{x}w^2\Big\ra \nonumber\\
&= \big\la\mathcal{F}_{R}(\hat{P}, \hat{\mathbf{S}}, \hat{\Phi}, \zeta, p, \phi, S), \, \mathbf{q}_{x}w^2\big\ra
 + O(1)   \v^{\f12 N_0} \big\{[\Phi]_{3,\hat{w}} + [\Phi]_{4,\hat{w}}\big\},
\end{align}
where we have used similar arguments as in \eqref{3.18-1}. We divide the rest proof into several steps. 

\smallskip

{\it Step 1. Estimate on  Rayleigh term.}  Integrating by parts in $y$, one has
\begin{align}\label{3.140}
\big\la\big(b_1 m u_s^2 \mathbf{q}_{xy}\big)_y,\, \mathbf{q}_{x}w^2\big\ra &= -\big\la b_1m u_s^2 \mathbf{q}_{xy},\, \mathbf{q}_{xy}w^2\big\ra - \big\la b_1 m u_s^2 \mathbf{q}_{xy},\, \mathbf{q}_{x}2ww_y\big\ra \nonumber\\
&= - \big(1-L\big) \|\sqrt{b_1}u_s \mathbf{q}_{xy}w\|^2 - O(1)L \|\mathbf{q}_{xx}w_y\|^2.
\end{align}
Integrating by parts in $x$ and using Lemma \ref{lemA.4}, one gets
\begin{align}\label{3.129-1}
\big\la \v \big(b_2 mu_s^2 \mathbf{q}_{xx}\big)_x ,\, \mathbf{q}_{x} w^2\big\ra 
&= - \|\sqrt{m b_2}\sqrt{\v}u_s \mathbf{q}_{xx}w\|^2 + O(1) \|\sqrt{\v}\mathbf{q}_xw\|_{x=L} \|u_s\sqrt{\v}\mathbf{q}_{xx}w\|_{x=L} \nonumber\\
&\quad + O(1) \|\sqrt{\v}\mathbf{q}_xw\|_{x=0}\cdot \|\sqrt{\v}\mathbf{q}_{xx}w\|_{x=0}\nonumber\\
&= - \big(1-L\big) \|\sqrt{b_2}\sqrt{\v}u_s \mathbf{q}_{xx}w\|^2 - CL^{1-} [\Phi]_{3,\hat{w}}.
\end{align}

\medskip

{\it Step 2. Estimate on the viscous term.} Integrating by parts in $y$, one obtains that
\begin{align}\label{3.141}
\Big\la\Big(\frac{b_1}{\rho_s} \Phi_{yyy}\Big)_y,\, \mathbf{q}_{x}w^2\Big\ra 
&= -\Big\la\frac{b_1}{\rho_s} \big[ u_s\mathbf{q}_{yyy} + 3u_{sy}\mathbf{q}_{yy} + 3 u_{syy} \mathbf{q}_y + u_{syyy} \mathbf{q}\big],\, \mathbf{q}_{xy}w^2\Big\ra \nonumber\\
&\quad  + \Big\la \frac{b_1}{\rho_s} \Phi_{yy},\, \mathbf{q}_{xy} 2ww_y\Big\ra + \Big\la \big(\frac{b_1}{\rho_s}2ww_y\big)_y \Phi_{yy},\, \mathbf{q}_{xy} \Big\ra \nonumber\\
&= \Big\la\frac{b_1}{\rho_s} u_s\mathbf{q}_{yy},\, \mathbf{q}_{xyy}w^2\Big\ra  + O(1)\sqrt{L} [\Phi]_{3,\hat{w}} 
= O(1)\sqrt{L} [\Phi]_{3,\hat{w}},
\end{align}
and
\begin{align}\label{3.142}
\Big\la\v \Big(\frac{b_2}{\rho_s} \Phi_{xyy}\Big)_x + \v \Big(\frac{b_2}{\rho_s} \Phi_{xxy}\Big)_y,\, \mathbf{q}_{x}w^2\Big\ra  
&\lesssim L^{1-} [\Phi]_{3,\hat{w}} .
\end{align}
We integrate by parts in $x$ to have
\begin{align}\label{3.144}
 \Big\la\v^2 \Big(\frac{b_2}{\rho_s} \Phi_{xxx}\Big)_x,\, \mathbf{q}_{x}w^2\Big\ra  
 &=-\Big\la\v^2 \frac{b_2}{\rho_s} \Phi_{xxx},\, \mathbf{q}_{xx}w^2\Big\ra - \Big\la\v^2 \frac{b_2}{\rho_s} \Phi_{xxx},\, \mathbf{q}_{x}w^2\Big\ra_{x=0} \nonumber\\
	&\lesssim \sqrt{\v} \|\v\Phi_{xxx}w\|\cdot \|\sqrt{\v}\mathbf{q}_{xx}w\| + \|\sqrt{\v}\mathbf{q}_xw\|_{x=0} \|\v^{\f32}\Phi_{xxx}w\|_{x=0} \nonumber\\
	& \lesssim  L^{\f12-} [\Phi]_{3,\hat{w}} + L^{\f12-} [\Phi]_{4,\hat{w}},
\end{align}
where we have used Lemmas \ref{lemA.3} \& \ref{lemA.4}, and  the following  trace estimate
\begin{align}\label{3.145}
\|\v^{\f32}\Phi_{xxx}w\|_{x=0}^2 &\leq 2\|\v^2\Phi_{xxxx}w\|\cdot \|\v \Phi_{xxx}w\| \lesssim [\Phi]_{4,w}^{\f12}\cdot  [\Phi]_{3,\hat{w}}^{\f12}.
\end{align}

Combining \eqref{3.141}-\eqref{3.144}, we obtain that 
\begin{align}\label{3.146-1}
\Big\la\Big(\frac{b_1}{\rho_s} \Delta_{\v}\Phi_y\Big)_y,\, \mathbf{q}_{x}w^2\Big\ra
+\Big\la\v \Big(\frac{b_2}{\rho_s} \Delta_{\v}\Phi_{x}\Big)_x,\, \mathbf{q}_{x}w^2\Big\ra 
&\lesssim  L^{\f12-} [\Phi]_{3,w} +  L^{\f12-} [\Phi]_{4,w}.
\end{align}

\smallskip

\noindent{\it Step 3. Estimate on $\big\la\mathcal{F}_{R}(\hat{P}, \hat{\mathbf{S}}, \hat{\Phi}, \zeta, p, \phi, S), \, \mathbf{q}_{x}w^2\big\ra$.} We divide it into several substeps.

\noindent{\it Step 3.1.} For $\mathcal{F}_{R1}, \mathcal{F}_{R2}$, it is clear that 
\begin{align}\label{3.147}
\big\la \mathcal{F}_{R1}, \, \mathbf{q}_x w^2 \big\ra &
= \big\la b_1\hat{d}_{11}S_x, \, \mathbf{q}_{xy} w^2 \big\ra + \big\la b_1\hat{d}_{11}S_x, \, \mathbf{q}_{x}2 ww_y \big\ra\nonumber\\
&\lesssim \f1N \|u_s\mathbf{q}_{xy}w\|^2 + \|\mathbf{q}_{x}w_y\|^2  + N \big\{\|u_s S_{x}w\|^2 + \|\v S_{x}w\|^2 \big\}\nonumber\\
&\lesssim \f1N \|u_s\mathbf{q}_{xy}w\|^2 +L \|\mathbf{q}_{xx}w_y\|^2 + N L^{\f12} [[[S]]]_{2,\hat{w}},
\end{align}
and
\begin{align}\label{3.148-1}
\big\la \mathcal{F}_{R2}, \, \mathbf{q}_x w^2 \big\ra &= \big\la b_1\mathfrak{g}_1, \, \mathbf{q}_{xy} w^2 \big\ra + \big\la b_1\mathfrak{g}_1, \, \mathbf{q}_{x}ww_y \big\ra + \v \big\la (b_2\mathfrak{g}_2)_x, \, \mathbf{q}_{x} w^2 \big\ra\nonumber\\
&\lesssim  L\|\mathbf{q}_{xx}w_y\|^2 + L[\Phi]_{3,\hat{w}}  + L[\phi]_{3,\hat{w}},
\end{align}
where we have used \eqref{A.9} in \eqref{3.147}-\eqref{3.148-1}.

\smallskip

\noindent{\it Step 3.2.} 
Noting \eqref{3.67-0}, we have 
\begin{align}\label{3.149}
\begin{split}
&\big\la(b_1\,I_{u}(\phi))_y, \, \mathbf{q}_{x}w^2 \big\ra
=\big\la b_1\,I_{u}(\phi), \, (\mathbf{q}_{x}w^2)_y \big\ra 
\lesssim L \|\mathbf{q}_{xx}w_y\|^2 +  L [\Phi]_{3,\hat{w}} + L [\phi]_{3,\hat{w}},\\
&\big\la\big(\v  b_2 I_{v}(\phi)\big)_x,\, \mathbf{q}_{x}w^2\big\ra 
\lesssim  \|\sqrt{\v}\mathbf{q}_xw\|\cdot \sqrt{\v}\|\big(I_{v}(\phi), \pa_xI_{v}(\phi)\big)w\| 
\lesssim L [\Phi]_{3,\hat{w}} + L [\phi]_{3,\hat{w}}.
\end{split}
\end{align}

Using   Lemma \ref{lemA.3},  we have 
\begin{align}
	\big\la (b_1[\tilde{d}_{12}+\cdots] p_y)_y,\, \mathbf{q}_{x}w^2 \big\ra
	&=-\big\la b_1[\tilde{d}_{12}+\cdots]\,p_{y},\, \mathbf{q}_{xy}w^2 \big\ra -  \big\la b_1[\tilde{d}_{12}+\cdots] \,p_y,\, \mathbf{q}_{x}2ww_y \big\ra\nonumber\\
	&\lesssim L\|\mathbf{q}_{xx}w_y\|^2 +   \sqrt{\v} [\Phi]_{3,\hat{w}} + \f{1}{\sqrt{\v}} [[p]]_{2,\hat{w}}, \label{3.150}\\
	\big\la \v (b_2 [\tilde{d}_{21} + \cdots] p_x)_x,\, \mathbf{q}_{x}w^2 \big\ra & =\big\la \v  b_2[\tilde{d}_{21} + \cdots] p_{xx},\, \mathbf{q}_{x}w^2 \big\ra + \big\la \v (b_2[\tilde{d}_{21} + \cdots])_x \, p_x,\, \mathbf{q}_{x}w^2 \big\ra \nonumber\\
	&\lesssim  L [\Phi]_{3,\hat{w}} + L \frac{1}{\v} [[p]]_{2,\hat{w}}, \label{3.151}
\end{align}
and
\begin{align}\label{3.152}
	\big\la \v (b_2 [\tilde{d}_{23}+ \cdots] p)_x,\, \mathbf{q}_{x}w^2 \big\ra
	&\lesssim \|\sqrt{\v}\mathbf{q}_{x}w\| \big(a\sqrt{\v}\|p_xw\| + \sqrt{\v}\|(b_2\tilde{d}_{23})_x\, pw\|\big)\nonumber\\
	&\lesssim L  [\Phi]_{3,\hat{w}} + L^2 \f{1}{\v}[[p]]_{2,\hat{w}} + L \f{1}{\v^2} [[p]]_{2,1}\cdot \mathbf{1}_{\{w=w_0\}} .
\end{align}

Noting \eqref{3.65-1} and  Lemmas \ref{lemA.2} $\&$ \ref{lemA.5-1}, we obtain
\begin{align}\label{3.153}
&\big\la (b_1[\tilde{d}_{13}+ \cdots] p)_y,\, \mathbf{q}_{x}w^2 \big\ra  
=\big\la b_1 [\tilde{d}_{13}+ \cdots]p_y,\, \mathbf{q}_{x}w^2 \big\ra + \big\la (b_1[\tilde{d}_{13}+ \cdots])_y\, p,\, \mathbf{q}_{x}w^2 \big\ra \nonumber\\
&\,\,\lesssim \|\f{\mathbf{q}_{x}}{y}w\|\big\{\|p_yw\| +  \|\mathbf{d}_2pw\| \big\} \lesssim L [\Phi]_{3,\hat{w}} +  L \f1{\v}[[p]]_{2,w}  + L \f1{\v^2}[[p]]_{2,1}  \mathbf{1}_{\{w=w_0\}}.
\end{align}
 
Integrating by parts, one can obtain
\begin{align}\label{3.154}
-\big\la (b_1\hat{d}_{12} S_{y})_y ,\, \mathbf{q}_{x}w^2 \big\ra & =\big\la b_1\hat{d}_{12} S_{y}  ,\, \mathbf{q}_{xy}w^2 \big\ra + \big\la b_1\hat{d}_{12} S_{y}  ,\, \mathbf{q}_{x}2ww_y \big\ra\nonumber\\
&\lesssim L \|\mathbf{q}_{xx}w_y\|^2 + L [\Phi]_{3,\hat{w}}  + L [[[S]]]_{2,w},
\end{align}
\begin{align}
- \big\la (b_1\hat{d}_{13} S)_y ,\, \mathbf{q}_{x}w^2 \big\ra 
&=  \big\la b_1\hat{d}_{13} S ,\, \mathbf{q}_{xy}w^2 \big\ra +  \big\la b_1\hat{d}_{13} S ,\, \mathbf{q}_{x} 2ww_y \big\ra\nonumber\\
&\lesssim  \big(\|\mathbf{q}_{xy}w\| + L \|\mathbf{q}_{xx}w_y\|\big) \|\mathbf{d}_2S w\| \nonumber\\
&\lesssim  L \|\mathbf{q}_{xx}w_y\|^2 + L[\Phi]_{3,\hat{w}} +  L [[[S]]]_{2,1} +  L \f{1}{\v} [[[S]]]_{2,1}  \mathbf{1}_{\{w=w_0\}},
\end{align}
and
\begin{align}\label{3.155-1}
	&\big\la \v (b_2\hat{d}_{21} S_x)_x +\v (b_2\hat{d}_{22} S_y)_x + \v (b_2\hat{d}_{23} S)_x,\, \mathbf{q}_{x}w^2 \big\ra \nonumber\\
	&\lesssim L  [\Phi]_{3,\hat{w}} + L [[[S]]]_{2,\hat{w}}  + L\f{1}{\v} [[[S]]]_{2,1}  \mathbf{1}_{\{w=w_0\}}.
\end{align}

For terms involving $q$, we have 
\begin{align}\label{3.156-1} 
&\v \Big|\Big\la \Big(b_2 \Big[\v \hat{d}_{21}\, (\f{T_{sy}}{p_s} q)_x  +  \v \hat{d}_{22}\, (\f{T_{sy}}{p_s} q)_y +  \v \hat{d}_{23}\, \f{T_{sy}}{p_s} q\Big] \Big)_x,\, \mathbf{q}_{x}w^2 \Big\ra\Big|\nonumber\\
&+ \Big|\Big\la \Big(b_1 \Big[\hat{d}_{12}\big(\f{T_{sy}}{p_s} q\big)_y + \hat{d}_{13}\f{T_{sy}}{p_s} q\Big] \Big)_y,\, \mathbf{q}_{x}w^2 \Big\ra\Big| \lesssim L\|\mathbf{q}_{xx}w_y\|^2 + L [\Phi]_{3,\hat{w}}.
\end{align}

Recalling \eqref{3.66-0}-\eqref{3.66-1} and combining \eqref{3.149}-\eqref{3.156-1}, one obtains 
\begin{align}\label{3.158}
\big\la \mathcal{F}_{R3},\, \mathbf{q}_{x}w^2\big\ra 
&\lesssim  L \|\mathbf{q}_{xx}w_y\|^2  + L [\Phi]_{3,\hat{w}} +  L \f{1}{\v} [[p]]_{2,\hat{w}}  + L [\phi]_{3,\hat{w}} +  L  [[[S]]]_{2,\hat{w}}  \nonumber\\
&\quad+ L \f{1}{\v}\Big\{[[[S]]]_{2,1} + \f{1}{\v} [[p]]_{2,1} \Big\} \mathbf{1}_{\{w=w_0\}} + \v^{2N_0-1}.
\end{align}

\smallskip

\noindent{\it Step 3.3.} Recall  $\mathcal{F}_{R4}$ in \eqref{3.92-1}. 
Integrating by parts in $y$, then using \eqref{3.93} and \eqref{6.176},  one obtains that 
\begin{align}\label{3.159-1}
&\Big\la  \Big(\frac{b_1 }{T^{\v}} u^{\v} \Delta_{\v}\tilde{T}\Big)_y,\, \mathbf{q}_x w^2 \Big\ra 
= -  \Big\la \frac{b_1 }{T^{\v}} u^{\v} \Delta_{\v}\tilde{T} ,\, \mathbf{q}_{xy} w^2 \Big\ra - \Big\la  \frac{b_1 }{T^{\v}} u^{\v} \Delta_{\v}\tilde{T} ,\, \mathbf{q}_x 2ww_y \Big\ra\nonumber\\
&\lesssim \f1N \|u_s\mathbf{q}_{xy}w\|^2 + L\|\mathbf{q}_{xx}w_y\|^2 +   L\f{1}{\v}[[p]]_{2,\hat{w}}  + L^{1-} [\phi]_{3,1}  +  L^{\f12}[[[S]]]_{2,w}.
\end{align}
Using \eqref{3.137-1}, one gets
\begin{align}\label{3.160-1}
\Big\la  \lambda\v\frac{b_1-b_2}{T^{\v}} v^{\v} \tilde{T}_{xyy},\, \mathbf{q}_x w^2 \Big\ra 
&\lesssim L\|\sqrt{\v}\mathbf{q}_{xx} w\|\cdot \|u_s\sqrt{\v}\tilde{T}_{xyy}w\| \nonumber\\
&\lesssim L [\Phi]_{3,\hat{w}} + L [[[\zeta]]] _{3,w}  + L [[p]]_{3,w} + L\f{1}{\v}[[p]]_{2,\hat{w}} + L[\phi]_{3,\hat{w}}.
\end{align}
It follows from  \eqref{3.93}-\eqref{3.94} and \eqref{6.176} that
\begin{align}\label{3.161-2}
-\Big\la  \v \Big(\frac{b_2v^{\v}}{T^{\v}}  \Delta_{\v}\tilde{T}\Big)_x ,\, \mathbf{q}_{x}w^2\Big\ra 
& \lesssim L[\Phi]_{3,\hat{w}} + \v [[[\zeta]]]_{2,w}   + L[[p]]_{2,\hat{w}}  + \v [\phi]_{3,1} + L\|u_s \sqrt{\v} \Delta_{\v}\tilde{T}_xw\|^2,
\end{align}
and
\begin{align}\label{3.161-1}
\v \Big\la  \big(\frac{b_1v^{\v}}{T^{\v}} \big)_y \tilde{T}_{xy} -  \big(\frac{b_2v^{\v}}{T^{\v}}\big)_x \tilde{T}_{yy},\, \mathbf{q}_x w^2 \Big\ra 
&\lesssim  \sqrt{\v} [\Phi]_{3,\hat{w}} + \sqrt{\v}[[[\zeta]]]_{2,w} + L\frac{1}{\v} [[p]]_{2,\hat{w}}  + \sqrt{\v}[\phi]_{3,1}.
\end{align}

Combining \eqref{3.92-1} and \eqref{3.159-1}-\eqref{3.161-1}, one has that 
\begin{align}\label{3.162-1}
\big\la \mathcal{F}_{R4},\, \mathbf{q}_{x}w^2\big\ra 
&\lesssim \f1N \|u_s\mathbf{q}_{xy}w\|^2 + L\|\mathbf{q}_{xx}w_y\|^2 + L [\Phi]_{3,\hat{w}} +   L\f{1}{\v}[[p]]_{2,\hat{w}}  + L [[p]]_{3,w}  \nonumber\\
&\quad +  L^{\f12}[[[S]]]_{2,w}  + L [[[\zeta]]] _{3,w}   + L^{1-}[\phi]_{3,\hat{w}} + L\|u_s \sqrt{\v} \Delta_{\v}\tilde{T}_xw\|^2.
\end{align}

\

\noindent{\it Step 3.4.} We consider $\big\la \mathcal{F}_{R5},\, \mathbf{q}_{x}w^2\big\ra$. Integrating by parts in $y$, and using  Lemmas \ref{lemA.3} \& \ref{lemA.5-1}, one gets
\begin{align}\label{3.192-1}
	\Big\la \Big(\f{b_1}{T^{\v}} \big(\f{1}{2\rho_s}\big)_{yy} u^{\v} p\Big)_y,\, \mathbf{q}_{x}w^2 \Big\ra
	&=-\Big\la \f{b_1}{T^{\v}} \big(\f{1}{2\rho_s}\big)_{yy} u^{\v} p,\, \mathbf{q}_{xy}w^2 \Big\ra + \Big\la  \f{b_1}{T^{\v}} \big(\f{1}{2\rho_s}\big)_{yy} u^{\v} p ,\, \mathbf{q}_{x} 2ww_y \Big\ra\nonumber\\
	&\lesssim L \|\mathbf{q}_{xx}w_y\|^2 + L [\Phi]_{3,\hat{w}} +  L \f1{\v}[[p]]_{2,w}  + L \f1{\v^2}[[p]]_{2,1}  \mathbf{1}_{\{w=w_0\}},
\end{align} 
\begin{align}
\Big\la \Big(\f{b_1}{T^{\v}} \big(\f{1}{\rho_s}\big)_{y} u^{\v} p_y\Big)_y,\, \mathbf{q}_{x}w^2 \Big\ra
&=- \Big\la \f{b_1}{T^{\v}}\big(\f{1}{\rho_s}\big)_{y} u^{\v} p_y,\, \mathbf{q}_{xy}w^2 \Big\ra - \Big\la \f{b_1}{T^{\v}} \big(\f{1}{\rho_s}\big)_{y} u^{\v} p_y,\, \mathbf{q}_{x} 2ww_y \Big\ra \nonumber\\
& \lesssim \sqrt{\v} \|\mathbf{q}_{xx}w_y\|^2 + \sqrt{\v} [\Phi]_{3,\hat{w}} +   \f1{\sqrt{\v}}[[p]]_{2,w},
\end{align}
and
\begin{align}
\Big\la \Big(\f{b_1u^{\v}}{T^{\v}} \big(\f{T_{sy}}{p_s}q\big)_{yy}\Big)_y,\, \mathbf{q}_{x}w^2 \Big\ra 
&=-\Big\la \f{b_1u^{\v}}{T^{\v}} \big(\f{T_{sy}}{p_s}q\big)_{yy} ,\, \mathbf{q}_{xy}w^2 \Big\ra - \Big\la \f{b_1u^{\v}}{T^{\v}} \big(\f{T_{sy}}{p_s}q\big)_{yy} ,\, \mathbf{q}_{x} 2ww_y \Big\ra \nonumber\\
&\lesssim \f1N \|u_s\mathbf{q}_{xy}w\|^2 + L \|\mathbf{q}_{xx}w_y\|^2 + L^{1-} [\phi]_{3,\hat{w}}.
\end{align}
Also, using  Lemmas \ref{lemA.3} \& \ref{lemA.5-1}, one can obtain
\begin{align}\label{3.195-1}
	\big\la \v \big(b_2G_{22}\big)_x, \mathbf{q}_{x}w^2\big\ra
	&\lesssim L [\Phi]_{3,\hat{w}} + L [\phi]_{3,\hat{w}} + L\f{1}{\v}[[p]]_{2,\hat{w}}.
\end{align}

Noting \eqref{2.98} and \eqref{3.92-1}, we have from  \eqref{3.192-1}-\eqref{3.195-1} that
\begin{align} 
	\big\la \mathcal{F}_{R5},\, \mathbf{q}_{x}w^2\big\ra
	&\lesssim \f1N \|u_s\mathbf{q}_{xy}w\|^2 + L \|\mathbf{q}_{xx}w_y\|^2  + L [\Phi]_{3,\hat{w}} + L^{1-} [\phi]_{3,\hat{w}} \nonumber\\
	&\quad +  L \f1{\v}[[p]]_{2,\hat{w}}  + L \f1{\v^2}[[p]]_{2,1} \mathbf{1}_{\{w=w_0\}}.
\end{align}

\

\noindent{\it Step 3.5.} We consider  $\big\la \mathcal{F}_{R6},\, \mathbf{q}_{x}w^2\big\ra$. Integrating by parts, then using \eqref{3.93}-\eqref{3.94},  one has
\begin{align}\label{3.199}
 \Big\la \big(\f{b_1}{T^{\v}}u^{\v} \hat{\mathbf{S}}_{yy}\big)_y ,\, \mathbf{q}_{x}w^2\Big\ra 
	&=- \Big\la \f{b_1}{T^{\v}}u^{\v} \hat{\mathbf{S}}_{yy} ,\, \mathbf{q}_{xy}w^2\Big\ra -  \Big\la \f{b_1}{T^{\v}}u^{\v} \hat{\mathbf{S}}_{yy} ,\, \mathbf{q}_{x} 2ww_y\Big\ra\nonumber\\
	&\lesssim \f1N \|u_s\mathbf{q}_{xy}w\|^2 + L\|\mathbf{q}_{xx}w_y\|^2  + L^{\f12-} [[[\hat{\mathbf{S}}]]]_{2,w},
\end{align}
and
\begin{align}\label{3.199-1}
  \v \Big\la \big(\f{b_1u^{\v}}{T^{\v}} \hat{\mathbf{S}}_{xy}\big)_x ,\, \mathbf{q}_{x}w^2\Big\ra  
	&=- \v \Big\la \big(\f{b_1u^{\v}}{T^{\v}} \hat{\mathbf{S}}_{xx}  ,\, \mathbf{q}_{xy}w^2\Big\ra  - \v \Big\la \big(\f{b_1u^{\v}}{T^{\v}} \hat{\mathbf{S}}_{xx}  ,\, \mathbf{q}_{x}2ww_y\Big\ra\nonumber\\
	&\quad +  \v \Big\la \big(\big(\f{b_1u^{\v}}{T^{\v}}\big)_x  \hat{\mathbf{S}}_{xy},\, \mathbf{q}_{x}w^2\Big\ra - \v \Big\la \big(\big(\f{b_1u^{\v}}{T^{\v}}\big)_y \hat{\mathbf{S}}_{xx}  ,\, \mathbf{q}_{x}w^2\Big\ra \nonumber\\
	&\lesssim  \f{1}{N} \|u_s\mathbf{q}_{xy}w\|^2 + \sqrt{\v} \|\mathbf{q}_{xx}w_y\|^2 +\sqrt{\v} [\Phi]_{3,\hat{w}} + \sqrt{\v} [[[\hat{\mathbf{S}}]]]_{2,\hat{w}}.
\end{align}
Hence it follows from $\eqref{2.98-1}_1$ and \eqref{3.199}-\eqref{3.199-1} that
\begin{align}\label{3.199-2}
	\big\la \mathcal{F}_{R6}(\hat{\mathbf{S}}),\, \mathbf{q}_{x}w^2\big\ra
	 &\lesssim  \f{1}{N} \|u_s\mathbf{q}_{xy}w\|^2 + L\|\mathbf{q}_{xx}w_y\|^2 + \sqrt{\v}[\Phi]_{3,w} + L^{\f12-} [[[\hat{\mathbf{S}}]]]_{2,w}.
\end{align}

\

\noindent{\it Step 3.6.} We consider $\big\la \mathcal{F}_{R7},\, \mathbf{q}_{x}w^2\big\ra$. Integrating by parts in $x$, then using \eqref{3.93}, \eqref{3.122} and Lemma \ref{lemA.3}, one gets
\begin{align}\label{3.200}
 \Big\la \v \Big(\frac{b_2v^{\v}}{p^{\v}}  \Delta_{\v}\hat{P}\Big)_x,\, \mathbf{q}_{x}w^2\Big\ra
	&\lesssim \sqrt{\v}  \|\sqrt{\v} \mathbf{q}_{xx}w\|\cdot \|\Delta_{\v}\hat{P}w\| +  \sqrt{\v} \| \mathbf{q}_{x}w\|_{x=L}\cdot \sqrt{\v}\|u_s\Delta_{\v}\hat{P} w\|_{x=L} \nonumber\\
	&\quad  + \|\sqrt{\v}\mathbf{q}_{x}w\|_{x=0}\cdot \sqrt{\v}\|u_s \Delta_{\v}\hat{P} w\|_{x=0} \nonumber\\
	&\lesssim \v^{\f12-} [\Phi]_{3,\hat{w}} +  \f{1}{\sqrt{\v}}\mathcal{K}_{2,w}(\hat{P}) 
	+ \sqrt{\v}\bar{\mathcal{K}}_{3,w}(\hat{P}).
\end{align}
and
\begin{align}
	&\Big\la \Big\{- \mu \frac{b_1}{p^{\v}} u^{\v} \Delta_{\v}\hat{P}  - \lambda \frac{b_1\rho^{\v}}{2\rho_sp^{\v}} u^{\v} \hat{P}_{yy}  + \lambda\v \frac{b_2-b_1}{p^{\v}} (U^{\v}\cdot\nabla)\hat{P}_x -  \lambda \v \frac{b_2\rho^{\v}}{2\rho_sp^{\v}} u^{\v} \hat{P}_{xx}\Big\}_y ,\, \mathbf{q}_{x}w^2\Big\ra \nonumber\\
	&= \Big\la \mu \frac{b_1}{p^{\v}} u^{\v} \Delta_{\v}\hat{P} + \lambda \frac{b_1\rho^{\v}}{2\rho_sp^{\v}} u^{\v} \hat{P}_{yy} - \lambda\v \frac{b_2-b_1}{p^{\v}} (U^{\v}\cdot\nabla)\hat{P}_x + \lambda \v \frac{b_2\rho^{\v}}{2\rho_sp^{\v}} u^{\v} \hat{P}_{xx},\, (\mathbf{q}_{x}w^2)_y\Big\ra  \nonumber\\
	&\lesssim \sqrt{\v}[\Phi]_{3,w} + \sqrt{\v}\|\mathbf{q}_{xx}w_y\|^2  + \f{1}{\sqrt{\v}} \mathcal{K}_{2,w}(\hat{P}).
\end{align}
It is clear that 
\begin{align}\label{3.202}
	&\Big\la \v \big(\frac{b_2}{p^{\v}} U^{\v} \big)_x  \cdot\nabla\hat{P}_y -  \v \big(\frac{b_2}{p^{\v}} U^{\v}\big)_y \cdot\nabla\hat{P}_x +  \v \big(\frac{b_2\rho^{\v}}{2\rho_s p^{\v}} u^{\v}\big)_x \hat{P}_{xy}  - \v  \Big(\frac{b_2\rho^{\v}}{2\rho_sp^{\v}} u^{\v}\Big)_y \hat{P}_{xx},\, \mathbf{q}_{x}w^2\Big\ra\nonumber\\
	& \lesssim \sqrt{\v}[\Phi]_{3,\hat{w}} +  \f{1}{\sqrt{\v}} \mathcal{K}_{2,w}(\hat{P}).
\end{align}
Hence it follows from \eqref{3.98-1} and \eqref{3.200}-\eqref{3.202} that 
\begin{align}\label{3.203-0}
\big\la \mathcal{F}_{R7},\, \mathbf{q}_{x}w^2\big\ra \lesssim \sqrt{\v} [\Phi]_{3,w} + \sqrt{\v}\|\mathbf{q}_{xx}w_y\|^2 + \f{1}{\sqrt{\v}}\mathcal{K}_{2,w}(\hat{P}) + \sqrt{\v}\bar{\mathcal{K}}_{3,w}(\hat{P}). 
\end{align}

\noindent{\it Step 3.7.}  Substituting \eqref{3.140}-\eqref{3.129-1} and \eqref{3.146-1} into \eqref{3.139}, and using the estimates in steps 3.1-3.6, we conclude \eqref{3.138-0}.
Therefore the proof of Lemma \ref{lem3.3} is completed. $\hfill\Box$

\medskip

\subsection{Estimate on $\|||\Phi||\|_w$}\label{sec5}
Unlike the weighted norm  $\|\sqrt{u_s}\phi_{xyyy}\|$ estimated in Guo-Iyer \cite{Guo-Iyer-2024},  we need to establish  estimate on $\|\nabla_{\v}^4\Phi\|$ which is important for us to control the pressure. As explained in the introduction, the subtle $H^{1/2}_{00}$ argument plays a key role in this subsection.

\begin{lemma}\label{lem3.4}
It holds that 
\begin{align}\label{3.Tr142}
	&\|\sqrt{\v}\sqrt{u_s}\Phi_{xyy}w\|_{x=0}^2 + \|\v \sqrt{u_s}\Phi_{xxy}w\|_{x=L}^2  + \|\big(\sqrt{\v}\Phi_{xyyy}, \v \Phi_{xxyy}, \v^{\f32}\Phi_{xxxy}\big)w\|^2 \nonumber\\
	&\lesssim \|\Phi_{yyyy}\|^2  +  L^{\f12\alpha-} \big\{ [\Phi]_{3,\hat{w}} + [\Phi]_{4,w} \big\}  + \|\mathcal{F}_{R}w\|^2.
\end{align}
\end{lemma}

\noindent{\bf Proof.} Multiplying \eqref{2.121} by $\v \Phi_{xxyy}w^2$, one has 
\begin{align} \label{3.Tr143}
	&\big\la\big(b_1 m u_s^2 \mathbf{q}_{xy}\big)_y,\, \v \Phi_{xxyy}w^2\big\ra + \big\la \v \big(b_2 mu_s^2 \mathbf{q}_{xx}\big)_x ,\, \v \Phi_{xxyy} w^2\big\ra  \nonumber\\
	&-\mu \Big\la\Big(\frac{b_1}{\rho_s} \Delta_{\v}\Phi_y\Big)_y,\, \v \Phi_{xxyy} w^2\Big\ra- \mu \Big\la\v \Big(\frac{b_2}{\rho_s}   \Delta_{\v}\Phi_x\Big)_x,\, \v \Phi_{xxyy} w^2\Big\ra \nonumber\\
	&= \big\la \mathcal{F}_{R}(\hat{P}, \hat{\mathbf{S}}, \zeta, p, \phi, S), \, \v \Phi_{xxyy} w^2\big\ra  + O(1)   \v^{\f12 N_0} \big\{[\Phi]_{3,\hat{w}} + [\Phi]_{4,w}\big\},
\end{align}
where we have used similar arguments as in \eqref{3.18-1}. Since the proof is long, we divide it into several steps. 

\smallskip

\noindent{\it Step 1. Estimate on Rayleigh term.} First we rewrite the Rayleigh terms as
\begin{align}\label{5.3}
\begin{split}
\big(b_1m u_s^2 \mathbf{q}_{xy}\big)_y  
&= b_1  u_s \big\{\Phi_{xyy}-\bar{u}_{sxyy} \mathbf{q} - \bar{u}_{syy}\mathbf{q}_{x} - \bar{u}_{sx} \mathbf{q}_{yy} - 2 \bar{u}_{sxy} \mathbf{q}_y - 2\bar{u}_{sy} \mathbf{q}_{xy}\big\} \\
&\quad + \big(b_1m u_s^2\big)_y \mathbf{q}_{xy},\\
\v \big(b_2mu_s^2 \mathbf{q}_{xx}\big)_x
&=\v b_2u_s \big\{\Phi_{xxx} - 3\bar{u}_{sx}\mathbf{q}_{xx} - 3 \bar{u}_{sxx} \mathbf{q}_x - \bar{u}_{sxxx} \mathbf{q}\big\} + \v \big(b_2mu_s^2\big)_x \mathbf{q}_{xx}.
\end{split}
\end{align} 
Then we have from $\eqref{5.3}_1$ that
\begin{align}\label{3.Tr144}
 \big\la\big(b_1\,m u_s^2 \mathbf{q}_{xy}\big)_y,\, \v \Phi_{xxyy}w^2\big\ra 
&=\big\la b_1 u_s \Phi_{xyy} ,\, \v  \Phi_{xxyy}w^2\big\ra
  + C L^{\f12\alpha-} [\Phi]_{3,\hat{w}}  + C L^{\f12\alpha-} [\Phi]_{4,w}\nonumber\\
&= -\f12 \|\sqrt{b_1u_s}  \sqrt{\v}\Phi_{xyy}w\|_{x=0}^2   + C L^{\f12\alpha-} [\Phi]_{3,\hat{w}}  + C L^{\f12\alpha-} [\Phi]_{4,w}.
\end{align}

Integrating by parts in $y$ and then in $x$, we have
\begin{align}\label{3.Tr149}
 \big\la \v \big(b_2\, m u_s^2 \mathbf{q}_{xx}\big)_x,\, \v \Phi_{xxyy}w^2\big\ra 
	&=\big\la \v b_2u_s \Phi_{xxx} ,\, \v \Phi_{xxyy}w^2\big\ra + O(1) \sqrt{\v} [\Phi]_{3,\hat{w}}  + O(1) \sqrt{\v} [\Phi]_{4,w}\nonumber\\
	&= - \big\la \v  b_2 u_s \Phi_{xxxy} ,\, \v  \Phi_{xxy}w^2\big\ra - \big\la \v  (b_2u_s w^2)_y \Phi_{xxx} ,\, \v  \Phi_{xxy}\big\ra\nonumber\\
	&\quad   + O(1) \sqrt{\v} [\Phi]_{3,\hat{w}}  + O(1) \sqrt{\v} [\Phi]_{4,w}\nonumber\\
	&=  -\f12 \|\sqrt{b_2 u_s}\v\Phi_{xxy}w\|_{x=L}^2  + O(1)\sqrt{\v} \big\{ [\Phi]_{3,\hat{w}}  +  [\Phi]_{4,w}\big\}.
\end{align}

\

{\it Step 2. Estimate on viscous term.} It is obvious that 
\begin{align}\label{3.Tr150}
	&-\Big\la\Big(\frac{b_1}{\rho_s} \Delta_{\v}\Phi_y\Big)_y + \v \Big(\frac{b_2}{\rho_s} \Delta_{\v}\Phi_x\Big)_x,\, \v \Phi_{xxyy}w^2\Big\ra\nonumber\\
	&=- \|\sqrt{\frac{b_1+b_2}{\rho_s}} \v \Phi_{xxyy}w\|^2 -\Big\la\Big(\frac{b_1}{\rho_s} \Phi_{yyy}\Big)_y,\, \v  \Phi_{xxyy}w^2\Big\ra  - \Big\la \v^2 \Big(\frac{b_2}{\rho_s} \Phi_{xxx}\Big)_x,\, \v \Phi_{xxyy}w^2\Big\ra\nonumber\\
	&\quad + O(1)\sqrt{\v} \big\{[\Phi]_{3,\hat{w}} + [\Phi]_{4,w} \big\},
\end{align}
where we have used the following estimate
\begin{align*}
 \v \Big|\Big\la \big(\frac{b_2}{\rho_s}\big)_x \Phi_{xyy} +  \big(\frac{b_1}{\rho_s}\big)_y \Phi_{xxy} ,\, \v \Phi_{xxyy}w^2\Big\ra \Big|  
&\lesssim \sqrt{\v} [\Phi]_{3,\hat{w}} + \sqrt{\v} [\Phi]_{4,w}.
\end{align*}
Then we need only to estimate 2th \& 3th terms on RHS of \eqref{3.Tr150}.

\smallskip

{\it Step 2.1.} Integrating by parts in $x$ and then $y$, one has
\begin{align}\label{3.Tr152}
 -\Big\la\Big(\frac{b_1}{\rho_s} \Phi_{yyy}\Big)_y,\, \v \Phi_{xxyy}w^2\Big\ra 
	&=- \Big\la \frac{b_1}{\rho_s} \Phi_{yyyy} ,\, \v  \Phi_{xxyy}w^2\Big\ra + O(1) L^{\f18} \big\{ [\Phi]_{3,\hat{w}} + \|||\Phi||\|^2_{w}\big\}\nonumber\\
	&=\Big\la \frac{b_1}{\rho_s} \Phi_{xyyyy} ,\, \v \Phi_{xyy}w^2\Big\ra + O(1) L^{\f18} \big\{ [\Phi]_{3,\hat{w}} + [\Phi]_{4,w}\big\} \nonumber\\
	&=- \|\sqrt{\f{b_1}{\rho_s}} \sqrt{\v}\Phi_{xyyy}w\|^2  - \Big\la \frac{b_1}{\rho_s}  \Phi_{xyyy} ,\, \v   \Phi_{xyy}\Big\ra_{y=0} \nonumber\\
	&\quad + O(1) L^{\f18} \big\{ [\Phi]_{3,\hat{w}} + [\Phi]_{4,w}\big\}.
\end{align}

For the boundary term on RHS of \eqref{3.Tr152}, we need to use the $H^{1/2}_{00}$ argument. Applying  Lemma \ref{lemA.9}, it is clear to have
\begin{align}\label{3.Tr154-1}
&\v \int_0^L \frac{b_1}{\rho_s}  \Phi_{xyyy}  \Phi_{xyy} dx \Big|_{y=0}  
= \v \int_0^L \frac{b_1}{\rho_s}  \Phi_{xyyy}  \Phi_{xyy} \chi(\frac{4x}{L}) dx \Big|_{y=0} +  \v \int_0^L \frac{b_1}{\rho_s} \bar{\chi}(\frac{4x}{L}) \Phi_{xyyy}  \Phi_{xyy}  dx \Big|_{y=0} \nonumber\\
&=\v \int_0^L   \Big(\chi(\f{2x}{L})\Phi_{yyy}\Big)_x  \Phi_{xyy} \frac{b_1}{\rho_s}\chi(\frac{4x}{L}) dx \Big|_{y=0} + \int_0^L \Big(\bar{\chi}(\f{8x}{L}) \Phi_{yyy}\Big)_x  \frac{b_1}{\rho_s} \bar{\chi}(\frac{4x}{L}) \Phi_{xyy}  dx \Big|_{y=0}  \nonumber\\
&\lesssim \v \|\frac{b_1}{\rho_s} \chi\, \Phi_{xyy}(\cdot,0)\|_{H^{1/2}}\cdot \|\chi\,\Phi_{yyy}(\cdot,0)\|_{H^{1/2}_{00}} + \v \|\frac{b_1}{\rho_s}\bar{\chi}\,\Phi_{xyy}(\cdot,0)\|_{H^{1/2}_{00}}\cdot \| \bar{\chi}\Phi_{yyy}(\cdot,0)\|_{H^{1/2}}.
\end{align}
In the following we shall control the terms on RHS of \eqref{3.Tr154-1}.

\smallskip

\noindent\underline{$\bullet$ \it Estimate on 1th term on RHS of \eqref{3.Tr154-1}.} 
Noting $\Phi|_{x=0}=0$, we apply Lemma \ref{lemAH.1} to obtain that 
\begin{align}\label{3.Tr19}
\|\chi \Phi_{yyy}(\cdot,0)\|^2_{H^{1/2}_{00}} &\lesssim \big\{L^{-1}\|\Phi_{yyy}\| + \|\Phi_{xyyy}\|\big\} \big\{\|\Phi_{yyy}\| + \|\Phi_{yyyy}\|\big\} \nonumber\\
&\lesssim \|\Phi_{xyyy}\| \big\{ \|\Phi_{yyy}\| +  \|\Phi_{yyyy}\|\big\},
\end{align}
which, together with \eqref{Z.11-2}, yields that 
\begin{align}\label{3.Tr20}
\eqref{3.Tr154-1}_1
&\lesssim \v \big\{\|\Phi_{xxyy} \| + L^{-1}\| \Phi_{xyy} \|\big\}^{\f12}\big\{ \|\Phi_{xyyy}\| +  \|\Phi_{xyy}\|\big\}^{\f12}  \|\Phi_{xyyy}\|^{\f12} \big\{ \|\Phi_{yyy}\| +  \|\Phi_{yyyy}\|\big\}^{\f12} \nonumber\\
&\lesssim \f1N \|(\v \Phi_{xxyy}, \sqrt{\v}\Phi_{xyyy})\|^2 + N \|\Phi_{yyyy}\|^2 + L^{\f14}[\Phi]_{3,1},
\end{align}
where we have used \eqref{A.12-2} and \eqref{A.12-4}.

\

\noindent\underline{$\bullet$ \it Estimate on 2th term on RHS of \eqref{3.Tr154-1}.} Noting $\Phi_{xyy}\big|_{x=L}=0$,  we also have from Lemma \ref{lemAH.1} that 
\begin{align}\label{3.Tr26-1}
\|\frac{b_1}{\rho_s}\bar{\chi}\,\Phi_{xyy}(\cdot,0)\|_{H^{1/2}_{00}} 
&\lesssim \big\{L^{-1} \|\Phi_{xyy}\| + \|\Phi_{xxyy}\| \big\} \big\{\|\Phi_{xyyy}\| +  \|\Phi_{xyy}\| \big\}\nonumber\\
	&\lesssim \|\Phi_{xxyy}\| \big\{\|\Phi_{xyyy}\| +  \|\Phi_{xyy}\| \big\},
\end{align}
which yields that 
\begin{align}\label{3.Tr26}
\eqref{3.Tr154-1}_2
&\lesssim \v \|\Phi_{xxyy}\|^{\f12}\big\{ \|\Phi_{xyyy}\| + \|\Phi_{xyy}\|\big\}^{\f12}  \|\Phi_{xyyy}\|^{\f12} \big\{ \|\Phi_{yyy}\| +  \|\Phi_{yyyy}\|\big\}^{\f12} \nonumber\\
&\lesssim \f1N \|(\v \Phi_{xxyy}, \sqrt{\v}\Phi_{xyyy})\|^2 + N \|\Phi_{yyyy}\|^2 + L^{\f14} [\Phi]_{3,1}.
\end{align}

\smallskip

$\bullet$ Substituting \eqref{3.Tr20} and \eqref{3.Tr26} into \eqref{3.Tr154-1}, one has
\begin{align}\label{3.Tr27}
\v \int_0^L \frac{b_1}{\rho_s}  \Phi_{xyyy}  \Phi_{xyy} dx \Big|_{y=0}   &\lesssim \f1N \|(\v \Phi_{xxyy}, \sqrt{\v}\Phi_{xyyy})\|^2 + N \|\Phi_{yyyy}\|^2 + L^{\f14} [\Phi]_{3,1},
\end{align}
which, together with  \eqref{3.Tr152}, yields that 
\begin{align}\label{3.155}
-\Big\la\Big(\frac{b_1}{\rho_s} \Phi_{yyy}\Big)_y,\, \v \Phi_{xxyy}w^2\Big\ra 
	&= - (1-\f1N) \|\sqrt{\f{b_1}{\rho_s}}\sqrt{\v}\Phi_{xyyy}w\|^2  + \f1N \|\v \Phi_{xxyy}\|^2  + N \|\Phi_{yyyy}\|^2  \nonumber\\
	&\quad + CL^{\f18} \big\{ [\Phi]_{3,\hat{w}} + [\Phi]_{4,w}\big\}.
\end{align}

\smallskip

{\it Step 2.2.} Integrating by parts in $x$ and $y$, we get
\begin{align}\label{3.156}
 - \Big\la \v^2 \Big(\frac{b_2}{\rho_s} \Phi_{xxx}\Big)_x,\, \v \Phi_{xxyy}w^2\Big\ra  
 &=\Big\la \v^2  \frac{b_2}{\rho_s} \Phi_{xxx} ,\, \v \Phi_{xxxyy}w^2\Big\ra  \nonumber\\
	&=- \|\sqrt{\f{b_2}{\rho_s}} \v^{\f32}\Phi_{xxxy}w\|^2  + C \sqrt{\v} \big\{ [\Phi]_{3,\hat{w}} + [\Phi]_{4,w}\big\}.
\end{align}

\smallskip

{\it Step 2.3.} Substituting \eqref{3.155}-\eqref{3.156}   into \eqref{3.Tr150}, one has
\begin{align}\label{3.157}
	&-\Big\la\Big(\frac{b_1}{\rho_s} \Delta_{\v}\Phi_y\Big)_y + \v \Big(\frac{b_2}{\rho_s} \Delta_{\v}\Phi_x\Big)_x,\, \v \Phi_{xxyy}w^2\Big\ra\nonumber\\
	&= - (1-\f1N)  \bigg\{ \|\sqrt{\f{b_1}{\rho_s}} \sqrt{\v}\Phi_{xyyy}w\|^2+ \|\sqrt{\frac{b_1+b_2}{\rho_s}}  \v \Phi_{xxyy}w\|^2 
	+ \|\sqrt{\f{b_2}{\rho_s}} \v^{\f32}\Phi_{xxxy}w\|^2\bigg\}\nonumber\\
	&\quad + N \|\Phi_{yyyy}\|^2 + CL^{\f18} \big\{ [\Phi]_{3,\hat{w}} + [\Phi]_{4,w}\big\}.
\end{align}

\smallskip

\noindent{\it Step 3.}
Substituting \eqref{3.157} and \eqref{3.Tr144}-\eqref{3.Tr149} into \eqref{3.Tr143}, we  conclude  \eqref{3.Tr142}. Therefore the proof of Lemma \ref{lem3.4} is completed. $\hfill\Box$

\medskip

\begin{lemma}\label{lem3.5}
It holds that  
\begin{align}\label{3.172-0}
&\|\sqrt{u_s}\v^{\f32}\Phi_{xxx}w\|_{x=0}^2   + \|\big(\v\Phi_{xxyy}, \v^{\f32} \Phi_{xxxy}, \v^2\Phi_{xxxx}, \Phi_{yyyy}\big)w\|^2 \nonumber\\
&\lesssim  a^3 [\Phi]_{3,\hat{w}} +  L^{\f18} [\Phi]_{4,w}  + \|\mathcal{F}_{R}w\|^2.
\end{align}
\end{lemma}

\noindent{\bf Proof.} Multiplying  \eqref{2.121} by $\v^2 \Phi_{xxxx}w^2$, one has 
\begin{align}  \label{3.172}
	&\big\la\big(b_1 m u_s^2 \mathbf{q}_{xy}\big)_y,\, \v^2 \Phi_{xxxx}w^2\big\ra + \big\la \v \big(b_2 mu_s^2 \mathbf{q}_{xx}\big)_x ,\, \v^2 \Phi_{xxxx}w^2\big\ra  \nonumber\\
	&-\Big\la\mu\,\Big(\frac{b_1}{\rho_s} \Delta_{\v}\Phi_y\Big)_y,\, \v^2 \Phi_{xxxx}w^2\Big\ra-\Big\la\mu\,\v \Big(\frac{b_2}{\rho_s}   \Delta_{\v}\Phi_x\Big)_x,\, \v^2 \Phi_{xxxx}w^2\Big\ra \nonumber\\
	&= \big\la \mathcal{F}_{R}(\hat{P}, \hat{\mathbf{S}}, \zeta, p, \phi, S), \, \v^2 \Phi_{xxxx}w^2\big\ra +  O(1)   \v^{\f12 N_0} \big\{[\Phi]_{3,\hat{w}} + [\Phi]_{4,w}\big\},
\end{align}
where we have used similar arguments as in \eqref{3.18-1}. 
 
\smallskip
 
 Using \eqref{5.3}, one has
\begin{align}\label{5.32-1}
\big\la \v \big(b_2 mu_s^2 \mathbf{q}_{xx}\big)_x ,\, \v^2 \Phi_{xxxx}w^2\big\ra 
&=\Big\la \v b_2u_s \Phi_{xxx},\, \v^2  \Phi_{xxxx}w^2\Big\ra   + C\sqrt{\v} \big\{[\Phi]_{3,\hat{w}} + [\Phi]_{4,w} \big\} \nonumber\\
&= -\f12 \|\sqrt{b_2}\sqrt{u_s}\v^{\f32}\Phi_{xxx}w\|^2_{x=0}  + C\sqrt{\v} \big\{[\Phi]_{3,\hat{w}} + [\Phi]_{4,w} \big\},\\
\big\la\big(b_1 m u_s^2 \mathbf{q}_{xy}\big)_y,\, \v^2 \Phi_{xxxx}w^2\big\ra
&= \f{1}{N}\|\v^2 \Phi_{xxxx}w\|^2 +  N a^3 [\Phi]_{3,\hat{w}}.\label{3.175}
\end{align}
A direct calculation shows that 
\begin{align}\label{3.203}
	&- \Big\la\Big(\frac{b_1}{\rho_s} \Delta_{\v}\Phi_y\Big)_y + \v \Big(\frac{b_2}{\rho_s} \Delta_{\v}\Phi_x\Big)_x,\, \v^2 \Phi_{xxxx}w^2\Big\ra\nonumber\\
	&=-  \|\sqrt{\frac{b_2}{\rho_s}}\v^2 \Phi_{xxxx}w\|^2 + \Big\la \frac{b_1+b_2}{\rho_s} \v \Phi_{xxy} ,\, \v^2  \Phi_{xxxxy}w^2\Big\ra +   \Big\la \frac{b_1}{\rho_s} \Phi_{yyy} ,\, \v^2 \Phi_{xxxxy}w^2\Big\ra  \nonumber\\
	&\quad   + CL^{\f18}\big\{ [\Phi]_{3,\hat{w}} + [\Phi]_{4,w}\big\}\nonumber\\
	&=-  \|\sqrt{\frac{b_2}{\rho_s}}\v^2 \Phi_{xxxx}w\|^2  - \|\sqrt{\frac{b_1+b_2}{\rho_s}}\v^{\f32} \Phi_{xxxy}w\|^2 - \|\sqrt{\frac{b_1}{\rho_s}}\v  \Phi_{xxyy}w\|^2 \nonumber\\
	&\quad  + CL^{\f18}\big\{ [\Phi]_{3,\hat{w}} + [\Phi]_{4,w}\big\}
\end{align}
Substituting \eqref{3.203} and \eqref{5.32-1}-\eqref{3.175} into \eqref{3.172}, we obtain
\begin{align}\label{3.251}
\|\sqrt{u_s}\v^{\f32}\Phi_{xxx}w\|_{x=0}^2   + \|\big(\v\Phi_{xxyy}, \v^{\f32} \Phi_{xxxy}, \v^2\Phi_{xxxx}\big)w\|^2
&\lesssim a^3 [\Phi]_{3,\hat{w}} +  L^{\f18} [\Phi]_{4,w}  + \|\mathcal{F}_{R}w\|^2.
\end{align}

\smallskip

A direct calculation shows that 
\begin{align*}
	\|\big(b_1 m u_s^2 \mathbf{q}_{xy}\big)_yw\|^2 + \|\v \big(b_2 mu_s^2 \mathbf{q}_{xx}\big)_xw\|^2 \lesssim  a^3 [\Phi]_{3,\hat{w}},
\end{align*}
which, together with  \eqref{2.121}  and \eqref{3.251}, yields that
\begin{align}\label{5.41}
\|\Phi_{yyyy}w\|^2
&\lesssim \|(\v\Phi_{xxyy}, \v^2 \Phi_{xxxx})w\|^2 + a^3 [\Phi]_{3, \hat{w}} + \|\mathcal{F}_{R}w\|^2  \nonumber\\
&\lesssim a^3 [\Phi]_{3,\hat{w}} +  L^{\f18} [\Phi]_{4,w}  + \|\mathcal{F}_{R}w\|^2.
\end{align}

Finally, combining  \eqref{3.251}-\eqref{5.41}, one concludes \eqref{3.172-0}.
 Therefore the proof of Lemma \ref{lem3.5} is completed. $\hfill\Box$
 
 \medskip

\subsection{Proof of Theorem \ref{thm4.1}} 
 Combining Lemmas \ref{lem3.1}, \ref{lem3.2} $\&$ \ref{lem3.3}, we have 
\begin{align}\label{3.202-0}
[\Phi]_{3,w}
&\lesssim   \f1N \big\{\|\mathbf{q}_{xx}w_y\|^2  + [\Phi]_{3,1}\big\}   +  L^{\f14-} [\Phi]_{4,w} + (\f{1}{N} + N a ) \|(\phi, p, S(\zeta))\hat{w}\|^2_{\mathbf{X}^{\v}} \nonumber\\
&\quad  + L^{\f14-} \f{1}{\v}\|(0, p, S)\|^2_{\mathbf{X}^{\v}} \cdot \mathbf{1}_{\{w=w_0\}}  + Na  \|(0, \hat{P}, \hat{\mathbf{S}})\hat{w}\|^2_{\mathbf{Y}^{\v}} +  L^{\f12\alpha-} \|u_s\sqrt{\v}\Delta_{\v}\tilde{T}_{x}w\|^2 \nonumber\\
&\quad  + \f{1}{\v} \|(\sqrt{u_s}\mathfrak{N}_3, \mathcal{F}_{R8}) w\|^2  +  \v^{N_0 },
\end{align}
which immediately concludes \eqref{18.3}.

\smallskip
 
Combining  \eqref{3.172-0}, \eqref{3.Tr142} and  \eqref{3.103-0}, we  obtain
\begin{align}\label{3.Tr87}
[\Phi]_{4,w}
&\lesssim a^3 [\Phi]_{3,\hat{w}} +  [[p]]_{3,w} + [[[\zeta]]] _{3,w}  + \mathcal{K}_{3,w}(\hat{P})  + [[[\hat{\mathbf{S}}]]]_{3,w}  +  a^2  \|(\phi, p, S(\zeta)) \hat{w}\|^2_{\mathbf{X}^{\v}}  \nonumber\\
&\quad +  L^{\f12}  \|(0, \hat{P}, \hat{\mathbf{S}})\hat{w}\|^2_{\mathbf{Y}^{\v}} + \v^{1-} \|u_s\sqrt{\v} \Delta_{\v}\tilde{T}_x w\|^2 + \|\mathcal{F}_{R8}w\|^2,
\end{align}
which immediately concludes \eqref{18.4}. Therefore the proof of Theorem \ref{thm4.1} is completed. $\hfill\Box$

\medskip

\section{Existence and Uniform Estimate for   Pressure} \label{sec7}
In  this section, we aim to establish the  existence and uniform estimates for pressure $P$. 
\begin{theorem}\label{thm5.1}
	There exists a unique strong solution $P$ to the following boundary value problem
	\begin{align}\label{18.39-1}
		\begin{cases}
			\dis  \Delta_{\v}P +  \v  (\mu+\lambda)\f{1}{p^{\v}}(U^{\v}\cdot\nabla)\Delta_{\v}P  = \mathcal{G}(\hat{P}, \hat{\mathbf{S}},\Phi, p,S, \zeta, \phi),\quad (x,y)\in (0,L)\times \R_+,\\
			\dis (\mu+\lambda) \f{1}{p_s} u_s\Delta_{\v}P \big|_{x=0}=-\sigma \mathbf{d}_{11} B,
		\end{cases}
	\end{align}
	with $P(0,0)=0,\,\, \lim_{y\to+\infty}P_y=0$, and also  the following boundary conditions
	\begin{align} \label{18.39-2}
		\begin{cases}
			\dis B\big|_{y=0} = \f{\mu}{\rho^{\v}} \Phi_{yyy},\\ 
			\dis \big\{ (\sigma-1)  p_s \mathbf{d}_{11} B - \f12 \v \lambda(1+\bm{\chi}) u_s  B_x\big\} \big|_{x=0} 
			= - \f12 \lambda(1+\bm{\chi}) u_s \Delta_{\v}P +   \mathbf{g}_{b1}(\hat{P},\Phi,p,S, \zeta,\phi),\\
			\dis \big\{\v B_{x} - \lambda\v\frac{b_2}{2p_s} u_s B_{yy}\big\} \big|_{x=L} = \Delta_{\v}P  + \mathbf{g}_{b2}(\hat{P}, \hat{\mathbf{S}}, \Phi, p, S, \zeta, \phi),
		\end{cases}
	\end{align}
	where $B=P_x + \bar{g}_1$, and the definitions of $\mathbf{g}_{b1}, \mathbf{g}_{b2}$ are given in \eqref{17.19-7}-\eqref{17.19-8}.

	Moreover the following estimates hold:
	\begin{align}
		&\mathscr{K}(P)+  \f{1}{\v}\mathcal{K}_{2,\hat{w}_0^{\v}}(P) + \|\f{B\hat{w}^{\v}_0}{\sqrt{u_s}}\|_{x=0}^2  \lesssim  a [\Phi]_{3,\hat{w}_0^{\v}}  + \|||\Phi||\|^2_1 + a^2 \|(0, \hat{P}, \hat{\mathbf{S}})\hat{w}_0^{\v}\|^2_{\mathbf{Y}^{\v}}   + \{\cdots\}_{1,\hat{w}_0^{\v}} \nonumber\\
		&\qquad\qquad\qquad\qquad\qquad\qquad\qquad\quad\,\, + \|(\mathbf{d}_2 \bar{g}_1, \bar{g}_1\mathbf{1}_{\{y\in[0,2]\}})\|^2+   \v \|\sqrt{u_s} \hat{\mathbf{S}}_{xy}\hat{w}_0^{\v}\|^2_{x=L}\mathbf{1_{D}} ,\label{18.41}\\
		& \mathcal{K}_{3,\hat{w}_0^{\v}}(P)
		\lesssim (\xi a + C_{\xi}a^3) [\Phi]_{3,\hat{w}_0^{\v}} + \xi \|||\Phi||\|^2_1 + C_{\xi}a^2 \|(0, \hat{P}, \hat{\mathbf{S}})\hat{w}_0^{\v}\|^2_{\mathbf{Y}^{\v}} +  \xi \{\cdots\}_{1,\hat{w}_0^{\v}}\nonumber\\
		&\qquad \qquad \qquad  + C_{\xi} \{\cdots\}_{2,\hat{w}_0^{\v}}  + \xi\v \|\sqrt{u_s} \hat{\mathbf{S}}_{xy}\hat{w}_0^{\v}\|^2_{x=L}\mathbf{1_{D}} + C_{\xi}\|u_s \nabla_{\v}^2\hat{\mathbf{S}}_{y}\hat{w}_0^{\v}\|^2\mathbf{1_D},\label{18.42}\\
		& \bar{\mathcal{K}}_{3,\hat{w}_0^{\v}}(P) 
		\lesssim a [\Phi]_{3,\hat{w}_0^{\v}} +  \|||\Phi||\|^2_1 +  a^2 \|(0, \hat{P}, \hat{\mathbf{S}})\hat{w}_0^{\v}\|^2_{\mathbf{Y}^{\v}}  + \{\cdots\}_{1,\hat{w}_0^{\v}} + \{\cdots\}_{2,\hat{w}_0^{\v}} \nonumber\\
		&\qquad\qquad\quad   +  \v \|\sqrt{u_s} \hat{\mathbf{S}}_{xy}\hat{w}_0^{\v}\|^2_{x=L}\mathbf{1_{D}} +  \|u_s \nabla_{\v}^2\hat{\mathbf{S}}_{y}\hat{w}_0^{\v}\|^2\mathbf{1_D},\label{18.43}\\
		&\|u_s  \Delta_{\v}P_y\hat{w}^{\v}_0\|^2
		\lesssim  a^3 \|(\Phi,0,0)\|^2_{\mathbf{X}^{\v}} + a^2 \|(0, \hat{P}, \hat{\mathbf{S}})\hat{w}_0^{\v}\|^2_{\mathbf{Y}^{\v}}  + \sqrt{\v}\|(\mathbf{d}_2 \bar{g}_1, \bar{g}_1\mathbf{1}_{\{y\in[0,2]\}})\|^2 \nonumber\\
		&\qquad\qquad\qquad\,\, + \sqrt{\v}\{\cdots\}_{1,\hat{w}_0^{\v}} +  \{\cdots\}_{2,\hat{w}_0^{\v}},\label{C17.223-1}
	\end{align}
	where we have used the following notations
	\begin{align}\label{17.217}
		\{\cdots\}_{1,w}:&=   a^2\|(\phi,p, \zeta)\hat{w}\|^2_{\mathbf{X}^{\v}} + \sqrt{\v} \|u_s\nabla_{\v}^2\bar{g}_1w\|^2   + \|\nabla_{\v}\bar{g}_{1}w\|^2  + \f{1}{\v}\|g_{2y}w\|^2  +  \v \|\sqrt{u_s} \bar{g}_{1x}w\|_{x=L}^2 \nonumber\\
		&\quad   + \f{1}{\v} \|\sqrt{u_s} (b_2g_2)_{y}w\|^2_{x=L}  + \v \|u_s \nabla\Delta_{\v}\tilde{T}w\|^2  +  \v \|\nabla_{\v}^2\tilde{T}\hat{w}\|^2   + \|\sqrt{u_s}\Delta_{\v}\tilde{T}w\|^2_{x=0} \nonumber\\
		&\quad   + \v  \|\sqrt{u_s}\sqrt{\v}\tilde{T}_{xy}w\|^2_{x=0}   + \v \|\f{v_s}{\sqrt{u_s}}(\Delta_{\v}\tilde{T}, \tilde{T}_{yy})w\|^2_{x=L}  +  \|\f{\mathfrak{N}_{12}}{\sqrt{u_s}}w\|^2_{x=0}  + \v \|\f{\mathfrak{N}_{22}}{\sqrt{u_s}} w\|^2_{x=L}    \nonumber\\
		&\quad  + \v  \|(\mathcal{N}_{1x}, \mathcal{N}_{2y})w\|^2,
	\end{align}
	and
	\begin{align}\label{17.218-1}
		\{\cdots\}_{2,w}:&=   a^4\|(\phi,p, \zeta)\hat{w}\|^2_{\mathbf{X}^{\v}}
		+ \sqrt{\v} \|\nabla_{\v}\bar{g}_{1}\hat{w}\|^2  + \sqrt{\v} \|u_{s}\nabla^2_{\v}\bar{g}_1 \hat{w}\|^2  + \|u_s\pa^2_{y}g_{2}w\|^2   \nonumber\\
		&\quad  + \|\sqrt{u_{s}}(b_2g_2)_{y}w\|_{x=L}^2   + a^2\|\nabla_{\v}^2\tilde{T}\hat{w}\|^2   + a^2 \|\nabla_{\v}\Delta_{\v}\tilde{T}w\|^2  + \v \|u_s\nabla_{\v}^2\tilde{T}_yw\|^2  \nonumber\\
		&\quad  + a^2\v^2 \|u_s \nabla_{\v}\Delta_{\v}\tilde{T}_yw\|^2  + \v^2 \|u_s (\mathcal{N}_{1xy}, \mathcal{N}_{2yy})w\|^2  + \|(\mathfrak{N}_{12}, \pa_{\v}\mathfrak{N}_{12})w\|^2 \nonumber\\
		&\quad  +   \v \|(\mathfrak{N}_{22}, \pa_{\v}\mathfrak{N}_{22})\hat{w}\|^2  + \v^2 \|\sqrt{u_s}(b_2\mathfrak{h}_2, b_2\mathfrak{N}_{21})_{y}w\|_{x=L}  +a^2\v \|\f{\mathfrak{N}_{12}}{\sqrt{u_s}}w\|^2_{x=0}.
	\end{align}
	We remark that $\eqref{18.39-2}_3$ holds in the sense of distribution, and the uniqueness is under the constraints \eqref{18.41}-\eqref{18.43}.
\end{theorem}

It is very difficult to solve the pressure BVP \eqref{18.39-1}-\eqref{18.39-2}
with high regularity. 
In this section, we shall first establish the existence and uniform estimates for a $\delta$-approximation problem. Then, by passing to the limit $\delta \to 0+$ to obtain the solution for \eqref{18.39-1}-\eqref{18.39-2}.

\subsection{Construction of $\mathbf{w}^{\delta}$ for a linearized problem}\label{Sec5.1}
For later use, we denote 
\begin{align}\label{17.20-0}
	\begin{split}
		&U_s^{\delta}:= (u_s^{\delta}, v_s)^t \equiv ( u_s+ \delta, v_s)^t ,\\
		& U^{\v,\delta}:= (u^{\v,\delta} , v^{\v})^{t} \equiv (u^{\v}+\delta, v^{\v})^{t},\\
		&(p^{\delta}, S^{\delta}, \zeta^{\delta}):=(p,S,\zeta)(x,y+\delta),
	\end{split}
\end{align}
and
\begin{align}\label{17.20-5}
	\begin{split}
		&T^{\delta}:=T(S^{\delta}, p^{\delta}, q)\quad \&\quad  \rho^{\delta}:=\rho(p^{\delta}, T^{\delta}),\\
		&(p^{\v,\delta},T^{\v,\delta},\rho^{\v,\delta}):=(p_s,T_s, \rho_s) + \v^{N_0} (p^{\delta},T^{\delta}, \rho^{\delta}).
	\end{split}
\end{align}

It follows directly from Lemma \ref{lemA.8} that $p^{\delta}, S^{\delta}, \zeta^{\delta}\in H^3((0,L)\times \R_+)$ and 
\begin{align}\label{17.20-2}
	\begin{split}
		\|\nabla^{k}_{\v}(p^{\delta}, S^{\delta}, \zeta^{\delta})\|&\lesssim \|\nabla^{k}_{\v}(p, S, \zeta)\|,\,\,k=1,2,\\
		\|u^{\d}_s \nabla^{3}_{\v}(p^{\delta}, S^{\delta}, \zeta^{\delta})\|&\lesssim \|u_s \nabla^{3}_{\v}(p, S, \zeta)\|.
	\end{split}
\end{align}
These estimates are essential for rigorously taking the limit $\delta\to0+$ for the highest-order derivatives. Moreover, analogous to \eqref{3.94}, we have
\begin{align}\label{17.20-6}
	\v^{N_0}|u_s^{\delta} \nabla( p^{\delta}, T^{\delta},  \rho^{\delta}, S^{\delta}, \zeta^{\delta})| &\lesssim \v^{N_0-1}  + \v^{N_0} \|u_s \nabla( p, T,  \rho,  S,  \zeta)\|_{L^\infty}\lesssim \v^{N_0-2},
\end{align}
which will be used frequently in the subsequent analysis.

\smallskip

Let $0<\delta\ll\v$, and $\dot{B}$ be a given function.  We first consider the following linear $\delta$-approximation BVP:
\begin{align}\label{17.20}
	\begin{cases}
		\dis -\delta \Delta_{\v}\mathbf{w}^{\delta} + \mathbf{w}^{\delta}  + \v (\mu+\lambda)  \f{1}{p^{\v,\delta}}(U^{\v,\delta} \cdot\nabla)\mathbf{w}^{\delta} = \hat{\mathcal{G}}^{\delta},\,\, (x,y)\in (0,L)\times \mathbb{R}_+,\\
		\dis -\delta \mathbf{w}^{\delta}_x + (\mu+\lambda) \f{1}{p_s} u_s^{\delta}\mathbf{w}^{\delta} \big|_{x=0}=-\sigma \mathbf{d}_{11}\dot{B},\\
		\dis \mathbf{w}^{\delta}_x\big|_{x=L}=\mathbf{w}^{\delta}_y\big|_{y=0}=0,
	\end{cases}
\end{align}
where 
\begin{align}\label{17.20-1}
	\hat{\mathcal{G}}^{\delta}:=\mathcal{G}(\hat{P}^{\d},\hat{\mathbf{S}}^{\delta}, \Phi, p^{\delta}, S^{\delta}, \zeta^{\delta}, \phi, q, q^{\delta}),
\end{align}
with $\hat{P}^{\d}=\hat{P}(x,y+\d)$ and  $\hat{\mathbf{S}}^{\d}=\hat{\mathbf{S}}(x,y+\d)$. For $q$, we only modify one of $q$  in  \eqref{17.17-2} to be $q^{\delta}=q(x,y+\delta)$, while keeping all others $q$  unchanged. Specifically, the term $\v \big(\f{1}{p_s}T_{sy} q\big)_{xy}$ in \eqref{17.17-2} should be modified  to $\v \big(\f{1}{p_s}T_{sy} q^{\delta}\big)_{xy}$ with $q^{\delta}:=q(x,y+\delta)$ in the case of DT. Then the following estimate holds
\begin{align}
	\v^2 \|u_s^{\delta}\big(p_s^{-1}T_{sy} q^{\delta}\big)_{xyy}w\|^2\lesssim \v^2 \|u_s^{\delta}q^{\delta}_{xyy}w\|^2 + \cdots\lesssim \v^2 \|\sqrt{u_s}q_{xyy}\hat{w}\|^2 + \cdots.
\end{align}
which will be utilized in $\eqref{17.181-1}_4$. For the case of NT,  we note that such a modification is not necessary due to $T_{sy}|_{y=0}=0$, however,  we still keep such modification on $q$ for notational consistency.

\smallskip

\smallskip

\begin{lemma}\label{lemP7.2} 
	Recall $w$ in \eqref{6.0-0}. There exists a unique weak solution $\mathbf{w}^{\delta}$ to BVP \eqref{17.20} satisfying 
	\begin{align}\label{17.25}
		&(1-\xi)\Big\{\delta \|\nabla_{\v}\mathbf{w}^{\delta}w\|^2 +  \|\mathbf{w}^{\delta}w\|^2 +    (\mu+\lambda) \v \int_0^\infty \f{u_s^{\delta}}{2p_s} |\mathbf{w}^{\delta}|^2w^2 dy \Big|_{x=L} \Big\} \nonumber\\
		&\leq  \f12 \f{\sigma^2}{\mu+\lambda} \v  \int_0^\infty \mathbf{d}_{11}^2\f{p_s}{u_s^{\delta}} |\dot{B}|^2 w^2 dy \Big|_{x=0}  + C_{\xi} \|\hat{\mathcal{G}}^{\delta} w\|^2 +  C_{\xi}\v \|\mathbf{w}^{\delta}w_y\|^2, 
	\end{align}
	where $\xi>0$ is a small positive constant with $\v\ll \xi\leq 1$.
\end{lemma}

\noindent{\bf Proof.} 
Multiplying \eqref{17.20} by a smooth test function $\varphi$, one obtains 
\begin{align}\label{17.21}
	\mathcal{B}^{\delta}[\mathbf{w}^{\delta},\varphi]
	&=\iint \hat{\mathcal{G}}^{\delta}\,\varphi dydx + \v \int_0^\infty \big[-\delta \mathbf{w}^{\delta}_x + (\mu+\lambda) \f{u^{\delta}_s}{p_s}  \mathbf{w}^{\delta}\big]\varphi  dy \Big|_{x=0} \nonumber\\
	&=\iint \hat{\mathcal{G}}^{\delta}\,\varphi dydx - \sigma \v \int_0^\infty \mathbf{d}_{11} \dot{B}\, \varphi  dy \Big|_{x=0}, 
\end{align}
where we have used the following notation of bilinear form:
\begin{align*}
	\mathcal{B}^{\delta}[\mathbf{w}^{\delta},\varphi]:&=\delta \iint \nabla_{\v}\mathbf{w}^{\delta}\cdot \nabla_{\v}\varphi dydx + \iint \mathbf{w}^{\delta}\varphi dydx  +(\mu+\lambda) \v \int_0^\infty \f{u^{\v,\delta}}{p^{\v,\delta}}  \mathbf{w}^{\delta} \varphi dy \Big|_{x=L} \nonumber\\
	&\quad - (\mu+\lambda) \v \iint \mathbf{w}^{\delta}\big[\varphi \mbox{div} \big(\f{1}{p^{\v,\delta}}U^{\v,\delta}\big) + \f{1}{p^{\v,\delta}} U^{\v,\delta}\cdot \nabla \varphi\big] dydx \nonumber\\
	&\quad  - (\mu+\lambda) \v \int_0^\infty  \big[\f{u^{\v,\delta}}{p^{\v,\delta}}- \f{u_s^{\delta}}{p_s} \big] \mathbf{w}^{\delta} \varphi dy \Big|_{x=0}.
\end{align*}
Using \eqref{17.20-6} and \eqref{3.92}--\eqref{3.94}, one can show that 
\begin{align}\label{17.23}
	\mathcal{B}^{\delta}[\mathbf{w}^{\delta},\mathbf{w}^{\delta}]
	&\geq \delta \|\nabla_{\v}\mathbf{w}^{\delta}\|^2 + (1-Ca\v) \|\mathbf{w}^{\delta}\|^2 + (\mu+\lambda) \v\int \f{u^{\v,\delta}}{2p^{\v,\delta}} |\mathbf{w}^{\delta}|^2 dy \Big|_{x=L} \nonumber\\
	&\quad  + [1-C\v^{N_0}] (\mu+\lambda) \v\int \f{u^{\v,\delta}}{2p^{\v,\delta}} |\mathbf{w}^{\delta}|^2 dy \Big|_{x=0}.
\end{align}

Noting \eqref{17.23} and applying the Lax-Milgram theorem, there exists a unique weak solution $\mathbf{w}^{\delta} \in H^1$ to BVP \eqref{17.20} in the sense of \eqref{17.21}.
Taking $\varphi=\mathbf{w}^{\delta} w^2$ in \eqref{17.21}, 
one  easily concludes \eqref{17.25}. Therefore the proof of Lemma \ref{lemP7.2} is completed. $\hfill\Box$

\smallskip
 
For later use, we still need uniform-in-$\delta$ estimate on $\|u_s\nabla_{\v}\mathbf{w}^{\delta}w\|$. 
\begin{lemma}\label{lemP7.3}
	Let  $\v\ll \xi$. We have 
	\begin{align}\label{17.126}
		&(1-\xi)\Big\{\delta \|u_s^{\delta} \nabla_{\v}\mathbf{w}^{\delta}_y w\|^2 +  \|u_s^{\delta}\mathbf{w}^{\delta}_yw\|^2 +  (\mu+\lambda) \v \int_0^\infty  \f{(u_s^{\delta})^3}{2p_s}   |\mathbf{w}^{\delta}_y|^2 w^2 dy \Big|_{x=L} \Big\}\nonumber\\
		&\leq \Big[\f{\sigma^2}{2(\mu+\lambda)}+\xi\Big]\v \int_0^\infty p_s\mathbf{d}_{11}^2 u_s^{\delta} |\dot{B}_{y}|^2 w^2  dy \Big|_{x=0} + \v^{\f32} \|(U_{s}^{\delta}\cdot \nabla) \mathbf{w}^{\delta}w\|^2  + \v \|u_s^{\delta} \mathbf{w}^{\delta}_y\|^2 \nonumber\\
		&\quad   + \sqrt{\v}\|\mathbf{w}^{\delta}w\|^2  + C_{\xi} a^2 \delta \|\mathbf{w}^{\delta}_y\hat{w}\|^2  +   C_{\xi} a^4 \v \|\f{\dot{B}w}{\sqrt{u^{\delta}_s}}\|^2_{x=0}  +  C_{\xi} \|u_s^{\delta}\hat{\mathcal{G}}^{\delta}_yw\|^2,
	\end{align}
	and 
	\begin{align}\label{17.99}
		&\f{\delta^2}{\v} \|\Delta_{\v}\mathbf{w}^{\delta} w\|^2 + \v\|(U_s^{\delta}\cdot \nabla)\mathbf{w}^{\delta}w\|^2  +  \delta \int_0^\infty  u_s^{\delta}  |\sqrt{\v}\mathbf{w}^{\delta}_x|^2w^2 dy\Big|_{x=0}  +  \delta \int_0^\infty u_s^{\delta}  |\mathbf{w}^{\delta}_{y}|^2 w^2\, dy \Big|_{x=L} \nonumber\\ 
		&\lesssim  \delta \|u_s^{\delta} \nabla_{\v}\mathbf{w}^{\delta}_y w\|^2  + \v \|u_s^{\delta}\mathbf{w}^{\delta}_{y}w\|^2  +  \f{\delta}{\v}\|\nabla_{\v}\mathbf{w}^{\delta} \hat{w}\|^2   +  \f{1}{\v}\|\mathbf{w}^{\delta}w\|^2+ \f{1}{\v}\|\hat{\mathcal{G}}^{\delta}w\|^2.
	\end{align}
\end{lemma} 

\noindent{\bf Proof.}  Applying $\pa_y$ to \eqref{17.119}, one obtains that 
\begin{align}\label{17.122}
-\delta \Delta_{\v}\mathbf{w}^{\delta}_y  + \mathbf{w}^{\delta}_y + (\mu+\lambda) \v \f{1}{p^{\v,\delta}} (U^{\v,\delta}\cdot\nabla)\mathbf{w}^{\delta}_y
 = \hat{\mathcal{G}}^{\delta}_y - (\mu+\lambda) \v \Big(\f{1}{p^{\v,\delta}} U^{\v,\delta}\Big)_y \cdot\nabla \mathbf{w}^{\delta}.
\end{align}
In the following, we divide the proof into two steps. 

{\it Step 1.} Multiplying \eqref{17.122} by $(u_s^{\delta})^2\mathbf{w}^{\delta}_yw^2$, one gets
\begin{align}\label{17.123}
	&\delta \|u_s^{\delta} \nabla_{\v}\mathbf{w}^{\delta}_yw\|^2 + \|u_s^{\delta}\mathbf{w}^{\delta}_yw\|^2  -  (\mu+\lambda) \v \iint \mbox{div}\Big(\f{(u_s^{\delta}w)^2}{2p^{\v,\delta}} U^{\v,\delta}\Big) |\mathbf{w}^{\delta}_y|^2 dydx  \nonumber\\
	& +  (\mu+\lambda) \v \int_0^\infty  \f{u^{\v,\delta}}{2p^{\v,\delta}} (u_s^{\delta})^2  |\mathbf{w}^{\delta}_y|^2 w^2 dy \Big|_{x=0}^{x=L} 
	 + \delta  \v  \int_0^\infty (u_s^{\delta}w)^2 \mathbf{w}^{\delta}_{xy} \mathbf{w}^{\delta}_y dy \Big|_{x=0} \nonumber\\
	&= -   (\mu+\lambda) \v \Big\la  \Big(\f{1}{p^{\v,\delta}} U^{\v,\delta}\Big)_y \cdot\nabla \mathbf{w}^{\delta},\, (u_s^{\delta})^2 \mathbf{w}^{\delta}_y w^2\Big\ra - \big\la \delta \nabla_{\v}\mathbf{w}^{\delta}_y, \, \nabla_{\v}\big((u_s^{\delta} w)^2\big) \mathbf{w}^{\delta}_y\big\ra  \nonumber\\
	&\quad  + \big\la \tilde{\mathcal{G}}^{\delta}_y,\, (u_s^{\delta})^2 \mathbf{w}^{\delta}_y w^2\big\ra,
\end{align}
where we have used $\mathbf{w}^{\delta}_{xy}\big|_{x=L}=0$ and $\mathbf{w}_y^{\delta}\big|_{y=0}=0$.

It is clear that 
\begin{align}\label{17.124}
	\begin{split}
	 \Big|\v \iint \mbox{div}\Big(\f{1}{p^{\v,\delta}}(u_s^{\delta}w)^2 U^{\v,\delta}\Big) |\mathbf{w}^{\delta}_y|^2 dydx \Big| &\lesssim  \v \|u_s^{\delta} \mathbf{w}^{\delta}_y \hat{w}\|^2,\\
	\Big|\v \Big\la  \Big(\f{1}{p^{\v,\delta}} U^{\v,\delta}\Big)_y \cdot\nabla \mathbf{w}^{\delta},\, (u_s^{\delta})^2 \mathbf{w}^{\delta}_y w^2 \Big\ra \Big|
	&\lesssim \sqrt{\v} \|u_s^{\delta} \mathbf{w}^{\delta}_yw\|^2 +  \sqrt{\v} \|u_s^{\delta}\sqrt{\v}\mathbf{w}^{\delta}_xw\|^2,\\
	\big|\big\la \delta \nabla_{\v}\mathbf{w}^{\delta}_y, \, \nabla_{\v}\big((u_s^{\delta} w)^2\big) \mathbf{w}^{\delta}_y\big\ra 
	& \lesssim \xi \delta \|u_s^{\delta} \nabla_{\v}\mathbf{w}^{\delta}_yw\|^2  + C_{\xi} a^2 \delta \|\mathbf{w}^{\delta}_y\hat{w}\|^2,\\
	\big\la \hat{\mathcal{G}}^{\delta}_y,\, (u_s^{\delta})^2 \mathbf{w}^{\delta}_y w^2 \big\ra
	&\lesssim \xi \|u_s^{\delta} \mathbf{w}^{\delta}_yw\|^2 + C_{\xi} \|u_s^{\delta} \hat{\mathcal{G}}^{\delta}_yw\|^2,
	\end{split}
\end{align}
where we have used \eqref{17.20-6} in $\eqref{17.124}_{1,2}$.

For boundary terms at $x=0$,  using $\eqref{17.20}_2$ and Holder's inequality, one has
\begin{align}\label{17.125}
	& \v  \int_0^\infty (u_s^{\delta}w)^2 \mathbf{w}^{\delta}_y \Big[-\delta \mathbf{w}^{\delta}_{xy} +  (\mu+\lambda)  \f{u^{\v,\delta}}{2p^{\v,\delta}}  \mathbf{w}^{\delta}_y\Big] dy \Big|_{x=0}  \nonumber\\
	&=-(\mu+\lambda) \v \int_0^\infty\f{1+O(\v)}{2p_s} (u_s^{\delta})^3  |\mathbf{w}^{\delta}_y|^2 w^2 dy \Big|_{x=0}  -(\mu+\lambda)  \v  \int_0^\infty (u_s^{\delta}w)^2 (\f{u_s^{\delta}}{p_s})_y \mathbf{w}^{\delta}_y \mathbf{w}^{\delta} dy \Big|_{x=0} \nonumber\\
	&\quad - \sigma \int_0^\infty \v (u_s^{\delta}w)^2 \mathbf{d}_{11} \mathbf{w}^{\delta}_y \dot{B}_{y} dy \Big|_{x=0} - \sigma \int_0^\infty \v (u_s^{\delta}w)^2 \mathbf{d}_{11y}  \mathbf{w}^{\delta}_y  \dot{B} dy \Big|_{x=0} \nonumber\\
	&\leq \Big[\f{\sigma^2}{2(\mu+\lambda)}+\xi\Big]\v \int_0^\infty p_s \mathbf{d}_{11}^2u_s^{\delta} |\dot{B}_{y}|^2 w^2 dy \Big|_{x=0} + C_{\xi} \v \int_0^\infty u_s^{\delta} |\mathbf{d}_2\mathbf{w}^{\delta}|^2w^2 dy\Big|_{x=0} \nonumber\\
	&\quad + C_{\xi} \v \int_0^\infty u_s^{\delta} |\mathbf{d}_2|^2  \dot{B}^2 w^2 dy\Big|_{x=0}.
\end{align}

A direct calculation shows that 
\begin{align}\label{17.99-1}
	\|u_s^{\d}\sqrt{\v}\mathbf{w}^{\delta}_xw\|^2 &\cong \v \|(U_{s}^{\delta}\cdot \nabla) \mathbf{w}^{\delta}w\|^2 + O(1)\v \|u_s^{\delta}\mathbf{w}^{\delta}_{y}w\|^2,
\end{align}
and
\begin{align}\label{17.99-2}
  \v \int_0^\infty u_s^{\delta} |\mathbf{w}^{\delta}|^2 w^2 dy\Big|_{x=0} 
 &\lesssim \f{a}{L} \v \|\mathbf{w}^{\delta}w\|^2 + a \sqrt{\v} \|\sqrt{\v}u_s^{\delta}\mathbf{w}^{\delta}_xw\|\cdot \|\mathbf{w}^{\delta}w\| \nonumber\\
	&\lesssim \v^{\f32} \|(U_{s}^{\delta}\cdot \nabla) \mathbf{w}^{\delta}w\|^2 + \sqrt{\v} \|u_s^{\delta}\mathbf{w}^{\delta}_{y}w\|^2 + \sqrt{\v}\|\mathbf{w}^{\delta}w\|^2.
\end{align}
Substituting \eqref{17.124}-\eqref{17.125} into \eqref{17.123}, then using \eqref{17.99-1}-\eqref{17.99-2}, one  concludes \eqref{17.126}.

\smallskip

{\it Step 2.} We control  $\v \|(U_{s}^{\delta}\cdot \nabla) \mathbf{w}^{\delta} w\|^2$.  It is clear from \eqref{17.20} that 
\begin{align}\label{17.5.20}
&\delta^2 \|\Delta_{\v}\mathbf{w}^{\delta}w\|^2 + (\mu+\lambda)^2\v^2\|\f{w}{p^{\v,\delta}}(U^{\v,\delta}\cdot \nabla)\mathbf{w}^{\delta}\|^2 - 2(\mu+\lambda)\v \delta \iint \f{u^{\v,\delta}}{p^{\v,\delta}} \mathbf{w}^{\delta}_x  \Delta_{\v}\mathbf{w}^{\delta}   w^2\, dydx \nonumber\\
&\quad - 2(\mu+\lambda)\v \delta \iint \f{v^{\v}}{p^{\v,\delta}}   \mathbf{w}^{\delta}_y  \Delta_{\v}\mathbf{w}^{\delta}  \, w^2 dydx 
\leq \|\hat{\mathcal{G}}^{\delta} w - \mathbf{w}^{\delta}w\|^2.
\end{align}
Using \eqref{17.99-1} and \eqref{3.92}, we have 
\begin{align}\label{17.26-1}
	\|\f{w}{p^{\v,\delta}}(U^{\v,\delta}\cdot \nabla)\mathbf{w}^{\delta}\|^2
	&\cong (1-\v)\|\f{w}{p_s}(U^{\delta}_s\cdot \nabla)\mathbf{w}^{\delta}\|^2 + O(1)\v^{2N_0-1} \|U^{(\v)}\cdot \nabla \mathbf{w}^{\delta}w\|^2\nonumber\\
	&\cong \|(U^{\delta}_s\cdot \nabla)\mathbf{w}^{\delta} w\|^2- C\v^{2N_0-1} \|u_s^{\delta}\mathbf{w}^{\delta}_{y}w\|^2.
\end{align}

Integrating by parts in $x$,  noting $\mathbf{w}^{\delta}_x|_{x=L}=0$  and \eqref{17.20-6}, one can obtain
\begin{align}
- 2 \v^2 \delta \iint \f{u^{\v,\delta}}{p^{\v,\delta}} \mathbf{w}^{\delta}_x  \mathbf{w}^{\delta}_{xx} w^2\, dydx 
	&= \v^2 \delta \int_0^\infty \f{u^{\v,\delta}}{p^{\v,\delta}} |\mathbf{w}^{\delta}_x|^2 w^2 dy\Big|_{x=0}  + O(1)\v \delta \|\sqrt{\v}\mathbf{w}^{\delta}_xw\|^2.
\end{align}

Noting $\mathbf{w}^{\delta}_y\big|_{y=0}=0$, \eqref{17.20-6} and \eqref{3.93}, we integrate by parts to get
\begin{align}
&- 2\v \delta \iint \f{u^{\v,\delta}}{p^{\v,\delta}} \mathbf{w}^{\delta}_x \mathbf{w}^{\delta}_{yy} \, w^2 dydx \nonumber\\
	&= \v \delta \iint \Big\{2\big(\f{u^{\v,\delta}}{p^{\v,\delta}} w^2\big)_y \mathbf{w}^{\delta}_x \mathbf{w}^{\delta}_{y}- \big(\f{u^{\v,\delta}}{p^{\v,\delta}}\big)_x |\mathbf{w}^{\delta}_{y}|^2 w^2 \Big\} dydx + \v \delta \int \f{u^{\v,\delta}}{p^{\v,\delta}} |\mathbf{w}^{\delta}_{y}|^2 w^2 \, dy \Big|_{x=0}^{x=L} \nonumber\\
	&=\v \delta \int_0^\infty \f{u^{\v,\delta}}{p^{\v,\delta}} |\mathbf{w}^{\delta}_{y}|^2 w^2\, dy \Big|_{x=L} 
	+\v \delta \|u_s^{\delta} \nabla_{\v}\mathbf{w}^{\delta}_yw\|^2  +   \delta \|\mathbf{w}^{\delta}_y\hat{w}\|^2,
\end{align}
where we have used 
\begin{align}\label{17.23-1}
	\v \delta \int_0^\infty u_s^{\delta} |\mathbf{w}^{\delta}_{y}|^2 w^2\, dy \Big|_{x=0}
	&\lesssim \v \delta \Big\{\|u_s^{\delta} \mathbf{w}^{\delta}_{xy}w\|\cdot \|\mathbf{w}^{\delta}_yw\| + \f{1}{L} \|\mathbf{w}^{\delta}_yw\|^2 \Big\}
	\nonumber\\
	&\lesssim \v \delta \|u_s^{\delta} \nabla_{\v}\mathbf{w}^{\delta}_yw\|^2  +   \delta \|\mathbf{w}^{\delta}_yw\|^2.
\end{align}

Integrating by parts, and  using $\mathbf{w}^{\delta}_x|_{x=L}=0$, \eqref{17.20-6}, \eqref{17.23-1} and \eqref{3.93}, one obtains that 
\begin{align}
	& 2\v^2 \delta \iint \f{v^{\v}}{p^{\v,\delta}} w^2 \mathbf{w}^{\delta}_y \mathbf{w}^{\delta}_{xx}  dydx\nonumber\\
	&= 2\v^2 \delta \int_0^\infty \f{v^{\v}}{p^{\v,\delta}} w^2 \mathbf{w}^{\delta}_y \mathbf{w}^{\delta}_{x}  dy \Big|_{x=0}  + 2\v^2 \delta \iint \Big\{\f{v^{\v}}{p^{\v,\delta}} \mathbf{w}^{\delta}_{xy} + \big(\f{v^{\v}}{p^{\v,\delta}}\big)_x \mathbf{w}^{\delta}_y \Big\} \mathbf{w}^{\delta}_{x}  w^2  dydx \nonumber\\
	&\lesssim \v^{\f32}\delta \int_0^\infty u_s |\sqrt{\v}\mathbf{w}^{\delta}_x|^2 w^2 dy \Big|_{x=0} + \v^{\f32} \delta \|u_s^{\delta} \nabla_{\v}\mathbf{w}^{\delta}_yw\|^2  + \sqrt{\v} \delta \|\mathbf{w}^{\delta}_yw\|^2
\end{align}
and
\begin{align}\label{17.5.26}
	- 2\v \delta \iint \f{v^{\v}}{p^{\v,\delta}} \mathbf{w}^{\delta}_y \mathbf{w}^{\delta}_{yy} w^2\, dydx 
	&=\v \delta \iint \big(\f{v^{\v}}{p^{\v,\delta}} w^2\big)_y  |\mathbf{w}^{\delta}_y|^2 dydx 
	  =O(1) \v\delta \|\mathbf{w}^{\delta}_y\hat{w}\|^2.
\end{align}

Finally, combining \eqref{17.5.20}-\eqref{17.5.26},  
we conclude \eqref{17.99}. Therefore the proof of Lemma \ref{lemP7.3} is completed. $\hfill\Box$

\subsection{The construction of $B^{\delta}$.}  
The task now is to recover the pressure function $P^{\delta}$ from the relation $\Delta_{\v}P^{\delta} \backsim  \mathbf{w}^{\delta}$. As mentioned in Section \ref{Sec1.2}, we work with the {\it good unknown} $B$ instead of the pressure itself. 
\subsubsection{Reduction of necessary boundary condition \eqref{2.104}}
To solve the elliptic problem $\Delta_{\v}P^{\delta}\cong \mathbf{w}^{\delta}$, we can only utilize the necessary boundary condition \eqref{2.104} since no additional conditions on pressure  are available  except \eqref{17.15}. 

\smallskip

For later use, we denote 
\begin{align}\label{17.28}
	\bm{\chi}(y)=
	\begin{cases}
		\chi(y), \quad \mbox{for NT},\\
		0, \qquad\,\,\, \mbox{for DT}.
	\end{cases}
\end{align}

\begin{lemma}\label{lemAF-1}
The necessary boundary condition \eqref{2.104} can be rewritten as
\begin{align}
B \big|_{y=0} &=\frac{\mu}{\rho^{\v}} \Phi_{yyy},\label{17.17}\\
\Big\{(\sigma-1)p_s \mathbf{d}_{11}  B - \v \frac{1}{2}\lambda  [1+\bm{\chi}] u_s B_{x} \Big\} \Big|_{x=0} & =-\frac{1}{2}\lambda  [1+\bm{\chi}] u_s \Delta_{\v}P + \mathbf{g}_{b1}(P,\Phi,p,S, \zeta,\phi),\label{17.19-5}\\
\Big\{\v B_{x} - \lambda\v\frac{b_2}{2p_s} u_s B_{yy}\Big\} \Big|_{x=L} & = \Delta_{\v}P + \mathbf{g}_{b2}(P,\mathbf{S},\Phi,p,S, \zeta,\phi),\label{17.19-3}
\end{align}
where 
\begin{align}
\mathbf{g}_{b1}:&= (\mu+\lambda)p_s \big(\frac{u^{\v} }{p^{\v}}-\frac{u_s}{p_s}\big) \Delta_{\v}P  - \f{1}{2}\lambda p_s\Big[\big(\frac{ u^{\v}}{p^{\v}}-\frac{ u_s}{p_s}\big) + \f{\bm{\chi}}{\rho_s} \big(\f{u^{\v}}{T^{\v}}-\f{u_s}{T_s}\big) -  \v^{N_0} \f{u^{\v}\rho^{(\v)}}{\rho_sp^{\v}}\Big]P_{yy}  \nonumber\\
&\quad  + \lambda\v  \frac{p_s}{p^{\v}} v^{\v} P_{xy}  -  \lambda  p_s  \frac{u^{\v}}{T^{\v}} \Big\{\big(\f{\bm{\chi}}{2\rho_s}\big)_y P_y + \big(\f{\bm{\chi}}{2\rho_s} \big)_{yy} P\Big\}   - p_s G_{11}(\tilde{T}(\zeta,p,q))  - p_s G_{12}(p,q)  \nonumber\\
&\quad  - \frac{1}{2}\v \lambda  [1+\bm{\chi}] u_s \bar{g}_{1x} + p_s\hat{d}_{11} \zeta_x   + p_s\mathfrak{g}_1(q)    + \rho_s u_s^2 \mathbf{q}_{xy}   - p_s\mathfrak{N}_{12},\label{17.19-7}\\
\mathbf{g}_{b2}:&= \lambda\v \big(\frac{b_2}{2p_s} u_s\big)_y P_{xy} - \pa_y G_3 + \big(b_2 [ g_2+ \v \mathfrak{g}_2 - \v\mathfrak{h}_2 - \v\mathfrak{N}_{21}]\big)_y   + \v \bar{g}_{1x} - \lambda\v\frac{b_2}{2p_s} u_s \bar{g}_{1yy}.\label{17.19-8}
\end{align}
with
\begin{align}\label{17.19-10}
G_3(\Phi,P, \mathbf{S},  p, \zeta, \phi):&= - b_2\frac{\mu\v}{p^{\v}} v^{\v} \Delta_{\v}P   - b_2\frac{\lambda\v}{p^{\v}} v^{\v} P_{yy} - \f12 \lambda \v b_2 P_{xy} \Big[\big(\frac{ u^{\v}}{p^{\v}}-\frac{ u_s}{p_s}\big) + \v^{N_0} \f{u^{\v}\rho^{(\v)}}{\rho_sp^{\v}}\Big]\nonumber \\
&\hspace{-1.3cm}  + \mathbf{1_D}   \lambda\v\frac{b_2}{T^{\v}} u^{\v} \mathbf{S}_{xy} + \v b_2 m u_s^2 \mathbf{q}_{xx} + \v b_2 G_{21}(\tilde{T}(\zeta,p,q)) + \v b_2 G_{22}(p,q)  + \v b_2 \mathfrak{N}_{22}.
\end{align}
The definitions of $\mathfrak{N}_{12}, \mathfrak{N}_{22}$ are given in \eqref{17.18-1} and \eqref{17.19-4} below, respectively.
\end{lemma}

\noindent{\bf Proof.}  In the following, we divide the proof into three steps. 

{\it Step 1.} It is clear that \eqref{17.17} follows directly  from \eqref{2.113}, \eqref{17.14} and $\eqref{2.104}_1$.

{\it Step 2.} Noting \eqref{17.28}, we represent  the boundary conditions of pseudo entropy for both  NT and DT 
\begin{align}
	\dis \mathbf{S}|_{x=0}=-\f{\bm\chi}{2\rho_s}P,
\end{align}
then we can rewrite $\eqref{2.104}_2$, i.e., $\mathbf{A}_1\big|_{x=0}=0$,  as
\begin{align}\label{17.18}
	&\Big\{(\sigma-1)p_s \mathbf{d}_{11}  B + \frac{1}{2}\lambda  [1+\bm{\chi}] u_s P_{yy}\Big\} \Big|_{x=0}\nonumber\\
	&= - \f{1}{2}\lambda p_s\Big[\big(\frac{ u^{\v}}{p^{\v}}-\frac{ u_s}{p_s}\big) + \f{\bm{\chi}}{\rho_s} \big(\f{u^{\v}}{T^{\v}}-\f{u_s}{T_s}\big) -  \v^{N_0} \f{u^{\v}\rho^{(\v)}}{\rho_sp^{\v}}\Big]P_{yy} + (\mu+\lambda)p_s \big(\frac{u^{\v} }{p^{\v}}-\frac{u_s}{p_s}\big) \Delta_{\v}P   \nonumber\\
	&\quad  + \lambda\v  \frac{p_s}{p^{\v}} v^{\v} P_{xy}  -  \lambda  p_s  \frac{u^{\v}}{T^{\v}} \Big\{\big(\f{\bm{\chi}}{2\rho_s}\big)_y P_y + \big(\f{\bm{\chi}}{2\rho_s} \big)_{yy} P\Big\}  - p_s G_{11}(\tilde{T}(\zeta,p,q))  - p_s G_{12}(p,q) \nonumber\\
	&\quad + p_s\hat{d}_{11} \zeta_x   + p_s\mathfrak{g}_1(q)  + \rho_s u_s^2 \mathbf{q}_{xy}   - p_s\mathfrak{N}_{12}.
\end{align}
where we have used 
\begin{align}\label{17.18-1}
	\mathfrak{N}_{12}:=\mathfrak{N}_{1} - \mathfrak{N}_{11}\quad  \mbox{with}\,\,\, \mathfrak{N}_{11} \,\,\&\,\,\mathfrak{N}_{1} \,\, \mbox{defined in} \,\,  \eqref{17.14-1}\,\, \& \,\, \eqref{9.26}\,\, \mbox{respectively}.
\end{align}  
Noting \eqref{17.14}, we have 
\begin{equation}\label{17.18-2}
	P_{yy} 
	=\Delta_{\v}P-\v B_{x}+ \v \bar{g}_{1x},
\end{equation} 
then we can obtain \eqref{17.19-5} from \eqref{17.18}. 

\smallskip 

{\it Step 3.} Similarly, using \eqref{2.113}, $\eqref{2.111}_2$ and $\eqref{2.111-6D}_2$, one deduces from $\eqref{2.104}_3$, i.e.,  $\mathbf{A}_2\big|_{x=L}=0$, that
\begin{align}\label{17.19-2}
	\Big\{P_y + \lambda\v\frac{b_2}{2p_s} u_s P_{xy}\Big\} \Big|_{x=L} 
	&=G_3 - b_2 [ g_2(p,\phi, T(S,p,q))+ \v \mathfrak{g}_2(q) - \v\mathfrak{h}_2 - \v\mathfrak{N}_{21}],
\end{align}
where we have used the notation (see \eqref{9.28} for the definition of $\mathfrak{N}_2$)
\begin{align}\label{17.19-4}
	\begin{split}
		\mathfrak{N}_{22}:&=\mathfrak{N}_{2} - \mathfrak{N}_{21}\quad \mbox{with}\,\, \mathfrak{N}_{21}:=-\lambda \v^{N_0} \f{1}{\rho^{\v}} u^{(\v)}_y \rho^{(\v)}_y.
	\end{split}
\end{align}
Applying $\pa_y$ to \eqref{17.19-2}, one obtains
\begin{align*}
	P_{yy}  \Big|_{x=L} &  =  - \lambda\v\frac{b_2}{2p_s} u_s P_{xyy} - \lambda\v \big(\frac{b_2}{2p_s} u_s\big)_y P_{xy}  + \pa_y G_3 - \big(b_2 [ g_2+ \v \mathfrak{g}_2 - \v\mathfrak{h}_2 - \v\mathfrak{N}_{21}]\big)_y,
\end{align*}
which, together with \eqref{17.18-2} and \eqref{17.14},  concludes \eqref{17.19-3}. Therefore the proof of Lemma \ref{lemAF-1} is completed. $\hfill\Box$

\subsubsection{An auxiliary function}
Since $\mathbf{w}^{\delta}$ is a viscous approximation of original problem, in the following, we can also make some approximations on the boundary conditions \eqref{17.17}--\eqref{17.19-3}.

We first consider an odd extension of $\Phi$, i.e.,
\begin{align}
\tilde{\Phi}(x,y)=
\begin{cases}
\Phi(x,y),\qquad\,\,\, (x,y) \in [0,L]\times \R_+,\\
-\Phi(-x,y),\quad (x,y) \in [-L,0]\times \R_+.
\end{cases}
\end{align}
Noting $\Phi_{x}|_{x=0}=\Phi_{xxx}|_{x=0}=0$, then  we further extend   $\tilde{\Phi}$ periodically in $x$ so that it is well-defined in $(x,y)\in \R\times \R_+$. It is clear that $\tilde{\Phi} \in H^4([-L,L]\times \R_+)$. 
We smooth $\tilde{\Phi}$ in $x$-direction,  denote it as $\tilde{\Phi}^{\delta}$ with $\delta^{\tau} (0<\tau\ll 1)$ as the mollifier parameter, such that
\begin{align}\label{5.28}
\begin{split}
&\|||\tilde{\Phi}^{\delta}||\|_1 + \sqrt{\delta} \|\tilde{\Phi}^{\delta}_{xxyyy}\| \lesssim \|||\Phi||\|_1 \,\,\,\, \mbox{and}\,\,\,\,
 \|\tilde{\Phi}^{\delta}_{yyy}-\Phi_{yyy}\| +  \|||\tilde{\Phi}^{\delta}-\Phi||\| \to 0\,\,\mbox{as}\,\, \delta\to0.
\end{split}
\end{align}
Also, we extend $\rho^{\v}$ evenly   to the domain $[-L, L] \times \mathbb{R}_+$,  then extend it periodically in the $x$-direction to $\mathbb{R} \times \mathbb{R}_+$, and finally smoothed along $x$ as above, and denoted it as $\tilde{\rho}^{\v,\d}$.

To deal with the boundary condition  \eqref{17.17}  at $y=0$, we  introduce an auxiliary function $\theta$, i.e.,
\begin{align}\label{17.26}
	\begin{cases}
	\dis 	-\Delta_{\v} \theta^{\d} + \theta^{\d} =0, \,\, (x,y) \in (0,L)\times \mathbb{R}_+,\\
	\dis 	\theta^{\d}|_{x=0}=0,\quad \theta^{\d}_x|_{x=L}=0,\\
	\dis 	\theta^{\d}|_{y=0}=\f{\mu}{\tilde{\rho}^{\v,\d}} \tilde{\Phi}^{\delta}_{yyy} \in H^{3/2}_{y=0}.
	\end{cases}
\end{align}
We remark that $\f{\mu}{\tilde{\rho}^{\v,\d}} \tilde{\Phi}^{\delta}_{yyy}\to \f{\mu}{\rho^{\v}} \Phi_{yyy} $ in $H^{1/2}_{y=0}$. 
With the help of the above $\theta^{\d}$, then we can extend the elliptic problem \eqref{17.35}-\eqref{17.36} from $(0,L)\times \mathbb{R}_+$ to $(0,L)\times \mathbb{R}$.  

\begin{proposition}[Extension lemma]\label{lemAF}
There exist a unique solution to \eqref{17.26} satisfying
\begin{align}\label{17.27}
		\begin{split}
			\|(1+y)^l (\sqrt{\v}\theta^{\d}_x, \theta^{\d}_y, \theta^{\d})\|^2 
			&\lesssim C_l \Big(\|\Phi_{yyy}\|^2 + \|||\Phi||\|^2_1\Big),\\
			\|u_s\, (1+y)^l (\v\theta^{\d}_{xx}, \sqrt{\v}\theta^{\d}_{xy}, \theta^{\d}_{yy})\|^2 &\lesssim C_{l} a^2 \Big(\|\Phi_{yyy}\|^2 + \|||\Phi||\|^2_1\Big),\\
			\delta \|\nabla_{\v}^2\theta^{\d}\|^2 &\lesssim C_{\v}\Big(\delta \|\Phi_{yyy}\| +  \|||\Phi||\|^2_1\Big).
		\end{split}
\end{align}
\end{proposition}

\noindent{\bf Proof.} In the following, we divide the proof into two steps. 

{\it Step 1. Existence.} Noting the construction of $\f{\mu}{\tilde{\rho}^{\v,\d}} \tilde{\Phi}^{\delta}_{yyy}$, and $\theta^{\d}|_{x=0}=\theta^{\d}_x|_{x=L}=0$, we can take odd extension of $\theta^{\d}$ so that \eqref{17.26} is defined in $(-L,L)\times \mathbb{R}_+$, hence it becomes a $x$-periodic problem.  We consider the solution in the following form
\begin{align}\label{AF3.2}
	\theta^{\d}(x,y) = \sum_{i\geq1} \theta_i(y) \sin(\f{i}{2L}\pi x).
\end{align} 
We remark that the summation is only for  odd $i$ so that 
$$\big(\sin(\f{i}{2L}\pi x)\big)'\Big|_{x=L}=\cos(\f{i}{2}\pi)=0.$$

For simplicity, we denote 
\begin{align}\label{AF3.3}
	\dis f(x):= \f{\mu}{\tilde{\rho}^{\v,\d}} \tilde{\Phi}^{\delta}_{yyy}\big|_{y=0} \quad \mbox{and}\quad 	f(x)=\sum_{i} c_{i}(f) \sin(\f{i\pi}{2L}x),
\end{align}
then it is clear that 
\begin{align}\label{AF3.4}
	\|f\|^2=\f{1}{L}\sum_{i}|c_{i}(f)|^2,\quad \|f_x\|^2 =\f{1}{L}\sum_{i} |c_{i}(f)|^2 (\f{i\pi}{2L})^2,
\end{align}
where we have used the fact 
\begin{align*}
	\int_{-L}^{L} \sin^2(\f{i}{2L}\pi x)dx =\f{1}{L}.
\end{align*}
It follows from $\eqref{17.26}_1$ that
\begin{align}\label{AF3.5}
	\theta_{i}(0) =  c_{i}(f).
\end{align}
We also assume $\theta_{i}(\infty)=0$.

Now we obtain the following ODEs
\begin{align}\label{AF3.6}
	\begin{cases}
		\dis \theta_{i}''(y) - \big\{1+ \v [\f{i\pi}{2L}]^2\} \theta_i(y) =0,\\
		\dis \theta_{i}(0) =  c_{i}(f)\quad \mbox{and}\quad 
		\theta_{i}(\infty)=0.
	\end{cases}
\end{align}
It is easy to solve above problem
\begin{align}\label{AF3.7}
	\theta_{i}(y)=\theta_i(0) \exp \Big(-y\sqrt{1+ \v [\f{i\pi}{2L}]^2}\Big) = c_{i}(f) \exp \Big(-y\sqrt{1+ \v [\f{i\pi}{2L}]^2}\Big),
\end{align}
which yields that 
\begin{align}\label{AF3.8}
	\theta_{i}'(y)  = - \sqrt{1+ \v [\f{i\pi}{2L}]^2} c_{i}(f)  \exp \Big(-y\sqrt{1+ \v [\f{i\pi}{2L}]^2}\Big).
\end{align}

\

Multiplying $\eqref{AF3.6}_1$ by $\theta_i(y)$, one obtains that 
\begin{align*}
	|\theta_i'|^2_{L^2_y} + \big\{1+\v [\f{i\pi}{2L}]^2\big\} |\theta_i|_{L^2_y} = -\theta_i'(0) \theta_i(0) = \sqrt{1+ \v [\f{i\pi}{2L}]^2} |c_{i}(f)|^2 ,
\end{align*}
which yields that 
\begin{align}\label{AF3.10}
	\|\theta^{\d}_y\|^2 + \|\sqrt{\v}\theta^{\d}_x\|^2 + \|\theta^{\d}\|^2
	&= \f{1}{L} \sum_i \Big\{ |\theta_i'|^2_{L^2_y} + \v [\f{i\pi}{2L}]^2 |\theta_i|_{L^2_y} \Big\} 
	=\f{1}{L} \sqrt{1+ \v [\f{i\pi}{2L}]^2} |c_{i}(f)|^2.
\end{align}
A direct calculation shows that 
\begin{align}\label{AF3.11}
	\f{1}{L} \sum_i \sqrt{1+ \v [\f{i\pi}{2L}]^2}  |c_{i}(f)|^2 
	&\lesssim \f{1}{L} \sum_i \big\{1+\sqrt{\v} \f{i\pi}{2L}\big\}  |c_{i}(f)|^2 \lesssim \|f\|^2   + \sqrt{\v}\|f\|_{H^{1/2}}^2 \nonumber\\
	&\lesssim \|\Phi_{yyy}\|^2 +\|(\Phi_{yyyy}, \sqrt{\v}\Phi_{xyyy})\|^2.
\end{align}
Hence it holds from \eqref{AF3.10}-\eqref{AF3.11} that 
\begin{align}\label{AF3.13}
	\|(\theta^{\d}, \theta^{\d}_y,\sqrt{\v}\theta^{\d}_x)\|^2  \lesssim \|\Phi_{yyy}\|^2 +\|(\Phi_{yyyy}, \sqrt{\v}\Phi_{xyyy})\|^2.
\end{align}
Therefore, if $f\in H^{1/2}_{y=0}$,  we establish the existence of  a unique weak solution $\theta^{\d}\in H^1$ of \eqref{17.26}.

\smallskip

For $\dis f \in H^{3/2}_{y=0}$,  we know $\theta^{\d}\in H^2$ by the regularity theory of elliptic equation with
\begin{align}\label{AF.41}
	\|\nabla^2_{\v}\theta^{\d}\|\leq C_{\v} \|\tilde{\Phi}^{\delta}_{yyy}\|^2_{H^{3/2}}\leq C_{\v} \Big(\|\Phi_{yyy}\| + \f{1}{\sqrt{\delta}} \|||\Phi||\|_1\Big).
\end{align}

\medskip

{\it Step 2. Weighted estimate.} Multiplying \eqref{17.26} by $y^2\theta^{\d}$, one has 
\begin{align}\label{AF3.39}
	\|y\theta^{\d}\|^2 + \|y(\sqrt{\v}\theta^{\d}_x,\theta^{\d}_y)\|^2   \lesssim \|\theta^{\d}\|^2 
	\lesssim \|\Phi_{yyy}\|^2 +\|(\Phi_{yyyy}, \sqrt{\v}\Phi_{xyyy})\|^2,
\end{align}
where we have used \eqref{AF3.13} and the following two estimates
\begin{align*}
	\begin{split}
		- \la \v \theta^{\d}_{xx},\, y^2 \theta^{\d} \ra 
		&= \|\sqrt{\v}y\theta^{\d}_x\|^2,\\
		-  \la \theta^{\d}_{yy},\, y^2 \theta^{\d} \ra& =\|y\theta^{\d}_y\|^2 + 2\la \theta^{\d}_{y},\, y \theta^{\d} \ra\geq \f12 \|y\theta^{\d}_y\|^2 - C\|\theta^{\d}\|^2.
	\end{split}
\end{align*}

\medskip

Multiplying \eqref{17.26} by $y^{2l} \theta^{\d}$, one has
\begin{align}\label{AF.42}
	\|y^l \theta^{\d}\|^2 - \la \v \theta^{\d}_{xx},\, y^{2l} \theta^{\d} \ra -  \la \theta^{\d}_{yy},\, y^{2l} \theta^{\d} \ra =0.
\end{align}
A direct calculation shows 
\begin{align*}
	\begin{split}
		- \la \v \theta^{\d}_{xx},\, y^{2l} \theta^{\d} \ra 
		&= \|\sqrt{\v}y^{l}\theta^{\d}_x\|^2,\\
		-  \la \theta^{\d}_{yy},\, y^{2l} \theta^{\d} \ra&=\|y^{l}\theta^{\d}_y\|^2 + 2\la \theta^{\d}_{y},\, 2ly^{2l-1} \theta^{\d} \ra\geq \f12 \|y^{l}\theta^{\d}_y\|^2 - Cl^2\|y^{l-1}\theta^{\d}\|^2,
	\end{split}
\end{align*}
which, together with \eqref{AF.42}, yields that 
\begin{align}\label{AF.43}
	\|y^{l}\theta^{\d}\|^2 + \|(\theta^{\d}_y, \sqrt{\v}\theta^{\d}_x)y^{l}\|^2 \leq Cl^2\|y^{l-1}\theta^{\d}\|^2.
\end{align}
By induction arguments on $l$, one obtains 
\begin{align}\label{AF.45-1}
	\|y^l\theta^{\d}\|^2 + \|y^l (\theta^{\d}_y, \sqrt{\v}\theta^{\d}_x)\|^2 &\leq C_{l}\|y^{l-1}\theta^{\d}\|^2\leq \cdots \leq C_{l} \|\theta^{\d}\|^2\nonumber\\
	&\leq C_{l} \big\{ \|\Phi_{yyy}\|^2 +\|(\Phi_{yyyy}, \sqrt{\v}\Phi_{xyyy})\|^2\big\}.
\end{align}

It follows from \eqref{17.26} that 
\begin{align}\label{AF.45}
	\|y^{l} \v \theta^{\d}_{xx}\|^2 + \|y^{l}\theta^{\d}_{yy}\|^2 + \la y^{2l} \v \theta^{\d}_{xx}, \theta^{\d}_{yy} \ra = \|y^{l} \theta^{\d}\|^2.
\end{align}
Integrating by parts, one gets 
\begin{align*}
	\la y^{2l} \v \theta^{\d}_{xx}, \theta^{\d}_{yy} \ra & = - \la y^{2l} \v \theta^{\d}_{x}, \theta^{\d}_{xyy} \ra  =\|\sqrt{\v} \theta^{\d}_{xy}y^{l}\|^2 -  2l\la y^{2l-1} \v \theta^{\d}_{x}, \theta^{\d}_{xy} \ra\nonumber\\
	&\geq \f12 \|\sqrt{\v} \theta^{\d}_{xy}y^{l}\|^2 - Cl^2\|\sqrt{\v}\theta^{\d}_x y^{l-1}\|^2,
\end{align*}
which, together with \eqref{AF.45}, yields that 
\begin{align}\label{AF.46}
	\|y^{l} (\v \theta^{\d}_{xx}, \theta^{\d}_{yy}, \sqrt{\v}\theta^{\d}_{xy})\|^2 
	&\lesssim \|\Phi_{yyy}\|^2 +\|(\Phi_{yyyy}, \sqrt{\v}\Phi_{xyyy})\|^2.
\end{align}

Combining \eqref{AF3.13}-\eqref{AF.41}, \eqref{AF.45-1} and \eqref{AF.46}, one concludes \eqref{17.27}. Therefore the proof of Proposition \ref{lemAF} is completed. $\hfill\Box$


\smallskip
 
\subsubsection{Approximations and extension} Recall \eqref{17.14}. As an approximation of $B$, we define
\begin{align}\label{17.35}
	B^{\delta} :&= P_x^{\delta}  + \bar{g}_1(p^{\delta},\phi, T(S^{\delta},p^{\delta},q))\equiv P_x^{\delta} + \bar{g}_1^{\delta}.
\end{align}
Then we consider the BVP
\begin{align}\label{17.36-0}
-\Delta_{\v}B^{\delta} = -\mathbf{w}^{\delta}_x - \Delta_{\v}\bar{g}^{\delta}_1,\,\, (x,y) \in (0,L)\times \mathbb{R}_+,
\end{align}
with the  $\delta$-approximation boundary conditions 
\begin{align}\label{17.36}
	\begin{cases}
		\dis B^{\delta} \big|_{y=0} = \f{\mu}{\tilde{\rho}^{\v,\d}} \tilde{\Phi}^{\delta}_{yyy},\\
		\dis \big\{ (\sigma-1) p_s \mathbf{d}_{11} B^{\delta}  - \v\frac{1}{2}\lambda  [1+\bm{\chi}] u_s^{\delta} B^{\delta}_{x} \big\} \big|_{x=0} = -\frac{1}{2}\lambda  [1+\bm{\chi}] \big(u_s^{\delta} \mathbf{w}^{\delta} + \delta\v \theta^{\d}_x\big) + \mathbf{g}^{\delta}_{b1},\\
		\dis \big\{\v B^{\delta}_{x} - \lambda\v\frac{b_2}{2p_s} u_s^{\delta} B^{\delta}_{yy}\big\} \big|_{x=L}  =\mathbf{w}^{\delta} - \lambda \delta \v \f{b_2}{2p_s}  \theta^{\d}_{yy} + \mathbf{g}^{\delta}_{b2},
	\end{cases}
\end{align}
where we have denoted 
\begin{align}\label{17.31}
	\begin{split}
		\mathbf{g}^{\delta}_{b1}:&=\mathbf{g}_{b1}(\hat{P},\Phi,p^{\delta},S^{\delta}, \zeta^{\delta},\phi),\\
		\mathbf{g}^{\delta}_{b2}:&=  \mathbf{g}_{b2}(\hat{P}, \hat{\mathbf{S}}, \Phi,p^{\delta},S^{\delta}, \zeta^{\delta},\phi),
	\end{split}
\end{align}
which means that we lift all functions $(p,S,\zeta)$ in $\mathbf{g}_{b1},\mathbf{g}_{b2}$ to their corresponding 
$\delta$-dependent versions. Similar notational conventions for such lifted quantities may also be used elsewhere without further detailed definition.

Since \eqref{17.36-0}-\eqref{17.36} is not a standard elliptic problem,  it remains unclear how to solve it directly. To address the lower boundary, we introduce a new function
\begin{align}\label{17.37}
\bar{B}^{\delta}:=B^{\delta} - \theta^{\d}.
\end{align}
where $\theta^{\d}$ is the solution of \eqref{17.26} established in Proposition \ref{lemAF}, and it will help us to handle the lower boundary.
Thus we can rewrite \eqref{17.35}-\eqref{17.36} as
\begin{align} \label{17.38}
	\begin{split}
		-\Delta_{\v}\bar{B}^{\delta} = -\mathbf{w}^{\delta}_x - \Delta_{\v}\bar{g}^{\delta}_1 + \theta^{\d}, \,\, (x,y) \in (0,L)\times \mathbb{R}_+,
	\end{split}
\end{align}
with boundary conditions
\begin{align}\label{17.39}
	\begin{cases}
		\dis \bar{B}^{\delta} \big|_{y=0}= 0,\\
		\dis \big\{ (\sigma-1) p_s \mathbf{d}_{11} \bar{B}^{\delta}  - \v\frac{1}{2}\lambda  [1+\bm{\chi}] u_s^{\delta} \bar{B}^{\delta}_{x} \big\} \big|_{x=0} = -\frac{1}{2}\lambda  [1+\bm{\chi}] \big(u_s^{\delta} \mathbf{w}^{\delta} - \v u_s \theta^{\d}_x\big) + \mathbf{g}^{\delta}_{b1},\\[2mm]
		\dis \big\{\v \bar{B}^{\delta}_{x} - \lambda\v\frac{b_2}{2p_s} u_s^{\delta} \bar{B}^{\delta}_{yy}\big\} \big|_{x=L}  = \mathbf{w}^{\delta} + \lambda  \v \f{b_2}{2p_s} u_s  \theta^{\d}_{yy} + \mathbf{g}^{\delta}_{b2}.
	\end{cases}
\end{align}

\

However, solving \eqref{17.38}-\eqref{17.39} remains nontrivial. To construct $B$ with sufficient regularity, we extend the original problem from $(0,L)\times \mathbb{R}_+$ to $(0,L)\times \mathbb{R}$, then a further approximation is necessary. 
For strictness, we further introduce  an $\gamma$-approximate problem (cut-off and smooth technique) with $\gamma=\gamma(\delta)=\sqrt{\delta}\to 0+$ as $\delta\to0$, i.e.,   
\begin{align}\label{17.40}
	-\Delta_{\v} \bar{B}^{\gamma,\delta} = - \Big(\mathbf{w}^{\delta}\bar{\chi}(\f{y}{\gamma})\Big)_x - \bar{\chi}(\f{y}{\gamma})\Delta_{\v}\bar{g}^{\delta}_1 + \bar{\chi}(\f{y}{\gamma}) \theta^{\d}, \,\, (x,y) \in (0,L)\times \mathbb{R}_+,
\end{align}
with the $\gamma$-approximate boundary conditions
\begin{align}\label{17.41}
\begin{cases}
	\dis \bar{B}^{\gamma,\d} \big|_{y=0}= 0,\\
	\dis \big\{ (\sigma-1) p_s \mathbf{d}_{11} \bar{B}^{\gamma,\d}  - \v\frac{1}{2}\lambda  [1+\bm{\chi}] u_s^{\delta} \bar{B}^{\gamma,\d}_{x} \big\} \big|_{x=0} = -\frac{1}{2}\lambda  [1+\bm{\chi}] \big(u_s^{\delta} \mathbf{w}^{\delta} - \v u_s \theta^{\d}_x\big)\bar{\chi}(\f{y}{\gamma}) + \mathbf{g}^{\gamma,\delta}_{b1},\\[2mm]
	\dis \big\{\v \bar{B}^{\gamma,\d}_{x} - \lambda\v\frac{b_2}{2p_s} u_s^{\delta} \bar{B}^{\gamma,\d}_{yy} \big\} \big|_{x=L}  = \mathbf{w}^{\delta}\bar{\chi}(\f{y}{\gamma}) + \lambda  \v \f{b_2}{2p_s} u_s  \theta^{\d}_{yy}\bar{\chi}(\f{y}{\gamma}) + \mathbf{g}^{\gamma,\delta}_{b2},
\end{cases}
\end{align}
where 
\begin{align}\label{17.51}
	\begin{split}
		\mathbf{g}_{b1}^{\gamma,\d}:&= \bar{\chi}(\f{y}{\gamma})\, \mathbf{g}^{\delta}_{b1},\\
		\mathbf{g}_{b2}^{\gamma,\d}:&=\bar{\chi}(\f{y}{\gamma})\Big[\lambda\v \big(\frac{b_2}{2p_s} u_s\big)_y \hat{P}^{\d}_{xy}  + \v \bar{g}^{\delta}_{1x} + \big(b_2 [ g^{\delta}_2+ \v \mathfrak{g}_2 - \v\mathfrak{h}^{\delta}_2 - \v\mathfrak{N}^{\delta}_{21}]\big)_y\Big] \\
		&\quad    - \Big(\bar{\chi}(\f{y}{\gamma})\hat{G}^{\delta}_3\Big)_y  - \lambda\v \f{b_2}{2p_s} u_s  \Big(\bar{\chi}(\f{y}{\gamma})\bar{g}^{\delta}_{1y}\Big)_y,\\
		\hat{G}^{\delta}_3:&= G_3(\Phi, \hat{P}^{\d}, \hat{\mathbf{S}}^{\d},  p^{\delta}, \zeta^{\delta}, \phi).
	\end{split}
\end{align}

Noting $\eqref{17.41}_1$,  we consider an odd extension in $y$ for   \eqref{17.40}-\eqref{17.41}, and define
\begin{align}\label{17.52-0}
	\tilde{B}^{\delta}(x,y):=
	\begin{cases}
	\dis	\bar{B}_1^{\gamma,\d}(x,y),\qquad  y\geq 0,\\
	\dis	-\bar{B}_1^{\gamma,\d}(x,-y),\,\,\,\,  y<0,
	\end{cases}
\&\quad 
	\tilde{\mathbf{w}}^{\delta}(x,y):=
	\begin{cases}
	\dis	\big(\mathbf{w}^{\delta}\bar{\chi}(\f{y}{\gamma})\big)(x,y),\qquad y\geq 0,\\
	\dis	-\big(\mathbf{w}^{\delta}\bar{\chi}(\f{y}{\gamma})\big)(x,-y),\,\,\,\,  y< 0,
	\end{cases}
\end{align}

\begin{align}\label{17.57-1}
	\tilde{\theta}^{\delta}(x,y):=
	\begin{cases}
	\dis	\big(\theta^{\d}\bar{\chi}(\f{y}{\gamma})\big)(x,y),\quad y\geq 0,\\
	\dis	-\big(\theta^{\d}\bar{\chi}(\f{y}{\gamma})\big)(x,-y),\,\, y< 0,
	\end{cases}
\& \quad 
	\tilde{\Xi}^{\delta}(x,y):=
	\begin{cases}
	\dis	\big(\theta^{\d}_{yy}\bar{\chi}(\f{y}{\gamma})\big)(x,y),\quad y\geq 0,\\
	\dis	-\big(\theta^{\d}_{yy}\bar{\chi}(\f{y}{\gamma})\big)(x,-y),\,\,  y< 0,
	\end{cases}
\end{align}

\begin{align}\label{17.57}
	\tilde{\mathbf{g}}^{\delta}_{b1}(x,y):=
	\begin{cases}
	\dis	 \mathbf{g}^{\gamma,\delta}_{b1}(x,y),\qquad y\geq 0,\\
	\dis	- \mathbf{g}_{b1}^{\gamma,\delta}(x,-y),\,\,\,\,  y< 0,
	\end{cases}
\&	\quad 
 \tilde{\mathbf{g}}^{\delta}_{b2}(x,y):=
	\begin{cases}
	\dis \mathbf{g}_{b2}^{\gamma,\delta}(x,y),\,\, y\geq 0\\
    \dis - \mathbf{g}_{b2}^{\gamma,\delta}(x,-y),\,\, y< 0,
	\end{cases}
\end{align}
and 
\begin{align}\label{17.48}
	\bar{\mathbf{g}}^{\delta}_1(x,y):=
	\begin{cases}
		\dis	\big(\bar{\chi}(\f{y}{\gamma})\Delta_{\v}\bar{g}^{\delta}_1\big)(x,y),\qquad y\geq 0,\\
		\dis	-\big(\bar{\chi}(\f{y}{\gamma})\Delta_{\v}\bar{g}^{\delta}_1\big)(x,-y),\,\,\,\, y< 0.
	\end{cases}
\end{align}
For the coefficient $u_s$, we apply an even extension
\begin{align}
	\tilde{u}_s(x,y):=
	\begin{cases}
		u_s(x,y),\quad y\geq 0\\
		u_s(x,-y),\,\, y< 0,
	\end{cases}
	\quad \mbox{and}\quad 
	\tilde{u}^{\delta}_{s}:=\delta + \tilde{u}_s.
\end{align}
Analogous extensions are applied to other relevant terms, the details of which are omitted for conciseness.

\medskip

Therefore, to construct the approximate pressure  $P^{\delta}$, we turn to the following approximate BVP:
\begin{align}\label{17.50}
\begin{cases}
\dis -\Delta_{\v} \tilde{B}^{\d} = - \tilde{\mathbf{w}}^{\d}_x - \bar{\mathbf{g}}^{\d}_1 +\tilde{\theta}^{\delta},\quad (x,y)\in (0,L)\times \mathbb{R},\\[1mm]
\dis \Big\{ (\sigma-1) p_s \mathbf{d}_{11} \tilde{B}^{\d}  - \v\frac{1}{2}\lambda  [1+\bm{\chi}] \tilde{u}_s^{\delta} \tilde{B}^{\d}_{x} \Big\} \Big|_{x=0} = -\frac{1}{2}\lambda  [1+\bm{\chi}] \big(\tilde{u}_s^{\delta} \tilde{\mathbf{w}}^{\d} - \v \tilde{u}_s \tilde{\theta}^{\delta}_x\big) + \tilde{\mathbf{g}}^{\d}_{b1},\\[2mm]
\dis \Big\{ \v \tilde{B}^{\d}_{x} - \lambda\v\frac{b_2}{2p_s} \tilde{u}_s^{\delta} \bar{B}^{\d}_{yy}\Big\} \Big|_{x=L}  = \tilde{\mathbf{w}}^{\d} + \lambda  \v \f{b_2}{2p_s} \tilde{u}_s \tilde{\Xi}^{\delta} + \tilde{\mathbf{g}}^{\d}_{b2}.
\end{cases}
\end{align}
We can further smooth the boundary data if necessary.

\subsubsection{Weak solution of \eqref{17.50}.} For later use, we define $\mathcal{H}^1$ as:
\begin{align}\label{17.124-1}
\mathcal{H}^1=\big\{ \varphi\in H^1([0,L]\times\mathbb{R}) \, : \,  \varphi_y|_{x=L}\in L^2(\mathbb{R})\big\},\quad \|\varphi\|_{\mathcal{H}^1}^2 : = \|\varphi\|_{H^1}^2 + \|\varphi_y\|_{x=L}^2.
\end{align}
\begin{lemma}\label{Blem7.4}
There exists a unique weak solution $\tilde{B}^{\delta} \in \mathcal{H}^1$ to \eqref{17.50}. Moreover, $\tilde{B}^{\delta}$ is odd with respect to $y \in \mathbb{R}$.
\end{lemma}

\noindent{\bf Proof.} Multiplying \eqref{17.50} by smooth function $\varphi$, one has
\begin{align}\label{17.93}
	&\int_0^L\int_{\mathbb{R}} \nabla_{\v}\tilde{B}^{\d}\cdot \nabla_{\v}\varphi dydx - \int_{\mathbb{R}} \v  \tilde{B}^{\d}_{x} \varphi dy \Big|_{x=L} + \int_{\mathbb{R}} \v  \tilde{B}^{\d}_{x} \varphi dy \Big|_{x=0}  \nonumber\\
	&=\int_0^L\int_{\mathbb{R}} \tilde{\mathbf{w}}^{\d} \varphi_x dydx - \int_{\mathbb{R}}\tilde{\mathbf{w}}^{\d}  \varphi dy \Big|_{x=L} + \int_{\mathbb{R}} \tilde{\mathbf{w}}^{\d}  \varphi dy \Big|_{x=0} - \int_0^L\int_{\mathbb{R}}\big(\bar{\mathbf{g}}^{\d}_1 -\tilde{\theta}^{\delta}\big) \varphi dydx.
\end{align}
For the boundary terms at $x=0,L$ in \eqref{17.93}, we have
\begin{align}\label{17.94}
&\int_{\mathbb{R}} \{\v\tilde{B}^{\d}_{x} - \tilde{\mathbf{w}}^{\d}\}  \varphi dy \Big|_{x=L} =\int_{\mathbb{R}}  \Big\{\lambda \v \f{b_2}{2p_s} \tilde{u}_s^{\delta} \tilde{B}^{\d}_{yy}  + \lambda\v\f{b_2}{2p_s} \tilde{u}_s \tilde{\Xi}^{\delta} + \tilde{\mathbf{g}}^{\d}_{b2}\Big\}\, \varphi dy \Big|_{x=L} \nonumber\\
&= - \lambda \v\int_{\mathbb{R}}   \tilde{B}^{\d}_{y} \Big(\f{b_2}{2p_s}\tilde{u}_s^{\delta} \varphi\Big)_y  dy \Big|_{x=L} + \int_{\mathbb{R}}  \Big\{\lambda\v \f{b_2}{2p_s}\tilde{u}_s \tilde{\Xi}^{\delta} + \tilde{\mathbf{g}}^{\d}_{b2} \Big\}\, \varphi dy \Big|_{x=L},
\end{align}
and
\begin{align}\label{17.95}
&\int_{\mathbb{R}} \big\{ \v  \tilde{B}^{\d}_{x}  -  \tilde{\mathbf{w}}^{\d} \big\} \varphi dy \Big|_{x=0} \nonumber\\
&= \int_{\mathbb{R}} \f{2}{\lambda [1+\bm{\chi}]} \f{1}{\tilde{u}^{\delta}_s} \Big\{(\sigma-1)p_s \mathbf{d}_{11}\tilde{B}^{\d}  - \f12  \lambda \v [1+\bm{\chi}] \tilde{u}_s \tilde{\theta}^{\delta}_x - \tilde{\mathbf{g}}^{\d}_{b1}\Big\}\varphi dy \Big|_{x=0}.
\end{align}

Combining\eqref{17.93}-\eqref{17.95}, we can express the weak formulation of \eqref{17.50} as  
\begin{align}\label{17.96}
	\mathfrak{B}^{\delta}[\tilde{B}^{\d},\varphi] &= \iint \tilde{\mathbf{w}}^{\d} \varphi_x dydx - \iint  \big(\bar{\mathbf{g}}^{\d}_1 -\tilde{\theta}^{\delta}\big) \varphi dydx + \int_{\mathbb{R}}  \Big\{\lambda\v \f{b_2}{2p_s}\tilde{u}_s \tilde{\Xi}^{\delta} + \tilde{\mathbf{g}}^{\d}_{b2} \Big\}\, \varphi dy \Big|_{x=L} \nonumber\\
	&\quad  + \int_{\mathbb{R}} \f{2}{\lambda [1+\bm{\chi}]} \f{\varphi}{\tilde{u}^{\delta}_s} \Big\{ \f12  \lambda \v [1+\bm{\chi}] \tilde{u}_s \tilde{\theta}^{\delta}_x + \tilde{\mathbf{g}}^{\d}_{b1}\Big\} dy \Big|_{x=0},
\end{align}
where we have used the following notation 
\begin{align}\label{17.125-1}
	\mathfrak{B}^{\delta}[\tilde{B}^{\d},\varphi]:&=\iint \nabla_{\v}\tilde{B}^{\d}\cdot \nabla_{\v}\varphi dydx + \lambda \v\int_{\mathbb{R}}   \tilde{B}^{\d}_{y} \,  \Big(\f{b_2}{2p_s}\tilde{u}_s^{\delta} \varphi\Big)_y  dy \Big|_{x=L} \nonumber\\
	&\quad   +  \int_{\mathbb{R}}  p_s \mathbf{d}_{11} \f{2(\sigma-1)}{\lambda[1+\bm{\chi}]}\f{1}{\tilde{u}^{\delta}_s} \tilde{B}^{\d}  \varphi dy \Big|_{x=0}.
\end{align}

Noting $\sigma>1$ and $0<\v\ll 1$, we have 
\begin{align}\label{17.66}
\mathfrak{B}^{\delta}[\tilde{B}^{\d}, \tilde{B}^{\d}]
&\geq  \f34\|\nabla_{\v}\tilde{B}^{\d}\|^2+ \f34\lambda \v\int_{\mathbb{R}} \f{b_2}{2p_s}\tilde{u}_s^{\delta}  |\tilde{B}^{\d}_{y}|^2  dy \Big|_{x=L} +  \f34\int_{\mathbb{R}} \f{2(\sigma-1)}{\lambda[1+\bm{\chi}]} \f{ p_s \mathbf{d}_{11}}{\tilde{u}^{\delta}_s} |\tilde{B}^{\d}|^2 dy \Big|_{x=0} .
\end{align}
It is direct to know that 
\begin{align}\label{17.68-1}
\|\tilde{B}^{\d}w\|&\lesssim \sqrt{L} \|\tilde{B}^{\d}w\|_{x=0} + L  \|\tilde{B}^{\d}_{x}w\| \lesssim \sqrt{La} \|\f{\tilde{B}^{\d}w}{\sqrt{\tilde{u}_s^{\delta}}}\|_{x=0} + L  \|\tilde{B}^{\d}_{x}w\| .
\end{align}
which, together with \eqref{17.66}, yields that 
\begin{align}\label{17.78-0}
	\mathfrak{B}^{\delta}[\tilde{B}^{\d}, \tilde{B}^{\d}]
	&\geq c_1(\v,\delta) \|\tilde{B}^{\d}\|^2_{\mathcal{H}^1}>0,
\end{align} 
where the positive constant $c_1(\v,\delta)$ depends on $\v,\delta$.

We control the terms on RHS of \eqref{17.96}. We integrate by parts in $y$ to obtain
\begin{align}\label{17.73}
& \lambda\v \int_{\mathbb{R}} \f{b_2}{2p_s} \tilde{u}_s \tilde{\Xi}^{\delta} \varphi dy \Big|_{x=L}
= \lambda\v \int_{0}^{\infty} \f{b_2u_s}{2p_s}   \theta^{\d}_{yy}\bar{\chi}(\f{y}{\gamma}) \varphi dy \Big|_{x=L} + \int_{-\infty}^0 (\cdots) dy \Big|_{x=L}  \nonumber\\
&=- \lambda\v \int_{0}^{\infty} \theta^{\d}_{y} \Big\{\big[\f{b_2u_s}{2p_s} \varphi\big]_y\bar{\chi}(\f{y}{\gamma}) + \f{b_2 u_s}{2p_s}  \varphi  \, \f{1}{\gamma}\bar{\chi}'(\f{y}{\gamma})\Big\} dy \Big|_{x=L}   + \int_{-\infty}^0(\cdots) dy \Big|_{x=L} \nonumber\\
&\lesssim  \v \|\theta^{\d}_y\|_{x=L} \big\{\|\varphi_y\|_{x=L} + \|\varphi\|_{x=L}\big\},
\end{align}
and
\begin{align}
\int_{\mathbb{R}} \tilde{\mathbf{g}}^{\d}_{b2} \varphi dy \Big|_{x=L} 
&=  \int_{0}^{\infty} \bar{\chi}(\f{y}{\gamma})\Big[\lambda\v \big(\frac{b_2u_s}{2p_s} \big)_y \hat{P}_{xy}  + \v \bar{g}^{\delta}_{1x} + \big(b_2 [ g^{\delta}_2+ \v \mathfrak{g}_2 - \v\mathfrak{h}^{\delta}_2 - \v\mathfrak{N}^{\delta}_{21}]\big)_y\Big] \, \varphi dy \Big|_{x=L} \nonumber\\
&\quad+  \int_{0}^{\infty} \Big\{\lambda \v \big(\f{b_2u_s}{2p_s}  \varphi\big)_y  \bar{\chi}(\f{y}{\gamma})\bar{g}^{\delta}_{1y} + \bar{\chi}(\f{y}{\gamma})\, \hat{G}^{\delta}_3\, \varphi_y\Big\} dy \Big|_{x=L}   + \int_{-\infty}^0 (\cdots) dy \Big|_{x=L} \nonumber\\
&\lesssim \Big\{\|(\v\hat{P}_{xy},\v\bar{g}^{\delta}_{1x}, \v \bar{g}^{\delta}_{1y}, \hat{G}^{\delta}_3)\|_{x=L} + \|\big(b_2 [ g^{\delta}_2+ \v \mathfrak{g}_2 - \v\mathfrak{h}^{\delta}_2 - \v\mathfrak{N}^{\delta}_{21}]\big)_y\|_{x=L} \Big\}\nonumber\\
&\quad \times\Big(\|\varphi\|_{x=L} + \|\varphi_y\|_{x=L}\Big).
\end{align}
A direct calculation shows that 
\begin{align}
	\begin{split}
		\left|\iint \nabla_{\v}\tilde{B}^{\d}\cdot \nabla_{\v}\varphi dydx \right| &\lesssim \|\nabla_{\v}\tilde{B}^{\d}\|\cdot \|\nabla_{\v}\varphi\|,\\
		\left|\iint \tilde{\mathbf{w}}^{\d} \varphi_x dydx\right| &\lesssim \|\tilde{\mathbf{w}}^{\d}\|\cdot \|\varphi_x\|,\\
		\left|\iint  \big(\bar{\mathbf{g}}^{\d}_1 -\tilde{\theta}^{\delta}\big) \varphi dydx \right|
		&\lesssim \|(\Delta_{\v}\bar{g}^{\delta}_1,\theta^{\d})\chi(\f{y}{\gamma})\|\cdot \|\varphi\|.
	\end{split}
\end{align}
For the boundary terms, we have 
\begin{align}
\begin{split}
\lambda \v\int_{\mathbb{R}}   \tilde{B}^{\d}_{y} \,  \big(\f{b_2}{2p_s}\tilde{u}_s^{\delta} \varphi\big)_y  dy \Big|_{x=L} &\lesssim   \|\tilde{B}^{\d}_{y}\|_{x=L} \big\{\|\varphi\|_{H^1} + \|\varphi_y\|_{x=L}\},\\
\int_{\mathbb{R}} p_s \mathbf{d}_{11}\f{2(\sigma-1)}{\lambda[1+\bm{\chi}]}  \f{1}{\tilde{u}^{\delta}_s} \tilde{B}^{\d}  \varphi dy \Big|_{x=0} 
&\lesssim C_{\delta} \|\tilde{B}^{\d}\|_{H^1}\cdot \|\varphi\|_{H^1},
\end{split}
\end{align}
\begin{align}\label{17.78-1}
\f{2}{\lambda} \int_{\mathbb{R}} \Big\{\lambda \v\tilde{u}_s \tilde{\theta}^{\delta}_x + \f{\tilde{\mathbf{g}}^{\d}_{b1}}{1+\bm{\chi}} \Big\}\f{\varphi}{\tilde{u}^{\delta}_s} dy \Big|_{x=0}
& \lesssim \Big\{\|\f{\mathbf{g}^{\delta}_{b1}}{\sqrt{u^{\delta}_s}}\|_{x=0} + \|\v\sqrt{u_s} \theta^{\d}_x\|_{x=0} \Big\}  \|\f{\varphi}{\sqrt{\tilde{u}^{\delta}_s}}\|_{x=0}.
\end{align}

Recall $\mathcal{H}^1$ in \eqref{17.124-1}, the coercivity of $\mathfrak{B}^{\delta}[\cdot,\cdot]$ in \eqref{17.78-0}, and the estimates \eqref{17.73}-\eqref{17.78-1}. Then, by applying the Lax-Milgram theorem, we can find a unique weak solution $\tilde{B}^{\d}$ of \eqref{17.50} in the sense of weak formula \eqref{17.96}. Due to above construction, it is clear that $\tilde{B}^{\d}$ is odd  in $y\in \mathbb{R}$. Therefore the proof of Lemma \ref{Blem7.4} is completed. $\hfill\Box$

\medskip


\subsubsection{Uniform a priori estimate on $\tilde{B}^{\d}$.} 
Although the construction in Lemma \ref{Blem7.4} provides the bound $\|\tilde{B}^{\d}\|_{\mathcal{H}^1} \lesssim C(1/\v, 1/\delta)$, this estimate is insufficient for taking the limit $\delta \to 0+$ due to its $\delta$-dependence. To address this, we will exploit the oddness of $\tilde{B}^{\d}$ to derive some uniform-in-$\delta$ estimates.

\begin{lemma}\label{Blem7.4-1}
	Let $\tilde{B}^{\d}$ be the solution constructed in Lemma \ref{Blem7.4}, then it holds that 
	\begin{align}\label{17.70-1}
		\v  \|\f{\tilde{B}^{\d}w}{\sqrt{\tilde{u}_s}}\|^2_{x=L}
		&\lesssim \sqrt{L} \|\nabla_{\v}\tilde{B}^{\d}w\|^2 + \v  \|\f{\tilde{B}^{\d}w}{\sqrt{\tilde{u}_s^{\delta}}}\|^2_{x=0} .
	\end{align}
We remark that \eqref{17.70-1} is useful for later uniform estimates.
\end{lemma}

\noindent{\bf Proof.} Recall the oddness of $\tilde{B}^{\d}$ on $y$.
Using $\eqref{A.5}_1$ and \eqref{17.68-1}, one has that 
\begin{align}\label{17.68}
	\begin{split}
		\|\f{\tilde{B}^{\d}}{\tilde{u}_s}w\|&\lesssim  \f{1}{a} \big\{ \|\tilde{B}^{\d}w\| + \|\f{\tilde{B}^{\d}}{y}w\| \big\} \lesssim \f{1}{a} \big\{ \|\tilde{B}^{\d}w\| + \|\tilde{B}^{\d}_yw\| \big\}\\
		&\lesssim \f{1}{a}\Big\{L^{\f12} \|\f{\tilde{B}^{\d}w}{\sqrt{\tilde{u}_s^{\delta}}}\|_{x=0} + L\|\tilde{B}^{\d}_{x}w\| + \|\tilde{B}^{\d}_{y}w\|\Big\},
	\end{split}
\end{align}
which yields that 
\begin{align}\label{17.69}
	\v  \|\f{\tilde{B}^{\d}w}{\sqrt{\tilde{u}_s^{\delta}}}\|^2_{x=L} 
	&\lesssim \v  \|\f{\tilde{B}^{\d}w}{\sqrt{\tilde{u}_s^{\delta}}}\|^2_{x=0} + \v \|\tilde{B}^{\d}_{x}w\|\cdot \|\f{\tilde{B}^{\d}}{\tilde{u}_s^{\delta}}w\|
	\lesssim \sqrt{L} \|\nabla_{\v}\tilde{B}^{\d}w\|^2 + \v  \|\f{\tilde{B}^{\d}w}{\sqrt{\tilde{u}_s^{\delta}}}\|^2_{x=0}.
\end{align}

It is clear that 
\begin{align}\label{17.70}
	\int_0^\infty \f{|\tilde{B}^{\d}|^2}{u_s} w^2 dy \Big|_{x=L} \lesssim \f{1}{a} \int_0^\infty \f{|\tilde{B}^{\d}|^2}{y} w^2 dy \Big|_{x=L} +  \|\f{\tilde{B}^{\d}w}{\sqrt{\tilde{u}_s^{\delta}}}\|^2_{x=L}.
\end{align}
Using \eqref{17.68} and $\eqref{A.5}_1$, we have 
\begin{align*}
	\f{1}{a}\int_0^\infty \f{|\tilde{B}^{\d}|^2}{y}w^2 dy \Big|_{x=L} 
	&\lesssim \f{1}{La} \int_0^L\int_0^\infty \f{|\tilde{B}^{\d}|^2}{y}w^2 dydx + \f{1}{a}\int_0^L\int_0^\infty \f{2}{y} |\tilde{B}^{\d}_{x}\tilde{B}^{\d}| w^2 dydx\nonumber\\
	&\lesssim \f{1}{La} \big\{\|\tilde{B}^{\d}_{y}w\| + \|\sqrt{\v}\tilde{B}^{\d}w\|\big\} \big\{ \|\tilde{B}^{\d}w\| +  L\|\tilde{B}^{\d}_{x}w\|\big\} \nonumber\\
	&\lesssim  \f{1}{La}\big\{\|\nabla_{\v}\tilde{B}^{\d}w\| + \sqrt{L\v} \|\f{\tilde{B}^{\d}w}{\sqrt{\tilde{u}_s^{\delta}}}\|_{x=0} \big\} \big\{\|\f{\tilde{B}^{\d}w}{\sqrt{\tilde{u}_s^{\delta}}}\|_{x=0} + L  \|\tilde{B}^{\d}_{x}w\|\big\},
\end{align*}
which, together with  \eqref{17.69}-\eqref{17.70}, concludes \eqref{17.70-1}
Therefore the proof of Lemma \ref{Blem7.4-1} is completed. $\hfill\Box$

\begin{lemma}\label{Blem7.5}
For the weak solution $\tilde{B}^{\d}\in \mathcal{H}^1$ of \eqref{17.50} constructed in Lemma \ref{Blem7.4}, we have the following uniform estimate
\begin{align}\label{17.107}
	&(1-\xi) \Big\{\|\nabla_{\v}\tilde{B}^{\d}w\|^2 + \lambda \v\int_{\mathbb{R}}  \f{b_2}{2p_s}\tilde{u}_s^{\delta} |\tilde{B}^{\d}_{y}|^2 w^2  dy \Big|_{x=L}  +  \int_{\mathbb{R}}  p_s \mathbf{d}_{11} \f{2(\sigma-1)}{\lambda[1+\bm{\chi}]}\f{|\tilde{B}^{\d}|^2}{\tilde{u}^{\delta}_s} w^2 dy \Big|_{x=0} \Big\}\nonumber\\
	&\leq  \|\tilde{\mathbf{w}}^{\d}w\|\cdot \|\tilde{B}^{\d}_xw\|   +  C_{\xi}\v \|\f{\tilde{B}^{\d}}{\sqrt{\tilde{u}_s^{\delta}}} w_y\|^2_{x=L} +  C_{\xi} \|\tilde{B}^{\d} w_y\|^2  + C_{\xi} (\cdots)_{1,w},
\end{align}
where we have denoted
\begin{align}\label{17.71-1}
	(\cdots)_{1,w}:&=  a^2 \v \|\sqrt{u_s} \hat{P}_{xy} w\|^2_{x=L} + \|\nabla_{\v}\bar{g}^{\d}_{1}w\|^2 +  \v^2 \|\sqrt{u_s} \bar{g}^{\d}_{1x}w\|_{x=0}^2   +  \v \|\sqrt{u_s} (\bar{g}^{\d}_{1x}, \bar{g}^{\d}_{1y})\hat{w}\|_{x=L}^2  \nonumber\\
	&\,\,\, +  \v \|\sqrt{u_s}(b_2\mathfrak{g}_2)_{y}w\|^2_{x=L}   + \f{1}{\v} \|\sqrt{u_s} (b_2g^{\d}_2)_{y}w\|^2_{x=L}   + \v \|\sqrt{u_s}(b_2\mathfrak{h}^{\delta}_2+b_2\mathfrak{N}^{\delta}_{21})_{y}w\|_{x=L} 
	 \nonumber\\
	&\,\,\, + \f{1}{\v}\|\f{\hat{G}^{\d}_3}{\sqrt{u^{\delta}_s}}w\|^2_{x=L} +   \|\f{\mathbf{g}^{\delta}_{b1} w}{\sqrt{u_s^{\delta}}} \|^2_{x=0}  +  \|y\theta^{\d} w\|^2  + \v \|\sqrt{u_s}(\sqrt{\v}\theta^{\d}_x, \theta^{\d}_y)\hat{w}\|^2_{x=0} .
\end{align}
\end{lemma}

\noindent{\bf Proof.} To derive the uniform in $\delta,\v$ estimate, we need to use the oddness of $\tilde{B}^{\d}$ in $y$. 
Taking $\varphi=\tilde{B}^{\d} w^2$ in \eqref{17.96}, one has
\begin{align}\label{17.52}
	&\|\nabla_{\v}\tilde{B}^{\d}w\|^2 + \lambda \v\int_{\mathbb{R}}  \f{b_2}{2p_s}\tilde{u}_s^{\delta} |\tilde{B}^{\d}_{y}|^2 w^2  dy \Big|_{x=L}  +  \int_{\mathbb{R}}  p_s \mathbf{d}_{11} \f{2(\sigma-1)}{\lambda[1+\bm{\chi}]}\f{|\tilde{B}^{\d}|^2}{\tilde{u}^{\delta}_s} w^2 dy \Big|_{x=0} \nonumber\\
	&= - 2\big\la \tilde{B}^{\d}_y, \tilde{B}^{\d}  ww_y \big\ra +  \big\la \tilde{\mathbf{w}}^{\d},\, \tilde{B}^{\d}_x w^2\big\ra + \big\la \tilde{\theta}^{\delta},\, \tilde{B}^{\d} w^2\big\ra  - \big\la \bar{\mathbf{g}}^{\d}_1,\, \tilde{B}^{\d} w^2 \big\ra \nonumber\\
	&\quad  -  \lambda \v\int_{\mathbb{R}}  \Big(\f{b_2}{2p_s}\tilde{u}_s^{\delta} w^2\Big)_y   \tilde{B}^{\d}_{y}  \tilde{B}^{\d}   dy \Big|_{x=L}  + \int_{\mathbb{R}}  \Big\{\lambda\v \f{b_2}{2p_s}\tilde{u}_s \tilde{\Xi}^{\delta} + \tilde{\mathbf{g}}^{\d}_{b2} \Big\} \tilde{B}^{\d} w^2 dy \Big|_{x=L}    \nonumber\\
	&\quad + \int_{\mathbb{R}} \f{2}{\lambda} \f{\tilde{B}^{\d} }{\tilde{u}^{\delta}_s} \Big\{ \f12  \lambda \v \tilde{u}_s \tilde{\theta}^{\delta}_x + \f{\tilde{\mathbf{g}}^{\d}_{b1}}{1+\bm{\chi}}\Big\} w^2 dy \Big|_{x=0}.
\end{align}

For terms on RHS of \eqref{17.52}, it is clear to have that 
\begin{align}\label{17.53}
\begin{split}
\eqref{17.52}_{1,2}
&\lesssim \|\tilde{\mathbf{w}}^{\d}w\|\cdot \|\tilde{B}^{\d}_xw\| + 2 \|\tilde{B}^{\d}_yw\|\cdot \|\tilde{B}^{\d} w_y\|,\\
\eqref{17.52}_{3}
&\lesssim \|y\tilde{\theta}^{\d} w\|\cdot \|\f{1}{y}\tilde{B}^{\d}w\| \lesssim \xi\|\tilde{B}^{\d}_{y} w\|^2 + \xi \|\tilde{B}^{\d} w_y\|^2 + C_{\xi} \|y\theta^{\d} w\|^2,\\
\eqref{17.52}_{5}	&\leq \xi \lambda \v\int_{\mathbb{R}}  \f{b_2}{2p_s}\tilde{u}_s^{\delta} |\tilde{B}^{\d}_{y}|^2 w^2  dy \Big|_{x=L} +  C_{\xi}\v \|\f{\tilde{B}^{\d}w}{\sqrt{\tilde{u}_s^{\delta}}}\|_{x=L}   +  C_{\xi}\v \|\f{\tilde{B}^{\d}w_y}{\sqrt{\tilde{u}_s^{\delta}}} \|^2_{x=L}.
\end{split}
\end{align}

Noting $\gamma=\sqrt{\delta}$, we always have the following useful estimates
\begin{align}\label{17.84-0}
u_s\bar{\chi}(\f{y}{\gamma})\cong u^{\d}_s\,\,\mbox{for}\,\,\, y\geq \gamma, \quad \mbox{and}\quad  \|\sqrt{\tilde{u}_s^{\delta}}f \bar{\chi}(\f{y}{\gamma})\|\lesssim \|\sqrt{\tilde{u}_s}f\|,
\end{align}
which will be used frequently. We integrate by parts to obtain
\begin{align}\label{17.61}
\eqref{17.52}_{4}	&=  \int_0^L\int_0^\infty \pa_{y}\bar{g}^{\d}_1  \Big\{\bar{\chi}(\f{y}{\gamma}) \tilde{B}^{\d}_y w^2 + \bar{\chi}(\f{y}{\gamma}) \tilde{B}^{\d}  w w_y + \f{1}{\gamma} \bar{\chi}'(\f{y}{\gamma})  \tilde{B}^{\d} w^2\Big\}  dydx \nonumber\\
	&\quad + \v \int_0^L\int_0^\infty  \bar{\chi}(\f{y}{\gamma})\pa_{x}\bar{g}^{\d}_1\tilde{B}^{\d}_x w^2  - \v  \int_0^\infty \bar{\chi}(\f{y}{\gamma})\pa_{x}\bar{g}^{\d}_1\tilde{B}^{\d} w^2 dy \Big|^{x=L}_{x=0}  - \int_0^L\int_{\infty}^0 (\cdots) \nonumber\\
	&\lesssim \xi \|\nabla_{\v}\tilde{B}^{\d}w\|^2 + \xi \|\tilde{B}^{\d} w_y\|^2 + \xi  \|\f{\tilde{B}^{\d}}{\sqrt{\tilde{u}_s^{\delta}}}w\|^2_{x=0} + \xi \v \|\f{\tilde{B}^{\d}}{\sqrt{\tilde{u}_s^{\delta}}}w\|_{x=L}  \nonumber\\
	&\quad  + C_{\xi} \|\nabla_{\v}\bar{g}^{\d}_{1} w\|^2 + C_{\xi} \v \|\sqrt{u_s} \bar{g}^{\d}_{1x}\|_{x=L}^2 + C_{\xi} \v^2 \|\sqrt{u_s} \bar{g}^{\d}_{1x}\|_{x=0}^2,
\end{align}
where we have used \eqref{A.5} to control $\|\f{\tilde{B}^{\d}}{y}w\|\lesssim \|\tilde{B}^{\d}_yw\| + \|\tilde{B}^{\d}w_y\|$.

%

 It is clear that 
\begin{align}\label{17.64}
\eqref{17.52}_{7}
&\lesssim \xi  \|\f{\tilde{B}^{\d}w}{\sqrt{\tilde{u}_s^{\d}}}\|^2_{x=0}  + C_{\xi} \v^2 \|\sqrt{u_s}\theta^{\d}_xw\|^2_{x=0}   +  C_{\xi} \|\f{\mathbf{g}^{\delta}_{b1}}{\sqrt{u_s^{\delta}}}w\|^2_{x=0} .
\end{align}

For $\eqref{17.52}_6$, we recall $\tilde{\Xi}^{\d}$ in \eqref{17.57-1}. Integrating by parts in $y$, one obtains
\begin{align}\label{17.63}
	&\lambda\v  \int_{\mathbb{R}} \f{b_2 \tilde{u}_s }{2p_s}\tilde{B}^{\d} \tilde{\Xi}^{\delta} w^2 dy \Big|_{x=L}
	\nonumber\\
	&=- \lambda\v  \int_0^\infty \f{b_2u_s}{2p_s}  \tilde{B}^{\d}_y\,  \theta^{\d}_y\bar{\chi}(\f{y}{\gamma}) w^2 dy \Big|_{x=L} -   \lambda\v  \int_0^\infty \big(\f{b_2u_{s}}{2p_s}w^2\big)_y \tilde{B}^{\d}\,   \theta^{\d}_y\bar{\chi}(\f{y}{\gamma})  dy \Big|_{x=L} \nonumber\\
	&\quad - \lambda\v  \int_0^\infty  \f{b_2u_s}{2\gamma p_s}  \tilde{B}^{\d}\,  \theta^{\d}_{y}\bar{\chi}'(\f{y}{\gamma}) w^2 dy \Big|_{x=L}  + \int_{-\infty}^{0} (\cdots) dy \Big|_{x=L}\nonumber\\
	&\leq \xi \v \int_{\mathbb{R}} \tilde{u}_s |\tilde{B}^{\d}_y|^2 w^2 dy \Big|_{x=L} +  \xi \v  \|\f{\tilde{B}^{\d}w}{\sqrt{\tilde{u}_s}}\|^2_{x=L}  + C_{\xi} \v \|\sqrt{u_s}\theta^{\d}_y\hat{w}\|^2_{x=0}.
\end{align}

\

Recall the definition of  $\tilde{\mathbf{g}}^{\d}_{b2}$ in \eqref{17.57}.
 It is clear that 
\begin{align}\label{17.63-2}
&\int_{0}^\infty \bar{\chi}(\f{y}{\gamma})\big[ \lambda\v (\f{b_2}{p_s} u_s)_y \hat{P}_{xy} + \v \bar{g}^{\delta}_{1x} + \big(b_2 [ g^{\delta}_2+ \v \mathfrak{g}_2 - \v\mathfrak{h}^{\delta}_2 - \v\mathfrak{N}^{\delta}_{21}]\big)_y\big]  \tilde{B}^{\d} w^2 dy \Big|_{x=L}\nonumber\\
&\lesssim  \xi \v  \|\f{\tilde{B}^{\d}w}{\sqrt{\tilde{u}_s}}\|^2_{x=L} + C_{\xi}  a^2 \v \|\sqrt{u_s} \hat{P}_{xy} w\|^2_{x=L} + C_{\xi} \v \|\sqrt{u_s}(b_2\mathfrak{h}^{\delta}_2+b_2\mathfrak{N}^{\delta}_{21})_{y}w\|_{x=L}  \nonumber\\
&\quad + C_{\xi} \v \|\sqrt{u_s}\bar{g}^{\delta}_{1x}w\|^2_{x=L}  +  C_{\xi} \v \|\sqrt{u_s}(b_2\mathfrak{g}_2)_{y}w\|^2_{x=L}   + C_{\xi}\f{1}{\v} \|\sqrt{u_s} (b_2g^{\d}_2)_{y}w\|^2_{x=L} .
\end{align}
Integrating by parts in $y$, one has
\begin{align}
&-  \lambda\v \int_{0}^\infty\f{b_2}{p_s} u_s  \Big(\bar{\chi}(\f{y}{\gamma})\bar{g}^{\d}_{1y}\Big)_y \tilde{B}^{\d} w^2 dy \Big|_{x=L}\nonumber\\
&\lesssim \xi \v \int_{\mathbb{R}} \tilde{u}_s |\tilde{B}^{\d}_y|^2 w^2 dy \Big|_{x=L} +  \xi \v  \|\f{\tilde{B}^{\d}w}{\sqrt{\tilde{u}_s}}\|^2_{x=L} + C_{\xi} \v \|\sqrt{u_s}\bar{g}^{\d}_{1y}\hat{w}\|^2_{x=L},
\end{align}
and
\begin{align}\label{17.63-1}
&-\int_{0}^\infty \Big(\bar{\chi}(\f{y}{\gamma})\hat{G}^{\d}_3\Big)_y   \tilde{B}^{\d} w^2 dy \Big|_{x=L} 
\nonumber\\
&\lesssim \xi \v \int_{\mathbb{R}} \tilde{u}^{\delta}_s |\tilde{B}^{\d}_y|^2w^2 dy \Big|_{x=L}  + \xi \v \|\sqrt{\tilde{u}^{\delta}_s} \tilde{B}^{\d}w_y\|^2_{x=L}
+ C_{\xi}  \f{1}{\v}\|\f{\hat{G}^{\d}_3}{\sqrt{\tilde{u}^{\delta}_s}}w\|^2_{x=L}. 
\end{align}

Noting \eqref{17.57-1}-\eqref{17.57}, then we have from \eqref{17.63}-\eqref{17.63-1} that 
\begin{align}\label{17.65}
\eqref{17.52}_6
&\leq \xi \v \int_{\mathbb{R}} \tilde{u}_s |\tilde{B}^{\d}_y|^2 w^2 dy \Big|_{x=L}  +  \xi \v  \|\f{\tilde{B}^{\d}w}{\sqrt{\tilde{u}_s}}\|^2_{x=L} + \xi \v \|\sqrt{\tilde{u}^{\delta}_s} \tilde{B}^{\d}w_y\|^2_{x=L}    + C_{\xi} \v \|\sqrt{u_s}\theta^{\d}_y\hat{w}\|^2_{x=0}\nonumber\\
&\quad + C_{\xi}  a^2 \v \|\sqrt{u_s} \hat{P}_{xy} w\|^2_{x=L} + C_{\xi} \v \|\sqrt{u_s}\bar{g}^{\d}_{1x}w\|^2_{x=L}  +  C_{\xi} \v \|\sqrt{u_s}(b_2\mathfrak{g}_2)_{y}w\|^2_{x=L} \nonumber\\
&\quad   + C_{\xi}\f{1}{\v} \|\sqrt{u_s} (b_2g^{\d}_2)_{y}w\|^2_{x=L}  + C_{\xi} \v \|\sqrt{u_s}\bar{g}^{\d}_{1y}\hat{w}\|^2_{x=L} 
+ C_{\xi}  \f{1}{\v}\|\f{\hat{G}^{\d}_3}{\sqrt{\tilde{u}^{\delta}_s}}w\|^2_{x=L}\nonumber\\
&\quad + C_{\xi} \v \|\sqrt{u_s}(b_2\mathfrak{h}^{\delta}_2+b_2\mathfrak{N}^{\delta}_{21})_{y}w\|_{x=L}.
\end{align}

Substituting \eqref{17.53}-\eqref{17.64} and \eqref{17.65} into \eqref{17.52}, then using \eqref{17.70-1}, we immediately conclude \eqref{17.107} by taking  $\sqrt{L}\lesssim \xi$. Therefore the proof of Lemma \ref{Blem7.5} is completed. $\hfill\Box$

\begin{lemma}\label{Blem7.6}
It holds that 
\begin{align}\label{17.92}
	&(1-\xi) \bigg\{ \v\|\tilde{u}_s^{\delta}\nabla_{\v}\tilde{B}^{\d}_yw\|^2  +   \int_{\mathbb{R}} \f{2p_s}{\lambda b_2} \tilde{u}_s^{\delta} \big[\v\tilde{B}^{\d}_{x}-\tilde{\mathbf{w}}^{\d}\big]^2w^2 dy \Big|_{x=L} \nonumber\\
	&\qquad\qquad\qquad   + \v \int_{\mathbb{R}} \f{2(\sigma-1)}{\lambda(1+\bm{\chi})} p_s \mathbf{d}_{11} \tilde{u}_s^{\delta}|\tilde{B}^{\d}_{y}|^2w^2 dy \Big|_{x=0}\bigg\} \nonumber\\
	&\leq \|\tilde{u}^{\delta}_sw\tilde{\mathbf{w}}^{\d}_y\|\cdot \sqrt{\v} \|\tilde{u}^{\delta}_s\sqrt{\v}\tilde{B}^{\d}_{xy}w\|   + \xi \|\tilde{u}_s^{\delta}\nabla_{\v} \tilde{\mathbf{w}}^{\d}w\|^2  + C_{\xi}  a^2 \Big\{\|(\sqrt{\v}\nabla_{\v}\tilde{B}^{\d}, \tilde{\mathbf{w}}^{\d})\hat{w}\|^2   \nonumber\\
	&\quad +  \v \|\f{\tilde{B}^{\d}}{\sqrt{\tilde{u}^{\delta}_{s}}} \hat{w}\|^2_{x=0} +  \v^2 \|\sqrt{\tilde{u}_s^{\delta}}\tilde{B}^{\d}_y\hat{w}\|^2_{x=L}\Big\} + C_{\xi} \, (\cdots)_{2,w},
\end{align}
where we have used  the notation
\begin{align}\label{17.99-10}
	(\cdots)_{2,w}:&= \v \|u_{s}\Delta_{\v}\bar{g}^{\d}_1 \hat{w}\|^2  + \v \|(u_s\nabla_{\v}\theta^{\d}_y,a\theta^{\d}_y, a\theta^{\d})\hat{w}\|^2  + a^2 \v^2 \|\sqrt{u_{s}}\theta^{\d}_y\hat{w}\|_{x=L}^2   +  a^2\v^2 \|\sqrt{u_{s}} \hat{P}_{xy}w\|_{x=L}^2\nonumber\\
	&\,\,\, + \|\sqrt{u_s}\hat{G}^{\d}_{3}w_y\|^2_{x=L}+  \v^2 \|\sqrt{u_{s}}\bar{g}^{\d}_{1x}w\|_{x=L}^2  +  a^2 \v^2 \|\sqrt{u_s}\bar{g}^{\d}_{1y}\hat{w}\|^2_{x=L}  + \|\sqrt{u_{s}}(b_2g^{\d}_2)_{y}w\|_{x=L}^2  \nonumber\\
	&\,\,\,  + \v^2 \|\sqrt{u_{s}}(b_2 \mathfrak{g}_2)_{y}w\|_{x=L}^2  +  a^2 \v \|\f{\hat{w}}{\sqrt{\tilde{u}^{\delta}_s}}\big(\v\tilde{u}_s\tilde{\theta}^{\delta}_x, \tilde{\mathbf{g}}^{\d}_{b1}\big)\|^2_{x=0}  + \v^2 \|\sqrt{u_s}(b_2\mathfrak{h}^{\delta}_2, b_2\mathfrak{N}^{\delta}_{21})_{y}w\|_{x=L} \nonumber\\
	&\,\,\,  + \sqrt{\v} \|(\v\tilde{u}_s\tilde{\theta}^{\delta}_x, \f{\tilde{\mathbf{g}}^{\d}_{b1}}{1+\bm{\chi}})w\|^2_{H^{1/2}_{x=0}} + \v^{\f32} \|\tilde{u}_sw\bar{g}^{\d}_{1y}\bar{\chi}(\f{y}{\gamma})\|^2_{H^{1/2}_{x=L}} + \f{1}{\sqrt{\v}} \|\bar{\chi}(\f{y}{\gamma})\hat{G}^{\d}_{3} w\|^2_{H^{1/2}_{x=L}}  .
\end{align}
\end{lemma}

\noindent{\bf Proof.} In the following, we divide the proof into four steps. 

{\it Step 1.} Applying $\pa_y$ to $\eqref{17.50}_1$, one has
\begin{align}\label{17.74}
	-\Delta_{\v} \tilde{B}^{\d}_y  = - \tilde{\mathbf{w}}^{\d}_{xy}  - \bar{\mathbf{g}}^{\d}_{1y} +\tilde{\theta}^{\delta}_y.
\end{align}
Multiplying \eqref{17.74} by  $\v \,(\tilde{u}_s^{\delta} w)^2\tilde{B}^{\d}_y$, one obtains that 
\begin{align}\label{17.75}
	&\v\|\tilde{u}_s^{\delta}\nabla_{\v}\tilde{B}^{\d}_yw\|^2 + \v \big\la \nabla_{\v}\tilde{B}^{\d}_y\cdot \nabla_{\v}(\tilde{u}_s^{\delta}w) ,\, 2\tilde{u}_s^{\delta}w \tilde{B}^{\d}_y \big\ra 
	- \v^2 \int_{\mathbb{R}}  (\tilde{u}_s^{\delta}w)^2\tilde{B}^{\d}_{xy} \tilde{B}^{\d}_y dy \Big|_{x=L} \nonumber\\
	&\quad + \v^2 \int_{\mathbb{R}}  (\tilde{u}_s^{\delta}w)^2\tilde{B}^{\d}_{xy} \tilde{B}^{\d}_y dy \Big|_{x=0}  \nonumber\\
	&= - \v \big\la  \tilde{\mathbf{w}}^{\d}_{xy},\, (\tilde{u}_s^{\delta}w)^2\tilde{B}^{\d}_y \big\ra - \v \big\la \bar{\mathbf{g}}^{\d}_{1y},\, (\tilde{u}_s^{\delta}w)^2\tilde{B}^{\d}_y \big\ra + \v \big\la \tilde{\theta}^{\delta}_{y},\, (\tilde{u}_s^{\delta}w)^2\tilde{B}^{\d}_y \big\ra.
\end{align}
A direct calculation shows that 
\begin{align}\label{17.76}
\begin{split}
\v \big|\big\la \nabla_{\v}\tilde{B}^{\d}_y\cdot \nabla_{\v}(\tilde{u}_s^{\delta}w) ,\, 2\tilde{u}_s^{\delta}w \tilde{B}^{\d}_y \big\ra \big|
&\leq \xi \v\|\tilde{u}_s^{\delta}\nabla_{\v}\tilde{B}^{\d}_y w\|^2 + C_{\xi} a^2 \v \|\tilde{B}^{\d}_y\hat{w}\|^2, 
\end{split}
\end{align}
\begin{align}
	\v \big\la \tilde{\theta}^{\delta}_{y},\, (\tilde{u}_s^{\delta}w)^2\tilde{B}^{\d}_y \big\ra
	& \leq a \v \|\tilde{B}^{\d}_yw\|\big\{ \|\tilde{u}_s^{\delta}\theta^{\d}_y \bar{\chi}(\f{y}{\gamma}) w\| + \|\tilde{u}_s^{\delta}\theta^{\d} \f{1}{\gamma} \bar{\chi}'(\f{y}{\gamma})w\| \big\} \nonumber\\
	&\lesssim a^2\v \|\tilde{B}^{\d}_yw\|^2 + a^2\v \|(\theta^{\d}_y, \theta^{\d} )w\|^2,
\end{align}
and
\begin{align}\label{17.78}
- \v \big\la \bar{\mathbf{g}}^{\d}_{1y},\, (\tilde{u}_s^{\delta}w)^2\tilde{B}^{\d}_y \big\ra
& = \v \big\la \bar{\mathbf{g}}^{\d}_{1},\, (\tilde{u}_s^{\delta}w)^2\tilde{B}^{\d}_{yy} \big\ra + \v \big\la \bar{\mathbf{g}}^{\d}_{1},\, 2\tilde{u}^{\delta}_sw (\tilde{u}^{\delta}_{s}w)_y \tilde{B}^{\d}_y \big\ra \nonumber\\
&\leq \xi \v\|\tilde{u}_s^{\delta} \tilde{B}^{\d}_{yy}w\|^2  + C_{\xi} \v \|u_{s}\Delta_{\v}\bar{g}^{\d}_1 \hat{w}\|^2 +  C a^2\v\|\tilde{B}^{\d}_yw\|^2 .
\end{align}

We integrate by parts in $x$ and then in $y$ to obtain that 
\begin{align}\label{17.79}
 - \v\big\la  \tilde{\mathbf{w}}^{\d}_{xy},\, (\tilde{u}_s^{\delta}w)^2\tilde{B}^{\d}_y \big\ra 
	&=\v\big\la  \tilde{\mathbf{w}}^{\d}_{y},\, (\tilde{u}_s^{\delta}w)^2\tilde{B}^{\d}_{xy} \big\ra  +  \v\big\la  \tilde{\mathbf{w}}^{\d}_{y},\, 2\tilde{u}_s^{\delta}\tilde{u}_{sx}^{\delta}w^2 \tilde{B}^{\d}_y \big\ra 
	+ \v \int_{\mathbb{R}} (\tilde{u}_s^{\delta}w)^2\tilde{\mathbf{w}}^{\d} \tilde{B}^{\d}_{yy} dy \Big|_{x=L}\nonumber\\
	&\quad  + \v \int_{\mathbb{R}} 2\tilde{u}_s^{\delta} w (\tilde{u}^{\delta}_{s}w)_y \tilde{\mathbf{w}}^{\d} \tilde{B}^{\d}_y dy \Big|_{x=L}  -  \v\int_{\mathbb{R}} (\tilde{u}_s^{\delta}w)^2\tilde{\mathbf{w}}^{\d}  \tilde{B}^{\d}_{yy} dy \Big|_{x=0} \nonumber\\
	&\quad  - \v \int_{\mathbb{R}} 2\tilde{u}_s^{\delta} w ( \tilde{u}^{\delta}_{s}w)_y \tilde{\mathbf{w}}^{\d} \tilde{B}^{\d}_y dy \Big|_{x=0}.
\end{align}

For the boundary terms on LHS of \eqref{17.75}, we integrate by parts in $y$ to have 
\begin{align}\label{17.81}
\begin{split}
- \v^2 \int_{\mathbb{R}}  (\tilde{u}_s^{\delta}w)^2\tilde{B}^{\d}_{xy} \tilde{B}^{\d}_y dy \Big|_{x=L} 
& =  \v^2 \int_{\mathbb{R}} \tilde{B}^{\d}_{x}  \big\{ (\tilde{u}_s^{\delta}w)^2\tilde{B}^{\d}_{y} \big\}_y dy \Big|_{x=L} ,\\
\v^2 \int_{\mathbb{R}}  (\tilde{u}_s^{\delta}w)^2\tilde{B}^{\d}_{xy} \tilde{B}^{\d}_y dy \Big|_{x=0}  
& = - \v^2 \int_{\mathbb{R}}  \tilde{B}^{\d}_{x}  \big\{ (\tilde{u}_s^{\delta}w)^2\tilde{B}^{\d}_{y} \big\}_y dy \Big|_{x=0} . 
\end{split}
\end{align}
 
Substituting \eqref{17.76}-\eqref{17.81} into \eqref{17.75} to obtain
\begin{align}\label{17.82}
	&(1-\xi)\v\|\tilde{u}_s^{\delta}\nabla_{\v}\tilde{B}^{\d}_yw\|^2  - \v \int_{\mathbb{R}}  \big\{(\tilde{u}_s^{\delta} w)^2\tilde{B}^{\d}_{y}\big\}_y \big[ \v \tilde{B}^{\d}_{x} - \tilde{\mathbf{w}}^{\d} \big] dy \Big|_{x=0} \nonumber\\
	&\,\,\, + \v \int_{\mathbb{R}} \big\{(\tilde{u}_s^{\delta} w)^2\tilde{B}^{\d}_{y}\big\}_y \big[\v\tilde{B}^{\d}_{x}-\tilde{\mathbf{w}}^{\d}\big] dy \Big|_{x=L}  \nonumber\\
	&\leq \v\big\la  \tilde{\mathbf{w}}^{\d}_{y},\, (\tilde{u}_s^{\delta}w)^2\tilde{B}^{\d}_{xy} \big\ra   +  \v\big\la  \tilde{\mathbf{w}}^{\d}_{y},\, 2\tilde{u}_s^{\delta}\tilde{u}_{sx}^{\delta}w^2 \tilde{B}^{\d}_y \big\ra  +   C_{\xi} a^2 \v \|\tilde{B}^{\d}_y\hat{w}\|^2 \nonumber\\
	&\quad  + C_{\xi} \v \|u_{s}\Delta_{\v}\bar{g}^{\d}_1 \hat{w}\|^2  + C a^2\v \|(\theta^{\d}_y, \theta^{\d} )w\|^2 .
\end{align}

{\it Step 2.}
Noting $\eqref{17.50}_2$, for the boundary terms at $x=0$ on LHS of \eqref{17.82}, we have 
\begin{align}\label{17.83}
	& - \v \int_{\mathbb{R}}  \big\{(\tilde{u}_s^{\delta} w)^2\tilde{B}^{\d}_{y}\big\}_y \big[ \v \tilde{B}^{\d}_{x} - \tilde{\mathbf{w}}^{\d} \big] dy \Big|_{x=0} \nonumber\\
	&=- \v \int_{\mathbb{R}} \f{2(\sigma-1)}{\lambda(1+\bm{\chi})} p_s \mathbf{d}_{11} \tilde{u}_s^{\delta}\tilde{B}^{\d}_{yy} \tilde{B}^{\d} w^2 dy \Big|_{x=0} +  \v \int_{\mathbb{R}} \tilde{u}_s^{\delta}\tilde{B}^{\d}_{yy} \Big[ \v \tilde{u}_s \tilde{\theta}^{\delta}_x + \f{2\tilde{\mathbf{g}}^{\d}_{b1}}{\lambda(1+\bm{\chi})}\Big] w^2 dy \Big|_{x=0} \nonumber\\
	&\quad -   \v\int_{\mathbb{R}} \tilde{B}^{\d}_y   w (\tilde{u}^{\delta}_{s}w)_y \Big\{\f{4(\sigma-1)}{\lambda(1+\bm{\chi})}  p_s \mathbf{d}_{11}  \tilde{B}^{\d}- \f{4\tilde{\mathbf{g}}^{\d}_{b1}}{\lambda(1+\bm{\chi})}  - 2\v \tilde{u}_s \tilde{\theta}^{\delta}_x \Big\} dy \Big|_{x=0} \nonumber\\
	&\geq (1-\xi)\v \int_{\mathbb{R}} \f{2(\sigma-1)}{\lambda(1+\bm{\chi})} p_s \mathbf{d}_{11} \tilde{u}_s^{\delta}|\tilde{B}^{\d}_{y}|^2w^2 dy \Big|_{x=0} - C_{\xi} a^2 \v \|\f{\tilde{B}^{\d}}{\sqrt{\tilde{u}^{\delta}_{s}}} \hat{w}\|^2_{x=0} \nonumber\\
	&\quad  - C_{\xi}a^2 \v \|\f{\hat{w}}{\sqrt{\tilde{u}^{\delta}_s}}\big(\v\tilde{u}_s\tilde{\theta}^{\delta}_x, \tilde{\mathbf{g}}^{\d}_{b1}\big)\|^2_{x=0}  - C \v \|\tilde{u}_s^{\delta}\tilde{B}^{\d}_{y}w\|_{H^{1/2}_{x=0}} \cdot \|(\v\tilde{u}_s\tilde{\theta}^{\delta}_x, \f{\tilde{\mathbf{g}}^{\d}_{b1}}{1+\bm{\chi}})w\|_{H^{1/2}_{x=0}}.
\end{align}
We remark that we have to apply $H^{1/2}_{x=0}(\R)$ argument for $\v\tilde{u}_s\tilde{\theta}^{\delta}_x$, otherwise we may meet $\|\nabla^3\theta^{\d}\|$ which is  hard to control.
 
\medskip

{\it Step 3.} We consider the boundary term at $x=L$ on LHS of \eqref{17.82}. It is direct to see
\begin{align}\label{17.91}
	&\v \int_{\mathbb{R}} 2\tilde{u}_s^{\delta} w ( \tilde{u}^{\delta}_{s}w)_y \tilde{B}^{\d}_y \big[\v \tilde{B}^{\d}_{x}-\tilde{\mathbf{w}}^{\d}\big] dy \Big|_{x=L} \nonumber\\
	&\lesssim \xi \|\sqrt{\tilde{u}^{\delta}_{s}} [\v\tilde{B}^{\d}_x-\tilde{\mathbf{w}}^{\d}]w\|_{x=L}^2 + C_{\xi} a^2 \v^2 \|\sqrt{\tilde{u}_s^{\delta}}\tilde{B}^{\d}_y\hat{w}\|^2_{x=L}.
\end{align}

It follows from $\eqref{17.50}_3$ that
\begin{align}\label{17.84}
  \v \int_{\mathbb{R}}  (\tilde{u}_s^{\delta}w)^2 \tilde{B}^{\d}_{yy} \big[\v\tilde{B}^{\d}_{x}-\tilde{\mathbf{w}}^{\d}\big] dy \Big|_{x=L} 
 &= - \int_{\mathbb{R}} \f{2p_s}{\lambda b_2}\tilde{u}_s^{\delta} w^2 \big[\v\tilde{B}^{\d}_{x}-\tilde{\mathbf{w}}\big]\big[\lambda\v\f{b_2}{2p_s} \tilde{u}_s \tilde{\Xi}^{\delta} + \tilde{\mathbf{g}}^{\d}_{b2}\big] dy \Big|_{x=L} \nonumber\\
 &\quad + \int_{\mathbb{R}} \f{2p_s}{\lambda b_2} \tilde{u}_s^{\delta} \big|\v\tilde{B}^{\d}_{x}-\tilde{\mathbf{w}}^{\d}\big|^2 w^2 dy \Big|_{x=L}.
\end{align}
To deal with the terms on RHS of \eqref{17.84}, we apply $H^{1/2}_{x=L}$ estimate to obtain
\begin{align}\label{17.85-0}
\begin{split}
\sqrt{\v}\|\tilde{u}_s^{\delta}[\v\tilde{B}^{\d}_x-\tilde{\mathbf{w}}^{\d}] w\|^2_{H^{1/2}_{x=L}} 
&\lesssim \v \|\tilde{u}_s^{\delta}\nabla_{\v}^2\tilde{B}^{\d}w\|^2 + \|\tilde{u}_s^{\delta}\nabla_{\v} \tilde{\mathbf{w}}^{\d}w\|^2  + a^2 \|(\sqrt{\v}\nabla_{\v}\tilde{B}^{\d}, \tilde{\mathbf{w}}^{\d})\hat{w}\|^2 ,\\
\sqrt{\v}\|\tilde{u}_s \theta^{\d}_{y} \bar{\chi}(\f{y}{\gamma})w\|^2_{H^{1/2}_{x=L}}
&\lesssim \|u_s \nabla_{\v}\theta^{\d}_yw\|^2 + a^2 \|\theta^{\d}_y\hat{w}\|^2,
\end{split}
\end{align}
which yields that 
\begin{align}\label{17.85}
	&- \v \int_{\mathbb{R}} \tilde{u}_s^{\delta}  \tilde{u}_s \big[\v\tilde{B}^{\d}_{x}-\tilde{\mathbf{w}}^{\d}\big] \tilde{\Xi}^{\delta} w^2 dy \Big|_{x=L} 
	\nonumber\\
	&=- \v\int_{\mathbb{R}} \tilde{u}^{\delta}_s w \big[\v\tilde{B}^{\d}_{x}-\tilde{\mathbf{w}}^{\d}\big]  \big(\tilde{u}_s w \theta^{\d}_{y} \bar{\chi}(\f{y}{\gamma})\big)_y dy \Big|_{x=L} + \v\int_{\mathbb{R}} \big((\tilde{u}^{\delta}_s w)^2\big)_y \big[\v\tilde{B}^{\d}_{x}-\tilde{\mathbf{w}}^{\d}\big]   \theta^{\d}_{y} \bar{\chi}(\f{y}{\gamma}) dy \Big|_{x=L}\nonumber\\
	&\quad + \v\f{1}{\gamma}\int_{\mathbb{R}} \tilde{u}_s^{\delta} \tilde{u}_s w^2 \big[\v\tilde{B}^{\d}_{x}-\tilde{\mathbf{w}}^{\d}\big]     \theta^{\d}_{y}\bar{\chi}'(\f{y}{\gamma}) dy \Big|_{x=L}\nonumber\\
	&\lesssim \v \|\tilde{u}_s^{\delta}[\v\tilde{B}^{\d}_x-\tilde{\mathbf{w}}^{\d}]w\|_{H^{1/2}_{x=L}}\|u_s \theta^{\d}_{y} \bar{\chi}(\f{y}{\gamma})w\|_{H^{1/2}_{x=L}} + a\v \|\sqrt{\tilde{u}^{\delta}_{s}} [\v\tilde{B}^{\d}_x-\tilde{\mathbf{w}}^{\d}]w\|_{x=L} \|\sqrt{\tilde{u}_{s}}\theta^{\d}_y\hat{w}\|_{x=L}\nonumber\\
	&\lesssim \xi\v \|\tilde{u}_s^{\delta}\nabla_{\v}^2\tilde{B}^{\d}w\|^2 + \xi \|\sqrt{\tilde{u}^{\delta}_{s}} [\v\tilde{B}^{\d}_x-\tilde{\mathbf{w}}^{\d}]w\|_{x=L}^2 + \xi \|u_s^{\delta}\nabla_{\v}\tilde{\mathbf{w}}^{\d}w\|^2 + \xi  a^2 \v \|\nabla_{\v}\tilde{B}^{\d}\hat{w}\|^2 \nonumber\\
	&\quad + \xi a^2 \|\tilde{\mathbf{w}}^{\d}\hat{w}\|^2 + C_{\xi} \v \|u_s\nabla_{\v}\theta^{\d}_yw\|^2  + C_{\xi} a^2 \v \|\theta^{\d}_y\hat{w}\|^2  + C_{\x} a^2 \v^2 \|\sqrt{u_{s}}\theta^{\d}_y\hat{w}\|_{x=L}^2.
\end{align}

We also note that 
\begin{align}\label{17.91-1}
&\int_{\mathbb{R}} \f{p_sw^2}{\lambda b_2}\tilde{u}_s^{\delta} \big[\v\tilde{B}^{\d}_{x}-\tilde{\mathbf{w}}^{\d}\big] \Big\{\bar{\chi}(\f{y}{\gamma})\big[ \lambda\v (\f{2b_2u_s}{p_s} )_y \hat{P}_{xy}  + \v \bar{g}^{\delta}_{1x} + \big(b_2 [ g^{\delta}_2+ \v \mathfrak{g}_2 - \v\mathfrak{h}^{\delta}_2 - \v\mathfrak{N}^{\delta}_{21}]\big)_y\big]\Big\}  dy \Big|_{x=L}\nonumber\\
&\lesssim \xi \|\sqrt{\tilde{u}^{\delta}_{s}} [\v\tilde{B}^{\d}_x-\tilde{\mathbf{w}}^{\d}]w\|_{x=L}^2 + C_{\xi}  a^2\v^2 \|\sqrt{u_{s}} \hat{P}_{xy}w\|_{x=L}^2  + C_{\xi}\|\sqrt{u_{s}}(b_2g^{\d}_2)_{y}w\|_{x=L}^2  \nonumber\\
&\quad+ C_{\xi}\v^2 \|\sqrt{u_{s}}\bar{g}^{\d}_{1x}w\|_{x=L}^2 + C_{\xi}\v^2 \|\sqrt{u_{s}}(b_2 \mathfrak{g}_2)_{y}w\|_{x=L}^2 + C_{\xi} \v^2 \|\sqrt{u_s}(b_2\mathfrak{h}^{\delta}_2+b_2\mathfrak{N}^{\delta}_{21})_{y}w\|_{x=L},
\end{align}
and
\begin{align}\label{17.92-1}
&-  \int_{\mathbb{R}}  \tilde{u}_s^{\delta}  w^2  \big[\v\tilde{B}^{\d}_{x}-\tilde{\mathbf{w}}^{\d}\big] \Big\{ \v \tilde{u}_s  \Big(\bar{\chi}(\f{y}{\gamma})\bar{g}_{1y}\Big)_y + \f{2p_s}{\lambda b_2}\Big(\hat{G}^{\d}_3\bar{\chi}(\f{y}{\gamma})\Big)_y\Big\} dy \Big|_{x=L} \nonumber\\
&\lesssim \|\tilde{u}_s^{\delta} w\big[\v\tilde{B}^{\d}_{x}-\tilde{\mathbf{w}}^{\d}\big]\|_{H^{1/2}_{x=L}} \Big\{\|\bar{\chi}(\f{y}{\gamma}) \hat{G}^{\d}_{3}w\|_{H^{1/2}_{x=L}} + \v \|\tilde{u}_sw\bar{g}^{\delta}_{1y}\bar{\chi}(\f{y}{\gamma})\|_{H^{1/2}_{x=L}}\Big\}  + C_{\xi} \|\sqrt{u_s}\hat{G}^{\d}_3 w_y\|^2_{x=L}\nonumber\\
&\quad + \xi \|\sqrt{\tilde{u}^{\delta}_{s}} w [\v\tilde{B}^{\d}_x-\tilde{\mathbf{w}}^{\d}]\|_{x=L}^2 + C_{\x}a^2 \v^2 \|\sqrt{u_s}\bar{g}^{\delta}_{1y}\hat{w}\|^2_{x=L}.
\end{align}

Recall $\tilde{\mathbf{g}}^{\d}_{b2}$ in \eqref{17.57}, then we have  from \eqref{17.84} and \eqref{17.85}-\eqref{17.92-1}  that 
\begin{align}\label{17.90}
& \v \int_{\mathbb{R}}  (\tilde{u}_s^{\delta}w)^2 \tilde{B}^{\d}_{yy} \big[\v\tilde{B}^{\d}_{x}-\tilde{\mathbf{w}}^{\d}\big] dy \Big|_{x=L}  \nonumber\\
&\geq (1-\xi) \int_{\mathbb{R}} \f{2p_s}{\lambda b_2} \tilde{u}_s^{\delta} \big|\v\tilde{B}^{\d}_{x}-\tilde{\mathbf{w}}^{\d}\big|^2w^2 dy \Big|_{x=L} - \xi\v \|\tilde{u}_s^{\delta}\nabla_{\v}^2\tilde{B}^{\d}w\|^2 - \xi \|\tilde{u}_s^{\delta}\nabla_{\v} \tilde{\mathbf{w}}^{\d}w\|^2\nonumber\\
&\quad  - C_{\xi} \Big\{ a^2 \|(\sqrt{\v}\nabla_{\v}\tilde{B}^{\d}, \tilde{\mathbf{w}}^{\d})\hat{w}\|^2 + \f{1}{\sqrt{\v}} \|\bar{\chi}(\f{y}{\gamma}) \hat{G}^{\d}_3 w\|^2_{H^{1/2}_{x=L}}   + \v^{\f32} \|\tilde{u}_sw\bar{g}^{\d}_{1y}\bar{\chi}(\f{y}{\gamma})\|^2_{H^{1/2}_{x=L}} \nonumber\\
&\quad  +  a^2\v^2 \|\sqrt{u_{s}} \hat{P}_{xy}w\|_{x=L}^2  + \|\sqrt{u_s} \hat{G}^{\d}_3 w_y\|^2_{x=L} +  a^2 \v^2 \|\sqrt{u_s}\bar{g}^{\d}_{1y}\hat{w}\|^2_{x=L}  +  \v^2 \|\sqrt{u_{s}}\bar{g}^{\d}_{1x}w\|_{x=L}^2\nonumber\\
&\quad  + \|\sqrt{u_{s}}(b_2g^{\d}_2)_{y}w\|_{x=L}^2 + \v^2 \|\sqrt{u_{s}}(b_2 \mathfrak{g}_2)_{y}w\|_{x=L}^2 + \v \|(u_s\nabla_{\v}\theta^{\d}_y,a\theta^{\d}_y)\hat{w}\|^2 \nonumber\\
&\quad + a^2 \v^2 \|\sqrt{u_{s}}\theta^{\d}_y\hat{w}\|_{x=L}^2 + \v^2 \|\sqrt{u_s}(b_2\mathfrak{h}^{\delta}_2, b_2\mathfrak{N}^{\delta}_{21})_{y}w\|_{x=L}\Big\}.
\end{align}

{\it Step 4.} We  have from \eqref{17.50} that 
\begin{align}\label{17.111-1}
\v \|\tilde{u}_s^{\delta}\v \tilde{B}^{\d}_{xx}w\|^2 \lesssim \v \|\tilde{u}_s^{\delta}\tilde{B}^{\d}_{yy}w\|^2 + \|\tilde{u}_s\sqrt{\v}\tilde{\mathbf{w}}^{\d}_xw\|^2 + \v  \|u_s \Delta_{\v}\bar{g}^{\d}_1w\|^2 + a^2 \v \|\theta^{\d} w\|^2,
\end{align}
Using  \eqref{Z.11-2}, one gets
\begin{align}\label{17.116-1}
&\v \|\tilde{u}_s^{\delta}\tilde{B}^{\d}_{y}w\|_{H^{1/2}_{x=0}} \cdot \|(\v\tilde{u}_s\tilde{\theta}^{\delta}_x, \f{\tilde{\mathbf{g}}^{\d}_{b1}}{1+\bm{\chi}})w\|_{H^{1/2}_{x=0}} \nonumber\\
&\lesssim \xi \v\|\tilde{u}_s^{\delta}\nabla_{\v}\tilde{B}^{\d}_y w\|^2 + \xi a^2\v \|\nabla_{\v}\tilde{B}^{\d}\hat{w}\|^2   + C_{\x} \sqrt{\v} \|(\v\tilde{u}_s\tilde{\theta}^{\delta}_x, \f{\tilde{\mathbf{g}}^{\d}_{b1}}{1+\bm{\chi}})w\|^2_{H^{1/2}_{x=0}}.
\end{align} 
Substituting \eqref{17.90}-\eqref{17.91} and \eqref{17.83} into \eqref{17.82} and using \eqref{17.111-1}-\eqref{17.116-1}, one concludes \eqref{17.92}. Therefore the proof of Lemma \ref{Blem7.6} is completed. $\hfill\Box$

\medskip

\subsubsection{ Iteration scheme and uniform estimates.} 
Recall the oddness of  $\tilde{B}^{\d}$  in $y\in \mathbb{R}$. Taking $w=1$,  we have from \eqref{17.25} and \eqref{17.107} that 
\begin{align} \label{17.110}
	&(1 -\xi_0 ) \Big\{(2-\xi_1) \|\nabla_{\v}\tilde{B}^{\d}\|^2 + \f{\delta}{\v} \|\nabla_{\v}\mathbf{w}^{\delta}\|^2  +  \f{1}{\v}\|\mathbf{w}^{\delta}\|^2\Big\}  +  (\mu+\lambda) \int_0^\infty \f{u^{\delta}_s}{2p_s} |\mathbf{w}^{\delta}|^2 dy \Big|_{x=L}\nonumber\\
	&\quad  + (2-\xi_1) (1-\xi_0) \Big\{ \v\int_0^\infty \f{\lambda b_2}{2p_s}u_s^{\delta}  |\tilde{B}^{\d}_{y}|^2  dy \Big|_{x=L}  +   \int_{0}^{\infty}  p_s \mathbf{d}_{11} \f{2(\sigma-1)}{\lambda[1+\bm{\chi}]}\f{|\tilde{B}^{\d}|^2}{\tilde{u}^{\delta}_s} dy \Big|_{x=0}\Big\}  \nonumber\\
	&\leq  \f12 \f{\sigma^2}{\mu+\lambda} \int_0^\infty p_s\mathbf{d}_{11}^2\f{|\dot{B}_1|^2}{u_s^{\delta}}  dy \Big|_{x=0} + C_{\xi_0} \f{1}{\v}\|\hat{\mathcal{G}}^{\d}\|^2  +  (2-\xi_1) \|\mathbf{w}^{\delta}\|\cdot \|\tilde{B}^{\d}_x\|   + 2C_{\xi_0} (\cdots)_{1,1},
\end{align}
where all the integrations and norms are on $(0,L)\times\mathbb{R}_+$ or $\mathbb{R}_+$, and $\xi_0>0$ is a small constant determined later.

A direct calculation shows that 
\begin{align}\label{17.111}
	&(1 -\xi_0 ) \Big\{(2-\xi_1) \|\nabla_{\v}\tilde{B}^{\d}\|^2 +  \f{1}{\v}\|\mathbf{w}^{\delta}\|^2\Big\} - (2-\xi_1) \|\mathbf{w}^{\delta}\|\cdot \|\tilde{B}^{\d}_x\|\nonumber\\
	&=-\xi_0 \Big\{(2-\xi_1) \|\nabla_{\v}\tilde{B}^{\d}\|^2 +  \f{1}{\v}\|\mathbf{w}^{\delta}\|^2\Big\} + (2-\xi_1) \|\nabla_{\v}\tilde{B}^{\d}\|^2 +  \f{1}{\v}\|\mathbf{w}^{\delta}\|^2 - (2-\xi_1) \|\mathbf{w}\|\cdot \|\tilde{B}^{\d}_x\|\nonumber\\
	&\geq  \f14 \Big\{ \|\nabla_{\v}\tilde{B}^{\d}\|^2 +  \f{1}{\v}\|\mathbf{w}^{\delta}\|^2\Big\},
\end{align}
where we have taken $\xi_0, \xi_1>0$ suitably small.

Noting $\chi\in[0,1]$, by taking $\dis\sigma=2+2\f{\mu}{\lambda}$, we have that 
\begin{align*} 
\f12 \f{\sigma^2}{\mu+\lambda}   -   \f{4(\sigma-1)}{\lambda(1+\bm{\chi})} &\leq \f12 \f{\sigma^2}{\mu+\lambda}   -   \f{2(\sigma-1)}{\lambda}  = \f1{2(\mu+\lambda)} \Big\{\sigma^2 - 4(\f{\mu}{\lambda}+1)(\sigma-1)\Big\} \nonumber\\
&= \f1{2(\mu+\lambda)} \Big\{ \cancel{|\sigma - 2(\f{\mu}{\lambda}+1)|^2} - 4\f{\mu}{\lambda}(\f{\mu}{\lambda}+1)\Big\} =-2\f{\mu}{\lambda^2} <0.
\end{align*}
Then there exist $\xi_0, \xi_1>0$ such that 
\begin{align*}
 \f12 \f{\sigma^2}{\mu+\lambda} -	(2-\xi_1) \f{1}{\lambda}(\sigma-1)(1-\xi_0) < -\f32 \f{\mu}{\lambda^2} <0,
\end{align*}
which is equivalent to 
\begin{align}\label{17.113}
	\Lambda := \f12 \f{\sigma^2}{\mu+\lambda}\Big[(2-\xi_1) \f{1}{\lambda}(\sigma-1)(1-\xi_0)\Big]^{-1} <1.
\end{align}
Therefore we have from \eqref{17.110}, \eqref{17.111} and \eqref{17.113} that
\begin{align}\label{17.114}
	&c_3\Big[\f{\delta}{\v} \|\nabla_{\v}\mathbf{w}^{\delta}\|^2  +  \f{1}{\v}\|\mathbf{w}^{\delta}\|^2  +  \|\nabla_{\v}\tilde{B}^{\d}\|^2+ \|\sqrt{u_s^{\delta}}\mathbf{w}^{\delta}\|^2_{x=L} + \v \|\sqrt{u_s^{\delta}} \tilde{B}^{\d}_y\|^2_{x=L} \Big] \nonumber\\
	&\,\,\, + \int_0^\infty  \mathbf{d}_{11} \f{p_s}{u^{\delta}_s} |\tilde{B}^{\d}|^2 dy \Big|_{x=0}  \leq \Lambda \int_0^\infty \mathbf{d}_{11} \f{p_s}{u_s^{\delta}} |\dot{B}_1|^2 dy \Big|_{x=0}  + C \f{1}{\v}\|\hat{\mathcal{G}}^{\d}\|^2  +  (\cdots)_{1,1},
\end{align}
where we have used the fact that $\mathbf{d}_{11}=1 - \f{u_s^2}{2T_s} + O(1)\v u_s < 1 + O(1)\v u_s$.

\smallskip

For the higher order derivative estimates, by noting $\gamma=\sqrt{\delta}$, we know that 
\begin{align}\label{17.126-0}
\|u_s^{\delta}\tilde{\mathbf{w}}^{\d}_yw\| \leq \|u_s^{\delta} \mathbf{w}^{\delta}_yw\| + \f{1}{\gamma}\|u_s^{\delta}\mathbf{w}^{\delta} \chi'(\f{y}{\gamma})\|\leq \|u_s^{\delta} \mathbf{w}^{\delta}_yw\| + Ca \|\mathbf{w}^{\delta}\|.
\end{align}

Taking $w=1$, using \eqref{17.126-0}, \eqref{17.99-1} and \eqref{17.99}, and by similar argument as \eqref{17.114},  we have from \eqref{17.126} and \eqref{17.92}  that
\begin{align}\label{17.107-1}
	&c_4\Big\{\v\|u_s^{\delta}\nabla_{\v}\tilde{B}^{\d}_y\|^2  + \delta \|u_s^{\delta} \nabla_{\v}\mathbf{w}^{\delta}_y\|^2 + \|u_s^{\delta}\mathbf{w}^{\delta}_y\|^2 + \|\sqrt{\tilde{u}_s^{\delta}} \big[\v\tilde{B}^{\d}_{x}-\tilde{\mathbf{w}}^{\d}\big]\|^2_{x=L} \nonumber\\
	&\qquad + \v \|(u_s^{\delta})^{\f32} \mathbf{w}^{\delta}_y\|^2_{x=L} \Big\}  +  \v \int_{0}^\infty p_s \mathbf{d}_{11} u_s^{\delta} |\tilde{B}^{\d}_{y}|^2 dy \Big|_{x=0} \nonumber\\
	&\leq \Lambda \v \int_0^\infty p_s\mathbf{d}_{11} u_s^{\delta} |\dot{B}_{y}|^2 dy \Big|_{x=0}    +  \xi\Big\{\f{\delta}{\v}\|\nabla_{\v}\mathbf{w}^{\delta}\|^2  +   \f{1}{\v}\|\mathbf{w}^{\delta}\|^2 +  \f{1}{\v}\|\hat{\mathcal{G}}^{\d}\|^2\Big\} +  C_{\xi} \|u_s^{\delta}\hat{\mathcal{G}}^{\d}_y\|^2\nonumber\\
	&\quad  + C_{\xi} a^2\v\Big\{\|\nabla_{\v}\tilde{B}^{\d}\|^2 +  \|\f{\tilde{B}^{\d}}{\sqrt{\tilde{u}^{\delta}_{s}}}\|^2_{x=0}   + \|\f{\dot{B}}{\sqrt{u^{\delta}_s}}\|^2_{x=0} +  \v  \|\sqrt{\tilde{u}_s^{\delta}}\tilde{B}^{\d}_y\|^2_{x=L} \Big\}   +  C_{\xi}(\cdots)_{2,1},
\end{align}
where all the integrations and norms are also on $(0,L)\times\mathbb{R}_+$ or $\mathbb{R}_+$.

\medskip

  Now we consider the iteration scheme for $n\geq 0,1,2,\cdots$. 
\begin{align}\label{17.115}
	\begin{cases}
		\dis -\delta \Delta_{\v}\mathbf{w}^{\delta, n+1} + \mathbf{w}^{\delta,n+1}  + (\mu+\lambda) \v \f{1}{p^{\v,\d}}(U^{\v,\d}\cdot\nabla)\mathbf{w}^{\delta,n+1} = \hat{\mathcal{G}}^{\d},\,\, (x,y)\in (0,L)\times \mathbb{R}_+,\\
		\dis -\delta \mathbf{w}_x^{\delta,n+1} + (\mu+\lambda) \f{1}{p_s}u_s^{\delta}\mathbf{w}^{\delta,n+1} \big|_{x=0}=-\sigma \mathbf{d}_{11}\tilde{B}^{n},\\[2.5mm]
		\mathbf{w}^{\delta,n+1}_x\big|_{x=L}=\mathbf{w}^{\delta,n+1}_y\big|_{y=0}=0,
	\end{cases}
\end{align}
and 
\begin{align}\label{17.116}
	\begin{cases}
		\dis -\Delta_{\v} \tilde{B}^{\d, n+1} = - \tilde{\mathbf{w}}^{\d, n+1}_x - \bar{\mathbf{g}}^{\d}_1 +\tilde{\theta}^{\delta}, \quad (x,y)\in (0,L)\times \mathbb{R},\\[1mm]
		\dis (\sigma-1) p_s \mathbf{d}_{11}\tilde{B}^{\d, n+1}  - \f12 \lambda \v [1+\bm{\chi}]\tilde{u}_s^{\delta} \tilde{B}^{\d,n+1}_{x} \big|_{x=0} \\
		\dis\hspace{5cm}  =- \f12 \lambda [1+\bm{\chi}] \big\{\tilde{u}_s^{\delta} \tilde{\mathbf{w}}^{\d,n+1} - \v \tilde{u}_s  \tilde{\theta}^{\delta}_x\big\}  + \tilde{\mathbf{g}}^{\d}_{b1},\\[1mm]
		\dis \v \tilde{B}^{\d,n+1}_{x} - \lambda \v \f{b_2}{2p_s} \tilde{u}_s^{\delta} \tilde{B}^{\d,n+1}_{yy} \big|_{x=L} = \tilde{\mathbf{w}}^{\d,n+1} + \lambda\v\f{b_2}{2p_s} \tilde{u}_s \tilde{\Xi}^{\delta} + \tilde{\mathbf{g}}^{\d}_{b2},
	\end{cases}
\end{align}
with a given $\tilde{B}^{0}=0$.

By induction arguments,  we can obtain a sequence of solutions $(\mathbf{w}^{\delta,n+1}, \tilde{B}^{\d, n+1}), n=0,1,2,\cdots$ of  \eqref{17.115}-\eqref{17.116} with
\begin{align}\label{17.117}
	&c_3\Big\{\f{\delta}{\v} \|\nabla_{\v}\mathbf{w}^{\delta,n+1}\|^2  +  \f{1}{\v}\|\mathbf{w}^{\delta,n+1}\|^2  +  \|\nabla_{\v}\tilde{B}^{\d,n+1}\|^2+ \|\sqrt{u_s^{\delta}}\mathbf{w}^{\delta,n+1}\|^2_{x=L}  + \v \|\sqrt{u_s^{\delta}} \tilde{B}^{\d,n+1}_y\|^2_{x=L} \Big\} \nonumber\\
	&\quad + \int_0^\infty  \mathbf{d}_{11} \f{p_s}{u^{\delta}_s} |\tilde{B}^{\d,n+1}|^2 dy \Big|_{x=0}  \leq \Lambda \int_0^\infty \mathbf{d}_{11} \f{p_s}{u_s^{\delta}} |\tilde{B}^{\d,n}|^2 dy \Big|_{x=0}  + C \f{1}{\v}\|\hat{\mathcal{G}}^{\d}\|^2  +  (\cdots)_{1,1}.
\end{align}
Similarly, we can obtain
\begin{align}\label{17.118}
	&c_3\Big\{\f{\delta}{\v} \|\nabla_{\v}(\mathbf{w}^{\delta,n+2}-\mathbf{w}^{\delta,n+1})\|^2  +  \f{1}{\v}\|(\mathbf{w}^{\delta,n+2}-\mathbf{w}^{\delta,n+1})\|^2 +  \|\nabla_{\v}(\tilde{B}^{\d,n+2}-\tilde{B}^{\d,n+1})\|^2 \nonumber\\
	&\quad  +   \|\sqrt{u_s^{\delta}}(\mathbf{w}^{\delta,n+2}-\mathbf{w}^{\delta,n+1})\|^2_{x=L} + \v \|\sqrt{u_s^{\delta}} (\tilde{B}^{\d,n+2}-\tilde{B}^{\d,n+1})_y\|^2_{x=L}  \Big\}\nonumber\\
	&\quad  + \int_0^\infty  \mathbf{d}_{11} \f{p_s}{u^{\delta}_s} |\tilde{B}^{\d,n+2}-\tilde{B}^{\d,n+1}|^2 dy \Big|_{x=0} \nonumber\\
	&  \leq \Lambda \int_0^\infty  \mathbf{d}_{11} \f{p_s}{u^{\delta}_s} |\tilde{B}^{\d,n+1}-\tilde{B}^{\d,n}|^2 dy \Big|_{x=0}  
	\leq \cdots \leq \Lambda^{n}\int_0^\infty  \mathbf{d}_{11} \f{p_s}{u^{\delta}_s} |\tilde{B}^{\d,1}-\tilde{B}^{\d,0}|^2 dy \Big|_{x=0} \nonumber\\
	&= C\Lambda^{n}\|\f{\tilde{B}^{\d,1}}{\sqrt{u_s^{\delta}}}\|_{x=0}^2.
\end{align}
and
\begin{align}\label{17.132} 
	&c_4\Big\{\v\|u_s^{\delta}(\nabla_{\v}\tilde{B}^{\d,n+2}_y-\nabla_{\v}\tilde{B}^{\d,n+1}_y)\|^2  + \delta \|u_s^{\delta} (\nabla_{\v}\mathbf{w}^{\delta,n+2}_y-\nabla_{\v}\mathbf{w}^{\delta,n+1}_y)\|^2\nonumber\\
	&\,\,\,  + \|u_s^{\delta}(\mathbf{w}^{\delta,n+2}_y- \mathbf{w}^{\delta,n+1}_y)\|^2  + \|\sqrt{\tilde{u}_s^{\delta}} \big[(\v\tilde{B}^{\d,n+2}_{x}-\tilde{\mathbf{w}}^{\d,n+2})-(\v\tilde{B}^{\d,n+1}_{x}-\tilde{\mathbf{w}}^{\d,n+1})\big]\|^2_{x=L} \nonumber\\
	&\,\,\,  + \v \|(u_s^{\delta})^{\f32} (\mathbf{w}^{\delta,n+2}_y-\mathbf{w}^{\delta,n+1}_y)\|^2_{x=L}\Big\}   +  \v \int_{0}^\infty p_s \mathbf{d}_{11} u_s^{\delta} |\tilde{B}^{\d,n+2}_{y}-\tilde{B}^{\d,n+1}_{y}|^2 dy \Big|_{x=0} \nonumber\\
	&\leq \Lambda \v \int_0^\infty p_s\mathbf{d}_{11} u_s^{\delta} |\tilde{B}^{\d,n+1}_{y}-\tilde{B}^{\d,n}_{y}|^2 dy \Big|_{x=0}    +  \xi \f{\delta}{\v}\|\nabla_{\v}\mathbf{w}^{\delta,n+2} - \nabla_{\v}\mathbf{w}^{\delta,n+1}\|^2 \nonumber\\
	&\,\,\,  +   \xi \f{1}{\v}\|\mathbf{w}^{\delta,n+2} - \mathbf{w}^{\delta,n+1}\|^2   + C_{\xi} a^2\v\Big\{\|\nabla_{\v}(\tilde{B}^{\d,n+2}- \tilde{B}^{\d,n+1})\|^2 +  \|\f{1}{\sqrt{\tilde{u}^{\delta}_{s}}}(\tilde{B}^{\delta,n+2} - \tilde{B}^{\delta,n+1})\|^2_{x=0} \nonumber\\
	&\,\,\, +  \v  \|\sqrt{\tilde{u}_s^{\delta}}(\tilde{B}^{\d,n+2}_y- \tilde{B}^{\d,n+1}_y)\|^2_{x=L}  + \|\f{1}{\sqrt{u^{\delta}_s}}(\tilde{B}^{\d,n+1} - \tilde{B}^{\d,n})\|^2_{x=0} \Big\} .
\end{align}

For later use, we define the notations
\begin{align}
	\mathbb{K}^{\delta}_{2,w}
	:&=c_3\, \Big[\delta \|\nabla_{\v}\mathbf{w}^{\delta}w\|^2  +  \|\mathbf{w}^{\delta}w\|^2  +  \v \|\nabla_{\v}\tilde{B}^{\d}w\|^2+ \v\|\sqrt{u_s^{\delta}}\mathbf{w}^{\delta}w\|^2_{x=L}\nonumber\\
	&\qquad  + \v^2 \|\sqrt{u_s^{\delta}} \tilde{B}^{\d}_yw\|^2_{x=L} \Big]  + C_0 (1-\Lambda) \v \|\f{\tilde{B}^{\d}w}{\sqrt{\tilde{u}_s^{\delta}}}\|_{x=0}^2,\label{17.137}\\
 \mathbb{K}^{\delta}_{3,w}
	:&= c_4\Big[ \delta \|u_s^{\delta} \nabla_{\v}\mathbf{w}^{\delta}_yw\|^2 + \|u_s^{\delta}\mathbf{w}^{\delta}_yw\|^2 + \v\|u_s^{\delta}\nabla_{\v}\tilde{B}^{\d}_yw\|^2  + \|\sqrt{u_s^{\delta}}[\v\tilde{B}^{\d}_{x}-\tilde{\mathbf{w}}^{\d}]w\|^2_{x=L}  \nonumber\\
	&\qquad  + \v \|(u_s^{\delta})^{\f32} \mathbf{w}^{\delta}_{y}w\|^2_{x=L}  \Big] + C_0 (1-\Lambda) \v \|\sqrt{u_s^{\delta}}\tilde{B}^{\d}_yw\|^2_{x=0},\label{17.137-1}\\
	\bar{\mathbb{K}}^{\delta}_{3,w}
	:&=\f{\delta^2}{\v} \|\Delta_{\v}\mathbf{w}^{\delta}w\|^2  + \v\|(U_s^{\delta}\cdot \nabla)\mathbf{w}^{\delta}w\|^2 + \|u_s^{\d}\nabla_{\v}\mathbf{w}^{\delta} w\|^2  + 	\v \|u_s^{\delta}\v \tilde{B}^{\d}_{xx}w\|^2 \nonumber\\
	&\quad + \delta \|\sqrt{u_s^{\delta}}\sqrt{\v} \mathbf{w}^{\delta}_xw\|^2_{x=0}   + \delta \|\sqrt{u_s^{\delta}}  \mathbf{w}^{\delta}_yw\|^2_{x=L}.\label{17.137-2}
\end{align}
With \eqref{17.118}-\eqref{17.132} and contraction mapping arguments, we obtain the following lemma.
\begin{lemma}\label{Blem7.7}
	There exist $(\mathbf{w}^{\delta}, \tilde{B}^{\d})\in H^2\times H^2$ such that 
	\begin{align}
		(\mathbf{w}^{\delta,n}, \tilde{B}^{\d, n})\,\, \to \,\, (\mathbf{w}^{\delta}, \tilde{B}^{\d}) \quad \mbox{strongly in }\,\,\, H^2\times H^2\,\, \mbox{as}\,\, n\to\infty.
	\end{align}
	Moreover, $(\mathbf{w}^{\delta},\tilde{B}^{\d})$ solve the following boundary value equations:
	\begin{align}\label{17.119}
		\begin{cases}
			\dis -\delta \Delta_{\v}\mathbf{w}^{\delta}  + \mathbf{w}^{\delta} + (\mu+\lambda) \v \f{1}{p^{\v,\d}} (U^{\v,\delta}\cdot\nabla)\mathbf{w}^{\delta}  = \hat{\mathcal{G}}^{\d}, \,\, (x,y)\in (0,L)\times \mathbb{R}_+,\\
			\dis -\delta \mathbf{w}^{\delta}_x  + (\mu+\lambda) \f{1}{p_s} u_s^{\delta}\mathbf{w}^{\delta}  \big|_{x=0}=-\sigma \mathbf{d}_{11} \tilde{B}^{\d} ,\\
			\mathbf{w}^{\delta}_x\big|_{x=L}=\mathbf{w}^{\delta}_y\big|_{y=0}=0,
		\end{cases}
	\end{align}
	and 
	\begin{align}\label{17.120}
		\begin{cases}
			\dis -\Delta_{\v} \tilde{B}^{\d} = - \tilde{\mathbf{w}}^{\d}_x - \bar{\mathbf{g}}^{\d}_1 +\tilde{\theta}^{\delta},\quad (x,y)\in (0,L)\times \mathbb{R},\\[1mm]
			\dis  (\sigma-1) p_s \mathbf{d}_{11} \tilde{B}^{\d}  - \v\frac{1}{2}\lambda  [1+\bm{\chi}] \tilde{u}_s^{\delta} \tilde{B}^{\d}_{x}  \Big|_{x=0} = -\frac{1}{2}\lambda  [1+\bm{\chi}] \big(\tilde{u}_s^{\delta} \tilde{\mathbf{w}}^{\d} - \v \tilde{u}_s \tilde{\theta}^{\delta}_x\big) + \tilde{\mathbf{g}}^{\d}_{b1},\\[2mm]
			\dis \v \tilde{B}^{\d}_{x} - \lambda\v\frac{b_2}{2p_s} \tilde{u}_s^{\delta} \tilde{B}^{\delta}_{yy} \Big|_{x=L}  = \tilde{\mathbf{w}}^{\d} + \lambda  \v \f{b_2}{2p_s} \tilde{u}_s \tilde{\Xi}^{\delta} + \tilde{\mathbf{g}}^{\d}_{b2},
		\end{cases}
	\end{align}
	with following uniform estimates:
	\begin{align}
		\f{1}{\v}\mathbb{K}^{\delta}_{2,1}
		&\lesssim  (\cdots)_{1,1} + \f{1}{\v}\|\hat{\mathcal{G}}^{\d}\|^2,\label{17.121}\\
		\mathbb{K}^{\delta}_{3,1}
		&\lesssim \xi \Big\{(\cdots)_{1,1} + \f{1}{\v}\|\hat{\mathcal{G}}^{\d}\|^2\Big\} + \Big\{(\cdots)_{2,1}  + C_{\xi} \|u_s^{\delta}\hat{\mathcal{G}}^{\d}_y\|^2\Big\},\label{17.107-2}
	\end{align}
	where we have taken $\v$ small such that  $\v\ll \xi$.  We remark that the small parameter $\xi$ of $(\cdots)_{1,1}$ is important for us. Since $\tilde{B}^{\d}$ is an odd extension of $\bar{B}^{\g,\d}$, so we have also solved \eqref{17.40}-\eqref{17.41}.
\end{lemma}

\begin{lemma}\label{Blem7.8}
It holds that 
\begin{align}\label{17.111-3}
\bar{\mathbb{K}}^{\delta}_{3,1}
&\lesssim  \Big\{(\cdots)_{1,1}  + \f{1}{\v}\|\hat{\mathcal{G}}^{\d}\|^2\Big\}   +  \Big\{(\cdots)_{2,1} + \|u_s^{\delta}\hat{\mathcal{G}}^{\d}_y\|^2\Big\}.
\end{align}
\end{lemma}

\noindent{\bf Proof.} Using \eqref{17.121}-\eqref{17.107-2},  one gets from \eqref{17.99} that 
\begin{align}\label{17.107-3}
&\f{\delta^2}{\v} \|\Delta_{\v}\mathbf{w}^{\delta}\|^2 + \v\|(U_s^{\delta}\cdot \nabla)\mathbf{w}^{\delta}\|^2 + \delta \|\sqrt{u_s^{\delta}}\sqrt{\v} \mathbf{w}^{\delta}_x\|^2_{x=0}   + \delta \|\sqrt{u_s^{\delta}}  \mathbf{w}^{\delta}_y\|^2_{x=L}\nonumber\\ 
&\lesssim (\cdots)_{1,1} +  \f{1}{\v}\|\hat{\mathcal{G}}^{\d}\|^2   +  (\cdots)_{2,1} +   \|u_s^{\delta}\hat{\mathcal{G}}^{\d}_y\|^2.
\end{align}
Combining \eqref{17.99-1} and \eqref{17.107-3}, one can also obtain
\begin{align}\label{17.110-1}
\|u_s^{\d}\sqrt{\v}\mathbf{w}^{\delta}_x\|^2 &\leq \v \|(U_{s}^{\delta}\cdot \nabla) \mathbf{w}^{\delta}\|^2 + C\v \|u^{\delta}_s \mathbf{w}^{\delta}_{y}\|^2\nonumber\\
&\lesssim  (\cdots)_{1,1} +  \f{1}{\v}\|\hat{\mathcal{G}}^{\d}\|^2   +  (\cdots)_{2,1} +   \|u_s^{\delta}\hat{\mathcal{G}}^{\d}_y\|^2.
\end{align}
We have from \eqref{17.110-1} and \eqref{17.111-1} that 
\begin{align}\label{17.111-2}
\v \|u_s^{\delta}\v \tilde{B}^{\d}_{xx}\|^2 
&\lesssim \v \|u_s^{\delta}\tilde{B}^{\d}_{yy}\|^2 +  \v \|(U_{s}^{\delta}\cdot \nabla) \mathbf{w}^{\delta}\|^2 + \v \|u^{\delta}_s \mathbf{w}^{\delta}_{y}\|^2 + \v  \|u_s \Delta_{\v}\bar{g}_1\|^2 + a^2\v \|\theta^{\d}\|^2 \nonumber\\
&\lesssim  (\cdots)_{1,1} +  \f{1}{\v}\|\hat{\mathcal{G}}^{\d}\|^2   +  (\cdots)_{2,1} +   \|u_s^{\delta}\hat{\mathcal{G}}^{\d}_y\|^2.
\end{align} 
Combining \eqref{17.107-3}-\eqref{17.111-2}, one concludes \eqref{17.111-3}. Therefore the proof of Lemma \ref{Blem7.8} is completed. $\hfill\Box$

\

Next, we consider the weighted estimates. We point out that, for strictness, we should replace the weight $w_0$ by $w=w_0\chi(\f{y}{k})$, then take $k\to +\infty$. We omit the process for simplicity of presentation since it is standard in some sense. 

\begin{lemma}\label{Blem7.9}
Recall $w_0$ in \eqref{6.0-0}-\eqref{6.0}. It holds that 
\begin{align}
\mathbb{K}^{\delta}_{2,w_0}
&\lesssim (\cdots)_{1,\hat{w}_0^{\v}} + \f{1}{\v}\|\hat{\mathcal{G}}^{\d}\hat{w}_0^{\v}\|^2 ,\label{17.145}\\
\v \mathbb{K}^{\delta}_{3,w_0}
&\lesssim \xi \Big\{ (\cdots)_{1,\hat{w}_0^{\v}} + \f{1}{\v}\|\hat{\mathcal{G}}^{\d} \hat{w}_0^{\v}\|^2 \Big\} + C_{\xi} \Big\{ (\cdots)_{2,\hat{w}_0^{\v}} + \|u_s^{\delta}\hat{\mathcal{G}}^{\d}_y\hat{w}_0^{\v}\|^2\Big\}, \label{17.146}\\
\v \bar{\mathbb{K}}^{\delta}_{3,w_0}
&\lesssim \Big\{ (\cdots)_{1,\hat{w}_0^{\v}} + \f{1}{\v}\|\hat{\mathcal{G}}^{\d}\hat{w}_0^{\v}\|^2 \Big\} +   \Big\{ (\cdots)_{2,\hat{w}_0^{\v}} + \|u_s^{\delta}\hat{\mathcal{G}}^{\d}_y\hat{w}_0^{\v}\|^2\Big\}.\label{17.147}
\end{align}
\end{lemma}
 
\noindent{\bf Proof.} Similar arguments as in \eqref{17.110}-\eqref{17.114}, we have from \eqref{17.25}, \eqref{17.107} and \eqref{17.68}-\eqref{17.69} that
\begin{align*}
\mathbb{K}^{\delta}_{2,w_0}  & \lesssim  \|\hat{\mathcal{G}}^{\d}w_0\|^2 + \v \mathbb{K}^{\delta}_{2,1} +  \v^2 \|\f{\tilde{B}^{\d}}{\sqrt{\tilde{u}_s^{\delta}}} w_{0y}\|^2_{x=L} +   \v \|\tilde{B}^{\d} w_{0y}\|^2  + \v (\cdots)_{1,w_0}\nonumber\\
&\lesssim \sqrt{L}\mathbb{K}^{\delta}_{2,w_0}  +L^2 \f{1}{\v} \mathbb{K}^{\delta}_{2,1}   +  \|\hat{\mathcal{G}}^{\d} w_0\|^2 +  \v (\cdots)_{1,w_0},
\end{align*}
which, together with \eqref{17.121} and the smallness of $L\ll1$, yields \eqref{17.145}. 

It follows from \eqref{17.126} and \eqref{17.92} that
\begin{align*}
\v \mathbb{K}^{\delta}_{3,w_0} &\lesssim    \v^2  \mathbb{K}^{\delta}_{3,1} + \xi  \mathbb{K}^{\delta}_{2,w_0}    + \xi \sqrt{\v}  \mathbb{K}^{\delta}_{2,1}  + \xi \|\hat{\mathcal{G}}^{\d} w_0\|^2 +  C_{\xi} \, \v  \|u_s^{\delta}\hat{\mathcal{G}}^{\d}_yw_0\|^2   + C_{\xi} \, \v  (\cdots)_{2,w_0}\nonumber\\
&\lesssim  L^2 \f{1}{\v}\|\hat{\mathcal{G}}^{\d}\|^2  + \xi \|\hat{\mathcal{G}}^{\d} w_0\|^2  + C_{\xi} \v^2 \|u_s^{\delta}\hat{\mathcal{G}}^{\d}_y\|^2+  C_{\xi} \, \v  \|u_s^{\delta}\hat{\mathcal{G}}^{\d}_yw_0\|^2  \nonumber\\
&\quad +  L^2 (\cdots)_{1,1} +  \xi\v (\cdots)_{1,w_0}  +  \v^2(\cdots)_{2,1} + C_{\xi} \, \v  (\cdots)_{2,w_0},
\end{align*}
where we have used \eqref{17.145} and \eqref{17.121}-\eqref{17.107-2} in the last inequality. Thus we obtain \eqref{17.146}.

Using \eqref{17.99}, \eqref{17.99-1} and \eqref{17.111-1}, one can obtain
\begin{align*} 
\v \bar{\mathbb{K}}^{\delta}_{3,w_0}
 &\lesssim \v  \mathbb{K}^{\delta}_{3,w_0} + \mathbb{K}^{\delta}_{2,w_0}    + \sqrt{\v} \mathbb{K}^{\delta}_{2,1}  +  \|\hat{\mathcal{G}}^{\d} w_0\|^2  + \v  (\cdots)_{2,w_0},
 \end{align*}
which, together with \eqref{17.145}-\eqref{17.146}, concludes \eqref{17.147}. Therefore the proof of Lemma \ref{Blem7.9} is completed. $\hfill\Box$

\subsection{Construction of $P^{\delta}$.}\label{ss7.6}
Since we have constructed $\mathbf{w}^{\delta}$ and $\tilde{B}^{\d}$, now we  define a unique $P^{\delta}$  in $(x,y)\in [0,L]\times \mathbb{R}_+$ by solving
\begin{align}\label{17.182-00}
\begin{split}
&P^{\delta}_x = \tilde{B}^{\d} + \theta^{\d} - \bar{g}^{\d}_1,\\
&P^{\delta}_{yy}= \tilde{\mathbf{w}}^{\d}- \v\tilde{B}^{\d}_{x}-\v \theta^{\d}_x + \v \bar{g}^{\d}_{1x},\\
&P^{\delta}(0,0)=0, \quad \lim_{y\to +\infty} P^{\d}=0
\end{split}
\end{align}
Then it is clear to have 
\begin{align}\label{17.182-0}
\begin{cases}
P^{\delta}_{xy}= \tilde{B}^{\d}_y + \theta^{\d}_y - \bar{g}^{\d}_{1y},\\
\sqrt{\v}P^{\delta}_{xx} = \sqrt{\v}\tilde{B}^{\d}_x + \sqrt{\v}\theta^{\d}_x - \sqrt{\v}\bar{g}^{\d}_{1x},\\
\Delta_{\v}P^{\delta}=\tilde{\mathbf{w}}^{\d},
\end{cases}
(x,y)\in [0,L]\times \mathbb{R}_+.
\end{align}

We integrate \eqref{17.182-00} from far field to obtain
\begin{align}\label{17.182-2}
P^{\delta}_y(x,y) = - \int_y^\infty P^{\delta}_{yy}(x,z) dz = - \int_y^\infty \big[\tilde{\mathbf{w}}^{\d}- \v\tilde{B}^{\d}_{x}-\v \theta^{\d}_x + \v \bar{g}^{\d}_{1x}\big](x,z) dz.
\end{align}
Due to the weighted estimates in Lemma \ref{Blem7.9}, \eqref{17.182-2} is well-defined.

Noting $P^{\delta}(0,0)=0$,  we define
\begin{align}\label{17.182-3}
	\begin{split}
		P^{\delta}(x,0)&= \int_0^x P^{\delta}_x(s,0)ds = \int_0^x \big[\theta^{\d} - \bar{g}^{\d}_1\big](s,0)ds,\\
		P^{\delta}(0,y)&=\int_0^y P^{\delta}_y(0,z) dz,
	\end{split}
\end{align}
where we have used the fact $\tilde{B}^{\d}(x,0)\equiv0$ in $\eqref{17.182-3}_1$. Then it is direct to define
\begin{align}
P^{\delta}(x,y)&=P^{\delta}(0,y) + \int_0^x P^{\delta}_x(s,y)ds \nonumber\\
&=\int_0^y P^{\delta}_y(0,z) dz + \int_0^x \big[\tilde{B}^{\d} + \theta^{\d} - \bar{g}^{\delta}_1\big](s,y)ds,
\end{align}
or equivalent
\begin{align}\label{17.168-1}
	P^{\delta}(x,y)&=P^{\delta}(x,0) + \int_0^y P^{\delta}_y(x,z)dz\nonumber\\
	&=\int_0^x \big[\theta^{\d} - \bar{g}^{\delta}_1\big](s,0)ds - \int_0^y \Big[\int_z^\infty \big[\tilde{\mathbf{w}}^{\d}- \v\tilde{B}^{\d}_{x}-\v \theta^{\d}_x + \v \bar{g}^{\delta}_{1x}\big](x,s) ds\Big]dz.
\end{align}

\medskip
 
\begin{lemma}\label{Blem7.10}
Recall \eqref{6.3-2} and \eqref{3.0}.	 For $P^{\d}$ defined in \eqref{17.182-00}-\eqref{17.168-1}, it holds that 
	\begin{align}\label{17.218}
		\begin{split}
			\mathcal{K}^{\delta}_{2,\hat{w}_0^{\v}}(P^{\delta})
			&\lesssim \mathbb{K}^{\delta}_{2,\hat{w}_0^{\v}} +  \v\|||\Phi||\|^2_1 +  \v \|\nabla_{\v}\bar{g}^{\delta}_1\hat{w}_0^{\v}\|^2 + \v^2 \|\sqrt{u_s^{\delta}}\theta^{\d}_y\hat{w}_0^{\v}\|^2_{x=L} + \v^2 \|\sqrt{u_s^{\delta}}\bar{g}^{\delta}_{1y}\hat{w}_0^{\v}\|^2_{x=L},\\
			 \mathcal{K}^{\delta}_{3,\hat{w}_0^{\v}}(P^{\delta})
			&\lesssim  \mathbb{K}^{\delta}_{3,\hat{w}_0^{\v}} + a^2 \mathbb{K}^{\delta}_{2,1}  + \v \|u_s^{\delta}\nabla^2_{\v}\theta^{\d} \hat{w}_0^{\v}\|^2 + \v \|u_s^{\delta}\nabla^2_{\v}\bar{g}^{\delta}_{1}\hat{w}_0^{\v}\|^2 + \v \|\sqrt{u_s^{\delta}}\bar{g}^{\delta}_{1y} \hat{w}_0^{\v}\|^2_{x=0}\\
			&\quad  + \v \|\sqrt{u_s^{\delta}}\, \sqrt{\v} \bar{g}^{\delta}_{1x} w\|^2_{x=L},\\
			\bar{\mathcal{K}}^{\delta}_{3,\hat{w}_0^{\v}}(P^{\delta}) &\lesssim \bar{\mathbb{K}}^{\delta}_{3,\hat{w}_0^{\v}} +  a^2 \mathbb{K}^{\delta}_{2,1},
		\end{split}
	\end{align}
	and
	\begin{align}\label{17.218-0}
		\frac{1}{\v} \mathfrak{K}_{2,w}(P^{\delta})
		&\lesssim \f{1}{\v}\mathcal{K}_{2,w}(P^{\delta}) +   \mathcal{K}_{3,w}(P^{\delta}).
	\end{align}
\end{lemma}

\noindent{\bf Proof.} In the following, we divide the proof into three steps. 

{\it Step 1.} It follows from  \eqref{17.182-00}-\eqref{17.182-0} that 
\begin{align}
\f{1}{\v}\|\nabla_{\v}^2P^{\delta}w\|^2 &\lesssim \|\nabla_{\v}\tilde{B}^{\d}w\|^2 + \f{1}{\v}\|\mathbf{w}^{\delta}w\|^2 + \|||\Phi||\|^2_1 +  \|\nabla_{\v}\bar{g}^{\delta}_1w\|^2 \nonumber\\
&\lesssim \f{1}{\v}\mathbb{K}^{\delta}_{2,w} +  \|||\Phi||\|^2_1 +  \|\nabla_{\v}\bar{g}^{\delta}_1w\|^2.\label{17.182-1}
\end{align}

For the third order derivatives, we note from \eqref{17.182-0}  that 
\begin{align}
	\begin{split}
		\nabla_{\v}P^{\delta}_{xy}&= \nabla_{\v}\tilde{B}^{\d}_y + \nabla_{\v}\theta^{\d}_y - \nabla_{\v}\bar{g}^{\delta}_{1y},\\
		\nabla_{\v}P^{\delta}_{yy}&= \nabla_{\v}\tilde{\mathbf{w}}^{\d} - \v \nabla_{\v}\tilde{B}^{\d}_{x}-\v  \nabla_{\v}\theta^{\d}_{x} + \v \nabla_{\v}\bar{g}^{\delta}_{1x},
	\end{split}
\end{align}
which yields that 
\begin{align}\label{17.212}
	\|u_s^{\delta}\nabla^{2}_{\v}P^{\delta}_{y}w\|^2 
	&\lesssim \|u_s^{\delta} \tilde{\mathbf{w}}^{\delta}_yw\|^2  + \v \|u_s^{\delta}\nabla_{\v}\tilde{B}^{\d}_yw\|^2  + \v \|u_s^{\delta}\nabla_{\v}\theta^{\d}_yw\|^2 + \v \|u_s^{\delta}\nabla^2_{\v}\bar{g}^{\delta}_{1}w\|^2 \nonumber\\
	&\lesssim \mathbb{K}^{\delta}_{3,w} + a^2\mathbb{K}^{\delta}_{2,1} + \v \|u_s^{\delta}\nabla^2_{\v}\theta^{\d} w\|^2 + \v \|u_s^{\delta}\nabla^2_{\v}\bar{g}^{\delta}_{1}w\|^2,
\end{align}
and
\begin{align}\label{17.212-1}
	\|u_s^{\delta}\nabla_{\v}\Delta_{\v}P^{\delta} w\|^2& = \|u_s^{\delta}\nabla_{\v}\tilde{\mathbf{w}}^{\d} w\|^2 \lesssim \|u_s^{\delta}\nabla_{\v}\mathbf{w}^{\delta} w\|^2+  a^2  \|\mathbf{w}^{\delta}\| \lesssim  a^2 \mathbb{K}^{\delta}_{2,1}.
\end{align}

{\it Step 2.} For trace at $x=0,L$, it follows from  \eqref{17.182-0} that 
\begin{align}
	\begin{split}
		\|\sqrt{u_s^{\delta}}\Delta_{\v}P^{\delta}w\|^2_{x=L}&=\|\sqrt{u_s^{\delta}}\tilde{\mathbf{w}}^{\delta}w\|^2_{x=L}\leq  \|\sqrt{u_s^{\delta}}\mathbf{w}^{\delta}w\|^2_{x=L},\\
		\v \|\sqrt{u_s^{\delta}}P^{\delta}_{xy}w\|^2_{x=L} &\lesssim \v \|\sqrt{u_s^{\delta}}\tilde{B}^{\d}_yw\|^2_{x=L} + \v \|\sqrt{u_s^{\delta}}\theta^{\d}_yw\|^2_{x=L} + \v \|\sqrt{u_s^{\delta}}\bar{g}^{\delta}_{1y}w\|^2_{x=L},\\
		\v \|\sqrt{u_s^{\delta}}P^{\delta}_{xy}w\|^2_{x=0} &\lesssim \v \|\sqrt{u_s^{\delta}}\tilde{B}^{\d}_yw\|^2_{x=0}   + \v \|\sqrt{u_s^{\delta}}\bar{g}^{\delta}_{1y}w\|^2_{x=0},\\
		\|\sqrt{u_s^{\delta}}P^{\delta}_{yy}w\|^2_{x=L} &\lesssim\|\sqrt{u_s^{\delta}}[\v\tilde{B}^{\d}_{x}-\tilde{\mathbf{w}}^{\d}]w\|^2_{x=L}+ \v \|\sqrt{u_s^{\delta}}\, \sqrt{\v} \bar{g}^{\delta}_{1x}w\|^2_{x=L},
	\end{split}
\end{align}
which  yields that 
	\begin{align}\label{17.214}
		\|\sqrt{u_s^{\delta}}\Delta_{\v}P^{\delta}w\|^2_{x=L} + \v \|\sqrt{u_s^{\delta}}P^{\delta}_{xy}w\|^2_{x=L}
		&\lesssim \f{1}{\v}\mathbb{K}^{\delta}_{2,w} + \v \|\sqrt{u_s^{\delta}}\theta^{\d}_yw\|^2_{x=L} + \v \|\sqrt{u_s^{\delta}}\bar{g}^{\delta}_{1y}w\|^2_{x=L},
	\end{align}
	and
	\begin{align}\label{17.214-1}
		& \|\sqrt{u_s^{\delta}}P^{\delta}_{yy}w\|^2_{x=L}  + \v \|\sqrt{u_s^{\delta}}P^{\delta}_{xy}w\|^2_{x=0} + \v \|u_s^{\delta}\sqrt{u_s^{\delta}} \Delta_{\v}P^{\delta}_y w\|^2_{x=L}\nonumber\\
		&\lesssim \mathbb{K}^{\delta}_{3,w} + a^2   \mathbb{K}^{\delta}_{2,1}  + \v \|\sqrt{u_s^{\delta}}\bar{g}^{\delta}_{1y}w\|^2_{x=0} + \v \|\sqrt{u_s^{\delta}}\, \sqrt{\v} \bar{g}^{\delta}_{1x}w\|^2_{x=L}.
\end{align}

{\it Step 3.} Combining  \eqref{17.182-1}, \eqref{17.212} and \eqref{17.214}-\eqref{17.214-1}, we can conclude \eqref{17.218}.  For \eqref{17.218-0}, we need only to use Newtonian-Leibniz formula to control $\|\sqrt{u_s} P^{\delta}_{yy}w\|_{x=0}$. Therefore the proof of Lemma \ref{Blem7.10} is completed. $\hfill\Box$

\medskip

\begin{lemma}\label{Blem7.11}
	We denote 
	\begin{align*}
		\omega_i(y)=\la y\ra \la Y\ra^{\mathfrak{l}_i}=(1+y) (1+\sqrt{\v}y)^{\mathfrak{l}_i} \,\,\, \mbox{for}\,\, \mathfrak{l}_i\geq 0,\,\, i=1,2,\cdots.
	\end{align*} 
	For the first order derivatives of $P^{\delta}$, we have
	\begin{align}\label{17.205}
		\begin{split}
			\v\|\sqrt{\v}P^{\delta}_xw\|^2 &\lesssim L^2  \mathbb{K}^{\delta}_{2,w} + \v^2 \|||\Phi||\|^2_1 + \v^2 \|\bar{g}_1w\|^2,\\
			\v \|P^{\delta}_y\omega_1\|^2 &\lesssim \mathbb{K}^{\delta}_{2,\hat{w}_0} + \v \|||\Phi||\|^2_1 + \v \|\nabla_{\v}\bar{g}_1\hat{w}_0\|^2,\quad \mathfrak{l}_1\leq \mathfrak{l}-1,\\
			\v^{\beta_1}\|P^{\delta}_y\|^2 &\lesssim  \mathbb{K}^{\delta}_{2,\hat{w}_0} + \v\|||\Phi||\|^2_1 + \v\|\nabla_{\v}\bar{g}_1\hat{w}_0\|^2 
			,\,\, \forall \, \beta_1 \in (0, 1).
		\end{split}
	\end{align}
	Recall \eqref{3.0-1}, it holds that 
	\begin{align}\label{17.156-0} 
		\mathscr{K}(P^{\delta})
		&\lesssim \f{1}{\v}\mathbb{K}^{\delta}_{2,\hat{w}^{\v}_0} 
		+\|||\Phi||\|^2_1  + \|\nabla_{\v}\bar{g}_1\hat{w}^{\v}_0\|^2+ \|(\mathbf{d}_2 \bar{g}_1, \bar{g}_1\mathbf{1}_{\{y\in[0,2]\}})\|^2.
	\end{align}
	We remark that \eqref{17.156-0} plays a key role in the iteration.
\end{lemma}

\noindent{\bf Proof.} In the following, we divide the proof into two steps. 

{\it Step 1.} We have from \eqref{17.68-1} that
\begin{align}\label{17.168-2}
	\v\|\tilde{B}^{\d}w\|^2& \lesssim  L\v \|\f{\tilde{B}^{\delta}}{\sqrt{\tilde{u}_s^{\delta}}}w\|^2_{x=0} + L^2  \|\sqrt{\v}\tilde{B}^{\delta}_{x}w\|^2 \lesssim L^2 \f{1}{\v} \mathbb{K}^{\delta}_{2,w}.
\end{align}
It follows from \eqref{17.168-2} and \eqref{17.182-00} - \eqref{17.182-0} that 
\begin{align} \label{17.182}
	\|\sqrt{\v}P^{\delta}_xw\|^2 &\lesssim \|\sqrt{\v}\tilde{B}^{\d}w\|^2 + \v \|(\theta^{\d}, \bar{g}^{\delta}_1)w\|^2  \lesssim L^2 \f{1}{\v}\mathbb{K}^{\delta}_{2,w} + \v \|||\Phi||\|^2_1 + \v \|\bar{g}^{\delta}_1w\|^2,
\end{align}
which yields $\eqref{17.205}_1$.

It follows from \eqref{17.182-2} that 
\begin{align}\label{17.180}
	\|\omega_1 P^{\delta}_y\|^2
	&=\int_0^L\int_0^\infty |\omega_1(y) |^2 \Big[\int_y^\infty \big[\tilde{\mathbf{w}}^{\d}- \v \tilde{B}^{\d}_x - \v \theta^{\d}_x + \v \bar{g}^{\delta}_{1x}\big](x,z) dz\Big]^2 dydx 
	\nonumber\\
	&\lesssim \int_0^\infty |\omega_1(y)|^2 \int_y^\infty |\omega_2(z)|^{-2}dz \, dy \, \|(\tilde{\mathbf{w}}^{\d}, \v\tilde{B}^{\d}_x, \v \theta^{\d}_x, \v \bar{g}^{\delta}_{1x})\omega_2\|^2 \nonumber\\
	&\lesssim  \f{1}{\v}\mathbb{K}^{\delta}_{2,\omega_2} +  \|||\Phi||\|^2_1 + \|\nabla_{\v}\bar{g}^{\delta}_1\omega_2\|^2 ,
\end{align}
where we have demanded $\mathfrak{l}_2>\mathfrak{l}_1 +1$ and  used the facts
\begin{align}
	\begin{split}
		\int_y^\infty (1+z)^{-2} (1+\sqrt{\v}z)^{-2\mathfrak{l}} dz &\lesssim (1+y)^{-1} (1+\sqrt{\v}y)^{-2\mathfrak{l}} ,\\
		\int_y^\infty (1+\sqrt{\v}z)^{-2\mathfrak{l}} dz &\lesssim \v^{-\f12} (1+\sqrt{\v}y)^{-2\mathfrak{l}+1}.
	\end{split}
\end{align}
Thus $\eqref{17.205}_2$ is proved.
We point out that, at this step, it is hard to have estimate on $\|P^{\delta}_y w_0\|$ due to \eqref{17.180} since we loss some space weight.

On the other hand, we have from \eqref{17.182-2} that 
\begin{align*} 
	\|P^{\delta}_y\|^2 &=\int_0^L\int_0^\infty\Big[\int_y^\infty \big[\tilde{\mathbf{w}}^{\d} - \v \tilde{B}^{\d}_x - \v \theta^{\d}_x + \v \bar{g}^{\delta}_{1x}\big](x,z) dz\Big]^2 dydx \nonumber\\
	&\lesssim \int_0^\infty \int_y^\infty |\omega_1(z)|^{-2}dz \, dy \, \|(\mathbf{w}^{\delta},\v\tilde{B}^{\d}_x, \v \theta^{\d}_x, \v \bar{g}^{\delta}_{1x})\omega_1\|^2 \nonumber\\
	&\lesssim \v^{-\beta_1} \|(\mathbf{w}^{\delta},\v\tilde{B}^{\d}_x, \v \theta^{\d}_x, \v \bar{g}^{\delta}_{1x})\omega_1\|^2 \nonumber\\
	&\lesssim \v^{1-\beta_1} \Big[\f{1}{\v}\mathbb{K}^{\delta}_{2,1} + \f{1}{\v}\mathbb{K}^{\delta}_{2,w_0} + \|||\Phi||\|^2_1 + \|\nabla_{\v}\bar{g}^{\delta}_1\hat{w}_0\|^2\Big]
	,\quad \forall \, \beta_1 \in (0, 1),
\end{align*}
which yields that  $\eqref{17.205}_3$.

\smallskip

{\it Step 2.} From \eqref{17.182}-\eqref{17.180},  we may not have good enough $\v$-decay estimate for $P^{\delta}_x, P^{\delta}_y$. Luckily, one may only need $\|\mathbf{d}_2P^{\delta}_x\|^2 $ and $\|\mathbf{d}_{2}P^{\delta}_y\|^2 $ to close the iteration from $\hat{P}^{\delta}\to P^{\delta}$.
In fact, a direct calculation shows that 
\begin{align}\label{17.197}
	\|\mathbf{d}_2P^{\delta}_x\|^2 
	\lesssim  \|\f{\tilde{B}^{\d}}{y}\|^2 + \|\mathbf{d}_2\theta^{\d}\|^2 +  \|\mathbf{d}_2 \bar{g}^{\delta}_1\|^2  
	&\lesssim  \f{1}{\v}\mathbb{K}^{\delta}_{2,1} + a  \|||\Phi||\|^2_1 +  \|\mathbf{d}_2 \bar{g}^{\delta}_1\|^2 ,
\end{align}
and
\begin{align}\label{17.198}
	\|\mathbf{d}_2P^{\delta}_y\|^2 &\lesssim \int_0^L\int_0^\infty |\mathbf{d}_2(y)|^2 \Big[\int_y^\infty \big[\tilde{\mathbf{w}}^{\d}- \v \tilde{B}^{\d}_x - \v \theta^{\d}_x + \v \bar{g}^{\delta}_{1x}\big](x,z) dz\Big]^2 dydx \nonumber\\
	&\lesssim  \int_0^\infty   |\mathbf{d}_2(y)|^2 \int_y^\infty (1+z)^{-2\beta} dz  dy  \|(\tilde{\mathbf{w}}^{\d} - \v \tilde{B}^{\d}_x - \v \theta^{\d}_x + \v \bar{g}^{\delta}_{1x})\omega_1^{\beta}\|^2\nonumber\\
	&\lesssim \v^{1-\beta} \|(\mathbf{w}^{\delta}, \v \tilde{B}^{\d}_x, \v \theta^{\d}_x, \v \bar{g}^{\delta}_{1x})w\|^2 + \v^{1-\beta}\|(\f{\mathbf{w}^{\delta}}{\sqrt{\v}}, \sqrt{\v} \tilde{B}^{\d}_x, \sqrt{\v} \theta^{\d}_x, \sqrt{\v} \bar{g}^{\delta}_{1x})\|^2 \nonumber\\
	&\lesssim \v^{1-\beta} \big[\mathbb{K}^{\delta}_{2,w} + \f{1}{\v}\mathbb{K}^{\delta}_{2,1}\big] +  \v^{1-\beta} \|||\Phi||\|^2_1 + \v^{1-\beta} \|\nabla_{\v}\bar{g}^{\delta}_1 \hat{w}^{\v}\|^2 , 
\end{align}
where we choose $\beta \in(\f12,1),\, \mathfrak{l}_1\geq 0$. That means we can gain some $\v$ decay for $\|\mathbf{d}_2P^{\delta}_y\|^2$.
Hence we have from \eqref{17.197}-\eqref{17.198} that 
\begin{align}\label{17.199}
\|\mathbf{d}_2 P^{\delta}_x\|^2 + \v^{-1+\beta}\|\mathbf{d}_2P^{\delta}_y\|^2  
&\lesssim \f{1}{\v}\mathbb{K}^{\delta}_{2,1} + \mathbb{K}^{\delta}_{2,w_0} +\|||\Phi||\|^2_1 +  \|\nabla_{\v}\bar{g}^{\delta}_1 \hat{w}^{\v}_0\|^2  + \|\mathbf{d}_2 \bar{g}^{\delta}_1\|^2 .
\end{align}

Similarly one can obtain
\begin{align}\label{17.199-0}
&\|P^{\delta}_x\mathbf{1}_{\{y\in[0,2]\}}\|^2  + \v^{-1+\beta} \|P^{\delta}_y\mathbf{1}_{\{y\in[0,2]\}}\|^2 \nonumber\\
&\lesssim  \f{1}{\v}\mathbb{K}^{\delta}_{2,1} + \mathbb{K}^{\delta}_{2,w_0} +\|||\Phi||\|^2_1 +  \|\nabla_{\v}\bar{g}^{\delta}_1 \hat{w}^{\v}_0\|^2  + \|\bar{g}^{\delta}_1\mathbf{1}_{\{y\in[0,2]\}}\|^2.
\end{align}

Recall  \eqref{3.0-1}. Then combining \eqref{17.199}-\eqref{17.199-0}, we conclude \eqref{17.156-0}. Therefore the proof of Lemma \ref{Blem7.11} is completed. $\hfill\Box$

\medskip


\begin{lemma}\label{Blem7.12}
It holds that 
\begin{align}\label{17.165}
	(\cdots)_{1,w}
	&\lesssim a [\Phi]_{3,w} + \|||\Phi||\|^2_1 
	+ a^2\|(0,\hat{P}, 0)\hat{w}\|^2_{\mathbf{Y}^{\v}} + a^2\|(\phi,p, \zeta)\hat{w}\|^2_{\mathbf{X}^{\v}} + \sqrt{\v} \|u_s\nabla_{\v}^2\bar{g}_1w\|^2\nonumber\\
	&\quad    + \|\nabla_{\v}\bar{g}_{1}w\|^2  +  \v \|\sqrt{u_s} \bar{g}_{1x}w\|_{x=L}^2  + \f{1}{\v} \|\sqrt{u_s} (b_2g_2)_{y}w\|^2_{x=L} + \v \|\sqrt{u_s} \hat{\mathbf{S}}_{xy}w\|^2_{x=L}\mathbf{1_{D}} \nonumber\\
	&\quad  + \v  \|\sqrt{u_s}\sqrt{\v}\tilde{T}_{xy}w\|^2_{x=0}   + \|\sqrt{u_s}\Delta_{\v}\tilde{T}w\|^2_{x=0} + \v \|u_s(\Delta_{\v}\tilde{T}, \tilde{T}_{yy})w\|^2_{x=L} \nonumber\\
	&\quad  +  \|\f{\mathfrak{N}_{12}}{\sqrt{u_s}}w\|^2_{x=0}  + \v \|\f{\mathfrak{N}_{22}}{\sqrt{u_s}} w\|^2_{x=L},
\end{align} 
and
\begin{align}\label{17.166-7}
(\cdots)_{2,w}&\lesssim  a^3[\Phi]_{3,\hat{w}}  + a^2\v \|||\Phi||\|^2_1 + L^2\|(0,\hat{P}, 0)\hat{w}\|^2_{\mathbf{Y}^{\v}} + a^4\|(\phi,p, \zeta)\hat{w}\|^2_{\mathbf{X}^{\v}}
+ \sqrt{\v} \|\nabla_{\v}\bar{g}_{1}\hat{w}\|^2   \nonumber\\
&\quad + \sqrt{\v} \|u_{s}\nabla^2_{\v}\bar{g}_1 \hat{w}\|^2    + \|\sqrt{u_{s}}(b_2g_2)_{y}w\|_{x=L}^2  + a^2\|\nabla_{\v}^2\tilde{T}\hat{w}\|^2   + \|u_s\nabla_{\v}\Delta_{\v}\tilde{T}w\|^2 \nonumber\\
&\quad + \v \|u_s\nabla_{\v}^2\tilde{T}_yw\|^2 + \|u_s \nabla_{\v}^2\hat{\mathbf{S}}_{y}\hat{w}\|^2\mathbf{1_D} +\v \|\hat{\mathbf{S}}_{xy}\hat{w}\|^2\mathbf{1_D} + \|(\mathfrak{N}_{12}, \pa_{\v}\mathfrak{N}_{12})w\|^2 \nonumber\\
&\quad   +   \v \|(\mathfrak{N}_{22}, \pa_{\v}\mathfrak{N}_{22})\hat{w}\|^2   + \v^2 \|\sqrt{u_s}(b_2\mathfrak{h}_2, b_2\mathfrak{N}_{21})_{y}w\|_{x=L} +a^2\v \|\f{\mathfrak{N}_{12}}{\sqrt{u_s}}w\|^2_{x=0}.
\end{align}

\end{lemma}
 
\noindent{\bf Proof.} We divide the proof into three steps.

{\it Step 1. Estimates on $\theta^{\d}\, \&\, \tilde{\theta}^{\d}$.} 
We note that 
\begin{align}\label{17.155}
\begin{split}
\sqrt{\v}\|\sqrt{u_s}\theta^{\d}_yw\|^2_{x=L}
&\lesssim  \|\sqrt{\v} u_s\theta^{\d}_{xy}w\|\cdot \|\theta^{\d}_yw\| + a\sqrt{\v}\|\theta^{\d}_yw\|^2
\lesssim a\big\{\|||\Phi||\|^2_1 + L^{\f14}[\Phi]_{3,1}\big\},\\
\sqrt{\v} \|\sqrt{\v}\sqrt{u_s}\theta^{\d}_xw\|^2_{x=0} 
&\lesssim  \|\v u_s\theta^{\d}_{xx}w\|\cdot \|\sqrt{\v}\theta^{\d}_xw\| + a\sqrt{\v}\|\sqrt{\v}\theta^{\d}_xw\|^2
\lesssim a\big\{\|||\Phi||\|^2_1 + L^{\f14}[\Phi]_{3,1}\big\},\\
\sqrt{\v} \|\sqrt{u^{\delta}_{s}}\theta^{\d}_yw\|_{x=L}^2 
&\lesssim \big\{\|\sqrt{\v} u_s\theta^{\d}_{xy}w\| +  \delta \|\sqrt{\v} \theta^{\d}_{xy}w\| \big\}\|\theta^{\d}_yw\| + a\sqrt{\v}\|\theta^{\d}_yw\|^2\\
&\lesssim a\big\{\|||\Phi||\|^2_1 + L^{\f14}[\Phi]_{3,1}\big\},
\end{split}
\end{align}
and
\begin{align}\label{17.158}
	\sqrt{\v} \|\v\tilde{u}_s\tilde{\theta}^{\delta}_xw\|^2_{H^{1/2}_{x=0}} &\lesssim \sqrt{\v} \|\big( (\v\tilde{u}_s\tilde{\theta}^{\delta}_xw)_y, \v\tilde{u}_s\tilde{\theta}^{\delta}_xw\big)\| \big\{ \|\v\pa_x(\tilde{u}_s\tilde{\theta}^{\delta}_xw)\| + \f{1}{L} \|\v\tilde{u}_s\tilde{\theta}^{\delta}_xw\|
	\big\} \nonumber\\
	&\lesssim \v \|u_s\nabla_{\v}^2\theta w\|^2 + a^2 \v \|\nabla_{\v}\theta \hat{w}\|^2
	\lesssim a^2\v \big\{\|||\Phi||\|^2_1 + L^{\f14}[\Phi]_{3,1}\big\}.
\end{align}
Hence we have from \eqref{17.155}-\eqref{17.158} and \eqref{17.27} that
\begin{align}\label{17.159}
\begin{split}
\|y\theta^{\d} \hat{w}\|^2 + \v^{\f32} \|\sqrt{u_s}\theta^{\d}_x \hat{w}\|^2_{x=0}  + \sqrt{\v} \|\sqrt{u_s}\theta^{\d}_y\hat{w}\|^2_{x=L} 
&\lesssim \|||\Phi||\|^2_1 + L^{\f14}[\Phi]_{3,1},\\
\v \|(a\theta^{\d},\,  a\theta^{\d}_y,\, u_s\nabla_{\v}^2\theta^{\d})\hat{w}\|^2 +  \v^2 \|\sqrt{u^{\delta}_{s}}\theta^{\d}_y\hat{w}\|_{x=L}^2 
&\lesssim a^2\v \|||\Phi||\|^2_1 + \v L^{\f14}[\Phi]_{3,1}.
\end{split}
\end{align}

\smallskip

{\it Step 2. Estimate on $(\cdots)_{1,w}$.} Recall $(\cdots)_{1,w}$ in \eqref{17.71-1}.  A direct calculation shows that 
\begin{align}
\v \|\sqrt{u^{\delta}_s}(b_2\mathfrak{g}_2)_{y}w\|^2_{x=L} &\lesssim L \v [\phi]_{3,w}.
\end{align}
Recall \eqref{17.19-7} and  \eqref{17.31},
then we have 
\begin{align}\label{17.161}
\|\f{\mathbf{g}^{\d}_{b1}}{\sqrt{\tilde{u}_s^{\delta}}} w\|^2_{x=0} 
&\lesssim a [\Phi]_{3,w} + \v^{\f{1-\beta}{2}}\|(0,\hat{P}, 0)\hat{w}\|^2_{\mathbf{Y}^{\v}} + a^2\|(\phi,p, \zeta)\hat{w}\|^2_{\mathbf{X}^{\v}} + \|\sqrt{u_s}\Delta_{\v}\tilde{T}w\|^2_{x=0}  \nonumber\\
&\quad + \v  \|\sqrt{u_s}\sqrt{\v}\tilde{T}_{xy}w\|^2_{x=0}  + \v \|\sqrt{u_s}\sqrt{\v}\bar{g}_{1x}w\|^2_{x=0}   +  \|\f{\mathfrak{N}_{12}}{\sqrt{u_s}}w\|^2_{x=0},
\end{align}
where we have used 
\begin{align}\label{17.162}
\begin{split}
\|\f{\mathfrak{g}_1w}{\sqrt{u^{\delta}_s}}\|^2_{x=0}&\lesssim  \v \|\sqrt{\v}q_xw\|^2_{x=0} \lesssim \sqrt{L}\v [\phi]_{3,w},\\
\|\f{G_{11}(\tilde{T})}{\sqrt{u_s^{\delta}}}w\|^2_{x=0} &\lesssim \|\sqrt{u_s}\Delta_{\v}\tilde{T}w\|^2_{x=0} + \v  \|\sqrt{u_s}\sqrt{\v}\tilde{T}_{xy}w\|^2_{x=0},\\
\|\f{G_{12}(p,q)}{\sqrt{u_s^{\delta}}}w\|^2_{x=0}&\lesssim \f{1}{\sqrt{\v}} [[p]]_{2,1}.
\end{split}
\end{align}

Recall  \eqref{17.19-4} and \eqref{17.51}, 
we have  
\begin{align}\label{17.163}
\f{1}{\v}\|\f{\hat{G}^{\d}_3}{\sqrt{u^{\delta}_s}}\|^2_{x=L}
&\lesssim a [\Phi]_{3,w} +  \v  \mathcal{K}_{3,w}(\hat{P}) + \mathcal{K}_{2,w}(\hat{P}) + \f{1}{\sqrt{\v}} [[p]]_{2,1}  + L \v [\phi]_{3,w} \nonumber\\
&\quad + \v \|\sqrt{u_s} \hat{\mathbf{S}}_{xy}w\|^2_{x=L}\mathbf{1_{D}} + \v \|u_s(\Delta_{\v}\tilde{T}, \tilde{T}_{yy})w\|^2_{x=L} + \v \|\f{\mathfrak{N}_{22}}{\sqrt{u_s}} w\|^2_{x=L} ,
\end{align}
where we have used the following facts
\begin{align}\label{17.164}
\begin{split}
\v \|\f{G_{21}(\tilde{T})}{\sqrt{u^{\delta}_s}}w\|^2_{x=L} 
 &\lesssim \v \|u_s (\Delta_{\v}\tilde{T}, \tilde{T}_{yy})w\|^2_{x=L},\\
 \v \|\f{G_{22}(p,q)}{\sqrt{u^{\delta}_s}}w\|^2_{x=L} 
 & \lesssim \f{1}{\sqrt{\v}} [[p]]_{2,1}  +  L [\phi]_{3,w}.
\end{split} 
\end{align}
 
We also note that 
\begin{align}\label{17.164-0}
	\v  \|\sqrt{u_s}\nabla_{\v}\bar{g}^{\d}_{1}w\|^2_{x=0,L}  &\lesssim \sqrt{\v} \|u_s\nabla_{\v}^2\bar{g}_1w\|^2 + \sqrt{\v} \|\nabla_{\v}\bar{g}_1w\|^2.
\end{align}

Substituting \eqref{17.159}-\eqref{17.163}   into  \eqref{17.71-1} and noting \eqref{17.164-0},  one concludes \eqref{17.165}. Here we have used similar estimates as in $\eqref{17.20-2}$ for $\bar{g}^{\d}_1$, $g^{\d}_2$, $\cdots$.

\smallskip

{\it Step 3. Estimate on $(\cdots)_{2,w}$.} Recall  $(\cdots)_{2,w}$ in  \eqref{17.99-10}.
 It follows from \eqref{17.31} and \eqref{17.19-7} that 
\begin{align}
\|\mathbf{g}^{\d}_{b1}w\|^2 
&\lesssim L^{\f14} [\Phi]_{3,w} + L^2\|(0,\hat{P}, 0)\hat{w}\|^2_{\mathbf{Y}^{\v}} + L^{\f18}\|(\phi,p, \zeta)\hat{w}\|^2_{\mathbf{X}^{\v}}  + a^2 \v \|\nabla_{\v}\bar{g}_1w\|^2  \nonumber\\
&\quad  + a^2\|\nabla_{\v}^2\tilde{T}w\|^2  + \|\mathfrak{N}_{12}w\|^2,\label{17.166-2}\\
\|\nabla_{\v}\mathbf{g}^{\d}_{b1}w\|^2 
&\lesssim  a^3 [\Phi]_{3,w}+  L^2\|(0,\hat{P}, 0)\hat{w}\|^2_{\mathbf{Y}^{\v}} + a^4\|(\phi,p, \zeta)\hat{w}\|^2_{\mathbf{X}^{\v}} + \v \|u_s\nabla_{\v}^2\bar{g}_1w\|^2 \nonumber\\
&\quad  + a^2\v \|\nabla_{\v}\bar{g}_1w\|^2   + \|u_s\nabla_{\v}\Delta_{\v}\tilde{T}w\|^2 + \v \|u_s\nabla_{\v}^2\tilde{T}_yw\|^2  + a^2\|\nabla_{\v}^2\tilde{T}w\|^2 \nonumber\\
&\quad + \|(\mathfrak{N}_{12}, \pa_{\v}\mathfrak{N}_{12})w\|^2,\label{17.166-1}
\end{align}
and
\begin{align}\label{17.166-0}
\f{1}{\gamma^2} \|\bar{\chi}'(\f{y}{\gamma}) \mathbf{g}^{\d}_{b1}w\|^2
&\lesssim  L^{\f18} [\Phi]_{3,1} + L^2\|(0,\hat{P}, 0)\|^2_{\mathbf{Y}^{\v}} + L^{\f18}\|(\phi,p, \zeta)\|^2_{\mathbf{X}^{\v}} + a^2\|\nabla_{\v}^2\tilde{T}\|^2  \nonumber\\
&\quad + a^2 \v \|\nabla_{\v}\bar{g}_1\|^2  +  \|\mathfrak{N}_{12y}\|^2,
\end{align}
where we have used the Hardy inequality and the facts $\mathfrak{N}_{12}\big|_{y=0}=0$ in \eqref{17.166-0}.

Noting  \eqref{17.51} and using \eqref{Z.11-2}, we have from \eqref{17.166-2}-\eqref{17.166-0} that  
\begin{align}\label{17.166-4}
\sqrt{\v} \|\tilde{\mathbf{g}}^{\d}_{b1}w\|^2_{H^{1/2}_{x=0}}
&\lesssim \|\nabla_{\v}\mathbf{g}^{\d}_{b1}w\|^2 + \f{1}{\gamma^2} \|\bar{\chi}'(\f{y}{\gamma}) \mathbf{g}^{\d}_{b1}w\|^2 + \|\mathbf{g}^{\d}_{b1}\|^2 + \f{\v}{L^2} \|\mathbf{g}^{\d}_{b1}w\|^2 \nonumber\\
&\lesssim    a^3 [\Phi]_{3,w}+  L^2\|(0,\hat{P}, 0)\hat{w}\|^2_{\mathbf{Y}^{\v}} + a^4\|(\phi,p, \zeta)\hat{w}\|^2_{\mathbf{X}^{\v}} + \v \|u_s\nabla_{\v}^2\bar{g}_1w\|^2 \nonumber\\
&\quad  + a^2\v \|\nabla_{\v}\bar{g}_1\hat{w}\|^2   + \|u_s\nabla_{\v}\Delta_{\v}\tilde{T}w\|^2 + \v \|u_s\nabla_{\v}^2\tilde{T}_yw\|^2  + a^2\|\nabla_{\v}^2\tilde{T}\hat{w}\|^2 \nonumber\\
&\quad   + \|(\mathfrak{N}_{12}, \pa_{\v}\mathfrak{N}_{12})w\|^2.
\end{align}
 
\smallskip

A direct calculation shows that 
\begin{align}\label{17.170-1}
\begin{split}
\|G_{21}w\|^2&\lesssim   a^2\|\nabla_{\v}^2 \tilde{T}w\|^2 , \\
\|\nabla_{\v}G_{21}w\|^2 &\lesssim \|v_s\nabla_{\v}(\Delta_{\v}\tilde{T}, \tilde{T}_{yy})w\|^2 
 + a^2\|\nabla_{\v}^2 \tilde{T}w\|^2,\\
\|G_{22}w\|^2&\lesssim a^2 \f{1}{\v} [[p]]_{2,1} + a^2 L^{\alpha} [\phi]_{3,1},\\
\|\nabla_{\v}G_{22}w\|^2 &\lesssim a^2 \f{1}{\v} [[p]]_{2,1} + a^{2} [\phi]_{3,\hat{w}},
\end{split}
\end{align}
which, together with \eqref{17.19-10} and \eqref{17.51}, yields that 
\begin{align}
\|\hat{G}^{\d}_3 w\|^2 
&\lesssim L^{\f18} \v [\Phi]_{3,w} + \v^2\|(0,\hat{P}, 0)\hat{w}\|^2_{\mathbf{Y}^{\v}} + a^2\v^2\|(\phi,p, 0)\hat{w}\|^2_{\mathbf{X}^{\v}}  +  a^2\v^2 \|\nabla_{\v}^2 \tilde{T}w\|^2\nonumber\\
&\,\,\, + \v^2 \|\mathfrak{N}_{22}w\|^2 +a^2\v^2\|\hat{\mathbf{S}}_{xy}w\|^2\mathbf{1_D} , \label{17.170-2}\\
\|\nabla_{\v}\hat{G}^{\d}_3w\|^2
&\lesssim a^3\v [\Phi]_{3,w}  + \v^2\|(0,\hat{P}, 0)\hat{w}\|^2_{\mathbf{Y}^{\v}} + a^2\v^2\|(\phi,p, 0)\hat{w}\|^2_{\mathbf{X}^{\v}}  + \v^2 \|v_s\nabla_{\v}(\Delta_{\v}\tilde{T}, \tilde{T}_{yy})w\|^2 \nonumber\\
&\,\,\, + a^2\v^2 \|\nabla_{\v}^2 \tilde{T}w\|^2  +   \v^2 \|(\mathfrak{N}_{22}, \pa_{\v}\mathfrak{N}_{22})w\|^2 + \v\|u_s \nabla_{\v}^2\hat{\mathbf{S}}_{y}w\|^2\mathbf{1_D} +\v^2 \|\hat{\mathbf{S}}_{xy}w\|^2\mathbf{1_D}, \label{17.170-3}
\end{align}
and
\begin{align}\label{17.170-5}
\|\f{1}{\gamma}\bar{\chi}'(\f{y}{\gamma})\hat{G}^{\d}_3\|^2 
&\lesssim  L^{\f18} \v [\Phi]_{3,1} + \v^2\|(0,\hat{P}, 0)\|^2_{\mathbf{Y}^{\v}} + a^2\v^2\|(\phi,p, 0)\|^2_{\mathbf{X}^{\v}}   +  a^2\v^2 \|\nabla_{\v}^2 \tilde{T}\|^2\nonumber\\
&\quad  + \v^2 \|\mathfrak{N}_{22y}\|^2 +a^2\v^2\|\hat{\mathbf{S}}_{xy}w\|^2\mathbf{1_D},
\end{align}
where we  used the Hardy inequality and the fact $\mathfrak{N}_{22}\big|_{y=0}=0$ in \eqref{17.170-5}.

Hence it follows from \eqref{17.170-2}-\eqref{17.170-5} that
\begin{align}\label{17.166-5}
\f{1}{\sqrt{\v}} \|\bar{\chi}(\f{y}{\gamma})\hat{G}^{\d}_3 w\|^2_{H^{1/2}_{x=L}}
& \lesssim \f{1}{\v} \|\nabla_{\v}\hat{G}^{\d}_3w\|^2 + \f{1}{\v} \|\f{1}{\g}\bar{\chi}'(\f{y}{\g})\hat{G}^{\d}_3w\|^2 + \f{1}{\v}\|\hat{G}^{\d}_3 w\|^2 
\nonumber\\
&\lesssim a^3  [\Phi]_{3,w}  + \v\|(0,\hat{P}, 0)\hat{w}\|^2_{\mathbf{Y}^{\v}} + a^2\v\|(\phi,p, 0)\hat{w}\|^2_{\mathbf{X}^{\v}}  + a^2\v \|\nabla_{\v}^2 \tilde{T}\hat{w}\|^2 \nonumber\\
&\,\,\, + \v \|v_s\nabla_{\v}(\Delta_{\v}\tilde{T}, \tilde{T}_{yy})\hat{w}\|^2  + \|u_s \nabla_{\v}^2\hat{\mathbf{S}}_{y}\hat{w}\|^2\mathbf{1_D} +\v \|\hat{\mathbf{S}}_{xy}\hat{w}\|^2\mathbf{1_D}\nonumber\\
&\,\,\, +   \v \|(\mathfrak{N}_{22}, \pa_{\v}\mathfrak{N}_{22})\hat{w}\|^2.
\end{align}

It is clear to have  
\begin{align}\label{17.166-6}
\sqrt{\v} \|u_s\bar{g}^{\d}_{1y}\bar{\chi}(\f{y}{\gamma}) w\|^2_{H^{1/2}_{x=L}} 
&\lesssim \sqrt{\v} \|\big(u_s\bar{g}^{\d}_{1y}\bar{\chi}(\f{y}{\gamma})w, \big(u_s\bar{g}^{\d}_{1y}\bar{\chi}(\f{y}{\gamma})w\big)_y\big)\| \Big\{\|\big(u_s\bar{g}^{\d}_{1y}w\big)_x\| + \f{1}{L}\|u_s\bar{g}^{\d}_{1y}w\|\Big\}\nonumber\\
& \lesssim \|u_s\nabla_{\v}^2\bar{g}_1w\|^2 + a^2 \|\bar{g}_{1y}\hat{w}\|^2.
\end{align}
We also note
\begin{align}\label{17.166-8}
\begin{split}
\v\|\sqrt{u_s}\Delta_{\v}\tilde{T}w\|^2_{x=0} 
&\lesssim \|\nabla_{\v}^2\tilde{T}w\|^2 +   \v\|u_s\nabla_{\v}\Delta_{\v}\tilde{T}w\|^2,\\
\v^2  \|\sqrt{u_s}\sqrt{\v}\tilde{T}_{xy}w\|^2_{x=0}
&\lesssim \v \|\nabla_{\v}^2\tilde{T}w\|^2 + \v^2\|u_s \nabla_{\v}^2\tilde{T}_yw\|^2.
\end{split}
\end{align}

Combining  above estimates and \eqref{17.99-10}, we conclude \eqref{17.166-7}. Here we also used similar  estimates as in $\eqref{17.20-2}$ for $\bar{g}^{\d}_1$, $g^{\d}_2$, $\cdots$.
Therefore the proof of Lemma \ref{Blem7.12} is completed. $\hfill\Box$

\medskip

\begin{lemma}\label{Blem7.13}
It holds that 
\begin{align}\label{17.198-0}
	\f{1}{\v}\|\hat{\mathcal{G}}^{\d}w\|^2 &\lesssim  L^{\f18} [\Phi]_{3,w}  + a^2 \|(0,\hat{P}, \hat{\mathbf{S}})\hat{w}\|^2_{\mathbf{Y}^{\v}} 
	+ a^3\|(\phi,p, \zeta)\hat{w}\|^2_{\mathbf{X}^{\v}} +  \|\nabla_{\v}g_{1}w\|^2 + \f{1}{\v}\|g_{2y}w\|^2 \nonumber\\
	&\quad    + \v \|u_s \nabla\Delta_{\v}\tilde{T}w\|^2  +  \v \|\nabla_{\v}^2\tilde{T}\hat{w}\|^2  + \v  \|(\mathcal{N}_{1x}, \mathcal{N}_{2y})w\|^2,
\end{align}
and
\begin{align}\label{17.200}
\|u_s^{\delta}\hat{\mathcal{G}}^{\d}_yw\|^2 
&\lesssim a^4 \v [\Phi]_{3,1}  + a^3 \v  \|||\Phi||\|^2_1 +  a^2 \|(0,\hat{P}, \hat{\mathbf{S}})\hat{w}\|^2_{\mathbf{Y}^{\v}} + a^4\|(\phi,p, \zeta)\hat{w}\|^2_{\mathbf{X}^{\v}}\nonumber\\
&\quad + a^2\v^2 \|u_s \nabla_{\v}\Delta_{\v}\tilde{T}_yw\|^2 + a^4 \v \|\nabla_{\v}\Delta_{\v}\tilde{T}w\|^2 + a^2 \v \|u_s\nabla^2_{\v}\tilde{T}_yw\|^2  + \v^2 \|\nabla^2_{\v}\tilde{T}w\|^2  \nonumber\\
&\quad + \v \|u_s\sqrt{\v} \pa_{xy}g_{1}w\|^2 + \|u_s\pa_{yy}g_{2}w\|^2 + \v^2 \|u_s (\mathcal{N}_{1xy}, \mathcal{N}_{2yy})w\|^2 .
\end{align}
\end{lemma}

\noindent{\bf Proof.} We divide the proof into two steps. 

{\it Step 1.} Noting the definition of $\tilde{d}_{22}$ in \eqref{7.33-1}, we have $\tilde{d}_{22y}=\v \big(\f{d_{22}}{T_s}\big)_y\cong \v\mathbf{d}_2$.  Then we have from \eqref{17.17-1} and \eqref{17.20-1} that 
\begin{align}\label{17.170}
\begin{split}
\|\mathcal{G}_1(\hat{P}^{\d})w\|^2
&\lesssim \v \|\mathbf{d}_2(\sqrt{\v}\hat{P}_x,\sqrt{\v}\hat{P}_y)w\|^2 + a^2 \|\nabla_{\v}^2\hat{P}w\|^2,\\
\|\mathcal{G}_2(\tilde{T})w\|^2 & \lesssim  \v^2 \|u_s \nabla\Delta_{\v}\tilde{T}w\|^2 +\v^2 \|\nabla_{\v}^2\tilde{T}\hat{w}\|^2, \\
\|\mathcal{G}_3(\Phi)w\|^2 &
\lesssim L^{\f1{8}} \v [\Phi]_{3,w}, 
\end{split}
\end{align}
and
\begin{align}\label{17.175}
\begin{split}
\|\mathcal{G}_4(p^{\d},S^{\d},\phi)w\|^2 &\lesssim a^2 \v [[p]]_{2,1}   + \v^2 [\phi]_{3,1} + \v \|\sqrt{\v}g_{1x}w\|^2 + \|g_{2y}w\|^2,\\
\|\mathcal{G}_5(\hat{\mathbf{S}}^{\d})w\|^2&\lesssim \v^2 [[[\hat{\mathbf{S}}]]]_{2,1}, \\
 \|\mathcal{G}_6(\zeta^{\d})w\|^2 
&\lesssim  a^3\v [[[\zeta]]]_{2,w},\\
\|\mathcal{G}_7(q)w\|^2 
&\lesssim  L^{\f1{8}} \v [\Phi]_{3,w},\\
\|\mathcal{G}_8(p^{\d},S^{\d},\phi) w\|^2&\lesssim \v^2 \|(\mathcal{N}_{1x}, \mathcal{N}_{2y})w\|^2.
\end{split}
\end{align}
From \eqref{17.170}-\eqref{17.175} and \eqref{17.17-1}, we prove \eqref{17.198-0}. Here we also used similar  estimates as in $\eqref{17.20-2}$ for $\bar{g}^{\d}_1$, $g^{\d}_2$, $\cdots$.

\smallskip

{\it Step 2.} Recall \eqref{17.17-1} and \eqref{17.20-1}.  For $\|u_s^{\delta}\hat{\mathcal{G}}^{\d}_{y}\|$, we notice that 
\begin{align}\label{17.181-1}
\begin{split}
\|u_s^{\delta}\pa_y\mathcal{G}_1(\hat{P}^{\d})w\|^2 
&\lesssim a^2\|u_s \nabla_{\v}^2\hat{P}_yw\|^2 + a^2 \|\nabla_{\v}^2\hat{P}\|^2 + \v  \|\mathbf{d}_2(\sqrt{\v}\hat{P}_x,\sqrt{\v}\hat{P}_y)\|^2 ,\\
\|u_s^{\delta}\pa_y\mathcal{G}_2(\tilde{T})w\|^2 &\lesssim \v^2 \|u_s^{\delta}u_s \nabla\Delta_{\v}\tilde{T}_yw\|^2 + a^4 \v \|\nabla_{\v}\Delta_{\v}\tilde{T}w\|^2 + a^2 \v^2 \|u_s\nabla^2_{\v}\tilde{T}_yw\|^2 \\
&\quad  + a^4 \v^2 \|\nabla^2_{\v}\tilde{T}w\|^2, \\
\|u_s^{\delta}\pa_y\mathcal{G}_{3}(\Phi)w\|^2 &\lesssim a^3 \v  \|||\Phi||\|^2_1 + a^4 \v [\Phi]_{3,1},\\
\|u_s^{\delta}\pa_y\mathcal{G}_4(p^{\d},S^{\d}, \phi)w\|^2 &\lesssim  \v [[p]]_{3,1} +  \v [[p]]_{2,1} + \v [\phi]_{3,1} +  \v \|u_s\sqrt{\v} \pa_{xy}g_{1}w\|^2 + \|u_s \pa_{yy}g_{2}w\|^2,
\end{split}
\end{align}
and
\begin{align}\label{17.181-2}                                                                                    
\begin{split}     
\|u_s^{\delta}\pa_y\mathcal{G}_5(\hat{\mathbf{S}}^{\d})w\|^2 &\lesssim \v [[[\hat{\mathbf{S}}]]]_{3,1} + a^2\v^2 [[[\hat{\mathbf{S}}]]]_{2,1}  ,\\
\|u_s^{\delta}\pa_y\mathcal{G}_6(\zeta^{\d})w\|^2&\lesssim  a^4 [[[\zeta]]] _{3,w}  + \v  [[[\zeta]]]_{2,1},\\
\|u_s^{\delta}\pa_y\mathcal{G}_{7}(q)w\|^2 &\lesssim  \v^2 [\phi]_{3,1},\\
\|u_s^{\delta}\pa_y\mathcal{G}_8(p^{\d},S^{\d},\phi) w\|^2 &\lesssim \v^2 \|u_s  (\mathcal{N}_{1xy}, \mathcal{N}_{2yy})w\|^2.
\end{split}                                            
\end{align}
We remark that the lift up  $q^{\delta}$ in \eqref{17.17-2} makes the estimate of $\eqref{17.181-1}_4$ rigorously.

Combining \eqref{17.181-1}-\eqref{17.181-2} and \eqref{17.17-1}, we conclude \eqref{17.200}. Here we have used similar  estimates as in $\eqref{17.20-2}$ for $\bar{g}^{\d}_1$, $g^{\d}_2$, $\cdots$. Therefore the proof of Lemma \ref{Blem7.13} is completed.   $\hfill\Box$

\begin{remark}
	We remark that it is in the proof of \eqref{17.181-1}-\eqref{17.181-2} that we need to introduce the lift up  $(p^{\d},S^{\d}, \zeta^{\d})$, otherwise $\|u_s^{\delta}\sqrt{\v} \pa_{xy}g_{1}\|\cong \|u^{\d}_s \sqrt{\v} p_{xyy}\| +\cdots$ and some other highest  order derivatives may not make sense in our framework.
\end{remark}

\begin{proposition}\label{Bprop7.11}
Recall \eqref{17.137}-\eqref{17.137-2} and \eqref{3.0}. For $P^{\delta}$ constructed in Section \ref{ss7.6},  we have
\begin{align}
&\f{1}{\v}\mathbb{K}^{\delta}_{2,\hat{w}_0^{\v}}+ \f{1}{\v}\mathcal{K}^{\delta}_{2,\hat{w}_0^{\v}}(P^{\delta}) \lesssim  a [\Phi]_{3,\hat{w}_0^{\v}}  + \|||\Phi||\|^2_1 + a^2 \|(0, \hat{P}, \hat{\mathbf{S}})\hat{w}_0^{\v}\|^2_{\mathbf{Y}^{\v}}  + \{\cdots\}_{1,\hat{w}_0^{\v}} \nonumber\\
&\qquad\qquad\qquad \qquad \qquad  + \v \|\sqrt{u_s} \hat{\mathbf{S}}_{xy}\hat{w}_0^{\v}\|^2_{x=L}\mathbf{1_{D}},\label{17.204}\\
&\mathbb{K}^{\delta}_{3,\hat{w}_0^{\v}}  +  \mathcal{K}^{\delta}_{3,\hat{w}_0^{\v}}(P^{\delta})
\lesssim (\xi a + C_{\xi}a^3) [\Phi]_{3,\hat{w}_0^{\v}} + \xi \|||\Phi||\|^2_1 + C_{\xi}a^2 \|(0, \hat{P}, \hat{\mathbf{S}})\hat{w}_0^{\v}\|^2_{\mathbf{Y}^{\v}} +  \xi \{\cdots\}_{1,\hat{w}_0^{\v}} \nonumber\\
&\qquad\qquad\qquad\qquad\quad + C_{\xi} \{\cdots\}_{2,\hat{w}_0^{\v}} + \xi\v \|\sqrt{u_s} \hat{\mathbf{S}}_{xy}\hat{w}_0^{\v}\|^2_{x=L}\mathbf{1_{D}} + C_{\xi}\|u_s \nabla_{\v}^2\hat{\mathbf{S}}_{y}\hat{w}_0^{\v}\|^2\mathbf{1_D},\label{17.206}\\
&\bar{\mathbb{K}}^{\delta}_{3,\hat{w}_0^{\v}}  + \bar{\mathcal{K}}^{\delta}_{3,\hat{w}_0^{\v}}(P^{\delta}) 
 \lesssim a [\Phi]_{3,\hat{w}_0^{\v}} +  \|||\Phi||\|^2_1 +  a^2 \|(0, \hat{P}, \hat{\mathbf{S}})\hat{w}_0^{\v}\|^2_{\mathbf{Y}^{\v}}  + \{\cdots\}_{1,\hat{w}_0^{\v}}  + \{\cdots\}_{2,\hat{w}_0^{\v}}\nonumber\\
 &\qquad\qquad \qquad\qquad \quad   + \v \|\sqrt{u_s} \hat{\mathbf{S}}_{xy}\hat{w}_0^{\v}\|^2_{x=L}\mathbf{1_{D}} +  \|u_s \nabla_{\v}^2\hat{\mathbf{S}}_{y}\hat{w}_0^{\v}\|^2\mathbf{1_D},
\label{17.205-1} 
\end{align}
where $\{\cdots\}_{1,w}$ and $\{\cdots\}_{2,w}$ are the ones defined in \eqref{17.217}-\eqref{17.218-1}.
\end{proposition}

\noindent{\bf Proof.} We have  Lemmas \ref{Blem7.7}-\ref{Blem7.9} that 
\begin{align}\label{17.215}
\begin{split}
	\f{1}{\v}\mathbb{K}^{\delta}_{2,\hat{w}_0^{\v}} 
&\lesssim (\cdots)_{1,{\hat{w}_0^{\v}}} + \f{1}{\v}\|\hat{\mathcal{G}}^{\d}\hat{w}_0^{\v}\|^2 ,\\
\tilde{\mathbb{K}}^{\delta}_{3,\hat{w}_0^{\v}} 
&\lesssim \xi \big\{(\cdots)_{1,{\hat{w}_0^{\v}}} + \f{1}{\v}\|\hat{\mathcal{G}}^{\d}\hat{w}_0^{\v}\|^2 \big\} + C_{\xi} \big\{(\cdots)_{2,{\hat{w}_0^{\v}}} + \|u_s^{\delta}\hat{\mathcal{G}}^{\d}_y\hat{w}_0^{\v}\|^2\big\}, \\
\bar{\mathbb{K}}^{\delta}_{3,\hat{w}_0^{\v}} 
&\lesssim \big\{(\cdots)_{1,{\hat{w}_0^{\v}}} + \f{1}{\v}\|\hat{\mathcal{G}}^{\d}\hat{w}_0^{\v}\|^2 \big\} + C_{\xi} \big\{(\cdots)_{2,{\hat{w}_0^{\v}}} + \|u_s^{\delta}\hat{\mathcal{G}}^{\d}_y\hat{w}_0^{\v}\|^2\big\}.
\end{split}
\end{align}
It follows from Lemmas \ref{Blem7.12} \& \ref{Blem7.13} that 
\begin{align}\label{17.216}
\begin{split}
(\cdots)_{1,{\hat{w}_0^{\v}}} + \f{1}{\v}\|\hat{\mathcal{G}}^{\d}\hat{w}_0^{\v}\|^2 
&\lesssim  a [\Phi]_{3,\hat{w}_0^{\v}}  + \|||\Phi^{\delta}||\|^2_1 + a^2 \|(0, \hat{P}, \hat{\mathbf{S}})\hat{w}_0^{\v}\|^2_{\mathbf{Y}^{\v}}  + \{\cdots\}_{1,\hat{w}_0^{\v}},\\
(\cdots)_{2,{\hat{w}_0^{\v}}} + \|u_s^{\delta}\hat{\mathcal{G}}^{\d}_y\hat{w}_0^{\v}\|^2
&\lesssim a^3[\Phi]_{3,\hat{w}_0^{\v}}  + a^2\v \|||\Phi^{\delta}||\|^2_1 + a^2 \|(0, \hat{P}, \hat{\mathbf{S}})\hat{w}_0^{\v}\|^2_{\mathbf{Y}^{\v}}    + \{\cdots\}_{2,\hat{w}_0^{\v}}.
\end{split}
\end{align} 

By similar estimates as in $\eqref{17.20-2}$, then we have from  \eqref{17.218} that
\begin{align}\label{17.219}
\begin{split}
\mathcal{K}^{\delta}_{2,\hat{w}_0^{\v}}(P^{\delta})
&\lesssim \f{1}{\v}\mathbb{K}^{\delta}_{2,\hat{w}_0^{\v}} +  \|||\Phi||\|^2_1 + L^{\f14} [\Phi]_{3,1}  + \{\cdots\}_{1,{\hat{w}_0^{\v}}},\\
\mathcal{K}^{\delta}_{3,\hat{w}_0^{\v}}(P^{\delta})
&\lesssim \mathbb{K}^{\delta}_{3,\hat{w}_0^{\v}} + a^2 \mathbb{K}^{\delta}_{2,1} + a^2\v \|||\Phi||\|^2_1 +  L^{\f14} \v  [\Phi]_{3,1} + \{\cdots\}_{2,{\hat{w}_0^{\v}}},\\
\bar{\mathcal{K}}^{\delta}_{3,\hat{w}_0^{\v}}(P^{\delta}) &\lesssim \bar{\mathbb{K}}^{\delta}_{3,\hat{w}_0^{\v}} +  a^2 \mathbb{K}^{\delta}_{2,1},
\end{split}
\end{align}
which, together with \eqref{17.215}-\eqref{17.216}, concludes \eqref{17.204}-\eqref{17.205-1}.
Therefore the proof of Proposition \ref{Bprop7.11} is completed.   
$\hfill\Box$

\smallskip
\subsection{A refined estimate for $\|u_s\Delta_{\v}P_y^{\delta} w\|$} 
In this subsection we aim to derive a refined estimate for $\|u_s\Delta_{\v}P_y^{\delta} w\|$, which is crucial to close the estimate for  DT case.

\begin{lemma}\label{CBlem7.6}
	Let $(\mathbf{w}^{\delta}, \tilde{B}^{\d})$ be the solution of \eqref{17.119}-\eqref{17.120}. It holds that 
	\begin{align}\label{C17.223}
		\|u_s^{\d} \Delta_{\v}P^{\d}_yw\|&\lesssim \v^{-\f18}\mathbb{K}^{\d}_{2,\hat{w}} + \sqrt{\v} \bar{\mathbb{K}}^{\d}_{3,\hat{w}} + a^3 \|(\Phi,0,0)\|^2_{\mathbf{X}^{\v}} + a^2 \|(0,\hat{P}, \hat{\mathbf{S}})w\|^2_{\mathbf{Y}^{\v}} +  \{\cdots\}_{2,\hat{w}}.
	\end{align}
	The key point is that small coefficients are available for terms on RHS of \eqref{C17.223}.
\end{lemma}

\noindent{\bf Proof.} In the following, we divide the proof into three steps. 

{\it Step 1.} Similar as Lemma \ref{lemP7.3}, one has
\begin{align}\label{C17.126}
	&(1-\xi)\Big\{\delta \|u_s^{\delta} \nabla_{\v}\mathbf{w}^{\delta}_y w\|^2 +  \|u_s^{\delta}\mathbf{w}^{\delta}_yw\|^2 +  (\mu+\lambda) \v \int_0^\infty  \f{1}{2p_s} ( u_s^{\delta})^3  |\mathbf{w}^{\delta}_y|^2 w^2 dy \Big|_{x=L} \Big\}\nonumber\\
	&\leq \Big[\f{\sigma^2}{2(\mu+\lambda)}+\xi\Big]\v \int_0^\infty p_s\mathbf{d}_{11}^2 u_s^{\delta} |\tilde{B}^{\delta}_{y}|^2 w^2  dy \Big|_{x=0} + \v^{\f32} \|(U_{s}^{\delta}\cdot \nabla) \mathbf{w}^{\delta}w\|^2  + \v \|u_s^{\delta} \mathbf{w}^{\delta}_y\|^2 \nonumber\\
	&\quad   + \sqrt{\v}\|\mathbf{w}^{\delta}w\|^2  + C_{\xi} a^2 \delta \|\mathbf{w}^{\delta}_y\hat{w}\|^2  +   C_{\xi} a^4 \v \|\f{\tilde{B}^{\d}w}{\sqrt{u^{\delta}_s}}\|^2_{x=0}  +  C_{\xi} \|u_s^{\delta}\hat{\mathcal{G}}^{\delta}_yw\|^2.
\end{align}

\smallskip

{\it Step 2.} Denote $\chi_{L}:=\chi(\f{4x}{L})$.	
Applying $\pa_y$ to $\eqref{17.120}_1$, then multiplying the resultant equation  by $\v \,(\tilde{u}_s^{\delta}  w \chi_L)^2\tilde{B}^{\d}_y$, and using similar arguments as in Lemma \ref{Blem7.6}, one obtains that 
\begin{align}\label{C17.82}
	&(1-\xi)\v\|\tilde{u}_s^{\delta}\nabla_{\v}\tilde{B}^{\d}_yw  \chi_{L} \|^2   
	+	(1-\xi)\v \int_{\mathbb{R}} \f{2(\sigma-1)}{\lambda(1+\bm{\chi})} p_s \mathbf{d}_{11} \tilde{u}_s^{\delta}|\tilde{B}^{\d}_{y}|^2w^2 dy \Big|_{x=0} \nonumber\\
	&\leq \|\tilde{u}_s^{\delta} \tilde{\mathbf{w}}^{\d}_{y} w\|\cdot \sqrt{\v}\|\sqrt{\v}\tilde{u}_s^{\delta} \tilde{B}^{\d}_{xy} w \chi_L\| + \v  \|\tilde{u}_s^{\delta} \tilde{\mathbf{w}}^{\d}_{y} w\|   +  \v^{1-} \|\tilde{B}^{\d}_y\hat{w}\|^2 +  C_{\xi} a^2 \v \|\f{\tilde{B}^{\d}\hat{w}}{\sqrt{\tilde{u}^{\delta}_{s}}} \|^2_{x=0}\nonumber\\
	&\quad  + C_{\xi} \v \|u_{s}\Delta_{\v}\bar{g}^{\d}_1 \hat{w}\|^2   + C a^2\v \|(\theta^{\d}_y, \theta^{\d} )w\|^2  + C_{\xi}a^2 \v \|\f{\hat{w}}{\sqrt{\tilde{u}^{\delta}_s}}\big(\v\tilde{u}_s\tilde{\theta}^{\delta}_x, \tilde{\mathbf{g}}^{\d}_{b1}\big)\|^2_{x=0} \nonumber\\
	&\quad + C_{\xi} \sqrt{\v} \|(\v\tilde{u}_s\tilde{\theta}^{\delta}_x, \f{\tilde{\mathbf{g}}^{\d}_{b1}}{1+\bm{\chi}})w\|^2_{H^{1/2}_{x=0}}.
\end{align}
The key point is that the boundary condition $\eqref{17.120}_3$ at $x=L$ can be ignored. Then we can skip the effect of $\v u_s \hat{\mathbf{S}}_{xy}|_{x=L}$ in  DT case.

\smallskip

{\it Step 3.} From \eqref{C17.126}-\eqref{C17.82}, taking $\xi>0$ suitably small and by similar arguments as in \eqref{17.110}-\eqref{17.113}, one has
\begin{align}\label{C17.226}
	&\v\|\tilde{u}_s^{\delta}\nabla_{\v}\tilde{B}^{\d}_y\hat{w} \chi_{L} \|^2 
	+ \delta \|u_s^{\delta} \nabla_{\v}\mathbf{w}^{\delta}_y \hat{w}\|^2 +  \|u_s^{\delta}\mathbf{w}^{\delta}_y\hat{w}\|^2 
	\nonumber\\
	&\lesssim \v^{-\f18}\mathbb{K}^{\d}_{2,\hat{w}} + \sqrt{\v} \bar{\mathbb{K}}^{\d}_{3,\hat{w}}  + a^3 \|(\Phi,0,0)\|^2_{\mathbf{X}^{\v}} +  \{\cdots\}_{3,\hat{w}}.
\end{align}
		It follows from $\eqref{17.182-0}_3$  that 
		\begin{align*}
			\|u_s^{\d} \Delta_{\v}P^{\d}_y \hat{w}\|^2=\|u_s^{\d} \tilde{\mathbf{w}}^{\d}_y \hat{w}\|^2 \lesssim \|u_s^{\d}\mathbf{w}^{\d}_y \hat{w}\|^2 + Ca^2 \|\mathbf{w}^{\d}\|^2,
		\end{align*}
		which, with \eqref{C17.226}, concludes \eqref{C17.223}.
		Therefore the proof of Lemma \ref{CBlem7.6} is completed. $\hfill\Box$

\subsection{Proof of Theorem \ref{thm5.1}.} Since the proof is long, we divide it into four steps.

{\it Step 1.}  Recall the definitions of $\Phi^{\d}, \hat{P}^{\d}, \hat{\mathbf{S}}^{\d}$ and $p^{\d}, S^{\d}, q^{\d}$.   For $ \mathbf{w}^{\delta}, P^{\delta}, \tilde{B}^{\delta}$ constructed above, we shall consider  the limit  $\delta\to 0$.   Firstly it is clear to have the following convergences
\begin{align}\label{18.23}
	\begin{split}
		&P^{\delta} \to P, \quad \mbox{weakly in}\,\, H^2_{loc}([0,L]\times\mathbb{R}_{+}),	\\
		&\mathbf{w}^{\delta} \to \mathbf{w}\equiv \Delta_{\v}P,\,\,\quad \mbox{weakly in}\,\, L^2_{loc}([0,L]\times\mathbb{R}_{+}),\\
		&\tilde{B}^{\delta} \to \tilde{B},\,\,\quad \mbox{weakly in}\,\, H^1_{loc}([0,L]\times\mathbb{R}),
	\end{split}
	\quad \mbox{as}\,\, \delta \to 0+,
\end{align}
and
\begin{align}\label{18.23-1}
	\begin{split}
		&\Phi^{\delta} \to \Phi, \quad \mbox{weakly in}\,\, H^4_{loc}([0,L]\times\mathbb{R}_{+}),\\
		&\hat{\mathbf{P}}^{\delta} \to \hat{\mathbf{P}} ,\,\,\quad \mbox{weakly in}\,\, H^1_{loc}([0,L]\times\mathbb{R}_+),\\
		&\hat{\mathbf{S}}^{\delta} \to \hat{\mathbf{S}} ,\,\,\quad \mbox{weakly in}\,\, H^1_{loc}([0,L]\times\mathbb{R}_+),\\
		&\theta^{\d} \to \theta, \quad \mbox{weakly in}\,\, H^1_{loc}([0,L]\times\mathbb{R}_{+}),\\	
		&u_s\theta^{\d} \to u_s\theta, \quad \mbox{weakly in}\,\, H^2_{loc}([0,L]\times\mathbb{R}_{+}),
	\end{split}
	\quad \mbox{as}\,\, \delta \to 0+.
\end{align}
We also have similar convergence for $p^{\d}, S^{\d}, q^{\d}$ as $\d\to0+$.
Hereafter the local convergence always means local in $y$. 
We remark that $\theta$ is the solution  of \eqref{17.26} but replacing $\f{\mu}{\tilde{\rho}^{\v,\d}} \tilde{\Phi}^{\delta}_{yyy}$ by $\f{\mu}{\rho^{\v}} \Phi_{yyy}$.	 

It follows from \eqref{18.23}-\eqref{18.23-1} and \eqref{17.182-00}-\eqref{17.182-0} that 
\begin{align}\label{18.24}
	\begin{cases}
		P_x = \tilde{B} + \theta - \bar{g}_1,\\
		P_{xy}= \tilde{B}_y + \theta_y - \bar{g}_{1y},\\
		\sqrt{\v}P_{xx} = \sqrt{\v}\tilde{B}_x + \sqrt{\v}\theta_x - \sqrt{\v}\bar{g}_{1x},	\\
		P_{yy}=\mathbf{w} - \v\tilde{B}_{x}-\v \theta_x + \v \bar{g}_{1x},\\
		\mathbf{w}=\Delta_{\v}P,
	\end{cases}
	(x,y)\in [0,L]\times \mathbb{R}_+.
\end{align}
Noting 
\begin{align}
	\begin{split}
		\|u_s\mathbf{w}^{\delta}\|^2_{H^1} &\lesssim \|\mathbf{w}^{\delta}\|^2 + \|u_s\nabla \mathbf{w}^{\delta}\|^2 ,\\
		\|u_s\nabla \tilde{B}^{\d}\|^2_{H^1} &\lesssim \|\nabla \tilde{B}^{\delta}\|^2 + \|u_s\nabla^2 \tilde{B}^{\delta}\|^2,\\
		\|u_s\nabla^2_{\v}P^{\delta}\|^2_{H^1}&\lesssim \|\nabla_{\v}^2P^{\delta}\|^2 + \|u_s\nabla^3_{\v}P^{\delta}\|^2,\\
		\|u_s\nabla^2_{\v}\hat{P}^{\delta}\|^2_{H^1}&\lesssim \|\nabla_{\v}^2\hat{P}^{\delta}\|^2 + \|u_s\nabla^3_{\v}\hat{P}^{\delta}\|^2,\\	\|u_s\nabla^2_{\v}\hat{\mathbf{S}}^{\delta}\|^2_{H^1}&\lesssim \|\nabla_{\v}^2\hat{\mathbf{S}}^{\delta}\|^2 + \|u_s\nabla^3_{\v}\hat{\mathbf{S}}^{\delta}\|^2,
	\end{split}
\end{align}
we have 
\begin{align}\label{17.232-1}
	\begin{split}
		&u_s\mathbf{w}^{\delta} \to u_s\mathbf{w},\quad \mbox{strongly in }\,\, L^2_{loc}([0,L]\times\mathbb{R}_+)\,\, as\,\, \delta \to0,\\
		&u_s\nabla_{\v}\tilde{B}^{\delta} \to u_s\nabla_{\v}\tilde{B},\quad \mbox{strongly in }\,\, L^2_{loc}([0,L]\times\mathbb{R}_+)\,\, as\,\, \delta \to0,\\
		&u_s\nabla^2_{\v}P^{\delta} \to u_s\nabla^2_{\v}P,\quad \mbox{strongly in }\,\, L^2_{loc}([0,L]\times\mathbb{R}_+)\,\, as\,\, \delta \to0,\\
		&u_s\nabla^2_{\v}\hat{P}^{\delta} \to u_s\nabla^2_{\v}\hat{P},\quad \mbox{strongly in }\,\, L^2_{loc}([0,L]\times\mathbb{R}_+)\,\, as\,\, \delta \to0,\\
		&u_s\nabla^2_{\v}\hat{\mathbf{S}}^{\delta} \to u_s\nabla^2_{\v}\hat{\mathbf{S}},\quad \mbox{strongly in }\,\, L^2_{loc}([0,L]\times\mathbb{R}_+)\,\, as\,\, \delta \to0,
	\end{split}
\end{align}
and
\begin{align}\label{17.232}
	\begin{split}
		& u_s\mathbf{w}^{\delta} \big|_{x=0,L} \to u_s\mathbf{w}\big|_{x=0,L},\quad  \mbox{strongly in} \,\, L^2_{loc}(\mathbb{R}_+)\,\, as\,\, \delta \to0,\\
		& \tilde{B}^{\delta} \big|_{x=0,L} \to \tilde{B}\big|_{x=0,L},\quad  \mbox{strongly in} \,\, L^2_{loc}(\mathbb{R})\,\, as\,\, \delta \to0,\\
		& \tilde{u}_s \nabla \tilde{B}^{\delta}\big|_{x=0,L} \to \tilde{u}_s \nabla \tilde{B} \big|_{x=0,L},\quad  \mbox{strongly in} \,\, L^2_{loc}(\mathbb{R})\,\, as\,\, \delta \to0,\\
		&u_s\nabla^2_{\v}P^{\delta}\big|_{x=0,L} \to u_s\nabla^2_{\v}P\big|_{x=0,L},\quad  \mbox{strongly in} \,\, L^2_{loc}(\mathbb{R})\,\, as\,\, \delta \to0,\\
		&u_s\nabla^2_{\v}\hat{P}^{\delta}\big|_{x=0,L} \to u_s\nabla^2_{\v}\hat{P}\big|_{x=0,L},\quad  \mbox{strongly in} \,\, L^2_{loc}(\mathbb{R})\,\, as\,\, \delta \to0,\\
		&u_s\nabla^2_{\v}\hat{\mathbf{S}}^{\delta}\big|_{x=0,L} \to u_s\nabla^2_{\v}\hat{\mathbf{S}}\big|_{x=0,L},\quad  \mbox{strongly in} \,\, L^2_{loc}(\mathbb{R})\,\, as\,\, \delta \to0.
	\end{split}
\end{align}
We remark that the weight $u_s$ in \eqref{17.232} is crucial for us to take the limit strongly.

\smallskip

{\it Step 2.} We consider the limit $\delta\to0$ in \eqref{17.119}.  For any given test function $\varphi \in H^1((0,L)\times\mathbb{R}_+)$ with compact support, a direct calculation shows that 
\begin{align}\label{17.233}
	\Big|\int_0^{L}\int_0^\infty \delta \nabla_{\v}\mathbf{w}^{\delta}\cdot \nabla_{\v}\varphi dydx\Big| \lesssim \delta \|\nabla_{\v}\mathbf{w}^{\delta}\|\cdot \|\nabla_{\v}\varphi\|\lesssim \sqrt{\delta} \to 0,\,\, \mbox{as}\,\, \delta \to0,
\end{align}
and
\begin{align}
	&\lim_{\delta\to0}\int_0^{L}\int_0^\infty \mathbf{w}^{\delta}\varphi dydx - (\mu+\lambda) \v \lim_{\delta\to0}\int_0^{L}\int_0^\infty \mathbf{w}^{\delta} \big[\varphi \mbox{div} \big(\f{1}{p^{\v,\d}}U^{\v,\d}\big) + \f{1}{p^{\v,\d}}U^{\v,\d}\cdot \nabla \varphi\big] dydx\nonumber\\
	&=\int_0^{L}\int_0^\infty \mathbf{w} \varphi dydx - (\mu+\lambda) \v\int_0^{L}\int_0^\infty \mathbf{w} \big[\varphi \mbox{div} \big(\f{1}{p^{\v}}U^{\v}\big) + \f{1}{p^{\v}}U^{\v}\cdot \nabla \varphi\big] dydx.
\end{align}
Similarly, it holds that 
\begin{align}
	\lim_{\delta\to0}\int_0^L\int_0^\infty  \hat{\mathcal{G}}^{\delta}\,\varphi dydx = \int_0^L\int_0^\infty  \hat{\mathcal{G}}\,\varphi dydx.
\end{align}

For the boundary term, we have from \eqref{17.232} that
\begin{align}
	\lim_{\delta\to0}\int_0^\infty \f{1}{p^{\v,\d}} u^{\v,\d} \mathbf{w}^{\delta} \varphi dy\Big|_{x=L} 
	&= \lim_{\delta\to0}\int_0^\infty \f{1}{p^{\v,\d}} u^{\v} \mathbf{w}^{\delta} \varphi dy\Big|_{x=L} + \lim_{\delta\to0}\int_0^\infty \f{\delta}{p^{\v,\d}}  \mathbf{w}^{\delta} \varphi dy\Big|_{x=L}\nonumber\\
	&=\int_0^\infty \f{1}{p^{\v}} u^{\v} \mathbf{w}^{\delta} \varphi dy\Big|_{x=L}.
\end{align}
Similarly,
\begin{align}\label{17.237}
	\begin{split}
		\lim_{\delta\to0}\int_0^\infty  \big[\f{u^{\v,\delta}}{p^{\v,\delta}}- \f{u_s^{\delta}}{p_s} \big] \mathbf{w}^{\delta} \varphi dy \Big|_{x=0}&= \int_0^\infty  \big[\f{u^{\v}}{p^{\v}}- \f{u_s}{p_s} \big] \mathbf{w} \varphi dy \Big|_{x=0},\\
		\lim_{\delta\to0}\int_0^\infty \mathbf{d}_{11}\tilde{B}^{\delta} \varphi  dy \Big|_{x=0} & =  \int_0^\infty \mathbf{d}_{11}\tilde{B} \varphi  dy \Big|_{x=0}.
	\end{split}	
\end{align}

Taking $\delta\to 0$ and using \eqref{17.233}-\eqref{17.237}, we have from \eqref{17.21} that 
\begin{align*}
	\mathcal{B}[\mathbf{w},\varphi]:&=\int_0^{L}\int_0^\infty \mathbf{w} \varphi dydx  - (\mu+\lambda) \v\int_0^{L}\int_0^\infty \mathbf{w} \big[\varphi \mbox{div} (\f{1}{p^{\v}}U^{\v}) + \f{1}{p^{\v}} U^{\v}\cdot \nabla \varphi\big] dydx \nonumber\\
	&\quad +(\mu+\lambda) \v \int_0^\infty \f{u_s }{p_s}  \mathbf{w} \varphi dy \Big|_{x=L} - (\mu+\lambda) \v \int_0^\infty  \big[\f{u^{\v}}{p^{\v}}- \f{u_s}{p_s} \big] \mathbf{w} \varphi dy \Big|_{x=0}\nonumber\\
	&=\int_0^{L}\int_0^\infty \hat{\mathcal{G}}\,\varphi dydx - \sigma \v \int_0^\infty \mathbf{d}_{11}\tilde{B} \varphi  dy \Big|_{x=0},
\end{align*}
which implies the weak formulation of $\mathbf{w}\equiv \Delta_{\v}P$, i.e.,
\begin{align}\label{17.238} 
	\begin{cases}
		\dis \mathbf{w}  + (\mu+\lambda) \v \f{1}{p^{\v}}(U^{\v}\cdot\nabla)\mathbf{w} = \hat{\mathcal{G}},\quad (x,y)\in (0,L)\times \mathbb{R}_+,\\
		\dis(\mu+\lambda) \f{1}{p_s} u_s\mathbf{w} \big|_{x=0}=-\sigma \mathbf{d}_{11} \tilde{B}.
	\end{cases}
\end{align}

{\it Step 3.} We take the limit $\delta\to0+$ in \eqref{17.120}. Since $\tilde{B}^{\d}$ is an odd function in $y$, so we only consider the test function  $\varphi\in C^2([0,L]\times \mathbb{R})$  which is also odd in $y$  and has compact support in $y$, thus it holds that $\varphi(x,0)\equiv 0$. 

For later use, we note 
\begin{align*}
	\begin{split}
		\big|\f{\varphi}{\tilde{u}_s^{\delta}}-\f{\varphi}{\tilde{u}_s}\big|^2\leq \f{\varphi^2}{\tilde{u}_s^2}\quad \mbox{and}\quad \big|\f{\varphi}{\tilde{u}_s^{\delta}}-\f{\varphi}{\tilde{u}_s}\big|^2\to0,\quad  \delta \to 0,\\
	\end{split}
\end{align*}
which, together with Lebesgue dominated convergence theorem, yields that 
\begin{align*}
	\lim_{\delta\to0}\iint |\f{\varphi}{\tilde{u}_s^{\delta}}-\f{\varphi}{\tilde{u}_s}\big|^2 dydx=0\quad \mbox{and}\quad \lim_{\delta\to0}\int_0^\infty |\f{\varphi}{\tilde{u}_s^{\delta}}-\f{\varphi}{\tilde{u}_s}\big|^2 dy\Big|_{x=0,L}=0.
\end{align*}

We need to take the limit $\delta\to0$ in the weak formulation \eqref{17.96}.
It is direct to see that 
\begin{align}\label{18.38}
	\begin{split}
		\lim_{\delta\to0}\iint \nabla_{\v}\tilde{B}^{\d}\cdot \nabla_{\v}\varphi dydx &= 2\int_0^L\int_0^\infty \nabla_{\v}\tilde{B}\cdot \nabla_{\v}\varphi dydx,\\
		\lim_{\delta\to0} \iint \tilde{\mathbf{w}}^{\d} \varphi_x dydx  &= 2\int_0^L\int_0^\infty \mathbf{w} \varphi_x dydx,\\
		\lim_{\delta\to0} \iint \bar{\mathbf{g}}_1^{\d} \varphi  dydx  &=  2\int_0^L\int_0^\infty \Delta_{\v}\bar{g}_1\, \varphi  dydx,\\
		\lim_{\delta\to0} \iint \tilde{\theta}^{\delta} \varphi dydx&=2\int_0^L\int_0^\infty \theta   \varphi dydx,
	\end{split}
\end{align}
and
\begin{align}\label{18.39}
	\begin{split}
		\lim_{\delta\to0}\int_{\mathbb{R}}  p_s \mathbf{d}_{11} \tilde{B}^{\d} \f{\varphi}{\tilde{u}^{\delta}_s} \f{1}{1+\bm{\chi}} dy \Big|_{x=0} & =\int_{\mathbb{R}}  p_s \mathbf{d}_{11} \tilde{B} \f{\varphi}{\tilde{u}_s}\f{1}{1+\bm{\chi}}  dy \Big|_{x=0},\\
		\lim_{\delta\to0}\int_{\mathbb{R}}  \tilde{u}_s\tilde{\theta}^{\delta}_x  \f{\varphi}{\tilde{u}^{\delta}_s}  dy \Big|_{x=0}
		&= \int_{\mathbb{R}} \tilde{\theta}_x  \tilde{\varphi}  dy \Big|_{x=0}.
	\end{split}
\end{align}

For terms involving $\tilde{\mathbf{g}}^{\d}_{b1}$, we have from \eqref{17.19-7}, \eqref{17.31}, \eqref{17.51} and \eqref{17.57} that 
\begin{align}\label{18.40}
	&\lim_{\delta\to0}\int_{\mathbb{R}}  \tilde{\mathbf{g}}^{\d}_{b1} \f{\varphi}{\tilde{u}^{\delta}_s} \f{1}{1+\bm{\chi}}dy \Big|_{x=0}
	\nonumber\\
	&=\lim_{\delta\to0}\int_{0}^\infty  \f{1}{1+\bm{\chi}} \bar{\chi}(\f{y}{\gamma}) \f{\varphi}{\tilde{u}^{\delta}_s} \bigg\{  - \f{1}{2}\lambda p_s\Big[\big(\frac{ u^{\v}}{p^{\v}}-\frac{ u_s}{p_s}\big) + \f{\bm{\chi}}{\rho_s} \big(\f{u^{\v}}{T^{\v}}-\f{u_s}{T_s}\big) -  \v^{N_0} \f{u^{\v}\rho^{(\v)}}{\rho_sp^{\v}}\Big]\hat{P}^{\d}_{yy}  \nonumber\\
	&\quad + (\mu+\lambda)p_s \big(\frac{u^{\v} }{p^{\v}}-\frac{u_s}{p_s}\big) \Delta_{\v}\hat{P}^{\d}  + \lambda\v  \frac{p_s}{p^{\v}} v^{\v} \hat{P}^{\d}_{xy}  -  \lambda  p_s  \frac{u^{\v}}{T^{\v}} \Big[\big(\f{\bm{\chi}}{2\rho_s}\big)_y \hat{P}^{\d}_y + \big(\f{\bm{\chi}}{2\rho_s} \big)_{yy} \hat{P}^{\d}\Big] \nonumber\\
	&\quad  - p_s G_{11}(\tilde{T}(\zeta^{\d},p^{\d},q))  - p_s G_{12}(p^{\d},q)   - \frac{1}{2}\v \lambda  [1+\bm{\chi}] u_s \bar{g}_{1x} + p_s\hat{d}_{11} \zeta^{\d}_x  + p_s\mathfrak{g}_1(q)  \nonumber\\
	&\quad    + \rho_s u_s^2 \mathbf{q}^{\d}_{xy} - p_s\mathfrak{N}^{\d}_{12}\bigg\}dy \Big|_{x=0}   + \lim_{\delta\to0}\int_{-\infty}^0 (\cdots) \, dy \Big|_{x=0} =:\int_{\mathbb{R}} \tilde{\mathbf{g}}_{b1} \f{\varphi}{\tilde{u}_s} \f{1}{1+\bm{\chi}} dy \Big|_{x=0}.
\end{align}

We notice from \eqref{17.232} that 
\begin{align}
	&\lim_{\delta\to0}\int_{\mathbb{R}}   \tilde{B}^{\d}_{y} \,  \big(\f{b_2}{2 p_s}\tilde{u}_s^{\delta} \varphi\big)_y  dy \Big|_{x=L} \nonumber\\
	&=\lim_{\delta\to0}\int_{\mathbb{R}} \tilde{u}_s \tilde{B}^{\d}_{y} \, \Big[\f{b_2}{2p_s}\varphi_y + \big(\f{b_2}{2p_s}\tilde{u}_s\big)_y \f{\varphi}{\tilde{u}_s}\Big]  dy \Big|_{x=L}  + \lim_{\delta\to0}\int_{\mathbb{R}} \delta  \tilde{B}^{\d}_{y} \,  \big(\f{b_2}{2p_s}  \varphi\big)_y  dy\Big|_{x=L} \nonumber\\
	&= \int_{\mathbb{R}} \tilde{u}_s \tilde{B}_{y} \, \Big[\f{b_2}{p_s}\varphi_y + \big(\f{b_2}{p_s}\tilde{u}_s\big)_y \f{\varphi}{\tilde{u}_s}\Big]  dy \Big|_{x=L}\nonumber\\
	&\equiv \int_{\mathbb{R}}   \tilde{B}_{y} \,  \big(\f{b_2}{p_s}\tilde{u}_s \varphi\big)_y  dy \Big|_{x=L},
\end{align}
where we have used the fact
\begin{align}
	\Big|\lim_{\delta\to0}\int_{\mathbb{R}} \delta  \tilde{B}^{\d}_{y} \,  \big(\f{b_2}{p_s}  \varphi\big)_y  dy\Big|_{x=L} \Big| \lesssim \sqrt{\delta} \|\sqrt{\tilde{u}^{\delta}_{s}}\tilde{B}_y^{\d}\|_{x=L}  \|\varphi_{y}\|_{\mathcal{H}^1} \lesssim \sqrt{\delta} \to 0\,\quad \mbox{as}\,\, \delta\to0.
\end{align}

\smallskip

For term involving $\tilde{\Xi}^{\d}$, due to $u_{s}|_{y=0}=\varphi|_{y=0}=0$, we have 
\begin{align}
	&\lambda\v \int_{0}^{\infty}   \f{b_2}{2p_s}\tilde{u}_s \tilde{\Xi}^{\delta}  \varphi dy \Big|_{x=L} \nonumber\\
	&= - \lambda\v \int_{0}^{\infty} \theta^{\delta}_{y} \big[\f{b_2}{2p_s} \tilde{u}_s\varphi\big]_y\bar{\chi}(\f{y}{\gamma}) dy \Big|_{x=L}  -  \lambda\v \int_{0}^{\infty} \f{b_2}{2p_s} \tilde{u}_s  \theta^{\delta}_{y}\, \varphi  \, \f{1}{\g}\bar{\chi}'(\f{y}{\gamma})dy \Big|_{x=L} \nonumber\\
	&\to -  \lambda\v \int_{0}^{\infty} \theta_{y} \big[\f{b_2}{2 p_s} u_s\varphi\big]_y dy \Big|_{x=L} \quad \mbox{as}\,\, \delta\to 0,
\end{align}
where we have used the fact
\begin{align}
	\int_{0}^{\infty} \f{b_2}{p_s} \tilde{u}_s  \theta^{\delta}_{y}\, \varphi  \, \f{1}{\g}\bar{\chi}'(\f{y}{\gamma})dy \Big|_{x=L}
	& \lesssim \sqrt{\g(\delta)} \|\sqrt{u_s}\theta^{\delta}_y\|_{x=L}\cdot \|\frac{\varphi}{y}\|_{x=L} \to 0\,\,\mbox{as}\,\, \delta\to0.
\end{align}

\smallskip

Next, we consider the boundary term in \eqref{17.96} involving $\tilde{\mathbf{g}}^{\d}_{b2}$. Recall \eqref{17.51},  \eqref{17.19-8} and \eqref{17.31}.	It follows from \eqref{17.232} that 
\begin{align}\label{17.252}
	&\int_{0}^{\infty} \bar{\chi}(\f{y}{\gamma})\Big[ \lambda\v \big(\frac{b_2}{2p_s} u_s\big)_y \hat{P}^{\d}_{xy}   + \v \bar{g}^{\d}_{1x}  + \big(b_2 [ g^{\d}_2+ \v \mathfrak{g}_2 - \v\mathfrak{h}^{\d}_2 - \v\mathfrak{N}^{\d}_{21}]\big)_y  \Big] \, \varphi dy \Big|_{x=L}\nonumber\\
	&\to \int_{0}^{\infty}\big[\lambda\v (\f{b_2}{2p_s} \tilde{u}_s)_y\hat{P}_{xy}  + \v \bar{g}_{1x} +  \big(b_2 [ g_2+ \v \mathfrak{g}_2 - \v\mathfrak{h}_2 - \v\mathfrak{N}_{21}]\big)_y\big] \, \varphi dy \Big|_{x=L}\quad \mbox{as}\,\, \delta\to 0,
\end{align}
and
\begin{align}
	\lambda \v  \int_{0}^{\infty} \big(\f{b_2}{2p_s} \tilde{u}_s\varphi\big)_y  \bar{\chi}(\f{y}{\gamma})\bar{g}^{\d}_{1y} dy \Big|_{x=L}
	\to \lambda \v  \int_{0}^{\infty} \big(\f{b_2}{2 p_s} \tilde{u}_s\varphi\big)_y \bar{g}_{1y} dy  \Big|_{x=L}\,\,\,\mbox{as}\,\,\delta\to 0.
\end{align}

\smallskip

Noting \eqref{17.19-4} and \eqref{17.51}, we have   
\begin{align}\label{17.254}
	&\int_{0}^{\infty} \bar{\chi}(\f{y}{\gamma}) \hat{G}^{\d}_{3}\, \varphi_y dy \Big|_{x=L}\nonumber\\
	&=\int_{0}^{\infty} \varphi_y \bar{\chi}(\f{y}{\gamma}) \bigg\{- b_2\frac{\mu\v}{p^{\v}} v^{\v} \Delta_{\v}\hat{P}^{\d}   - b_2\frac{\lambda\v}{p^{\v}} v^{\v} \hat{P}^{\d}_{yy} - \f12 \lambda \v b_2 \hat{P}^{\d}_{xy} \Big[\big(\frac{ u^{\v}}{p^{\v}}-\frac{ u_s}{p_s}\big) + \v^{N_0} \f{u^{\v}\rho^{(\v)}}{\rho_sp^{\v}}\Big] \nonumber \\
	&\,\,\,+ \lambda \v u_s^{\v} \hat{\mathbf{S}}^{\d}_{xy} \mathbf{1_{D}} + \v b_2 m u_s^2 \mathbf{q}^{\d}_{xx} + \v b_2 G_{21}(\tilde{T}(\zeta^{\d},p^{\d},q)) + \v b_2 G_{22}(p^{\d},q)   + \v b_2 [h^{\d}_2 + \mathfrak{N}^{\d}_{22}]\bigg\}\,  dy \Big|_{x=L}\nonumber\\
	&\to  \int_{0}^{\infty} G_{3}(\Phi, \hat{P}, \hat{\mathbf{S}}, p, \zeta, \phi) \, \varphi_y dy \Big|_{x=L} \,\,  \mbox{as}\,\, \delta\to0.
\end{align}

Combining \eqref{18.38}-\eqref{17.254}, we obtain the weak formulation
\begin{align} \label{17.258}
	&\int_0^L\int_0^\infty \big[-\Delta_{\v}\tilde{B} + \mathbf{w}_x + \Delta_{\v}\bar{g}_1 -\theta\big]\, \varphi dydx   + \int_0^\infty \big[\v \tilde{B}_x -  \mathbf{w}\big]\varphi dy \Big|_{x=L} \nonumber\\
	&= \int_0^\infty \Big\{\big[\v \tilde{B}_x- \mathbf{w}\big] + \f{2}{\lambda(1+\bm{\chi})} \f{1}{u_s}\big\{\f12 \lambda \v [1+\bm{\chi}] u_s \theta_x +  \mathbf{g}_{b1}\big\} -  (\sigma-1) p_s \mathbf{d}_{11}   \f{\tilde{B}}{u_s} \Big\}\varphi dy \Big|_{x=0} \nonumber\\
	&\quad  - \lambda \v\int_0^\infty\tilde{B}_{y} \,  \big(\f{b_2}{2p_s}u_s  \varphi\big)_y  dy \Big|_{x=L} - \lambda\v \int_{0}^{\infty} \theta_{y} \big(\f{b_2}{2p_s} u_s\varphi\big)_y  dy \Big|_{x=L} \nonumber\\
	&\quad  +\int_{0}^{\infty} \Big[\lambda\v (\f{b_2}{2p_s} u_s)_y \hat{P}_{xy}  + \v \bar{g}_{1x} + \big(b_2 [ g_2+ \v \mathfrak{g}_2 - \v\mathfrak{h}_2 - \v\mathfrak{N}_{21}]\big)_y\Big] \, \varphi dy \Big|_{x=L} \nonumber\\
	&\quad  + \lambda \v  \int_{0}^{\infty} \big(\f{b_2}{2p_s} u_s\varphi\big)_y\bar{g}_{1y} dy \Big|_{x=L} + \int_{0}^{\infty} G_{3} \, \varphi_y dy \Big|_{x=L}.
\end{align}

From \eqref{18.24}, it holds that (in the sense of distribution)
\begin{align}
	\Delta_{\v}\tilde{B} = \mathbf{w}_x + \Delta_{\v}\bar{g}_1 -\theta,\quad (x,y)\in (0,L)\times\R_+.
\end{align}
Since $\tilde{B}$ is odd in $y$, we have  
\begin{align}\label{17.263-1}
	\tilde{B}\big|_{y=0}=0 \quad  \Longrightarrow \quad B\big|_{y=0} = \f{\mu}{\rho^{\v}} \Phi_{yyy} \quad \mbox{with}\,\,\, B:=P_x + \bar{g}_1.
\end{align}
For the boundary condition at $x=0$, we have from \eqref{17.258} that 
\begin{align*}
	\big\{ (\sigma-1)  p_s \mathbf{d}_{11} \tilde{B} - \f12 \v \lambda(1+\bm{\chi}) u_s  \tilde{B}_x\big\} \big|_{x=0} 
	&= - \f12 \lambda(1+\bm{\chi}) u_s \mathbf{w} +\big\{\f12 \lambda \v [1+\bm{\chi}] u_s \theta_x +  \mathbf{g}_{b1}\big\}, 
\end{align*}
which, together with  $\tilde{B}=B -\theta$ and $\theta|_{x=0}=0$, yields that 
\begin{align}\label{17.261}
	\big\{ (\sigma-1)  p_s \mathbf{d}_{11} B - \f12 \v \lambda(1+\bm{\chi}) u_s  B_x\big\} \big|_{x=0} 
	&= - \f12 \lambda(1+\bm{\chi}) u_s \Delta_{\v}P +   \mathbf{g}_{b1}(\hat{P},   \Phi, p, S, \zeta, \phi).
\end{align}
For the boundary condition at $x=L$, noting \eqref{17.258}, one gets that (in the sense of distribution)
\begin{align*}
	\big\{\v \tilde{B}_{x} - \lambda\v\frac{b_2}{2p_s} u_s \tilde{B}_{yy}\big\} \big|_{x=L} &= \mathbf{w} + \lambda  \v \f{b_2}{2p_s} u_s  \theta_{yy} + \mathbf{g}_{b2}, 
\end{align*}
which, together with $\theta_x|_{x=L}=0$, yields 
\begin{align}\label{17.263}
	\big\{\v B_{x} - \lambda\v\frac{b_2}{2p_s} u_s B_{yy}\big\} \big|_{x=L} &= \Delta_{\v}P  + \mathbf{g}_{b2}(\hat{P}, \hat{\mathbf{S}}, \Phi, p, S, \zeta, \phi).
\end{align}

\smallskip

{\it Step 4.} Noting $\mathbf{w}=\Delta_{\v}P$, \eqref{17.238} and \eqref{17.263-1}-\eqref{17.263}, we prove \eqref{18.39-1}-\eqref{18.39-2}. Using \eqref{17.156-0}, Proposition \ref{Bprop7.11} and Lemma \ref{CBlem7.6}, we conclude \eqref{18.41}-\eqref{C17.223-1}. The uniqueness can be derived by using the contraction property like \eqref{17.114}.   Therefore  the proof of Theorem \ref{thm5.1} is completed. $\hfill\Box$

\medskip
 

\section{Existence and Uniform Estimate for Pseudo Entropy}\label{sec6}
In this section, we aim to establish the existence and uniform estimates for solutions of pseudo entropy BVPs \eqref{2.111},  \eqref{2.111-6D} and \eqref{2.111-1}.

\subsection{Existence and uniform estimate for Neumann BVP  \eqref{2.111}} \label{secT.1}

\begin{theorem}\label{thm6.7}
There exists a unique strong solution $\mathbf{S}$ to \eqref{2.111} with the following uniform estimates
\begin{align}
	&[[[\mathbf{S}]]]_{2,\hat{w}_0^{\v}} + 	[[[\mathbf{S}]]]_{b,\hat{w}^{\v}_0}
	\lesssim_{C_0}  \mathscr{K}(P) +  \f{1}{\v} \mathcal{K}_{2,1}(P)   +  \mathcal{K}_{3,1}(P)    +  \f{1}{\v} [[p]]_{2,\hat{w}_0^{\v}}  + [[p]]_{3,\hat{w}^{\v}_0}  \nonumber\\
	&\qquad\qquad\qquad\qquad\qquad  +  a^2 \|(\phi,p, S)\hat{w}_0^{\v}\|^2_{\mathbf{X}^{\v}}  + \f{1}{\v}\|\mathfrak{N}_3 \hat{w}_0^{\v}\|^2  + \|u_s \pa_y\mathfrak{N}_3 \hat{w}^{\v}_0\|^2,\\
	&[[[\mathbf{S}]]]_{3,\hat{w}^{\v}_0}
	\lesssim_{C_0}  \sqrt{\v} [[[\mathbf{S}]]]_{b,1} + L^{\f12} [[[\mathbf{S}]]]_{2,\hat{w}^{\v}_0}  +  L^2 \mathscr{K}(P)  + a^2 \f{1}{\v}\mathcal{K}_{2,1}(P)  +  \mathcal{K}_{3,1}(P)  \nonumber\\
	&\qquad\qquad\quad  + [[p]]_{3,\hat{w}^{\v}_0} +  L^{\f18} \|(\phi,p,S) \hat{w}^{\v}_0\|^2_{\mathbf{X}^{\v}} + \|u_s \pa_y\mathfrak{N}_3 \hat{w}^{\v}_0\|^2.
\end{align}
\end{theorem}

Let $P^{\d}$ be the one  constructed in Section \ref{ss7.6}. We consider the following approximate BVP
\begin{align}\label{2.111-6}
\begin{cases}
\dis	2 \rho_s u_s \mathbf{S}^{\d}_x + 2\rho_s v_s \mathbf{S}^{\d}_y - \kappa \Delta_{\v}\mathbf{S}^{\d} 
 = \kappa \Delta_{\v}\big(\frac{1}{2\rho_s}[\chi P^{\d} + \bar{\chi}p]\big) + \kappa \Delta_{\v}\big(\f{1}{p_s}T_{sy}q\big)  + \mathfrak{J},\\
\dis \mathbf{S}^{\d} \big|_{x=0}= -  \frac{\chi}{2\rho_s}P^{\d},\,\,\, \mathbf{S}^{\d}_{x}\big|_{x=L}=0 \,\, \,\mbox{and}\,\,\,
 \mathbf{S}^{\d}_y \big|_{y=0} = - \big(\frac{1}{2\rho_s}P^{\d}\big)_y.
\end{cases}
\end{align}
We define
\begin{align*}
	T^{\d} := \mathbf{S}^{\d}  + \frac{1}{2\rho_s}[\chi P^{\d} + \bar{\chi} p] + \f{1}{p_s}T_{sy}q,
\end{align*}
then it is clear to see that $T^{\d}$ satisfies a type of elliptic equation  \eqref{C.01}. Applying Proposition \ref{Bprop6.16} (see Section \ref{SecB} for more details), we  obtain the existence of solution $T^{\d}$ with $T^{\d}-\bar{\chi}p\in H^3 $. Noting $\chi P^{\delta} \in H^3, T_{sy}q\in H^3$, we conclude  $\mathbf{S}^{\d} \in H^3$. In the following of this subsection, we aim to establish some uniform estimates,  then take $\delta\to 0+$ to prove Theorem \ref{thm6.7}.

\begin{lemma}\label{lemTH.1-0}
	It holds that 
	\begin{align}\label{1TC.21}
		L^{-1}\|\mathbf{S}^{\d}_yw\|^2_{x=0,L} + L^{-2}\|\mathbf{S}^{\d}_yw\|^2&\lesssim \|\mathbf{S}^{\d}_{xy}w\|^2 +  \v^{\f{1-\b}{2}-} \big\{\mathscr{K}(P^{\d}) + \f{1}{\v} \mathfrak{K}_{2,1}(P^{\d}) \big\}.
	\end{align}
\end{lemma}

\noindent{\bf Proof.} Noting \eqref{2.111-6}, we have 
\begin{align*}
	\|\mathbf{S}^{\d}_yw\|^2_{x=0}&\lesssim \|\chi P^{\d}_y\|^2_{x=0} \lesssim \f{1}{L} \|\chi P^{\d}_y\|^2 +\|\chi P^{\d}_y\|\cdot \|\chi P^{\d}_{xy}\|\lesssim \v^{\f{1-\b}{2}} \big\{\mathscr{K}(P^{\d}) + \f{1}{\v} \mathfrak{K}_{2,1}(P^{\d}) \big\}.
\end{align*}
The cut-off function is crucial in above estimate  due to the lack of good enough control on $\|P^{\d}_y\|$.

It is clear that 
\begin{align*}
	\begin{split}
		\|\mathbf{S}^{\d}_yw\|^2 &\lesssim L \|\mathbf{S}^{\d}_yw\|^2_{x=0} + L^2\|\mathbf{S}^{\d}_{xy}w\|^2 
		\lesssim L^2\|\mathbf{S}^{\d}_{xy}w\|^2  + L \v^{\f{1-\b}{2}} \big\{\mathscr{K}(P^{\d}) + \f{1}{\v} \mathfrak{K}_{2,1}(P^{\d}) \big\},
	\end{split}
\end{align*}
and
\begin{align*}
	\|\mathbf{S}^{\d}_y w\|^2_{x=L}&\lesssim \|\mathbf{S}^{\d}_{xy}w\|\cdot \|\mathbf{S}^{\d}_yw\| +  \|\mathbf{S}^{\d}_y w\|^2_{x=0} \lesssim L\|\mathbf{S}^{\d}_{xy}w\|^2 + \|\mathbf{S}^{\d}_y w\|^2_{x=0} \nonumber\\
	&\lesssim L \|\mathbf{S}^{\d}_{xy}w\|^2  +  \v^{\f{1-\b}{2}} \big\{\mathscr{K}(P^{\d}) + \f{1}{\v} \mathfrak{K}_{2,1}(P^{\d}) \big\}.
\end{align*}
Therefore the proof of Lemma \ref{lemTH.1-0} is completed. $\hfill\Box$

\begin{lemma}\label{lemTH6.3}
	It holds that 
	\begin{align}\label{1TC.51-0}
		&\|(\sqrt{\v}\mathbf{S}^{\d}_{xx}, \mathbf{S}^{\d}_{xy})w\|^2 + \|\sqrt{u_s}\mathbf{S}^{\d}_xw\|^2_{x=0}\nonumber\\
		&\lesssim \sqrt{L} \|\mathbf{S}^{\d}_{xx}w_y\|^2 +  (\xi + L^{\f38-}) [[[\mathbf{S}^{\d}]]]_{2,\hat{w}}+ C_{\xi} \mathscr{K}(P^{\d})+ C_{\xi}\f{1}{\v} \mathfrak{K}_{2,1}(P^{\d}) + C_{\xi}\f{1}{\v} [[p]]_{2,\hat{w}}\nonumber\\
		&\quad  + C_{\xi} a^2\|(\phi, 0, S)\hat{w}\|^2_{\mathbf{X}^{\v}} + C_{\xi} L \f{1}{\v}\|(0,p, S)\|^2_{\mathbf{X}^{\v}} \mathbf{1}_{\{w=w_0\}} + C_{\xi}  \f{1}{\v}\|\mathfrak{N}_3 w\|^2.
	\end{align}
\end{lemma}

\noindent{\bf Proof.} Multiplying \eqref{T.1} by $\mathbf{S}^{\d}_{xx}w^2$, one obtains that
\begin{align}\label{1TC.12}
	&\kappa \|\sqrt{\v}\mathbf{S}^{\d}_{xx}w\|^2 + \big\la \kappa \mathbf{S}^{\d}_{yy},\, \mathbf{S}^{\d}_{xx} w^2\big\ra + \|\sqrt{\rho_s u_s}\mathbf{S}^{\d}_xw\|^2_{x=0} - \cancel{\|\sqrt{\rho_s u_s}\mathbf{S}^{\d}_xw\|^2_{x=L}}\nonumber\\
	& \leq C \|\sqrt{u_s}\mathbf{S}^{\d}_xw\|^2 + \big|\big\la 2\rho_sv_s \mathbf{S}^{\d}_y,\, \mathbf{S}^{\d}_{xx} w^2\big\ra \big|  + \Big|\Big\la \kappa \Delta_{\v}\big(\frac{1}{2\rho_s}\chi P^{\d}\big),\, \mathbf{S}^{\d}_{xx}  w^2 \Big\ra\Big|\nonumber\\
	&\quad  + \Big|\Big\la \kappa \Delta_{\v}\big(\frac{\bar{\chi}}{2\rho_s} p\big),\, \mathbf{S}^{\d}_{xx} w^2\Big\ra \Big| + \Big|\Big\la \kappa \Delta_{\v}\big(\f{1}{p_s}T_{sy}q\big),\, \mathbf{S}^{\d}_{xx} w^2 \Big\ra \Big| + \Big|\Big\la  \mathfrak{J},\, \mathbf{S}^{\d}_{xx} w^2 \Big\ra \Big|.
\end{align}
Since the proof is  complicate, we divide it into five steps.

\smallskip

{\it Step 1.} Integrating by parts, one has 
\begin{align}\label{1TC.13}
	\la \kappa \mathbf{S}^{\d}_{yy},\, \mathbf{S}^{\d}_{xx} w^2\ra 
	&=\kappa \|\mathbf{S}^{\d}_{xy}w\|^2  - \la \kappa \mathbf{S}^{\d}_{y},\, \mathbf{S}^{\d}_{xx} 2ww_y \ra + \kappa\la \mathbf{S}^{\d}_{y},\, \mathbf{S}^{\d}_{xy}w^2 \ra_{x=0} - \kappa \la \mathbf{S}^{\d}_{y},\, \mathbf{S}^{\d}_{xx} w^2\ra_{y=0},
\end{align}
where we have used $\mathbf{S}^{\d}_x|_{x=L}=0$.
Using \eqref{1TC.21}, we know that 
\begin{align}\label{1TC.14-0}
\eqref{1TC.13}_2\lesssim L\|\mathbf{S}^{\d}_{xx}w_y\|^2 + L \|\mathbf{S}^{\d}_{xy}w\|^2  + \v^{\f{1-\b}{2}} \big\{\mathscr{K}(P^{\d}) + \f{1}{\v} \mathfrak{K}_{2,1}(P^{\d}) \big\}.
\end{align}
Noting  $\eqref{2.111-6}_2$, one has that
\begin{align}\label{1TC.14}
\eqref{1TC.13}_3= - \Big\la \kappa \mathbf{S}^{\d}_{xy},\, (\f{\chi}{2\rho_s} P^{\delta})_y w^2\Big\ra_{x=0},
\end{align}
and
\begin{align}\label{1TC.15}
\eqref{1TC.13}_4	&=-\kappa \Big\la  \mathbf{S}^{\d}_{xxy},\, (\f{\chi}{2\rho_s}P^{\d})_y w^2\Big\ra  - \kappa \Big\la  \mathbf{S}^{\d}_{xx},\, (\f{\chi}{2\rho_s}P^{\d})_{yy} w^2\Big\ra - \kappa \Big\la  \mathbf{S}^{\d}_{xx},\, (\f{\chi}{2\rho_s}P^{\d})_y 2ww_y\Big\ra\nonumber\\
	&=\kappa \Big\la  \mathbf{S}^{\d}_{xy},\, (\f{\chi}{2\rho_s}P^{\d})_{xy} w^2\Big\ra  - \kappa \Big\la  \mathbf{S}^{\d}_{xx},\, (\f{\chi}{2\rho_s}P^{\d})_{yy}w^2\Big\ra - \kappa \Big\la  \mathbf{S}^{\d}_{xx},\, (\f{\chi}{2\rho_s}P^{\d})_y 2ww_y\Big\ra\nonumber\\
	&\quad   +  \Big\la  \kappa \mathbf{S}^{\d}_{xy},\, (\f{\chi}{2\rho_s}P^{\d})_y w^2\Big\ra_{x=0}  .
\end{align}

Next we shall control terms on RHS of \eqref{1TC.15}.
Noting \eqref{6.174}, one obtains that 
\begin{align}\label{1TC.29-1}
 \eqref{1TC.15}_3&\lesssim \|\mathbf{S}^{\d}_{xx}w_y\|\cdot \|(P^{\d}_y\chi, P^{\d}\chi)\|  
 \lesssim L \|\mathbf{S}^{\d}_{xx}w_y\|^2 + L\mathscr{K}(P^{\d}) + L\f{1}{\v} \mathfrak{K}_{2,1}(P^{\d}),\\
\eqref{1TC.15}_1	&\lesssim \|\mathbf{S}^{\d}_{xy}\|\cdot \|(P^{\d}_{xy}, P^{\d}_x, P^{\d}_y, P^{\d})\chi\|  
\lesssim \xi \|\mathbf{S}^{\d}_{xy}\|^2 + C_{\xi} \mathscr{K}(P^{\d}) + C_{\xi} \f{1}{\v} \mathfrak{K}_{2,1}(P^{\d}).\label{1TC.29}
\end{align}

For $\eqref{1TC.15}_2$, using  , we get that 
\begin{align}\label{1TC.30}
	\eqref{1TC.15}_2	&\lesssim \big|\big\la  \mathbf{S}^{\d}_{xx},\, \f{\chi}{2\rho_s} P^{\d}_{yy} w^2 \big\ra\big| + \big|\big\la  \mathbf{S}^{\d}_{xx},\, (\f{\chi}{\rho_s})_y P^{\d}_{y} w^2 \big\ra\big| + \big|\big\la  \mathbf{S}^{\d}_{xx},\, (\f{\chi}{2\rho_s})_{yy} P^{\d} w^2 \big\ra\big| \nonumber\\
	&\leq \xi \|(\sqrt{\v}\mathbf{S}^{\d}_{xx}, \mathbf{S}^{\d}_{xy})\|^2 + \xi \|\sqrt{u_s}\mathbf{S}^{\d}_x\|^2_{x=0}  + \f{1}{a^2} \|u_s\mathbf{S}^{\d}_x\|^2 + \|\sqrt{\v}\mathbf{S}^{\d}_x\|^2_{x=0}  \nonumber\\
	&\quad  + C_{\xi} \mathscr{K}(P^{\d})+ C_{\xi}\f{1}{\v} \mathfrak{K}_{2,1}(P^{\d}).
\end{align}
where we have used $P^{\d}(0,0)=0$, \eqref{2TC.3-1} and integration by parts to derive the following two estimates
\begin{align}\label{1TC.12-2}
	\Big|\Big\la \mathbf{S}^{\d}_{xx},\,  (\f{\chi}{\rho_s})_{yy} P^{\d} w^2\Big\ra\Big|
	&\leq  \big|\big\la \mathbf{S}^{\d}_{x},\,  \big\{(\f{\chi}{\rho_s})_{yy} P^{\d}_x + (\f{\chi}{\rho_s})_{xyy} P^{\d}\big\} w^2\big\ra\big|  +  \big|\big\la \mathbf{S}^{\d}_{x},\,  (\f{\chi}{\rho_s})_{yy} P^{\d} w^2\big\ra_{x=0}\big| \nonumber\\
	&\lesssim \|\mathbf{S}^{\d}_x\| \big(\|P^{\d}_x\chi\| + \|P^{\d}\chi \|\big) + \f{1}{a} \|P^{\d}_y\chi\|_{x=0}\cdot \|\sqrt{u_s}\mathbf{S}^{\d}_{x}\|_{x=0}\nonumber\\
	&\leq \xi \|\mathbf{S}^{\d}_{xy}\|^2 + \xi \|\sqrt{u_s}\mathbf{S}^{\d}_x\|^2_{x=0} +  \f{1}{a^2} \|u_s\mathbf{S}^{\d}_x\|^2 + C_{\xi}  \f{1}{\v} \mathfrak{K}_{2,1}(P^{\d})  + C_{\xi} \mathscr{K}(P^{\d}) ,
\end{align}
and
\begin{align}\label{1TC.12-3}
	\Big|\Big\la  \mathbf{S}^{\d}_{xx},\, (\f{\chi}{\rho_s})_y P^{\d}_{y} w^2\Big\ra\Big| 
	&\lesssim \Big|\Big\la  \mathbf{S}^{\d}_{x},\, (\f{\chi}{\rho_s})_y P^{\d}_{xy}w^2\Big\ra\Big| + \Big|\Big\la  \mathbf{S}^{\d}_{x},\, (\f{\chi}{\rho_s})_{xy} P^{\d}_{y}w^2\Big\ra\Big| + \Big|\Big\la \mathbf{S}^{\d}_{x},\, (\f{\chi}{\rho_s})_y P^{\d}_{y}w^2\Big\ra_{x=0}\Big| \nonumber\\
	&\lesssim \xi \|\mathbf{S}^{\d}_{xy}\|^2 + \xi \|\sqrt{u_s}\mathbf{S}^{\d}_x\|^2_{x=0} + \|\sqrt{\v}\mathbf{S}^{\d}_x\|^2_{x=0}  +  \f{1}{a^2} \|u_s\mathbf{S}^{\d}_x\|^2 + C_{\xi}\f{1}{\v} \mathfrak{K}_{2,1}(P^{\d})\nonumber\\
	&\quad  + C_{\xi} \v^{\f{1-\beta}{2}-} \mathscr{K}(P^{\d}) ,
\end{align}
We also remark that  $\rho_{sy}=\f{1}{T_s}\big\{p_{sy} -\rho_s T_{sy} \big\}=O(1)T_{sy} + O(1) \sqrt{\v}$ and $T_{sy}|_{y=0}=0$ are used  in \eqref{1TC.12-3}.   

Combining \eqref{1TC.13}-\eqref{1TC.30}, and noting the {\it cancellation} of boundary terms in \eqref{1TC.14}-\eqref{1TC.15}, one obtains that
\begin{align}\label{1TC.31}
	\la \kappa \mathbf{S}^{\d}_{yy},\, \mathbf{S}^{\d}_{xx} w^2\ra 
	&\geq \f34\kappa \|\mathbf{S}^{\d}_{xy}w\|^2  -  \xi \|(\sqrt{\v}\mathbf{S}^{\d}_{xx}, \mathbf{S}^{\d}_{xy})\|^2 - C\xi \|\sqrt{u_s}\mathbf{S}^{\d}_x\|^2_{x=0}   -  L\|\mathbf{S}^{\d}_{xx}w_y\|^2  \nonumber\\
	&\quad  - L^{\f12-} [[[\mathbf{S}^{\d}]]]_{2,\hat{w}} - C_{\xi}\f{1}{\v} \mathfrak{K}_{2,1}(P^{\d})  - C_{\xi} \mathscr{K}(P^{\d}),
\end{align}
where we have used $\eqref{2TC.3}_2$.
We remark that, by noting $\mathbf{S}^{\d} \in H^3, \chi P^{\d} \in H^3$, the term $\dis \big\la  \mathbf{S}^{\d}_{xy},\, (\f{\chi}{\rho_s}P^{\d})_y \big\ra_{x=0}$ make sense, i.e., it is finite. This is the main reason that we use $P^{\d}$ in \eqref{2.111-6}.

\smallskip

{\it Step 2.} For 2th term on RHS of \eqref{1TC.12}, one gets
\begin{align}\label{1TC.33}
\eqref{1TC.12}_2	&\leq  \big|\big\la 2\rho_sv_s \mathbf{S}^{\d}_{xy},\, \mathbf{S}^{\d}_{x} w^2\big\ra \big| + \big|\big\la 2(\rho_sv_s)_x \mathbf{S}^{\d}_y,\, \mathbf{S}^{\d}_{x} w^2 \big\ra \big| +  \big|\big\la 2\rho_sv_s \mathbf{S}^{\d}_y,\, \mathbf{S}^{\d}_{x} w^2\big\ra_{x=0} \big| \nonumber\\
&\leq \f18\kappa\|\mathbf{S}^{\d}_{xy}w\|^2 + \f18\|\sqrt{\rho_su_s}\mathbf{S}^{\d}_xw\|^2_{x=0} + L^{\f12-} [[[\mathbf{S}^{\d}]]]_{2,w} + L\big\{\mathscr{K}(P^{\d}) + \f{1}{\v} \mathfrak{K}(P^{\d}) \big\},
\end{align}
where we have used the fact \eqref{1TC.21}.

\smallskip

For 3th term on RHS of \eqref{1TC.12}, noting \eqref{1TC.30}, we have 
\begin{align}\label{1TC.30-0}
\eqref{1TC.12}_3	&\lesssim \xi \|(\sqrt{\v}\mathbf{S}^{\d}_{xx}, \mathbf{S}^{\d}_{xy})\|^2 + \xi \|\sqrt{u_s}\mathbf{S}^{\d}_x\|^2_{x=0}  + L^{\f12-} [[[\mathbf{S}^{\d}]]]_{2,1}  + C_{\xi} \mathscr{K}(P^{\d})+ C_{\xi}\f{1}{\v} \mathfrak{K}_{2,1}(P^{\d}).
\end{align}
For 4th term on RHS of \eqref{1TC.12}, it holds that 
\begin{align}\label{1TC.12-1}
\eqref{1TC.12}_4& = \Big\la \v (\f{\bar{\chi}}{\rho_s} p)_{xx} ,\, \mathbf{S}^{\d}_{xx} w^2 \Big\ra + \Big\la \f{\bar{\chi}}{\rho_s} p_{yy} +  2(\f{\bar{\chi}}{\rho_s})_y p_{y} + (\f{\bar{\chi}}{\rho_s})_{yy} p,\, \mathbf{S}^{\d}_{xx} w^2\Big\ra \nonumber\\
	&\lesssim \xi \|(\mathbf{S}^{\d}_{xy},\sqrt{\v}\mathbf{S}^{\d}_{xx}) w\|^2 + \xi \|\mathbf{S}^{\d}_{xy}\|^2  + \xi \|\sqrt{u_s}\mathbf{S}^{\d}_xw\|^2_{x=0}  + L^{\f12-} [[[\mathbf{S}^{\d}]]]_{2,\hat{w}}  + C_{\xi}\f{1}{\v} [[p]]_{2,\hat{w}},
\end{align}
where we have used the fact
\begin{align}\label{1TC.12-5}
 \big|\big\la \mathbf{S}^{\d}_{xx} ,\, (\f{\bar{\chi}}{\rho_s})_{yy} p\, w^2\big\ra\big| 
	&\leq  \big|\big\la \mathbf{S}^{\d}_{x},\,  \big\{(\f{\bar{\chi}}{\rho_s})_{yy} p_x +  (\f{\bar{\chi}}{\rho_s})_{xyy} p\big\} \,  w^2\big\ra\big|  +  \big|\big\la \mathbf{S}^{\d}_{x},\,  (\f{\bar{\chi}}{\rho_s})_{yy} p\, w^2 \big\ra_{x=0}\big| \nonumber\\
	&\lesssim a^{-1}\|u_s \mathbf{S}^{\d}_x w\| \big(\|p_x\| + \|\mathbf{d}_3p w\| + \|p\chi\|\big) + a^{-\f12} \|p_y\|_{x=0}\cdot \|\sqrt{u_s}\mathbf{S}^{\d}_{x}w\|_{x=0}\nonumber\\
	&\leq \xi \|\sqrt{u_s}\mathbf{S}^{\d}_xw\|^2_{x=0}  + L^{\f12-} [[[\mathbf{S}^{\d}]]]_{2,w} + C_{\xi} \f{1}{\v} [[p]]_{2,1}.
\end{align}

\smallskip

{\it Step 3.}  We consider 5th term on RHS of \eqref{1TC.12}. A direct calculation shows that 
\begin{align}\label{1TC.36}
	\Big\la \v \big(\f{1}{p_s}T_{sy}q\big)_{xx}, \, \mathbf{S}^{\d}_{xx}\Big\ra \lesssim a\|\sqrt{\v}\mathbf{S}^{\d}_{xx}w\|  \|(\sqrt{\v}q_{xx}, \sqrt{\v}q_x)w\| 
	\lesssim \xi \|\sqrt{\v}\mathbf{S}^{\d}_{xx}w\|^2 + C_{\xi} L^{\alpha} [\phi]_{3,w}.
\end{align}
Integrating by parts in $x$, and noting $q|_{x=0}=0$, $\mathbf{S}^{\d}_{x}|_{x=L}=0$ and $T_{sy}|_{y=0}=0$,   one obtains that 
\begin{align}\label{1TC.37}
	\Big\la \big(\f{T_{sy}}{p_s}q\big)_{yy}, \, \mathbf{S}^{\d}_{xx} w^2\Big\ra 
	&=- \Big\la \f{T_{sy}}{p_s}q_{xyy}, \, \mathbf{S}^{\d}_{x} w^2\Big\ra - \Big\la \big(\f{T_{sy}}{p_s}\big)_x q_{yy}, \, \mathbf{S}^{\d}_{x}\Big\ra -  \Big\la 2\big(\big(\f{T_{sy}}{p_s}\big)_{y}q_{y}\big)_x, \, \mathbf{S}^{\d}_{x}\Big\ra \nonumber\\
	&\quad -  \Big\la \big(\big(\f{T_{sy}}{p_s}\big)_{yy} q\big)_x, \, \mathbf{S}^{\d}_{x}w^2\Big\ra \nonumber\\
	&= \Big\la \f{T_{sy}}{p_s}q_{xy}, \, \mathbf{S}^{\d}_{xy} w^2\Big\ra + \Big\la \big(\f{T_{sy}}{p_s}w^2\big)_y q_{xy}, \, \mathbf{S}^{\d}_{x} \Big\ra - \Big\la (\f{T_{sy}}{p_s})_x q_{yy}, \, \mathbf{S}^{\d}_{x}w^2\Big\ra\nonumber\\
	&\quad  -  \Big\la 2\big(\f{T_{sy}}{p_s}\big)_{y}q_{y}\big)_x, \, \mathbf{S}^{\d}_{x} w^2\Big\ra  -  \Big\la \big(\big(\f{T_{sy}}{p_s}\big)_{yy} q\big)_x, \, \mathbf{S}^{\d}_{x}w^2\Big\ra \nonumber\\
	&\lesssim \|\mathbf{S}^{\d}_{xy}w\|\cdot \|q_{xy}\| + \|\mathbf{S}^{\d}_{x}w\| \cdot \|(q_{xy}, q_{yy})\| \nonumber\\
	&\lesssim \xi \|\mathbf{S}^{\d}_{xy}\hat{w}\|^2   + L^{\f12-} [[[\mathbf{S}^{\d}]]]_{2,\hat{w}} + C_{\xi} L^{\alpha} [\phi]_{3,1},
\end{align}
which, together with \eqref{1TC.36},  yields immediately that
\begin{align}\label{1TC.38}
\eqref{1TC.12}_5	&\lesssim \xi \|(\sqrt{\v}\mathbf{S}^{\d}_{xx}, \mathbf{S}^{\d}_{xy})w\|^2  +  \xi \|\mathbf{S}^{\d}_{xy}\|^2  + L^{\f12-} [[[\mathbf{S}^{\d}]]]_{2,\hat{w}} + C_{\xi} L^{\alpha} [\phi]_{3,\hat{w}}.
\end{align}

\smallskip

{\it Step 4.} We consider 6th term on RHS of \eqref{1TC.12}, that is the term involving $\mathfrak{J}(\phi,p,S, T(p, q, S))$, whose precise definition is given in \eqref{2.84}. Integrating by parts in $x$, then using Lemmas \ref{lemA.3} \& \ref{lemA.5-1}, one has
\begin{align}\label{1TC.39}
	\Big\la  2\rho_sv_s \big(\f{1}{p_s}T_{sy}q\big)_y -  2\big[\rho_su_s \big(\f{1}{p_s}T_{sy}\big)_x-  T_{sy} \bar{u}_{sx}\big]q,\, \mathbf{S}^{\d}_{xx} w^2\Big\ra \lesssim \|u_s\mathbf{S}^{\d}_{x}w\|^2 + L^{\alpha} [\phi]_{3,w},
\end{align}
and
\begin{align}\label{1TC.40}
	&\Big\la \rho_s (U_s\cdot \nabla) (\frac{1}{\rho_s})\, p, \, \mathbf{S}^{\d}_{xx} w^2\Big\ra\nonumber\\
	& = - \Big\la \rho_s (U_s\cdot \nabla) (\frac{1}{\rho_s})\, p_x + \big(\rho_s (U_s\cdot \nabla) (\frac{1}{\rho_s})\big)_x\, p, \, \mathbf{S}^{\d}_{x} w^2\Big\ra- \Big\la \rho_s (U_s\cdot \nabla) (\frac{1}{\rho_s})\, p, \, \mathbf{S}^{\d}_{x} w^2 \Big\ra_{x=0}\nonumber\\
	&\lesssim  \|u_s\mathbf{S}^{\d}_{x}w\|\cdot \|\mathbf{d}_2(p_x, p)w\| + \|\sqrt{u_s}\mathbf{S}^{\d}_xw\|_{x=0} \cdot \|\mathbf{d}_2pw\|_{x=0} \nonumber\\
	&\lesssim \sqrt{\v} \|\sqrt{u_s}\mathbf{S}^{\d}_{x}w\|^2_{x=0} +  \|u_s\mathbf{S}^{\d}_{x}w\|^2 +  \f{1}{\v} [[p]]_{2,1} + L^2\f{1}{\v^2} [[p]]_{2,1}\cdot\mathbf{1}_{\{w=w_0\}}.
\end{align}

Noting $p_{sy}\cong \sqrt{\v}$, we integrate by parts in $x$ to obtain
\begin{align}\label{1TC.41}
 \Big\la \f{p_{sy}}{\rho_s}  \phi_x , \mathbf{S}^{\d}_{xx} w^2 \Big\ra
&	= \Big\la \f{p_{sy}}{\rho_s}\bar{u}_s q_x , \mathbf{S}^{\d}_{xx} w^2 \Big\ra + \Big\la \f{p_{sy}}{\rho_s} \bar{u}_{sx} q , \mathbf{S}^{\d}_{xx} w^2 \Big\ra \nonumber\\
	&=-\Big\la \f{p_{sy}}{\rho_s} \bar{u}_s q_{xx} , \mathbf{S}^{\d}_{x} w^2\Big\ra - \Big\la (\f{p_{sy}}{\rho_s}\bar{u}_s)_x q_x , \mathbf{S}^{\d}_{x} w^2\Big\ra - \Big\la \f{p_{sy}}{\rho_s}\bar{u}_s q_x , \mathbf{S}^{\d}_{x} w^2\Big\ra_{x=0} \nonumber\\
	&\quad - \Big\la \f{p_{sy}}{\rho_s} \bar{u}_{sx} q_x , \mathbf{S}^{\d}_{x} w^2\Big\ra - \Big\la (\f{p_{sy}}{\rho_s}\bar{u}_{sx})_x q , \mathbf{S}^{\d}_{x} w^2\Big\ra\nonumber \\
	&\lesssim a \|u_s \mathbf{S}^{\d}_{x}w\|\cdot \|(\sqrt{\v}q_{xx}, \sqrt{\v}q_{x})w\| + a \|\sqrt{u_s}\mathbf{S}^{\d}_{x}w\|_{x=0} \cdot \|\sqrt{\v}q_xw\|_{x=0} \nonumber\\
	&\lesssim \xi \|\sqrt{u_s}\mathbf{S}^{\d}_{x}w\|^2_{x=0} + \|u_s\mathbf{S}^{\d}_xw\|^2 + C_{\xi} L^{\alpha} [\phi]_{3,w}.
\end{align}
It is  clear to see that 
\begin{align}\label{1TC.42}
\big\la \mathfrak{M} +\mathfrak{N}_3, \, \mathbf{S}^{\d}_{xx}w^2\big\ra 
& = \big\la \mathfrak{M}_0, \, \mathbf{S}^{\d}_{xx} w^2 \big\ra + \big\la \v (\mathfrak{M}_1+\mathfrak{M}_2) , \, \mathbf{S}^{\d}_{xx} w^2\big\ra + \big\la \mathfrak{N}_3, \, \mathbf{S}^{\d}_{xx} w^2 \big\ra\nonumber\\
&= - \big\la \pa_x\mathfrak{M}_0, \, \mathbf{S}^{\d}_{x}w^2\big\ra  - \big\la \mathfrak{M}_0, \, \mathbf{S}^{\d}_{x}w^2\big\ra_{x=0} + \xi \|\sqrt{\v}\mathbf{S}^{\d}_{xx}w\|^2\nonumber\\
&\quad   + C_{\xi} \v\|(\mathfrak{M}_1,\mathfrak{M}_2)w\|^2 + C_{\xi}  \f{1}{\v}\|\mathfrak{N}_3 w\|^2.
\end{align}
It follows from  \eqref{1C.43-1}-\eqref{1TC.44} that 
\begin{align}\label{1TC.45}
	|\big\la \pa_x\mathfrak{M}_0, \, \mathbf{S}^{\d}_{x}w^2\big\ra| 
	&\lesssim  \xi \|\mathbf{S}^{\d}_{xy}\|^2 + \f{1}{a^2}\|u_s\mathbf{S}^{\d}_xw\|^2+ C_{\xi} \Big\{ a^2 \f{1}{\v} [[p]]_{2,1} + a^2 [[[S]]]_{2,1}\Big\}\nonumber\\
	&\quad + C_{\xi} La\Big\{\f{1}{\v^2} [[p]]_{2,1} + \f{1}{\v} [[[S]]]_{2,1}\Big\} \mathbf{1}_{\{w=w_0\}}  + C_{\xi} a^2 [\phi]_{3,w},
\end{align}
and
\begin{align}\label{1TC.46}
	\big|\big\la \mathfrak{M}_0, \, \mathbf{S}^{\d}_{x} w^2 \big\ra_{x=0}\big| 
	&\lesssim \|\sqrt{u_s}\mathbf{S}^{\d}_xw\|_{x=0} \|p_y\|_{x=0} + a\int_0^\infty \mathbf{d}_2 |\mathbf{S}^{\d}_x p|w^2 dy \big|_{x=0} \nonumber\\
	& \lesssim \xi \|\sqrt{u_s}\mathbf{S}^{\d}_xw\|^2_{x=0} + C_{\xi} \f{1}{\sqrt{\v}} [[p]]_{2,1} + C_{\xi} \f{1}{a\v^{\f32}} [[p]]_{2,1} \cdot \mathbf{1}_{\{w=w_0\}},
\end{align} 
where we have used the compatibility condition $\dis S|_{x=0}=-\f{\chi}{2\rho_s}p$ and  $p(0,0)=0$.

Noting $\mathfrak{J}$ in \eqref{2.84}, then  combining \eqref{1TC.39}-\eqref{1TC.46} and using \eqref{1TC.48}, one obtains  that 
\begin{align} \label{1TC.49}
	\eqref{1TC.12}_6&\lesssim  \xi \|(\sqrt{\v}\mathbf{S}^{\d}_{xx}, \mathbf{S}^{\d}_{xy})\hat{w}\|^2  + \xi \|\sqrt{u_s}\mathbf{S}^{\d}_xw\|^2_{x=0}  + L^{\f12-} [[[\mathbf{S}^{\d}]]]_{2,\hat{w}} + C_{\xi} a^2 [\phi]_{3,\hat{w}}  + C_{\xi} a^2 [[[S]]]_{2,1} \nonumber\\
	&\quad  + C_{\xi}\f{1}{\v} [[p]]_{2,1} + C_{\xi} L\Big\{\f{1}{\v^2} [[p]]_{2,1} + \f{1}{\v} [[[S]]]_{2,1}\Big\} \mathbf{1}_{\{w=w_0\}}    + C_{\xi}  \f{1}{\v}\|\mathfrak{N}_3 w\|^2.
\end{align}

{\it Step 5.} Substituting \eqref{1TC.31}-\eqref{1TC.12-1}, \eqref{1TC.38} and \eqref{1TC.49} into \eqref{1TC.12}, 
one obtains \eqref{1TC.51-0}.
Therefore the proof of Lemma \ref{lemTH6.3} is completed. $\hfill\Box$

\begin{lemma}\label{lemTH6.6}
	It holds that 
	\begin{align} \label{1TC.55-1}
		L^{-\f12} \big\{\|u_s \mathbf{S}^{\d}_{x}w\|^2 + \|\mathbf{S}^{\d}_{yy}w\|^2 \big\}
		&\lesssim  L^{\f12} [[[\mathbf{S}^{\d}]]]_{2,\hat{w}}  +  L \mathscr{K}(P^{\d}) + L  \f{1}{\v}\mathfrak{K}_{2,1}(P^{\d}) + L^{\f12-} \|(\phi,p,S)\hat{w}\|^2_{\mathbf{X}^{\v}}\nonumber\\
		&\quad  + L^{\f12}\f{1}{\v}\|(0,p,S)\|^2_{\mathbf{X}^{\v}} \mathbf{1}_{\{w=w_0\}}  +  \f{1}{\v^{0+}}\|\mathfrak{N}_3 w\|^2.
	\end{align}
\end{lemma}

\noindent{\bf Proof.} It is clear to have 
\begin{align}\label{1TC.10-1}
	\begin{split}
		\|\Delta_{\v}(\frac{\chi}{\rho_s}p)w\|^2 &\lesssim \|\Delta_{\v}p \chi \|^2 + \|\nabla_{\v}p \chi \|^2 + \|\chi p\|^2 \lesssim L^2  \f{1}{\v} [[p]]_{2,1},\\
		\|\Delta_{\v}(\frac{1}{\rho_s}p)w\|^2  
		&\lesssim \|\Delta_{\v}pw\|^2 + \|\mathbf{d}_2\nabla_{\v}pw\|^2 + a\|\mathbf{d}_3pw\|^2 \lesssim [[p]]_{2,w}+ L^2 \f{1}{\v} [[p]]_{2,1},\\
		\|\Delta_{\v}(\frac{\chi}{\rho_s}P^{\d})w\|^2 &\lesssim \|\Delta_{\v}P^{\d} \chi \|^2 + \|\nabla_{\v}P^{\d} \chi \|^2 + \|\chi P^{\d}\|^2 \lesssim L^2 \mathscr{K}(P^{\d}) + L^2 \f{1}{\v}\mathfrak{K}_{2,1}(P^{\d}),
	\end{split}
\end{align}
and 
\begin{align}\label{1TC.18}
	\|\Delta_{\v}\big(\f{1}{p_s}T_{sy}q\big) w\|^2 \lesssim \|T_{sy}\Delta_{\v}q w\|^2 + a\|\mathbf{d}_3\nabla_{\v}q w\|^2 + a \|\mathbf{d}_4 q w\|^2 \lesssim L^{1-} [\phi]_{3,1},
\end{align}
Then it follows from \eqref{2.111-6} and \eqref{1TC.10-1}-\eqref{1TC.18} that
\begin{align}\label{1TC.58-1}
 \|\big(2\rho_su_s \mathbf{S}^{\d}_x - \kappa \Delta_{\v}\mathbf{S}^{\d}\big)w\|^2  
	&\lesssim  L^2  \f{1}{\v} [[p]]_{2,\hat{w}} +  L^2 \mathscr{K}(P^{\d}) + L^2 \f{1}{\v}\mathfrak{K}_{2,1}(P^{\d})  +  L^{1-} [\phi]_{3,1}\nonumber\\ 
	&\quad +  \|\mathfrak{J}w\|^2+  \|v_s \mathbf{S}^{\d}_yw\|^2.
\end{align}

It is clear that 
\begin{align}\label{1TC.59}
	\|\big(2\rho_su_s \mathbf{S}^{\d}_x - \kappa \Delta_{\v}\mathbf{S}^{\d}\big)w\|^2
	&=4\|\rho_s u_s \mathbf{S}^{\d}_{x}w\|^2 + \kappa^2 \|(\v \mathbf{S}^{\d}_{xx}, \mathbf{S}^{\d}_{yy})w\|^2 + 2\kappa^2 \big\la \v \mathbf{S}^{\d}_{xx},\, \mathbf{S}^{\d}_{yy}w^2\big\ra \nonumber\\
	&\quad  - 4 \kappa \big\la \rho_s u_s \mathbf{S}^{\d}_{x},\, \v \mathbf{S}^{\d}_{xx} w^2\big\ra   - 4\kappa  \big\la \rho_s u_s \mathbf{S}^{\d}_{x},\, \mathbf{S}^{\d}_{yy}w^2 \big\ra.
\end{align}
A direct calculation shows that 
\begin{align}\label{1TC.60}
2\kappa^2 \big|\la \v \mathbf{S}^{\d}_{xx},\, \mathbf{S}^{\d}_{yy}w^2\big\ra\big| +  4 \kappa \big|\big\la \rho_s u_s \mathbf{S}^{\d}_{x},\, \v \mathbf{S}^{\d}_{xx} w^2\big\ra \big|
&\lesssim \sqrt{\v} [[[\mathbf{S}^{\d}]]]_{2,w},
\end{align}
and
\begin{align}\label{1TC.61}
	- 4 \big\la \rho_s u_s \mathbf{S}^{\d}_{x},\, \mathbf{S}^{\d}_{yy} w^2\big\ra 
	&=  4 \big\la \rho_s u_s \mathbf{S}^{\d}_{xy},\, \mathbf{S}^{\d}_{y} w^2\big\ra + 4 \big\la (\rho_s u_s)_y \mathbf{S}^{\d}_{x},\, \mathbf{S}^{\d}_{y} w^2 \big\ra + 8 \big\la \rho_s u_s \mathbf{S}^{\d}_{x},\, \mathbf{S}^{\d}_{y} w w_y\big\ra\nonumber \\
	&\lesssim \big\{\|\mathbf{S}^{\d}_{xy}w\| + \|\mathbf{S}^{\d}_x\hat{w}\| \big\} \|\mathbf{S}^{\d}_yw\| 
	\lesssim L [[[\mathbf{S}^{\d}]]]_{2,\hat{w}},
\end{align}
where we have used $\eqref{2TC.3}_1$ and \eqref{1TC.21} in \eqref{1TC.61}.

Hence it follows from \eqref{1TC.59}-\eqref{1TC.61} that 
\begin{align}\label{1TC.62}
4\|\rho_s u_s \mathbf{S}^{\d}_{x}w\|^2 +  \kappa^2 \|\mathbf{S}^{\d}_{yy}w\|^2 
	&\leq  \|2\rho_su_s \mathbf{S}^{\d}_x - \kappa \Delta_{\v}\mathbf{S}^{\d}\|^2  +   L [[[\mathbf{S}^{\d}]]]_{2,\hat{w}}.
\end{align}

Recall the precise definition $\mathfrak{J}$  in \eqref{2.84}, one has that 
\begin{align}\label{1TC.59-1}
	\|w\mathfrak{J}\|^2 
	& \lesssim 
	L  \|(\phi,p,S)\hat{w}\|^2_{\mathbf{X}^{\v}}+  L\Big\{\f{1}{\v^2} [[p]]_{2,1} + \f{1}{\v} [[[S]]]_{2,1}\Big\} \mathbf{1}_{\{w=w_0\}}  + \|\mathfrak{N}_3 w\|^2,
\end{align}
which, together with \eqref{1TC.58-1}-\eqref{1TC.59} and \eqref{1TC.62}-\eqref{1TC.59-1}, yields that 
\begin{align}\label{1TC.55-2}
\|u_s \mathbf{S}^{\d}_{x}w\|^2 + \|\mathbf{S}^{\d}_{yy}w\|^2  
	&\lesssim L [[[\mathbf{S}^{\d}]]]_{2,\hat{w}} + L^{1-}  \|(\phi,p,S)\hat{w}\|^2_{\mathbf{X}^{\v}}
	+    L^2 \mathscr{K}(P^{\d}) + L^2 \f{1}{\v}\mathfrak{K}_{2,1}(P^{\d}) \nonumber\\ 
	&\quad   +  L\Big\{\f{1}{\v^2} [[p]]_{2,1} + \f{1}{\v} [[[S]]]_{2,1}\Big\} \mathbf{1}_{\{w=w_0\}}  + \|\mathfrak{N}_3 w\|^2.
\end{align}
Thus we conclude \eqref{1TC.55-1}. 
%
Therefore the proof of Lemma \ref{lemTH6.6} is completed. $\hfill\Box$

\begin{corollary}\label{cor6.6} 
It holds that 
	\begin{align}\label{1TC.56}
		[[[\mathbf{S}^{\d}]]]_{2,\hat{w}_0^{\v}}
		&\lesssim \mathscr{K}(P^{\d})+  \f{1}{\v} \mathfrak{K}_{2,1}(P^{\d}) +  \f{1}{\v} [[p]]_{2,\hat{w}_0^{\v}} +  a^2 \|(\phi,p, S)\hat{w}_0^{\v}\|^2_{\mathbf{X}^{\v}} + \f{1}{\v}\|\mathfrak{N}_3 \hat{w}_0^{\v}\|^2.
	\end{align} 
\end{corollary}

\noindent{\bf Proof.}   
Combining \eqref{1TC.21}, \eqref{1TC.51-0} and \eqref{1TC.55-1}, one obtains that 
\begin{align*}
	&\|(\sqrt{\v}\mathbf{S}^{\d}_{xx}, \mathbf{S}^{\d}_{xy})\hat{w}_0^{\v}\|^2 + L^{-\f12} \big\{\|u_s \mathbf{S}^{\d}_{x}\hat{w}_0^{\v}\|^2 + \|\mathbf{S}^{\d}_{yy}\hat{w}_0^{\v}\|^2 \big\} + \|\sqrt{u_s}\mathbf{S}^{\d}_x\hat{w}_0^{\v}\|^2_{x=0, L}\nonumber\\
	&\quad + 	L^{-1}\|\mathbf{S}^{\d}_y\hat{w}_0^{\v}\|^2_{x=0,L} + L^{-2}\|\mathbf{S}^{\d}_y\hat{w}_0^{\v}\|^2 \nonumber\\
	&\lesssim  (\xi + L^{\f38-}) [[[\mathbf{S}^{\d}]]]_{2,\hat{w}_0^{\v}}+ C_{\xi} \mathscr{K}(P^{\d})+ C_{\xi}\f{1}{\v} \mathfrak{K}_{2,1}(P^{\d}) + C_{\xi}\f{1}{\v} [[p]]_{2,\hat{w}_0^{\v}}\nonumber\\
	&\quad  + C_{\xi} a^2\|(\phi, 0, S) \hat{w}_0^{\v}\|^2_{\mathbf{X}^{\v}} + C_{\xi}  \f{1}{\v}\|\mathfrak{N}_3 \hat{w}_0^{\v}\|^2.
\end{align*}
Taking $\xi>0$ suitably small in above estimates, one concludes \eqref{1TC.56}. Therefore the proof of Corollary \ref{cor6.6} is completed.  $\hfill\Box$

\begin{lemma}\label{lemTH6.8}
	It holds that 
	\begin{align}\label{1TC.62-0}
		[[[\mathbf{S}^{\d}]]]_{b,\hat{w}^{\v}_0}
		&\lesssim  \mathscr{K}(P^{\d}) +  \f{1}{\v} \mathfrak{K}_{2,1}(P^{\d})   + [[P^{\d}]]_{3,1}   +  \f{1}{\v} [[p]]_{2,\hat{w}_0^{\v}}  + [[p]]_{3,\hat{w}^{\v}_0}  \nonumber\\
		&\quad  +  a^2 \|(\phi,p, S)\hat{w}_0^{\v}\|^2_{\mathbf{X}^{\v}}  + \f{1}{\v}\|\mathfrak{N}_3 \hat{w}_0^{\v}\|^2  + \|u_s \pa_y\mathfrak{N}_3 \hat{w}^{\v}_0\|^2.
	\end{align}
\end{lemma} 
\noindent{\bf Proof.} Applying $\pa_y$ to \eqref{2.111-6}, one gets
\begin{align}\label{1TC.58}
	& 2 \rho_s u_s \mathbf{S}^{\d}_{xy} + 2\rho_s v_s \mathbf{S}^{\d}_{yy} + 2 (\rho_s u_s)_y \mathbf{S}^{\d}_x + 2(\rho_s v_s)_y \mathbf{S}^{\d}_y  - \kappa \Delta_{\v}\mathbf{S}^{\d}_y   \nonumber\\
	&= \kappa \Delta_{\v}\Big(\frac{\chi}{2\rho_s}P^{\d}\Big)_y + \kappa \Delta_{\v}\Big(\frac{\bar{\chi}}{2\rho_s}p\Big)_y + \kappa \Delta_{\v}\Big(\f{1}{p_s}T_{sy}q\Big)_y + \mathfrak{J}_y.
\end{align}

Multiplying \eqref{1TC.58} by $u_s \mathbf{S}^{\d}_{xy}w^2$, one has 
\begin{align}\label{1TC.23}
	&2\|\sqrt{\rho_s}u_s \mathbf{S}^{\d}_{xy}w\|^2  - \kappa \la  \v \mathbf{S}^{\d}_{xxy},\, u_s \mathbf{S}^{\d}_{xy}w^2\ra - \kappa \la \mathbf{S}^{\d}_{yyy},\, u_s \mathbf{S}^{\d}_{xy} w^2\ra\nonumber\\
	&= -\Big\la 2\rho_s v_s \mathbf{S}^{\d}_{yy} + 2 (\rho_s u_s)_y \mathbf{S}^{\d}_x + 2(\rho_s v_s)_y \mathbf{S}^{\d}_y,\, u_s \mathbf{S}^{\d}_{xy}w^2\Big\ra + \Big\la \kappa \Delta_{\v}\Big(\frac{\chi}{2\rho_s}P^{\d}\Big)_y,\, u_s \mathbf{S}^{\d}_{xy}w^2\Big\ra\nonumber\\
	& + \Big\la \kappa \Delta_{\v}\Big(\frac{\bar{\chi}}{2\rho_s}p\Big)_y,\, u_s \mathbf{S}^{\d}_{xy}w^2\Big\ra
	+ \Big\la \kappa \Delta_{\v}\Big(\f{1}{p_s}T_{sy}q\Big)_y,\, u_s \mathbf{S}^{\d}_{xy}w^2 \Big\ra  + \Big\la \mathfrak{J}_y ,\, u_s \mathbf{S}^{\d}_{xy}w^2\Big\ra.
\end{align} 
In the following, we divide the proof into two steps. 

{\it Step 1.} We consider the terms on LHS of \eqref{1TC.23}. We integrate by parts to obtain
\begin{align}\label{1TC.25}
	- \kappa \la  \v \mathbf{S}^{\d}_{xxy},\, u_s \mathbf{S}^{\d}_{xy} w^2\ra 
	&=\f12 \kappa \|\sqrt{u_s}\sqrt{\v}\mathbf{S}^{\d}_{xy}w\|^2_{x=0} + O(1)  \v [[[\mathbf{S}^{\d}]]]_{2,w},
\end{align}
and 
\begin{align}\label{1TC.26}
	&- \kappa \la \mathbf{S}^{\d}_{yyy},\, u_s \mathbf{S}^{\d}_{xy}w^2 \ra \nonumber\\
	& = \kappa \la \mathbf{S}^{\d}_{yy},\, u_s \mathbf{S}^{\d}_{xyy} w^2\ra + \kappa \la \mathbf{S}^{\d}_{yy},\, u_{sy} \mathbf{S}^{\d}_{xy}w^2\ra +  \kappa \la \mathbf{S}^{\d}_{yy},\, u_s \mathbf{S}^{\d}_{xy} 2w w_y\ra  \nonumber\\
	&=\f12 \kappa \|\sqrt{u_s}\mathbf{S}^{\d}_{yy}w\|^2_{x=L} - \f12 \kappa \|\sqrt{u_s}\mathbf{S}^{\d}_{yy}w\|^2_{x=0}   + O(1) a \|\mathbf{S}^{\d}_{yy}w\|\cdot \|(\mathbf{S}^{\d}_{xy}, \mathbf{S}^{\d}_{yy})\hat{w}\| \nonumber\\
	&=\f12 \kappa \|\sqrt{u_s}\mathbf{S}^{\d}_{yy}\|^2_{x=L}  +  O(1) \Big\{\f{1}{\v}\mathfrak{K}_{2,1}(P^{\d}) + \v^{\f{1-\beta}{2}} \mathscr{K}(P^{\d}) +  L^{\f14}  [[[\mathbf{S}^{\d}]]]_{2,\hat{w}}\Big\},
\end{align}
where we have used 
\begin{align}\label{1TC.67}
	\|\sqrt{u_s}\mathbf{S}^{\d}_{yy}w\|^2_{x=0}
	&=	\|\sqrt{u_s}(\f{1}{2\rho_s}\chi P^{\d})_{yy}w\|^2_{x=0}
	 \lesssim  \v^{\f{1-\beta}{2}} \mathscr{K}(P^{\d}) + \f{1}{\v}\mathfrak{K}_{2,1}(P^{\d}) .
\end{align}

{\it Step 2.} We consider terms on RHS of \eqref{1TC.23}. For the 1th term on RHS of \eqref{1TC.23}, it holds that 
\begin{align}\label{1TC.24}
\eqref{1TC.23}_1	&\lesssim  a\|(\mathbf{S}^{\d}_{yy},\mathbf{S}^{\d}_{y})w\|\cdot \|\mathbf{S}^{\d}_{xy}w\| +  \|u_s\mathbf{S}^{\d}_xw\|\cdot \|\mathbf{S}^{\d}_{xy}w\| \lesssim L^{\f14} [[[\mathbf{S}^{\d}]]]_{2,\hat{w}} .
\end{align}

A direct calculation shows that 
\begin{align}
	\|u_s\Delta_{\v}\big(\chi \rho_s^{-1} P^{\d}\big)_yw\|^2 
	& \lesssim \|u_s \Delta_{\v}P^{\d}_y\chi\|^2 + a^2\|\chi \nabla^2_{\v}P^{\d}\|^2 + a^2\|\chi\nabla_{\v}P^{\d}\|^2 + a^2 \|\chi P^{\d}\|^2\nonumber\\
	&\lesssim [[P^{\d}]]_{3,1}  + a^2 L^2 \mathscr{K}(P^{\d}) + a^2L^2 \f{1}{\v}\mathfrak{K}_{2,1}(P^{\d}),\label{1TC.69}\\
	\|u_s\Delta_{\v}\big(\bar{\chi}\rho_s^{-1} p\big)_yw\|^2 
	& \lesssim \|u_s \Delta_{\v}p_yw\|^2 + a^2\|\big(\mathbf{d}_2\nabla^2_{\v}p, \mathbf{d}_3\nabla_{\v}p, \mathbf{d}_4 p\big)w\|^2 \nonumber\\
	&\lesssim [[p]]_{3,w} + L^2 \f{1}{\v} [[p]]_{2,1}.\label{1TC.69-1}
\end{align} 
Then it follows from \eqref{1TC.69}-\eqref{1TC.69-1} that 
\begin{align}
\eqref{1TC.23}_2&\lesssim   [[[\mathbf{S}^{\d}]]]_{2,w} + [[P^{\d}]]_{3,1}  +  L^2 \mathscr{K}(P^{\d}) + L^2 \f{1}{\v}\mathfrak{K}_{2,1}(P^{\d}), \label{1TC.71}\\
\eqref{1TC.23}_3&\lesssim   [[[\mathbf{S}^{\d}]]]_{2,w} + [[p]]_{3,w} + L^2 \f{1}{\v} [[p]]_{2,1}. \label{1TC.71-1}
\end{align}
Also it is clear that 
\begin{align}\label{1TC.72}
\eqref{1TC.23}_4 &\lesssim  [[[\mathbf{S}^{\d}]]]_{2,w} + a^2 L^{\f14} [\phi]_{3,w}.
\end{align}

Next we consider 5th on RHS of \eqref{1TC.23}. Recall $\mathfrak{J}$ in \eqref{2.84}. A direct calculation shows that  
\begin{align}\label{1TC.70}
&	\|\big(\rho_sv_s \big(p_s^{-1}T_{sy}q\big)_y\big)_yw\|^2 + \|\big(\big[\rho_su_s \big(p_s^{-1}T_{sy}\big)_x-T_{sy}\bar{u}_{sx}\big]q\big)w\|^2 \nonumber\\
& \quad + \|\big(\rho_s^{-1}p_{sy}  \phi_x\big)_yw\|^2  + \|\big(\rho_s (U_s\cdot \nabla) (\rho_s^{-1})\, p\big)_yw\|^2  \lesssim L [\phi]_{3,w}+ L^2\f{1}{\v} [[p]]_{2,1}, 
\end{align}
which implies that 
\begin{align}\label{1TC.77}
	&\big\la  \big(\rho_sv_s \big(p_s^{-1}T_{sy}q\big)_y\big)_y + \big(\big[\rho_su_s \big(p_s^{-1}T_{sy}\big)_x-  T_{sy} \tilde{u}_{sx}\big]q\big)_y + \big(\rho_s^{-1}p_{sy}  \phi_x\big)_y ,\, u_s \mathbf{S}^{\d}_{xy} w^2 \big\ra \nonumber\\
	&\quad + \big\la \big(\rho_s (U_s\cdot \nabla) (\frac{1}{\rho_s})\, p\big)_y , \, u_s\mathbf{S}^{\d}_{xy} \big\ra  \lesssim L^{\f12} [[[\mathbf{S}^{\d}]]]_{2,w} + L^{\f12} [\phi]_{3,w}  + L \f{1}{\v} [[p]]_{2,1},
\end{align}

It is clear that 
\begin{align}\label{1TC.66-0}
	\Big\la \pa_y\mathfrak{M}  + \pa_y\mathfrak{N}_3 , \, u_s\mathbf{S}^{\d}_{xy} \Big\ra
	& \lesssim  [[[\mathbf{S}^{\d}]]]_{2,w} + \|u_s(\pa_y\mathfrak{M}_0, \v \pa_y \mathfrak{M}_1, \v \pa_y \mathfrak{M}_2)w\|^2    + \|u_s \pa_y\mathfrak{N}_3 w\|^2 + \v^{2N_0}.
\end{align}
It follows from \eqref{7.47} that (for both NT \& DT)
\begin{align}\label{D3.76S}
	|\pa_y\mathfrak{M}_0 |&\lesssim \mathbf{d}_2u_s |(p_{yy}, S_{yy})| + a \mathbf{d}_2 |(p_y,S_y)| + a\mathbf{d}_2|(p,S)|  + \mathbf{d}_2 |(\phi_{yyy}, \phi_{yy}, \phi_y)| \nonumber\\
	&\quad  + u_s \mathbf{d}_2^2 |q_{yy}|+ a \mathbf{d}_2^2 |(q_y,q)|,
\end{align}
which yields that (for both NT \& DT)
\begin{align}\label{1TC.79}
	\|\pa_y\mathfrak{M}_0 w\|^2 &\lesssim L^{\f18} [\phi]_{3,1}   +  L \f{1}{\v} [[p]]_{2,1} + L^{\f18} [[[S]]]_{2,1}   +  L \Big\{\f{1}{\v^2} [[p]]_{2,1} + \f{1}{\v}[[[S]]]_{2,1}\Big\} \mathbf{1}_{\{w=w_0\}}. 
\end{align}
For $\mathfrak{M}_1, \mathfrak{M}_2$, we have from \eqref{7.35-2}-\eqref{7.36-2} that (for both NT \& DT)
\begin{align*}
|\pa_y\big(\mathfrak{M}_1 , \mathfrak{M}_2 \big) \big| 
	&\lesssim   u_s \mathbf{d}_1 |(p_{yy},p_{xy}, S_{xy}, S_{yy})| + u_s \mathbf{d}_2 |(p_x,p_y,S_x,S_y)| + u_s\mathbf{d}_2 |(p,S)| \nonumber\\
	&\quad + \mathbf{d}_1 |(\phi_{yyy}, \phi_{xyy}, \phi_{xxy}, \nabla^2\phi, \nabla \phi)| + a u_s\mathbf{d}_2 |(q_{xy}, q_{yy}, q_x,q_y,q)|,
\end{align*} 
which yields that (for both NT \& DT)
\begin{align}\label{1TC.81}
	\v \|w\pa_y\big(\mathfrak{M}_1 , \mathfrak{M}_2 \big)\|^2 
	&\lesssim  
L \|(\phi,p,S)\hat{w}\|^2_{\mathbf{X}^{\v}} +  L \Big\{\f{1}{\v} [[p]]_{2,1} + [[[S]]]_{2,1}\Big\} \mathbf{1}_{\{w=w_0\}}.
\end{align} 

Recall again $\mathfrak{J}$ in \eqref{2.84}, we have from \eqref{1TC.77}-\eqref{1TC.81} that 
\begin{align}\label{1TC.66}
	\eqref{1TC.23}_5	&\lesssim [[[\mathbf{S}^{\d}]]]_{2,w} +  L^{\f18} \|(\phi,p,S)\hat{w}\|^2_{\mathbf{X}^{\v}} 
	+  L \Big\{\f{1}{\v^2} [[p]]_{2,1} + \f{1}{\v}[[[S]]]_{2,1}\Big\} \mathbf{1}_{\{w=w_0\}}   \nonumber\\
	&\quad+ \|u_s \pa_y\mathfrak{N}_3 w\|^2 .
\end{align}

Substituting \eqref{1TC.25}-\eqref{1TC.24}, \eqref{1TC.71}-\eqref{1TC.72}, and \eqref{1TC.66} into \eqref{1TC.23}  one obtains that 
\begin{align}\label{1C.67}
	&  \|\sqrt{u_s}\sqrt{\v}\mathbf{S}^{\d}_{xy}w\|^2_{x=0}  + \|\sqrt{u_s}\mathbf{S}^{\d}_{yy}w\|^2_{x=L}   \nonumber\\
	&\lesssim [[[\mathbf{S}^{\d}]]]_{2,\hat{w}}   + [[P^{\d}]]_{3,1}  + \f{1}{\v}\mathfrak{K}_{2,1}(P^{\d}) + L^2 \mathscr{K}(P^{\d})   + [[p]]_{3,w}  +  a^2 \|(\phi,p,S)\hat{w}\|^2_{\mathbf{X}^{\v}}   \nonumber\\
	&\quad + L\f{1}{\v} \|(0,p,S)\|^2_{\mathbf{X}^{\v}} \mathbf{1}_{\{w=w_0\}} + \|u_s \pa_y\mathfrak{N}_3 w\|^2,
\end{align}  
which, together with \eqref{1TC.56}, concludes \eqref{1TC.62-0}.
Therefore the proof of Lemma \ref{lemTH6.8} is completed.  $\hfill\Box$

\begin{lemma}\label{lemTH6.9}
	It holds that 
	\begin{align} \label{1TC.85}
		[[[\mathbf{S}^{\d}]]]_{3,\hat{w}^{\v}_0}
		&\lesssim  \sqrt{\v} [[[\mathbf{S}^{\d}]]]_{b,1} + L^{\f12}  [[[\mathbf{S}^{\d}]]]_{2,\hat{w}^{\v}_0} +  [[P^{\d}]]_{3,1}  +  L^2 \mathscr{K}(P^{\d})  + a^2 \f{1}{\v}\mathfrak{K}_{2,1}(P^{\d}) \nonumber\\
		&\quad + [[p]]_{3,\hat{w}^{\v}_0}  +  L^{\f18} \|(\phi,p,S) \hat{w}^{\v}_0\|^2_{\mathbf{X}^{\v}} + \|u_s \pa_y\mathfrak{N}_3 \hat{w}^{\v}_0\|^2 .
	\end{align}
\end{lemma}

\noindent{\bf Proof.} It follows from \eqref{1TC.58} and \eqref{1TC.69}-\eqref{1TC.69-1}  that 
\begin{align}\label{1TC.68}
	& 4\|\rho_su_s^2 \mathbf{S}^{\d}_{xy}w\|^2  + \kappa^2 \|u_s(\v  \mathbf{S}^{\d}_{xxy}, \mathbf{S}^{\d}_{yyy})w\|^2   - 4 \kappa  \big\la \rho_s u_s^3 \mathbf{S}^{\d}_{xy}, \, \Delta_{\v}\mathbf{S}^{\d}_{y} w^2 \big\ra  + 2\kappa^2 \big\la \v u_s^2 \mathbf{S}^{\d}_{xxy}, \mathbf{S}^{\d}_{yyy} w^2 \big\ra\nonumber\\
	&\lesssim L^{\f38} a^2  [[[\mathbf{S}^{\d}]]]_{2,w} +  [[P^{\d}]]_{3,1}  + a^2 L^2 \mathscr{K}(P^{\d}) + a^2L^2 \f{1}{\v}\mathfrak{K}_{2,1}(P^{\d})  + L^{\f18} [\phi]_{3,\hat{w}}  + [[p]]_{3,w} \nonumber\\
	&\quad   + L  \f{1}{\v} [[p]]_{2,\hat{w}}  + L^{\f18} [[[S]]]_{2,\hat{w}}  +  L \Big\{\f{1}{\v^2} [[p]]_{2,1} + \f{1}{\v}[[[S]]]_{2,1}\Big\} \mathbf{1}_{\{w=w_0\}} + \|u_s \pa_y\mathfrak{N}_3 w\|^2 .
\end{align}

We integrate by parts to obtain
\begin{align}\label{1TC.75}
\big\la \v u_s^2 \mathbf{S}^{\d}_{xxy}, \mathbf{S}^{\d}_{yyy} w^2\big\ra 
	&=   - \big\la \v u_s^2 \mathbf{S}^{\d}_{xy}, \mathbf{S}^{\d}_{xyyy} w^2\big\ra -  \big\la \v 2u_su_{sx} \mathbf{S}^{\d}_{xy}, \mathbf{S}^{\d}_{yyy}w^2 \big\ra - \big\la \v u_s^2 \mathbf{S}^{\d}_{xy}, \mathbf{S}^{\d}_{yyy} w^2\big\ra_{x=0} \nonumber\\
	&=\|u_s\sqrt{\v}\mathbf{S}^{\d}_{xyy}w\|^2  + O(1) a\sqrt{\v} \|\mathbf{S}^{\d}_{xy}\hat{w}\|\cdot \|u_s(\sqrt{\v} \mathbf{S}^{\d}_{xyy}, \sqrt{\v} \mathbf{S}^{\d}_{yyy})w\| \nonumber\\
	&\quad   - \big\la \v u_s^2 \mathbf{S}^{\d}_{xy}, \mathbf{S}^{\d}_{yyy} w^2\big\ra_{x=0}.
\end{align}
To control the boundary term on RHS of \eqref{1TC.75}, we have from \eqref{1TC.67} that 
\begin{align*}
	\int_0^\infty \f{1}{y^{\ast}}\big|u_s (\f{\chi}{\rho_s}  P^{\d})_{yy}w\big|^2 dy \Big|_{x=0} \lesssim  a\|\sqrt{u_s} \big(\f{\chi}{\rho_s} P^{\d}\big)_{yy} w\|^2_{x=0}
	\lesssim a\v^{\f{1-\beta}{2}} \mathscr{K}(P^{\d}) + a\f{1}{\v}\mathfrak{K}_{2,1}(P^{\d}),
\end{align*}
which shows that $u_sw\mathbf{S}^{\d}_{yy}\big|_{x=0} \in H^{1/2}_{00}(\mathbb{R}_+)$. Then one obtains  that 
\begin{align}\label{1TC.76}
	&\big\la \v u_s^2 \mathbf{S}^{\d}_{xy}, \mathbf{S}^{\d}_{yyy} w^2\big\ra_{x=0}  = \big\la \v u_sw\mathbf{S}^{\d}_{xy}, \big(u_sw\mathbf{S}^{\d}_{yy}\big)_y  \big\ra_{x=0} + \big\la \v u_sw \mathbf{S}^{\d}_{xy}, (u_sw)_y \mathbf{S}^{\d}_{yy} \big\ra_{x=0} \nonumber\\
	&= - \big\la \v u_s\mathbf{S}^{\d}_{xy}, \big(u_sw (\f{\chi}{2\rho_s} P^{\d})_{yy}\big)_y \big\ra_{x=0} - \big\la \v u_sw \mathbf{S}^{\d}_{xy},   (u_sw)_y \big(\f{\chi}{2\rho_s}  P^{\d}\big)_{yy}\big\ra_{x=0} \nonumber\\
	&\lesssim \v\|(u_sw \chi \mathbf{S}^{\d}_{xy})(0,\cdot)\|_{H^{1/2}}  \|\big(u_sw\big(\f{\chi}{\rho_s} P^{\d}\big)_{yy}\big)(0,\cdot)\|_{H^{1/2}_{00}} + a \v \|\sqrt{u_s}\mathbf{S}^{\d}_{xy}\|_{x=0}  \|\sqrt{u_s}\big(\f{\chi}{\rho_s} P^{\d}\big)_{yy}\|_{x=0}\nonumber\\
	&\lesssim \xi  [[[\mathbf{S}^{\d}]]]_{3,1}  + \sqrt{\v} [[[\mathbf{S}^{\d}]]]_{2,1}+ C_{\xi} [[P^{\d}]]_{3,1} + C_{\xi}L^2 \mathscr{K}(P^{\d}) + C_{\xi}L^2 \f{1}{\v}\mathfrak{K}_{2,1}(P^{\d}),
\end{align}
where we have used \eqref{D3.30-0} and the following two estimates
\begin{align}
	\begin{split}
		\v^{\f32} \|(u_sw\chi\mathbf{S}^{\d}_{xy})(0,\cdot)\|^2_{H^{1/2}}&\lesssim \|u_s(\sqrt{\v}\mathbf{S}^{\d}_{xyy}, \v \mathbf{S}^{\d}_{xxy})\|^2 +  \v^{1-} \|\mathbf{S}^{\d}_{xy}\|^2, \\
		\sqrt{\v}\|(u_sw (\f{\chi}{\rho_s}P^{\d})_{yy})(0,\cdot)\|^2_{H^{1/2}_{00}}
		&\lesssim [[P^{\d}]]_{3,1}  + L^2 \mathscr{K}(P^{\d}) + L^2 \f{1}{\v}\mathfrak{K}_{2,1}(P^{\d}).
	\end{split}
\end{align}
Hence it follows from \eqref{1TC.75}-\eqref{1TC.76} that 
\begin{align}\label{1TC.75-1}
	\big\la \v u_s^2 \mathbf{S}^{\d}_{xxy}, \mathbf{S}^{\d}_{yyy} w^2\big\ra 
	&\geq \|u_s\sqrt{\v}\mathbf{S}^{\d}_{xyy}w\|^2 - \xi  [[[\mathbf{S}^{\d}]]]_{3,\hat{w}}  - \sqrt{\v} [[[\mathbf{S}^{\d}]]]_{2,\hat{w}} - C_{\xi} [[P^{\d}]]_{3,1}\nonumber\\
	&\quad -  C_{\xi}L^2 \mathscr{K}(P^{\d})- C_{\xi}L^2 \f{1}{\v}\mathfrak{K}_{2,1}(P^{\d}).
\end{align}

Integrating by parts, one has that 
\begin{align}\label{1C.79-1}
	- 4 \kappa  \big\la \rho_s u_s^3 \mathbf{S}^{\d}_{xy}, \v \mathbf{S}^{\d}_{xxy}w^2\big\ra = 2\kappa \|\sqrt{\rho_su_s^3}\sqrt{\v} \mathbf{S}^{\d}_{xy}w\|^2_{x=0} + O(1)\v [[[\mathbf{S}^{\d}]]]_{2,w},
\end{align}
and
\begin{align}\label{1TC.79-2}
	&- 4\kappa \big\la \rho_s u_s^3 \mathbf{S}^{\d}_{xy}, \mathbf{S}^{\d}_{yyy} w^2\big\ra 
	=  4\kappa \big\la \rho_s u_s^3 \mathbf{S}^{\d}_{xyy}, \mathbf{S}^{\d}_{yy} w^2\big\ra + \kappa \big\la (\rho_s u_s^3 w^2)_{y} \mathbf{S}^{\d}_{xy}, \mathbf{S}^{\d}_{yy}\big\ra \nonumber\\
	&\geq 2\kappa \|\sqrt{\rho_s u_s^3} \mathbf{S}^{\d}_{yy}w\|^2_{x=L} - 2\kappa \|\sqrt{\rho_s u_s^3} \mathbf{S}^{\d}_{yy}\|^2_{x=0} - C a^3 \|\mathbf{S}^{\d}_{yy}w\|^2 - C a\|u_s^2\mathbf{S}^{\d}_{xy}w\|\cdot \|\mathbf{S}^{\d}_{yy}w_y\|  \nonumber\\
	&\geq 2\kappa \|\sqrt{\rho_s u_s^3} \mathbf{S}^{\d}_{yy}w\|^2_{x=L} - \f18 \|\sqrt{\rho_s}u_s^2\mathbf{S}^{\d}_{xy}w\|^2 -  a^2 \f{1}{\v}\mathfrak{K}_{2,1}(P^{\d}) - \v^{\f{1-\beta}{2}} \mathscr{K}(P^{\d})  - a^2 \|\mathbf{S}^{\d}_{yy}\hat{w}\|^2,
\end{align}
where we have used \eqref{1TC.67} in the last inequality of \eqref{1TC.79-2}.

Substituting \eqref{1TC.75-1}-\eqref{1TC.79-2} into \eqref{1TC.68}, then using \eqref{1TC.79}, \eqref{1TC.81} , one gets that 
\begin{align*}
	&\|u_s^2 \mathbf{S}^{\d}_{xy}w\|^2   + \|u_s(\v \mathbf{S}^{\d}_{xxy}, \sqrt{\v} \mathbf{S}^{\d}_{xyy}, \mathbf{S}^{\d}_{yyy})w\|^2 + \|\sqrt{u_s^3} \mathbf{S}^{\d}_{yy}w\|^2_{x=L} + \|\sqrt{u_s^3}\sqrt{\v} \mathbf{S}^{\d}_{xy}w\|^2_{x=0}  \nonumber\\
	&\lesssim  \xi  [[[\mathbf{S}^{\d}]]]_{3,1}  + L^{\f12}  [[[\mathbf{S}^{\d}]]]_{2,\hat{w}} + C_{\xi}[[P^{\d}]]_{3,1}  + C_{\xi}L^2 \mathscr{K}(P^{\d})   + a^2 \f{1}{\v}\mathfrak{K}_{2,1}(P^{\d}) + [[p]]_{3,\hat{w}} \nonumber\\
	&\quad +  L^{\f18} \|(\phi,p,S) \hat{w} \|^2_{\mathbf{X}^{\v}} +  L \f{1}{\v}\Big\{\f{1}{\v} [[p]]_{2,1} + [[[S]]]_{2,1}\Big\} \mathbf{1}_{\{w=w_0\}}  + \|u_s \pa_y\mathfrak{N}_3 w\|^2 ,
\end{align*}
which immediately concludes \eqref{1TC.85}. Therefore the proof of Lemma \ref{lemTH6.9} is completed.  $\hfill\Box$

\medskip

\noindent{\bf Proof of Theorem \ref{thm6.7}.}  By taking the limit $\delta\to0+$ in \eqref{2.111-6}, and noting \eqref{1TC.56}, \eqref{1TC.62-0} and \eqref{1TC.85},  we conclude Theorem \ref{thm6.7} directly. Therefore the proof is completed. $\hfill\Box$

\subsection{Existence and uniform estimate for Neumann BVP \eqref{2.111-1}} \label{secT.2}
\begin{theorem}\label{THlem6.10}
There exist a unique strong solution $\zeta$ of \eqref{2.111-1} with the following uniform estimates 
	\begin{align}
		[[[\zeta]]]_{2,\hat{w}_0^{\v}}
		&\lesssim  \f{1}{\v} [[p]]_{2,\hat{w}_0^{\v}}  +  a^2 \|(\phi,p, S)\hat{w}_0^{\v}\|^2_{\mathbf{X}^{\v}} + \f{1}{\v}\|\mathfrak{N}_3 \hat{w}_0^{\v}\|^2,\label{1ZC.56}\\
		[[[\zeta]]]_{3,\hat{w}^{\v}_0}
		&\lesssim [[p]]_{3,\hat{w}^{\v}_0} +  a^2  \|(\phi,p,S (\zeta)) \hat{w}^{\v}_0\|^2_{\mathbf{X}^{\v}} + \|u_s \pa_y\mathfrak{N}_3 \hat{w}^{\v}_0\|^2,\label{1ZC.85}\\
		[[[\zeta]]]_{b,\hat{w}_0^{\v}} 
		&\lesssim \f{1}{\v} [[p]]_{2,\hat{w}_0^{\v}}  + [[p]]_{3,\hat{w}^{\v}_0}   +  a^2 \|(\phi,p, S)\hat{w}_0^{\v}\|^2_{\mathbf{X}^{\v}} + \f{1}{\v}\|\mathfrak{N}_3 \hat{w}_0^{\v}\|^2  + \|u_s \pa_y\mathfrak{N}_3 \hat{w}^{\v}_0\|^2.\label{1ZC.62-0}
	\end{align}
\end{theorem}

\noindent{\bf Proof.} Since $\nabla^3p\notin L^2$, for strictness, we need first to study the following approximate problem
\begin{align}\label{1ZC.94}
\begin{cases}
\dis	2 \rho_s u_s \zeta^{\epsilon}_x + 2\rho_s v_s \zeta^{\e}_y - \kappa \Delta_{\v}\zeta^{\e}  = \kappa \Delta_{\v}\Big(\frac{1}{2\rho_s}\big[\chi p^{\e} + \bar{\chi} p\big]\Big) + \kappa \Delta_{\v}\Big(\f{1}{p_s}T_{sy}q\Big) + \mathfrak{J},\\
\dis \zeta \big|_{x=0}= -  \frac{\chi}{2\rho_s}p^{\epsilon},\,\,\, \zeta_{x}\big|_{x=L}=0 \,\, \,\mbox{and}\,\,\,
\zeta_y \big|_{y=0} = - \big(\frac{1}{2\rho_s}p^{\epsilon}\big)_y,
\end{cases}
\end{align}
where $p^{\e}:=p(x,y+\e)$ with $\epsilon\ll \v$. Although $p^{\epsilon}(0,0) \neq 0$, we observe that $p^{\epsilon}(0,0) = p(0,\epsilon) = O(\sqrt{\epsilon}) \ll \v$, which does not introduce  additional difficulty in the proof.  Now it is clear that $\chi p^{\e}\in H^3$. Then we apply the existence and uniform estimates in Section  \ref{secT.1} to \eqref{1ZC.94}, and finally take the limit $\e\to 0+$ to complete the proof of Theorem \ref{THlem6.10}. $\hfill\Box$

\begin{lemma}\label{THlem6.11}
	Recall $\dis\tilde{T}(\zeta,p,q)=\zeta + \f{1}{2\rho_s}p + \f{1}{p_s}T_{sy}q$,  we have 
	\begin{align}
		\|\sqrt{u_s}\Delta_{\v}\tilde{T}\hat{w}_0^{\v}\|^2_{x=0} 
		&\lesssim   a^2 \|(\phi,p,S(\zeta))\hat{w}_0^{\v}\|^2_{\mathbf{X}^{\v}} +  \|\sqrt{u_s}\mathfrak{N}_3 \hat{w}_0^{\v}\|^2_{x=0} + \v^{2N_0},\label{1ZC.107}\\
		\|\nabla_{\v}\Delta_{\v}\tilde{T}\hat{w}_0^{\v}\|^2 
		&\lesssim  a^2 \|(\phi,p,S(\zeta))\hat{w}_0^{\v}\|^2_{\mathbf{X}^{\v}}  +  \|\nabla_{\v}\mathfrak{N}_3 \hat{w}_0^{\v}\|^2 + \v^{2N_0},\label{1ZC.107-1}\\
		\|u_s \nabla_{\v}\Delta_{\v}\tilde{T}_y \hat{w}_0^{\v}\|^2 
		&\lesssim   a^2 \|(\phi,p,S(\zeta))\hat{w}_0^{\v}\|^2_{\mathbf{X}^{\v}}  +  \|u_s \nabla_{\v}\pa_y\mathfrak{N}_3 {\hat{w}_0^{\v}}\|^2 + \v^{2N_0}.\label{1ZC.107-2}
	\end{align} 
\end{lemma}

\noindent{\bf Proof.} It follows from \eqref{2.111-1} that 
\begin{align}\label{1ZC.96}
	\kappa \Delta_{\v}\tilde{T} = 2 \rho_s u_s \zeta_x + 2\rho_s v_s \zeta_y  - \mathfrak{J}(\phi,p,T(S,p,q)),
\end{align}
which yields that 
\begin{align}\label{1ZC.97}
	\begin{split}
		\|\sqrt{u_s}\Delta_{\v}\tilde{T}w\|^2_{x=0} 
		&\lesssim  a^2 [[[\zeta]]]_{2,w} + \f{1}{\sqrt{\v}} [[p]]_{2,1} + \|\sqrt{u_s}\mathfrak{J}w\|^2_{x=0},\\
		\|\nabla_{\v}\Delta_{\v}\tilde{T}w\|^2 
		&\lesssim a^2 [[[\zeta]]]_{2,w} + \|\nabla_{\v}\mathfrak{J}w\|^2,\\
		\|u_s \nabla_{\v}\Delta_{\v}\tilde{T}_y w\|^2 
		&\lesssim a^2 [[[\zeta]]] _{3,w}  +  a^2 [[[\zeta]]]_{2,w} + \|u_s\nabla_{\v}\mathfrak{J}_y  w\|^2. 
	\end{split}
\end{align}

Recall the definition of $\mathfrak{J}$ in \eqref{2.84}. Noting \eqref{1C.43-1} and \eqref{1TC.47},  one gets that 
\begin{align}\label{1ZC.111}
&\|\sqrt{u_s}\mathfrak{M}_0w\|^2_{x=0}
\lesssim  \|p_y\|^2_{x=0} +  \f{1}{\v} \|p_y\|^2_{x=0}\cdot \mathbf{1}_{\{w=w_0\}}  \lesssim  \f{1}{\sqrt{\v}} [[p]]_{2,1} +  \f{1}{\v^{\f32}} [[p]]_{2,1} \mathbf{1}_{\{w=w_0\}},\\
&\|\sqrt{u_s}(\v \mathfrak{M}_1, \v \mathfrak{M}_2) w\|^2_{x=0}  
 \lesssim \v \|p_y\|^2_{x=0} + \v \|\sqrt{\v}p_x\|^2_{x=0} + a^2 \v^2 \|\sqrt{u_s}S_x\|^2_{x=0} + \v [\phi]_{3,w} \nonumber\\
&\qquad \qquad\qquad \qquad\qquad\, \lesssim  \sqrt{\v} [[p]]_{2,1} + \v^2 [[[S]]]_{2,1}  + \v [\phi]_{3,w}.\label{1ZC.112}
\end{align}
Then we have from \eqref{2.84} and \eqref{1ZC.111}-\eqref{1ZC.112} that 
\begin{align*}
	\|\sqrt{u_s}\mathfrak{J} w\|^2_{x=0}
	&\lesssim L \|(\phi,p,S) \hat{w}\|^2_{\mathbf{X}^{\v}}  +  \f{1}{\v^{\f32}} [[p]]_{2,1} \mathbf{1}_{\{w=w_0\}}  + \|\sqrt{u_s}\mathfrak{N}_3 w\|^2 + \v^{2N_0},
\end{align*}
which, together with $\eqref{1ZC.97}_1$, yields that
\begin{align*}
	\|\sqrt{u_s}\Delta_{\v}\tilde{T}w\|^2_{x=0} 
	&\lesssim   
	a^2 \|(\phi,p,S\,(\zeta)) \hat{w}\|^2_{\mathbf{X}^{\v}} +  \f{1}{\v^{\f32}} [[p]]_{2,1}\mathbf{1}_{\{w=w_0\}}   + \|\sqrt{u_s}\mathfrak{N}_3 w\|^2_{x=0} + \v^{2N_0}.
\end{align*} 
Then it is easy to conclude \eqref{1ZC.107}.

\smallskip

It follows from \eqref{7.47} that
\begin{align}
	|\nabla_{\v}\mathfrak{M}_0|&\lesssim \mathbf{d}_2u_s |\nabla_{\v}(p_{y}, S_{y})| + a \mathbf{d}_2 |\nabla_{\v}(p, S)| + a\mathbf{d}_2|(p,S)|  + \mathbf{d}_2 |\nabla_{\v}(\phi_{yy},  \phi_{y}, \phi)| \nonumber\\
	&\quad  + u_s \mathbf{d}_2^2 |\nabla_{\v}q_{y}|+ a \mathbf{d}_2^2 |(\nabla_{\v}q,q)|,\\
	|\nabla_{\v}\big(\mathfrak{M}_1, \mathfrak{M}_2\big)|&\lesssim   u_s \mathbf{d}_1 |\nabla_{\v}(p_{y}, p_{x}, S_{x}, S_{y})| + u_s \mathbf{d}_2 |(p_x,p_y,S_x,S_y)| + u_s\mathbf{d}_2 |(p,S)| \nonumber\\
	&\quad + \mathbf{d}_1 |(\nabla_{\v}\nabla^2\phi, \mathbf{d}_2\nabla^2\phi, \mathbf{d}_2\nabla \phi)| + a u_s\mathbf{d}_2 |(\nabla_{\v}\nabla q, \nabla q, q)|,
\end{align} 
which immediately  yields that 
\begin{align}\label{1ZC.114}
	\|\nabla_{\v} (\mathfrak{M}_0, \v \mathfrak{M}_1, \v \mathfrak{M}_2) w\|^2 
	&\lesssim 
	L^{\f18} \|(\phi,p,S) \hat{w}\|^2_{\mathbf{X}^{\v}} +  L \Big\{\f{1}{\v^2} [[p]]_{2,1} + \f{1}{\v}[[[S]]]_{2,1}\Big\} \cdot \mathbf{1}_{\{w=w_0\}}, 
\end{align}
Hence we have from \eqref{2.84} and \eqref{1ZC.114} that 
\begin{align}\label{1ZC.115}
	\|\nabla_{\v}\mathfrak{J}w\|^2
	& \lesssim 
	L^{\f18} \|(\phi,p,S) \hat{w}\|^2_{\mathbf{X}^{\v}} +  L \Big\{\f{1}{\v^2} [[p]]_{2,1} + \f{1}{\v}[[[S]]]_{2,1}\Big\}  \mathbf{1}_{\{w=w_0\}} + \|\nabla_{\v}\mathfrak{N}_3 w\|^2 + \v^{2N_0}.
\end{align}
Similarly, we can have 
\begin{align}\label{1ZC.116}
	\|u_s\nabla_{\v}\mathfrak{J}_y  w\|^2
	&\lesssim 
	[\phi]_{4,w}  + a^2 \|(\phi,p,S) \hat{w}\|^2_{\mathbf{X}^{\v}}  +  \|u_s \nabla_{\v}\pa_y\mathfrak{N}_3 w\|^2 + \v^{2N_0}.
\end{align}

Substituting \eqref{1ZC.115}-\eqref{1ZC.116} into \eqref{1ZC.97}, one has 
\begin{align*}
	\|\nabla_{\v}\Delta_{\v}\tilde{T}w\|^2 
	&\lesssim  
	a^2 \|(\phi,p,S\,(\zeta)) \hat{w}\|^2_{\mathbf{X}^{\v}} + L \Big\{\f{1}{\v^2} [[p]]_{2,1} + \f{1}{\v}[[[S]]]_{2,1}\Big\} \mathbf{1}_{\{w=w_0\}}  
	+  \|\nabla_{\v}\mathfrak{N}_3 w\|^2 + \v^{2N_0},\\
	\|u_s \nabla_{\v}\Delta_{\v}\tilde{T}_y w\|^2 
	&\lesssim  
	a^2 \|(\phi,p,S\,(\zeta)) \hat{w}\|^2_{\mathbf{X}^{\v}} +  \|u_s \nabla_{\v}\pa_y\mathfrak{N}_3 w\|^2 + \v^{2N_0}.
\end{align*}
Then it is direct to conclude \eqref{1ZC.107-1}-\eqref{1ZC.107-2}. Therefore the proof of Lemma \ref{THlem6.11} is completed. $\hfill\Box$

\smallskip
\subsection{Existence and uniform estimate for  Dirichlet BVP \eqref{2.111-6D}}\label{SecD7}
\begin{theorem}\label{thm7.1}
There exists a unique  strong solution $\mathbf{S}$ to \eqref{2.111-6D} with the following uniform estimates
\begin{align}
	[[[\mathbf{S}]]]_{b,\hat{w}^{\v}_0}&\lesssim_{C_0} [[[\mathbf{S}]]]_{3,\hat{w}^{\v}_0} + [[[\mathbf{S}]]]_{2,\hat{w}^{\v}_0}  + L^{\f14\alpha-} \|(\phi,p,S) \hat{w}^{\v}_0\|^2_{\mathbf{X}^{\v}},\label{18.9-1}\\
	[[[\mathbf{S}]]]_{2,\hat{w}^{\v}_0} &\lesssim_{C_0}    [\phi]_{4,1} + \f1{\v}[[p]]_{2,\hat{w}^{\v}_0}	+ a\|(\phi, p, S)\hat{w}^{\v}_0\|^2_{\mathbf{X}^{\v}}     + \f{1}{\v}\|\mathfrak{N}_3 \hat{w}^{\v}_0\|^2   +  \f{1}{\v}\|(\mathfrak{h}_1, \mathfrak{N}_{11})\|^2_{H^1} +  \v^{N_0},\\
	[[[\mathbf{S}]]]_{3,\hat{w}^{\v}_0} 
	&\lesssim_{C_0} a^4  [[[\mathbf{S}]]]_{2,\hat{w}^{\v}_0} +  \|u_s \Delta_{\v}p_y \hat{w}^{\v}_0\|^2  
	+ L^{\f18} \|(\phi, p, S)\hat{w}^{\v}_0\|^2_{\mathbf{X}^{\v}} +  \|u_s \pa_y\mathfrak{N}_3 \hat{w}^{\v}_0\|^2. \label{18.9-3}
\end{align}
\end{theorem}

Applying Proposition \ref{propC.1} (see Appendix \ref{secC} for more details), we  obtain the existence of solution $\mathbf{S}$ for \eqref{2.111-6D}.  In the following, we need only to establish some uniform estimates for $\mathbf{S}$.

 
\begin{lemma}\label{lemAD.0}
	It holds that 
	\begin{align}\label{D3.2-0}
		L^{-1}\|\mathbf{S}_y w\|^2_{x=0, L} +L^{-2}\|\mathbf{S}_yw\|^2 &\lesssim \|\mathbf{S}_{xy}w\|^2.
	\end{align}
\end{lemma}

\noindent{\bf Proof.} By noting $\mathbf{S}|_{x=0}=0$, one has
\begin{align*}
	\begin{split}
		\|\mathbf{S}_yw\|^2 &\lesssim 
		L^2\|\mathbf{S}_{xy}w\|^2,\\ 
		\|\mathbf{S}_y w\|^2_{x=L}&\lesssim \|\mathbf{S}_{xy}w\|\cdot \|\mathbf{S}_yw\|  \lesssim L\|\mathbf{S}_{xy}w\|^2.
	\end{split}
\end{align*}
Therefore the proof of Lemma \ref{lemAD.0} is completed. $\hfill\Box$

\smallskip

\begin{lemma}\label{lemD3.8.2}
	It holds that 
	\begin{align}\label{D3.30-1}
		[[[\mathbf{S}]]]_{b,\hat{w}^{\v}_0}&\lesssim [[[\mathbf{S}]]]_{3,\hat{w}^{\v}_0} + [[[\mathbf{S}]]]_{2,\hat{w}^{\v}_0}  + L^{\f14\alpha-} \|(\phi,p,S) \hat{w}^{\v}_0\|^2_{\mathbf{X}^{\v}}.
	\end{align}
\end{lemma}

\noindent{\bf Proof.}  
Noting $\eqref{1C.43-1}_1$, one has that 
\begin{align*}
	\|\mathfrak{M}_0 w\|^2_{x=L}&\lesssim a^4\|(p_y, S_y)w\|^2_{x=L} + a^2\|\mathbf{d}_2 (p, S)w\|^2_{x=L} + L^{\f12} [\phi]_{3,\hat{w}} \nonumber\\
	&\lesssim La^2 [[[S]]]_{2,\hat{w}}  + L\f{1}{\v} [[p]]_{2,w}  + L\f{1}{\v^2} [[p]]_{2,1} \mathbf{1}_{\{w=w_0\}}  + L^{\f12} [\phi]_{3,\hat{w}},
\end{align*}
which, together with $\eqref{2.111-6D}_2$ and \eqref{AD.1}, yields that
\begin{align}\label{D3.30}
	\|\mathbf{S}_{yy}w\|^2_{x=L}&\lesssim a\|\sqrt{u_s}\mathbf{S}_xw\|^2_{x=L}  
	+ L^{\f14\alpha-} \|(\phi,p,S) \hat{w}\|^2_{\mathbf{X}^{\v}}
	+ L\f{1}{\v^2} [[p]]_{2,1} \mathbf{1}_{\{w=w_0\}}.
\end{align}
From \eqref{D3.30} and \eqref{D3.30-0}, one
concludes \eqref{D3.30-1}. Therefore the proof of Lemma \ref{lemD3.8.2} is completed. $\hfill\Box$

\begin{lemma}\label{lemD3.1}
It holds that 
\begin{align} \label{D3.8}
& \|(\sqrt{\v}\mathbf{S}_{xx}, \mathbf{S}_{xy})w\|^2 + \|\sqrt{u_s} \mathbf{S}_x w\|^2_{x=L} + \|\sqrt{u_s} \mathbf{S}_xw\|^2_{x=0}  \nonumber\\
&\lesssim \xi \|(\sqrt{\v}\mathbf{S}_{xx}, \mathbf{S}_{xy} )\|^2   + \xi \|\mathbf{S}_{yy}\|^2_{x=L}    + C_{\xi} \|\mathbf{S}_xw_y\|^2  +  L^{\f14-} [[[\mathbf{S}]]]_{2,\hat{w}}  + C_{\xi} [\phi]_{4,1}  \nonumber\\
&\quad 	+  C_{\xi} \f1{\v}[[p]]_{2,\hat{w}} + C_{\xi}a \|(\phi, p, S) \hat{w}\|^2_{\mathbf{X}^{\v}} + C_{\xi} L\f{1}{\v} \Big\{\f{1}{\v} [[p]]_{2,1} + [[[S]]]_{2,1}\Big\} \mathbf{1}_{\{w=w_0\}}   
 \nonumber\\
 &\quad  + C_{\xi}  \f{1}{\v}\|\mathfrak{N}_3 w\|^2   +  C_{\xi}\f{1}{\v}\|(\mathfrak{h}_1, \mathfrak{N}_{11})\|^2_{H^1} +  \v^{N_0}.
\end{align}
\end{lemma}

\noindent{\bf Proof.} 
Multiplying $\eqref{2.111-6D}_1$ by $\mathbf{S}_{xx} w^2$, we have 
\begin{align}\label{D3.2}
	\kappa \|\sqrt{\v}\mathbf{S}_{xx}\|^2 + \big\la \mathbf{S}_{xx} w^2, \, \kappa \mathbf{S}_{yy}- 2\rho_s u_s \mathbf{S}_x\big\ra  
	&= \big\la \mathbf{S}_{xx} w^2, \,  2\rho_s v_s \mathbf{S}_y\big\ra - \f{1}{2}\big\la \mathbf{S}_{xx}, \,  \kappa \Delta_{\v}\big(\rho_s^{-1}p\big)\big\ra \nonumber\\
	&\quad  - \big\la \mathbf{S}_{xx} w^2, \, \kappa \Delta_{\v}\big(p_s^{-1}T_{sy}q\big)\big\ra  - \big\la \mathbf{S}_{xx}w^2, \, \mathfrak{J}\big\ra.
\end{align}
In the following, we divide the proof into seven steps. 

{\it Step 1.} For the 2th term on LHS of \eqref{D3.2}, integrating by parts and noting $\mathbf{S}|_{x=0}=0$, one obtains
\begin{align}\label{D3.3}
&\big\la \mathbf{S}_{xx} w^2, \, \kappa \mathbf{S}_{yy}- 2\rho_s u_s \mathbf{S}_x\big\ra \nonumber\\
	&= \big\la \mathbf{S}_{x}w^2, \, \kappa \mathbf{S}_{yy}- 2\rho_s u_s \mathbf{S}_x\big\ra_{x=L} + 2\|\sqrt{\rho_su_s} \mathbf{S}_xw\|^2_{x=L}  +\big\la  \rho_s u_s w^2, \,  (|\mathbf{S}_{x}|^2)_x \big\ra \nonumber\\
	&\quad   + \big\la 2 (\rho_s u_s)_x w^2,   \,   |\mathbf{S}_{x}|^2\big\ra  - \kappa \big\la \mathbf{S}_{x} w^2, \, \mathbf{S}_{xyy} \big\ra  \nonumber\\
	&=- \big\la \mathbf{S}_{x}w^2, \, \kappa \mathbf{S}_{xyy} \big\ra + \big\la \mathbf{S}_{x} w^2, \, \kappa \mathbf{S}_{yy}- \rho_s u_s \mathbf{S}_x\big\ra_{x=L} +  \|\sqrt{\rho_su_s} \mathbf{S}_xw\|^2_{x=0}   + \big\la (\rho_s u_s)_x , \, |\mathbf{S}_{x}|^2w^2\big\ra  \nonumber\\
	&=\kappa \|\mathbf{S}_{xy}w\|^2  + \kappa \big\la \mathbf{S}_{x} w^2, \, \mathbf{S}_{xy} \big\ra_{y=0}  + \|\sqrt{\rho_su_s}\mathbf{S}_xw\|^2_{x=0} + \big\la \mathbf{S}_{x}, \, \kappa \mathbf{S}_{yy}- \rho_s u_s \mathbf{S}_x\big\ra_{x=L}  \nonumber\\
	&\quad + O(1)\|u_s \mathbf{S}_x w\|\cdot  \|\mathbf{S}_x w\| +O(1) \|\mathbf{S}_{xy}w\|\cdot \|\mathbf{S}_xw_y\| ,
\end{align}

\smallskip

{\it Step 2.}  For the 1th term on RHS of \eqref{D3.2},  we integrate by parts to obtain
\begin{align}\label{D3.11}
\eqref{D3.2}_1&=  - \big\la \mathbf{S}_{x} w^2, \,  2\rho_s v_s \mathbf{S}_{xy}\big\ra - \big\la \mathbf{S}_{x}  w^2, \,  2 (\rho_s v_s)_x \mathbf{S}_{y}\big\ra 
 +  \big\la \mathbf{S}_{x} w^2, \,  2\rho_s v_s \mathbf{S}_y\big\ra_{x=L}\nonumber\\
&\lesssim \xi\|\sqrt{u_s}\mathbf{S}_xw\|^2_{x=L} + (\xi + C_{\xi}L) \|\mathbf{S}_{xy}w\|^2 + C_{\xi} \|u_s \mathbf{S}_xw\|^2,
\end{align}
where we have used  $\mathbf{S}|_{x=0}=0$ and the following estimates
\begin{align*}
	- \big\la \mathbf{S}_{x}, \,  2\rho_s v_s \mathbf{S}_{xy} w^2\big\ra &\lesssim  \xi \|\mathbf{S}_{xy}w\|^2 + C_{\xi} \|u_s \mathbf{S}_xw\|^2,\\
	-\big\la \mathbf{S}_{x}, \,  2 (\rho_s v_s)_x \mathbf{S}_{y}  w^2\big\ra &\lesssim \xi\|\mathbf{S}_{xy} w\|^2 + C_{\xi} L^2 \|u_s \mathbf{S}_xw\|^2,\\
	\big\la \mathbf{S}_{x}, \,  2\rho_s v_s \mathbf{S}_y  w^2\big\ra_{x=L}&
	\lesssim \xi \|\sqrt{u_s}\mathbf{S}_xw\|^2_{x=L} + C_{\xi} L \|\mathbf{S}_{xy}w\|^2.
\end{align*}
 
\smallskip

{\it Step 3.}  For the 2th term on RHS of \eqref{D3.2}, it is clear that 
\begin{align}\label{D3.13-1}
	\big\la \mathbf{S}_{xx} w^2, \,  \v \big(\rho_s^{-1}p\big)_{xx} + \rho_s^{-1}p_{yy} + (2\rho_s^{-1})_y\,p_{y}\big\ra
	&\lesssim \xi \|\sqrt{\v}\mathbf{S}_{xx}w\|^2 + C_{\xi} \f1{\v}[[p]]_{2,\hat{w}}.
\end{align}
Integrating by parts in $x$, then using \eqref{2TC.3-1}, we  obtain
\begin{align}\label{D3.13-3}
	\big\la \mathbf{S}_{xx} w^2, \,  \big(\rho_s^{-1}\big)_{yy}p\big\ra
	&= \big\la \mathbf{S}_{x} w^2, \,  \big(\rho_s^{-1}\big)_{yy}p\big\ra_{x=L} + O(1)\|\mathbf{S}_xw\| \big\{\f{a}{\sqrt{\v}}\|\sqrt{\v}p_xw\| + \|\mathbf{d}_3 p w\|\big\}\nonumber\\
	&\quad  + O(1)a^{-\f12}\|\sqrt{u_s}\mathbf{S}_xw\|_{x=0}\cdot \|(\mathbf{d}_3p w, p_y)\|_{x=0}\nonumber\\
	&= \big\la \mathbf{S}_{x} w^2, \,  \big(\rho_s^{-1}\big)_{yy}p\big\ra_{x=L} + \xi \|\mathbf{S}_{xy}\|^2 + \xi \|\sqrt{u_s}\mathbf{S}_xw\|^2_{x=0} 
	 \nonumber\\
	&\quad + C_{\xi} L^{\f12} [[[\mathbf{S}]]]_{2,\hat{w}} + C_{\xi} a^2 \f{1}{\v} [[p]]_{2,\hat{w}}. 
\end{align}
Combining \eqref{D3.13-1}-\eqref{D3.13-3}, we get
\begin{align}\label{D3.14-0}
\eqref{D3.2}_2 &=\big\la \mathbf{S}_{x} w^2, \,  \kappa\big(\rho_s^{-1}\big)_{yy}p\big\ra_{x=L} + \xi \|(\sqrt{\v}\mathbf{S}_{xx}, \mathbf{S}_{xy})\hat{w}\|^2 + \xi \|\sqrt{u_s}\mathbf{S}w\|^2_{x=0} \nonumber\\
&\quad + C_{\xi} L^{\f12} [[[\mathbf{S}]]]_{2,\hat{w}} + C_{\xi} \f1{\v}[[p]]_{2,\hat{w}}.
\end{align}

\smallskip

{\it Step 4.}  For the 3th term on RHS of \eqref{D3.2}, it is noted that 
\begin{align}\label{D1TC.36}
\big\la \v \big(p_s^{-1}T_{sy}q\big)_{xx}, \, \mathbf{S}_{xx} w^2 \big\ra 
&\lesssim \|\sqrt{\v}\mathbf{S}_{xx}w\|\cdot  \|(\sqrt{\v}q_{xx}, \sqrt{\v}q_x)\| 
\lesssim \xi \|\sqrt{\v}\mathbf{S}_{xx}w\|^2 + C_{\xi} L^{\alpha} [\phi]_{3,1}.
\end{align}
Integrating by parts in $x$ and then $y$, and noting $q|_{x=0}=0$,    one obtains that 
\begin{align}\label{D1TC.37}
&\big\la \big(p_s^{-1}T_{sy}q\big)_{yy}, \, \mathbf{S}_{xx}w^2\big\ra \nonumber\\
&= \big\la p_s^{-1}T_{sy}q_{xy}, \, \mathbf{S}_{xy} w^2\big\ra + \big\la \big(p_s^{-1}T_{sy} w^2\big)_y q_{xy}, \, \mathbf{S}_{x} \big\ra- \big\la (p_s^{-1}T_{sy})_x q_{yy}, \, \mathbf{S}_{x} w^2\big\ra\nonumber\\
&\quad  -  \big\la 2\big(\big(p_s^{-1}T_{sy}\big)_{y}q_{y}\big)_x, \, \mathbf{S}_{x} w^2\big\ra  -  \big\la \big(\big(p_s^{-1}T_{sy}\big)_{yy} q\big)_x, \, \mathbf{S}_{x} w^2\big\ra  + \big\la p_s^{-1}T_{sy}q_{xy}, \, \mathbf{S}_{x} w^2\big\ra_{y=0} \nonumber\\
&\quad + \big\la p_s^{-1}T_{sy}q_{yy}, \, \mathbf{S}_{x} w^2\big\ra_{x=L} + \big\la 2\big(p_s^{-1}T_{sy}\big)_{y}q_y, \, \mathbf{S}_{x} w^2\big\ra_{x=L}  +  \big\la \big(p_s^{-1}T_{sy}\big)_{yy} q, \, \mathbf{S}_{x}w^2\big\ra_{x=L} \nonumber\\
&= O(a) \|\mathbf{S}_{xy}w\|\cdot \|q_{xy}\| + O(a)\|\mathbf{S}_{x}\hat{w}\| \cdot \|(q_{xy}, q_{yy})\|  + \big\la p_s^{-1}T_{sy}q_{xy}, \, \mathbf{S}_{x} w^2\big\ra_{y=0} \nonumber\\
&\quad  +  \big\la \big(p_s^{-1}T_{sy}\big)_{yy} q, \, \mathbf{S}_{x} w^2\big\ra_{x=L} + \big\la p_s^{-1}T_{sy}q_{yy}, \, \mathbf{S}_{x} w^2\big\ra_{x=L} + \big\la 2\big(p_s^{-1}T_{sy}\big)_{y}q_y, \, \mathbf{S}_{x} w^2\big\ra_{x=L}\nonumber\\
&= \big\la p_s^{-1}T_{sy}q_{yy}, \, \mathbf{S}_{x} w^2\big\ra_{x=L} + \big\la 2\big(p_s^{-1}T_{sy}\big)_{y}q_y, \, \mathbf{S}_{x} w^2\big\ra_{x=L} +  \xi \|\mathbf{S}_{xy} w\|^2 + + C_{\xi} L^{\f12} [[[\mathbf{S}]]]_{2,\hat{w}} \nonumber\\
&\quad  + L^{\f12-} \|\sqrt{u_s}\mathbf{S}_x w\|^2_{x=L}   + C_{\xi} L^{\alpha} [\phi]_{3,1} + C_{\xi} \|T_{sy}q_{xy}\|_{y=0}^2,
\end{align}
where we have used the following two estimates
\begin{align}\label{D1TC.001}
\begin{split}
	\big\la p_s^{-1}T_{sy}q_{xy}, \, \mathbf{S}_{x}\big\ra_{y=0}
&\lesssim \xi \|\mathbf{S}_{xy}\|^2 + C_{\xi} \|T_{sy}q_{xy}\|_{y=0}^2  + \f1{a^2}\|u_s\mathbf{S}_x\|^2,\\
\big\la \big(p_s^{-1}T_{sy}\big)_{yy} q, \, \mathbf{S}_{x} w^2 \big\ra_{x=L}&\lesssim \xi \|\sqrt{u_s} \mathbf{S}_x\hat{w}\|^2_{x=L} + C_{\xi} L^{\f12-} \|q_{xy}\|^2.
\end{split}
\end{align}
We remark that the trace theorem and \eqref{B.14} (or \eqref{2TC.3-1}) are used in $\eqref{D1TC.001}_1$.

Noting from $\eqref{1.7-3}_2$, i.e., $|T_{sy}(x,0)|\lesssim |u_{sy}(x,0)|^2$, 
one gets
\begin{align}
	\|T_{sy}q_{xy}\|_{y=0}^2\lesssim a^3 \|\sqrt{u_{sy}}q_{xy}\|^2_{y=0} \lesssim a^3 [\phi]_{3,1},
\end{align}
which, together with \eqref{D1TC.36}-\eqref{D1TC.37}, yields that  
\begin{align}\label{D3.19-0}
	\eqref{D3.2}_3&= - \big\la \kappa p_s^{-1}T_{sy}q_{yy}, \, \mathbf{S}_{x} w^2\big\ra_{x=L} - \big\la 2\kappa \big(p_s^{-1}T_{sy}\big)_{y}q_y, \, \mathbf{S}_{x} w^2\big\ra_{x=L} + \xi \|(\mathbf{S}_{xy}, \sqrt{\v}\mathbf{S}_{xx})w\|^2  \nonumber\\
	&\quad + C_{\xi} L^{\f12} [[[\mathbf{S}]]]_{2,\hat{w}} + L^{\f12-} \|\sqrt{u_s}\mathbf{S}_xw\|^2_{x=L} 
	+ C_{\xi} a^3 [\phi]_{3,1}.
\end{align}

\smallskip

{\it Step 5.} We consider 4th term on RHS of \eqref{D3.2}. Recall  $\mathfrak{J}$   in \eqref{2.84}.
Integrating by parts in $x$, one has
\begin{align}\label{D1TC.39}
\big\la  2\rho_sv_s \big(p_s^{-1}T_{sy}q\big)_y -  2\big[\rho_su_s \big(p_s^{-1}T_{sy}\big)_x-  T_{sy} \bar{u}_{sx}\big]q,\, \mathbf{S}_{xx} w^2\big\ra 
&\lesssim L^{\alpha} [[[\mathbf{S}]]]_{2,\hat{w}} + L^{\alpha} [\phi]_{3,\hat{w}}
\end{align}
and
\begin{align}\label{D1TC.40}
	&\big\la \rho_s (U_s\cdot \nabla) (\rho_s^{-1})\, p, \, \mathbf{S}_{xx} w^2\big\ra\nonumber\\
	&\lesssim  \|u_s\mathbf{S}_{x}w\|\cdot \|\mathbf{d}_2(p_x, p)w\| + \|\sqrt{u_s}\mathbf{S}_xw\|_{x=0} \cdot \|\mathbf{d}_2pw\|_{x=0}  + \|\sqrt{u_s}\mathbf{S}_xw\|_{x=L} \cdot \|\mathbf{d}_2pw\|_{x=L} \nonumber\\
	&\lesssim \xi \|\sqrt{u_s}\mathbf{S}_{x}w\|^2_{x=0,L}  + C_{\xi} L^{\f12} [[[\mathbf{S}]]]_{2,\hat{w}} + C_{\xi} \f{1}{\v} [[p]]_{2,1} + C_{\xi} L^2\f{1}{\v^2} [[p]]_{2,1} \mathbf{1}_{\{w=w_0\}}.
\end{align}

Noting $p_{sy}\cong \sqrt{\v}$, we integrate by parts in $x$ to obtain
\begin{align}\label{D1TC.41}
	&\big\la \rho_s^{-1}p_{sy}  \phi_x , \mathbf{S}_{xx} w^2 \big\ra
	\nonumber\\
	&=-\big\la \rho_s^{-1}p_{sy}\bar{u}_s q_{xx} , \mathbf{S}_{x} w^2\big\ra - \big\la (\rho_s^{-1}p_{sy}\bar{u}_s)_x q_x , \mathbf{S}_{x} w^2\big\ra + \big\la \rho_s^{-1}p_{sy}\bar{u}_s q_x , \mathbf{S}_{x} w^2\big\ra\big|^{x=L}_{x=0} \nonumber\\
	&\quad  - \big\la \rho_s^{-1}p_{sy} \bar{u}_{sx} q_x , \mathbf{S}_{x} w^2\big\ra - \big\la (\rho_s^{-1}p_{sy} \bar{u}_{sx})_x q , \mathbf{S}_{x} w^2\big\ra +  \big\la \rho_s^{-1}p_{sy} \bar{u}_{sx} q , \mathbf{S}_{x} w^2 \big\ra_{x=L} \nonumber\\
	&\lesssim \|u_s\mathbf{S}_{x}w\| \big\{ \|\sqrt{\v}q_{xx}w\| + \|\sqrt{\v} q_xw\|\big\} +  \|\sqrt{u_s}\mathbf{S}_{x}w\|_{x=0} \cdot \|\sqrt{\v}q_xw\|_{x=0} \nonumber\\
	&\quad +   \|\sqrt{u_s}\mathbf{S}_{x}w\|_{x=L}  \|\sqrt{\v}q_x w\|_{x=L}  \nonumber\\
	&\lesssim \xi \|\sqrt{u_s}\mathbf{S}_{x}w\|^2_{x=0, L} +  C_{\xi} L^{\f12} [[[\mathbf{S}]]]_{2,\hat{w}}+ C_{\xi} L^{\alpha-} [\phi]_{3,w}.
\end{align}
It is also clear to see that 
\begin{align}\label{D1TC.42}
	\big\la \mathfrak{M} +\mathfrak{N}_3, \, \mathbf{S}_{xx}w^2\big\ra 
	& = \big\la \mathfrak{M}_0, \, \mathbf{S}_{xx} w^2 \big\ra + \big\la \v (\mathfrak{M}_1+\mathfrak{M}_2) , \, \mathbf{S}_{xx} w^2\big\ra + \big\la \mathfrak{N}_3, \, \mathbf{S}_{xx} w^2 \big\ra\nonumber\\
	&= - \big\la \pa_x\mathfrak{M}_0, \, \mathbf{S}_{x}w^2\big\ra +\big\la \mathfrak{M}_0, \, \mathbf{S}_{x}w^2\big\ra_{x=L} - \big\la \mathfrak{M}_0, \, \mathbf{S}_{x}w^2\big\ra_{x=0}\nonumber\\
	&\quad  + \xi \|\sqrt{\v}\mathbf{S}_{xx}w\|^2   + C_{\xi} \v\|(\mathfrak{M}_1,\mathfrak{M}_2)w\|^2 + C_{\xi}  \f{1}{\v}\|\mathfrak{N}_3 w\|^2.
\end{align}

It follows from  \eqref{1C.43-1}-\eqref{1TC.44} that 
\begin{align}
	|\big\la \pa_x\mathfrak{M}_0, \, \mathbf{S}_{x}w^2\big\ra| 
	&\lesssim  \xi \|\mathbf{S}_{xy}\|^2 +  C_{\xi} L^{\f12} [[[\mathbf{S}]]]_{2,\hat{w}} + C_{\xi} a^2 \|(\phi,p,S) \hat{w}\|^2_{\mathbf{X}^{\v}}\nonumber\\
	&\quad + C_{\xi} La\Big\{\f{1}{\v^2} [[p]]_{2,1} + \f{1}{\v} [[[S]]]_{2,1}\Big\} \mathbf{1}_{\{w=w_0\}},\label{D1TC.45}\\
	\big|\big\la \mathfrak{M}_0, \, \mathbf{S}_{x} w^2 \big\ra_{x=0}\big| 
	&\lesssim \|\sqrt{u_s}\mathbf{S}_xw\|_{x=0} \|p_y\|_{x=0} + a\int_0^\infty \mathbf{d}_2 |\mathbf{S}_x p|w^2 dy \big|_{x=0} \nonumber\\
	& \lesssim \xi \|\sqrt{u_s}\mathbf{S}_xw\|^2_{x=0} + C_{\xi} \v^{-\f12} [[p]]_{2,1} + C_{\xi}  \v^{-\f32-} [[p]]_{2,1} \mathbf{1}_{\{w=w_0\}},\label{D1TC.46}
\end{align} 
where we have used the compatibility condition $\dis S|_{x=0}=0$ and  $p(0,0)=0$.

Recall again $\mathfrak{J}$  in \eqref{2.84}. Then combining \eqref{D1TC.39}-\eqref{D1TC.46}  and \eqref{1TC.48}, one obtains that 
\begin{align}\label{D3.31}
	\eqref{D3.2}_4&=-\big\la \mathfrak{M}_0, \, \mathbf{S}_{x}w^2\big\ra_{x=L} +  \xi  \|\sqrt{u_s}\mathbf{S}_{x}w\|^2_{x=0,L}  + \xi \|(\sqrt{\v}\mathbf{S}_{xx}, \mathbf{S}_{xy})\hat{w}\|^2   + C_{\xi} L^{\f12} [[[\mathbf{S}]]]_{2,\hat{w}}\nonumber\\
	&\quad   + C_{\xi} \f{1}{\v} [[p]]_{2,1} + C_{\xi} a^2 \|(\phi,p,S) \hat{w}\|^2_{\mathbf{X}^{\v}}+ C_{\xi} L\Big\{\f{1}{\v^2} [[p]]_{2,1} + \f{1}{\v} [[[S]]]_{2,1}\Big\} \mathbf{1}_{\{w=w_0\}} \nonumber\\
	&\quad   + C_{\xi}  \f{1}{\v}\|\mathfrak{N}_3 w\|^2.
\end{align}

\smallskip

{\it Step 6.} Substituting \eqref{D3.3}, \eqref{D3.11}, \eqref{D3.14-0},  \eqref{D3.19-0} and \eqref{D3.31} into \eqref{D3.2}, then using $\eqref{2.111-6D}_2$,  one gets
\begin{align}\label{D3.5}
	&\f12 \kappa \|(\sqrt{\v}\mathbf{S}_{xx}, \mathbf{S}_{xy})w\|^2 + \|\sqrt{\rho_s u_s} \mathbf{S}_x w\|^2_{x=L} + \|\sqrt{\rho_s u_s} \mathbf{S}_xw\|^2_{x=0} +\kappa \big\la \mathbf{S}_{x} w^2, \, \mathbf{S}_{xy} \big\ra_{y=0}  \nonumber\\
	&\lesssim \|\mathbf{S}_xw_y\|^2 +  \xi \|(\sqrt{\v}\mathbf{S}_{xx}, \mathbf{S}_{xy})\|^2 + L^{\f14-} [[[\mathbf{S}]]]_{2,\hat{w}}    
	 	+  \f1{\v}[[p]]_{2,\hat{w}} + C_{\xi} a^2 \|(\phi,p,S) \hat{w}\|^2_{\mathbf{X}^{\v}}\nonumber\\
	&\quad   + C_{\xi} La\Big\{\f{1}{\v^2} [[p]]_{2,1} + \f{1}{\v} [[[S]]]_{2,1}\Big\} \mathbf{1}_{\{w=w_0\}} + C_{\xi}  \f{1}{\v}\|\mathfrak{N}_3 w\|^2 .
\end{align}

For the boundary term $\big\la \mathbf{S}_{x}, \, \mathbf{S}_{xy} \big\ra_{y=0}$  on LHS of \eqref{D3.5}, by noting $\eqref{2.111-6D}_2$ and $\eqref{D3.1-1}_1$, one gets 
\begin{align}\label{D3.6}
	\big\la \mathbf{S}_{x}, \, \mathbf{S}_{xy} \big\ra_{y=0}
	&= \big\la (2\rho_s)^{-1} g_1, \, \mathbf{S}_{xy} \big\ra_{y=0} - \mu\big\la (2\rho_s^2)^{-1} \phi_{yyy}, \, \mathbf{S}_{xy} \big\ra_{y=0}   - \big\la \big((2\rho_s)^{-1}\big)_x \mathfrak{P}, \, \mathbf{S}_{xy} \big\ra_{y=0}  \nonumber\\
	&\quad - \big\la (2\rho_s)^{-1} \big[ \mathfrak{h}_1 + \mathfrak{N}_{11}\big], \, \mathbf{S}_{xy} \big\ra_{y=0} +  \mu \v^{N_0} \big\la \rho^{(\v)} (\rho^{\v}\rho_s)^{-1} \phi_{yyy}, \, \mathbf{S}_{xy} \big\ra_{y=0} .
\end{align}

{\it Step 6.1.} For 1th term on RHS of \eqref{D3.6},
it is noted that
\begin{align}\label{D3.24}
	g_1\big|_{y=0} 
	&= \mu p_s^{-1}  u_{sy}\, p_y -  2\mu u_{sy} T_s^{-1}\, S_y  - 2\mu T_s^{-1} u_{sy}\,  \big[\big((2\rho_s)^{-1}\big)_yp + p_s^{-1}T_{sy} q_y\big]\nonumber\\
	&\quad + \tilde{d}_{13} \, p  + \hat{d}_{13}\, T + 2\mu\rho_s^{-2}\rho_{sy} \phi_{yy} - \mathfrak{R}_{s1}.
\end{align}
Noting $\mathbf{S}|_{x=0}=0$ and integrating by parts, one has that
\begin{align}
&\big\la  \big\{\mu (2\rho_sp_s)^{-1}  u_{sy}\, p_y\big\}, \, \mathbf{S}_{xy} \big\ra_{y=0} 
\nonumber\\
& =\big\la \big\{\cdots\big\}_x, \, \mathbf{S}_{yy} \big\ra - \big\la \big\{\cdots\big\}_y, \, \mathbf{S}_{xy} \big\ra  - \big\la \big\{\cdots\big\}, \, \mathbf{S}_{yy} \big\ra_{x=L} \nonumber\\
&\lesssim a\|\mathbf{S}_{yy}\|\cdot \|(p_y, p_{xy})\|  + a \|\mathbf{S}_{xy}\| \|(p_y, p_{yy})\| + a \|\mathbf{S}_{yy}\|_{x=L}\cdot \|p_y\|_{x=L}\nonumber\\
&\lesssim L^{\f14} [[[\mathbf{S}]]]_{2,\hat{w}} + \v^{\f12-} \|\mathbf{S}_{yy}\|^2_{x=L} + L^{\f14} \f1{\v}[[p]]_{2,1},
\end{align}
and
\begin{align}
	&-\big\la \big\{\mu p_s^{-1} u_{sy} \big[\big((2\rho_s)^{-1}\big)_yp + p_s^{-1}T_{sy} q_y\big]\big\} , \, \mathbf{S}_{xy} \big\ra_{y=0}\nonumber\\
	&=\big\la \big\{\cdots \big\}_x , \, \mathbf{S}_{yy} \big\ra- \big\la \big\{\cdots \big\}_y , \, \mathbf{S}_{xy} \big\ra  -\big\la \big\{\cdots \big\}, \, \mathbf{S}_{yy} \big\ra_{x=L} \nonumber\\
	&\lesssim L^{\f14} [[[\mathbf{S}]]]_{2,\hat{w}}   + L^{\f14} \|\mathbf{S}_{yy}\|^2_{x=L} + L^{\f14}\f1{\v}[[p]]_{2,1} + L^{\f14}[\phi]_{3,1}.
\end{align}

Also it holds that 
\begin{align}
	&\big\la \big((2\rho_s)^{-1}\big\{\big[\tilde{d}_{13}   + \, \hat{d}_{13}(2\rho_s)^{-1}\big] p+\hat{d}_{13}p_s^{-1}T_{sy}q  + 2\mu \rho_s^{-2}\rho_{sy} \phi_{yy} - \mathfrak{R}_{s1}\big\}\big), \, \mathbf{S}_{xy} \big\ra_{y=0}\nonumber\\
	&=\big\la \big\{\cdots\big\}_x, \, \mathbf{S}_{yy} \big\ra  - \big\la \big\{\cdots\big\}_y, \, \mathbf{S}_{xy} \big\ra - \big\la \big\{\cdots\big\} , \, \mathbf{S}_{yy} \big\ra_{x=L}\nonumber\\
	&\lesssim L^{\f18} [[[\mathbf{S}]]]_{2,\hat{w}}   + L^{\f12-} \|\mathbf{S}_{yy}\|^2_{x=L} + L^{\f18}\f1{\v}[[p]]_{2,1} + L^{\f18}[\phi]_{3,1} + \v^{N_0},
\end{align}

\begin{align}
	& \big\la (2\rho_s)^{-1} \hat{d}_{13}S, \, \mathbf{S}_{xy} \big\ra_{y=0}
	\nonumber\\
	&=\big\la \big( (2\rho_s)^{-1} \hat{d}_{13}S\big)_x, \, \mathbf{S}_{yy} \big\ra  - \big\la \big( (2\rho_s)^{-1} \hat{d}_{13}S\big)_y, \, \mathbf{S}_{xy} \big\ra + \big\la (2\rho_s)^{-1} \hat{d}_{13}S, \, \mathbf{S}_{yy} \big\ra_{x=L}\nonumber\\
	&\lesssim a \|\mathbf{S}_{yy}\|\cdot \|S_x\| + \sqrt{L} \|\mathbf{S}_{xy}\|^2  + \sqrt{L} \|\mathbf{S}_{yy}\|^2_{x=L} + \sqrt{L} \|S_{xy}\|^2  + \sqrt{L} \|S_x\|^2\nonumber\\
	&\lesssim L^{\f14-} [[[\mathbf{S}]]]_{2,1} +  L^{\f14-} [[[S]]]_{2,1}  + \sqrt{L} \|\mathbf{S}_{yy}\|^2_{x=L},
\end{align}
and
\begin{align}\label{D3.28-1}
	&-\big\la \mu p_s^{-1} u_{sy}\, S_y , \, \mathbf{S}_{xy} \big\ra_{y=0} 
	\nonumber\\
	&=\big\la \mu p_s^{-1} u_{sy}\, S_y , \, \mathbf{S}_{yy} \big\ra_{x=L} - \big\la \big(\mu p_s^{-1} u_{sy}\, S_y\big)_x , \, \mathbf{S}_{yy} \big\ra +  \big\la \big(\mu p_s^{-1} u_{sy}\, S_y \big)_y, \, \mathbf{S}_{xy} \big\ra \nonumber\\
	&\lesssim  L^{\f14-} [[[\mathbf{S}]]]_{2,1} +  L^{\f14-} [[[S]]]_{2,1}  + \sqrt{L} \|\mathbf{S}_{yy}\|^2_{x=L}.
\end{align}


\smallskip


Combining \eqref{D3.24}-\eqref{D3.28-1}, then one obtains that 
\begin{align}\label{D3.33}
\eqref{D3.6}_1 &\lesssim  L^{\f18} [[[\mathbf{S}]]]_{2,\hat{w}} + L^{\f14} \|\mathbf{S}_{yy}\|^2_{x=L} + L^{\f18} \|(\phi, p, S) \hat{w}\|^2_{\mathbf{X}^{\v}} + \v^{N_0}.
\end{align}

\smallskip

{\it Step 6.2.} For the 2th term on RHS of \eqref{D3.6}, noting $\mathbf{S}|_{x=0}=0$, one integrates by parts to have
\begin{align}\label{D3.7}
	- \big\la \rho_s^{-2} \phi_{yyy}, \, \mathbf{S}_{xy} \big\ra_{y=0}
	&= \big\la \rho_s^{-2} \phi_{yyyy}, \, \mathbf{S}_{xy} \big\ra  - \big\la \rho_s^{-2} \phi_{xyyy}, \, \mathbf{S}_{yy} \big\ra 
	 - \big\la \big(\rho_s^{-2}\big)_x \phi_{yyy}, \, \mathbf{S}_{yy} \big\ra \nonumber\\      
	&\quad +  \big\la \big(\rho_s^{-2}\big)_y \phi_{yyy}, \, \mathbf{S}_{xy} \big\ra  +  \big\la \rho_s^{-2} \phi_{yyy}, \, \mathbf{S}_{yy} \big\ra_{x=L}  .
\end{align}

It is clear that 
\begin{align}\label{D3.35}
\eqref{D3.7}_{1,3,4}	&\lesssim \xi \|(\mathbf{S}_{xy}, \mathbf{S}_{yy})\|^2 + C_{\xi} L^{\f18-} [\phi]_{3,1}  + C_{\xi} \|||\phi||\|^2_1.
\end{align}
For the boundary term on RHS of \eqref{D3.7}, it follows from $\eqref{2.111-6D}_2$ that
\begin{align}\label{D3.9}
\eqref{D3.7}_5	&\lesssim \xi \|\mathbf{S}_{yy}\|^2_{x=L}  +  C_{\xi} a [\phi]_{3,1}  + C_{\xi}\|||\phi||\|^2_1,
\end{align}
where we have used  the following estimates
\begin{align}\label{D3.9-0}
	\|\phi_{yyy}\|^2_{x=L}
	&\lesssim \iint |\phi_{yyyy} \phi_{xyy}| dydx + \int_0^L  |\phi_{yyy}\phi_{xyy}|  dx\Big|_{y=0} \nonumber\\
	&  \lesssim \sqrt{a} \|||\phi||\|_1 \cdot [\phi]_{3,1}^{\f12}  + L^{\f18-}  [\phi]_{3,1}.
\end{align}

Applying the equation \eqref{2.111-6D}, one has
\begin{align}\label{D3.9-1}
\eqref{D3.7}_2	&=   \big\la \rho_s^{-2} \phi_{xyyy}, \, \v \mathbf{S}_{xx} \big\ra - \big\la \rho_s^{-1} \phi_{xyyy}, \, 2 u_s \mathbf{S}_x \big\ra  - \big\la \rho_s^{-1}  \phi_{xyyy}, \, 2 v_s \mathbf{S}_y\big\ra\nonumber\\
&\quad + \f12 \big\la \rho_s^{-2} \phi_{xyyy}, \,   \Delta_{\v}\big(\rho_s^{-1}p\big)\big\ra + \big\la \rho_s^{-2} \phi_{xyyy}, \,  \Delta_{\v}\big(p_s^{-1}T_{sy}q\big)\big\ra +  \big\la \rho_s^{-2}  \phi_{xyyy}, \,  \mathfrak{J}\big\ra.
\end{align}
A direct calculation shows that 
\begin{align}\label{D3.10}
\eqref{D3.9-1}_1&\lesssim \xi \|\sqrt{\v}\mathbf{S}_{xx}\|^2 + C_{\xi} \|||\phi||\|^2_1, 
\end{align}
\begin{align}\label{D3.10-1}
\eqref{D3.9-1}_2& = \big\la \rho_s^{-1} \phi_{xyy}, \, 2 u_s \mathbf{S}_{xy} \big\ra + \big\la  \phi_{xyy}, \, (\rho_s^{-1}u_s)_y \mathbf{S}_x \big\ra\nonumber\\
	&\lesssim \xi \|\mathbf{S}_{xy}\|^2 + C_{\xi} a^2 [\phi]_{3,1} +  a^{-2}\|u_s \mathbf{S}_x\|^2,
\end{align}
\begin{align}
\eqref{D3.9-1}_3&=  \big\la \rho_s^{-1} \phi_{xyy}, \, 2 v_s \mathbf{S}_{yy}\big\ra  - \big\la  \phi_{xyy}, \, (\rho_s^{-1} v_s)_y \mathbf{S}_y\big\ra\nonumber\\
	&\lesssim L^{\f14-} [[[\mathbf{S}]]]_{2,1} + L^{\f14-}a^2 [\phi]_{3,1},
\end{align}
where we have used \eqref{B.14} (or \eqref{2TC.3-1}) in the proof of
 \eqref{D3.10-1}.

Integrating by parts in $y$, and noting , we get 
\begin{align}\label{D3.12}
\eqref{D3.9-1}_4&=\f12 \big\la  \rho_s^{-3} \phi_{xyyy}, \,  \Delta_{\v}p \big\ra  +   \big\la  \rho_s^{-2} \phi_{xyyy}, \,  \nabla_{\v}(\rho_s^{-1})\cdot \nabla_{\v}p \big\ra   \nonumber\\
	&\quad - \f12 \big\la  \phi_{xyy}, \,  \big\{\rho_s^{-2}\Delta_{\v}(\rho_s^{-1}) p\big\}_y \big\ra  - \f12 \big\la  \phi_{xyy}, \,  \rho_s^{-2}\Delta_{\v}(\rho_s^{-1})\, p \big\ra_{y=0} \nonumber\\
	&\lesssim  \|||\phi||\|^2_1 + \sqrt{L}[\phi]_{3,1} +  \f{1}{\v} [[p]]_{2,1}.
\end{align}
where we have used the fact $\|\phi_{xyy}\|^2_{y=0}\lesssim a \|\sqrt{u_{sy}} q_{xy}\|^2_{y=0} $.

For the 5th term on RHS of \eqref{D3.9-1}, we note
\begin{align}\label{D3.14}
\eqref{D3.9-1}_5 &= \big\la  \rho_s^{-2} \phi_{xyyy}, \, \nabla_{\v}\big(p_s^{-1}T_{sy}\big)\cdot \nabla_{\v}q \big\ra + \big\la  \rho_s^{-2}  \phi_{xyyy}, \,  q \Delta_{\v}\big(p_s^{-1}T_{sy}\big) \big\ra \nonumber\\
	&\quad + \big\la  \phi_{xyyy}, \,   \rho_s^{-2} p_s^{-1}T_{sy}\Delta_{\v}q \big\ra.
\end{align}
Noting $q|_{y=0}=0$ and $\phi|_{y=0}=0$, we integrate by parts in $y$ to obtain
\begin{align}
\eqref{D3.14}_{1,2}&= -  \big\la   \phi_{xyy}, \, \rho_s^{-2} \big(p_s^{-1}T_{sy}\big)_y  q_y \big\ra_{y=0} -\big\la  \phi_{xyy}, \, \big\{ \rho_s^{-2}  \nabla_{\v}\big(p_s^{-1}T_{sy}\big)\cdot \nabla_{\v}q\big\}_y \big\ra \nonumber\\
&\quad - \big\la   \phi_{xyy}, \,  \big\{q \rho_s^{-2} \Delta_{\v}\big(p_s^{-1}T_{sy}\big)\big\}_y \big\ra    \lesssim   L^{\f18-} [\phi]_{3,1} .
\end{align}
For the 3th term  on RHS of \eqref{D3.14}, we have
\begin{align}\label{D3.45}
\eqref{D3.14}_3&=\big\la \pa_{xy}\big\{\bar{u}_sq_{yy}\big\}, \,  \rho_s^{-2} p_s^{-1} T_{sy}\Delta_{\v}q \big\ra + \big\la  \bar{u}_{sy}q_{xyy}  , \,  \rho_s^{-2} p_s^{-1}T_{sy}\Delta_{\v}q \big\ra + L^{\f18-} [\phi]_{3,1}\nonumber\\
&\lesssim L^{\f14\a-} [\phi]_{3,1},
\end{align}
where we have used the following estimates
\begin{align}
	\big\la \pa_{xy}\big\{\bar{u}_sq_{yy}\big\}, \,  \rho_s^{-2} p_s^{-1} T_{sy}\Delta_{\v}q \big\ra
	&= - \big\la \pa_{x}\big\{\bar{u}_sq_{yy}\big\}, \,  \big\{\rho_s^{-2}p_s^{-1}T_{sy}\Delta_{\v}q\big\}_y \big\ra  \lesssim L^{\f18-} [\phi]_{3,1},
\end{align}
and
\begin{align}\label{D3.43}
&\big\la  \bar{u}_{sy}q_{xyy}  , \,  \rho_s^{-2}p_s^{-1}T_{sy}\Delta_{\v}q \big\ra\nonumber\\
&=  \f12 \big\la  \rho_s^{-2}p_s^{-1}T_{sy} u_{sy}  , \,   |q_{yy}|^2 \big\ra_{x=L}  -  \v  \big\la  \rho_s^{-2}p_s^{-1}T_{sy} u_{sy} q_{xy} , \,  q_{xxy} \big\ra   + L^{\f18-} [\phi]_{3,1}\nonumber\\
&=\f12 \big\la  \rho_s^{-2}p_s^{-1}T_{sy} u_{sy}  , \,   |q_{yy}|^2 \big\ra_{x=L} + \f12 \v  \big\la  \rho_s^{-2}p_s^{-1}T_{sy} u_{sy} , \,  |q_{xy}|^2 \big\ra_{x=0}  + L^{\f18-} [\phi]_{3,1}\nonumber\\
&\lesssim L^{\f14\a-} [\phi]_{3,1}.
\end{align}
We remark that Lemma \ref{lemAD.1} is used  in \eqref{D3.43}.

Combining \eqref{D3.14}-\eqref{D3.45}, we get
\begin{align}\label{D3.18}
\eqref{D3.9-1}_5	&\lesssim L^{\f14\a-} [\phi]_{3,1} + L^{\f14\a-} [\phi]_{4,1}.
\end{align}

\smallskip

We shall estimate 6th term on RHS of \eqref{D3.9-1}. Recall $\mathfrak{J}$ in \eqref{2.84}.  Integrating by parts in $y$, then we have 
\begin{align}\label{D3.57}
	&\Big\la \rho_s^{-2} \phi_{xyyy},\,  2\rho_sv_s \big(\f{1}{p_s}T_{sy}q\big)_y + 2\big[\rho_su_s \big(\f{1}{p_s}T_{sy}\big)_x-  T_{sy} \bar{u}_{sx}\big]q + \f{1}{\rho_s}p_{sy}\phi_x + \rho_s (U_s\cdot \nabla) (\frac{1}{\rho_s})\, p\Big\ra \nonumber\\
	&\lesssim   L^{\f12} [\phi]_{3,1} + L^{\f12}\f{1}{\v} [[p]]_{2,1},
\end{align}
Integrating by parts in $x$ and using \eqref{3.59-1}-\eqref{1TC.44}, one gets 
\begin{align}\label{D3.58}
\big\la \rho_s^{-2}  \phi_{xyyy}, \,  \mathfrak{M}_0\big\ra
&=- \big\la \rho_s^{-2} \phi_{yyy}, \,  \pa_x\mathfrak{M}_0\big\ra - \big\la \big(\rho_s^{-2}\big)_x  \phi_{yyy}, \,  \mathfrak{M}_0\big\ra + \big\la \rho_s^{-2} \phi_{yyy}, \,  \mathfrak{M}_0\big\ra_{x=L}\nonumber\\
&\lesssim a^2 \f{1}{\v} [[p]]_{2,1} + a^2 [[[S]]]_{2,1}  + a [\phi]_{3,1} + \|||\phi||\|^2_1.
\end{align}
We have from \eqref{D3.57}-\eqref{D3.58} and \eqref{1TC.48} that 
\begin{align}\label{D3.59}
\eqref{D3.9-1}_6 &\lesssim a^2 \f{1}{\v} [[p]]_{2,1} + a^2 [[[S]]]_{2,1}  + a [\phi]_{3,1} + \|||\phi||\|^2_1 + \f{1}{\v} \|\mathfrak{N}_{3}\|^2.
\end{align}

Substituting \eqref{D3.10}-\eqref{D3.12}, \eqref{D3.18} and \eqref{D3.59} into \eqref{D3.9-1}, one obtains
\begin{align}\label{D3.19}
\eqref{D3.7}_2	&\lesssim \xi \|(\sqrt{\v}\mathbf{S}_{xx}, \mathbf{S}_{xy})\|^2  +  L^{\f14-} [[[\mathbf{S}]]]_{2,1} +  C_{\xi} \f{1}{\v} [[p]]_{2,1}  + C_{\xi} \|||\phi||\|^2_1  + C_{\xi} a [\phi]_{3,1}.
\end{align}

Then it follows from \eqref{D3.7}-\eqref{D3.9} and \eqref{D3.19} that 
\begin{align}\label{D3.50}
\eqref{D3.6}_2	&\lesssim  \xi \|\sqrt{u_s}\mathbf{S}_x\|^2_{x=L}  + \xi \|(\sqrt{\v}\mathbf{S}_{xx}, \mathbf{S}_{xy})\|^2  + \xi \|\mathbf{S}_{yy}\|^2_{x=L}  +  L^{\f14-} [[[\mathbf{S}]]]_{2,1} \nonumber\\
&\quad + C_{\xi} a [\phi]_{3,1}  +  C_{\xi} \f{1}{\v} [[p]]_{2,1}  + C_{\xi} \|||\phi||\|^2_1.
\end{align}

\smallskip

{\it Step 6.3.} For the 3th term on RHS of \eqref{D3.6}, we note from $\eqref{D3.1-1}_1$ that 
\begin{align*}
\|\mathfrak{P}\chi\|^2 + \|\mathfrak{P}_x\chi\|^2&\lesssim \int_0^L \Big|\Big(\frac{\mu}{\rho_s} \phi_{yyy} - g_1 + \mathfrak{h}_1 + \mathfrak{N}_{11} - \mu \v^{N_0} \frac{\rho }{\rho^{\v}\rho_s}  \phi_{yyy}\Big)(x,0)\Big|^2 dx\nonumber\\
&\lesssim 
L^{\f14-} \|(\phi, p, S)\|^2_{\mathbf{X}^{\v}}+ \|(\mathfrak{h}_1, \mathfrak{N}_{11})\|^2_{H^1} + \v^{N_0},
\end{align*}
which, together with $\mathfrak{P}(0)=0$, yields that 
\begin{align}
\eqref{D3.6}_3&= \f12 \big\la \big[(\rho_s^{-1})_x \chi\big]_y \mathfrak{P} , \, \mathbf{S}_{xy} \big\ra  - \f12 \big\la \big[(\rho_s^{-1})_x \mathfrak{P}\big]_x  \chi , \, \mathbf{S}_{yy} \big\ra  + \big\la (\f{1}{2\rho_s})_x \mathfrak{P} \chi, \, \mathbf{S}_{yy} \big\ra_{x=L}  
	\nonumber\\
	&\lesssim \xi \|(\mathbf{S}_{xy}, \mathbf{S}_{yy})\|^2 + \xi \|\mathbf{S}_{yy}\|^2_{x=L}+  L^{\f14-} \|(\phi, p, S)\|^2_{\mathbf{X}^{\v}} +  C_{\xi} \|(\mathfrak{h}_1, \mathfrak{N}_{11})\|^2_{H^1}.
\end{align}

For the 4th and 5th terms on RHS of \eqref{D3.6}, we have 
\begin{align}
\eqref{D3.6}_4	&= \f12 \big\la \big\{\rho_s^{-1} \big[ \mathfrak{h}_1 + \mathfrak{N}_{11}\big]\big\}_y, \, \mathbf{S}_{xy} \big\ra - \f12 \big\la \big\{\rho_s^{-1} \big[ \mathfrak{h}_1 + \mathfrak{N}_{11}\big]\big\}_x, \, \mathbf{S}_{yy} \big\ra \nonumber\\
&\quad   + \f12 \big\la \rho_s^{-1} \big[ \mathfrak{h}_1 + \mathfrak{N}_{11}\big], \, \mathbf{S}_{yy} \big\ra_{x=L} \nonumber\\
	&\lesssim \xi\|(\mathbf{S}_{xy}, \mathbf{S}_{yy})\|^2 + \sqrt{L} \|\mathbf{S}_{yy}\|^2_{x=L}  + C_{\xi}\f{1}{\v}\|(\mathfrak{h}_1, \mathfrak{N}_{11})\|^2_{H^1},
\end{align}
and
\begin{align}\label{D3.51}
\eqref{D3.6}_5 &= \mu \v^{N_0} \big\la\big\{ \rho (\rho^{\v}\rho_s)^{-1} \Delta_{\v}\phi_{y}\big\}_x, \, \mathbf{S}_{yy} \big\ra - \mu \v^{N_0} \big\la \rho (\rho^{\v}\rho_s)^{-1} \Delta_{\v}\phi_{y}, \, \mathbf{S}_{yy} \big\ra_{x=L} \nonumber\\
	&\quad - \mu \v^{N_0} \big\la \big\{ \rho (\rho^{\v}\rho_s)^{-1} \Delta_{\v}\phi_{y}\big\}_y, \, \mathbf{S}_{xy} \big\ra\nonumber\\
	&\lesssim \v^{\f{N_0}{2}} \big\{\|(\mathbf{S}_{xy}, \mathbf{S}_{yy})\| + \|\mathbf{S}_{yy}\|_{x=L}^2 + [\phi]_{3,1} + [\phi]_{4,1}\big\}.
\end{align}

\smallskip

{\it Step 6.4.} Substituting \eqref{D3.33} and \eqref{D3.50}-\eqref{D3.51}  into \eqref{D3.6}, one has
\begin{align}\label{D3.21}
	\big|\big\la \mathbf{S}_{x}, \, \mathbf{S}_{xy} \big\ra_{y=0}\big|
	&\lesssim  \xi \|(\sqrt{\v}\mathbf{S}_{xx}, \mathbf{S}_{xy}, \mathbf{S}_{yy})\|^2  + \xi \|\sqrt{u_s}\mathbf{S}_x\|^2_{x=L}   + \xi \|\mathbf{S}_{yy}\|^2_{x=L} +  L^{\f14-} [[[\mathbf{S}]]]_{2,1}\nonumber\\
	&\quad  + C_{\xi} [\phi]_{4,1} + C_{\xi} \f{1}{\v} [[p]]_{2,1}  + C_{\xi}a \|(\phi, p, S)\|^2_{\mathbf{X}^{\v}} +  C_{\xi}\f{1}{\v}\|(\mathfrak{h}_1, \mathfrak{N}_{11})\|^2_{H^1} +  \v^{N_0}.
\end{align}

{\it Step 7.} Substituting \eqref{D3.21} into \eqref{D3.5}, 
we conclude \eqref{D3.8}.
Therefore the proof of Lemma \ref{lemD3.1} is completed. $\hfill\Box$

\begin{lemma}\label{lemD3.3}
	It holds that 
	\begin{align}\label{D3.61}
		\|u_s \mathbf{S}_{x}w\|^2 + \|\mathbf{S}_{yy}w\|^2 
		&\lesssim L^{-} [[[\mathbf{S}]]]_{2,\hat{w}}  + L^{1-} \|(\phi, p, S)\hat{w}\|^2_{\mathbf{X}^{\v}} 
		+  L\Big\{\f{1}{\v^2} [[p]]_{2,1} + \f{1}{\v} [[[S]]]_{2,1}\Big\} \mathbf{1}_{\{w=w_0\}}\nonumber\\
		&\quad     + \|\mathfrak{N}_3 w\|^2.
	\end{align}
\end{lemma}

\noindent{\bf Proof.} It follows from \eqref{2.111-6D} that 
\begin{align}\label{D3.62}
	\|\big(2\rho_su_s \mathbf{S}_x - \kappa \mathbf{S}_{yy}\big) w\|^2 
	&\lesssim \v \|\sqrt{\v} \mathbf{S}_{xx} w\|^2 + L^2 \|\mathbf{S}_{xy}w\|^2  +   L^{1-} [\phi]_{3,\hat{w}}  + a^2L [[[S]]]_{2,1} 
	\nonumber\\
	&\quad   +   L  \f{1}{\v} [[p]]_{2,\hat{w}}  +  La\Big\{\f{1}{\v^2} [[p]]_{2,1} + \f{1}{\v} [[[S]]]_{2,1}\Big\} \mathbf{1}_{\{w=w_0\}}  + \|\mathfrak{N}_3 w\|^2,
\end{align}
where we have used \eqref{1TC.59-1}, \eqref{1TC.18} and $\eqref{1TC.10-1}_2$.

It is clear that 
\begin{align}\label{D1TC.59}
	\|\big(2\rho_su_s \mathbf{S}_x - \kappa \mathbf{S}_{yy}\big) w\|^2
	&=4\|\rho_s u_s \mathbf{S}_{x} w\|^2 + \kappa^2 \|\mathbf{S}_{yy} w\|^2 - 4\kappa  \big\la \rho_s u_s \mathbf{S}_{x},\, \mathbf{S}_{yy} w^2\big\ra.
\end{align}
A direct calculation shows that 
\begin{align}\label{D1TC.61}
	- 4 \big\la \rho_s u_s \mathbf{S}_{x},\, \mathbf{S}_{yy} w^2\big\ra 
	&=  4 \big\la \rho_s u_s \mathbf{S}_{xy},\, \mathbf{S}_{y} w^2\big\ra + 4 \big\la (\rho_s u_s w^2)_y \mathbf{S}_{x},\, \mathbf{S}_{y}\big\ra \nonumber \\
	&\lesssim L^{-}\|\mathbf{S}_{xy} w\|^2 + L\|\mathbf{S}_{xy}w\| \big\{\|u_s \mathbf{S}_x w_y\| + \|\mathbf{S}_x\|\big\}\nonumber\\
	&\lesssim L^{-}\|\mathbf{S}_{xy} \hat{w}\|^2 + L^- \|u_s \mathbf{S}_x \hat{w}\|^2.
\end{align}
Hence we conclude \eqref{D3.61} from \eqref{D3.62}-\eqref{D1TC.61}.  
Therefore the proof of Lemma \ref{lemD3.3} is completed. $\hfill\Box$

\begin{corollary}
It holds that 
\begin{align}\label{D3.67-0}
[[[\mathbf{S}]]]_{2,\hat{w}^{\v}_0} &\lesssim    [\phi]_{4,1} + \f1{\v}[[p]]_{2,\hat{w}^{\v}_0}	+ a\|(\phi, p, S)\hat{w}^{\v}_0\|^2_{\mathbf{X}^{\v}}     + \f{1}{\v}\|\mathfrak{N}_3 \hat{w}^{\v}_0\|^2   +  \f{1}{\v}\|(\mathfrak{h}_1, \mathfrak{N}_{11})\|^2_{H^1} +  \v^{N_0}.
\end{align}
\end{corollary}
\noindent{\bf Proof.} We conclude \eqref{D3.67-0} directly from \eqref{D3.8}, \eqref{D3.61}, \eqref{D3.30} and \eqref{D3.2-0}.

\medskip

\begin{lemma}\label{lemD3.4}
It holds that 
\begin{align}\label{D3.66}
[[[\mathbf{S}]]]_{3,\hat{w}^{\v}_0}
&\lesssim a^4  [[[\mathbf{S}]]]_{2,\hat{w}^{\v}_0} +  \|u_s \Delta_{\v}p_y \hat{w}^{\v}_0\|^2  
+ L^{\f18} \|(\phi, p, S)\hat{w}^{\v}_0\|^2_{\mathbf{X}^{\v}} +  \|u_s \pa_y\mathfrak{N}_3 \hat{w}^{\v}_0\|^2. 
\end{align}
\end{lemma}

\noindent{\bf Proof.} Applying $\pa_y$ to \eqref{2.111-6D}, one gets
\begin{align}\label{D3.67}
	& 2 \rho_s u_s \mathbf{S}_{xy} + 2\rho_s v_s \mathbf{S}_{yy} + 2 (\rho_s u_s)_y \mathbf{S}_x + 2(\rho_s v_s)_y \mathbf{S}_y  - \kappa \Delta_{\v} \mathbf{S}_y   \nonumber\\
	&=  \f{1}{2}\kappa\Delta_{\v}\big(\rho_s^{-1} p\big)_y + \kappa \Delta_{\v}\big(p_s^{-1}T_{sy}q\big)_y +  \mathfrak{J}_y.
\end{align}
Then we have from \eqref{D3.67} that
\begin{align}\label{D3.68}
&  \kappa^2  \|u_s (\v \mathbf{S}_{xxy}, \mathbf{S}_{yyy}) w\|^2  +  2\kappa^2 \v  \big\la \mathbf{S}_{yyy},\, u_s^2 \mathbf{S}_{xxy} w^2\big\ra  \nonumber\\
&\lesssim  \|u_s \Delta_{\v}p_y w\|^2 +  a^4  [[[\mathbf{S}]]]_{2,\hat{w}}  
 +  L\Big\{\f{1}{\v^2} [[p]]_{2,1} + \f{1}{\v} [[[S]]]_{2,1}\Big\} \mathbf{1}_{\{w=w_0\}} \nonumber\\
&\quad + L^{\f18} \|(\phi, p, S)\hat{w}\|^2_{\mathbf{X}^{\v}}   +  \|u_s \pa_y\mathfrak{N}_3 w\|^2,
\end{align}
where we have used \eqref{2.84}, \eqref{1TC.70}   and \eqref{1TC.79}-\eqref{1TC.81}.

We integrate by parts to obtain
\begin{align}\label{D3.69}
\big\la \v u_s^2 \mathbf{S}_{xxy}, \mathbf{S}_{yyy} w^2\big\ra 
&=  - \big\la \v u_s^2 \mathbf{S}_{xy}, S_{xyyy} w^2\big\ra -\big\la \v 2u_su_{sx} \mathbf{S}_{xy}, \mathbf{S}_{yyy} w^2 \big\ra + \big\la \v u_s^2 \mathbf{S}_{xy}, \mathbf{S}_{yyy} w^2\big\ra_{x=L}  \nonumber\\
&=\|u_s\sqrt{\v} \mathbf{S}_{xyy}w\|^2 + \big\la \v u_s^2 \mathbf{S}_{xy}, \mathbf{S}_{yyy} w^2\big\ra_{x=L} +  O(1) \sqrt{\v} \|\mathbf{S}_{xy} \hat{w}\|^2\nonumber\\
&\quad  + O(1) \sqrt{\v} \|u_s(\sqrt{\v} \mathbf{S}_{xyy}, \mathbf{S}_{yyy})w\|^2.
\end{align}

Next we shall control the boundary term  on RHS of \eqref{D3.69}. Applying $\pa_y$ to the boundary condition of $\eqref{2.111-6D}_2$ at $x=L$, one has
\begin{align}\label{D3.69-1}
 \kappa \mathbf{S}_{yyy} \big|_{x=L}&=2\rho_s u_s  \mathbf{S}_{xy} + 2 (\rho_s u_s)_y  \mathbf{S}_x - \kappa \big\{p_s^{-1}T_{sy}q_{yy}\big\}_y - 2\kappa \big\{\big(p_s^{-1}T_{sy}\big)_{y}q_y\big\}_y \nonumber\\
 &\quad  - \kappa\big\{\big(\rho_s^{-1}\big)_{yy}p\big\}_y + \pa_y \mathfrak{M}_0.
\end{align}
Noting \eqref{D3.76S}, we get
\begin{align}\label{D3.77S}
\|\pa_y \mathfrak{M}_0w\|^2_{x=L}
&\lesssim a^4\|S_{yy}w\|_{x=L}^2 + a^3 [[p]]_{3,w}  + a^4 L \f{1}{\v}[[p]]_{2,w}  + a^2 L\f{1}{\v^2} [[p]]_{2,1}\mathbf{1}_{\{w=w_0\}} \nonumber\\
&\quad + L [[[S]]]_{2,\hat{w}} + a^2 [\phi]_{3,1} + a^2 [\phi]_{4,1},
\end{align}
which yields that 
\begin{align}\label{D3.73}
	&\v\big\la u_s^2 \mathbf{S}_{xy} w^2, \,  2\rho_s u_s  \mathbf{S}_{xy} + 2 (\rho_s u_s)_y  \mathbf{S}_x  - 2\kappa \big\{\big(p_s^{-1}T_{sy}\big)_{y}q_y\big\}_y  - \kappa\big\{\big(\rho_s^{-1}\big)_{yy}p\big\}_y + \pa_y \mathfrak{M}_0\big\ra_{x=L} \nonumber\\
	&\geq \f74 \v \|\sqrt{\rho_su_s^3} \mathbf{S}_{xy}w\|^2_{x=L} - \v \|\sqrt{u_s} \mathbf{S}_{x} w\|^2_{x=L} -  L \|(\phi, p, S)\hat{w}\|^2_{\mathbf{X}^{\v}}.
\end{align}

Recall $y^{\ast}$ in \eqref{D.7}, then we note
\begin{align*}
	\int_0^\infty \f{1}{y^*} \big|\f{u_s}{p_s}T_{sy}q_{yy} w\big|^2 dy \Big|_{x=L} 
	&\lesssim  \|\sqrt{u_s} q_{yy} w\|^2_{x=L} \lesssim L \|\sqrt{u_s}q_{xyy}w\|^2
	\lesssim L [\phi]_{3,w},
\end{align*}
which, together with Lemma \ref{lemAH.1}, yields that 
\begin{align}\label{D3.76}
	\|\big(\f{u_s}{p_s}T_{sy}q_{yy} w\big)(L,\cdot)\|^2_{H^{1/2}_{00}}
	&\lesssim La^3 [\phi]_{3,w} + \|\big(\f{u_s}{p_s}T_{sy}q_{yy} w\big)_y\| \Big\{\f{1}{L}\|u_sT_{sy}q_{yy}w\| + \|\big(\f{u_s}{p_s}T_{sy}q_{yy}w\big)_x\|\Big\} \nonumber\\
	&\lesssim L^{\f18}[\phi]_{3,\hat{w}}.
\end{align}
Also it holds that 
\begin{align}\label{D3.77}
	\v^{\f32}\|(u_s \mathbf{S}_{xy} w)(L,\cdot)\|^2_{H^{1/2}} 
	&\lesssim \v^{\f32}\|\big((u_s \mathbf{S}_{xy}w)_y, u_s \mathbf{S}_{xy}w\big)\| \Big\{\f{1}{L}\|u_s \mathbf{S}_{xy}w\| + \|(u_s \mathbf{S}_{xy} w)_x\|\Big\} \nonumber\\
	&\lesssim \| u_s(\v \mathbf{S}_{xxy}, \sqrt{\v}\mathbf{S}_{xyy})w\|^2 + \v \|\mathbf{S}_{xy}\hat{w}\|^2.
\end{align}

Hence it follows from \eqref{D3.30-0} and  \eqref{D3.76}-\eqref{D3.77}  that 
\begin{align}\label{D3.76-0}
&\v\big\la u_s^2 \mathbf{S}_{xy} w^2, \, \big\{p_s^{-1}T_{sy}q_{yy}\big\}_y\big\ra_{x=L} \nonumber\\
& =\v\big\la u_s \mathbf{S}_{xy} w, \, \big\{p_s^{-1} u_s T_{sy}q_{yy}\big\}_y\big\ra_{x=L}  - \v\big\la u_s w \mathbf{S}_{xy} , \, p_s^{-1}T_{sy}q_{yy} (u_sw)_y \big\ra_{x=L}  \nonumber\\
&\lesssim  \v \|(u_s \mathbf{S}_{xy} w)(L,\cdot)\|_{H^{1/2}} \|\big(u_s p_s^{-1}T_{sy}q_{yy} w\big)(L,\cdot)\|_{H^{1/2}_{00}} + \sqrt{\v} \|u_s \v \mathbf{S}_{xxy}\|^2 \nonumber\\
&\quad + \sqrt{\v} \|\mathbf{S}_{xy}\hat{w}\|^2 + \sqrt{\v} [\phi]_{3,\hat{w}}\nonumber\\
&\lesssim \sqrt{\v}  \| u_s(\v \mathbf{S}_{xxy}, \sqrt{\v}\mathbf{S}_{xyy}) w\|^2  + \sqrt{\v} \|\mathbf{S}_{xy}\hat{w}\|^2 + L^{\f18} [\phi]_{3,\hat{w}}.
\end{align}
which, together with \eqref{D3.69-1} and \eqref{D3.73},  that 
\begin{align}\label{D3.77-0}
	\v\big\la u_s^2 \mathbf{S}_{xy} w^2, \,  \kappa \mathbf{S}_{yyy} \big\ra_{x=L}
	&\geq  \f34 \v \|\sqrt{\rho_su_s^3} \mathbf{S}_{xy}w\|^2_{x=L} - \sqrt{\v}  \| u_s(\v \mathbf{S}_{xxy}, \sqrt{\v}\mathbf{S}_{xyy})w\|^2 \nonumber\\
	&\quad  
	- \sqrt{\v} [[[\mathbf{S}]]]_{2,\hat{w}}  - L^{\f18} \|(\phi, p, S)\hat{w}\|^2_{\mathbf{X}^{\v}}.
\end{align}

Finally, combining \eqref{D3.68}-\eqref{D3.69} and  \eqref{D3.77-0}, 
then one concludes \eqref{D3.66}. Therefore the proof of Lemma \ref{lemD3.4} is completed. $\hfill\Box$

\medskip

\noindent{\bf Proof of Theorem \ref{thm7.1}.}  We conclude Theorem \ref{thm7.1} directly from  \eqref{D3.30-1} and \eqref{D3.67-0}-\eqref{D3.66}. Therefore the proof is completed. $\hfill\Box$

\section{Proof of Main Theorem}\label{Sec7}

For fixed $(\phi, p, S)$, we aim to prove that the mapping
\begin{align}\label{18.1}
	(\hat{P}, \hat{\mathbf{S}}) \longrightarrow (\Phi, P, \mathbf{S})
\end{align}
is a contraction under certain constraints on $a$ and $L$. Regarding the pseudo-entropy in the Dirichlet case, we emphasize that $\hat{\mathbf{S}} \equiv \mathbf{S} \equiv \zeta$ which is independent of the aforementioned mapping \eqref{18.1}.

To derive the contraction, we define the notations $\mathbb{Y}^{\v}$ and $\mathbb{X}^{\v}$:
\begin{align}\label{18.10}
	\|(\phi, p, S)w\|^2_{\mathbb{Y}^{\v}} :&=[\phi]_{3,w} + n_2 [\phi]_{4,w} + n_1 \Big[\mathscr{K}(p) + \f{1}{\v}\mathcal{K}_{2,w}(p)  + \|\f{Bw}{\sqrt{u_s}}\|_{x=0}^2 \Big]   + C_1 n_3\mathcal{K}_{3,w}(p)\nonumber\\
	&\, + \xi \bar{\mathcal{K}}_{3,w}(p) + n_4 \|u_s\Delta_{\v}p_y w\|^2 \mathbf{1_D} + n_3 [[[S]]]_{3,w}  + [[[S]]]_{2,w} + [[[S]]]_{b,w},\\
	\|(\phi, p, S)w\|^2_{\mathbb{X}^{\v}}:&=[\phi]_{3,w} + n_2 [\phi]_{4,w}  + n_1 \f{1}{\v}[[p]]_{2,w}  + C_1 n_3 [[p]]_{3,w} + \xi \|u_s \sqrt{\v}\Delta_{\v}p_x w\|^2  \nonumber\\
	&\quad  + n_4 \|u_s\Delta_{\v}p_y w\|^2 \mathbf{1_D} +  n_3 [[[S]]]_{3,w}  + [[[S]]]_{2,w} + [[[S]]]_{b,w} ,\label{18.11}
\end{align}
where the positive constants $C_1, n_1,n_2, n_3, n_4\geq 1$ and $\xi\in(0,1)$ will be determined later. 

\smallskip

Combining Theorems \ref{thm4.1}, \ref{thm5.1}, \ref{thm6.7} \& \ref{thm7.1}, and using Lemma \ref{THlem6.11},  one has that 
\begin{align}\label{18.12}
&\|(\Phi, P, \mathbf{S})\hat{w}_0^{\v}\|^2_{\mathbb{Y}^{\v}} \nonumber\\
&\lesssim_{C_0} \tilde{C}_2 a \|(\Phi,0,0)\|^2_{\mathbf{X}^{\v}} + (\f{1}{N} + \tilde{C}_1 a ) \|(\phi, p, S(\zeta))\hat{w}_0^{\v}\|^2_{\mathbf{X}^{\v}}  +   \tilde{C}_2 a \|(0, \hat{P}, \hat{\mathbf{S}})\hat{w}_0^{\v}\|^2_{\mathbf{Y}^{\v}} \nonumber\\
&\quad  + n_2  \mathcal{K}_{3,\hat{w}_0^{\v}}(\hat{P})  + n_3  [[p]]_{3,\hat{w}_0^{\v}} + n_2  [[[\zeta]]] _{3,\hat{w}_0^{\v}}  + n_2  [[[\hat{\mathbf{S}}]]]_{3,\hat{w}_0^{\v}}  + [n_1 + C_1 n_3 \xi] \|||\Phi||\|^2_1  \nonumber\\
&\quad + n_3 \sqrt{\v} [[[\mathbf{S}]]]_{b,1} + n_3 L^{\f12} [[[\mathbf{S}]]]_{2,\hat{w}^{\v}_0}   +  \f{1}{\v}\mathcal{K}_{2,1}(P)  + n_3 \mathcal{K}_{3,1}(P) + \mathscr{K}(P)    +  \f{1}{\v} [[p]]_{2,\hat{w}_0^{\v}}\nonumber\\
&\quad  +  [\phi]_{4,1} + n_3  \|u_s \Delta_{\v}p_y \hat{w}^{\v}_0\|^2  \mathbf{1_D} + \tilde{C}_1\v \|\sqrt{u_s} \mathbf{S}_{xy}\hat{w}_0^{\v}\|^2_{x=L}\mathbf{1_{D}} + \tilde{C}_1\|u_s \nabla_{\v}^2\mathbf{S}_{y}\hat{w}_0^{\v}\|^2\mathbf{1_D} \nonumber\\
&\quad + \tilde{C}_2 \{\cdots\}_{1,\hat{w}_0^{\v}} + \tilde{C}_2 \{\cdots\}_{2,\hat{w}_0^{\v}}   + \tilde{C}_2 \|(\mathbf{d}_2 \bar{g}_1, \bar{g}_1\mathbf{1}_{\{y\in[0,2]\}})\|^2  + n_3 \|\nabla_{\v}\mathfrak{N}_3 \hat{w}^{\v}_0\|^2 \nonumber\\
&\quad + \f{1}{\v} \|(\mathfrak{N}_3, \mathcal{F}_{R8}) \hat{w}_0^{\v}\|^2   +  \f{1}{\v}\|(\mathfrak{h}_1, \mathfrak{N}_{11})\|^2_{H^1}
  +  \v^{N_0 }  
\end{align}
where we have denoted $\tilde{C}_1:=C(n_1,n_2, n_3, \xi, N)$, $\tilde{C}_2:=C(\tilde{C}_1, n_4)$, and $\hat{\mathbf{S}}=\mathbf{S}$ in the case of DT. 


Noting the notations in Section \ref{sec:appendixA}, using Theorem \ref{THlem6.10} and Lemma  \ref{THlem6.11}, and by tedious calculation, we can obtain
\begin{align}\label{18.15}
\begin{split}
\|(0,0, \zeta)\hat{w}_0^{\v}\|^2_{\mathbb{X}^{\v}} &\lesssim \|(\phi,p, S)\hat{w}_0^{\v}\|^2_{\mathbb{X}^{\v}} + \v^{N_0},\\
\{\cdots\}_{1,\hat{w}_0^{\v}} + \{\cdots\}_{2,\hat{w}_0^{\v}} &\lesssim a^2\|(\phi,p, S)\hat{w}_0^{\v}\|^2_{\mathbb{X}^{\v}} + \v^{N_0},\\
\|(\mathbf{d}_2 \bar{g}_1, \bar{g}_1\mathbf{1}_{\{y\in[0,2]\}})\|^2 &\lesssim  L \|(\phi, p, S)\hat{w}^{\v}_0\|^{2}_{\mathbb{X}^{\v}}  + \v^{N_0},\\
\f{1}{\v}\|(\mathfrak{h}_1, \mathfrak{N}_{11})\|^2_{H^1} + \f{1}{\v}\|(\mathfrak{N}_3,\mathcal{F}_{R8}) \hat{w}_0^{\v}\|^2  + \|\nabla_{\v}\mathfrak{N}_3 \hat{w}^{\v}_0\|^2 &\lesssim L \|(\phi, p, S)\hat{w}^{\v}_0\|^{2}_{\mathbb{X}^{\v}}  + \v^{N_0}.
\end{split}
\end{align}
Using \eqref{D3.30-0},    \eqref{18.9-3} and \eqref{18.15}, we have 
\begin{align}\label{18.15-1}
&\tilde{C}_1\v \|\sqrt{u_s} \mathbf{S}_{xy}\hat{w}_0^{\v}\|^2_{x=L}\mathbf{1_{D}} +  \tilde{C}_1\|u_s \nabla_{\v}^2\mathbf{S}_{y}\hat{w}_0^{\v}\|^2\mathbf{1_D} \leq  \f{1}{n_1} [[[\mathbf{S}]]]_{2,\hat{w}^{\v}_0} \mathbf{1_{D}}  + \tilde{C}_1 [[[\mathbf{S}]]]_{3,\hat{w}^{\v}_0} \mathbf{1_{D}}  \nonumber\\
&\lesssim \big\{\f{1}{n_1} + \tilde{C}_1a^4\big\} [[[\mathbf{S}]]]_{2,\hat{w}^{\v}_0}   + \tilde{C}_1 \|u_s \Delta_{\v}p_y \hat{w}^{\v}_0\|^2  \mathbf{1_{D}}
+ \tilde{C}_1 L^{\f18} \|(\phi, p, S)\hat{w}^{\v}_0\|^2_{\mathbb{X}^{\v}} + \tilde{C}_1 \v^{N_0}.
\end{align}
This is the main reason why we need the refined pressure estimate \eqref{C17.223-1}.

We first fix $C_1, n_1, n_2, n_3$ such that $n_3\gg n_2\gg n_1 \gg C_1\gg C_0\geq 1$ and $n_4\gg \tilde{C}_1$, then take $0<\v\ll L\ll a \ll \xi \ll  1$ and $N\gg \tilde{C}_1$.  Then we have from \eqref{18.12}-\eqref{18.15-1} that 
\begin{align}\label{18.18}
\|(\Phi, P, \mathbf{S})\hat{w}_0^{\v}\|^2_{\mathbb{Y}^{\v}}
&\lesssim_{C_0}  o(1) \|(0, \hat{P}, \hat{\mathbf{S}})\hat{w}_0^{\v}\|^2_{\mathbb{Y}^{\v}}  + o(1) \|(\phi, p, S )\hat{w}_0^{\v}\|^2_{\mathbb{X}^{\v}} + \tilde{C}_2\v^{N_0}.
\end{align}
We remark that above $o(1)$ depends on $\f{1}{n_1}, \f{n_1}{n_2}, \f{n_2}{n_3}, \f{1}{C_1}, \f{\tilde{C}_1}{n_4}$ and $a, \xi, \f{1}{N}$.

\begin{proposition}\label{Prop8.1}
For given $(\phi, p, S)$, there exists a fixed point $(\Phi, P, \mathbf{S})$ for the mapping \eqref{18.1} satisfying 
\begin{align}\label{18.18-1}
\|(\Phi, P, \mathbf{S})\hat{w}_0^{\v}\|^2_{\mathbb{Y}^{\v}}
&\lesssim_{C_0}  o(1) \|(\phi, p, S )\hat{w}_0^{\v}\|^2_{\mathbb{X}^{\v}} + \tilde{C}_2\v^{N_0}.
\end{align}
Specifically, $(\Phi, P, \mathbf{S})$ satisfies the following boundary value problems
\begin{equation}\label{18.27}
	\begin{cases}
		\dis \big(b_1 m u_s^2 \mathbf{q}_{xy}\big)_y  + \v \big(b_2 mu_s^2 \mathbf{q}_{xx}\big)_x   - \mu\Big(\frac{b_1}{\rho^{\v}} \Delta_{\v}\Phi_y\Big)_y-\mu\v \Big(\frac{b_2}{\rho^{\v}}   \Delta_{\v}\Phi_x\Big)_x \\
		 \dis \qquad \qquad\qquad\qquad\qquad\qquad\qquad\quad   =\mathcal{F}_{R}(P, \mathbf{S},  \zeta, p, \phi, S),\quad (x,y)\in (0,L)\times \R_+,\\
		\dis \Phi|_{x=0}=\Phi_{xx} |_{x=0}= \Phi_{x}|_{x=L}=\Phi_{xxx}|_{x=L}=0,
	\end{cases}
\end{equation}
and 
\begin{align}\label{18.28}
\begin{cases}
\dis \Delta_{\v}P + \v (\mu+\lambda)  \f{1}{p^{\v}}(U^{\v}\cdot\nabla)\Delta_{\v}P  = \mathcal{G}(P, \mathbf{S},\Phi, p,S, \zeta, \phi), \quad (x,y)\in (0,L)\times \R_+,\\
(\mu+\lambda) \f{1}{p_s} u_s\Delta_{\v}P \big|_{x=0}=-\sigma \mathbf{d}_{11} B,  \\
\mathbf{A}_1\big|_{x=0}=\mathbf{A}_2\big|_{x=L}=\mathbf{A}_1\big|_{y=0}=0, 
\end{cases}
\end{align}
along with either the Neumann BVP for the pseudo entropy
\begin{align} \label{18.29}
	\begin{cases}
		\dis 2 \rho_s u_s \mathbf{S}_x + 2\rho_s v_s \mathbf{S}_y - \kappa \Delta_{\v}\mathbf{S} = \kappa \Delta_{\v}\big(\frac{1}{2\rho_s}[\chi P + \bar{\chi}p]\big) + \kappa \Delta_{\v}\big(\f{1}{p_s}T_{sy}q\big)  \\
		\dis\qquad\qquad\qquad\qquad\qquad\qquad\quad  + \mathfrak{J}(\phi,p,T(S,p,q)), \quad (x,y)\in (0,L)\times \R_+,\\[2mm]
		\dis \mathbf{S} \big|_{x=0}= -  \frac{\chi}{2\rho_s}P,\,\,\, \mathbf{S}_{x}\big|_{x=L}=0 \,\, \,\mbox{and}\,\,\,
		\mathbf{S}_y \big|_{y=0} = - \big(\frac{1}{2\rho_s}P\big)_y, 
	\end{cases}
\end{align} 
or the Dirichlet BVP for the pseudo entropy
\begin{align}\label{18.29D}
\begin{cases}
\dis	2 \rho_s u_s \mathbf{S}_x + 2\rho_s v_s \mathbf{S}_y - \kappa \Delta_{\v}\mathbf{S}  = \kappa \Delta_{\v}\big(\frac{1}{2\rho_s} p \big) + \kappa \Delta_{\v}\big(\f{1}{p_s}T_{sy}q\big) \\
\dis\qquad\qquad\qquad\qquad\qquad\qquad\quad  + \mathfrak{J}(\phi,p,T(S,p,q)), \quad (x,y)\in (0,L)\times \R_+,\\
\dis\mathbf{S} |_{y=0}=-\f1{2\rho_s}\mathfrak{P}(x),\quad \mathbf{S} |_{x=0}=0 \,\,\,\, \mbox{and}\,\,\, \big\{2\rho_s u_s \mathbf{S}_x-\kappa \mathbf{S}_{yy}\big\}|_{x=L}= \mathfrak{G}.
\end{cases}
\end{align}
where $\mathbf{A}_1$ and $\mathbf{A}_2$ are the ones in \eqref{2.109}, and $B=P_x+\bar{g}_1$. 
\end{proposition}



\noindent{\bf Proof.} Noting \eqref{18.18}, it is clear that the mapping \eqref{18.1} is contract. Then, by the standard contraction argument, there is a fixed point  $(\Phi, \hat{P}, \hat{\mathbf{S}})=(\Phi, P, \mathbf{S})$ to the boundary value problems \eqref{2.121},  \eqref{18.39-1}-\eqref{18.39-2}, and \eqref{2.111} (or \eqref{2.111-6D}). By noting  \eqref{18.39-2} and Lemma \ref{lemAF-1}, we can derive the boundary condition $\eqref{18.28}_3$.  Finally, it is clear that \eqref{18.18-1} follows directly from \eqref{18.18}. Therefore the proof of Proposition \ref{Prop8.1} is completed. $\hfill\Box$

\smallskip


\begin{lemma}\label{lem9.1}
Recall the definition of $\mathbf{A}_{1}, \mathbf{A}_2$ in \eqref{2.109}. For $\Phi, P, \mathbf{S}$ constructed in Proposition \ref{Prop8.1}, it holds that 
\begin{align}\label{19.0}
\mathbf{A}_{1}=\mathbf{A}_{2}\equiv 0.
\end{align}
That means we solve the momentum equations \eqref{2.109}.
\end{lemma}

\noindent{\bf Proof.} From \eqref{18.27} and \eqref{18.28}, we know that $(\mathbf{A}_{1}, \mathbf{A}_{2})$ satisfies the following {\it div-curl} system
\begin{align}\label{19.1}
	\begin{cases}
		\v \pa_x \mathbf{A}_1 + \pa_y \mathbf{A}_2  =0,\\
		\pa_y\big(b_1\mathbf{A}_1\big) - \pa_{x} \big(b_2 \mathbf{A}_2\big) =0,\\
		\mathbf{A}_1\big|_{x=0}=\mathbf{A}_2\big|_{x=L}=\mathbf{A}_1\big|_{y=0}=0.
	\end{cases}
\end{align}
Applying the div-curl inversion Lemma \ref{lemR.8} to \eqref{19.1}, we immediately conclude \eqref{19.0}. Therefore the proof of Lemma \ref{lem9.1} is completed.  $\hfill\Box$

\smallskip

Noting \eqref{18.18-1}, we only have a weak estimate on $(P_x,P_y)$ which is not good enough for us to close  the mapping $(\phi,p,S) \to (\Phi,P,\mathbf{S})$. By using \eqref{19.0}, we can have a better $\v$-decay estimate on $(P_x,P_y)$, then derive control on $\|(\Phi, P, \mathbf{S})\hat{w}_0^{\v}\|^2_{\mathbb{X}^{\v}}$.
\begin{lemma}\label{lem9.2}
	It holds that 
	\begin{align}\label{7.15}
		\|(\Phi, P, \mathbf{S})\hat{w}_0^{\v}\|^2_{\mathbb{X}^{\v}}
		&\lesssim_{C_0}  o(1) \|(\phi, p, S )\hat{w}_0^{\v}\|^2_{\mathbb{X}^{\v}} + \tilde{C}_2\v^{N_0}.
	\end{align}
\end{lemma}

\noindent{\bf Proof.} We have from \eqref{19.0}  that 
\begin{align*}
	\|\sqrt{\v}P_x w\|^2 
	&\lesssim a L^{\f14} \v [\Phi]_{3,w}  + a^2 \v [[P]]_{2,w}+ L^{\f38}a^2 \v  [[[\mathbf{S}]]]_{2,w} + L^{\f38}a^2 \v [[[\zeta]]]_{2,w}  + \v \|g_1w\|^2 \nonumber\\
	&\quad  + L \v [\phi]_{3,1} + a^2 \v \|\nabla^2_{\v}\tilde{T}w\|^2 + L^2 [[p]]_{2,1} + \v \|\mathcal{N}_1w\|^2,\\
	\|P_yw\|^2 &\lesssim a^2 L^{\f14} \v [\Phi]_{3,w}  + a^2 \v [[P]]_{2,w} + a^2 \v^2  [[[\mathbf{S}]]]_{2,w} + \|g_2w\|^2 +   L \v [\phi]_{3,1}  + \v [[p]]_{2,1} \nonumber\\
	&\quad  + a\v^2 \|\nabla_{\v}^2 \tilde{T} w\|^2 + \v^2 \|\mathcal{N}_2w\|^2,
\end{align*}
which yield that  
\begin{align}\label{19.3}
	\v^{-1}\|(\sqrt{\v}P_x, P_y) w\|^2 
	&\lesssim a^2 \|(\phi, p, S)\hat{w}_0^{\v}\|^2_{\mathbb{X}^{\v}}  + \tilde{C}  \v^{N_0}.
\end{align}

Similar as \eqref{17.218-0}, it holds that
\begin{align*}
	\frac{1}{\v} \mathfrak{K}_{2,w}(p)
	&\lesssim \f{1}{\v}\mathcal{K}_{2,w}(p) +   \mathcal{K}_{3,w}(p),
\end{align*}
which, together with \eqref{19.3}, \eqref{18.18-1} and \eqref{18.10}-\eqref{18.11}, concludes \eqref{7.15}. Therefore the proof of Lemma \ref{lem9.2} is completed. $\hfill\Box$

\medskip

\noindent{\bf Proof of Theorem \ref{thm1.1}:}  Now it is time to consider the mapping  $(\phi, p, S) \to (\Phi, P, \mathbf{S})$. Given $(\phi, p, S)=(\phi^0, p^0, S^0)=(0,0,0)$, using Proposition \ref{Prop8.1} and Lemma \ref{lem9.2}, then we can construct a solution $(\Phi,P,\mathbf{S})=(\Phi^1,P^1, \mathbf{S}^1)$ and also $\zeta^1$. Inductively,   we can obtain $(\Phi^{k+1},P^{k+1}, \mathbf{S}^{k+1})$ and $\zeta^{k+1}$ for $k=0,1,2, \cdots$ satisfying
\begin{align*}
	\begin{split}
		\|(0,0, \zeta^{k+1})\hat{w}_0^{\v}\|^2_{\mathbb{X}^{\v}} \leq C_0 \|(\Phi^k,P^k, \mathbf{S}^k)\hat{w}_0^{\v}\|^2_{\mathbb{X}^{\v}} + \tilde{C}_2 \v^{N_0},\\
		\|(\Phi^{k+1}, P^{k+1}, \mathbf{S}^{k+1})\hat{w}_0^{\v}\|^2_{\mathbb{X}^{\v}} 
		\leq o(1) \|(\Phi^k,P^k, \mathbf{S}^k)\hat{w}_0^{\v}\|^2_{\mathbb{X}^{\v}}  + \tilde{C}_2  \v^{N_0}.
	\end{split}
\end{align*}
which immediately yields that
\begin{align}\label{7.70}
	\begin{split}
		\|(\Phi^{k+1}, P^{k+1}, \mathbf{S}^{k+1})\hat{w}_0^{\v}\|^2_{\mathbb{X}^{\v}}& \leq 2 \tilde{C}_2  \v^{N_0},\\
		\|(0,0, \zeta^{k+1})\hat{w}_0^{\v}\|^2_{\mathbb{X}^{\v}} & \leq C_0[2\tilde{C}_2+1] \v^{N_0}.
	\end{split}
\end{align}

Next, we would like to prove 
\begin{equation}\label{7.71}
	\|(\Phi^{k+1}-\Phi^{k}, P^{k+1}- P^{k}, \mathbf{S}^{k+1}-\mathbf{S}^{k})\|^2_{\mathbb{X}^{\v}} 
	\leq o(1) \|(\Phi^{k}-\Phi^{k-1}, P^{k}- P^{k-1}, \mathbf{S}^{k}-\mathbf{S}^{k-1})\|^2_{\mathbb{X}^{\v}}.
\end{equation}
The proof of \eqref{7.71} is very similar to \eqref{18.18}, but with some additional regularity due to the nonlinear coefficient $\f{1}{p^{\v}}U^{\v}$ in the highest order term of pressure transport equation \eqref{18.28}.  Let  
\begin{align*}
	\bar{\Phi}:= \Phi^{k+1} - \Phi^{k}, \quad \bar{P}:=P^{k+1} - P^{k} \quad  \mbox{and}\quad  \bar{\mathbf{S}}:=\mathbf{S}^{k+1} - \mathbf{S}^{k},
\end{align*}
then we have from \eqref{18.28} that 
\begin{align}\label{7.72}
	&\Delta_{\v} \bar{P} + (\mu+\lambda) \v \, \f{1}{P^{\v,k+1}} (U^{\v,k+1}\cdot \nabla) \Delta_{\v}\bar{P} \nonumber\\
	&= (\mu+\lambda)\v^{N_0+1}  \Big\{ \f{U^{k+1}}{p_s+ \v^{N_0}P^{k+1}} - \f{ U^{k}}{p_s+ \v^{N_0}P^{k}}  \Big\} \cdot \nabla \Delta_{\v}P^{k} + (\cdots).
\end{align}
Now we explain a little bit on the key different for the third order derivative term on RHS of \eqref{7.72}. Considering $\dis \int_{0}^L\int_0^\infty \pa_y\eqref{7.72}\cdot u_s^2 \Delta_{\v}\bar{P}_y $, then the new term is
\begin{align}
	&\v^{N_0+1} \Big|\int_0^{L}\int_{\R_+} \Big\{ \f{U^{k+1}}{p_s+ \v^{N_0}P^{k+1}} - \f{ U^{k}}{p_s+ \v^{N_0}P^{k}} \Big\} \cdot \nabla \Delta_{\v}P^{k}_y \cdot u_s^2 \Delta_{\v}\bar{P}_y  dydx\Big| \nonumber\\
	&\lesssim \v^{N_0+1}\|u_s^2\nabla\Delta_{\v}P^{k}_y\hat{w}^{\v}_0\|\cdot \|u_s \Delta_{\v}\bar{P}_y\| \Big\{\|\f{1}{\hat{w}^{\v}_0}(P^{k+1}-P^k)\|_{L^\infty} + \|\f{1}{u_s\hat{w}^{\v}_0} (U^{k+1}-U^k)\|_{L^\infty}\Big\}\nonumber\\
	&\lesssim \v^{N_0-2}\|u_s^2\nabla\Delta_{\v}P^{k}_y \hat{w}^{\v}_0\| \cdot \|(\bar{\Phi},\bar{P}, \bar{\mathbf{S}})\|^2_{\mathbb{X}^{\v}}.
\end{align}
Thus, if $\v^{N_0-2}\|u_s^2\nabla\Delta_{\v}P^{k}_y \hat{w}^{\v}_0\|\ll 1$, then we derive \eqref{7.71}. In fact, by noting the additional space weight $u_s$ and $\v^{N_0-2}$ decay, this can be done  by similar arguments as in  Sections \ref{sec4-0}-\ref{sec6}, for instance, we can get additional estimates
\begin{align}
	\v^{\f12N_0} \Big\{\|u_s \nabla_{\v}^{4}\Phi^{k+1}_{y} \hat{w}^{\v}_0\|+ \|u_s^2\nabla^3_{\v}P^{k+1}_y\| + \|u_s^2\nabla^3_{\v}S^{k+1}_y\|\Big\} \lesssim 1.
\end{align}
We omit the details of proof here for simplicity.

Finally, with the help of \eqref{7.71}, it is clear to known that there is a fixed point $(\phi,p, S)=(\Phi, P, \mathbf{S})$ and  $\zeta=\mathbf{S}$ with $\mathfrak{P}(x)=p(x,0)$, thus one solves \eqref{2.109}-\eqref{2.111-6D} and \eqref{17.15} (or equivalently  \eqref{T.1} and \eqref{2.92} with \eqref{pair1}, \eqref{0.3}, \eqref{pressure}, \eqref{2.17-0}, \eqref{T.2}-\eqref{3.61-11}, and \eqref{3.61-1}). Moreover, the following uniform estimate holds
\begin{align*}
	\|(\Phi, P, \mathbf{S})\hat{w}_0^{\v}\|^2_{\mathbb{X}^{\v}} \lesssim \v^{N_0},
\end{align*}
which  immediately  yields \eqref{1.51}-\eqref{1.52}. Therefore the proof of Theorem \ref{thm1.1} is completed. $\hfill\Box$

\small 
\appendix
\section{Reformulations and Notations}\label{sec:appendixA}
In this section, we reformulate the remainder equations \eqref{1.8-1}.
\subsection{Reformulation for momentum equations } \label{subsec10.1}
Substituting \eqref{7.1-1} into $\eqref{1.8-1}_{2,3}$ and using \eqref{1.19-1}, after tedious calculations, one can obtain that 
	\begin{align}\label{7.36}
	\begin{cases}
		\dis A_1:= [u_s \phi_{xy} - u_{sy}\phi_x] + I_{u}(\phi) + \tilde{d}_{11} p_x^{(\v)}
		+ \tilde{d}_{12} p_y^{(\v)} + \tilde{d}_{13} p^{(\v)} + \hat{d}_{11} T_x^{(\v)}
		+ \hat{d}_{12} T_y^{(\v)} + \hat{d}_{13} T^{(\v)} \\[2mm]
		\dis \qquad \,\, - \f{\mu}{\rho^{\v}}\Delta_{\v}\phi_y+ \mu\frac{u^{\v}}{p^{\v}}  \Delta_{\v}p^{(\v)} - \mu\frac{u^{\v}}{T^{\v}}  \Delta_{\v}T^{(\v)} + \f{\lambda \v}{p^{\v}}  (U^{\v}\cdot\nabla)p^{(\v)}_x - \f{\lambda \v}{T^{\v}} \, (U^{\v}\cdot\nabla) T^{(\v)}_x - h_1 - \mathfrak{h}_1\\
		\dis\qquad  \quad  - \mathfrak{N}_1 - \mathfrak{R}_{s1}=0,\\[2mm]
		\dis A_2:= - \v u_s \phi_{xx} + \v I_{v}(\phi) + \v \tilde{d}_{21} p_x^{(\v)}
		+   \tilde{d}_{22} p_y^{(\v)} + \v \tilde{d}_{23} p^{(\v)}  + \v \hat{d}_{21} T_x^{(\v)}
		+ \v \hat{d}_{22} T_y^{(\v)} + \v \hat{d}_{23} T^{(\v)}\\[2mm]
		\dis \qquad \quad + \f{\mu \v}{\rho^{\v}}\Delta_{\v}\phi_x + \mu\v \frac{v^{\v}}{p^{\v}}  \Delta_{\v}p^{(\v)} - \mu\v \frac{v^{\v}}{T^{\v}}  \Delta_{\v}T^{(\v)} + \f{\lambda \v}{p^{\v}}  (U^{\v}\cdot\nabla)p^{(\v)}_y - \f{\lambda \v}{T^{\v}} \, (U^{\v}\cdot\nabla) T^{(\v)}_y \\
		\dis \qquad \quad - \v  h_2 - \v \mathfrak{h}_2 - \v \mathfrak{N}_2 - \v \mathfrak{R}_{s2}=0,
	\end{cases}
\end{align}
where we have denoted the following notations
\begin{align}
	I_{u}(\phi)&:= v_{s}\phi_{yy} + u_{sx}\phi_y + \f{2\mu}{\rho_s^2} \nabla_{\v}\rho_s\cdot \nabla_{\v}\phi_{y} + \f{\lambda \v}{\rho_s^2} \nabla \rho_s\cdot \nabla^{\perp} \phi_x + \f{\lambda \v}{\rho_s^2} \Big\{ \f{2}{\rho_s} \rho_{sx}\rho_{sy} -  \rho_{sxy}\Big\} \phi_x \nonumber\\
	&\quad - \f{1}{\rho_s}\Big\{(U_s\cdot \nabla)\rho_s  - \f{\mu}{\rho_s} \Delta_{\v}\rho_s  + \f{2\mu}{\rho_s^2} |\nabla_{\v}\rho_s|^2 + \lambda \v \big[\f{2}{\rho_s^2} |\rho_{sx}|^2   -  \f{1}{\rho_s} \rho_{sxx}\big]\Big\} \phi_{y},\label{7.12-1}\\
	I_{v}(\phi)&:= - v_s \phi_{xy} +  \nabla v_s \cdot \nabla^{\perp}\phi  - \f{2 \mu  }{\rho_s^2}\nabla_{\v}\rho_s\cdot \nabla_{\v}\phi_{x} + \f{\lambda}{\rho_s^2} \nabla \rho_s\cdot \nabla^{\perp} \phi_y \nonumber\\
	&\quad + \f{1}{\rho_s} \Big\{(U_s\cdot \nabla)\rho_s   - \f{\mu}{\rho_s} \Delta_{\v}\rho_s +\f{2 \mu}{\rho_s^2} |\nabla_{\v}\rho_s|^2 + \f{\lambda}{\rho_s} \big[\f{2}{\rho_s} \rho_{sy}^2 -  \rho_{syy} \Big]\Big\}  \phi_x \nonumber\\
	&\quad - \f{\lambda}{\rho_s^2}\Big\{\f{2}{\rho_s} \rho_{sy}\rho_{sx}-   \rho_{sxy} \Big\} \phi_y,\label{7.13-1}
\end{align}
\begin{align}
	d_{11}& = - u_s^2 + \f{(2\mu+\lambda)\v}{\rho_s} u_{sx} + \f{\lambda \v}{\rho_s} \mbox{div}U_s - \f{(2\mu+\lambda)\v}{\rho_s^2}\rho_{sx} u_s  - \f{\lambda \v}{\rho_s^2} \nabla \rho_s\cdot U_s\nonumber\\
	&=:-u_s^2 + \v d_{11R},\label{A.4}\\
	d_{12}&= -u_sv_s +  \f{\lambda \v}{\rho_s} v_{sx}  + \f{2\mu}{\rho_s}   u_{sy}   - \f{2\mu}{\rho_s^2} \rho_{sy} u_s  - \f{\lambda \v}{\rho_s^2} \rho_{sx} v_s,\\
	d_{13}&= - (U_s\cdot \nabla)u_s - \f{1}{\rho_s^2} \Big[2\mu \nabla_{\v}\rho_s\cdot \nabla_{\v}u_{s} + \lambda \v \rho_{sx}  \mbox{div}U_s  +  \lambda \v  U_{sx}\cdot \nabla\rho_s\Big]  + \f{\lambda \v}{\rho_s} \mbox{div}U_{sx} \nonumber \\
	&\quad+ \f{\mu}{\rho_s}\Delta_{\v}u_s + \f{u_s}{\rho_s}\Big\{(U_s\cdot \nabla)\rho_s  - \f{\mu}{\rho_s} \Delta_{\v}\rho_s  + \f{2\mu}{\rho_s^2} |\nabla_{\v}\rho_s|^2 + \lambda \v \big[\f{2}{\rho_s^2} |\rho_{sx}|^2   -  \f{1}{\rho_s} \rho_{sxx}\big]\Big\} \nonumber \\
	&\quad +  \frac{\lambda \v  v_s}{\rho_s^2} \Big\{ \f{2}{\rho_s} \rho_{sx}\rho_{sy} -   \rho_{sxy}\Big\} ,\\
d_{21}&:=  - v_su_s + \f{\lambda}{\rho_s} u_{sy} +  \f{2\mu \v}{\rho_s} v_{sx}   - \f{2 \mu\v}{\rho_s^2} v_s \rho_{sx} - \f{\lambda}{\rho_s^2} \rho_{sy} u^{\v},  \\
d_{22}&:= - v_s^2 + \f{2\mu+\lambda}{\rho_s} v_{sy}  + \f{\lambda}{\rho_s}\mbox{div}U_s  - \f{2\mu+\lambda}{\rho_s^2}\rho_{sy}  v_s - \f{\lambda}{\rho_s^2} \nabla \rho_s\cdot U_s,\\
d_{23}&:= -(U_s\cdot \nabla)v_s  - \f{1}{\rho_s^2} \big[2 \mu \nabla_{\v}\rho_s\cdot \nabla_{\v}v_{s} +  \lambda \rho_{sy}  \mbox{div}U_s  +  \lambda  U_{sy}\cdot \nabla \rho_s \big] + \f{\lambda}{\rho_s} \mbox{div}U_{sy}   \nonumber\\
&\quad + \f{\mu}{\rho_s}\Delta_{\v}v_s + \f{v_s}{\rho_s}\Big\{(U_s\cdot \nabla)\rho_s   - \f{\mu}{\rho_s} \Delta_{\v}\rho_s +\f{2 \mu}{\rho_s^2} |\nabla_{\v}\rho_s|^2 + \f{\lambda}{\rho_s} \big[\f{2}{\rho_s} \rho_{sy}^2 -  \rho_{syy} \Big]\Big\}    \nonumber\\
&\quad + \f{\lambda u_s }{\rho_s^2}\Big\{\f{2}{\rho_s} \rho_{sy}\rho_{sx}-   \rho_{sxy} \Big\}.
\end{align}
\begin{align} \label{7.31-1}
\begin{split}
	\tilde{d}_{11}&=1 + d_{11} \frac{1}{T_s},\quad 
	\tilde{d}_{12}= d_{12}  \frac{1}{T_s} ,  \quad 
	\hat{d}_{11}=- d_{11} \frac{\rho_s}{T_s} , \quad 
	\hat{d}_{12} = - d_{12} \frac{\rho_s}{T_s} , \\
	\tilde{d}_{13}&= \Big\{\frac{1}{T_s} d_{13} - d_{11}\f{\pa_x T_s}{T_s^2} - d_{12}\f{\pa_y T_s}{T_s^2}\Big\},\\
	\hat{d}_{13}&= -  \Big\{\frac{\rho_s}{T_s} d_{13} + d_{11}\Big(\f{\rho_s\pa_x T_s}{T_s^2}  - \frac{\pa_x\rho_s}{T_s} \Big) + d_{12}\Big(\f{\rho_s\pa_y T_s}{T_s^2}  - \frac{\pa_y\rho_s}{T_s} \Big)  \Big\}  , 
\end{split}
\end{align}
\begin{align}\label{7.33-1}
\begin{split}
	\tilde{d}_{21}&= d_{21} \frac{1}{T_s},\quad 
	\tilde{d}_{22}= 1+  \v d_{22}  \frac{1}{T_s} ,
	\quad \hat{d}_{21}=- d_{21} \frac{\rho_s}{T_s} , \quad 
	\hat{d}_{22} = - d_{22}  \frac{\rho_s}{T_s},\\
	\tilde{d}_{23}&= \Big\{\frac{1}{T_s} d_{23} - d_{21}\f{\pa_x T_s}{T_s^2} - d_{22}\f{\pa_y T_s}{T_s^2}  \Big\},\\
	\hat{d}_{23}&= -\Big\{\frac{\rho_s}{T_s} d_{23} + d_{21}\Big(\f{\rho_s\pa_x T_s}{T_s^2}  - \frac{\pa_x\rho_s}{T_s} \Big) + d_{22}\Big(\f{\rho_s\pa_y T_s}{T_s^2}  - \frac{\pa_y\rho_s}{T_s} \Big)  \Big\}. 
\end{split}
\end{align}

\begin{align}\label{7.12}
	I^{\v}_{u}(\phi)&:= v_{s}\phi_{yy} + u_{sx}\phi_y + \f{2\mu}{(\rho^{\v})^2} \nabla_{\v}\rho_s\cdot \nabla_{\v}\phi_{y} + \f{\lambda \v}{(\rho^{\v})^2} \nabla \rho_s\cdot \nabla^{\perp} \phi_x + \f{\lambda \v}{(\rho^{\v})^2} \Big\{ \f{2}{\rho^{\v}} \rho_{sx}\rho_{sy} -  \rho_{sxy}\Big\} \phi_x \nonumber\\
	&\quad - \f{1}{\rho^{\v}}\Big\{(U_s\cdot \nabla)\rho_s  - \f{\mu}{\rho^{\v}} \Delta_{\v}\rho_s  + \f{2\mu}{(\rho^{\v})^2} |\nabla_{\v}\rho_s|^2 + \lambda \v \big[\f{2}{(\rho^{\v})^2} |\rho_{sx}|^2   -  \f{1}{\rho^{\v}} \rho_{sxx}\big]\Big\} \phi_{y},\\
	I^{\v}_{v}(\phi)&:= - v_s \phi_{xy} +  \nabla v_s \cdot \nabla^{\perp}\phi  - \f{2 \mu  }{(\rho^{\v})^2}\nabla_{\v}\rho_s\cdot \nabla_{\v}\phi_{x} + \f{\lambda}{(\rho^{\v})^2} \nabla \rho_s\cdot \nabla^{\perp} \phi_y \nonumber\\
	&\quad + \f{1}{\rho^{\v}} \Big\{(U_s\cdot \nabla)\rho_s   - \f{\mu}{\rho^{\v}} \Delta_{\v}\rho_s +\f{2 \mu}{(\rho^{\v})^2} |\nabla_{\v}\rho_s|^2 + \f{\lambda}{\rho^{\v}} \big[\f{2}{\rho^{\v}} \rho_{sy}^2 -  \rho_{syy} \Big]\Big\}  \phi_x \nonumber\\
	&\quad - \f{\lambda}{(\rho^{\v})^2}\Big\{\f{2}{\rho^{\v}} \rho_{sy}\rho_{sx}-   \rho_{sxy} \Big\} \phi_y,\label{7.13}
\end{align}

\begin{align}
		d^{\v}_{11}& = -u^{\v} u_s + \f{(2\mu+\lambda)\v}{\rho^{\v}} u_{sx} + \f{\lambda \v}{\rho^{\v}} \mbox{div}U_s - \f{(2\mu+\lambda)\v}{(\rho^{\v})^2}\rho_{sx} u^{\v}  - \f{\lambda \v}{(\rho^{\v})^2} \nabla \rho_s\cdot U^{\v},\\
		d^{\v}_{12}&= - u^{\v}v_s +  \f{\lambda \v}{\rho^{\v}} v_{sx}  + \f{2\mu}{\rho^{\v}}   u_{sy}   - \f{2\mu}{(\rho^{\v})^2} \rho_{sy} u^{\v}  - \f{\lambda \v}{(\rho^{\v})^2} \rho_{sx} v^{\v},\\
		d^{\v}_{13}&= - (U_s\cdot \nabla)u_s - \f{1}{(\rho^{\v})^2} \Big[2\mu \nabla_{\v}\rho_s\cdot \nabla_{\v}u_{s} + \lambda \v \rho_{sx}  \mbox{div}U_s  +  \lambda \v  U_{sx}\cdot \nabla\rho_s\Big]  + \f{\lambda \v}{\rho^{\v}} \mbox{div}U_{sx} \nonumber \\
		&\quad+ \f{\mu}{\rho^{\v}}\Delta_{\v}u_s + \f{u_s}{\rho^{\v}}\Big\{(U_s\cdot \nabla)\rho_s  - \f{\mu}{\rho^{\v}} \Delta_{\v}\rho_s  + \f{2\mu}{(\rho^{\v})^2} |\nabla_{\v}\rho_s|^2 + \lambda \v \big[\f{2}{(\rho^{\v})^2} |\rho_{sx}|^2   -  \f{1}{\rho^{\v}} \rho_{sxx}\big]\Big\} \nonumber\\
		&\quad +  \frac{\lambda \v  v_s}{(\rho^{\v})^2} \Big\{ \f{2}{\rho^{\v}} \rho_{sx}\rho_{sy} -   \rho_{sxy}\Big\} ,\\
	d^{\v}_{21}&:=  - v^{\v} u_s + \f{\lambda}{\rho^{\v}} u_{sy} +  \f{2\mu \v}{\rho^{\v}} v_{sx}   - \f{2 \mu\v}{(\rho^{\v})^2} v^{\v} \rho_{sx} - \f{\lambda}{(\rho^{\v})^2} \rho_{sy} u^{\v},  \\
	d^{\v}_{22}&:= - v^{\v} v_s + \f{2\mu+\lambda}{\rho^{\v}} v_{sy}  + \f{\lambda}{\rho^{\v}}\mbox{div}U_s  - \f{2\mu+\lambda}{(\rho^{\v})^2}\rho_{sy}  v^{\v}  - \f{\lambda}{(\rho^{\v})^2} \nabla \rho_s\cdot U^{\v},\\
	d^{\v}_{23}&:= -(U_s\cdot \nabla)v_s  - \f{1}{(\rho^{\v})^2} \big[2 \mu \nabla_{\v}\rho_s\cdot \nabla_{\v}v_{s} +  \lambda \rho_{sy}  \mbox{div}U_s  +  \lambda  U_{sy}\cdot \nabla \rho_s \big] + \f{\lambda}{\rho^{\v}} \mbox{div}U_{sy}   \nonumber\\
	&\quad + \f{\mu}{\rho^{\v}}\Delta_{\v}v_s + \f{v_s}{\rho^{\v}}\Big\{(U_s\cdot \nabla)\rho_s   - \f{\mu}{\rho^{\v}} \Delta_{\v}\rho_s +\f{2 \mu}{(\rho^{\v})^2} |\nabla_{\v}\rho_s|^2 + \f{\lambda}{\rho^{\v}} \big[\f{2}{\rho^{\v}} \rho_{sy}^2 -  \rho_{syy} \Big]\Big\}    \nonumber\\
	&\quad + \f{\lambda u_s }{(\rho^{\v})^2}\Big\{\f{2}{\rho^{\v}} \rho_{sy}\rho_{sx}-   \rho_{sxy} \Big\},
\end{align}
\begin{align} \label{7.31-10}
\begin{cases}
	\dis \tilde{d}^{\v}_{11}=1 + d^{\v}_{11} \frac{1}{T^{\v}},\quad 
	\tilde{d}^{\v}_{12}= d^{\v}_{12}  \frac{1}{T^{\v}} ,  \quad 
	\hat{d}^{\v}_{11}=- d^{\v}_{11} \frac{\rho^{\v}}{T^{\v}} , \quad 
	\hat{d}^{\v}_{12} = - d^{\v}_{12} \frac{\rho^{\v}}{T^{\v}} , \\
\dis 	\tilde{d}^{\v}_{13}= \Big\{\frac{1}{T^{\v}} d^{\v}_{13} - d^{\v}_{11}\f{\pa_x T_s}{(T^{\v})^2} - d^{\v}_{12}\f{\pa_y T_s}{(T^{\v})^2}\Big\},\\
\dis 	\hat{d}^{\v}_{13}= -  \Big\{\frac{\rho_s}{T^{\v}} d^{\v}_{13} + d^{\v}_{11}\Big(\f{\rho_s\pa_x T_s}{(T^{\v})^2}  - \frac{\pa_x\rho_s}{T^{\v}} \Big) + d^{\v}_{12}\Big(\f{\rho_s\pa_y T_s}{(T^{\v})^2}  - \frac{\pa_y\rho_s}{T^{\v}} \Big)  \Big\}  , 
\end{cases}
\end{align}
\begin{align}\label{7.33-10}
\begin{cases}
\dis	\tilde{d}^{\v}_{21}= d^{\v}_{21} \frac{1}{T^{\v}},\quad 
	\tilde{d}^{\v}_{22}= 1+  \v d^{\v}_{22}  \frac{1}{T^{\v}} ,
	\quad \hat{d}^{\v}_{21}=- d^{\v}_{21} \frac{\rho^{\v}}{T^{\v}} , \quad 
	\hat{d}^{\v}_{22} = - d^{\v}_{22}  \frac{\rho^{\v}}{T^{\v}},\\
\dis 	\tilde{d}^{\v}_{23}= \Big\{\frac{1}{T^{\v}} d^{\v}_{23} - d^{\v}_{21}\f{\pa_x T_s}{(T^{\v})^2} - d^{\v}_{22}\f{\pa_y T_s}{(T^{\v})^2}  \Big\},\\
\dis 	\hat{d}^{\v}_{23}= -\Big\{\frac{\rho_s}{T^{\v}} d^{\v}_{23} + d_{21}\Big(\f{\rho_s\pa_x T_s}{(T^{\v})^2}  - \frac{\pa_x\rho_s}{T^{\v}} \Big) + d^{\v}_{22}\Big(\f{\rho_s\pa_y T_s}{(T^{\v})^2}  - \frac{\pa_y\rho_s}{T^{\v}} \Big)  \Big\}. 
\end{cases}
\end{align}
and 
\begin{align}
\mathfrak{h}_1:&= \big[I^{\v}_{u}(\phi)-I_{u}(\phi)\big] + \big[\frac{1}{T^{\v}} d^{\v}_{11}  -  \frac{1}{T_s} d_{11}\big] p_x^{(\v)}
+ \big[\tilde{d}^{\v}_{12}-\tilde{d}_{12}\big] p_y^{(\v)} + \big[\tilde{d}^{\v}_{13}-\tilde{d}_{13}\big] p^{(\v)}\nonumber\\
&\quad  + \big[\hat{d}^{\v}_{11}-\hat{d}_{11}\big] T_x^{(\v)}
+ \big[\hat{d}^{\v}_{12}-\hat{d}_{12}\big] T_y^{(\v)} + \big[\hat{d}^{\v}_{13} - \hat{d}_{13}\big] T^{(\v)},\label{7.33-11}\\
\mathfrak{h}_2:&= \big[I^{\v}_{v}(\phi)-I_{v}(\phi)\big] + \big[\tilde{d}^{\v}_{21}-\tilde{d}_{21}\big] p_x^{(\v)}
+ \big[\frac{1}{T^{\v}_s} d^{\v}_{22} - \frac{1}{T_s} d_{22}\big] p_y^{(\v)} + \big[\tilde{d}^{\v}_{23}-\tilde{d}_{23}\big] p^{(\v)}\nonumber\\
&\quad  + \big[\hat{d}^{\v}_{21}-\hat{d}_{21}\big] T_x^{(\v)}
+ \big[\hat{d}^{\v}_{22}-\hat{d}_{22}\big] T_y^{(\v)} + \big[\hat{d}^{\v}_{23} - \hat{d}_{23}\big] T^{(\v)}.\label{7.33-12}
\end{align}
\begin{align}\label{7.36-1.1}
h^{\v}_1:&= 
-  \mu\f{u^{\v}}{\rho^{\v}}   \Big(\f{\rho_s\nabla_{\v} T_s}{(T^{\v})^2}  - \frac{\nabla_{\v}\rho_s}{T^{\v}} -  \nabla_{\v}\big(\frac{\rho^{\v}}{T^{\v}}\big) \Big) \cdot \nabla_{\v}T^{(\v)}  +  \mu   \frac{u^{\v}}{p^{\v}T^{\v}} \nabla_{\v}[T_s+T^{\v}]\cdot \nabla_{\v} p^{(\v)}\nonumber\\
&\quad + \mu \, \mbox{div}_{\v}\Big(\f{\nabla_{\v} T_s}{(T^{\v})^2}\Big)\,  \f{u^{\v} }{\rho^{\v}}  p^{(\v)} -  \mu\, \mbox{div}_{\v}\Big(\f{\rho_s\nabla_{\v}T_s}{(T^{\v})^2}  - \frac{\nabla_{\v}\rho_s}{T^{\v}}\Big)\f{u^{\v}  }{\rho^{\v}}  \,T^{(\v)}  - \f{\lambda \v}{\rho^{\v}} U^{\v}\cdot\nabla\big(\frac{1}{T^{\v}}\big)\, p^{(\v)}_x \nonumber\\
&\quad  + \lambda \v \f{ T_{sx}}{p^{\v}T^{\v}} U^{\v}\cdot \nabla p^{(\v)}  + \f{\lambda \v}{\rho^{\v}} U^{\v}\cdot \nabla\big(\frac{\rho^{\v}}{T^{\v}}\big) \, T^{(\v)}_x - \f{\lambda \v}{\rho^{\v}}\Big(\f{\rho_s T_{sx}}{(T^{\v})^2}  - \frac{\rho_{sx}}{T^{\v}} \Big)\, U^{\v}\cdot\nabla T^{(\v)}\nonumber\\
&\quad + \f{\lambda \v}{\rho^{\v}} U^{\v}\cdot \nabla\Big(\f{T_{sx}}{(T^{\v})^2}\Big)\, p^{(\v)} - \f{\lambda \v}{\rho^{\v}} U^{\v}\cdot\nabla\Big(\f{\rho_sT_{sx}}{(T^{\v})^2}  - \frac{\rho_{sx}}{T^{\v}} \Big)\, T^{(\v)},
\end{align}

\begin{align}\label{7.36-1.2}
\v h^{\v}_2:&=  -  \mu\v \f{v^{\v}}{\rho^{\v}}   \Big(\f{\rho_s\nabla_{\v} T_s}{(T^{\v})^2}  - \frac{\nabla_{\v}\rho_s}{T^{\v}} -  \nabla_{\v}\big(\frac{\rho^{\v}}{T^{\v}}\big) \Big) \cdot \nabla_{\v}T^{(\v)} +  \mu  \v \frac{v^{\v}}{p^{\v}T^{\v}} \nabla_{\v}[T_s+T^{\v}]\cdot \nabla_{\v} p^{(\v)} \nonumber\\
&\quad + \mu \v\, \mbox{div}_{\v}\Big(\f{\nabla_{\v} T_s}{(T^{\v})^2}\Big)\,  \f{v^{\v} }{\rho^{\v}}  p^{(\v)} -  \mu\v\, \mbox{div}_{\v}\Big(\f{\rho_s\nabla_{\v}T_s}{(T^{\v})^2}  - \frac{\nabla_{\v}\rho_s}{T^{\v}}\Big)\f{v^{\v}  }{\rho^{\v}}  \,T^{(\v)}  - \f{\lambda \v}{\rho^{\v}} U^{\v}\cdot\nabla\big(\frac{1}{T^{\v}}\big)\, p^{(\v)}_y \nonumber\\
&\quad   + \lambda \v \f{ T_{sy}}{p^{\v}T^{\v}} U^{\v}\cdot \nabla p^{(\v)}  + \f{\lambda \v}{\rho^{\v}} U^{\v}\cdot \nabla\big(\frac{\rho^{\v}}{T^{\v}}\big) \, T^{(\v)}_y - \f{\lambda \v}{\rho^{\v}}\Big(\f{\rho_s T_{sy}}{(T^{\v})^2}  - \frac{\rho_{sy}}{T^{\v}} \Big)\, U^{\v}\cdot\nabla T^{(\v)} \nonumber\\
&\quad+ \f{\lambda \v}{\rho^{\v}} U^{\v}\cdot \nabla\Big(\f{T_{sy}}{(T^{\v})^2}\Big)\, p^{(\v)} - \f{\lambda \v}{\rho^{\v}} U^{\v}\cdot\nabla\Big(\f{\rho_sT_{sy}}{(T^{\v})^2}  - \frac{\rho_{sy}}{T^{\v}} \Big)\, T^{(\v)}.
\end{align}

\begin{align}
	h_1:&=-  \mu\f{u_s}{\rho_s}   \Big(\f{\rho_s\nabla_{\v} T_s}{(T_s)^2}  - \frac{\nabla_{\v}\rho_s}{T_s} -  \nabla_{\v}\big(\frac{\rho_s}{T_s}\big) \Big) \cdot \nabla_{\v}T^{(\v)}  +  2\mu   \frac{u_s}{p_sT_s} \nabla_{\v}T_s\cdot \nabla_{\v} p^{(\v)}\nonumber\\
	&\quad + \mu \, \mbox{div}_{\v}\Big(\f{\nabla_{\v} T_s}{(T_s)^2}\Big)\,  \f{u_s}{\rho_s}  p^{(\v)} -  \mu\, \mbox{div}_{\v}\Big(\f{\rho_s\nabla_{\v}T_s}{(T_s)^2}  - \frac{\nabla_{\v}\rho_s}{T_s}\Big)\f{u_s}{\rho_s}  \,T^{(\v)}  - \f{\lambda \v}{\rho_s} U_s\cdot\nabla\big(\frac{1}{T_s}\big)\, p^{(\v)}_x \nonumber\\
	&\quad  + \lambda \v \f{ T_{sx}}{p_sT_s} U_s\cdot \nabla p^{(\v)}  + \f{\lambda \v}{\rho_s} U_s\cdot \nabla\big(\frac{\rho_s}{T_s}\big) \, T^{(\v)}_x - \f{\lambda \v}{\rho_s}\Big(\f{\rho_s T_{sx}}{(T_s)^2}  - \frac{\rho_{sx}}{T_s} \Big)\, U_s\cdot\nabla T^{(\v)}\nonumber\\
	&\quad + \f{\lambda \v}{\rho_s} U_s\cdot \nabla\Big(\f{T_{sx}}{(T_s)^2}\Big)\, p^{(\v)} - \f{\lambda \v}{\rho_s} U_s\cdot\nabla\Big(\f{\rho_sT_{sx}}{(T_s)^2}  - \frac{\rho_{sx}}{T_s} \Big)\, T^{(\v)},
\end{align}

\begin{align} 
	\v h_2:&=  -  \mu\v \f{v_s}{\rho_s}   \Big(\f{\rho_s\nabla_{\v} T_s}{(T_s)^2}  - \frac{\nabla_{\v}\rho_s}{T_s} -  \nabla_{\v}\big(\frac{\rho_s}{T_s}\big) \Big) \cdot \nabla_{\v}T^{(\v)} +  2\mu  \v \frac{v_s}{p_sT_s} \nabla_{\v}T_s\cdot \nabla_{\v} p^{(\v)} \nonumber\\
	&\quad + \mu \v\, \mbox{div}_{\v}\Big(\f{\nabla_{\v} T_s}{(T_s)^2}\Big)\,  \f{v_s }{\rho_s}  p^{(\v)} -  \mu\v\, \mbox{div}_{\v}\Big(\f{\rho_s\nabla_{\v}T_s}{(T_s)^2}  - \frac{\nabla_{\v}\rho_s}{T_s}\Big)\f{v_s}{\rho_s}  \,T^{(\v)}  - \f{\lambda \v}{\rho_s} U_s\cdot\nabla\big(\frac{1}{T_s}\big)\, p^{(\v)}_y \nonumber\\
	&\quad   + \lambda \v \f{ T_{sy}}{p_sT_s} U_s\cdot \nabla p^{(\v)}  + \f{\lambda \v}{\rho_s} U_s\cdot \nabla\big(\frac{\rho_s}{T_s}\big) \, T^{(\v)}_y - \f{\lambda \v}{\rho_s}\Big(\f{\rho_s T_{sy}}{(T_s)^2}  - \frac{\rho_{sy}}{T_s} \Big)\, U_s\cdot\nabla T^{(\v)} \nonumber\\
	&\quad+ \f{\lambda \v}{\rho_s} U_s\cdot \nabla\Big(\f{T_{sy}}{(T_s)^2}\Big)\, p^{(\v)} - \f{\lambda \v}{\rho_s} U_s\cdot\nabla\Big(\f{\rho_sT_{sy}}{(T_s)^2}  - \frac{\rho_{sy}}{T_s} \Big)\, T^{(\v)}.
\end{align}
%

\begin{align}
	\mathfrak{N}_{1} &:=  R_{u}^{(\v)}   - \v^{N_0} \Big\{\f{2\mu}{\rho^{\v}}\nabla_{\v}\rho^{(\v)}\cdot \nabla_{\v} u^{(\v)} + \f{\lambda \v }{\rho^{\v}} U^{(\v)}_x\cdot \nabla\rho^{(\v)} + \f{\lambda \v }{\rho^{\v}}\rho_x^{(\v)} \mbox{div}U^{(\v)} \Big\} + [h^{\v}_1- h_1],\label{9.26}\\
	\mathfrak{R}_{s1}&:=   - \v^{-N_0} B_{s1} + \v^{-N_0} \Big\{ \f{1}{(\rho^{\v})^2} \big[2\mu\nabla_{\v}\rho_s\cdot\nabla_{\v}A_{s1} + \lambda \v \nabla \rho_s\cdot A_{sx} + \lambda \v \rho_{sx}  \mbox{div}A_{s}\big]   -  \f{\mu}{\rho^{\v}}\Delta_{\v}A_{s1} \nonumber\\
	&\quad - \f{\lambda \v}{\rho^{\v}} \mbox{div}A_{sx}   + A_{s}\cdot \nabla u_s +  (U_s\cdot\nabla) A_{s1}  \Big\} - \f{\lambda\v^{1-N_0}}{(\rho^{\v})^2}  A_{s2} \Big\{ \f{2}{\rho^{\v}} \rho_{sx}\rho_{sy} - \rho_{sxy}\Big\} \nonumber\\
	&\quad - \f{\v^{-N_0}}{\rho^{\v}} A_{s1}\Big\{(U_s\cdot \nabla)\rho_s  - \f{\mu}{\rho^{\v}} \Delta_{\v}\rho_s  + \f{2\mu}{(\rho^{\v})^2} |\nabla_{\v}\rho_s|^2 + \lambda \v \big[\f{2}{(\rho^{\v})^2} |\rho_{sx}|^2   -  \f{1}{\rho^{\v}} \rho_{sxx}\big]\Big\}.
\end{align}
and
\begin{align}
	\mathfrak{N}_2&:= R_{v}^{(\v)} - \f{\v^{N_0}}{\rho^{\v}} \Big[2 \mu\nabla_{\v}\rho^{(\v)}\cdot \nabla_{\v} v^{(\v)} + \lambda  (U^{(\v)}_y\cdot \nabla)\rho^{(\v)} + \lambda  \rho_y^{(\v)} \mbox{div}U^{(\v)} \Big] + \v [h^{\v}_2-h_2],\label{9.28}\\
	\mathfrak{R}_{s2}&:=   -  \v^{-N_0}  B_{s2} + \v^{-N_0} \Big\{\f{1}{(\rho^{\v})^2} \big[ 2 \mu \nabla_{\v}\rho_s\cdot\nabla_{\v}A_{s2}  +  \lambda   \nabla \rho_s\cdot A_{sy} + \lambda \rho_{sy} \mbox{div}A_{s}\big]  - \f{\mu}{\rho^{\v}}\Delta_{\v}A_{s2} \nonumber\\
	&\quad - \f{\lambda}{\rho^{\v}}\mbox{div}A_{sy}    +  (U_s\cdot\nabla) A_{s2} + A_{s}\cdot \nabla v_s\Big\}   - \f{\lambda \v^{-N_0} }{(\rho^{\v})^2} A_{s1}  \Big\{\f{2}{\rho^{\v}} \rho_{sy}\rho_{sx}-   \rho_{sxy} \Big\}  \nonumber\\
	&\quad - \f{\v^{-N_0}}{\rho^{\v}}  A_{s2} \Big\{(U_s\cdot \nabla)\rho_s   - \f{\mu}{\rho^{\v}} \Delta_{\v}\rho_s +\f{2 \mu}{(\rho^{\v})^2} |\nabla_{\v}\rho_s|^2 + \f{\lambda}{\rho^{\v}} \big[\f{2}{\rho^{\v}} \rho_{sy}^2 -  \rho_{syy} \Big]\Big\},
\end{align}

\begin{remark}
We remark that the linear term $h_1$ can be absorbed into $\tilde{d}_{11} p_x^{(\v)}
+ \tilde{d}_{12} p_y^{(\v)} + \tilde{d}_{13} p^{(\v)} + \hat{d}_{11} T_x^{(\v)}
+ \hat{d}_{12} T_y^{(\v)} + \hat{d}_{13} T^{(\v)} $. And 
 $\v h_2$ can be absorbed into $\v \tilde{d}_{21} p_x^{(\v)}
 +   \tilde{d}_{22} p_y^{(\v)} + \v \tilde{d}_{23} p^{(\v)}  + \v \hat{d}_{21} T_x^{(\v)}
 + \v \hat{d}_{22} T_y^{(\v)} + \v \hat{d}_{23} T^{(\v)}$. Then one needs only add some terms into  $\tilde{d}_{ij}, \hat{d}_{i,j} i=1,2; j=1,2,3$, but the properties are completely similar, so we do not modify the notations  $\tilde{d}_{ij}, \hat{d}_{i,j} i=1,2; j=1,2,3$ for simplicity, and keep in mind that $h_1, \v h_2$ have already been absorbed.
\end{remark}

The term $\la (\hat{d}_{11}T_x)_y, q_{xx} \ra$ is hard to control, to deal with this term,  by introducing the pseudo entropy \eqref{s}, one obtains
\begin{align}\label{7.33-4}
	\hat{d}_{11} T_x = \hat{d}_{11} S_x + \f{1}{2\rho_s}\hat{d}_{11} p_x + \hat{d}_{11}\f{1}{p_s}T_{sy} q_x + \hat{d}_{11}(\f{1}{2\rho_s})_x p  + \hat{d}_{11} \big(\f{1}{p_s}T_{sy}\big)_x q.
\end{align}
Applying \eqref{7.33-4} to \eqref{7.36}, we have
\begin{align}\label{2.87}
	\begin{cases}
		\dis \big[u_s \phi_{xy} - u_{sy} \phi_x\big]   + \mathbf{d}_{11}\, p_x + u_s^2\f{1}{T_s^2}T_{sy} q_x+ g_1  + \hat{d}_{11} S_x - \v d_{11R}\f{1}{T_s^2}T_{sy} q_x    + \hat{d}_{11} \big(\f{1}{p_s}T_{sy}\big)_x q \\[2mm]
		\dis \quad  =\frac{\mu}{\rho^{\v}} \Delta_{\v}\phi_y - \frac{\mu+\lambda}{p^{\v}} u^{\v} \Delta_{\v}p  + \frac{\lambda}{p^{\v}} u^{\v} p_{yy} - \frac{\lambda\v}{p^{\v}} v^{\v} p_{xy}  - \frac{\lambda}{T^{\v}} u^{\v} T_{yy} + G_{11}(T) + \mathcal{N}_1,\\[3mm]
		\dis -\v u_s \phi_{xx} +   \tilde{d}_{22}\, p_y + g_2 \\
		\dis \quad  =-\frac{\mu\v }{\rho^{\v}}   \Delta_{\v}\phi_x - \frac{\mu\v}{p^{\v}} v^{\v} \Delta_{\v}p  - \frac{\lambda\v}{p^{\v}} u^{\v} p_{xy}  - \frac{\lambda\v}{p^{\v}} v^{\v} p_{yy}  + \frac{\lambda\v}{T^{\v}} u^{\v} T_{xy}  + \v G_{21}(T) + \v \mathcal{N}_2,
	\end{cases}
\end{align}
where we have used the notations  
\begin{align}\label{2.82}
	\begin{split}
		\mathbf{d}_{11}:&=\tilde{d}_{11} + \f{1}{2\rho_s} \hat{d}_{11} \equiv 1 + \f{1}{2T_s} d_{11} = 1 -\f{u_s^2}{2T_s}  + \v\f{1}{2T_s} d_{11R},\\
		g_1(p,T,\phi):&= \tilde{d}_{12}\, p_y  + \big[\tilde{d}_{13} + \big(\f{1}{2\rho_s}\hat{d}_{11}\big)_x \big]\, p + \hat{d}_{12}\, T_y  + \hat{d}_{13}\, T + I_u(\phi) - \mathfrak{R}_{s1},  
		\\
		g_2(p,T,\phi):&= \v \tilde{d}_{21}\, p_x + \v \tilde{d}_{23}\, p   + \v \hat{d}_{21}\, T_x 	+ \v \hat{d}_{22}\, T_y + \v \hat{d}_{23}\, T + \v I_{v}(\phi) - \v  \mathfrak{R}_{s2},  
	\end{split}
\end{align}
and
\begin{align}\label{2.83}
	\begin{split}
		G_{11}(T):&= \frac{\mu+\lambda}{T^{\v}} u^{\v} \Delta_{\v}T + \frac{\lambda\v}{T^{\v}} v^{\v} T_{xy},\\
		G_{21}(T):&= \frac{\mu}{T^{\v}} v^{\v} \Delta_{\v}T  + \frac{\lambda}{T^{\v}}v^{\v} T_{yy},\\
		\mathcal{N}_1:&= \mathfrak{h}_i + \mathfrak{N}_i,\,\,\, i=1,2. 
	\end{split}
\end{align}

In general, it is still very hard to control the term $\dis u_s^2\f{m}{T_s}T_{sy} q_x$ on LHS of \eqref{2.87} for $q_{xx}$ multiplier.  A direct calculation shows that 
\begin{align}\label{3.4}
	\begin{split}
		\big[u_s \phi_{xy} - u_{sy} \phi_x\big]&= m u_s^2 q_{xy} + u_s\bar{u}_{sx}q_y + [u_s\bar{u}_{sxy} - u_{sy}\bar{u}_{sx}]q +  u_s^2 m_y q_x,\\
		u_s \phi_{xx}&=m u_s^2 q_{xx} + 2u_s \bar{u}_{sx} q_x + u_s \bar{u}_{sxx} q.
	\end{split}
\end{align}
Fortunately, by taking $m:=\f{1}{T_s}$, the bad term on $q_x$ can be canceled, i.e., 
\begin{align}\label{2.91}
	u_s^2 m_y q_x + u_s^2\f{m}{T_s}T_{sy} q_x = u_s^2 q_x \big[m_y + m \f{T_{sy}}{T_s}\big]\equiv 0.
\end{align}

Therefore, based on above calculations, we further reduce the momentum equations  $\eqref{1.8-1}_{2,3}$ as 
\begin{align}\label{2.92-0}
	\begin{cases}
		\dis m u_s^2 q_{xy}  + \mathbf{d}_{11}\, p_x+ \hat{d}_{11} S_x  + g_1(p,\phi, T(S,p,q))  + \mathfrak{g}_1(q) \\
		\dis\,\, =\frac{\mu}{\rho^{\v}} \Delta_{\v}\phi_y - \frac{\mu+\lambda}{p^{\v}} u^{\v} \Delta_{\v}p  + \frac{\lambda}{p^{\v}} \big[1- \f{\rho^{\v}}{2\rho_s}\big] u^{\v} p_{yy} - \frac{\lambda\v}{p^{\v}} v^{\v} p_{xy}  - \frac{\lambda}{T^{\v}} u^{\v} S_{yy} + G_1 + \mathcal{N}_1,\\[3mm]
		\dis -\v m u_s^2 q_{xx} +   \tilde{d}_{22}\, p_y + g_2(p,\phi, T(S,p,q)) + \v \mathfrak{g}_2(q)\\
		\dis \,\, =-\frac{\mu\v }{\rho^{\v}}   \Delta_{\v}\phi_x - \frac{\mu\v}{p^{\v}} v^{\v} \Delta_{\v}p  - \frac{\lambda\v}{p^{\v}} \big[1-\f{\rho^{\v}}{2\rho_s}\big] u^{\v} p_{xy}  - \frac{\lambda\v}{p^{\v}} v^{\v} p_{yy} + \frac{\lambda\v}{T^{\v}} u^{\v} S_{xy}  + \v G_2 + \v \mathcal{N}_2,
	\end{cases}
\end{align}
where we have denoted 
\begin{align}\label{2.93}
	\begin{split}
		\mathfrak{g}_1(q):&= u_s\bar{u}_{sx}q_y + [u_s\bar{u}_{sxy} - u_{sy}\bar{u}_{sx}]q - \v d_{11R}\f{m}{T_s}T_{sy} q_x    + \hat{d}_{11} \big(\f{1}{p_s}T_{sy}\big)_x q, \\
		\mathfrak{g}_2(q):&=- 2  u_s \bar{u}_{sx} q_x -  u_s \bar{u}_{sxx} q.
	\end{split}
\end{align}
and $G_1:=G_{11}+G_{12},\,\,\, G_{2}:=G_{21} + G_{22}$ with
\begin{align}\label{2.96-1}
	\begin{split}
		G_{12}:&=- \frac{\lambda}{T^{\v}} u^{\v} \Big\{\big(\f1{2\rho_s}\big)_{yy}p + \big(\f1{\rho_s}\big)_yp_y + \big(\f{1}{p_s}T_{sy} q\big)_{yy}\Big\},\\
		G_{22}:&= \frac{\lambda }{T^{\v}} u^{\v}  \Big\{\big(\f1{2\rho_s}\big)_{xy}p + \big(\f1{2\rho_s}\big)_xp_y + \big(\f1{2\rho_s}\big)_yp_x +  \big(\f{1}{p_s}T_{sy} q\big)_{xy}\Big\}.
	\end{split}
\end{align}

\
 
\subsection{Reformulation for temperature equation}\label{subsec10.3}
For later use, we rewrite \eqref{7.1-1} as
\begin{align}\label{7.1}
	\begin{split}
		\rho_s u^{(\v)}&= \phi_y -  u_s \rho^{(\v)} + N_1,\\
		\rho_s v^{(\v)}&=  -\phi_x -  v_s \rho^{(\v)} + N_2,
	\end{split}
\end{align}
where 
\begin{align*}
	\begin{split}
		N_{1}&= - \v^{-N_0} A_{s1} - \v^{N_0} \rho^{(\v)} u^{(\v)},\\
		N_{2}&= - \v^{-N_0} A_{s2} - \v^{N_0} \rho^{(\v)} v^{(\v)}.
	\end{split}
\end{align*}
It follows from \eqref{7.1} that
\begin{align}\label{7.33-2.0}
	2\rho_s (U^{(\v)}\cdot \nabla)T_s 
	&= 2T_{sx} \phi_y - 2T_{sy} \phi_x - 2 (U_s\cdot\nabla)T_s \, \rho^{(\v)} + 2  [T_{sx}N_1 + T_{sy}N_2]\nonumber\\
	&=2T_{sx} \phi_y - 2T_{sy} \phi_x - 2 (U_s\cdot\nabla)T_s \, \rho^{(\v)} + \mathfrak{N}_{31},
\end{align}
\begin{align}\label{7.33-2}
	-(U^{(\v)}\cdot \nabla) p_s 
	& = \f{1}{\rho_s} \big[ p_{sy} \phi_x - p_{sx} \phi_y\big]  + \f{1}{\rho_s}(U_s\cdot \nabla)p_s \,\rho^{(\v)} - \f{1}{\rho_s}[p_{sx}N_1 + p_{sy}N_2]\nonumber\\
	&=\f{1}{\rho_s} \big[ p_{sy} \phi_x - p_{sx} \phi_y\big]  + \f{1}{\rho_s}(U_s\cdot \nabla)p_s \,[\f1{T_s}p^{(\v)}-\f{\rho_s}{T_s}T^{(\v)}] + \mathfrak{N}_{32},
\end{align}
\begin{align}\label{7.34-2}
	&4\mu \mathbb{S}^{\v}(U_s):\mathbb{S}^{\v}(U^{(\v)}) + 2 (\lambda-\mu) \v \mbox{div}U_s\cdot \mbox{div}U^{(\v)}\nonumber\\
	&= \f{2\mu}{\rho_s} u_{sy}[\phi_{yy} - \frac{\rho_{sy}}{\rho_s} \phi_y ]- \f{2\mu}{\rho_s} u_{sy} u_s\,  [\f1{T_s}p^{(\v)}_y -\f{\rho_s}{T_s}T^{(\v)}_y ] - \f{2\mu}{\rho_s} u_{sy} u_s\,  [(\f1{T_s})_y p^{(\v)}- (\f{\rho_s}{T_s})_y T^{(\v)}]\nonumber\\
	&\quad  - \f{2\mu}{\rho_s} u_{sy} [u_{sy}-\frac{\rho_{sy}}{\rho_s}u_s]\,[\f1{T_s}p^{(\v)}-\f{\rho_s}{T_s}T^{(\v)}]  
	+ \v \mathfrak{M}_1 + \v \mathfrak{M}_2 + \mathfrak{N}_{33}
\end{align}
where have used the following notations
\begin{align}
	\begin{split}
		\mathfrak{N}_{31}:&= 2  [T_{sx}N_1 + T_{sy}N_2],\\
		\mathfrak{N}_{32}:&= \f{1}{\rho_s}(U_s\cdot \nabla)p_s \,N_{\rho} - \f{1}{\rho_s}[p_{sx}N_1 + p_{sy}N_2],
	\end{split}
\end{align}
\begin{align}\label{7.35-2}
	\mathfrak{M}_1:&= \f{2\mu}{\rho_s}\Big\{ v_{sx}[\phi_{yy} - \frac{\rho_{sy}}{\rho_s} \phi_y ]  + 2u_{sx} [\phi_{xy} - \frac{\rho_{sx}}{\rho_s} \phi_y ] - 2 v_{sy}   [\phi_{xy} - \frac{\rho_{sy}}{\rho_s} \phi_x]\nonumber\\
	&\quad - [u_{sy}+\v v_{sx}] [\phi_{xx} - \frac{\rho_{sx}}{\rho_s} \phi_x] \Big\}  - \f{2\mu}{\rho_s} [v_{sx}u_s + 2v_{sy} v_s]  [\f1{T_s}p^{(\v)}-\f{\rho_s}{T_s}T^{(\v)}]_y \nonumber\\
	&\quad -\f{2\mu}{\rho_s} \Big\{2u_{sx} u_s  +  [u_{sy}+\v v_{sx}]  v_s \Big\} [\f1{T_s}p^{(\v)}-\f{\rho_s}{T_s}T^{(\v)}]_x  \nonumber\\
	&\quad - \f{2\mu}{\rho_s} \Big\{ v_{sx} [u_{sy}-\frac{\rho_{sy}}{\rho_s}u_s] + 2u_{sx}  [u_{sx}-\frac{\rho_{sx}}{\rho_s}u_s] + 2 v_{sy}  [v_{sy}-\frac{\rho_{sy}}{\rho_s}v_s] \nonumber\\
	&\qquad\qquad + [u_{sy}+\v v_{sx}]  [v_{sx}-\frac{\rho_{sx}}{\rho_s}v_s]\Big\} [\f1{T_s}p^{(\v)}-\f{\rho_s}{T_s}T^{(\v)}],
\end{align}
\begin{align}\label{7.36-2}
	\mathfrak{M}_2:&= 2 \f{\lambda-\mu}{\rho_s}  \mbox{div}U_s  \Big\{ \f{\rho_s}{T_s} (U_s\cdot \nabla)T^{(\v)}   - \f1{T_s}(U_s\cdot \nabla)p^{(\v)} - (U_s\cdot \nabla)(\f1{T_s})\, p^{(\v)}   + (U_s\cdot \nabla)(\f{\rho_s}{T_s}) \, T^{(\v)}   \nonumber\\
	&\quad + [\frac{1}{\rho_s} (U_s\cdot \nabla)\rho_s -\mbox{div}U_s]\,[\f1{T_s}p^{(\v)} - \f{\rho_s}{T_s}T^{(\v)}]  + [\frac{\rho_{sy}}{\rho_s} \phi_x- \frac{\rho_{sx}}{\rho_s} \phi_y]\Big\} 
\end{align}
and
\begin{align}\label{7.36-3}
	\mathfrak{N}_{33}:&=  \f{2\mu}{\rho_s} [u_{sy}+\v v_{sx}] \Big(-u_s \pa_y N_{\rho} - [u_{sy} -\f{\rho_{sy}}{\rho_s} u_s] N_{\rho} + [\pa_y N_1 - \frac{\rho_{sy}}{\rho_s} N_1]\Big) \nonumber\\
	&\quad +  \v \f{1}{\rho_s} [4\mu u_{sx} + 2(\lambda-\mu)\mbox{div}U_s] \Big(-u_s \pa_x N_{\rho} - [u_{sx} -\f{\rho_{sx}}{\rho_s} u_s] N_{\rho} + [\pa_x N_1 - \frac{\rho_{sx}}{\rho_s} N_1]\Big) \nonumber\\
	&\quad  + \v \f{1}{\rho_s}[4\mu v_{sy} + 2 (\lambda-\mu)   \mbox{div}U_s] \Big(- v_s \pa_y N_{\rho} -[v_{sy}-\frac{\rho_{sy}}{\rho_s}v_s] N_{\rho} + [\pa_y N_2 - \frac{\rho_{sy}}{\rho_s}N_2]\Big) \nonumber\\
	&\quad +  \v \f{2\mu }{\rho_s}[u_{sy}+\v v_{sx}]  \Big(- v_s \pa_x N_{\rho} -[v_{sx}-\frac{\rho_{sx}}{\rho_s}v_s] N_{\rho} + [\pa_x N_2 - \frac{\rho_{sx}}{\rho_s}N_2]\Big).
\end{align}
We also denote 
\begin{align}\label{7.47}
\mathfrak{M}_0:&= - \f{2\mu}{\rho_s} u_{sy} u_s\,  [\f1{T_s}p^{(\v)}_y -\f{\rho_s}{T_s}T^{(\v)}_y ] - \f{2\mu}{\rho_s} u_{sy} u_s\,  [(\f1{T_s})_y p^{(\v)}- (\f{\rho_s}{T_s})_y T^{(\v)}] \nonumber\\
&\quad - \Big\{\f{1}{\rho_s}(U_s\cdot \nabla)p_s + \f{2\mu}{\rho_s} u_{sy} [u_{sy}-\frac{\rho_{sy}}{\rho_s}u_s]\Big\}\,[\f1{T_s}p^{(\v)}-\f{\rho_s}{T_s}T^{(\v)}] \nonumber\\
&\quad + \f{2\mu}{\rho_s} u_{sy}[\phi_{yy} - \frac{\rho_{sy}}{\rho_s} \phi_y ] + \f{1}{\rho_s} [p_{sx}-2\rho_s T_{sx}] \phi_y.
\end{align}

Substituting \eqref{7.33-2}-\eqref{7.36-2} into $\eqref{1.8-1}_4$, one obtains
\begin{align}\label{7.46}
		&2\rho_s (U_s\cdot\nabla)T^{(\v)}  - (U_s\cdot\nabla)p^{(\v)}    + \f{1}{\rho_s}[p_{sy}-2\rho_s T_{sy}]\phi_x  = \kappa \Delta_{\v}T^{(\v)} + \mathfrak{M} + \mathfrak{N}_3,
\end{align}
where 
\begin{align}\label{7.46-1}
	\begin{split}
		\mathfrak{M}:&= \mathfrak{M}_0 + \v \mathfrak{M}_1 + \v \mathfrak{M}_2,\\
		\mathfrak{N}_3:&=\mathfrak{N}_{33}  - \mathfrak{N}_{31}  - \mathfrak{N}_{32}  + R^{(\v)}_{T} - \v^{-N_0}B_{s3}.
	\end{split}
\end{align}
We point out that all the nonlinear terms in \eqref{7.46} are lower order derivatives, and the term $\dis\f{1}{\rho_s}[p_{sy}-2\rho_s T_{sy}]\phi_x $ will cause serious difficulty.

Due to the lack of $\sqrt{\v}$, the terms $u_s p_x$ and $T_{sy}\phi_x$  are hard to control when we apply the multiplier $T_{xx}$, which are crucial to obtain uniform estimate on $\|(\sqrt{\v}T_{xx}, T_{xy})\|^2$. The key point is to  introduce the {\it new} quantity, i.e., {\it Pseudo entropy} (see \eqref{s}), to deal with these terms.


In fact, noting \eqref{s}, we have
\begin{align}\label{2.80}
	&2\rho_s (U_s\cdot \nabla) T - (U_s\cdot \nabla) p  -  2 T_{sy} \phi_x \nonumber\\
	&=2\rho_s u_s \Big(T-\frac{1}{2\rho_s}p - \f{1}{p_s}T_{sy}q\Big)_x + 2\rho_s v_s \big(T-\frac{1}{2\rho_s}p\big)_y  + \rho_s (U_s\cdot \nabla) (\frac{1}{\rho_s})\, p \nonumber\\
	&\quad + 2\big[\rho_su_s \big(\f{1}{p_s}T_{sy}\big)_x-  T_{sy} \bar{u}_{sx}\big]q.
\end{align}
Noting the {\it Pseudo entropy} in \eqref{s},  substituting \eqref{2.80} into \eqref{7.46}, one can rewrite the temperature equation $\eqref{1.8-1}_4$ as
\begin{align}\label{T.1-0}
	2 \rho_s u_s S_x + 2\rho_s v_s S_y - \kappa \Delta_{\v}S
	&= \kappa \Delta_{\v}\big(\frac{1}{2\rho_s}p\big) + \kappa \Delta_{\v}\big(\f{1}{p_s}T_{sy}q\big) + \mathfrak{J}(\phi,p,T),
\end{align}
where 
\begin{align}\label{2.84}
	\mathfrak{J}(\phi,p,T):&=-  2\rho_sv_s \big(\f{1}{p_s}T_{sy}q\big)_y -  2\big[\rho_su_s \big(\f{1}{p_s}T_{sy}\big)_x-  T_{sy} \bar{u}_{sx}\big]q  - \f{1}{\rho_s}p_{sy}  \phi_x  \nonumber\\
	&\quad  -  \rho_s (U_s\cdot \nabla) (\frac{1}{\rho_s})\, p + \mathfrak{M} + \mathfrak{N}_3. 
\end{align}

\subsection{Some notations}\label{appA.3}
We define following notations  for \eqref{2.112}:
\begin{align}\label{2.98}
	\begin{split}
		\mathcal{F}_{R1}:&=-\big(b_1 \,\hat{d}_{11} \zeta_x\big)_y,  \\
		\mathcal{F}_{R2}:&=  \v \big(b_2\, \mathfrak{g}_2\big)_x  - \big(b_1\mathfrak{g}_1\big)_y,\\
		\mathcal{F}_{R3}:&=\big(b_2\,g_2(p,\phi, T(S,p,q))\big)_x - \big( b_1 g_1(p,\phi, T(S,p,q))\big)_y,\\
		\mathcal{F}_{R4}:&=\pa_y\big(b_1G_{11}(\tilde{T}(\zeta,p,q))\big)  - \v \pa_x\big(b_2 G_{21}(\tilde{T}(\zeta,p,q))\big),\\
		\mathcal{F}_{R5}:&=\pa_y\big(b_1G_{12}(p,q)\big)  - \v \pa_x\big(b_2 G_{22}(p,q)\big),\\
		\mathcal{F}_{R8}:&=\pa_y\big(b_1\mathcal{N}_1\big)  - \v \pa_x\big(b_2 \mathcal{N}_2\big),
	\end{split}
\end{align} 
and 
\begin{align}\label{2.98-1}
	\begin{split}
		\mathcal{F}_{R6}:&= - \lambda\Big(\frac{b_1}{T^{\v}} u^{\v} \mathbf{S}_{yy}\Big)_y - \lambda\v\Big(\frac{b_2}{T^{\v}} u^{\v} \mathbf{S}_{xy}\Big)_x,\\
		\mathcal{F}_{R7}:&=\mu\v \Big(\frac{b_2}{p^{\v}} v^{\v} \Delta_{\v}P\Big)_x - \mu \Big(\frac{b_1}{p^{\v}} u^{\v} \Delta_{\v}P\Big)_y - \lambda \Big(\frac{b_1\rho^{\v}}{2\rho_sp^{\v}} u^{\v} P_{yy}\Big)_y +  \lambda \v \Big(\frac{b_2\rho^{\v}}{2\rho_s p^{\v}} u^{\v} P_{xy}\Big)_x\\
		&\quad  + \lambda\v\Big(\frac{b_2}{p^{\v}} (U^{\v}\cdot\nabla)P_y \Big)_x
		- \lambda\v \Big(\frac{b_1}{p^{\v}} (U^{\v}\cdot\nabla)P_x\Big)_y.
	\end{split}
\end{align}

On the other hand, , we define the following notations  for \eqref{2.115}:
\begin{align}\label{17.17-1}
	\begin{split}
		\mathcal{G}_1:&=  - \v \mathbf{d}_{11x} \, P_x - \tilde{d}_{22y}\, P_y - \f{1}{2T_s}d_{11}  \v P_{xx} - \v \f{1}{T_s}d_{22} P_{yy} - (\mu+\lambda)\v \big(\frac{u^{\v}}{p^{\v}}\big)_x  \Delta_{\v}P \\
		&\quad  - \mu\v \big(\frac{v^{\v}}{p^{\v}}\big)_y  \Delta_{\v}P + \lambda \v \big(\frac{u^{\v}}{p^{\v}}\big[1-\f{\rho^{\v}}{2\rho_s}\big]\big)_x P_{yy} - \lambda\v \big(\frac{v^{\v}}{p^{\v}}\big)_y  P_{yy}  - \lambda\v^2 \big(\frac{v^{\v}}{p^{\v}}\big)_x P_{xy} \\
		&\quad  - \lambda\v \big(\frac{u^{\v}}{p^{\v}}\big[1-\f{\rho^{\v}}{2\rho_s}\big]\big)_y P_{xy} ,\\
		\mathcal{G}_2:&=\v \pa_x G_{11}(\tilde{T}(\zeta,p,q)) + \v \pa_y G_{21}(\tilde{T}(\zeta,p,q)) \\
		&=\v (\mu+\lambda) (\frac{u^{\v}}{T^{\v}} \cdot \nabla)\Delta_{\v}T(\zeta,p,q) +  \Big\{\v (\mu+\lambda) \big(\frac{u^{\v}}{T^{\v}}\big)_x  \Delta_{\v} + \lambda\v^2 \big(\frac{v^{\v}}{T^{\v}}\big)_x \pa_{xy}\\
		&\quad   + \mu \v  \big(\frac{v^{\v}}{T^{\v}}\big)_y \Delta_{\v} + \lambda \v  \big(\frac{v^{\v}}{T^{\v}}\big)_y \pa_{yy} \Big\}T(\zeta,p,q),\\
		\mathcal{G}_3:&=(\f{\mu \v}{\rho^{\v}})_x \Delta_{\v}\Phi_{y} - (\f{\mu \v}{\rho^{\v}})_y \Delta_{\v}\Phi_{x}  + \v (mu_s^2)_y \mathbf{q}_{xx} - \v (m u_s^2)_x\mathbf{q}_{xy}, \\
		\mathcal{G}_4:&= \v \pa_x G_{12} + \v \pa_y G_{22} - \v \pa_xg_{1}(p,\phi, T(S,p,q)) -  \pa_yg_{2}(p,\phi, T(S,p,q)), \\
		\mathcal{G}_5:&= \big(\frac{\lambda\v}{T^{\v}} u^{\v}\big)_y \mathbf{S}_{xy} - \big(\frac{\lambda\v}{T^{\v}} u^{\v}\big)_x \mathbf{S}_{yy}, \\
		\mathcal{G}_6:&=- \v (\hat{d}_{11} \zeta_{x} )_x,\\
		\mathcal{G}_7:&=  - \v \mathfrak{g}_{1x} - \v \mathfrak{g}_{2y},\\ 
		\mathcal{G}_8:&= \v \mathcal{N}_{1x} + \v  \mathcal{N}_{2y},
	\end{split}
\end{align}
and
\begin{align}\label{17.17-2}
	G_{12x} + G_{22y}
	&= \frac{\lambda u^{\v}}{T^{\v}}   \Big\{  \big(\f1{2\rho_s}\big)_xp_{yy} - \big(\f1{2\rho_s}\big)_yp_{xy} \Big\} - \big( \frac{\lambda u^{\v}}{T^{\v}} \big)_x \Big\{\big(\f1{2\rho_s}\big)_{yy}p + \big(\f1{\rho_s}\big)_yp_y + \big(\f{T_{sy}}{p_s} q\big)_{yy}\Big\}\nonumber\\
	&\quad + \big(\frac{\lambda u^{\v}}{T^{\v}}\big)_y  \Big\{\big(\f1{2\rho_s}\big)_{xy}p + \big(\f1{2\rho_s}\big)_xp_y + \big(\f1{2\rho_s}\big)_yp_x +  \big(\f{T_{sy}}{p_s} q\big)_{xy}\Big\}.
\end{align}

\section{Existence of Solution for Temperature Equation with Neumann BC}\label{SecB}
In this section, we aim to establish the  well-posedness of temperature equation with Neumann BC at $y=0$, that is,
\begin{align}\label{C.01}
	\begin{cases}
		2\rho_s u_s T_{x} + 2\rho_s v_s T_{y} - \kappa \Delta_{\v}T = f,\quad (x,y)\in (0,L)\times \mathbb{R}_+,\\
		T_{y}|_{y=0}=0,\quad T|_{x=0} =\bar{\chi}(y)h,\quad T_{x}|_{x=L}=-g_x.
	\end{cases}
\end{align}

\begin{proposition}\label{Bprop6.16}
	Let $L>0$ be suitably small  $\rm($independent of $\v>0$$\rm)$. Then there is a unique solution $T$ of BVP \eqref{C.01} with 
	\begin{align}\label{B.69}
		\|T-\bar{\chi}(y) h\|^2_{H^{3}} &\leq C_{\v}\, \Big\{\|f\|^2_{H^1} + \|\nabla(\bar{\chi}h)\|^2_{H^2} + \|\sqrt{\v}g_x\|^2_{H^{3/2}_{x=L}} \Big\},
	\end{align}
	and
	\begin{align}\label{B.70}
		\|u_s\nabla_{\v}^2T_{yy}\|^2
		&\lesssim  \|u_s f_{yy}\|^2 + \|u_s(h_{yyyy}, \sqrt{\v} h_{xyyy}, \v h_{xxyy})\|^2   + \v \|u_s g_{xyy}\|^2_{H^{1/2}_{x=L}}\nonumber\\
		&\,\,\, + \|\sqrt{u_s}\sqrt{\v}g_{xyy}\|^2_{x=L}  + C_{\v}\, \Big\{\|f\|^2_{H^1} + \|\nabla(\bar{\chi}h)\|^2_{H^2} + \|\sqrt{\v}g_x\|^2_{H^{3/2}_{x=L}} \Big\}.
	\end{align}
\end{proposition}

\medskip


\begin{lemma}\label{Blem10.5}
Let $L>0$ be suitably small $\rm($independent of $\v>0$$\rm)$.  There exists a unique strong solution $T\in H^3$  to the  the following BVP
\begin{align}\label{C.51}
	\begin{cases}
		\dis	2\rho_s u_s T_{x} + 2\rho_s v_s T_{y} - \kappa \Delta_{\v}T = f,\quad (x,y)\in (0,L)\times \mathbb{R}_+,\\
		\dis T_{y}|_{y=0}=0,	\quad  T|_{x=0} =0, \quad T_{x}|_{x=L}=-g_x,
	\end{cases}
\end{align}	
satisfying
	\begin{align}\label{B.67}
		\|T\|^2_{H^{3}} &\leq C_{\v}\, \Big\{\|f\|^2_{H^1} + \|\sqrt{\v}g_x\|^2_{H^{3/2}_{x=L}} \Big\},
	\end{align}
	and
	\begin{align}\label{B.68}
		\|u_s(T_{yyyy}, \sqrt{\v} T_{xyyy}, \v T_{xxyy})\|^2
		&\lesssim  \|u_s f_{yy}\|^2 + \v \|u_s g_{xyy}\|^2_{H^{1/2}_{x=L}} + \|\sqrt{u_s}\sqrt{\v}g_{xyy}\|^2_{x=L} \nonumber\\
		&\quad + C_{\v}\, \Big\{\|f\|^2_{H^1} + \|\sqrt{\v}g_x\|^2_{H^{3/2}_{x=L}} \Big\}.
	\end{align}
\end{lemma}

\noindent{\bf Proof.} Noting $T_{y}|_{y=0}$, we consider the even extension
\begin{align*}
	\tilde{T}(x,y)=
	\begin{cases}
		T(x,y),\\
		T(x,-y).
	\end{cases}
\end{align*}
Similarly, we denote $\tilde{\rho}_s, \tilde{u}_s, \tilde{f}$ as the even extension of $\rho_s, u_s, f$ in $y$, and $\tilde{v}_s$ as the odd extension of $v_s$ in $y$. 

We point out that even though $f\in H^{k}((0,L)\times\mathbb{R}_+)$ and $g_x \in H^{k+\f12}(\mathbb{R}_+)$, in general,  $\tilde{f}\notin H^{k}((0,L)\times\mathbb{R})$ and $\tilde{g}_x \notin H^{k+\f12}(\mathbb{R})$. Let $\tilde{f}^{\e}, \tilde{g}^{\e}_x$ a smooth sequence, then we consider the following BVP
\begin{align}\label{C.52}
	\begin{cases}
		\dis 2\tilde{\rho}_s \tilde{u}_s \tilde{T}_{x} + 2\tilde{\rho}_s \tilde{v}_s \tilde{T}_{y} - \kappa \Delta_{\v}\tilde{T} = \tilde{f}^{\e},\quad (x,y)\in (0,L)\times \mathbb{R},\\
		\dis \tilde{T}|_{x=0}=0, \quad  \tilde{T}_{x} \big|_{x=L}=-\tilde{g}^{\e}_x.
	\end{cases}
\end{align}
Since the proof is very long, we divide it into several steps.

\smallskip

\noindent{\it Step 1. Auxiliary problem.} Denote 
\begin{align}\label{B.8-0}
\dis(\rho_+, u_+, v_+):=\lim_{y\to\infty} (\rho_s, u_s, v_s)\quad \mbox{and}\quad 
\tilde{\chi}_{\v}(y)=
\begin{cases}
	\bar{\chi}(\v^2 y), \qquad\,\,\, \mbox{for}\,\, y\geq 0,\\
	-\bar{\chi}(-\v^2 y), \quad \mbox{for}\,\, y<0.
\end{cases}
\end{align}
  We first consider the following auxiliary problem
\begin{align}\label{C.53}
	\begin{cases}
		\dis \mathcal{L}\tilde{T}:=2 \rho_{+}u_+ \tilde{T}_{x} + 2 \rho_{+}v_+ \tilde{\chi}_{\v}\,  \tilde{T}_{y} - \kappa \Delta_{\v}\tilde{T} = \tilde{f}^{\e},\quad (x,y)\in (0,L)\times \mathbb{R},\\
		\dis \tilde{T}|_{x=0}=0, \quad  \tilde{T}_{x} \big|_{x=L}=- \tilde{g}_x^{\e}.
	\end{cases}
\end{align}
 
For later use, we define the function space 
\begin{align}
	\mathcal{W}_0^1:=\big\{\varphi \in H^1((0,L)\times\mathbb{R}_+)\,|\, \varphi(0,y)=0 \big\}.
\end{align}
For the weak formulation of \eqref{C.53}, we multiply \eqref{C.53} by the test function $\varphi\in \mathcal{W}^1_0$ to obtain
\begin{align}\label{C.54}
	\mathbf{B}[\tilde{T},\varphi]
	&= \iint \tilde{f} \varphi dydx - \kappa \v \int \tilde{g}^{\e}_{x} \varphi dy \Big|_{x=L},
\end{align}
where 
\begin{align*}
	\mathbf{B}[\tilde{T},\varphi]:&=2\rho_+ u_+\iint \varphi \,  \tilde{T}_{x} dydx + 2\rho_+ v_+\iint \tilde{\chi}_{\v}\, \varphi \,  \tilde{T}_{y} dydx  + \kappa \iint \nabla_{\v}\tilde{T}\cdot \nabla_{\v}\varphi dydx .
\end{align*}
It is clear to see that 
\begin{align}\label{C.55}
	\mathbf{B}[\tilde{T},\tilde{T}]=\kappa \|\nabla_{\v}\tilde{T}\|^2 + \rho_+u_+ \|\tilde{T}\|^2_{x=L} +  O(1) \v^2 \rho_+v_+ \|\tilde{T}\|^2 \geq c_{\v} \|\tilde{T}\|^2_{H^1},
\end{align}
and
\begin{align}\label{C.56}
	\iint \tilde{f}^{\e} \tilde{T} dydx - \kappa \v \int \tilde{g}^{\e}_{x} \tilde{T} dy \Big|_{x=L} \lesssim \sqrt{L} \big(\|\tilde{f}^{\e}\| + \v \|\tilde{g}^{\e}_x\|_{x=L}\big) \|\tilde{T}_{x}\|.
\end{align}

Noting \eqref{C.54}-\eqref{C.56}, using Lax-Milgram theory, one gets an unique solution $\tilde{T}\in \mathcal{W}_0^1$ such that \eqref{C.53} holds. That means we  establish a weak solution $\tilde{T}$ for \eqref{C.53} with 
\begin{align}
	\|\tilde{T}\|_{H^{1}((0,L)\times\mathbb{R})} \leq C_{\v}\, \big\{\|\tilde{f}^{\e}\|_{L^2((0,L)\times\mathbb{R})} + \|\tilde{g}^{\e}_x\|_{x=L} \big\}.
\end{align} 
By the regularity theory of elliptic equation, one has that
\begin{align}
	\|\tilde{T}\|_{H^{k+2}((0,L)\times\mathbb{R})} \lesssim C_{\v} \big\{\|\tilde{f}^{\e}\|_{H^k((0,L)\times\mathbb{R})} + \|\tilde{g}^{\e}_x\|_{H^{k+\f12}(\mathbb{R})}\big\}.
\end{align}

\smallskip

\noindent{\it Step 2. Fredholm arguments.} 
We define 
\begin{align}\label{C.18}
	\tilde{\rho}^{\tau}_s=(\tilde{\rho}_{s}*J_{\tau}), \quad \tilde{u}^{\tau}_s=(\tilde{u}_{s}*J_{\tau}) \quad \mbox{and}\quad \tilde{v}^{\tau}_s=(\tilde{v}_{s}*J_{\tau}),
\end{align}
where the smooth is only for $y\in\mathbb{R}$.  
We consider an approximate problem \eqref{C.52}, i.e., 
\begin{align}\label{C.59}
	\begin{cases}
		\dis 2 \rho_+ u_+ \tilde{T}_{x} + 2 \rho_+ v_+ \tilde{\chi}_{\v}\,\tilde{T}_{y} - \kappa \Delta_{\v}\tilde{T} + 2[\tilde{\rho}^{\tau}_s\tilde{u}^{\tau}_s-\rho_+u_+] \tilde{T}_{x} + 2[\tilde{\rho}^{\tau}_s\tilde{v}^{\tau}_s-\rho_+v_+ \tilde{\chi}_{\v}] \tilde{T}_{y} = \tilde{f}^{\e},\\
		\dis \tilde{T}|_{x=0}=0, \quad  \tilde{T}_{x} \big|_{x=L}=-\tilde{g}^{\e}_x.
	\end{cases}
\end{align}
Denote $\bar{T} = \mathcal{L}\tilde{T} \Longleftrightarrow \tilde{T}=\mathcal{L}^{-1}\bar{T}$, then we can represent \eqref{C.59} as
\begin{align}
	\bar{T} + 2[\tilde{\rho}^{\tau}_s\tilde{u}^{\tau}_s-\rho_+u_+]\big(\mathcal{L}^{-1}\bar{T}\big)_x +  2[\tilde{\rho}^{\tau}_s\tilde{v}^{\tau}_s-\rho_+v_+\tilde{\chi}_{\v}]\big(\mathcal{L}^{-1}\bar{T}\big)_y = \tilde{f}^{\e}.
\end{align}
We know that
\begin{align*}
	2[\tilde{\rho}^{\tau}_s\tilde{u}^{\tau}_s-\rho_+u_+]\big(\mathcal{L}^{-1}\bar{T}\big)_x +  2[\tilde{\rho}^{\tau}_s\tilde{v}^{\tau}_s-\rho_+v_+\tilde{\chi}_{\v}]\big(\mathcal{L}^{-1}\bar{T}\big)_y:\,\, \bar{T}\in H^1 \longmapsto H^2.
\end{align*}
Due to the decay property $|\tilde{\rho}^{\tau}_s\tilde{u}^{\tau}_s-\rho_+u_+| + |\tilde{\rho}^{\tau}_s\tilde{v}^{\tau}_s-\rho_+v_+\tilde{\chi}_{\v}|\lesssim \la \sqrt{\v}y\ra^{-\f12 \mathfrak{l}_0}$, we see that 
$$2[\tilde{\rho}^{\tau}_s\tilde{u}^{\tau}_s-\rho_+u_+]\big(\mathcal{L}^{-1}\bar{T}\big)_x +  2[\tilde{\rho}^{\tau}_s\tilde{v}^{\tau}_s-\rho_+v_+\tilde{\chi}_{\v}]\big(\mathcal{L}^{-1}\bar{T}\big)_y,$$
is a {\it compact} operator from $H^1\longmapsto H^1$. Therefore, by Fredholm alternative, to establish the solvability of \eqref{C.59}, we need only to prove the uniqueness of the homogeneous problem
\begin{align}\label{C.62}
	\begin{cases}
		\dis 2\tilde{\rho}^{\tau}_s \tilde{u}^{\tau}_s \tilde{T}_{x} +  2\tilde{\rho}^{\tau}_s \tilde{v}^{\tau}_s \tilde{T}_{y} - \kappa \Delta_{\v}\tilde{T} = 0,\quad (x,y)\in (0,L)\times \mathbb{R},\\
		\dis \tilde{T}|_{x=0}=0, \quad  \tilde{T}_{x} \big|_{x=L}=0.
	\end{cases}
\end{align}

Multiplying \eqref{C.62} by $\tilde{T}_{x}$, one gets that 
\begin{align}\label{C.63}
	\|\sqrt{\tilde{\rho}_s\tilde{u}_s}\tilde{T}_{x}\|^2 + \f12 \kappa \v \|\tilde{T}_{x}\|^2_{x=0} +  \f12 \kappa  \|\tilde{T}_{y}\|^2_{x=L} \lesssim \|\sqrt{\tilde{u}_s}\tilde{T}_{y}\|^2.
\end{align}
On the other hand, multiplying \eqref{C.62} by $\tilde{T}$, one has
\begin{align}\label{C.64}
	\|\sqrt{\tilde{\rho}_s\tilde{u}_s} \tilde{T}\|^2_{x=L} + \kappa \|\nabla_{\v}\tilde{T}\|^2 \lesssim \|\sqrt{\tilde{u}_s}\tilde{T}\|^2\lesssim L^2 \|\sqrt{\tilde{u}_s} \tilde{T}_{x}\|^2.
\end{align}
Noting the smallness of $L>0$ (independent of $\v$ and $\tau$), combining \eqref{C.63}-\eqref{C.64}, one obtain $\|\nabla_{\v}\tilde{T}\|=0$ which shows the uniqueness of \eqref{C.62}.

Thus we  establish the existence and uniqueness of strong solution  $\tilde{T} \in H^2((0,T)\times \mathbb{R})$ of the following problem
\begin{align}\label{B.56}
	\begin{cases}
		\dis 2\tilde{\rho}^{\tau}_s \tilde{u}^{\tau}_s \tilde{T}_{x} +  2\tilde{\rho}^{\tau}_s \tilde{v}^{\tau}_s \tilde{T}_{y} - \kappa \Delta_{\v}\tilde{T} = \tilde{f}^{\e},\quad (x,y)\in (0,L)\times \mathbb{R},\\
		\dis \tilde{T}|_{x=0}=0, \quad  \tilde{T}_{x} \big|_{x=L}=- \tilde{g}^{\e}_x.
	\end{cases}
\end{align}
We point out that $\|\tilde{T}\|_{H^1}$ may depend on $\tau>0$. For the higher regularities, since $\|\pa_{yy}\tilde{u}_s^{\tau}\|$ depends on $\tau^{-1}$, then by the regularity theory of elliptic equation, we only have 
\begin{align}\label{B.57}
	\|\tilde{T}\|_{H^{k+2}}\lesssim C\Big\{\tau^{-1}, \|\tilde{f}^{\e}\|_{H^k} , \|\tilde{g}^{\e}_x\|_{H^{k+\f12}_{x=L}}, \|\tilde{T}\|_{H^1}\Big\}.
\end{align}

\smallskip

\noindent{\it Step 3. Some uniform-in-$\tau$ estimates.} 
To solve \eqref{C.52}, we need to take the limit $\tau\to0$ in \eqref{B.56}. Since the estimate in \eqref{B.57} depend on $\tau$,   we need to derive some uniform-in-$\tau$ estimates.

Multiplying \eqref{B.56} by $\tilde{T}$, one has 
\begin{align*}
	\|\sqrt{\tilde{\rho}^{\tau}_s \tilde{u}^{\tau}_s} \tilde{T}\|^2_{x=L} + \kappa \|\nabla_{\v}\tilde{T}\|^2 
	&=\big\la (\tilde{\rho}^{\tau}_s \tilde{u}^{\tau}_s)_x + (\tilde{\rho}^{\tau}_s \tilde{v}^{\tau}_s)_y,\, \tilde{T}^2 \big\ra + \kappa \big\la \v \tilde{T}_x, \, \tilde{T} \big\ra_{x=L} + \big\la \tilde{f}^{\e},\, \tilde{T} \big\ra\nonumber\\
	&\lesssim \|\sqrt{\tilde{u}^{\tau}_s} \tilde{T}\|^2 + \v \|\tilde{g}^{\e}_x\|_{x=L} \|\tilde{T}\|_{x=L} + L\|\tilde{f}^{\e}\|\cdot \|\tilde{T}_x\|,
\end{align*}
which yields immediately that 
\begin{align}\label{B.58}
	\|\sqrt{\tilde{\rho}^{\tau}_s \tilde{u}^{\tau}_s} \tilde{T}\|^2_{x=L} + \kappa \|\nabla_{\v}\tilde{T}\|^2  
	&\lesssim L\|\sqrt{\tilde{u}^{\tau}_s} \tilde{T}_x\|^2 + \v \|\tilde{g}^{\e}_x\|^2_{x=L} + \f{L^2}{\v} \|\tilde{f}^{\e}\|^2.
\end{align}
Multiplying \eqref{C.52} by $\tilde{T}_x$, one has 
\begin{align}\label{B.59}
	&\|\sqrt{\tilde{\rho}^{\tau}_s \tilde{u}^{\tau}_s} \tilde{T}_x\|^2 + \f12 \kappa \|\sqrt{\v}\tilde{T}_x\|^2_{x=0} + \f12 \kappa \|\tilde{T}_y\|^2_{x=L} \nonumber\\
	&\leq \f12 \kappa \|\sqrt{\v}\tilde{g}_x^{\e}\|^2_{x=L} + C\|\tilde{T}_y\|^2 +  C\|\sqrt{\v}\tilde{T}_x\|^2 +  C\f{1}{\v} \|\tilde{f}^{\e}\|^2.
\end{align}
Combining \eqref{B.58}-\eqref{B.59} and noting  $L\ll1$ $\rm($independent of $\v>0$$\rm)$, one obtains
\begin{align}\label{B.60}
	\|\sqrt{\tilde{u}^{\tau}_s} \tilde{T}\|^2_{x=L} + \|\nabla_{\v}\tilde{T}\|^2   + \|\sqrt{\tilde{u}^{\tau}_s} \tilde{T}_x\|^2 +  \|\sqrt{\v}\tilde{T}_x\|^2_{x=0} + \|\tilde{T}_y\|^2_{x=L} 
	\lesssim \|\sqrt{\v}\tilde{g}^{\e}_x\|^2_{x=L} + \f{1}{\v} \|\tilde{f}^{\e}\|^2.
\end{align}
The key point is that the constant in \eqref{B.60} is independent of $\tau>0$. Then, with  \eqref{B.60}, applying the regularity theory of elliptic problem, one has
\begin{align}\label{B.61}
	\|\tilde{T}\|^2_{H^{2}} \leq C_{\v} \, \Big\{\|\tilde{f}^{\e}\|^2 + \|\tilde{g}^{\e}_x\|^2_{H^{1/2}_{x=L}}\Big\}.
\end{align}

\smallskip

It is direct to know that    $\|(\pa_y \tilde{\rho}_{s}^{\tau}, \pa_y \tilde{u}_{s}^{\tau}, \pa_y \tilde{v}_{s}^{\tau})\|_{L^\infty}$ are independent of $\tau>0$. Then, using the regularity theory of elliptic problem and \eqref{B.61},  we can obtain 
\begin{align}\label{B.62}
	\|\tilde{T}\|^2_{H^{3}} &\leq C_{\v}\,\Big\{ \|\big(\tilde{f}^{\e}, \tilde{\rho}^{\tau}_s \tilde{u}^{\tau}_s \tilde{T}_{x}, \, \tilde{\rho}^{\tau}_s \tilde{v}^{\tau}_s \tilde{T}_{y}\big)\|^2_{H^1} + \|\tilde{g}^{\e}_x\|^2_{H^{3/2}_{x=L}}\Big\}  \leq C_{\v}\, \Big\{\|\tilde{f}^{\e}\|^2_{H^1} + \|\sqrt{\v}\tilde{g}^{\e}_x\|^2_{H^{3/2}_{x=L}} \Big\}.
\end{align}

\

For the fourth order derivatives, to obtain uniform-in-$\tau$ estimates, we have to consider a space weight $\tilde{u}_s$. Applying $\pa_y$ to \eqref{B.56}, one has
\begin{align*} 
	\dis \kappa \Delta_{\v}\tilde{T}_{yy} = 2\big(\tilde{\rho}^{\tau}_s \tilde{u}^{\tau}_s \tilde{T}_{x}\big)_{yy} + 2\big(\tilde{\rho}^{\tau}_s \tilde{v}^{\tau}_s \tilde{T}_{y}\big)_{yy} - \tilde{f}^{\e}_{yy},
\end{align*}
which yields immediately that 
\begin{align}\label{B.63}
	&\kappa^2 \|\tilde{u}_s(\tilde{T}_{yyyy}, \v \tilde{T}_{xxyy})\|^2 + 2\kappa^2 \big\la \v \tilde{u}_s^2 \tilde{T}_{xxyy}, \tilde{T}_{yyyy} \big\ra  \nonumber\\
	& \lesssim \|\tilde{u}_s \big(\tilde{\rho}^{\tau}_s \tilde{u}^{\tau}_s \tilde{T}_{x}\big)_{yy}\|^2 + \|\tilde{u}_s \big(\tilde{\rho}^{\tau}_s \tilde{v}^{\tau}_s \tilde{T}_{y}\big)_{yy}\|^2 + \|\tilde{u}_s\tilde{f}^{\e}_{yy}\|^2.
\end{align}
 A direct calculation shows that 
\begin{align}\label{B.64}
	\big\la \v \tilde{u}_s^2 \tilde{T}_{xxyy}, \tilde{T}_{yyyy} \big\ra
	& = - \big\la \v \tilde{u}_s^2 \tilde{T}_{xyy}, \tilde{T}_{xyyyy} \big\ra - \big\la 2\v \tilde{u}_s\tilde{u}_{sx} \tilde{T}_{xyy}, \tilde{T}_{yyyy} \big\ra + \big\la \v \tilde{u}_s^2 \tilde{T}_{xyy}, \tilde{T}_{yyyy} \big\ra_{x=L}\nonumber\\
	&=\|\tilde{u}_s \sqrt{\v} \tilde{T}_{xyyy}\|^2 + \big\la 2\v \tilde{u}_s \tilde{u}_{sy} \tilde{T}_{xyy}, \tilde{T}_{xyyy} \big\ra - \big\la 2\v \tilde{u}_s\tilde{u}_{sx} \tilde{T}_{xyy}, \tilde{T}_{yyyy} \big\ra\nonumber\\
	&\quad  + \big\la \v \tilde{u}_s^2 \tilde{T}_{xyy}, \tilde{T}_{yyyy} \big\ra_{x=L}\nonumber\\
	&\geq \f45 \|\tilde{u}_s \sqrt{\v} \tilde{T}_{xyyy}\|^2 - \xi \|\tilde{u}_s\tilde{T}_{yyyy}\|^2 - C_{\xi} \|\tilde{T}\|^2_{H^3} - \big\la \v \tilde{u}_s \tilde{g}^{\e}_{xyy}, \tilde{u}_{sy}T_{yyy} \big\ra_{x=L}\nonumber\\
	&\quad  + \big\la \v \tilde{u}_s \tilde{g}^{\e}_{xyy}, \big(\tilde{u}_s \tilde{T}_{yyy}\big)_y \big\ra_{x=L} \nonumber\\
	&\geq \f45 \|\tilde{u}_s \sqrt{\v} \tilde{T}_{xyyy}\|^2 - C \xi \|\tilde{u}_s\tilde{T}_{yyyy}\|^2 - Ca\v \|\sqrt{\tilde{u}_s}\tilde{g}^{\e}_{xyy}\|_{x=L}\cdot \|\sqrt{\tilde{u}_s}\tilde{T}_{yyy}\|_{x=L}\nonumber\\
	&\quad - C \v \|\tilde{u}_s\tilde{g}^{\e}_{xyy}\|_{H^{1/2}_{x=L}}\cdot \|\tilde{u}_s\tilde{T}_{yyy}\|_{H^{1/2}_{x=L}} - C_{\xi} \|\tilde{T}\|^2_{H^3}\nonumber\\
	&\geq \f34 \|\tilde{u}_s \sqrt{\v} \tilde{T}_{xyyy}\|^2 - C \xi \|\tilde{u}_s\tilde{T}_{yyyy}\|^2  - C_{\xi} \|\tilde{T}\|^2_{H^3} - C \v \|\tilde{u}_s\tilde{g}^{\e}_{xyy}\|^2_{H^{1/2}_{x=L}}\nonumber\\
	&\quad  - C\|\sqrt{\tilde{u}_s}\sqrt{\v}\tilde{g}^{\e}_{xyy}\|^2_{x=L}.
\end{align}

For the first and second terms on RHS of \eqref{B.63}, we have to be careful since $\pa_y \tilde{\rho}_s$ and $\pa_y \tilde{u}_s$ are discontinuous at $y=0$, then  $\|\pa_{yy}\tilde{\rho}_s^{\tau}\|$ and $\|\pa_{yy}\tilde{u}_s^{\tau}\|$ depend on $\tau^{-1}$. In fact, we note  that 
\begin{align} \label{B.36}
	|\tilde{\rho}^{\tau}_{syy}(x,y)|
	&=  \mathbf{1}_{\{|y|\leq 3\tau\}}\f{1}{\tau}\int \pa_y(\tilde{\rho}_{s1}^{\tau})(x,z)\,  J_{\tau}(y-z) dz + \mathbf{1}_{\{|y|\geq 3\tau\}} \int  (\pa_{yy}\tilde{\rho}_{s})(x,z)\,  J_{\tau}(y-z) dz\nonumber\\
	&\lesssim  \mathbf{1}_{\{|y|<3\tau\}}\f{1}{\tau} \|\rho_{sy}\|_{L^\infty} +  \mathbf{1}_{\{|y|\geq 3\tau\}} \|\pa_{yy}\rho_{s}\|_{L^\infty},
\end{align}
and similarly
\begin{align}\label{B.37}
	|\tilde{u}^{\tau}_{syy}(x,y)| \lesssim \mathbf{1}_{\{|y|\geq 3\tau\}} \|\pa_{yy} u_{s}\|_{L^\infty} + \mathbf{1}_{\{|y|<3\tau\}}\f{1}{\tau} \|u_{sy}\|_{L^\infty}.
\end{align}
Hence we have from \eqref{B.36}-\eqref{B.37} that 
\begin{align}\label{B.38}
	|\tilde{u}_s\big(\tilde{\rho}^{\tau}_s \tilde{u}^{\tau}_s\big)_{yy}|
	&\lesssim  |\tilde{u}_s\tilde{\rho}^{\tau}_{syy} \tilde{u}^{\tau}_s| + 2|\tilde{u}_s\tilde{\rho}^{\tau}_{sy} \tilde{u}^{\tau}_{sy}| + |\tilde{u}_s\tilde{\rho}^{\tau}_s \tilde{u}^{\tau}_{syy}|
	\lesssim 1.
\end{align}
Using \eqref{B.36}-\eqref{B.38}, one has
\begin{align}\label{B.65}
	\|\tilde{u}_s \big(\tilde{\rho}^{\tau}_s \tilde{u}^{\tau}_s \tilde{T}_{x}\big)_{yy}\|^2 + \|\tilde{u}_s \big(\tilde{\rho}^{\tau}_s \tilde{v}^{\tau}_s \tilde{T}_{y}\big)_{yy}\|^2
	&\lesssim \|\tilde{T}\|^2_{H^3}.
\end{align}
Substituting \eqref{B.64}-\eqref{B.65} into \eqref{B.63}, one has that 
\begin{align}\label{B.66}
	\|\tilde{u}_s(\tilde{T}_{yyyy}, \sqrt{\v} \tilde{T}_{xyyy}, \v \tilde{T}_{xxyy})\|^2
	& \lesssim  \|\tilde{u}_s\tilde{f}^{\e}_{yy}\|^2 + \v \|\tilde{u}_s\tilde{g}^{\e}_{xyy}\|^2_{H^{1/2}_{x=L}} + \|\sqrt{\tilde{u}_s}\sqrt{\v}\tilde{g}^{\e}_{xyy}\|^2_{x=L} \nonumber\\
	&\quad + C_{\v}\, \Big\{\|\tilde{f}^{\e}\|^2_{H^1} + \|\sqrt{\v}\tilde{g}^{\e}_x\|^2_{H^{3/2}_{x=L}} \Big\}.
\end{align}

\smallskip

\noindent{\it Step 4. Taking limit $\tau, \e\to 0$.} With the help of above estimates, we can take $\tau\to0$ and then $\e\to0$ in \eqref{B.56}, one obtains the existence of unique solution $\tilde{T}$ of \eqref{C.52}, hence the solution $T$ of \eqref{C.51}. The estimates \eqref{B.67}-\eqref{B.68} follows directly from \eqref{B.62} and \eqref{B.66}. Therefore the proof of Lemma \ref{Blem10.5} is completed. $\hfill\Box$

\

\noindent{\bf Proof of Proposition \ref{Bprop6.16}:}  To solve the original BVP \eqref{C.01}, we denote 
\begin{align*}
	\mathcal{T}:=T-\bar{\chi}(y) h,
\end{align*}
then  \eqref{C.01} is rewritten as
\begin{align}\label{C.01-0}
	\begin{cases}
		2\rho_s u_s \mathcal{T}_{x} + 2\rho_s v_s \mathcal{T}_{y} - \kappa \Delta_{\v}\mathcal{T} = F,\\
		T_{y}|_{y=0}=0,\quad \mathcal{T}|_{x=0} =0\quad \mbox{and}\quad  \mathcal{T}_{x}|_{x=L}=-g_x-\bar{\chi}(y)h_x
	\end{cases}
\end{align}
where 
\begin{align*}
	F:=f- \rho_s u_s \big(\bar{\chi}(y)h\big)_x - 2\rho_s v_s \big(\bar{\chi}(y) h\big)_y 
	+ \Delta_{\v} \big(\bar{\chi}(y) h\big).
\end{align*}
Applying Lemma \ref{Blem10.5} to \eqref{C.01-0}, we conclude  Proposition \ref{Bprop6.16}. $\hfill\Box$

%

\medskip
\section{Existence of Solution for Pseudo Entropy  with Dirichlet BC}\label{secC}
In this section, we aim to establish the existence of solution for the pseudo entropy equation with Dirichlet boundary condition at $y=0$. 
We consider 
\begin{align}\label{D4.16}
	\begin{cases}
		2 \rho_s u_s S_x + 2\rho_s v_s S_y - \kappa \Delta_{\v}S = f,\quad (x,y) \in (0,L)\times \mathbb{R}_+,\\
		S|_{y=0}=- h(x,0) \quad \&\quad 
		S|_{x=0}=0,\\
		\big\{2\rho_s u_s S_x- \kappa S_{yy}\big\}\big|_{x=L} =  g(L,y).
	\end{cases}
\end{align}
with $h(0,0)=0$.
\begin{proposition}\label{propC.1}
	Let $L>0$ be suitably small $\rm($independent of $\v>0$$\rm)$. There is a unique solution $S$ of BVP \eqref{D4.16} with 
	\begin{align}
		S\in H^2 ([0,L]\times \R_+) \quad \& \quad u_s \nabla^3 S \in L^2 ([0,L]\times \R_+).
	\end{align}
\end{proposition}

\smallskip

We first reformulate \eqref{D4.16}. We define 
\begin{align}\label{D4.3-0}
	\bar{S}(x,y)=S(x,y)+h(x,y) \chi(y).
\end{align}
then it holds that 
\begin{align}\label{D4.18}
	\begin{cases}
		2 \rho_s u_s \bar{S}_x + 2\rho_s v_s \bar{S}_y - \kappa \Delta_{\v}\bar{S} = \bar{f},\quad (x,y) \in (0,L)\times \mathbb{R}_+,\\
		\dis \bar{S}|_{y=0}=0 \quad \&\quad 
		\bar{S}|_{x=0}=h(0,y) \chi(y),\\
		\big[2\rho_s u_s \bar{S}_x-\kappa \bar{S}_{yy}\big]\big|_{x=L}
		=\bar{g}\big|_{x=L}.
	\end{cases}
\end{align}
where 
\begin{align*}
	\begin{split}
		\bar{f}:&= f + 2 \rho_s u_s  (h\chi)_x + 2\rho_s v_s (h \chi)_y - \kappa \Delta_{\v}  (h\chi ),\\
		\bar{g}:&=g + 2\rho_s u_s (h\chi)_x  - \kappa (h\chi)_{yy}.
	\end{split}
\end{align*}

Furthermore we consider 
\begin{align}
	\tilde{S}:=\bar{S} - h(0,y) \chi(y) \chi(\f{2x}{L}), 
\end{align}
then we can obtain
\begin{align}\label{D4.21}
	\begin{cases}
		2 \rho_s u_s \tilde{S}_x + 2\rho_s v_s \tilde{S}_y - \kappa \Delta_{\v}\tilde{S} = \tilde{f},\quad (x,y) \in (0,L)\times \mathbb{R}_+,\\
		\dis \tilde{S}|_{y=0}=0 \quad \&\quad 
		\tilde{S}|_{x=0}=0,\\
		\big[2\rho_s u_s\tilde{S}_x-\kappa \tilde{S}_{yy}\big]\big|_{x=L}
		=\bar{g}\big|_{x=L}.
	\end{cases}
\end{align}
where 
\begin{align*}
		\tilde{f}&:=\bar{f} -2 \rho_s u_s  \Big(h(0,y)\chi(y) \chi(\f{2x}{L}) \Big)_x - 2\rho_s v_s \Big(h(0,y)\chi(y) \chi(\f{2x}{L})\Big)_y   + \kappa \Delta_{\v}\Big(h(0,y) \chi(y) \chi(\f{2x}{L})\Big).
\end{align*}
Then we need only to establish the existence for \eqref{D4.21}.

\smallskip

 Noting $\tilde{S}|_{y=0}=0$, we shall use the odd extension argument.  For this, we first consider the following model problem:
\begin{align}\label{D4.17-0}
	\begin{cases}
		2 \rho_+ u_+ \tilde{\Theta}_x + 2\rho_+ v_+ \tilde{\chi}_{\v} \tilde{\Theta}_y - \kappa \Delta_{\v}\tilde{\Theta}=\tilde{\mathscr{F}},\quad (x,y) \in (0,L)\times \mathbb{R},\\
		\dis 
		\tilde{\Theta}|_{x=0}=0,\quad \big[2\rho_+ u_+ \tilde{\Theta}_x - \kappa \tilde{\Theta}_{yy}\big]\big|_{x=L}
		=\tilde{\mathscr{G}} \big|_{x=L}   \quad \&\quad \lim_{y\to \infty}\tilde{\Theta}=0.
	\end{cases}
\end{align}
where $\tilde{\chi}_{\v}$ is the odd function  in \eqref{B.8-0}.
\begin{lemma}\label{lemD.2}
Assume $\tilde{\mathscr{F}}\in L^2([0,L]\times \mathbb{R})$ and $\tilde{\mathscr{G}}\in L^2(\R)$, then there exists a  unique strong solution $\tilde{\Theta}\in H^2$  to  \eqref{D4.17-0} satisfying
\begin{align}\label{D4.8}
\|\tilde{\Theta}\|_{H^2} +  \|\tilde{\Theta}_{yy}\|_{x=L}
&\lesssim  C_{\v} \big\{\|\tilde{\mathscr{F}}\| + \|\tilde{\mathscr{G}}\|_{x=L}\big\}.
\end{align}
Moreover, we assume $\tilde{\mathscr{F}}\in H^1([0,L]\times \mathbb{R})$ and $\tilde{\mathscr{G}}\in  H^{1/2}(\R)$, then it holds that 
\begin{align}\label{D4.9}
\|\tilde{\Theta}\|_{H^3}  \lesssim C_{\v}\big\{ \|\tilde{\mathscr{F}}\|_{H^1}  + \|\tilde{\mathscr{G}}\|_{H^{1/2}}\big\}.
\end{align}
\end{lemma}

\noindent{\bf Proof.}  We divide the following proof into several steps since it is long.

\smallskip

\noindent{\it Step 1. Existence of weak solution to \eqref{D4.17-0}.} We define the function space
\begin{align}
	\mathfrak{H}^1:=\big\{ w \in H^1([0,L]\times\R)\,\big|\,\,  w|_{x=0}=0,\,\,  w_y|_{x=L}\in L^2(\R)\big\}.
\end{align}

Multiplying \eqref{D4.17-0} by a smooth test function $\varphi$ with $\varphi(0,y)=0$, one gets
\begin{align}\label{D4.26}
\mathbb{B}[\tilde{\Theta}, \varphi]	&:= \kappa \iint \big[\v \tilde{\Theta}_x \varphi_x + \tilde{\Theta}_y \varphi_y\big] dydx  -2\rho_+ u_+  \iint \tilde{\Theta} \varphi_x dydx - 2\rho_+ v_+ \iint\big(\tilde{\chi}_{\v}\varphi)_y  \tilde{\Theta} dydx  \nonumber\\
	&\quad  + 2\rho_+ u_+  \int  \tilde{\Theta} \varphi dy \Big|_{x=L}  + \f{\kappa^2 \v}{2\rho_+u_+}  \int \varphi_y  \tilde{\Theta}_{y}   dy \Big|_{x=L} \nonumber\\
	&=\iint \tilde{\mathscr{F}}\varphi dydx + \f{\kappa \v}{2\rho_+u_+}  \int_0^\infty \varphi \, \tilde{\mathscr{G}}  dy \Big|_{x=L}
\end{align}
where we have used the boundary condition in $\eqref{D4.17-0}$ to derive
\begin{align}\label{D4.27}
	- \kappa \v \int_0^\infty \tilde{\Theta}_x \varphi dy \Big|_{x=L}
	&= -\f{\kappa \v}{2\rho_+u_+}  \int_0^\infty \varphi \big[\kappa \tilde{\Theta}_{yy} + \tilde{\mathscr{G}}\big] dy \Big|_{x=L} \nonumber\\
	&=\f{\kappa^2 \v}{2\rho_+u_+}  \int_0^\infty \varphi_y  \tilde{\Theta}_{y}   dy \Big|_{x=L}  -\f{\kappa \v}{2\rho_+u_+}  \int_0^\infty \varphi  \tilde{\mathscr{G}}  dy \Big|_{x=L}.
\end{align}

 Noting 
\begin{align*}
	\Big|\iint \tilde{\chi}_{\v} \, \pa_y(|\tilde{\Theta}|^2)  dydx\Big|&=\Big|\iint \pa_y \tilde{\chi}_{\v} |\tilde{\Theta}|^2 dydx\Big| \lesssim L \v^2 \|\tilde{\Theta}_x\|^2,
\end{align*}
it is direct to get
\begin{align}\label{D4.25}
\mathbb{B}[\tilde{\Theta}, \tilde{\Theta}] \geq \rho_+u_+ \|\tilde{\Theta}\|^2_{x=L}  + \f12 \kappa \|(\sqrt{\v}\tilde{\Theta}_x, \tilde{\Theta}_y)\|^2 + \f{\kappa^2 \v}{2\rho_+u_+} \|\tilde{\Theta}_y\|^2_{x=L}.
\end{align}

Multiplying \eqref{D4.17-0} by $\tilde{\Theta}$, one obtains 
\begin{align}\label{D4.24}
\iint \tilde{\mathscr{F}}\varphi dydx + \f{\kappa \v}{2\rho_+u_+}  \int_0^\infty \varphi \, \tilde{\mathscr{G}}  dy \Big|_{x=L}	& \lesssim \|\varphi_x\| \big\{L \|\tilde{\mathscr{F}}\|  +  \v \|\tilde{\mathscr{G}}\|_{x=L} \big\}.
\end{align}


Noting \eqref{D4.25}-\eqref{D4.24}, then applying Lax-Milgram theorem to \eqref{D4.26}, there exists a unique weak solution $\tilde{\Theta}\in \mathfrak{H}^1$ of  \eqref{D4.17-0}  satisfying
\begin{align}\label{D4.25-0}
	\|(\tilde{\Theta}, \tilde{\Theta}_y)\|^2_{x=L} + \|(\tilde{\Theta}_x, \tilde{\Theta}_y)\|^2 &\lesssim  C_{\v} \big\{\|\tilde{\mathscr{F}}\|^2 + \|\tilde{\mathscr{G}}\|^2_{x=L}\big\}.
\end{align}

\smallskip

\noindent{\it Step 2. $H^2$-estimate}.  Multiplying \eqref{D4.17-0} by $\tilde{\Theta}_{xx}$, one obtains that 
\begin{align}\label{D4.29}
	&\rho_+u_+ \|\tilde{\Theta}_x\|^2 \big|^{x=L}_{x=0}  - \kappa \|\sqrt{\v}\tilde{\Theta}_{xx}\|^2 - \kappa \big\la \tilde{\Theta}_{xx},\, \tilde{\Theta}_{yy} \big\ra + 2\rho_+ v_+ \big\la \tilde{\Theta}_{xx},\, \tilde{\chi}_{\v} \tilde{\Theta}_y\big\ra =\big\la \tilde{\mathscr{F}}, \,  \tilde{\Theta}_{xx}\big\ra.
\end{align}
A direct calculation shows that 
\begin{align*}
\begin{split}
\big\la \tilde{\Theta}_{xx},\, \tilde{\chi}_{\v} \tilde{\Theta}_y\big\ra&=- \big\la \tilde{\Theta}_{x},\, \tilde{\chi}_{\v} \tilde{\Theta}_{xy}\big\ra +  \big\la \tilde{\Theta}_{x},\, \tilde{\chi}_{\v}\tilde{\Theta}_y\big\ra_{x=L},\\
	-  \big\la \tilde{\Theta}_{xx},\, \tilde{\Theta}_{yy} \big\ra
&=\big\la \tilde{\Theta}_{x},\, \tilde{\Theta}_{xyy} \big\ra -  \big\la \tilde{\Theta}_{x},\, \tilde{\Theta}_{yy} \big\ra_{x=L} =-\|\tilde{\Theta}_{xy}\|^2-  \big\la \tilde{\Theta}_{x},\, \tilde{\Theta}_{yy} \big\ra_{x=L},
\end{split}
\end{align*}
which, together with \eqref{D4.29}, yields that 
\begin{align}\label{D4.32}
	&\kappa \|(\sqrt{\v}\tilde{\Theta}_{xx}, \tilde{\Theta}_{xy})\|^2 + \rho_+u_+ \|\tilde{\Theta}_x\|^2_{x=0} +  \big\la \tilde{\Theta}_{x},\, \kappa \tilde{\Theta}_{yy} - \rho_+u_+\tilde{\Theta}_x\big\ra_{x=L}\nonumber\\
	&=-\big\la \tilde{\mathscr{F}}, \,  \tilde{\Theta}_{xx}\big\ra - 2\rho_+ v_+ \big\la \tilde{\Theta}_{x},\, \tilde{\chi}_{\v} \tilde{\Theta}_{xy}\big\ra +  2\rho_+ v_+ \big\la \tilde{\Theta}_{x},\, \tilde{\chi}_{\v} \tilde{\Theta}_y\big\ra_{x=L}.
\end{align}
Noting $\eqref{D4.17-0}$, we have 
\begin{align}\label{D4.33}
	\big\la \tilde{\Theta}_{x},\, \kappa \tilde{\Theta}_{yy} - \rho_+u_+\tilde{\Theta}_x\big\ra_{x=L}
	&=\rho_+u_+\|\tilde{\Theta}_x\|^2_{x=L} + O(1)\|\tilde{\Theta}_{x}\|_{x=L} \cdot \|\tilde{\mathscr{G}}\|_{x=L}\nonumber\\
	&\geq \f34 \rho_+u_+\|\tilde{\Theta}_x\|^2_{x=L} - C\|\tilde{\mathscr{G}}\|_{x=L}^2.
\end{align}
and 
\begin{align}\label{D4.34}
	&-\big\la \tilde{\mathscr{F}}, \,  \tilde{\Theta}_{xx}\big\ra - 2\rho_+ v_+ \big\la \tilde{\Theta}_{x},\, \tilde{\chi}_{\v} \tilde{\Theta}_{xy}\big\ra +  2\rho_+ v_+ \big\la \tilde{\Theta}_{x},\, \tilde{\chi}_{\v} \tilde{\Theta}_y\big\ra_{x=L} \nonumber\\
	&\leq \f18 \kappa \|(\sqrt{\v}\tilde{\Theta}_{xx}, \tilde{\Theta}_{xy})\|^2 + C_{\v} \|\tilde{\mathscr{F}}\|^2 + C_{\kappa} \|\tilde{\Theta}_x\|^2 + \|\tilde{\Theta}_x\|_{x=L}\cdot \|\tilde{\Theta}_y\|_{x=L}.
\end{align}
Substituting \eqref{D4.33}-\eqref{D4.34} into \eqref{D4.32} and using \eqref{D4.25-0}, one obtains that 
	\begin{align}\label{D4.35}
		\|(\sqrt{\v}\tilde{\Theta}_{xx}, \tilde{\Theta}_{xy})\|^2 +  \|\tilde{\Theta}_x\|^2_{x=0} +   \|\tilde{\Theta}_x\|^2_{x=L}  
		&\lesssim  C_{\v} \big\{\|\tilde{\mathscr{F}}\|^2 + \|\tilde{\mathscr{G}}\|^2_{x=L}\big\}.
	\end{align}
Also it follows directly from \eqref{D4.17-0}, \eqref{D4.25-0} and \eqref{D4.35} that 
	\begin{align}\label{D4.36}
		\|\tilde{\Theta}_{yy}\|^2 + \|\tilde{\Theta}_{yy}\|^2_{x=L}&\lesssim \|(\tilde{\Theta}_x, \tilde{\Theta}_y)\|^2 + \|\v \tilde{\Theta}_{xx}\|^2 + \|\tilde{\mathscr{F}}\|^2 + \|\tilde{\Theta}_x\|_{x=L}^2  + \|\tilde{\mathscr{G}}\|^2_{x=L}\nonumber\\
		&\lesssim C_{\v} \big\{\|\tilde{\mathscr{F}}\|^2 + \|\tilde{\mathscr{G}}\|^2_{x=L}\big\}.
	\end{align}
Hence we proved \eqref{D4.8} from \eqref{D4.25-0}, \eqref{D4.35}-\eqref{D4.36}. 
 
\smallskip

\noindent{\it Step 3. $H^3$-estimate.}  Noting  $\mathscr{F}(x,0)=0$ and $\mathscr{G}\in H_{00}^{1/2}(\R_+)$, we know that  $\tilde{\mathscr{F}} \in H^1((0,L)\times \R)$ and $\tilde{\mathscr{G}}\in H^{1/2}(\R)$ with 
\begin{align}\label{D.22}
 \|\tilde{\mathscr{F}}\|^2_{H^1}  + \| \tilde{\mathscr{G}}(L,\cdot)\|^2_{H^{1/2}} \lesssim \|\mathscr{F}\|^2_{H^1}  + \| \mathscr{G}(L,\cdot)\|^2_{H^{1/2}_{00}}.
\end{align}

Applying $\pa_y$ to $\eqref{D4.17-0}_1$, one gets
\begin{align}\label{D4.37}
	2 \rho_+ u_+ \tilde{\Theta}_{xy} + 2\rho_+ v_+ \tilde{\chi}_{\v} \tilde{\Theta}_{yy} + 2\rho_+ v_+ \pa_y\tilde{\chi}_{\v} \tilde{\Theta}_{y} - \kappa \Delta_{\v}\tilde{\Theta}_y = \tilde{\mathscr{F}}_y.
\end{align}
Multiplying \eqref{D4.37} by $\tilde{\Theta}_{xxy}$, we have 
\begin{align}\label{D4.39}
	&\rho_+ u_+ \|\tilde{\Theta}_{xy}\|^2_{x=L} - \kappa \big\la \tilde{\Theta}_{xy},\, \tilde{\Theta}_{yyy}\big\ra_{x=L}  - \rho_+ u_+ \|\tilde{\Theta}_{xy}\|^2_{x=0} - \kappa \| ( \sqrt{\v} \tilde{\Theta}_{xxy}, \tilde{\Theta}_{xyy})\|^2 \nonumber\\
	&= O(1) \|\tilde{\Theta}_{xxy}\| \big\{ \|\tilde{\mathscr{F}}_y\| +   \|(\tilde{\Theta}_{yy}, \tilde{\Theta}_{y})\| \big\},
\end{align}
where we have used 
\begin{align*}
\begin{split}
	- \kappa \big\la\tilde{\Theta}_{xxy},\, \tilde{\Theta}_{yyy}\big\ra 
= \kappa \big\la \tilde{\Theta}_{xy},\, \tilde{\Theta}_{xyyy}\big\ra - \kappa \big\la \tilde{\Theta}_{xy},\, \tilde{\Theta}_{yyy}\big\ra_{x=L}  &=- \kappa \|\tilde{\Theta}_{xyy}\|^2  - \kappa \big\la\tilde{\Theta}_{xy},\, \tilde{\Theta}_{yyy}\big\ra_{x=L},\\
2\rho_+ v_+\big\la \tilde{\Theta}_{xxy},\, \tilde{\chi}_{\v} \tilde{\Theta}_{yy} +  \pa_y\tilde{\chi}_{\v} \tilde{\Theta}_{y} \big\ra &\lesssim \|\tilde{\Theta}_{xxy}\|\cdot \|(\tilde{\Theta}_{yy}, \tilde{\Theta}_y)\|.
\end{split}
\end{align*}

For the boundary terms of \eqref{D4.39} at $x=L$, we use $\eqref{D4.17-0}_2$ to obtain
\begin{align}\label{D4.42}
	&\rho_+ u_+ \|\tilde{\Theta}_{xy}\|^2_{x=L} - \kappa \big\la\tilde{\Theta}_{xy},\, \tilde{\Theta}_{yyy}\big\ra_{x=L} 
	=-\rho_+ u_+ \|\tilde{\Theta}_{xy}\|^2_{x=L} +  \big\la \tilde{\Theta}_{xy},\,  \tilde{\mathscr{G}}_y \big\ra_{x=L} \nonumber\\
	&\geq - \f34\rho_+ u_+ \|\tilde{\Theta}_{xy}\|^2_{x=L}  - C \|\tilde{\Theta}_{xy}(L,\cdot)\|_{H^{1/2}}\cdot \| \tilde{\mathscr{G}}(L,\cdot)\|_{H^{1/2}}\nonumber\\
	&\geq - \f34\rho_+ u_+ \|\tilde{\Theta}_{xy}\|^2_{x=L} - \f18 \kappa \|( \sqrt{\v} \tilde{\Theta}_{xxy}, \tilde{\Theta}_{xyy})\|^2  - C_{L} \|\tilde{\Theta}_{xy}\|^2 - C_{L}\| \tilde{\mathscr{G}}(L,\cdot)\|_{H^{1/2}},
\end{align}
where we have used \eqref{Z.11-1} to derive the following estimate
\begin{align*}
	\|\tilde{\Theta}_{xy}(L,\cdot)\|^2_{H^{1/2}}&\lesssim \|\tilde{\Theta}_{xyy}\| \big\{\f{1}{L} \|\tilde{\Theta}_{xy}\| + \|\tilde{\Theta}_{xxy}\|\big\} \lesssim \|(\tilde{\Theta}_{xyy}, \tilde{\Theta}_{xxy})\|^2 + \f{1}{L^2} \|\tilde{\Theta}_{xy}\|^2.
\end{align*}

Substituting \eqref{D4.42} into \eqref{D4.39}, and using \eqref{D4.35}-\eqref{D4.36},  one has that 
\begin{align}\label{D4.45}
\|(\tilde{\Theta}_{xxy}, \tilde{\Theta}_{xyy})\|^2 +\|\tilde{\Theta}_{xy}\|^2_{x=0,L}   
&\lesssim C_{\v}\big\{ \|\tilde{\mathscr{F}}_y\|^2  + \| \tilde{\mathscr{G}}(L,\cdot)\|^2_{H^{1/2}}\big\}.
\end{align}
Also we have from  \eqref{D4.17-0}, \eqref{D4.35}-\eqref{D4.36} and \eqref{D4.45} that
\begin{align}\label{D4.46}
\|(\tilde{\Theta}_{yyy}, \tilde{\Theta}_{xxx})\|^2\lesssim C_{\v}\big\{ \|\tilde{\mathscr{F}}\|^2_{H^1}  + \| \tilde{\mathscr{G}}(L,\cdot)\|^2_{H^{1/2}}\big\}.
\end{align}
Thus we conclude  from \eqref{D4.45}-\eqref{D4.46} and \eqref{D.22}. Therefore the proof of Lemma \ref{lemD.2} is completed. $\hfill\Box$


\medskip

Denote $\tilde{\rho}_s, \tilde{u}_s$ as the even extension of $\rho_s, u_s$ in $y$, and $\tilde{v}_s$ as the odd extension of $v_s$ in $y$. We consider the following BVP
\begin{align}\label{D4.21-0}
	\begin{cases}
		2 \tilde{\rho}_s \tilde{u}_s \tilde{\Theta}_x + 2\tilde{\rho}_s \tilde{v}_s \tilde{\Theta}_y - \kappa \Delta_{\v}\tilde{\Theta} = \tilde{\mathscr{F}},\quad (x,y) \in (0,L)\times \mathbb{R},\\
		\dis  
		\tilde{\Theta}|_{x=0}=0,\quad 
		\big[2\tilde{\rho}_s \tilde{u}_s \tilde{\Theta}_x-\kappa \tilde{\Theta}_{yy}\big]\big|_{x=L}
		=\tilde{\mathscr{G}}\big|_{x=L} \quad \&\quad  \lim_{y\to \infty}\tilde{\Theta}=0.
	\end{cases}
\end{align}

\begin{lemma}\label{lemD.3}
Assume $\tilde{\mathscr{F}}\in H^{1/2}(\R)$ and $\tilde{\mathscr{G}}\in H^1((0,L)\times \mathbb{R})$, then there is a unique strong solution $\tilde{\Theta} \in H^3$ to \eqref{D4.21-0}.
\end{lemma}

\noindent{\bf Proof.} We divide the proof into two steps.

\noindent{\it Step 1. Uniqueness of \eqref{D4.21-0}.} To prove the uniqueness for solutions of \eqref{D4.21}, we need only to consider the following homogeneous problem
\begin{align}\label{D4.47}
	\begin{cases}
	2 \tilde{\rho}_s \tilde{u}_s \tilde{\Theta}_x + 2\tilde{\rho}_s \tilde{v}_s \tilde{\Theta}_y - \kappa \Delta_{\v}\tilde{\Theta} = 0,\quad (x,y) \in (0,L)\times \mathbb{R},\\
	\dis  
	\tilde{\Theta}|_{x=0}=0,\quad 
	\big[2\tilde{\rho}_s \tilde{u}_s\tilde{\Theta}_x-\kappa \tilde{\Theta}_{yy}\big]\big|_{x=L}
	=0 \quad \&\quad  \lim_{y\to \infty}\tilde{\Theta}=0.
\end{cases}
\end{align}
Multiplying \eqref{D4.47} by $\tilde{\Theta}_{xx}$, one obtains
\begin{align}\label{D4.48}
	\kappa \|(\sqrt{\v}\tilde{\Theta}_{xx}, \tilde{\Theta}_{xy})\|^2 + \|\sqrt{\tilde{\rho}_s \tilde{u}_s} \tilde{\Theta}_x\|^2_{x=0} 
	&\leq  \|\sqrt{\tilde{\rho}_s \tilde{u}_s} \tilde{\Theta}_x\|^2_{x=L} - \kappa \big\la \tilde{\Theta}_{x} ,\, \tilde{\Theta}_{yy}\big\ra_{x=L}   + L\|\sqrt{\tilde{u}_s}\tilde{\Theta}_x\|^2_{x=L}   \nonumber\\
	&\quad + \xi \|\tilde{\Theta}_{xy}\|^2  + C\|\tilde{u}_s \tilde{\Theta}_x\|\cdot \|\tilde{\Theta}_x\| + C_{\xi}\|\tilde{u}_s \tilde{\Theta}_x\|^2,
\end{align}
Using the boundary condition in $\eqref{D4.47}$, one has
\begin{align*}
	\|\sqrt{\tilde{\rho}_s \tilde{u}_s} \tilde{\Theta}_x\|^2_{x=L} - \kappa \big\la \tilde{\Theta}_{x} ,\, \tilde{\Theta}_{yy}\big\ra_{x=L}
	&= \big\la \tilde{\Theta}_{x} ,\, \tilde{\rho}_s \tilde{u}_s \tilde{\Theta}_{x} - \kappa \tilde{\Theta}_{yy}\big\ra_{x=L}
	=-\|\sqrt{\tilde{\rho}_s \tilde{u}_s} \tilde{\Theta}_x\|^2_{x=L},
\end{align*}
which, together  with \eqref{D4.48}, yields that 
\begin{align}\label{D4.50}
\kappa \|(\sqrt{\v}\tilde{\Theta}_{xx}, \tilde{\Theta}_{xy})\|^2 + \|\sqrt{\tilde{u}_s} \tilde{\Theta}_x\|^2_{x=0,L} 
&\lesssim \f{1}{a^2}\|\tilde{u}_s \tilde{\Theta}_x\|^2.
\end{align}

It follows from \eqref{D4.47} that 
\begin{align}\label{D4.51}
	\|\big(2\tilde{\rho}_s \tilde{u}_s \tilde{\Theta}_x - \kappa \Delta_{\v} \tilde{\Theta}\big)\|^2 
	&\lesssim   L \|\tilde{\Theta}_{xy}\|^2.
\end{align}
It is clear that 
\begin{align}\label{D4.52}
	\|\big(2\tilde{\rho}_s \tilde{u}_s \tilde{\Theta}_x - \kappa \Delta_{\v}\tilde{\Theta}\big)\|^2
	&=4\|\tilde{\rho}_s \tilde{u}_s \tilde{\Theta}_{x}\|^2 + \kappa^2 \|(\v \tilde{\Theta}_{xx}, \tilde{\Theta}_{yy})\|^2 + 2\kappa^2 \big\la \v \tilde{\Theta}_{xx},\, \tilde{\Theta}_{yy}\big\ra \nonumber\\
	&\quad  - 4 \kappa \big\la \tilde{\rho}_s \tilde{u}_s \tilde{\Theta}_{x},\, \v \tilde{\Theta}_{xx} \big\ra   - 4\kappa  \big\la \tilde{\rho}_s \tilde{u}_s \tilde{\Theta}_{x},\, \tilde{\Theta}_{yy}\big\ra.
\end{align}
A direct calculation shows that 
\begin{align}\label{D4.53}
	\begin{split}
		2\kappa^2 \big\la \v \tilde{\Theta}_{xx},\, \tilde{\Theta}_{yy} \big\ra&
		\lesssim \sqrt{\v} \|\tilde{\Theta}_{yy}\|^2 + \sqrt{\v}\|\sqrt{\v}\tilde{\Theta}_{xx}\|^2,\\
		-4 \big\la \tilde{\rho}_s \tilde{u}_s \tilde{\Theta}_{x},\, \v \tilde{\Theta}_{xx} w^2\big\ra
		&\lesssim \sqrt{\v} \|\tilde{u}_s \tilde{\Theta}_x\|^2  + \sqrt{\v}\|\sqrt{\v}\tilde{\Theta}_{xx}\|^2,
	\end{split}
\end{align}
and
\begin{align}\label{D4.55}
	- 4 \big\la \tilde{\rho}_s \tilde{u}_s \tilde{\Theta}_{x},\, \tilde{\Theta}_{yy}\big\ra 
	&=  4 \big\la \tilde{\rho}_s \tilde{u}_s \tilde{\Theta}_{xy},\, \tilde{\Theta}_{y}\big\ra + 4 \big\la (\tilde{\rho}_s \tilde{u}_s)_y \tilde{\Theta}_{x},\, \tilde{\Theta}_{y}\big\ra \nonumber \\
	&=2 \|\sqrt{\tilde{\rho}_s \tilde{u}_s} \tilde{\Theta}_{y}\|^2_{x=L} + L^{1-} \|u_s \tilde{\Theta}_x\|^2 + L^{-}\|\tilde{\Theta}_{xy}\|^2.
\end{align}

It follows from \eqref{D4.51}-\eqref{D4.55} that 
\begin{align} \label{D4.56}
\|\tilde{u}_s \tilde{\Theta}_{x}\|^2 + \|\tilde{\Theta}_{yy}\|^2 +  \|\sqrt{\tilde{u}_s} \tilde{\Theta}_{y}\|^2_{x=L} 
&\lesssim L^{-}\|(\sqrt{\v}\tilde{\Theta}_{xx}, \tilde{\Theta}_{xy})\|^2 .
\end{align}
Now combining \eqref{D4.50} and \eqref{D4.56}, and noting $L\ll 1$ (independent of $\v$), we obtain 
\begin{align*}
	\kappa \|(\sqrt{\v}\tilde{\Theta}_{xx}, \tilde{\Theta}_{xy}, \tilde{\Theta}_{yy})\|^2  + 	\|\tilde{u}_s \tilde{\Theta}_{x}\|^2 + \|\sqrt{\tilde{u}_s}  \tilde{\Theta}_x\|^2_{x=0,L}   +  \|\sqrt{\tilde{u}_s} \tilde{\Theta}_{y}\|^2_{x=L} =0,
\end{align*}
which yields that $\tilde{\Theta}\equiv 0.$ Thus the solution of \eqref{D4.47} is unique in $H^2$.

\smallskip

\noindent{\it Step 2. Fredholm argument.} 
We denote
\begin{align*}
	\begin{split}
		\mathcal{L}_1:&=2 \rho_+ u_+ \pa_x + 2\rho_+ v_+ \tilde{\chi}_{\v} \pa_y - \kappa \Delta_{\v} ,\\
		\mathcal{L}_2:&=2\rho_+ u_+ \pa_x - \kappa \pa_{yy}.
	\end{split}
\end{align*}
then we can rewrite \eqref{D4.17-0} as
\begin{align}\label{D4.58}
	\begin{cases}
		\mathcal{L}_1 \tilde{\Theta} = \tilde{\mathscr{F}}, \\
		\mathcal{L}_2 \tilde{\Theta}\big|_{x=L}
		= \tilde{\mathscr{G}}|_{x=L} 
		\quad \&\quad 
		\tilde{\Theta}|_{x=0}=0.
	\end{cases}
\end{align}
For later use, we shall denote the solution of \eqref{D4.58} as $\tilde{\Theta}:=\mathcal{L}^{-1}[\tilde{\mathscr{F}}, \tilde{\mathscr{G}}]$.


Let $\tilde{\Theta}$ be the solution of \eqref{D4.21-0} and denote
\begin{align}
	\tilde{\mathbf{F}}:=\mathcal{L}_1\tilde{\Theta}\quad \&\quad  \tilde{\mathbf{G}}|_{x=L}:=\mathcal{L}_2\tilde{\Theta} |_{x=L} \,\, \Longleftrightarrow\,\, \tilde{\Theta}=\mathcal{L}^{-1}[\tilde{\mathbf{F}}, \tilde{\mathbf{G}}], 
\end{align}
then we can represent \eqref{D4.21-0} as
\begin{align}
	\begin{cases}
		\tilde{\mathbf{F}} = \tilde{\mathscr{F}} -2 [\rho_s u_s-\rho_+u_+] \big(\mathcal{L}^{-1}[\tilde{\mathbf{F}}, \tilde{\mathbf{G}}]\big)_x -2 [\rho_s v_s-\rho_+v_+ \tilde{\chi}_{\v}] \big(\mathcal{L}^{-1}[\tilde{\mathbf{F}}, \tilde{\mathbf{G}}]\big)_y,\\
		\tilde{\mathbf{G}}\big|_{x=L}= - 2[\rho_s u_s- \rho_+ u_+] \big(\mathcal{L}^{-1}[\tilde{\mathbf{F}}, \tilde{\mathbf{G}}]\big)_x\big|_{x=L} + \tilde{\mathscr{G}}\big|_{x=L}
	\end{cases}
\end{align}
Noting \eqref{D4.25-0}, \eqref{D4.35}-\eqref{D4.36} and \eqref{D4.45}-\eqref{D4.46}, we know that
\begin{align}
	\big[\tilde{\mathbf{F}},\, \tilde{\mathbf{G}}\big] \in H^1([0,L]\times \R) \times H^{1/2}(\R)\,\, \mapsto \,\, H^3 ([0,L]\times \R).
\end{align}
Then it is clear that the following operator 
\begin{align}
	\begin{cases}
		-2 [\rho_s u_s-\rho_+u_+] \big(\mathcal{L}^{-1}[\tilde{\mathbf{F}}, \tilde{\mathbf{G}}]\big)_x -2 [\rho_s v_s-\rho_+v_+ \tilde{\chi}_{\v}] \big(\mathcal{L}^{-1}[\tilde{\mathbf{F}}, \tilde{\mathbf{G}}]\big)_y,\\
		- 2[\rho_s u_s- \rho_+ u_+] \big(\mathcal{L}^{-1}[\tilde{\mathbf{F}}, \tilde{\mathbf{G}}]\big)_x\big|_{x=L},
	\end{cases}
\end{align}
is compact in $H^1([0,L]\times \R) \times H^{1/2}(\R)$.

Finally, noting the uniqueness of homogeneous problem \eqref{D4.47}, then one can establish the existence and uniqueness of solution $\tilde{\Theta} \in H^3$ of \eqref{D4.21-0}  by applying the Fredholm alternative. Therefore the proof of Lemma \ref{lemD.3} is completed. $\hfill\Box$

\smallskip

\noindent{\bf Proof of Proposition \ref{propC.1}:} Noting \eqref{D4.3-0}-\eqref{D4.21}, we need only to establish the existence for \eqref{D4.21}. We consider the following approximation problem:
\begin{align}\label{D4.69}
	\begin{cases}
		2 \rho_s u_s \tilde{S}^{\d}_x + 2\rho_s v_s \tilde{S}^{\d}_y - \kappa \Delta_{\v}\tilde{S}^{\d} = \tilde{f} \bar{\chi}(\f{y}{\delta}),\quad (x,y) \in (0,L)\times \mathbb{R}_+,\\
		\dis \tilde{S}^{\d}|_{y=0}=0 \quad \&\quad 
		\tilde{S}^{\d}|_{x=0}=0,\\
		\big[2\rho_s u_s\tilde{S}^{\d}_x-\kappa \tilde{S}_{yy}^{\d}\big]\big|_{x=L}
		=\bar{g}\bar{\chi}(\f{y}{\delta})\big|_{x=L}.
	\end{cases}
\end{align}

Noting $\tilde{S}^{\d}|_{y=0}=0$, we apply an odd extension for $\tilde{S}^{\d}$,  we can reduce \eqref{D4.69} to be  an equation like \eqref{D4.21-0}, then we apply Lemma \ref{lemD.3} to obtain the existence of solution $\tilde{S}^{\d}\in H^3([0,L]\times \R_+)$ of \eqref{D4.69}. We remark the the cut-off in \eqref{D4.69} is necessary since  $\tilde{f}(x,0)\neq 0$ and $\bar{g}(L,0)\neq0$ in general.

Finally, for \eqref{D4.69}, applying the {\it a priori} uniform-in-$\delta$ estimates for $\tilde{S}^{\d}$ established in \eqref{D3.30-1},  \eqref{D3.67-0} and \eqref{D3.66}, we can take the limit $\delta\to0$ so that the original problem \eqref{D4.21} is solved. Therefore  the proof Proposition \ref{propC.1}  is completed. $\hfill\Box$

\medskip

\noindent{\bf Acknowledgments.} 
Yong Wang’s research  is partially supported by 
the National Natural
Science Foundation of China grants No. 12421001 \& 12288201, and CAS Project for Young
Scientists in Basic Research, grant No. YSBR-031. Yan Guo thanks Ian Tice for ref. on $H_{00}^{1/2}$. 


\end{document}